\documentclass[11pt]{article}
\usepackage[a4paper,margin=25mm]{geometry}
\usepackage[T1]{fontenc}
\usepackage{lmodern,amsmath,amssymb,amsthm,mathtools,mathrsfs,microtype}
\usepackage{booktabs,tabularx,float}
\usepackage{tikz}
\usetikzlibrary{shapes.geometric}
\usepackage[colorlinks=true,linkcolor=blue,citecolor=blue,urlcolor=blue]{hyperref}
\numberwithin{equation}{section}
\newtheorem{theorem}{Theorem}[section]
\newtheorem{definition}[theorem]{Definition}
\newtheorem{lemma}[theorem]{Lemma}
\newtheorem{proposition}[theorem]{Proposition}
\newtheorem{corollary}[theorem]{Corollary}
\theoremstyle{remark}
\newtheorem{remark}[theorem]{Remark}
\newcommand{\E}{\mathbb E}
\newcommand{\Prob}{\mathbb P}
\newcommand{\Z}{\mathbb Z}
\newcommand{\C}{\mathbb C}
\newcommand{\R}{\mathbb R}
\newcommand{\Tr}{\operatorname{Tr}}
\newcommand{\wt}{\widetilde}
\newcommand{\diag}{\operatorname{diag}}
\newcommand{\supp}{\operatorname{supp}}
\newcommand{\ellscale}{\ell_{W,g}}
\newcommand{\ellstar}{\ell_{W,g,L}}

\newcommand{\dist}{\operatorname{dist}}
\newcommand{\diam}{\operatorname{diam}}
\newcommand{\cL}{\mathcal L}
\newcommand{\cK}{\mathcal K}
\newcommand{\cQ}{\mathcal Q}
\newcommand{\op}{\mathrm{op}}
\newcommand{\rank}{\operatorname{rank}}
\newcommand{\ran}{\operatorname{ran}}
\newcommand{\HS}{\mathrm{HS}}
\newcommand{\one}{\mathbf 1}
\newcommand{\norm}[1]{\left\|#1\right\|}
\title{The Localization Scale of the Block Anderson Model\\
in One Dimension at Arbitrary Coupling and Volume}
\author{Jun Yin\thanks{The author acknowledges support from the National Science Foundation
under Grant No.~\href{https://www.nsf.gov/awardsearch/show-award?AWD_ID=2550762}{2550762}
and from MIT--Anthropic for the use of AI tools in this work.}\\
\small Department of Mathematics, University of California, Los Angeles\\
\small\href{mailto:jyin@math.ucla.edu}{\texttt{jyin@math.ucla.edu}}}
\date{}
\hypersetup{
  pdftitle={The Localization Scale of the Block Anderson Model in One Dimension at Arbitrary Coupling and Volume},
  pdfauthor={Jun Yin}
}
\begin{document}
\maketitle
\begin{abstract}
We study the one-dimensional block Anderson model with independent
\(W\times W\) Gaussian unitary ensemble potentials and deterministic
nearest-neighbor hopping, with hopping-to-disorder ratio \(g\).
The system consists of \(L\) blocks,
with total dimension \(N=WL\). As \(W\to\infty\), we identify the bulk localization
length, up to an arbitrarily small power of \(W\), with
\(\min\{N,\max\{W,W^2g^2\}\}\), uniformly for deterministic
couplings $g\ge0$ and cycle lengths $L\ge3$. For every fixed $D>0$ and
all sufficiently large $W$, the expected fraction of eigenvectors in the
selected energy interval satisfying these bounds is at least $1-W^{-D}$.
In polynomial volume the bounds hold simultaneously for all eigenvectors
in that interval with probability at least $1-W^{-D}$. The localization
length is the number of coordinates in the shortest cyclic interval
containing at least half of an eigenvector's mass.
\end{abstract}

\tableofcontents
\clearpage

\section{Introduction}
\label{sec:components}

The localization length describes the region carrying a substantial
fraction of an eigenvector's mass. Its dependence on disorder and on the
number of interacting degrees of freedom is a central question in the
theory of random Schr\"odinger operators. Anderson \cite{Anderson}
introduced a random lattice Hamiltonian to explain the suppression of
diffusion by disorder. Two fundamental
approaches to rigorous localization are the multiscale analysis of
Fr\"ohlich--Spencer \cite{FS} and the fractional-moment method of
Aizenman--Molchanov \cite{AM}. They establish localization at strong
disorder or suitable spectral edges. Quantifying the localization length
as the disorder weakens, or as the number of channels grows, requires
additional estimates beyond the existence of localization.

Random band matrices provide a finite-dimensional setting for this
question. Their independent entries are appreciable only within a band
of width \(W\), so the bandwidth controls the range of interaction.
In one dimension, the predicted localization scale is \(W^2\), and the
corresponding transition for a matrix of size \(N\) occurs near
\(W=\sqrt N\); the supersymmetric analysis of Fyodorov--Mirlin
\cite{FM} gives an early theoretical explanation of this scale.
On the localization side, for each fixed small \(\varepsilon>0\),
Schenker \cite{Schenker} proved localization for Gaussian band matrices
when \(W\le N^{1/8-\varepsilon}\). Chen--Smart
\cite{CS} and Cipolloni--Peled--Schenker--Shapiro \cite{CPSS} reached
\(W\le N^{1/4-\varepsilon}\), by different refinements of fluctuation arguments.
More recently, Drogin \cite{DroginRBM} proved an
eigenfunction-correlator bound with exponential decay length \(W^2\)
and a polynomial prefactor in \(W\), on any fixed compact energy interval,
for a class of real band matrices with independent entries, including
the real Wegner orbital model. On the delocalization side, Yau--Yin
\cite{YY} proved bulk delocalization and
universality for Gaussian band matrices when
\(W\ge N^{1/2+c}\), for any fixed \(c>0\).
These results concern ensembles with random off-diagonal entries;
the deterministic hopping in the present paper requires a different
localization argument.

The dimension also changes the picture substantially.
S.~Dubova, K.~Yang, H.-T.~Yau and J.~Yin \cite{DYY2D} proved bulk delocalization and
universality for two-dimensional Gaussian band matrices with
\(W\ge N^c\), for any fixed \(c>0\).
S.~Dubova, F.~Yang, H.-T.~Yau and J.~Yin \cite{DYYY}
established corresponding results in dimensions \(d\ge3\),
and also treated block Anderson operators in a range of couplings.
The stochastic analysis of resolvent loops developed in these works
connects the quantum problem to a deterministic diffusion profile.
The spatial estimates and coupling ranges depend on the dimension;
the high-dimensional delocalization theorems do not themselves give
the one-dimensional estimates considered here.

Block models retain a closer connection to the Anderson operator:
each lattice site carries several internal coordinates, while hopping
connects neighboring sites. This viewpoint goes back to Wegner's
orbital model \cite{Wegner}. In the block Anderson model the onsite (i.e.\ within the same block)
potential is random and the hopping is deterministic, whereas in the
Wegner orbital model the hopping blocks are random as well.
For normalized Gaussian blocks of size \(W\) and identity hopping
of strength \(g\), Peled--Schenker--Shamis--Sodin
\cite{PSSS} proved exponential decay of fractional moments
of the Green function for \(g\le cW^{-1/2}\) with \(c\) sufficiently
small. This gives the upper bound on the scale of a single block in
this regime, up to logarithmic factors. For small coupling above this
threshold, Truong--Yang--Yin \cite{TYY} proved the bulk
lower bound on the scale \(\min\{N,(Wg)^2\}\), up to arbitrarily small
power losses, in polynomial finite volume for
\(W^{-1/2+\epsilon}\le g\le W^{-\epsilon}\) with fixed small \(\epsilon>0\).
Their work also gives delocalization in two dimensions.

At order-one coupling, Drogin \cite{Drogin}
proved the \(W^2\) upper scale in coordinate distance for \(g=1\) and
GOE onsite blocks. His bound concerns the limiting exponential decay
rate of the Green function at each fixed energy and follows from a
lower bound of order \(W^{-1}\) on the gaps between Lyapunov exponents.

Here we study a cycle of \(L\) blocks with independent normalized
complex GUE potentials of size \(W\), coupled by \(gI_W\).
In coordinate distance, our main result establishes matching lower
and upper bounds on the scale $\min\{N,\max\{W,W^2g^2\}\}$, where
$N=WL$, up to arbitrarily small powers of $W$ and uniformly over
deterministic finite $L$ and $g$.

For the lower bound we follow the loop hierarchy of Yau--Yin
\cite{YY}, adapted to block Schr\"odinger operators by
Truong--Yang--Yin \cite{TYY}. It propagates estimates for traces of
cyclic products of Green functions to the spectral scale needed to bound
the mass of an eigenvector on a single block. Our main addition is
at large \(g\), where transfer matrices give the initial estimates
for this propagation and, for arbitrarily large \(g\), reduce the
spatial estimates to \(W^{O(1)}\) sample points.

For the upper bound we use the fractional-moment method of
Aizenman--Molchanov \cite{AM,ASFH}. Beyond the range covered by
\cite{PSSS}, we follow the boundary comparison of Macera--Sodin
\cite{MS} to deduce Green-function decay from the weakest
expansion of a transfer product on the subspace allowed by the boundary
condition. At small \(g\), we use Le Page's weighted-operator approach
to random matrix products \cite{LePage} and
adapt Drogin's quantitative nonconcentration argument
\cite{Drogin} to prove contraction. Our perturbation
of the onsite blocks varies the expansion in one direction by matching
the variation on the few directions that a short product strongly
amplifies. At large \(g\), we estimate the weakest
expansion directly by perturbing around the free transfer. A correction
to the logarithmic growth reveals a positive drift of order
\((Wg^2)^{-1}\) per block. This argument covers the band center and,
for \(g\ge W\), the entire bulk.

\subsection{The model}
\label{sec:model}

Let \(W,L\) be positive integers, and put \(N\coloneqq WL\).
We identify \(\Z_N\) and \(\Z_L\) with integers modulo \(N\) and
\(L\), respectively. Divide \(\Z_N\) into the blocks
\begin{equation}\label{eq:blocks}
 \mathcal I_a\coloneqq \{(a-1)W+1,\ldots,aW\},\qquad a\in\Z_L,
\end{equation}
where representatives \(a=1,\ldots,L\) are used in this formula.
For a coordinate \(i\in\Z_N\), write \([i]\) for the block containing
\(i\). Indices \(a,b,c,\ldots\) denote blocks, and
\(i,j,\ldots\) denote coordinates. We also identify
\(\C^N=\C^L\otimes\C^W\) in the order given by \eqref{eq:blocks}.

For each \(a\in\Z_L\), let \(V_a\) be a normalized complex GUE
matrix of size \(W\). The matrices \((V_a)_{a\in\Z_L}\) are
independent, and their covariance is
\begin{equation}\label{eq:covariance}
 \E(V_a)_{\alpha\beta}(V_b)_{\gamma\delta}
 =W^{-1}\delta_{ab}\delta_{\alpha\delta}\delta_{\beta\gamma},
 \qquad 1\le\alpha,\beta,\gamma,\delta\le W.
\end{equation}
Thus the diagonal entries are real Gaussians of variance \(W^{-1}\);
the real and imaginary parts of each upper off-diagonal entry are
independent Gaussians of variance \((2W)^{-1}\).

Let \(A_L\) be the adjacency matrix of the cycle \(\Z_L\) of \(L\) blocks, so
\((A_Lf)(a)=f(a-1)+f(a+1)\). For a coupling \(g\), define
\begin{equation}\label{eq:model}
 H=H_{g,W,L}\coloneqq \frac{gA_L\otimes I_W+V}{\max\{1,g\}},
 \qquad V\coloneqq \diag(V_a:a\in\Z_L).
\end{equation}
Equivalently, for a vector \(\psi=(\psi_a)_{a\in\Z_L}\),
\begin{equation}\label{eq:blockaction}
 (H\psi)_a=\frac{V_a\psi_a+g\psi_{a-1}+g\psi_{a+1}}{\max\{1,g\}},
 \qquad \psi_a\in\C^W.
\end{equation}
The parameter \(g\) is the ratio of the hopping strength to the onsite
disorder strength. The factor \(\max\{1,g\}\) is included only for later
convenience; it leaves the eigenvectors and their localization lengths
unchanged. For
\(g\le1\), the model is \(gA_L\otimes I_W+V\); for \(g\ge1\), it is
\(A_L\otimes I_W+V/g\). We use this normalization throughout.
We write
\begin{equation}\label{eq:eigenvectors}
 \lambda_1\le\cdots\le\lambda_N,\qquad
 H\psi_k=\lambda_k\psi_k,\qquad \|\psi_k\|_2=1,
\end{equation}
for the eigenvalues and normalized eigenvectors. The spectrum is
almost surely simple: the characteristic-polynomial discriminant is a
polynomial in the real Gaussian coordinates, is nonzero for sufficiently
separated diagonal onsite entries, and hence vanishes with probability zero.

Throughout the paper, $g\ge0$ and $L\ge3$ are finite deterministic
parameters. We consider $W\to\infty$, allowing $g=g(W)$ and $L=L(W)$
to vary without further restriction. Probability estimates hold uniformly
over the stated deterministic parameter range, with a separate event for
each choice of $g,L$.

\paragraph{Notation.}

We write \(\C_+\coloneqq \{z\in\C:\Im z>0\}\) and
\(G(z)\coloneqq (H-z)^{-1}\). For \(\psi\in\C^N\),
\[
 \psi_a\coloneqq (\psi(i))_{i\in\mathcal I_a}.
\]
Matrix operator and entrywise
maximum norms are denoted by \(\|A\|\) and
\(\|A\|_{\max}\coloneqq \max_{i,j}|A_{ij}|\), respectively.
We use \(\langle A\rangle\coloneqq \Tr A\) for the unnormalized trace.
For an $N\times N$ block matrix $A$, define its normalized block traces by
\[
 \langle A\rangle_{ab}\coloneqq W^{-1}\Tr A_{ab}.
\]
The trace identifies the orbital indices in blocks $a$ and $b$.
For block observables, we write
\[
 P_a\coloneqq \diag(\mathbf1(i\in\mathcal I_a)),\qquad E_a\coloneqq W^{-1}P_a
\]
for the orthogonal projection onto block $a$ and its normalized version.

Constants \(C,c>0\) may change from line to line and may depend on
the fixed parameters in a statement, but not on \(W,L,g\).
The notation \(X\lesssim Y\) means \(X\le CY\), and
\(X\asymp Y\) means both \(X\lesssim Y\) and \(Y\lesssim X\).

\medskip\noindent\textbf{Localization length.}

A cyclic interval of length \(r\) is a set
\(I=\{x,x+1,\ldots,x+r-1\}\subseteq\Z_N\), with
\(1\le r\le N\), where addition is modulo \(N\).
\begin{definition}[Localization length]\label{def:length}
For a normalized vector \(\psi\in\C^N\), define
\begin{equation}\label{eq:length}
 \ell(\psi)\coloneqq \min\left\{|I|:\ I\subseteq\Z_N\text{ is a cyclic interval},
                  \quad\sum_{i\in I}|\psi(i)|^2\ge\frac12\right\}.
\end{equation}
\end{definition}
In particular, \(1\le\ell(\psi)\le N\). The definition selects an
interval containing half of the mass; the interval and its location
may depend on \(\psi\). Its length counts coordinates, rather than
blocks. It measures finite-volume mass concentration and does not
define an asymptotic exponential decay rate.

Introduce the deterministic coordinate scales
\begin{equation}\label{eq:scales}
 \ellscale\coloneqq\max\{W,W^2g^2\},\qquad
 \ellstar\coloneqq\min\{WL,\ellscale\}.
\end{equation}
The uncapped scale is $\ellscale$, and the finite-volume scale is $\ellstar$.
The corresponding block scales are $\ellscale/W$ and $\ellstar/W$.
Equivalently,
\begin{equation}\label{eq:three-regimes}
 \ellstar=
 \begin{cases}
 W, & Wg^2\le1,\\
 W^2g^2, & 1<Wg^2<L,\\
 N, & Wg^2\ge L.
 \end{cases}
\end{equation}
The first regime is the single-block scale. In the second regime,
the characteristic number of occupied blocks is \(Wg^2\). In the
third regime, the finite-volume scale is reached.

\medskip\noindent\textbf{The deterministic bulk.}

The deterministic band edges depend on $g$. We therefore define the
bulk using the deterministic spectral distribution of the normalized
model \eqref{eq:model}.

Let \(m_g:\C_+\to\C_+\) be the solution of the infinite-chain
Dyson equation corresponding to \eqref{eq:model}:
\begin{equation}\label{eq:dyson}
 m_g(z)=(2\pi)^{-1}\int_{-\pi}^{\pi}
 \left(\dfrac{2g\cos p}{\max\{1,g\}}-z
                  -\dfrac{m_g(z)}{\max\{1,g\}^2}\right)^{-1}\mathrm dp,
 \qquad m_g(z)\sim-z^{-1}\quad(|z|\to\infty).
\end{equation}
This is \cite[(2.6)]{TYY} with the hopping and onsite (i.e.\ within the same block) covariance of
\eqref{eq:model}. Let $\nu_g$ be the probability measure with
Stieltjes transform $m_g$, and define its distribution function by
\begin{equation}\label{eq:distribution}
 F_g(x)\coloneqq \nu_g(({-\infty},x])
       =\pi^{-1}\int_{-\infty}^x\Im m_g(E+i0)\,\mathrm dE.
\end{equation}
Here $m_g(E+i0)$ denotes the boundary value from the upper half-plane.
For every finite $g$, the measure $\nu_g$ is even and has a continuous
density, positive inside a single spectral band.

As $g\to\infty$, $\nu_g$ tends to the arcsine law on $[-2,2]$,
with density
\begin{equation}\label{eq:arcsine}
 \rho_\infty(x)\coloneqq \frac1{\pi\sqrt{4-x^2}}\quad (|x|<2),
 \qquad \rho_\infty(x)\coloneqq 0\quad (|x|\ge2).
\end{equation}
\begin{lemma}[Deterministic spectral distribution and bulk windows]
\label{da:lem:distribution}
Let $\nu_g$ and $F_g$ be the measure and distribution function defined above.
\begin{enumerate}
\item[(i)] For every $g\ge0$,
\begin{equation}\label{eq:bounded-physical-bulk}
 \supp\nu_g\subseteq[-4,4].
\end{equation}
\item[(ii)] The measures are uniformly atomless:
\begin{equation}\label{rd:eq:noatoms}
 \lim_{h\downarrow0}\sup_{g\ge0,\ x\in\R}\nu_g([x-h,x+h])=0.
\end{equation}
\item[(iii)] For every $\kappa\in(0,1/2)$, there are
$\kappa_0\in(0,2)$ and $0<c_\kappa<C_\kappa<\infty$, depending only on
$\kappa$, such that
\begin{equation}\label{eq:bulk-bridge}
 g\in[0,c_\kappa]\cup[C_\kappa,\infty),\qquad F_g(-|E|)\ge\kappa
 \quad\Longrightarrow\quad |E|\le2-\kappa_0.
\end{equation}
\end{enumerate}
\end{lemma}
The proof, including the preceding distributional facts, is given in
Appendix~\ref{app:bulk-regularity}.

The following rough density law is uniform in both parameters. It will
determine the threshold for the cycle length $L$ in our choice of energy interval.
\begin{lemma}[A uniform rough density law]\label{lem:rough-density}
For every $\delta>0$ there are $W_\delta,L_\delta,C_\delta,c_\delta>0$
such that, for $W\ge W_\delta$, $L\ge L_\delta$, and every $g$,
\begin{equation}\label{eq:rough-counting}
 \Prob\left(\sup_{x\in\R}
 \left|N^{-1}\#\{k:\lambda_k\le x\}-F_g(x)\right|>\delta\right)
 \le C_\delta e^{-c_\delta W^2L}.
\end{equation}
\end{lemma}
Its proof is given in Appendix~\ref{sec:rough-density-proof}.
For $0<\kappa<1/2$, using the threshold $L_{\kappa/2}$
from Lemma~\ref{lem:rough-density}, define
\begin{equation}\label{eq:bulk-sampling}
 J_{g,\kappa,L}\coloneqq
 \begin{cases}
 \R, & L<L_{\kappa/2},\\
 \{E:F_g(-|E|)\ge\kappa/2\}, & L\ge L_{\kappa/2}.
 \end{cases}
\end{equation}
For $L\ge L_{\kappa/2}$, the interval excludes $\nu_g$-mass $\kappa/2$
at each edge, so its total mass is $1-\kappa$. The lemma, with
$\delta=\kappa/2$, shows that it contains at least $(1-2\kappa)WL$
eigenvalues except on an event of probability $C_\kappa e^{-c_\kappa W^2L}$,
for all sufficiently large $W$.  
 
\subsection{Main results}

\begin{theorem}[Localization length]\label{thm:main}
Fix $0<\kappa<1/2$ and $D,\tau>0$. For all sufficiently large $W$,
uniformly over $L,g$,
\begin{equation}\label{eq:main}
 \E\!\left[
 \frac{\#\{k:\lambda_k\in J_{g,\kappa,L},\ 
          W^{-\tau}\ellstar\le\ell(\psi_k)\le W^\tau\ellstar\}}
      {\#\{k:\lambda_k\in J_{g,\kappa,L}\}}
 \right]\ge 1-W^{-D}.
\end{equation}
Here $\ellstar$ is defined in \eqref{eq:scales}, and we set $0/0=0$.
\end{theorem}

It is easy to obtain the next corollary from Theorem~\ref{thm:main}
by Markov's inequality.

\begin{corollary}[Simultaneous bounds in polynomial volume]
\label{cor:simultaneous}
Fix $B_L,D,\tau>0$ and $0<\kappa<1/2$. For all sufficiently large $W$,
uniformly over $L\le W^{B_L}$ and $g$,
\begin{equation}\label{eq:simultaneous}
 \Prob\left\{
 W^{-\tau}\ellstar\le\ell(\psi_k)\le W^\tau\ellstar
 \text{ for every }k\text{ with }\lambda_k\in J_{g,\kappa,L}
 \right\}\ge1-W^{-D}.
\end{equation}
\end{corollary}

The power loss $W^\tau$ may be improvable. The unrestricted-volume
statement controls the expected proportion of eigenvectors satisfying
the bounds; it does not assert that every eigenvector in the selected
interval satisfies them when $L$ is arbitrarily large.

\paragraph{Organization of the paper.}
Sections~\ref{sec:intro-spreading-methods}
and~\ref{sec:intro-concentration-methods} explain the lower- and
upper-bound strategies, respectively, and state the component estimates.
The proofs in Sections~\ref{sec:extension-regions}
and~\ref{sec:upper-regions} assemble these estimates into
Theorem~\ref{thm:main}.
Section~\ref{sec:unrestricted} proves
the reduction from resolvent moments to spreading at arbitrary volume.
Section~\ref{sec:flow-overview} presents the loop hierarchy and the common
propagation theorem, and completes Theorems~\ref{thm:bounded-lower}
and~\ref{uvu:thm:terminal} using the inputs proved later.
Section~\ref{sec:lower-deterministic} collects the deterministic Dyson,
kernel and propagation-scale estimates.
Section~\ref{sec:large-proof} proves the transfer estimates and spatial
sampling, and verifies the stochastic inputs to the common propagation
Theorem~\ref{gc:thm:common}, uniformly for arbitrarily large coupling.

The upper-bound reduction is presented in Section~\ref{sec:upper-method}.
Section~\ref{sec:small-upper-proof} proves the small-coupling vector-field
estimates, Section~\ref{sec:compact-proof} constructs Gaussian fields from
the finite-product differential, and Section~\ref{sec:large-upper-proof}
proves the large-coupling transfer estimates.

Appendix~\ref{app:deterministic} proves the deterministic inputs,
Appendix~\ref{app:loop-identities} proves the loop identities and
primitive-loop estimates, and Appendix~\ref{app:gaussian-continuation}
proves the bounded-coupling initialization and the common propagation theorem.
Appendix~\ref{sec:elementary} contains orbital spreading and the rough
density law used by both bounds.
Appendix~\ref{app:strong-concentration} proves the fractional-moment
conversion and the finite-volume criterion;
Appendices~\ref{app:compact-circle}
and~\ref{app:upper-gue} contain the circle-operator estimates and the
auxiliary GUE estimates.

\subsection{Version note and AI disclosure}
\label{sec:version-ai}

\paragraph{Version note.}
This version gives the full proof of Theorem~\ref{thm:main}. The statement
of the theorem is final, although the exposition may undergo substantial
revision. The author will prepare a separate, shorter account for journal
submission. The Lean formalization proves Theorem~\ref{thm:main}
assuming no mathematical results beyond mathlib. Differences between
its arguments and those in the text are listed in
\texttt{PAPER-VS-LEAN.md} in the repository
(\url{https://github.com/JYin80/Lean-LSBA1d}). The author used
Claude Opus 5.5 and Claude Sonnet 5.5 to write the Lean proofs.

\paragraph{AI disclosure.}
The author wrote the first draft of the main text. In subsequent work with
OpenAI's ChatGPT 6 Astra, the author extended the results to the full range
of couplings and volumes in Theorem~\ref{thm:main}. The model participated
in mathematical discussions, explored proof strategies proposed by the
author, and helped revise the manuscript.

\clearpage

\section{Proof strategy for the lower bound}
\label{sec:intro-spreading-methods}

\subsection{Spatial cutoff and the loop hierarchy}

\paragraph{Spatial cutoff.}
To bound the localization length from below, we control the largest
block mass $\max_a\|\psi_a\|^2$, where $\psi_a$ is the $a$th block
component of $\psi$, through the diagonal blocks of the Green function
$G(z)=(H-z)^{-1}$. At the spectral heights used in the proof,
Theorem~\ref{uv:thm:smooth-cutoff} shows that $\Im G_{aa}$ on a cycle
with $L\ge\ellscale=\max\{W,W^2g^2\}$ is comparable to the same
quantity on a local periodic restriction of at most $C\ellscale$
blocks. Thus it suffices to establish the resolvent estimates for
$L\le C\ellscale$, instead of arbitrary $L$.
Theorem~\ref{uv:lem:lower-volume} combines this cutoff with the
conversion from resolvent moments to spreading.

\paragraph{The basic strategy: loop hierarchy.}
The main tool for the spreading bound is the loop hierarchy developed
in \cite[Sections~2 and~5]{YY} and extended to block Schr\"odinger
operators in \cite[Sections~3--7]{TYY}. We use a flow $(H_t,z_t)$ in which the onsite blocks
are independent Hermitian Brownian motions and the spectral height
decreases along a deterministic characteristic. At $t=1$, the matrix
$H_1$ has the law of the target matrix $H$.

Resolvent loops $\cL_k(t)$ are cyclic traces of products of
$G_t=(H_t-z_t)^{-1}$, with normalized block projections inserted
between the factors. Along the flow
\eqref{fo:eq:flow}--\eqref{fo:eq:characteristic}, It\^o's formula
expresses their drifts and quadratic variations through other
loops, giving the \emph{loop hierarchy}. The deterministic
primitive loops $\cK_k(t)$ retain the terms that split a loop into
two, with the same initial data. Following the five-step bootstrap
scheme of \cite{YY,TYY}, we control
$D_k=\cL_k-\cK_k$ up to times close to $1$
(Theorem~\ref{gc:thm:common}). Matching the parameters at the terminal
time to those of the target model then gives the required loop and
resolvent bounds.
The tree representation of primitive loops, their Ward identities,
and their cancellation properties come from this earlier development.
Fluctuation averaging uses the Dyson-residual approach of
He--Knowles--Rosenthal \cite{HKR}, as developed for block operators
in \cite[Section~6]{TYY}.

\subsection{Extremely small, small, and bounded coupling}

\paragraph{Extremely small coupling.}
Invariance under the same unitary rotation in every block gives a
localization length of at least $cW$, apart from an exponentially small
expected proportion of eigenvectors (Lemma~\ref{thm:orbital}).
Thus the spreading bound for $g\le W^{-1/2+o(1)}$ follows immediately.

\paragraph{Small coupling.}
For fixed small $\epsilon>0$ and polynomially small couplings
$W^{-1/2+\epsilon}\le g\le W^{-\epsilon}$, the one-dimensional spreading result and the corresponding
loop-hierarchy argument are already established in
\cite[Theorem~2.1]{TYY} in polynomial finite volume. The corresponding
resolvent moment bound here is included in Theorem~\ref{thm:bounded-lower}.

\paragraph{Bounded coupling.}
In the range $W^{-\epsilon}\le g\le O(1)$, the required resolvent
and stability estimates have a direct precedent in \cite[Section~8]{DYYY},
which studies block Anderson operators in dimensions $d\ge3$.
We adapt these inputs to the one-dimensional setting in
Lemma~\ref{fo:lem:bounded-initial}; propagation then gives the
resolvent moment bound of Theorem~\ref{thm:bounded-lower}.

\subsection{Large coupling: transfer matrices and initialization}
For a long time, we did not see how to initialize the loop hierarchy
at large coupling. The flow argument and the five-step bootstrap
scheme could proceed once the initial estimates were available.
The difficulty was obtaining these estimates: even at times bounded
away from $1$, the spectral height is already $\eta\asymp g^{-2}$.
In the normalized model, large $g$ means weak onsite randomness,
and our usual resolvent methods did not give the bounds needed to
start the flow.

The way forward was to reverse our perspective: weaker randomness
also brings the model closer to deterministic free propagation,
which is stable in the bulk. Taking inspiration from the
transfer-matrix methods used to prove localization, we therefore
sought direct control of Green-function amplitudes on the block
scale $g^2$. This is much shorter than the predicted localization
scale of $Wg^2$ blocks. Thus the required initial bounds on
Green-function block norms can already be obtained directly from
transfer-matrix estimates.

For the normalized model with $g\ge1$, fix a block $b$ and
$z=E+i\eta$ with $\eta>0$, and write $G(z)=(H-z)^{-1}$.
Away from the source block $b$, the resolvent equation gives
\[
 \binom{G_{j+1,b}}{G_{jb}}
 =T_j(z)\binom{G_{jb}}{G_{j-1,b}},\qquad
 T_j(z)=\begin{pmatrix}zI_W-V_j/g&-I_W\\I_W&0\end{pmatrix}.
\]
Number the blocks consecutively along an arc of the cycle. For $j\le j'$,
provided none of the blocks $j,\ldots,j'$ equals $b$ modulo $L$,
iteration gives
\begin{equation}\label{GGT}
\binom{G_{j'+1,b}}{G_{j'b}}
 =T_{j'}(z)\cdots T_j(z)\binom{G_{jb}}{G_{j-1,b}}.
\end{equation}
All Green blocks in this recurrence are evaluated at $z$.
For $|E|\le2-\kappa$, write $E=2\sin\theta$ with $|\theta|<\pi/2$. The change of
coordinates $S_E$ from Definition~\ref{def:transfer-matrices} separates
the free rotation from the random perturbation:
\[
 S_E^{-1}\cdot T_j(z)\cdot S_E
   =R_E+g^{-1}X_j+i\eta B_E,\qquad
 R_E=\diag(-ie^{i\theta},ie^{-i\theta})\otimes I_W.
\]
Here $X_j$ is a centered linear function of the onsite GUE matrix $V_j$,
while $B_E$ is deterministic and bounded. For fixed $p,C_0>0$ and
all sufficiently large $W$, Lemma~\ref{uvs:lem:transfer-moments} gives
\[
 \sup_j\E\max_{j\le j'<j+n}
 \left(\|T_{j'}(z)\cdots T_j(z)\|^p
       +\|(T_{j'}(z)\cdots T_j(z))^{-1}\|^p\right)\le C_p,
\]
for $n\le C_0g^2$ and $n\eta\le C_0$. Thus the norms of the stacked block
columns $\binom{G_{j+1,b}}{G_{jb}}$ remain comparable along at most
$n$ transfer steps avoiding $b$, with comparison
factors having bounded moments. For sufficiently long cycles,
summing the corresponding quadratic-form inequalities over $n\asymp g^2$
blocks and using the matrix Ward identity,
\[
 \eta\sum_j G_{jb}(z)^*G_{jb}(z)=\Im G_{bb}(z),
\]
we obtain the initial bound of Lemma~\ref{fo:lem:resolvent-moments}:
\[
 \E\|G_{ab}(E+i\eta)\|^p\le C_p,
 \qquad \eta\asymp g^{-2},
\]
uniformly in $g,L$.
The same transfer and Ward estimates give spreading on the coordinate
scale $W\min\{L,g^2\}$ (Lemma~\ref{uv:lem:short-spreading}), handling
short cycles directly.

The weighted Ward estimate of
Lemma~\ref{uvc:lem:weighted-ward} adds exponential spatial decay to the
initial resolvent-moment bound:
\[
 \E\|G_{ab}(z)\|^p\le C_p e^{-pc\|a-b\|_L/g^2},
\]
where $\|a-b\|_L$ is the block distance on the cycle.
Together with unitary invariance, Gaussian fluctuation averaging,
and Dyson stability, these bounds allow an early-time loop induction
that yields Lemma~\ref{fo:prop:inputs}, under its stated volume assumptions.
Propagation then gives the resolvent moment bound of
Theorem~\ref{uvu:thm:terminal}.

\subsection{Extremely large coupling: spatial sampling}
When $\log g\gg\log W$, even the reduced cycle can contain more
than polynomially many blocks. A union bound over all block tuples
would then destroy estimates whose failure probabilities are powers
of $W$. On the other hand, in this regime, the transfer-matrix product
$T_n(z)\cdots T_1(z)$ is extremely stable when $n\lesssim g^2W^{-B}$,
for $B>2$. More precisely, for $z=E+i\eta$, Lemma~\ref{uvc:lem:cells}
gives, with arbitrarily high polynomial probability,
\[
 T_n(z)\cdots T_1(z)=T(E)^n+O(W^{-B/4}),\qquad
 T(E)=\begin{pmatrix}EI_W&-I_W\\I_W&0\end{pmatrix},
 \qquad n\lesssim g^2W^{-B},\quad 0<\eta\lesssim g^{-2}.
\]
For $j\le j'$ such that none of the blocks $j,\ldots,j'$ equals $b$
modulo $L$, write $n=j'-j+1$.
Using \eqref{GGT}, we obtain
\[
 \binom{G_{j'+1,b}}{G_{j'b}}
 \approx T(E)^n\binom{G_{jb}}{G_{j-1,b}},
 \qquad n\lesssim g^2W^{-B}.
\]
Thus there are scalar coefficients $c_{j',j}$ and $c_{j',j-1}$ such that
\[
 G_{j'b}\approx c_{j',j}G_{jb}+c_{j',j-1}G_{j-1,b}.
\]
An analogous approximation holds in the index $b$. Recall that
$L\lesssim g^2W^2$ from the above spatial cutoff argument. Partition $\mathbb Z_L$
into $W^{O(1)}$ consecutive cells of length less than $g^2W^{-B}$,
and let $\mathcal J$ be the set of their first blocks. The set
$\mathcal J$ is a $g^2W^{-B}$-net in $\mathbb Z_L$. Then any
$G_{ab}$ can be approximated by the scalar source term $d_{ab}I_W$
plus a linear combination of the four blocks
\[
 G_{\alpha,\beta},\quad G_{\alpha-1,\beta},\quad
 G_{\alpha,\beta-1},\quad G_{\alpha-1,\beta-1},\quad \alpha,\; \beta \in {\cal J}
\] 
where $\alpha$ and $a$ lie in the same cell, as do $\beta$ and $b$.
Proposition~\ref{uvc:prop:green}
interpolates every Green block, with error $\mathcal E_{ab}$ and scalar source
term $d_{ab}I_W$, as follows:
\[
 G_{ab}
 =\sum_{\alpha,\beta}c_{a\alpha}c_{b\beta}G_{\alpha\beta}
       +d_{ab}I_W+\mathcal E_{ab},
\]
where $\alpha,\beta\in\mathcal J\cup(\mathcal J-1)$.
By Lemma~\ref{uvc:lem:space}, both the resolvent loops $\cL_k$ in
\eqref{fo:eq:loops} and the primitive loops $\cK_k$ in
\eqref{fo:eq:primitive} admit similar approximations.
By the maximum-volume argument of Lemma~\ref{uvc:lem:fekete},
Lemma~\ref{uvc:thm:ordinary-sampling} gives deterministic sets $\mathcal S_k$
with $|\mathcal S_k|\le W^C$ such that, with probability at least $1-W^{-D}$,
\[
 \max_{\boldsymbol a}|D_k(\boldsymbol a,\boldsymbol\sigma)|
 \le C_k\max_{\boldsymbol a\in\mathcal S_k}
       |D_k(\boldsymbol a,\boldsymbol\sigma)|+W^{-D},
 \qquad D_k=\cL_k-\cK_k.
\]

We can therefore restrict the estimates to these sampling points and
apply the same five-step scheme of \cite{YY,TYY} as in the other
coupling regimes to prove Theorem~\ref{uvu:thm:terminal} also in this regime.

\begin{figure}[H]
\centering
\hypersetup{allcolors=black}
\begin{tikzpicture}[
 method/.style={draw=black,rounded corners=2pt,align=center,
   font=\fontsize{8}{9.5}\selectfont,inner xsep=4pt,inner ysep=5pt,
   minimum height=0.9cm},
 framework/.style={method,text width=11.5cm},
 result/.style={framework,fill=black!3},
 route/.style={->,line width=0.5pt,draw=black},
 support/.style={route,dashed},
 note/.style={font=\fontsize{7.5}{9}\selectfont,fill=white,inner sep=2pt}]
 \node[framework] (cutoff) at (7.8,0)
   {$\Delta$: \textbf{Spatial cutoff}\\
    Reduce to a local cycle with $L\lesssim\ellscale$ blocks\\
    Theorem~\ref{uv:thm:smooth-cutoff}};
 \node[framework] (loops) at (7.8,-2.3)
   {\textbf{The basic strategy: loop hierarchy}\\
    Compare random loops $\cL_k$ with primitive loops $\cK_k$;
    the five-step scheme of \cite{YY,TYY}\\
    $\Delta$: Primitive-loop estimates (Lemma~\ref{fo:lem:deterministic});
    fluctuation averaging};
 \draw[route] (cutoff.south) --
   node[note,right] {Local resolvent estimates} (loops.north);

 \node[method,text width=3.8cm] (bounded) at (2.5,-5)
   {\textbf{Small and bounded coupling}\\
    Loop hierarchy: \cite{TYY}\\
    Stability method adapted from \cite{DYYY}\\
    Initial loop estimates: Lemma~\ref{fo:lem:bounded-initial}};
 \node[method,text width=3.8cm] (transfer) at (8,-5)
   {\textbf{Large coupling}\\
    $\Delta\Delta$: \textbf{Transfer matrices}\\
    Free rotation + small random perturbations\\
    Stability on $g^2$ blocks: Lemma~\ref{uvs:lem:transfer-moments}};
 \draw[route] ([xshift=-3cm]loops.south) -- (4.8,-3.5)
   -- (2.5,-3.5) -- (bounded.north);
 \draw[route] ([xshift=0.2cm]loops.south) -- (transfer.north);

 \node[method,text width=3.8cm] (seed) at (8,-7.7)
   {\textbf{Ward + weighted Ward}\\
    Initial moments at $\eta\asymp g^{-2}$\\
    Exponential spatial tails\\
    Lemmas~\ref{fo:lem:resolvent-moments}, \ref{uvc:lem:weighted-ward}};
 \draw[route] (transfer.south) -- (seed.north);
 \node[method,text width=3.8cm] (short) at (13.2,-5)
   {\textbf{Short cycles}\\
    Ward + Poisson smoothing\\
    $\eta\asymp1/\min\{L,g^2\}$\\
    Lemma~\ref{uv:lem:short-spreading}};
 \draw[route] (transfer.east) -- (short.west);
 \node[method,text width=3.8cm] (sampling) at (13.2,-8.3)
   {\textbf{Extremely large coupling}\\
    $\Delta\Delta$: \textbf{Spatial sampling}\\
    Retain the free oscillation\\
    Interpolate within each cell\\
    Only $W^{O(1)}$ sampled tuples\\
    Proposition~\ref{uvc:prop:green}; Lemma~\ref{uvc:thm:ordinary-sampling}};
 \draw[route] (transfer.east) -- (10.7,-5)
   -- (10.7,-8.3) -- (sampling.west);

 \node[method,text width=4.2cm] (initial) at (8,-10.3)
   {\textbf{Initial loop estimates}\\
    Lemma~\ref{fo:prop:inputs}};
 \draw[route] (seed.south) -- (initial.north);
 \draw[support] (sampling.south) -- (13.2,-10.3) -- (initial.east);

 \node[result] (conclusion) at (7.8,-12.7)
   {\textbf{Common loop propagation and terminal matching}\\
    Theorem~\ref{gc:thm:common}\textup{(ii)}\\
    Resolvent moments: Theorems~\ref{thm:bounded-lower}, \ref{uvu:thm:terminal}};
 \draw[route] (bounded.south) -- (2.5,-11.4)
   -- (4.8,-11.4) -- ([xshift=-3cm]conclusion.north);
 \draw[route] (initial.south) -- ([xshift=0.2cm]conclusion.north);
 \draw[support] (seed.west) -- (5.45,-7.7)
   -- node[note,left] {Coarse moments} (5.45,-11.2)
   -- (6.7,-11.2) -- ([xshift=-1.1cm]conclusion.north);
 \draw[support] (sampling.east) -- (15.5,-8.3)
   -- (15.5,-12.7) -- (conclusion.east);
 \node[result] (spreading) at (7.8,-15)
   {\textbf{Spreading lower bound}};
 \draw[route] (conclusion.south) --
   node[note,right] {Theorem~\ref{uv:lem:lower-volume}} (spreading.north);
 \draw[route] (short.east) -- (16,-5)
   -- (16,-15) -- (spreading.east);
\end{tikzpicture}
\caption{Lower-bound methods in the order of the preceding explanation.
$\Delta$: main contributions of this work; $\Delta\Delta$: the two central contributions.}
\label{fig:lower-methods}
\end{figure}
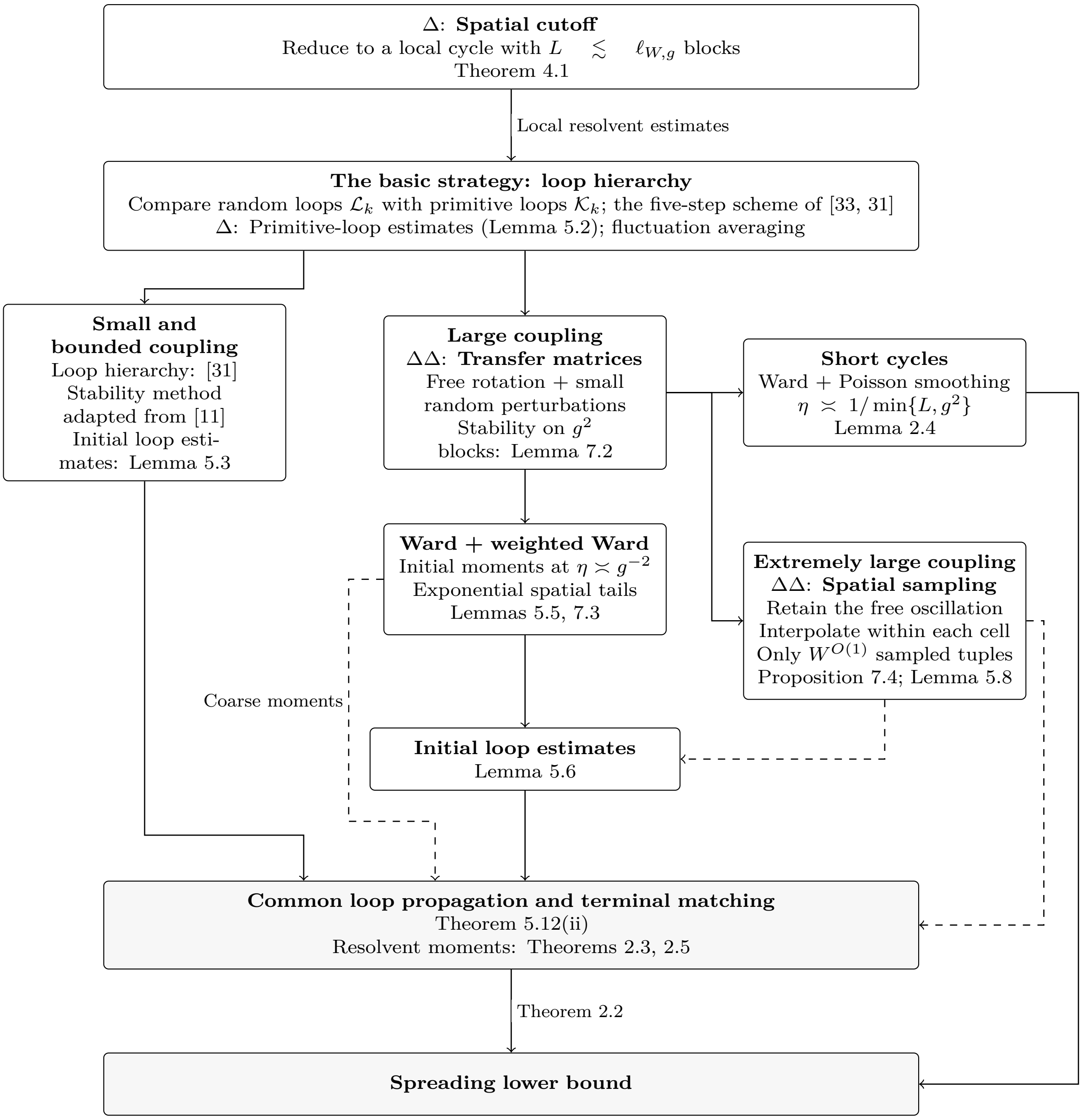

\subsection{Component estimates and proof of the lower bound}
\label{sec:extension-regions}

Lemma~\ref{thm:orbital} handles the single-block scale, and
Lemma~\ref{uv:lem:short-spreading} handles short cycles at large coupling.
For the remaining parameters, Theorem~\ref{uv:lem:lower-volume} reduces spreading to a
resolvent moment bound on a local cycle of $L_c<\ellscale$ blocks,
with $g$ unchanged. Theorems~\ref{thm:bounded-lower} and~\ref{uvu:thm:terminal}
supply this bound in the respective coupling ranges.
Figures~\ref{fig:lower-methods} and~\ref{fig:lower-ranges} summarize
the methods and coupling regimes in the order of the preceding explanation.

\begin{figure}[!t]
\centering
\renewcommand{\arraystretch}{1.35}
\fontsize{8.5}{10.5}\selectfont
\begin{tabularx}{\textwidth}{@{}p{0.30\textwidth}p{0.40\textwidth}X@{}}
\toprule
\textbf{Coupling} & \textbf{Main idea} & \textbf{Corresponding results}\\
\midrule
\textbf{Extremely small}\newline
$g\le W^{-1/2+o(1)}$
& Common unitary invariance already gives the single-block scale.
& Lemma~\ref{thm:orbital}\\[3pt]
\textbf{Small}\newline
$W^{-1/2+\epsilon}\le g\le W^{-\epsilon}$
& The loop-hierarchy method of \cite{TYY}.
& Theorem~\ref{thm:bounded-lower}\\[3pt]
\textbf{Bounded}\newline
$W^{-\epsilon}\le g\le C$
& Adapt the resolvent and stability method of \cite{DYYY}, then use the same loop hierarchy.
& Lemma~\ref{fo:lem:bounded-initial}; Theorem~\ref{thm:bounded-lower}\\[3pt]
\textbf{Large}\newline
$g\ge C_\kappa$
& Transfer and Ward estimates give short-cycle spreading directly; the remaining cycles use initialization and loop propagation.
& Short cycles: Lemma~\ref{uv:lem:short-spreading}.\newline
Remaining cycles: Theorem~\ref{uvu:thm:terminal}, using Lemmas~\ref{fo:lem:resolvent-moments}--\ref{fo:prop:inputs}.\\[3pt]
\textbf{Extremely large}\newline
$\log g\gg\log W$
& For the remaining cycles, retain free oscillations and interpolate within cells: only $W^{O(1)}$ sampled tuples are needed for loop propagation.
& Short cycles: Lemma~\ref{uv:lem:short-spreading}.\newline
Remaining cycles: Theorem~\ref{uvu:thm:terminal}, using Proposition~\ref{uvc:prop:green} and Lemma~\ref{uvc:thm:ordinary-sampling}.\\
\bottomrule
\end{tabularx}
\normalsize
\caption{Qualitative coupling regimes in the order of the preceding explanation,
with fixed $0<\epsilon<1/4$ and $C<\infty$.
The regimes may overlap; the last row extends the large-coupling argument.
The cited results give the precise energy and volume assumptions.}
\label{fig:lower-ranges}
\end{figure}

\begin{lemma}[Spreading from unitary invariance]
\label{thm:orbital}
For the model \eqref{eq:model}, there are absolute constants $c_0,c>0$
such that, for all sufficiently large $W$ and all $L,g$,
\begin{equation}\label{eq:haar-count}
 \frac1{WL}\E\#\{k:\ell(\psi_k)<c_0W\}\le4We^{-cW}.
\end{equation}
There is no energy restriction. Moreover, for every $D>0$ there is
$c_D>0$ such that, with probability at least $1-N^{-D}$,
simultaneously for all eigenvectors,
\begin{equation}\label{eq:orbital-length}
 \ell(\psi_k)\ge c_D(\log N)^{-1}W.
\end{equation}
\end{lemma}
The proof is in Appendix~\ref{subsec:pf_thm31}.
Lemma~\ref{thm:orbital} gives the required lower bound in the region
\begin{equation}\label{uv:eq:haar-lower-region}
 \ellstar\le c_0W^{1+\tau}.
\end{equation}

\subsubsection*{From resolvent moments to spreading at arbitrary volume.}
\begin{theorem}[Resolvent moments imply spreading, with volume extension]
\label{uv:lem:lower-volume}
Fix $E_0,D>0$, $0<\tau<1$, and $0<\nu<\tau/2$.
Let $J\subset[-E_0,E_0]$ be a deterministic interval. Put
\begin{equation}\label{uv:eq:lower-volume-scales}
 \eta=W^\nu/\ellscale,\qquad R=\lceil2/\eta\rceil.
\end{equation}
Choose the cycle length and moment height
\[
 L_c=\min\{L,2R+1\},\qquad
 \eta_c=\frac{W^\nu}{\ell_{W,g,L_c}},\qquad
 J_{\eta_c}=J+[-\eta_c,\eta_c].
\]
Let $G_{L_c}(z)=(H_{g,W,L_c}-z)^{-1}$ be the resolvent of the normalized
periodic model \eqref{eq:model} on this cycle.
Suppose that, for some $p>1+(D+1)/(\tau-\nu)$ and a constant $C_p$,
\begin{equation}\label{uv:eq:lower-volume-moment}
 \sup_{E\in J_{\eta_c}}\max_{1\le x\le WL_c}
 \E\left(\Im(G_{L_c}(E+i\eta_c))_{xx}\right)^p\le C_p.
\end{equation}
Then, for all sufficiently large $W$,
\begin{equation}\label{uv:eq:lower-volume-output}
 \frac1{WL}\E\#\left\{k:\lambda_k\in J,\
       \ell(\psi_k)<W^{-\tau}\ellstar\right\}\le W^{-D}.
\end{equation}
The width threshold depends only on $E_0,D,\tau,\nu,p,C_p$ and the
threshold in the moment assumption. It is independent of $L,g$
whenever that assumption holds uniformly in these parameters.
\end{theorem}
The proof is in Section~\ref{uv:sec:local-resolvents}.
The assumption requires only moments at deterministic energies;
in applications we keep the original interval $J=J_{g,\kappa,L}$
on the cycle of $L_c$ blocks.

\medskip\noindent\textbf{Bounded coupling.}
\begin{theorem}[Bounded-coupling resolvent moments]
\label{thm:bounded-lower}
Fix $C,p>0$, $0<\kappa<1/2$, and $0<\nu<1$.
Suppose that the deterministic coupling $g=g(W)$ and integer cycle
length $L=L(W)$ satisfy
\begin{equation}\label{eq:bounded-lower-range}
 W^{-1/2}\le g\le C,\qquad \log W\ll L\le W^2.
\end{equation}
For all sufficiently large $W$,
\begin{equation}\label{eq:bounded-lower-moment}
 \sup_{E:\,F_g(-|E|)\ge\kappa}\max_{1\le x\le WL}
 \E\bigl(\Im G_{xx}(E+i\eta_L)\bigr)^p\le C_p,
 \qquad \eta_L=\frac{W^\nu}{\ellstar}.
\end{equation}
The moment constant depends only on $C,p,\kappa,\nu$.
For every fixed positive function $\omega(W)\to\infty$, the width threshold
may additionally depend on $\omega$, and is uniform over deterministic
parameters satisfying $W^{-1/2}\le g\le C$ and
$\omega(W)\log W\le L\le W^2$.
\end{theorem}
The proof in Section~\ref{sec:fo-completion} follows the loop hierarchy
method of \cite{TYY}; we prove the bounded-coupling inputs here by adapting
the resolvent and stability method of
\cite[Section~8 and Appendix~B.1]{DYYY}.

\medskip\noindent\textbf{Large coupling.}
\begin{lemma}[Spreading on the scale $g^2$]\label{uv:lem:short-spreading}
Fix $0<\kappa<1/2$. There are constants $W_0,C_\kappa,c_\kappa>0$
such that, for every $W\ge W_0$, all $L$, and all $g\ge C_\kappa$,
\begin{equation}\label{uv:eq:constant-scattering-outline}
 \frac1{WL}\E\#\left\{k:\lambda_k\in J_{g,\kappa,L},\
   \ell(\psi_k)<c_\kappa W\min\{L,g^2\}\right\}
 \le C_\kappa e^{-c_\kappa W}.
\end{equation}
The constants are independent of $L,g$.
\end{lemma}
The proof in Section~\ref{sec:uniform-resolvent-seed} uses the transfer
estimates and Ward's identity, followed by Poisson smoothing.

\begin{theorem}[Large-coupling resolvent moments]
\label{uvu:thm:terminal}
Fix $0<\kappa<2$, $0<\nu<1$, and $p>0$.
There are constants $C_\kappa>1$ and $C<\infty$ such that, for the periodic
model \eqref{eq:model}, if deterministic $g=g(W)$ and $L=L(W)$ satisfy
\[
 g\ge C_\kappa,\qquad L\gg g^2\log W,\qquad L\le g^2W^2,
\]
then, for all sufficiently large $W$,
\begin{equation}\label{uvu:eq:local-moments}
 \sup_{|E|\le2-\kappa}\ \max_{1\le x\le WL}
 \E\bigl(\Im G_{xx}(E+i\eta_L)\bigr)^p\le C,
 \qquad \eta_L=\frac{W^\nu}{\ellstar}.
\end{equation}
Here $L\gg g^2\log W$ means $L/(g^2\log W)\to\infty$ as $W\to\infty$.
The moment constant $C$ depends only on $\kappa,\nu,p$, not on $g,L$.
For every fixed positive function $\omega(W)\to\infty$, the width threshold
may additionally depend on $\omega$, and is uniform over deterministic
$g\ge C_\kappa$ and $\omega(W)g^2\log W\le L\le g^2W^2$.
\end{theorem}

The proof of Theorem~\ref{uvu:thm:terminal} is completed in
Section~\ref{sec:fo-completion}.

\begin{proof}[Proof of the lower bound in Theorem~\ref{thm:main}]
For both bounds, the case $\tau\ge1$ follows from $\tau=1/2$ by monotonicity.
Fix $0<\tau<1$, $D>0$, and $0<\nu<\tau/2$, and choose
$p$ as in Theorem~\ref{uv:lem:lower-volume}. Define
\[
 B_-=\#\{k:\lambda_k\in J_{g,\kappa,L},\ \ell(\psi_k)<W^{-\tau}\ellstar\}.
\]
Lemma~\ref{thm:orbital} covers \eqref{uv:eq:haar-lower-region}.
Otherwise $\min\{WL,W^2g^2\}=\ellstar>c_0W^{1+\tau}$ for large $W$,
so $L>c_0W^\tau$ and $g\ge W^{-1/2}$.
Set $J:=J_{g,\kappa,L}$. By \eqref{eq:bulk-sampling} and
Lemma~\ref{da:lem:distribution}\textup{(i)}, $J\subseteq[-4,4]$ and
$F_g(-|E|)\ge\kappa/2$ on $J$.
Fix a sufficiently large $C_*$ depending only on $\kappa,\nu,p$.

If $g\ge C_*$ and $L\le W^{\tau/2}g^2$, then $\ellstar=WL$ and
$W^{-\tau}\ellstar\le W^{-\tau/2}W\min\{L,g^2\}$.
Lemma~\ref{uv:lem:short-spreading} therefore gives
$(WL)^{-1}\E B_-\le C_\kappa e^{-c_\kappa W}\le W^{-D}$
for sufficiently large $W$.

For the remaining parameters, select $L_c,\eta_c$ as in
Theorem~\ref{uv:lem:lower-volume} with $E_0=4$.
The exclusions above and \eqref{uv:eq:lower-volume-scales} give
\[
 \log W\ll L_c,\qquad
 \begin{cases}
 L_c\le W^2,&g<C_*,\\
 g^2W^{\tau/2}\le L_c\le g^2W^2,&g\ge C_*.
 \end{cases}
\]
For $g\ge C_*$, we thus have
$L_c/(g^2\log W)\ge W^{\tau/2}/\log W\to\infty$ uniformly.
Lemma~\ref{da:lem:distribution}\textup{(ii)} gives
\[
 \begin{aligned}
 \eta_c&\le W^{\nu-1}=o(1),\\
 F_g(-|E|)&\ge\frac\kappa2-\sup_{x\in\R}
 \nu_g([x-\eta_c,x+\eta_c])\ge\frac\kappa4,
 \qquad E\in J_{\eta_c},
 \end{aligned}
\]
uniformly in $g,L$. Apply Theorem~\ref{thm:bounded-lower} for $g<C_*$
and Theorem~\ref{uvu:thm:terminal} for $g\ge C_*$ at this reduced margin,
using Lemma~\ref{da:lem:distribution}\textup{(iii)} in the latter case.
This supplies \eqref{uv:eq:lower-volume-moment}, and
Theorem~\ref{uv:lem:lower-volume} yields
\begin{equation}\label{uv:eq:lower-regions-output}
 (WL)^{-1}\E B_-\le W^{-D}.
\end{equation}
\end{proof}

\section{Proof strategy for the upper bound}
\label{sec:intro-concentration-methods}

Our goal is to place at least half the mass of most eigenvectors
in the selected energy interval within $\max\{W,W^2g^2\}$ consecutive
coordinates, up to arbitrarily small powers of $W$. The total-volume
cap is $WL$, so we may assume
$L\ge Wg^2$.

\subsection{From fractional moments to transfer growth}
\label{sec:intro-upper-reductions}

\paragraph{The starting point: Fractional moment method.}
The classical fractional-moment method \cite{AM,ASFH} reduces the
upper bound to decay of Green-function moments. For the Hamiltonian
$H_\Lambda$ on a finite path or cycle $\Lambda$ of $L$ blocks,
write $G_\Lambda(x,y;E)$ for the $(x,y)$ block of
$(H_\Lambda-E)^{-1}$ and $\|x-y\|_L$ for the block distance.
The Gaussian regularization estimates of \cite{APSSS} control the
local Green-function factors used in this reduction.
For $g<W$, we use the fixed-energy input
\eqref{up:eq:fmr-pointwise} on the selected energy window $J$:
\[
 \sup_{E\in J}\E\|G_\Lambda(x,y;E)\|^q
 \le W^C e^{-c\|x-y\|_L/\mathfrak b},\qquad 0<q<1,
\]
with $\mathfrak b=\max\{1,Wg^2\}$ up to arbitrarily small powers of $W$.
Lemma~\ref{uva:thm:upper-audit}\textup{(i)} converts this into
concentration within $C\mathfrak b\log W$ blocks.

The fixed-energy argument applies for $g\le W^B$ for any fixed $B>0$,
since the factor $g$ from single-block spectral averaging can then be
absorbed into a power of $W$. The two methods overlap. We switch at
$g=W$ because the negative-moment estimate of Theorem~\ref{uvg:thm:transfer}
covers the entire bulk, including $E=0$, uniformly for $g\ge W$.
In this range we use the energy-averaged version
\eqref{uva:eq:conversion-input} on finite Dirichlet paths,
\[
 \int_J\E\|G_\Lambda(x,y;E)\|^q\,dE
 \le W^C e^{-c\|x-y\|_L/(Wg^2)},
\]
where $J$ is a fixed bulk window. To remove the factor $g$ for arbitrarily
large couplings, we average over groups of $m\asymp g^2$ blocks. Transfer products
and their inverses remain stable on this scale
(Lemma~\ref{uvs:lem:transfer-moments}); normalization reduces the
coupling-dependent cost to $g/\sqrt m\asymp1$.
Lemma~\ref{uva:thm:upper-audit}\textup{(ii)} then gives concentration
within $CWg^2\log W$ blocks.

For $g\le c_*W^{-1/2}$, the required fractional-moment decay is
already supplied by \cite[Theorem~1]{PSSS}, yielding
Theorem~\ref{thm:strong-upper}. In the remaining ranges,
we obtain it from transfer growth, as described next.

\paragraph{Reduction to transfer growth estimates.}
We follow the transfer-to-Green-function comparison of Macera--Sodin
\cite[Section~2.2]{MS}.
To illustrate the passage from transfer growth to Green-function decay,
take $g\ge1$ and a local interval $\Lambda=[0,n]$.
Let $L_0,R_0$ be Hermitian endpoint shifts and let $\iota_j$ insert
a vector into block $j$. As in \eqref{up:eq:shifted-resolvent}, set
\[
 H_\Lambda^{L_0,R_0}
 :=H_\Lambda+\iota_0L_0\iota_0^*
             +\iota_nR_0\iota_n^*.
\]
This is the standard boundary reduction: the shifts vanish for
Dirichlet truncation and encode the eliminated exterior of a larger
path. At a real energy $E$ away from poles, write
\[
 \widetilde G=(H_\Lambda^{L_0,R_0}-E)^{-1},
\]
and denote its blocks by $\widetilde G_{ab}$.
Using the transfer matrices $T_j(E)$ from
Section~\ref{sec:intro-spreading-methods}, iteration from $n$ to $1$,
away from the source at $0$, gives
\[
 \binom{\widetilde G_{00}}{\widetilde G_{10}}
 =T_1(E)\cdots T_n(E)
       \binom{\widetilde G_{n0}}{\widetilde G_{n+1,0}},\quad \widetilde G_{n+1,0}:=R_0\widetilde G_{n0}.
\]
Consequently,
\[
 \|\widetilde G_{n0}\|
 \le\left(\inf_{\|v\|=1}
   \left\|T_1(E)\cdots T_n(E)\binom{v}{R_0v}\right\|\right)^{-1}
   \bigl(\|\widetilde G_{00}\|+\|\widetilde G_{10}\|\bigr).
\]
The infimum measures the least growth allowed by the right boundary.
Use the fixed unitary change of coordinates $S$ from
Definition~\ref{def:transfer-matrices}, and let $\mathcal F_W$ denote
the set of orthonormal boundary frames in these coordinates:
\[
 S=\frac1{\sqrt2}\begin{pmatrix}1&1\\i&-i\end{pmatrix}\otimes I_W,
 \qquad
 \mathcal F_W=\left\{\Phi\in\C^{2W\times W}:
 \Phi^*\Phi=I_W,\quad
 \Phi^*\diag(I_W,-I_W)\Phi=0\right\}.
\]
For every Hermitian $R_0$, there exists $\Phi_{R_0}\in\mathcal F_W$
such that $S\Phi_{R_0}$ is an
orthonormal frame of $\{(v,R_0v):v\in\C^W\}$
(Definition~\ref{def:boundary-frames}). Normalizing the boundary data gives
\[
 \|\widetilde G_{n0}\|
 \le\left(\inf_{\|v\|=1}
   \left\|T_1(E)\cdots T_n(E)S\Phi_{R_0}v\right\|\right)^{-1}
   \bigl(\|\widetilde G_{00}\|+\|\widetilde G_{10}\|\bigr).
\]
To make $\widetilde G_{n0}$ small, we seek expansion by a power of $W$
with high probability at a suitable deterministic block length.
The bound must be uniform over deterministic $R_0$, but need not hold
simultaneously for all $R_0$.
Define the amplitude matrix as in \eqref{up:eq:D-definition}:
\[
 D_n(\Phi;E)=\sqrt2\begin{pmatrix}I_W&0\end{pmatrix}
                  S^*T_n(E)\cdots T_1(E)S\Phi.
\]
Here we switch the product order for convenience. The smallest singular
value of $D_n(\Phi;E)$ measures the weakest expansion:
\[
 \|D_n(\Phi;E)^{-1}\|^{-1}
 =\inf_{\|v\|=1}\|T_n(E)\cdots T_1(E)S\Phi v\|.
\]
For general $g>0$, we henceforth use the normalized transfers and
$D_n$ from Definition~\ref{def:boundary-frames}.
For $g<W$, the near-source Green-function factors
have fractional moments bounded by powers of $W$. Sufficiently strong
high-probability transfer growth absorbs this cost and yields the local
resolvent target \eqref{uv:eq:upper-volume-input}
(Proposition~\ref{up:prop:conversion}\textup{(i)}).
For fixed $D>0$, the high-probability input has the form below, with
$\mathfrak b$ supplied by the estimates for each coupling range:
\begin{equation}\label{DnWDWD}
 \sup_{E\in J}\ \sup_{\Phi\in\mathcal F_W}
 \Prob\!\left(\|D_{n_0}(\Phi;E)^{-1}\|>W^{-D}\right)
 \le W^{-D},\qquad
 n_0\lesssim\mathfrak b\log W.
\end{equation}

For $g\ge W$, we instead use the negative-moment bound from
Section~\ref{sec:intro-direct-transfer}: for some fixed $s\in(0,1)$
and a fixed bulk interval $J$,
\begin{equation}\label{intro:eq:negative-transfer}
 \sup_{E\in J}\sup_{\Phi\in\mathcal F_W}
 \E\|D_n(\Phi;E)^{-1}\|^s
 \le CW e^{-cn/(Wg^2)},\qquad n\ge1.
\end{equation}

\subsection{Reduction to operator contraction}
\label{sec:intro-operator-contraction}
In this subsection, we discuss the proof of \eqref{DnWDWD} for
$c_*W^{-1/2}\le g<W^\epsilon$, and away from zero energy for
$W^\epsilon\le g<W$, with any fixed $0<\epsilon<1$.
Using the projective and flag geometry of \cite{Furstenberg,Lessa},
and adapting Drogin's framework \cite[Sections~2 and~4]{Drogin}, we reduce
growth in the weakest direction of a long transfer product to
contraction of $P_{1/2}$.
Our target $\|D_n^{-1}\|^{-1}$ is the smallest singular value of $D_n$.
The total volume expansion records the product of all $W$ singular
values; dividing by the expansion on a $(W-1)$-dimensional subspace
isolates one transverse direction. For
$\mathcal T_n=S^*T_n\cdots T_1S$, $Y=\ran\Phi$ with $\Phi\in\mathcal F_W$, and
$U\subset Y$ of dimension $W-1$, define
\[
 \begin{aligned}
 c(\mathcal T_n;U,Y)&=\log\frac{\operatorname{vol}(\mathcal T_n|_Y)}
                          {\operatorname{vol}(\mathcal T_n|_U)},\qquad
 \|D_n(\Phi)^{-1}\|^{-1}
 &=\inf_{\substack{U\subset Y,\; \dim U=W-1}}e^{c(\mathcal T_n;U,Y)},\quad Y=\ran\Phi
 \end{aligned}
\]
Here $\operatorname{vol}(\mathcal T_n|_Z)$ denotes the volume expansion on $Z$.
Since the exponent $D$ in \eqref{DnWDWD} is arbitrary, polynomial
factors in $W$ can be absorbed. We can therefore use the $W$ row norms
of $D_n^{-1}$ instead of locating the random weakest direction.
For $0<s\le1$, the row identity in \eqref{eq:c-inverse-cocycle} and
$\|D_n^{-1}\|^s\le\sum_{j=1}^W\|e_j^*D_n^{-1}\|^s$ give
\[
 \E\|D_n^{-1}\|^s
 \le W\sup_{\substack{U\subset Y\\\dim U=W-1}}
       \E e^{-s\cdot  c(\mathcal T_n;U,Y)}, \quad Y =\ran\Phi
\]
We group a suitable fixed number $q=O(1)$ of blocks into one step
$\mathcal T_q=S^*T_q\cdots T_1S$ and define
\[
 (P_s\circ h)(U,Y)
 =\E\bigl[e^{-s c(\mathcal T_q;U,Y)}h(\mathcal T_qU,\mathcal T_qY)\bigr].
\]
Here $h$ is a bounded function of the pair. The advantage is that
estimating $P_s$ requires analyzing only $q=O(1)$ transfer matrices,
while $D_n$ contains a long product of $n$ such matrices.
Since $c$ adds under
successive transfers, iteration over independent batches gives,
as in Lemma~\ref{up:lem:operator-reduction},
\[
 \sup_{\Phi\in\mathcal F_W}\E\|D_{qk}(\Phi;E)^{-1}\|^s
 \le W\|(P_s)^k\circ 1\|_\infty,\qquad k\ge1,
\]
where the supremum norm ranges over the pairs $(U,Y)$ above.
The ideal input would be a one-step supremum-norm contraction
\[
 \|P_s\|_{\infty\to\infty}\le1-c/\mathfrak b,
\]
which would immediately give exponential decay of the iterates on
scale $\mathfrak b$.

The proof uses an angular $L^4$ bound. For fixed $U$, the admissible
spaces $Y=\ran\Phi$ containing $U$, with $\Phi\in\mathcal F_W$, form a
circle $\mathcal C_U$ with normalized angular measure $m_U$. Write
$\|h\|_{\mathrm{ang},4}:=\sup_U\|h(U,\cdot)\|_{L^4(\mathcal C_U,m_U)}$
(Definition~\ref{up:def:angular-norm}).
The angular Jacobian satisfies
\[
 \frac{dm_{\mathcal T_qU}(\mathcal T_qY)}{dm_U(Y)}
 =e^{-2c(\mathcal T_q;U,Y)}=\bigl(e^{-s c(\mathcal T_q;U,Y)}\bigr)^{2/s}.
\]
 Take $s=1/2$.
For each fixed $U$, Jensen's inequality and then the change of variables
$Y'=\mathcal T_qY$ give
\begin{equation}\label{intro:eq:angular-nonexpansion}
 \begin{aligned}
 \int_{\mathcal C_U}|(P_{1/2}\circ h)(U,Y)|^4\,dm_U(Y)
 &\le\E\int_{\mathcal C_U}e^{-2c(\mathcal T_q;U,Y)}
       |h(\mathcal T_qU,\mathcal T_qY)|^4\,dm_U(Y)\\
 &=\E\int_{\mathcal C_{\mathcal T_qU}}
       |h(\mathcal T_qU,Y')|^4\,dm_{\mathcal T_qU}(Y')\;\le\|h\|_{\mathrm{ang},4}^4.
 \end{aligned}
\end{equation}
The last step uses the supremum over all possible first subspaces
in the angular norm. Then the supremum over $U$ yields
\[
 \|P_{1/2}\circ h\|_{\mathrm{ang},4}\le\|h\|_{\mathrm{ang},4}.
\]
The constant $1$ comes from the exact match between the fourth power
of the weight $e^{-c/2}$ and the Jacobian $e^{-2c}$.
This bound holds for every bounded measurable $h$.
The same calculation uses $L^{2/s}$ for other small positive $s$;
we choose $s=1/2$ for convenience. To apply
Lemma~\ref{up:lem:operator-iteration}, we need the strict improvement
\[
 \|P_{1/2}\circ h\|_{\mathrm{ang},4}
 \le(1-c/\mathfrak b)\|h\|_{\mathrm{ang},4}
\]
for bounded nonnegative $h$ invariant under the same unitary rotation
in both $W$-dimensional components, as in \eqref{up:eq:invariant-function}.
The functions $h=(P_{1/2})^r\circ1$ have this invariance, so iteration
gives $\|(P_{1/2})^k\circ1\|_{\mathrm{ang},4}\le(1-c/\mathfrak b)^k$.

This angular estimate still allows narrow peaks, so we use smoothing
to obtain control at each fixed boundary pair. Following Drogin's
regularization of the transfer-product law \cite[Section~4]{Drogin},
we split $P_{1/2}=G_{1/2}+B_{1/2}$ using input weights $\chi_P$ and
$1-\chi_P$. The first part maps angular $L^4$ to $L^\infty$, while the
second contracts in $L^\infty$ (Lemma~\ref{up:lem:good-bad}).
Repeated angular contraction, followed by one smoothing step and a
geometric summation of the remainder terms, gives
\[
 \|(P_{1/2})^k\circ1\|_\infty\le W^C e^{-ck/\mathfrak b}
\]
by Lemma~\ref{up:lem:operator-iteration}. The smoothing cost $W^C$
is paid only once. Lemma~\ref{up:lem:operator-reduction} then gives
\eqref{DnWDWD} at $n_0=qk\asymp\mathfrak b\log W$ with a sufficiently
large constant.

\subsection{A criterion for contraction of \texorpdfstring{$P_{1/2}$}{P(1/2)}}
\label{sec:intro-contraction-criterion}
We adapt Drogin's quantitative nonconcentration criterion
\cite[Proposition~2.6 and Claims~2.7--2.8]{Drogin}, using Gaussian
integration by parts; see \cite[Section~3]{Hairer} for this mechanism.
Recall that equality in the Jensen step of
\eqref{intro:eq:angular-nonexpansion} requires, for fixed $(U,Y)$,
the quantity being averaged to be almost surely constant in the
Gaussian input $\omega=(V_1,\ldots,V_q)$.
Denote this quantity by
\[
 Z_h:=e^{-c(\mathcal T_q;U,Y)/2}h(\mathcal T_qU,\mathcal T_qY),
 \qquad (P_{1/2}\circ h)(U,Y)=\E Z_h.
\]
To obtain strict contraction, we therefore need to rule out near
constancy of $Z_h$ in $\omega$ quantitatively.

We vary  $(V_1,\ldots,V_q)$ in a direction
$X=(X_1,\ldots,X_q)$, with each $X_i$ Hermitian. Write $D_X$ for differentiation along this
perturbation. We aim to  choose $X$ so that
\[
 D_X\!\left[c(\mathcal T_q;U,Y)\right]=1,\qquad
 D_X\!\left[h(\mathcal T_qU,\mathcal T_qY)\right]=0.
\]
To ensure the second identity, we require the first-order changes of
the two terminal subspaces $\mathcal T_qU,\mathcal T_qY$ induced by $X$
to be a common unitary rotation acting identically on both
$W$-dimensional components, as in \eqref{up:eq:angular-geometry}.
Every invariant $h$ then has zero derivative. This condition controls
the motion of the subspaces while allowing stretching within them.
The task therefore has two parts:
\begin{enumerate}
\item Achieve $D_X[c(\mathcal T_q;U,Y)]=1$ under this geometric condition.
\item Control the loss of Gaussian probability weight along $X$.
\end{enumerate}
The choice of $X$ may depend on $(U,Y)$
but works for all such $h$.

The construction needs only work on a sufficiently large region of
the Gaussian inputs. A smooth cutoff $0\le\chi\le1$ restricts it
to that region, where the two derivative identities hold, and makes
$\chi X$ a regular field on the full input space. Its Gaussian
divergence is defined by
\[\delta_{\rm G}(\chi X):\quad 
 \E D_{\chi X}\varphi
 =\E[\varphi\,\delta_{\rm G}(\chi X)]
\]
for smooth test functions $\varphi$. This measures the cost of
integration by parts, including the variation of the Gaussian
density and the derivatives of the field and cutoff. The required
bounds are
\[
 \E\chi\ge\tfrac12,\qquad
 \E|\delta_{\rm G}(\chi X)|^2\le C\mathfrak b.
\]
The first prevents us from placing the entire perturbation on a
rare event. The second determines the contraction scale. To see the connection, test integration by parts with a bounded
primitive of the indicator of $[a-t,a+t]$, evaluated at $\log Z_h$.
The primitive has absolute value at most $t$, while
$D_{\chi X}\log Z_h=-\chi/2$. Cauchy--Schwarz therefore gives
\[
 \E\!\left[\chi\,\mathbf1_{\{|\log Z_h-a|\le t\}}\right]
 \le Ct\bigl(\E|\delta_{\rm G}(\chi X)|^2\bigr)^{1/2}
 \le Ct\sqrt{\mathfrak b}.
\]
Taking $t=c_0\mathfrak b^{-1/2}$ with $c_0$ sufficiently small,
we see that $\log Z_h$ cannot concentrate within distance $t$ of
any constant. This quantitative nonconcentration yields the strict
contraction of Lemma~\ref{up:lem:angular}:
\[
 \|P_{1/2}\circ h\|_{\mathrm{ang},4}
 \le(1-c/\mathfrak b)\|h\|_{\mathrm{ang},4}.
\]
The next two subsections explain how to construct $X$ and control
its Gaussian cost in the different coupling ranges.

\subsection{Small coupling: exploiting amplified directions}
\label{sec:intro-small-fields}

We first examine how the transfer matrices act on $U$ and $Y$,
using $T_j$ and the fixed unitary matrix $S$ from
Definition~\ref{def:boundary-frames}. Set $U_0=SU$, $Y_0=SY$, and write
$U_j=T_j\cdots T_1U_0$, $Y_j=T_j\cdots T_1Y_0$.
For each fixed incoming flag, almost surely for every $j\ge1$,
we represent $Y_j$ by a Hermitian matrix $R_j$, and $U_j$ by
$R_j$ and a nonzero vector $v_j$, as
 \[
 \begin{aligned}
 Y_j&=\left\{\binom{x}{-gR_jx}:x\in\C^W\right\},\qquad U_j&=\left\{\binom{x}{-gR_jx}:x\in\C^W,\ v_j^*x=0\right\}.
 \end{aligned}
\]
Choose a nonzero representative $v_1$ and fix the subsequent
representatives by the recurrence
\[
 R_j=(V_j-EI_W-g^2R_{j-1})^{-1},\qquad
 v_j=R_jv_{j-1},\qquad j\ge2.
\]
First consider a simultaneous unitary rotation. Fix $K^*=-K$,
set $\widehat K=\operatorname{diag}(K,K)$, and let
\[
 \begin{aligned}
 V_j(s)&=e^{sK}V_je^{-sK},\qquad 1\le j\le q,\qquad
 U_0(s)&=e^{s\widehat K}U_0,\qquad
 Y_0(s)=e^{s\widehat K}Y_0.
 \end{aligned}
\]
Every propagated pair then undergoes that same unitary rotation:
$U_j(s)=e^{s\widehat K}U_j$, $Y_j(s)=e^{s\widehat K}Y_j$.
The normal representatives may be chosen as $v_j(s)=e^{sK}v_j$.
Take this simultaneous variation of the matrices and the initial
flag as $X^{[0]}$, and write $D_{\rm rot}=D_{X^{[0]}}$. Then
\[
 D_{\rm rot} V_j= [K, V_j],\quad D_{\rm rot}R_j=[K,R_j],\qquad
 D_{\rm rot}v_j=Kv_j,\qquad
 D_{\rm rot}c=D_{\rm rot}h=0.
\]

\textbf{Keeping the initial subspaces and the preparation blocks fixed.}
We now keep $U_0,Y_0$ fixed, as required in the contraction
criterion. The first four blocks give the regularity and nondegeneracy
estimates for $R_4,v_4$ in Lemma~\ref{su:lem:preparation}.
We freeze these blocks as well. Thus
$R_4,v_4$ are fixed. Our first modification is the uncorrected
field $X^{[1]}$, which rotates only the later matrices:
\[
 X^{[1]}_j=0\quad(1\le j\le4),\qquad
X^{[1]}_j=[K,V_j]\quad(5\le j\le q).
\]
At the first active block this gives
\[
 \begin{aligned}
 D_{X^{[1]}}R_5&=[K,R_5]-g^2R_5[K,R_4]R_5,\qquad
 D_{X^{[1]}}v_5&=Kv_5-R_5Kv_4-g^2R_5[K,R_4]v_5.
 \end{aligned}
\]
The terms involving $[K, R_4]$ arise because $R_4$ no longer
rotates; even $Y_5$ therefore need not follow a pure unitary
rotation. The term $-R_5Kv_4$ arises because the incoming normal
$v_4$ is fixed.

\textbf{Correcting the first active block.}
To remove the terms involving $[K,R_4]$, define $X^{[2]}$ by
changing only the fifth component of $X^{[1]}$:
\[
 X^{[2]}_5=[K,V_5]-g^2[K,R_4],\qquad
 X^{[2]}_j=X^{[1]}_j\quad(j\ne5).
\]
Propagating this correction gives, for $5\le j\le q$,
\[
 \begin{aligned}
 D_{X^{[2]}}R_j&=[K,R_j],\qquad D_{X^{[2]}}v_j&=Kv_j-(R_j\cdots R_5)Kv_4.
 \end{aligned}
\]
Thus each $Y_j$ for $5\le j\le q$ follows the unitary rotation. The
variation of $v_j$, which specifies $U_j$ through $v_j^*x=0$, retains
the extra term caused by fixing $v_4$.

\textbf{The changes in $c$ and $h$ under $X^{[2]}$.}
Write
\[
 \xi_q=(R_q\cdots R_5)Kv_4,\qquad
 M_q=(I_W+g^2R_q^2)^{-1}.
\]
The common unitary rotation contributes nothing to $c$ or $h$.
After removing it, $R_q$ is fixed and the infinitesimal change
of $v_q$ is $-\xi_q$. With the first four blocks fixed,
\eqref{su:eq:graph-volume-ratio} expresses $e^{-2c}$ as a fixed
multiple of $v_q^*M_qv_q$. Including the rotation term in the
normal variation, we obtain
\[
 D_{X^{[2]}}c
 =\frac{\Re(v_q^*M_q\xi_q)}{v_q^*M_qv_q},\qquad
 D_{X^{[2]}}v_q=Kv_q-\xi_q.
\]
The ideal change $\xi_q=v_q$ would give
$D_{X^{[2]}}c=1$ and $D_{X^{[2]}}h=0$.

\textbf{Choosing $K$ at the penultimate block.}
Exact matching is impossible. We instead seek an approximate match
at $q-1$, choosing $K$ using only $V_1,\ldots,V_{q-1}$ and reserving
the last block for correcting the terminal motion. Here
\[
 D_{X^{[2]}}v_{q-1}=Kv_{q-1}-\xi_{q-1},\qquad
 \xi_{q-1}=(R_{q-1}\cdots R_5)Kv_4.
\]
We therefore seek a small skew-Hermitian $K$ such that
\[
 (R_{q-1}\cdots R_5)Kv_4\approx(R_{q-1}\cdots R_5)v_4=v_{q-1}.
\]
The error is measured relative to $\|v_{q-1}\|$. A small $Kv_4$
can reproduce this output if it matches $v_4$ on the few directions
that contribute substantially after propagation.

\textbf{Correcting the terminal motion.}
A scalar normalization of $K$ makes
$\Re(v_{q-1}^*\xi_{q-1})=\|v_{q-1}\|^2$.
The remaining matching error has no real radial component, so we
can choose a small skew-Hermitian $L$ satisfying
\[
 Lv_{q-1}=\xi_{q-1}-v_{q-1}.
\]
Define $X^{[3]}$ by changing only the last component:
\[
 X^{[3]}_q=X^{[2]}_q-[L,V_q-g^2R_{q-1}],\qquad
 X^{[3]}_j=X^{[2]}_j\quad(j<q).
\]
The terminal variations are
\[
 D_{X^{[3]}}R_q=[K-L,R_q],\qquad
 D_{X^{[3]}}v_q=(K-L)v_q-v_q.
\]
We take $X=X^{[3]}$. The unit radial contraction gives $D_Xc=1$.
Since the radial term $-v_q$ does not change the hyperplane
$v_q^*x=0$, both terminal subspaces undergo only the common
unitary rotation generated by $K-L$. Invariance of $h$ then
gives $D_Xh=0$ (Lemma~\ref{su:lem:field-geometry}).

It remains to achieve the match at $q-1$ with as small a
skew-Hermitian generator $K$ as possible.

\subsubsection{Constructing \texorpdfstring{$K$}{K}}
\label{sec:intro-small-generator}

Write $\mathcal R=R_{q-1}\cdots R_5$. We first look for a small
vector $y$ such that $\mathcal Ry\approx\mathcal Rv_4$, and then
realize $y$ as the image of $v_4$ under a small skew-Hermitian generator.
Skew-Hermiticity requires $\Re(v_4^*y)=0$. We will
impose the stronger condition $v_4^*y=0$.

\textbf{Isolating $V_5$ in the product.}
For two blocks we have
\[
 R_6R_5
 =(V_6-EI_W)^{-1}
   \bigl[R_5^{-1}-g^2(V_6-EI_W)^{-1}\bigr]^{-1}.
\]
The first factor depends only on $V_6$. To isolate $V_5$ in the
full product, define
\[
 \overleftarrow R_{q-1}=(V_{q-1}-EI_W)^{-1},\quad  
 \overleftarrow R_j
 =(V_j-EI_W-g^2\overleftarrow R_{j+1})^{-1},
 \qquad 6\le j<q-1.
\]
The same elimination gives
\begin{equation}\label{us:eq:cavity}
 \mathcal R=R_{q-1}\cdots R_5
 =\underbrace{\overleftarrow R_{q-1}\cdots\overleftarrow R_6}_{A}
  \underbrace{(V_5-EI_W-g^2R_4-g^2\overleftarrow R_6)^{-1}}_{\widetilde R}.
\end{equation}
The two shifts in $\widetilde R$ come from the blocks to the left
and right of block $5$, respectively (Section~\ref{su:sec:explicit-generator}).
Conditioning on all inputs except $V_5$ fixes both shifts, $A$, and
$v_4$, while $V_5$ retains its GUE law.

First take $W^{-\beta}\le g\le W^{-\alpha}$ with fixed
$0<\alpha<\beta<1/2$, so $g^2$ exceeds the bulk level spacing $1/W$
by a power of $W$. In this range, up to small powers of $W$, we need
to retain at most $Wg^2$ right singular directions of $A$.
For each such unit vector $u$,
\[
 u^*\widetilde Rv_4=(\widetilde Ru)^*v_4.
\]
For the fixed vector $v_4$, averaging over $V_5$ with the cutoff
weight $\chi^2$ bounds its squared projection onto this span by
$O(g^2)\|v_4\|^2$, up to small powers of $W$.
By choosing $y$ with the same projection, we can construct
a small skew-Hermitian $K$ such that $\mathcal R Kv_4\approx\mathcal Rv_4$
and $\E\|\chi K\|_{\HS}^2\lesssim W^\epsilon g^2$.

\textbf{Why $Wg^2$ directions suffice in one block.}
The factors of $A$ obey the same recursion from the right as $R_j$
does from the left. For two consecutive blocks we have
\[
 R_j=\Gamma_j^{-1},\qquad \Gamma_j=V_j-EI_W-g^2R_{j-1},
 \qquad \Gamma_{j+1}=V_{j+1}-EI_W-g^2R_j.
\]
If $\Gamma_j e=\lambda e$ and $\|e\|=1$, then $R_j e=\lambda^{-1}e$.
The next block therefore sees a shift $-g^2/\lambda$ along this direction.
This shift reaches the order-one scale of the GUE block when
$|\lambda|\asymp g^2$. Conditional on the preceding blocks,
$\Gamma_j$ is a GUE matrix with a fixed Hermitian shift, so the
one-block estimate \eqref{us:eq:wegner} gives
\[
 \E_{V_j}\#\{i:|\lambda_i(\Gamma_j)|\le g^2\}\lesssim Wg^2 .
\]
Block $j$ amplifies the squared norm in this direction by $1/\lambda^2$.
In block $j+1$, the shift $-g^2/\lambda$ along $e$ gives a spectral
density at zero, smoothed on scale $g^2$, of order
$[1+(g^2/\lambda)^2]^{-1}$ by \eqref{us:eq:relative-local-law} on
the regular chain. Multiplying this density by $O(g^2)$ bounds the
spectral mass of $e$ in the band $[-g^2,g^2]$ of the next block.
The squared amplification in block $j$ times the next-block density
is of order
\[
 \frac{1}{\lambda^2}\cdot\frac{1}{1+(g^2/\lambda)^2}
 =\frac{1}{\lambda^2+g^4}.
\]
Thus for $|\lambda|\lesssim g^2$ the weight per direction is of
order $g^{-4}$. For $k\ge0$ with $2^kg^2\ll1$, the band $|\lambda|\asymp 2^kg^2$
has about $2^k Wg^2$ directions. Its weight per direction is of
order $2^{-2k}$ relative to the band at $g^2$, so its
total contribution is only of order $2^{-k}$ times that of the
band at $g^2$.

The directions in the band $|\lambda|\asymp g^2$ change from block to
block. A product term passing through this band has rank at most
the rank of its spectral projection, hence at most $Wg^2$ up to
small powers of $W$ on the regular chain. For a vector fixed by the
preceding blocks, the band in the next block captures at most an
$O(g^2)$ fraction of its squared norm on average, up to the same
power loss, by \eqref{us:eq:projexp}. Inversion then amplifies this
part by order $g^{-4}$ in squared norm, giving a one-step gain in
mean squared norm of at most order $g^2\cdot g^{-4}=g^{-2}$.
The product moment estimate combines the contributions from all
spectral bands to bound $\Tr|A|^p$.

\textbf{A first choice of $K$.}
Recall that our target is $A\widetilde R K v_4\approx A\widetilde R v_4$.
We retain the right singular directions of $A$ with singular values
above $t_{\mathrm{cut}}$, and let $Q$ project orthogonally onto their
span after applying $\widetilde R$. To count these directions we use
\[
 \rank Q=\#\{i:\sigma_i(A)>t_{\mathrm{cut}}\}
 \le t_{\mathrm{cut}}^{-p}\Tr|A|^p .
\]
With $p$ as in \eqref{su:eq:packet2}, the moment bound
\eqref{su:eq:retained-rank-moment} gives the count \eqref{su:eq:rank}.
If $Qx=0$, then $\widetilde R x$
lies entirely in the discarded right singular subspace of $A$.
We choose $t_{\mathrm{cut}}$ to make the discarded output small
relative to $v_{q-1}=A\widetilde Rv_4$. On the cutoff support the
choice in \eqref{su:eq:construction-threshold} gives
\[
 \|A\widetilde R x\|\le t_{\mathrm{cut}}\|\widetilde R\|\,\|v_4\|
 \lesssim W^{-4}\|v_{q-1}\|
 \qquad (Qx=0,\ \|x\|\le\|v_4\|).
\]
Write $v_4=a+b$, with $a=Qv_4$
and $b=(I_W-Q)v_4$. The overlap estimate
\eqref{su:eq:construction-overlap-goals} controls $\|a\|^2/\|v_4\|^2$.
We might try to set
$Kv_4=a$. This is impossible for a skew-Hermitian $K$ when $a\ne0$,
since $v_4^*a=\|a\|^2$ whereas $\Re(v_4^*Kv_4)=0$.

When $\|a\|\ll\|v_4\|$, we have $\|b\|\asymp\|v_4\|$. We can
use this large component $b$ to produce $a$. As a zeroth choice, take
the rank-one map $K^{[0]}=ab^*/\|b\|^2$, so $K^{[0]}v_4=a$.
Subtracting its adjoint gives our first skew-Hermitian choice,
\[
 K^{[1]}=K^{[0]}-(K^{[0]})^*=\frac{ab^*-ba^*}{\|b\|^2},\qquad
 K^{[1]}v_4=a-\frac{\|a\|^2}{\|b\|^2}\,b .
\]
The correction lies along $b$, so $QK^{[1]}v_4=a$ exactly.
Since $a\perp b$, we have
$\|K^{[1]}\|_{\HS}^2=2\|a\|^2/\|b\|^2$.
Thus the small overlap gives a small generator, and
\eqref{su:eq:construction-error-goal} controls the relative error
from the discarded output.

\textbf{A spectral cutoff for $g^2R_4$.}
The Gaussian divergence at block $5$ contains the term
$-W\Tr(V_5[\chi K,g^2R_4])$ in \eqref{intro:eq:blockwise-divergence}.
To control the commutator, take $K^{[2]}=PK^{[1]}P$,
where $P$ is the spectral cutoff for $g^2R_4$ in \eqref{su:eq:source-projection}.
Here $Q$ selects the component we must match, whereas $P$ restricts
the directions on which $K$ may act. By construction,
\[
 K^{[1]}b=a,\qquad K^{[1]}a=-\frac{\|a\|^2}{\|b\|^2}\,b.
\]
Thus the component $a$ in $QK^{[1]}v_4$ comes entirely from
$K^{[1]}b$. Suppose, for example, that $Pa=a$ and $Pb=0$.
Applying $K^{[2]}$ then gives
\[
 v_4=a+b \ \xrightarrow{P}\ a
 \ \xrightarrow{K^{[1]}}\ -\frac{\|a\|^2}{\|b\|^2}\,b
 \ \xrightarrow{P}\ 0 .
\]
The cutoff preserves the target $a$ but removes the input $b$ that
the first generator needs to produce it.

For the third choice, $K^{[3]}$, we first choose its image
$y=K^{[3]}v_4$ inside $\ran P$ so that it still matches $a=Qv_4$.
In the range $W^{-\beta}\le g\le W^{-\alpha}$, up to small powers of
$W$, $Q$ retains at most $Wg^2$ input directions selected by the large
singular values of $A$, while $P$ discards at most $Wg^2$ directions
corresponding to large eigenvalues of $|g^2R_4|$. In the critical
window the same construction takes $P=I_W$ and uses the rank and
overlap bounds in \eqref{su:eq:construction-overlap-goals}.
In Section~\ref{su:sec:explicit-generator}
we construct $y$ by \eqref{us:eq:H0}. On the cutoff support $Pv_4\ne0$ and
\[
 Py=y,\qquad v_4^*y=0,\qquad Qy=Qv_4,\qquad \|y\|\le2\|Qv_4\|.
\]
As for $K^{[1]}$, take a rank-one map minus its adjoint,
\[
 K^{[3]}=\|Pv_4\|^{-2}
       \bigl(y(Pv_4)^*-(Pv_4)y^*\bigr).
\]
Both $y$ and $Pv_4$ lie in $\ran P$, so $PK^{[3]}P=K^{[3]}$.
The generator sends $v_4$ to $y$ and has squared norm
\[
 \|K^{[3]}\|_{\HS}^2
 =2\frac{\|y\|^2}{\|Pv_4\|^2}
 \lesssim\frac{\|Qv_4\|^2}{\|Pv_4\|^2}.
\]
Since $Q(y-v_4)=0$ and $\|y-v_4\|\lesssim\|v_4\|$, the cutoff bound
above gives $\|A\widetilde R(y-v_4)\|\lesssim W^{-4}\|v_{q-1}\|$.
Section~\ref{su:sec:explicit-generator} denotes $K^{[3]}$ by $K_0$.
To match the radial component at $q-1$ exactly, we use
\eqref{us:eq:normalize},
\[
 \xi=\frac{\Re\bigl(v_{q-1}^*\mathcal R K^{[3]}v_4\bigr)}
               {\|v_{q-1}\|^2}=1+O(W^{-4}),\qquad K=K^{[3]}/\xi .
\]
The terminal correction $L$ then absorbs the remaining error.

\textbf{Property of $K$.}
Recall the recursion for the block normals and the factorization
\eqref{us:eq:cavity},
\[
 R_j=(V_j-EI_W-g^2R_{j-1})^{-1},\quad v_j=R_jv_{j-1}\quad(j\ge2),
 \qquad \mathcal R=R_{q-1}\cdots R_5=A\widetilde R.
\]
At block $4$, $P=\mathbf1_{\{|g^2R_4|\le W^{\varepsilon/800}\}}$
selects the directions on which $K$ may act, with $P=I_W$ in the
critical window. The input directions seen by the retained output form
the range of $Q$, the orthogonal projection onto the image under
$\widetilde R$ of the right singular subspace of $A$ above $t_{\mathrm{cut}}$.
For fixed incoming data, $K=K(V_1,\ldots,V_{q-1})$. The smooth cutoff
$\chi$ from Section~\ref{sec:intro-contraction-criterion} selects inputs
for which $Pv_4\ne0$ and the discarded output is small. On its support,
with $\xi_{q-1}=\mathcal R Kv_4$ and $v_{q-1}=\mathcal Rv_4$, we have
\[
 \begin{aligned}
 \text{Generator}\quad & K^*=-K,\qquad K=PKP,\qquad \rank K\le2,\\
 \text{Size}\quad & \|K\|_{\HS}\lesssim \|Qv_4\|/\|Pv_4\|,\\
 \text{Block }5\quad & \|g^2[K,R_4]\|_{\HS}
     \le 2\|g^2PR_4P\|\,\|K\|_{\HS},\\
 \text{Matching error}\quad & \|\xi_{q-1}-v_{q-1}\|
     \lesssim W^{-4}\|v_{q-1}\|,\\
 \text{Radial component}\quad & \Re(v_{q-1}^*\xi_{q-1})=\|v_{q-1}\|^2.
 \end{aligned}
\]
With the terminal correction $Lv_{q-1}=\xi_{q-1}-v_{q-1}$, the field
satisfies $D_Xc=1$ and $D_Xh=0$ by Lemma~\ref{su:lem:field-geometry}.

\subsubsection{Controlling the Gaussian cost}
\label{sec:intro-small-gaussian-cost}

To estimate $\E|\delta_{\rm G}(\chi X)|^2$, write $\delta_{{\rm G},j}$
for the Gaussian divergence in $V_j$. The four preparation blocks do
not move, so
\[
 \delta_{\rm G}(\chi X)
 =\sum_{j=5}^{q}\delta_{{\rm G},j}(\chi X_j).
\]
Condition on the preparation blocks, so $R_4,v_4,P$ are fixed.
Put $\widetilde V_5=V_5-g^2R_4$ and $\widetilde V_j=V_j$ for
$6\le j<q$. For $5\le j<q$ we have $\chi X_j=[\chi K,\widetilde V_j]$,
and the Gaussian divergence formula \eqref{up:eq:gaussian-divergence} reads
\[
 \delta_{{\rm G},j}(\chi X_j)
 = W\Tr(V_j[\chi K,\widetilde V_j])
   -\operatorname{div}_{V_j}(\chi X_j).
\]
Varying $V_j$ in a Hermitian direction $B$ gives
\[
  D_{V_j}(\chi X_j)[B]
   =[D_{V_j}(\chi K)[B],\widetilde V_j]+[\chi K,B].
\]
To take the ordinary divergence, put $B=E_\alpha$ and sum the
$E_\alpha$ components over a real orthonormal basis of Hermitian
matrices, using $\langle A,B\rangle_{\HS}=\Re\Tr(A^*B)$.
The second commutator drops out since
$\langle E_\alpha,[\chi K,E_\alpha]\rangle_{\HS}=0$.
In the first, expand $D_{V_j}(\chi K)[E_\alpha]$ in a real orthonormal
basis $(H_k)$ of the skew-Hermitian matrices supported on $\ran P$.
The sum over $\alpha$ then reconstructs the direction
$\sum_\alpha\langle E_\alpha,[H_k,\widetilde V_j]\rangle_{\HS}E_\alpha
 =[H_k,\widetilde V_j]$. Since $\Tr(V_j[\chi K,V_j])=0$, this gives
\begin{equation}\label{intro:eq:blockwise-divergence}
 \delta_{{\rm G},j}(\chi X_j)
 =-\mathbf1_{\{j=5\}}Wg^2\Tr(V_5[\chi K,R_4])
  -\sum_k\big\langle H_k,
      D_{V_j}(\chi K)\bigl[[H_k,\widetilde V_j]\bigr]\big\rangle_{\HS}.
\end{equation}
\textbf{The shift at block $5$.}
The coefficient $[\chi K,g^2R_4]$ in the trace in
\eqref{intro:eq:blockwise-divergence} depends on $V_5$.
Let $M_0$ bound $\|g^2PR_4P\|$ on the cutoff support. Since $P$
commutes with $R_4$ and $\chi K=P(\chi K)P$, we have
\[
 W\|[\chi K,g^2R_4]\|_{\HS}^2
 \le 4WM_0^2\|\chi K\|_{\HS}^2.
\]
Write $\delta_{\rm act}=\sum_{j=5}^{q-1}\delta_{{\rm G},j}(\chi X_j)$.
Two Gaussian integrations by parts for $\E|\delta_{\rm act}|^2$ give
\eqref{su:eq:active-divergence-decomposition}. There the shift costs at most
$CW(1+M_0^2)\E\|\chi K\|_{\HS}^2+CWM_0e^{-cW}$, as in
\eqref{su:eq:active-matrix-derivative-bound}.

\textbf{Summing the blocks for each rotation direction.}
For the derivative contribution to $\delta_{\rm act}$, exchange the
$j$ and $k$ sums in \eqref{intro:eq:blockwise-divergence} using
\[
 \begin{aligned}
  \sum_{j=5}^{q-1}\sum_k\big\langle H_k,
        D_{V_j}(\chi K)\bigl[[H_k,\widetilde V_j]\bigr]\big\rangle_{\HS}
  &=\sum_k\left\langle H_k,
        \sum_{j=5}^{q-1}D_{V_j}(\chi K)\bigl[[H_k,\widetilde V_j]\bigr]
                          \right\rangle_{\HS}.
 \end{aligned}
\]

The inner sum is the derivative $D_{\rm mat}^{H_k}(\chi K)$ obtained
by conjugating every $\widetilde V_j$ by $\mathcal U_t=e^{tH_k}$
while holding $v_4$ fixed.
 
\textbf{Another key property of $K$.}
Under this rotation the matching target becomes $Q_t v_4$, where
$Q_t=\mathcal U_t Q\mathcal U_t^*$. If we also rotate $v_4$, as in
\eqref{su:eq:cut-generator-equivariance}, the rotated vector $\mathcal U_t y$ matches the new target since
\[
 Q_t(\mathcal U_t y)=\mathcal U_t Qv_4=Q_t(\mathcal U_t v_4).
\]
Both $y$ and $Pv_4$ in the rank-two formula rotate by $\mathcal U_t$,
while $\xi$ and $\chi$ stay fixed. Thus $\chi K$ becomes
$\mathcal U_t(\chi K)\mathcal U_t^*$. Differentiating at $t=0$ and
subtracting the contribution from $v_4$ gives the derivative in the inner sum,
\[
 \underbrace{D_{\rm mat}^{H_k}(\chi K)}_{\text{matrices only}}
 =\underbrace{[H_k,\chi K]}_{\text{matrices and vector}}
  -\underbrace{D_{v_4}(\chi K)[H_kv_4]}_{\text{vector only}}.
\]
Here $D_{v_4}$ varies only the vector argument in the construction of
$\chi K$, with all matrix inputs fixed.
In the scalar divergence, $\langle H_k,[H_k,\chi K]\rangle_{\HS}=0$
by cyclicity. Substitution therefore gives the exact identity
\[
 \delta_{\rm act}
 =-Wg^2\Tr(V_5[\chi K,R_4])
   +\sum_k\big\langle H_k,D_{v_4}(\chi K)[H_kv_4]\big\rangle_{\HS}.
\]
In the second-moment bound \eqref{su:eq:active-matrix-derivative-bound},
substitute the same derivative identity before estimating the squared norms.
By \eqref{su:eq:commutator-basis-bound}, the squared commutators sum to
at most $2W\|\chi K\|_{\HS}^2$, which is absorbed into the size term below.
For the vector derivatives, $H_kv_4=H_k(Pv_4)$. In coordinates with
$Pv_4$ along the first axis, only the first column of $H_k$ moves the
vector. Summing over $k$ as in \eqref{su:eq:matrix-source-derivative-comparison}
therefore gives
\[
 \begin{aligned}
 \E|\delta_{\rm act}|^2
 &\lesssim \E\!\left[\|Pv_4\|^2\sum_a
       \|D_{v_4}(\chi K)[e_a]\|_{\HS}^2\right] +W(1+M_0^2)\E\|\chi K\|_{\HS}^2+WM_0e^{-cW},
 \end{aligned}
\]
where $(e_a)$ is a real orthonormal basis of the tangent directions
$\{w\in\ran P:\Re(v_4^*w)=0\}$. Only the first term requires a derivative.

\textbf{Variations of $v_4$ that change or preserve $Qv_4$.}
Vary $v_4$ in a unit tangent direction $w$, keeping the matrices fixed.
Writing $\dot y=D_{v_4}y[w]$ and differentiating $Qy=Qv_4$ and $v_4^*y=0$ gives
\[
 Q\dot y=Qw,\qquad v_4^*\dot y=-w^*y.
\]
If $Qw$ is of order one, $y$ must follow the moving target at order-one
speed, however small $y$ itself is. If $Qw=0$, the target stays put.
The adjustment needed to keep $y$ orthogonal to $v_4$ is driven by
$w^*y$, which already carries the small factor $\|y\|$. For our rank-two
generator this gives $\|Pv_4\|\|D_{v_4}K[w]\|_{\HS}\lesssim\|K\|_{\HS}$
when both $Qw=0$ and $D_{v_4}\xi[w]=0$, by \eqref{su:eq:pair-source-regularity}.

Even when $Qw=0$, the scalar normalization $\xi$ can change.
Choose the tangent basis adapted to $Qw=0$ and $D_{v_4}\xi[w]=0$.
The first condition imposes at most $2\rank Q$ real constraints,
and the second adds at most one. Thus at most $2\rank Q+1$ basis directions
need the order-one bound. On the others the mean squared cost is bounded
by $C\E\|\chi K\|_{\HS}^2$. Differentiating $\chi$ adds only $O(1)$
in total by \eqref{su:eq:pair-cutoff-derivative}. For $W^{-\beta}\le g\le W^{-\alpha}$,
suppressing small powers of $W$, the construction gives
$\rank Q\lesssim Wg^2$ on the cutoff support and $\E\|\chi K\|_{\HS}^2\lesssim g^2$.
Thus the derivative bound in \eqref{su:eq:pair-source-energy} comes from
\[
 \underbrace{O(Wg^2)\times O(1)}_{\text{target or normalization changes}}
 \;+\;
 \underbrace{O(W)\times O(g^2)}_{\text{remaining directions, cost controlled by }\E\|\chi K\|_{\HS}^2}
 =O(Wg^2).
\]

\textbf{The resulting contraction scale.}
Conditional on $V_1,\ldots,V_{q-1}$, the last block contributes a
centered Gaussian with variance $Wg^4\|[\chi L,R_{q-1}]\|_{\HS}^2$
by \eqref{su:eq:last-divergence-moments}. Its cross term with
$\delta_{\rm act}$ vanishes. The $W^{-4}$ matching error gives
$\|L\|_{\HS}\lesssim W^{-4}$ on the cutoff support, where also
$\|R_{q-1}\|\lesssim W^3$, so this block costs only $O(W^{-1})$.

Averaging over the preparation blocks and restoring the small powers of
$W$, Lemma~\ref{up:lem:small-contraction} gives
$\E|\delta_{\rm G}(\chi X)|^2\lesssim\mathfrak b$ with
\[
 \mathfrak b=
 \begin{cases}
 W^{\varepsilon/2}Wg^2,& W^{-\beta}\le g\le W^{-\alpha},\\
 W^{18\theta},& c_*W^{-1/2}\le g\le W^{-1/2+\theta/2}.
 \end{cases}
\]
Here $0<\alpha<\beta<1/2$, $\varepsilon,c_*>0$ and $0<\theta\le1/100$
are fixed. Since $\E\chi\ge\tfrac12$, the criterion in
Section~\ref{sec:intro-contraction-criterion} gives the angular $L^4$
contraction of $P_{1/2}$ by $1-c/\mathfrak b$ on bounded nonnegative
functions invariant under simultaneous unitary rotations, uniformly
over incoming flags and energies in a fixed bulk interval. Iteration and smoothing as in
Section~\ref{sec:intro-operator-contraction} then give \eqref{DnWDWD} at length
$n_0\lesssim\mathfrak b\log W$.

\subsection{Perturbations from a fixed matrix direction}
\label{sec:intro-fixed-fields}

The range $g=W^{o(1)}$ lies between two regimes in which stability
has different origins. At small $g$, the product strongly amplifies
directions selected by the random blocks, but averaging over $V_5$
controls their overlap with the incoming normal. At large $g$, each
transfer is a small perturbation of the free transfer, so a fixed-length
product stays close to the free evolution.

It may be possible to obtain the main bound using just these two
approaches. Extending Section~\ref{sec:small-upper-proof} to $g=W^{o(1)}$
would require controlling the matching error and the final correction
even when $g\asymp1$. Extending Section~\ref{sec:large-upper-proof} to
the full bulk for $W^\epsilon\le g<C_\kappa\sqrt W$ would require a
drift estimate uniform in the energies outside its near-zero window.
For simplicity, we treat both ranges by adapting Drogin's construction
for the GOE model \cite[Section~4]{Drogin} to our GUE blocks. The Gaussian cost of realizing
the prescribed stretch is at most $W^{1+o(1)}$ for $g=W^{o(1)}$.
Away from zero energy we bound the cost by $O(Wg^2)$ for
$W^\epsilon\le g<W$, with fixed $0<\epsilon<1$.

We switch from the notation of
Section~\ref{sec:small-upper-proof} to that of Section~\ref{sec:compact-proof},
in which $P$ and $R_i$ denote blocks of test matrices. We fix $U$
pointwise and stretch the unit
direction $e\in Y\cap U^\perp$ by $e^t$ before applying the transfers.
This leaves both terminal subspaces fixed and adds $t$ to their
logarithmic volume ratio.

Write $F=T_1\cdots T_q$ in the product order of
Section~\ref{sec:compact-proof}. In the coordinates of $F$, the stretch
is $e^{\sqrt2 tM}$ with
\[
 M=\tfrac1{\sqrt2}\,S\bigl(ee^*-(J_{\rm sig}e)(J_{\rm sig}e)^*\bigr)S^*,
 \qquad J_{\rm sig}=\operatorname{diag}(I_W,-I_W).
\]

We seek Hermitian perturbations $X=(X_1,\ldots,X_q)$ of the GUE blocks
that realize the stretch to first order, so we require
\[
 \left.\frac{d}{dt}F(V_1+tX_1,\ldots,V_q+tX_q)\right|_{t=0}
 =\sqrt2\,FM .
\]
Multiplying on the left by $F^{-1}$ and using the product rule gives
the equivalent equation
\[
 \mathcal A(X):=F^{-1}D_XF
 =\frac1g\sum_{i=1}^q Q_i^{-1}
     \begin{pmatrix}0&0\\ X_i&0\end{pmatrix}Q_i=\sqrt2 M,
 \qquad Q_i=T_{i+1}\cdots T_q,\quad Q_q=I .
\]

\subsubsection{The general perturbation and its Gaussian cost}

Taking the real Hilbert--Schmidt inner product of both sides with a
$2W\times 2W$ test matrix $Z$ in the transfer Lie algebra, we need
\[
 \langle X,\mathcal A^*Z\rangle=\langle\sqrt2 M,Z\rangle
 \qquad\text{for every }Z.
\]
If we can recover $Z=\mathcal K_0\mathbf R$ from the data
$\mathbf R=(R_1,\ldots,R_q)=g\mathcal A^*Z$, then the right side equals
$\langle\sqrt2\,g\mathcal K_0^*M,\mathcal A^*Z\rangle$.
Thus $X=\sqrt2\,g\mathcal K_0^*M$ solves the equation.

To bound its size, Section~\ref{sec:compact-proof} uses Drogin's block
elimination \cite[Proposition~4.4]{Drogin} on the support of the
smoothing cutoff $\chi_P$ for $G_{1/2}$ from Section~\ref{sec:intro-operator-contraction}.
Subtracting two block equations gives the
coefficient $(V_{k_j}-V_{i_j})/g$. Normalizing this difference costs
a factor $g$. The GUE bound \eqref{eq:c-main-coercivity}, obtained
using \cite{CGV}, then controls the inverse uniformly in $W$. For $g\ge1$
this gives
\[
 \|\mathcal K_0\|\lesssim g,\qquad
 \|X\|\le\sqrt2\,g\|\mathcal K_0\|\lesssim g^2,
\]
since $\|M\|_{\rm HS}=1$.

\medskip\noindent\textbf{The Gaussian cost.}
The divergence also depends on how $X$ varies with the onsite blocks.
For any fixed $B>0$, Section~\ref{sec:compact-proof} obtains
throughout $W^{-B}\le g\le W^B$
\[
 \|X\|\lesssim\max\{g^{-1},g^2\},\qquad
 \|\nabla_V X\|_{\rm HS}\lesssim\sqrt W\max\{g^{-1},g^2\}
 \quad\text{on }\supp\chi_P .
\]
Here $\|X\|^2=\sum_i\|X_i\|_{\rm HS}^2$. The factor $\sqrt W$ comes
from summing squared derivatives over a Hilbert--Schmidt orthonormal
basis of Hermitian matrices, as in \eqref{eq:c-completeness}.
Substituting the bounds for $X$ and $\nabla_V X$ in the Gaussian
integration-by-parts estimate gives
\[
 \E|\delta_{\rm G}(\chi_P X)|^2
 \le W\E\|\chi_P X\|^2+\E\|\nabla_V(\chi_P X)\|_{\rm HS}^2
 \lesssim W\max\{g^{-2},g^4\}.
\]
This is the divergence bound of Lemma~\ref{up:lem:density-contraction}\textup{(i)}.
Since $\E\chi_P\ge1/2$ and the stretch fixes the terminal flag with
$D_Xc=1$, the criterion of Section~\ref{sec:intro-contraction-criterion}
(Lemma~\ref{up:lem:angular}) gives contraction by $1-c/\mathfrak b$
with $\mathfrak b=W\max\{g^{-2},g^4\}$. For intermediate coupling
$g=W^{o(1)}$, the ratio of this scale to the target $Wg^2$ is
$\max\{g^{-4},g^2\}=W^{o(1)}$. This loss fits within the arbitrarily
small power of $W$ allowed in the upper bound.

\subsubsection{Saving a factor of \texorpdfstring{$g$}{g} away from zero energy}

For $W^\epsilon\le g<W$, the general perturbation costs $Wg^4$ and we
need $Wg^2$. Retain one block equation \eqref{eq:c-block-equations}
before subtraction. With $E$ the energy and $P$ the upper-left block
of $Z$, it reads
\[
 E(P+P^*)=F_i(\mathbf R)+g^{-1}(V_iP+P^*V_i).
\]
On the cutoff, $\|F_i(\mathbf R)\|_{\rm HS}\lesssim\|\mathbf R\|$ and
$\|V_i\|\lesssim1$. The factor $g^{-1}$ therefore lets us use the
coarse bound $\|P\|_{\rm HS}\lesssim g\|\mathbf R\|$ to obtain
\[
 |E|\,\|P+P^*\|_{\rm HS}
 \lesssim \|\mathbf R\|+g^{-1}\|P\|_{\rm HS}
 \lesssim \|\mathbf R\|.
\]
For $|E|\ge\kappa_1>0$, division by $|E|$ saves a factor of $g$ on
$P+P^*$. This is exactly the diagonal combination tested by the
Hermitian target $M$ in $\langle M,Z\rangle$. Section~\ref{sec:compact-proof}
bounds the off-diagonal blocks by back substitution and estimates the
derivative of the minimum-norm perturbation. With
$\mathbf R=g\mathcal A^*Z$, the resulting size and Gaussian cost are
\[
 \|X\|\lesssim_{\kappa_1}g,\qquad
 \E|\delta_{\rm G}(\chi_P X)|^2\lesssim_{\kappa_1}Wg^2 .
\]
This gives the divergence input of Lemma~\ref{up:lem:density-contraction}\textup{(ii)}
at scale $\mathfrak b=Wg^2$, hence \eqref{DnWDWD} by contraction,
iteration and smoothing. Near zero, Section~\ref{sec:intro-direct-transfer}
estimates transfer growth directly.

\subsection{Direct transfer growth}
\label{sec:intro-direct-transfer}

We need a new method for two reasons. First, near $E=0$, the estimate
on $E(P+P^*)$ in Section~\ref{sec:intro-fixed-fields} gives no uniform
control of $P+P^*$. The general perturbation then costs $Wg^4$, too
large for the target scale $Wg^2$. Second, the smoothing estimate in
Section~\ref{sec:intro-operator-contraction} assumes $g\le W^B$ for a
fixed $B$, whereas our upper bound must also cover superpolynomial $g$.

In both ranges we treat the random term $V_j/g$ as a perturbation of
the free rotation, which preserves lengths. The larger $g$ is, the
closer each step is to this rotation and the more accurate the
expansion in $g^{-1}$ becomes. As $g$ decreases, the drift we seek
grows like $(Wg^2)^{-1}$, but the second-order error bound $O(g^{-4})$
grows faster. Its size relative to the drift is $O(W/g^2)$. Thus for
$g\ge W$ second order suffices even at $E=0$. The harder case is near
zero energy when $W^\epsilon\le g\ll\sqrt W$. There the error can mask
the drift, so we must expand further.  

\medskip\noindent\textbf{Following one direction through the product.}
For each fixed energy and initial frame, our targets are the two
bounds on $\|D_n^{-1}\|$ recalled below, uniform in both parameters.
\[
 \begin{aligned}
 \Prob\bigl(\|D_{n_0}^{-1}\|>W^{-D}\bigr)&\le W^{-D}
    &&\text{at } n_0\asymp Wg^2\log W,\\
 \E\|D_n^{-1}\|^s&\le CW e^{-cn/(Wg^2)}
    &&\text{for every } n\ge1.
 \end{aligned}
\]
The first is \eqref{DnWDWD} near zero energy for $W^\epsilon\le g<W$
and fixed $0<\epsilon<1$.
The second is \eqref{intro:eq:negative-transfer} for a fixed $s\in(0,1)$
and applies throughout the bulk for $g\ge W$, including at $E=0$.
Roughly speaking, we estimate the weakest growth of $D_n$ directly by
finding a gain of order $(Wg^2)^{-1}$ per block in the logarithmic
volume ratio.

\medskip\noindent\textbf{How one block moves the normal.}
As in Section~\ref{sec:intro-small-fields}, we follow a hyperplane
$U\subset Y=\ran\Phi$ through the transfer product, with
$\Phi\in\mathcal F_W$ and $\dim U=W-1$. In the free coordinates of
Definition~\ref{def:transfer-matrices}, the two spaces after $j$ blocks are
\[
 Y_j=S_E^{-1}T_j(E)\cdots T_1(E)S Y,\qquad
 U_j=S_E^{-1}T_j(E)\cdots T_1(E)S U,
\]
with the empty product at $j=0$. Section~\ref{sec:intro-small-fields}
used a Hermitian graph matrix $R_j$ and an unnormalized normal. Here
$C_j$ is unitary, and $v_j$ is a unit normal in the upper coordinate:
\[
 \begin{aligned}
 Y_j&=\{(u,C_ju):u\in\C^W\},\qquad C_j^*C_j=I_W,\\
 U_j&=\{(u,C_ju):u\in\C^W,\ v_j^*u=0\},\qquad \|v_j\|=1.
 \end{aligned}
\]
Write $E=2\sin\theta$, $h_\theta=(2\cos\theta)^{-1}$,
and $\lambda=-ie^{i\theta}$. The onsite block $V_j$ updates the pair
$(C_{j-1},v_{j-1})$ by
\[
 \begin{aligned}
 C_j&=-e^{-2i\theta}\bigl[C_{j-1}+\tfrac{ih_\theta}{g}V_j(I_W+C_{j-1})\bigr]F_j^{-1},\\
 v_j&=\lambda\frac{F_j^{-*}v_{j-1}}{\|F_j^{-*}v_{j-1}\|},
 \qquad F_j=I_W-\tfrac{ih_\theta}{g}V_j(I_W+C_{j-1}),
 \end{aligned}
\]
where $F_j^{-*}=(F_j^{-1})^*$. The upper coordinate evolves by
$u'=\lambda F_j u$. The normal transforms by the inverse adjoint because
\[
 v_{j-1}^*u=0
 \quad\Longleftrightarrow\quad
 v_{j-1}^*(\lambda F_j)^{-1}u'=0
 \quad\Longleftrightarrow\quad
 (\lambda F_j^{-*}v_{j-1})^*u'=0.
\]
For $g\gg1$, the term involving $V_j/g$
perturbs the free rotation $u\mapsto\lambda u$. After $n$ blocks the
transfer on the upper coordinate is $\lambda^nF_n\cdots F_1$. Since
$|\lambda|=1$, the coordinate comparison \eqref{lg:eq:coordinate-comparison}
gives
\[
 \|D_n(\Phi;E)^{-1}\|\asymp_\kappa\|(F_n\cdots F_1)^{-*}\|.
\]
The one-step logarithmic volume ratio in \eqref{eq:c-volume-cocycle}
is $c_j:=-\log\|F_j^{-*}v_{j-1}\|$. Starting from a unit normal
$v_0=w$ and multiplying the normalization factors gives
\[
 \ell_n(C_0,w):=\sum_{j=1}^n c_j\big|_{v_0=w}
 =-\log\|(F_n\cdots F_1)^{-*}w\|.
\]
The weakest growth is the minimum of this sum over all initial unit
normals. Thus, for each realization,
\[
 \log\|D_n(\Phi;E)^{-1}\|
 =-\min_{\|w\|=1}\sum_{j=1}^n c_j\big|_{v_0=w}+O_\kappa(1).
\]

\medskip\noindent\textbf{Martingalization of $c_i$.}
Write $\E_{i-1}$ for expectation conditional on
$V_1,\ldots,V_{i-1}$. Centering each increment gives
\[
 \sum_{i=1}^n c_i
 =\sum_{i=1}^n\E_{i-1}c_i
 +\underbrace{\sum_{i=1}^n(c_i-\E_{i-1}c_i)}_{\text{martingale}}.
\]
The martingale enters the final estimates as an error term
(Lemmas~\ref{uvg:lem:mgf} and~\ref{lg:lem:fluctuations}).

\medskip\noindent\textbf{Where the positive drift is hidden.}
For the incoming state, write
\[
 (C,v)=(C_{i-1},v_{i-1}),\qquad
 m_1:=W^{-1}\Tr C,\qquad f_1:=v^*Cv.
\]
The logarithmic gain from block $i$ is
\[
 c_i=-\tfrac12\log\bigl[v^*(F_i^*F_i)^{-1}v\bigr],\qquad
 F_i=I_W-\tfrac{ih_\theta}{g}V_i(I_W+C).
\]
The linear term is $(h_\theta/g)\Im(v^*V_iCv)$. It has mean zero,
but expanding the logarithm to second order brings in its square.
The GUE covariance gives
\[
 \E_{i-1}\!\left[\left(\tfrac{h_\theta}{g}\Im(v^*V_iCv)\right)^2\right]
 =\frac{h_\theta^2}{2Wg^2}\bigl(1-\Re(f_1^2)\bigr).
\]
The constant part of this variance is the gain we seek.
The full expansion \eqref{lg:eq:c2} gives,
uniformly in the incoming state,
\[
 \E_{i-1}c_i
 =\underbrace{\frac{h_\theta^2}{2Wg^2}}_{\text{desired gain}}
 +\frac{h_\theta^2}{g^2}
 \left[\Re f_1+\Re(m_1f_1)-\frac{\Re(f_1^2)}{2W}\right]+O(g^{-4}).
\]
At a single block, the first term carries a factor $1/W$, while the
terms involving $\Re f_1$ and $\Re(m_1f_1)$ can be of order $g^{-2}$.
Over $n$ blocks the first term has the same positive value each time,
so it adds up to $nh_\theta^2/(2Wg^2)$. The other terms change with
the state. For example, at $E=0$ the free rotation sends $f_1$ to
$-f_1$, so consecutive linear terms cancel in the free evolution.
The correction below uses this cancellation and controls the quadratic
terms to recover a positive drift on the scale $n/(Wg^2)$.

\medskip\noindent\textbf{Correcting the drift at a bounded endpoint cost.}
Choose a bounded function $Q$ of the state. Before taking conditional
expectations, replace each increment $c_i$ by
\[
 Y_i:=c_i-(Q_i-Q_{i-1}),\qquad Q_i:=Q(C_i,v_i).
\]
The differences $Q_i-Q_{i-1}$ telescope. Applying the martingale
decomposition to $Y_i$ therefore gives
\[
 \sum_{i=1}^n c_i
 =\underbrace{Q_n-Q_0}_{\text{endpoint difference}}
 +\underbrace{\sum_{i=1}^n\E_{i-1}Y_i}_{\text{corrected drift}}
 +\underbrace{\sum_{i=1}^n(Y_i-\E_{i-1}Y_i)}_{\text{martingale}}.
\]
The endpoint difference costs at most $2\|Q\|_\infty$ for any $n$.
For the choices of $Q$ below, the martingale again contributes an error term.
To see how $Q$ changes the drift, write
\[
 \E_{i-1}Y_i=\E_{i-1}c_i-\E_{i-1}Q(C_i,v_i)+Q(C_{i-1},v_{i-1}).
\]
For a bounded test function $F$ of the state, define the transition
operator $P$ by
\[
 (P\circ F)(C_{i-1},v_{i-1})=\E_{i-1}F(C_i,v_i).
\]
The same operator applies at every step because the update rule is fixed
and the onsite blocks $V_i$ are independent and identically distributed.
Thus
\[
 \E_{i-1}Y_i=\E_{i-1}c_i-[(P-I)\circ Q](C_{i-1},v_{i-1}).
\]
We choose $Q$ so that its expected change absorbs the state-dependent
terms in $\E_{i-1}c_i$ and reveals the positive drift.
This use of a correction follows
Schulz-Baldes \cite[Lemma~1 and Section~6.1]{SBanomaly} and
Sadel--Schulz-Baldes \cite[Sections~3--5]{SSB}.
In both ranges, we seek a bounded correction $Q$ such that
$\E_{i-1}Y_i\ge c_\kappa/(Wg^2)$ for most steps $i$.
Near zero energy, for $W^\epsilon\le g<W$, this bound holds with high
probability simultaneously for $n_*<i\le n_0$, where the initial
interval has length $n_*\asymp g^2\log W$ and the target length is
$n_0\asymp Wg^2\log W$.

\subsubsection{A quadratic drift correction for \texorpdfstring{$g\ge W$}{g at least W}}

Suppose $g\ge W$ and fix $\kappa\in(0,2)$. Throughout the bulk
$|E|\le 2-\kappa$, the drift bound holds for every incoming state
from the first block onward. Together with the fluctuation estimate, this gives
the negative-moment bound \eqref{intro:eq:negative-transfer} at every
length $n\ge1$. In this range the $O(g^{-4})$ remainder in the drift
expansion is already smaller than the target scale $(Wg^2)^{-1}$.

\medskip\noindent\textbf{Using the free rotation to cancel the linear term.}
The random block enters the update only through $V_i/g$. For the
smooth functions $F$ used here, we can therefore compute $P\circ F$
to the required order by expanding about the free step and averaging
over $V_i$. The term linear in $V_i$ has mean zero.
For example, the linear term in the conditional drift \eqref{lg:eq:c2}
is $h_\theta^2g^{-2}\Re f_1$. At energy $E=2\sin\theta$ the free step
multiplies $f_1$ by $-e^{-2i\theta}$, so its change is
$-(1+e^{-2i\theta})f_1$. This gives the first correction
\[
 Q^{[1]}(C,v)=-h_\theta^2g^{-2}\Re\bigl[f_1/(1+e^{-2i\theta})\bigr],\qquad
 (P-I)\circ Q^{[1]}=h_\theta^2g^{-2}\Re f_1+O(g^{-4}).
\]
 
\medskip\noindent\textbf{Using damping to absorb the two quadratic terms.}
For simplicity we describe the quadratic correction at $E=0$, where
$h_0=1/2$. Section~\ref{sec:uniform-large-g} gives the calculation
for general bulk energies with $h_\theta$ and the free multiplier
$-e^{-2i\theta}$. After the linear correction, the state-dependent part at order
$g^{-2}$ is
\[
 (4g^2)^{-1}\bigl[\Re(m_1f_1)-(2W)^{-1}\Re(f_1^2)\bigr].
\]
At $V_i=0$ the free step gives
\[
 (C_i,v_i)=(-C_{i-1},-iv_{i-1}),\qquad (m_1,f_1)\longmapsto(-m_1,-f_1).
\]
Thus $m_1f_1\mapsto(-m_1)(-f_1)=m_1f_1$ and
$f_1^2\mapsto(-f_1)^2=f_1^2$. For the random step, the terms linear
in $V_i/g$ average to zero, so
\[
 \begin{aligned}
 P\circ m_1&=-m_1+O(g^{-2}),&
 P\circ f_1&=-f_1+O(g^{-2}),\\
 P\circ(m_1f_1)&=m_1f_1+O(g^{-2}),&
 P\circ(f_1^2)&=f_1^2+O(g^{-2}).
 \end{aligned}
\]
For $f_1$, dividing by the change in the free multiplier, $-1-1=-2$,
gave $Q^{[1]}$. For the two quadratic terms this divisor is $1-1=0$.
Since we seek a lower bound on the drift, it suffices to offset their
negative contribution. For this we use the corrected moments
\eqref{uvg:eq:M}--\eqref{uvg:eq:F},
\[
 \begin{aligned}
 M_1(C)&=m_1-(8g^2)^{-1}\bigl(1+m_2+2m_1^2\bigr),\quad m_2=W^{-1}\Tr C^2,\\
 F_1(C,v)&=f_1-(8g^2)^{-1}\bigl(4m_1f_1+2f_2-f_1^2-|f_1|^2\bigr),\quad f_2=v^*C^2v.
 \end{aligned}
\]
For the original moments, the $O(g^{-2})$ errors in the one-step
averages are as large as the damping we seek. The corrections
$M_1=m_1+O(g^{-2})$ and $F_1=f_1+O(g^{-2})$ give the more precise
bounds \eqref{uvg:eq:Mdrift}--\eqref{uvg:eq:Fdrift},
\[
 \begin{aligned}
 -(P-I)\circ |M_1|^2
 &\ge g^{-2}\bigl[c_0|M_1|^2-C(W^{-2}+g^{-4})\bigr],\\
 -(P-I)\circ |F_1|^2
 &\ge g^{-2}\bigl[c_0|F_1|^2-C|M_1|^2-C(W^{-1}+g^{-4})\bigr].
 \end{aligned}
\]

To apply these bounds to the quadratic drift, replace $m_1,f_1$ by
$M_1,F_1$. Since all four moments are uniformly bounded,
\[
 m_1f_1-M_1F_1=(m_1-M_1)f_1+M_1(f_1-F_1)=O(g^{-2}),\qquad
 f_1^2-F_1^2=O(g^{-2}).
\]
The prefactor $(4g^2)^{-1}$ makes the resulting drift error $O(g^{-4})$.
Take $$Q=Q^{[1]}+A|M_1|^2+B|F_1|^2$$ First choose $B>0$ small enough
that the error $CB/(Wg^2)$ uses at most half of the gain $1/(8Wg^2)$.
By Young's inequality the product $|M_1F_1|$ is bounded by
$\delta |F_1|^2+(4\delta)^{-1}|M_1|^2$. Taking $\delta$ small compared
with $B$ and then $A$ large, the two damping bounds absorb the
quadratic loss at scale $g^{-2}$ and leave
\[
 \E_{i-1}Y_i
 \ge c_1/(Wg^2)-Cg^{-4}-C/(W^2g^2)
 \ge c/(Wg^2).
\]
Since $g\ge W$, we have $W/g^2\le 1/W$, so both errors are at most
$C/W$ times $(Wg^2)^{-1}$. For large $W$ the lower bound therefore
holds for every incoming state from the first block onward.

Lemma~\ref{uvg:lem:mgf} combines this drift with the fluctuation bound to give
$\E_{i-1}e^{-sY_i}\le e^{-c/(Wg^2)}$ for some fixed $s\in(0,1)$.
Iterating the conditional bound and using the bounded endpoint difference
and the basis estimate \eqref{lg:eq:scalar-reduction} gives
$\E\|D_n(\Phi;E)^{-1}\|^s\le CW e^{-cn/(Wg^2)}$ for every $n\ge1$,
as required in \eqref{intro:eq:negative-transfer}.

\subsubsection{A higher-order correction near zero energy}

Fix $0<\epsilon<1$ and suppose $W^\epsilon\le g<W$. We explain the
argument at $E=0$, where $h_0=1/2$ and the free step sends $C$ to $-C$.
Lemma~\ref{up:lem:center-transfer} proves the bound uniformly on a
fixed interval around zero. Here, for each fixed $D>0$, the target is
\[
 \sup_{\Phi\in\mathcal F_W}
 \Prob\bigl(\|D_{n_0}(\Phi;0)^{-1}\|>W^{-D}\bigr)\le W^{-D},
 \qquad n_0=\lceil C_1Wg^2\log W\rceil,
\]
for sufficiently large fixed $C_1$ and all large $W$.

\medskip\noindent\textbf{The initial interval.}
For $g\ge W$, the corrected drift bound in the preceding subsection
holds for every incoming state from the first block onward. Here the
early increments $c_i$ can have either sign, and losses can accumulate.
In the moment estimate \eqref{lg:eq:momentout}, the contribution of the
initial state is bounded by $Ce^{-ci/g^2}$. This bound is still of
order one at $i\asymp g^2$ and falls to a negative power of $W$ at
$n_*=\lceil C_*g^2\log W\rceil$. We therefore split the total growth as
\[
 \sum_{i=1}^{n_0}c_i
 =\underbrace{\sum_{i=1}^{n_*}c_i}_{\ell_{n_*}}
  +\sum_{i=n_*+1}^{n_0}c_i.
\]
For the first sum, the rough bound $|\E_{i-1}c_i|\le Cg^{-2}$ and the
martingale estimate give $|\ell_{n_*}|\le C_D\log W$ with probability
at least $1-W^{-D}$ by \eqref{lg:eq:initial-segment}. Once the moments
have stabilized, the corrected conditional drift is of order
$(Wg^2)^{-1}$ per block, or $W^{-1}$ over $g^2$ blocks. Over the
remaining $n_0-n_*$ blocks it therefore contributes order $C_1\log W$.
Choosing $C_1$ large makes this gain exceed the possible initial loss
$C_D\log W$.

\medskip\noindent\textbf{Further expansion.}
We first estimate the trace and directional moments that occur in the
conditional drift. The first two moments already show why an even
directional term needs a different correction from $Q^{[1]}$.

\medskip\noindent\textbf{Smallness of the first trace and directional moments.}
For $n_*\le i<n_0$ we estimate $c_{i+1}$ at the incoming state
$(C_i,v_i)$. To illustrate the first two corrections, take
$g\asymp W^\epsilon$ with $0<\epsilon<1/4$, so $W^{-1}=o(g^{-4})$.
Start with the first corrections to $m_1=W^{-1}\Tr C_i$
and $f_1=v_i^*C_i v_i$ from the preceding subsection, denoted here
by $M_1^{[1]}$ and $F_1^{[1]}$. The damping estimates
\eqref{uvg:eq:Mdrift}--\eqref{uvg:eq:Fdrift} show why we bound the
trace moment first. Write $u_i=\E|M_1^{[1]}(C_i)|^2$ and
$w_i=\E|F_1^{[1]}(C_i,v_i)|^2$. Taking expectations gives
\[
 \begin{aligned}
 u_{i+1}&\le(1-cg^{-2})u_i+Cg^{-6},\\
 w_{i+1}&\le(1-cg^{-2})w_i
               +Cg^{-2}u_i+Cg^{-6}.
 \end{aligned}
\]
In the first line the geometric sum of the damping factors is $O(g^2)$.
Thus $u_i\lesssim e^{-ci/g^2}+g^{-4}$. Inserting this in the
second line and summing again gives
$w_i\lesssim e^{-c'i/g^2}+g^{-4}$ for some $c'>0$.
For $i\ge n_*$, choosing $C_*$ sufficiently large gives
\[
 u_i\lesssim g^{-4},\qquad w_i\lesssim g^{-4}.
\]
For every fixed $1\le p<\infty$, iterating
\eqref{lg:eq:example-moment-M}--\eqref{lg:eq:example-moment-F} gives
$\|M_1^{[1]}(C_i)\|_{L^p}+\|F_1^{[1]}(C_i,v_i)\|_{L^p}\lesssim_p g^{-2}$
for $i\ge n_*$, with $C_*$ chosen for $p$.

\medskip\noindent\textbf{Estimating $M_2^{[1]}$.}
At $E=0$ the free step $C\mapsto-C$ leaves $m_2=W^{-1}\Tr C^2$
unchanged. Squaring the one-step update \eqref{lg:eq:matrix-update}
and averaging its normalized trace gives
\[
 (P-I)\circ m_2=-(2g)^{-2}
 \bigl[(1+2m_1+m_2)^2+2(1+m_1)(m_1+2m_2+m_3)\bigr]+O(g^{-4}).
\]
The odd part of the bracket is $2(3m_1+m_3+4m_1m_2)$. It changes
sign under the free step, so the correction
$M_2^{[1]}=m_2-(2g)^{-2}(3m_1+m_3+4m_1m_2)$ cancels it.
The even part is $1+6m_2+m_2^2+6m_1^2+2m_1m_3$, as in
\eqref{lg:eq:center-trace}. For $i\ge n_*$, the first-moment bound
on $m_1$ and $|m_3|\le1$ give, in each fixed $L^p$,
\[
 (P-I)\circ M_2^{[1]}
 =-(2g)^{-2}(1+6m_2+m_2^2)+O(g^{-4}).
\]
The unique zero of this leading polynomial in the
unit disk is $a_2^{(0)}=-3+\sqrt8$. Set $X_i=M_2^{[1]}(C_i)-a_2^{(0)}$.
Factoring the polynomial as $(m_2-a_2^{(0)})(6+m_2+a_2^{(0)})$
and using $m_2(C_i)-a_2^{(0)}=X_i+O(g^{-2})$ gives
\[
 \E_i X_{i+1}
 =\bigl[1-(2g)^{-2}(6+m_2(C_i)+a_2^{(0)})\bigr]X_i+O(g^{-4}).
\]
Since $|m_2(C_i)|\le1$ and $\Re(6+m_2(C_i)+a_2^{(0)})>4$, the
multiplier has modulus at most $1-cg^{-2}$ for large $g$.
The same moment iteration as for $M_1^{[1]}$ gives
\[
 \|M_2^{[1]}(C_i)-a_2^{(0)}\|_{L^p}\lesssim_p g^{-2},
 \qquad n_*\le i<n_0.
\]
Since $m_2-M_2^{[1]}=O(g^{-2})$, the same estimate holds for $m_2$.

\medskip\noindent\textbf{Estimating $F_2^{[1]}$.}
The free step also fixes $f_2=v^*C^2v$. Each term in the bracket
below changes sign under the free step. Thus the odd terms in
\eqref{lg:eq:example-raw-f} cancel in the drift of the first
correction \eqref{lg:eq:example-corrections},
\[
 F_2^{[1]}=f_2-(8g^2)^{-1}
 \bigl[(3+4m_2)f_1+3f_3+8m_1f_2-(f_1+\overline{f_1})f_2\bigr].
\]
In the even part of the order-$(2g)^{-2}$ coefficient, every term
except $-6f_2-2m_2f_2$ contains $m_1$ or $f_1$. Both first moments
are $O(g^{-2})$ for $i\ge n_*$. Using $f_2=F_2^{[1]}+O(g^{-2})$
in the remaining term gives, in each fixed $L^p$,
\[
 (P-I)\circ F_2^{[1]}
 =-(2g^2)^{-1}(3+m_2)F_2^{[1]}+O(g^{-4}).
\]
Since $\Re(3+m_2)\ge2$, the random step damps $F_2^{[1]}$ toward
zero at rate $g^{-2}$. Using the bounds above for
$M_1^{[1]},F_1^{[1]}$ in the moment recursion
\eqref{lg:eq:first-F2-moment} gives
\[
 \|F_2^{[1]}(C_i,v_i)\|_{L^p}\lesssim_p g^{-2},
 \qquad n_*\le i<n_0.
\]
The same bound holds for $f_2$.

\medskip\noindent\textbf{An even term in the corrected drift.}
Substitute these estimates into the drift after
$Q^{[1]}=-(8g^2)^{-1}\Re f_1$. By \eqref{lg:eq:c2},
\[
 \E_i c_{i+1}-[(P-I)\circ Q^{[1]}](C_i,v_i)
 =(8Wg^2)^{-1}+(4g^2)^{-1}
  \bigl[\Re(m_1f_1)-(2W)^{-1}\Re(f_1^2)\bigr]+O(g^{-4}).
\]
The two products are $O(g^{-4})$ in each fixed $L^p$ by the
first-moment bounds. As we expand further,
a term linear in an odd $f_j$ with a fixed coefficient can be cancelled
just as for $f_1$, since its free increment is $-2f_j$.
Consider instead a fixed multiple of $\Re f_2$. Its free increment is zero.

To see where such a term enters, keep only the coefficient of $f_2$,
taking real parts and omitting terms with an explicit $W^{-1}$.
At fourth order the logarithm in $c_{i+1}$ and the update of $Q^{[1]}$
give $(2g)^{-4}[-(1+m_1)^2+(1+m_1)]$. The constant parts cancel.
Substituting the first trace correction
$m_1=M_1^{[1]}+(8g^2)^{-1}(1+m_2+2m_1^2)$ leaves a contribution
$-\tfrac12(1+a_2^{(0)})$ to the coefficient of $(2g)^{-6}\Re f_2$
when the traces are set to their centers. The next Wick terms from the
logarithm and the update of $Q^{[1]}$ contribute
$(4+2a_2^{(0)})-(5+2a_2^{(0)})=-1$. Thus, by \eqref{lg:eq:first-even-drift}, the even term to be matched is
\[
 -\tfrac12(2g)^{-6}(3+a_2^{(0)})\Re f_2 .
\]
The $F_2^{[1]}$ calculation above explains how to match it.
At the trace centers, a unit coefficient of $f_2$ produces an $f_2$
coefficient $-(2g^2)^{-1}(3+a_2^{(0)})$ in $(P-I)\circ f_2$
by \eqref{lg:eq:example-raw-f}. To match the displayed term, we therefore
include $\tfrac14(2g)^{-4}\Re f_2$ in $Q^{[2]}-Q^{[1]}$. The factor
$g^2$ in this division is the cost of the even mode.

\medskip\noindent\textbf{Correcting the moments to higher order.}
For the full range $W^\epsilon\le g<W$ we carry out the moment corrections first.
Write $M_j^{[r]},F_j^{[r]}$ for moment $j$ after $r$ rounds.
Each round expands their one-step updates one order further in $g^{-2}$.
A monomial of total index $d$ acquires the free sign $(-1)^d$.
If its parity differs from that of $j$, the two free multipliers
differ by $\pm2$, so its coefficient can be cancelled as in the first
corrections above. The terms of the same parity remain in the update.
Consider the two even terms $m_2f_2$ and $m_1f_3$ in the example above.
The first combines with $3f_2$ to give the damping coefficient $3+m_2$,
whose real part is at least $2$. For the second, the bound on $m_1$
already suffices. Since $|f_3|\le1$, its contribution to the one-step
drift is $O(g^{-4})$ in each fixed $L^p$.

For general $j$, we first center the trace moments at $a_j$, as we
did for the second trace. The second-order drift formula for the
directional moments \eqref{lg:eq:center-corner} gives the terms of
the same parity as $j$. Its diagonal term gives damping at order
$g^{-2}$, and the terms involving lower directional moments are
controlled inductively. The only higher directional moment at this
order is $f_{j+1}$, which appears multiplied by $m_1$. Thus
$|f_{j+1}|\le1$ suffices at this step. Higher-order couplings carry
at least $g^{-4}$, so even after summing over the $O(g^2)$ damping
time they retain a factor $g^{-2}$. Iterating the moment recursion
\eqref{lg:eq:momentrec} uses this gain to reach the required accuracy.

Wick pairing restricts the indices at each
order. For instance, at order $g^{-2r}W^0$ in the conditional drift,
each index has modulus at most $r$ and the sum of their moduli is at
most $2r$ (see the proof of Lemma~\ref{lg:lem:ideal}).
Choose the expansion depth $R$ with $2\epsilon R>1$. A remaining
one-step error $O(g^{-2R-2})$ accumulates over the damping time
$O(g^2)$ to $O(g^{-2R})=o(W^{-1})$. Writing $M_j,F_j$ for the
resulting moments, \eqref{lg:eq:momentout} gives, for fixed small
$\eta>0$, with high probability simultaneously at the required indices
and times $n_*\le i<n_0$,
\[
 |M_j(C_i)-a_j|\le W^{-1+\eta},\qquad
 |F_j(C_i,v_i)|\le W^{-1/2+\eta}.
\]

\medskip\noindent\textbf{Choosing the drift correction.}
Now expand $\E_i c_{i+1}$ in these corrected moments. Set the traces
to their centers and collect the terms linear in $F_j$. The preceding
bounds control the products and the errors from replacing the traces.
At each round we add a linear combination of the $F_j$'s to $Q^{[s]}$
so that $(P-I)\circ(Q^{[s+1]}-Q^{[s]})$ matches the linear
coefficients still present at that order. The divisor found for $f_2$ persists for
every even index. For $j>0$, the two diagonal terms in the
second-order drift formula for the directional moments
\eqref{lg:eq:center-corner} are $-3jF_j$ and $-ja_2F_j$ at the
trace centers. Including the free sign, the coefficient of $F_j$ in
$(P-I)\circ F_j$ through order $g^{-2}$ at $E=0$ is
\[
 (-1)^j\bigl[1-j(2g)^{-2}(3+a_2)\bigr]-1 .
\]
It is of order one for odd $j$ and of order $g^{-2}$ for even $j$,
since $3+a_2=3+a_2^{(0)}+O(g^{-2})=\sqrt8+O(g^{-2})$
by \eqref{lg:eq:centers}. The other linear terms
at this order have lower indices of the same parity, so we can
match coefficients by triangular inversion. The higher-order part of
the update is $O(g^{-4})$, so even after the $O(g^2)$ inverse each
successive correction gains a factor $O(g^{-2})$.
Section~\ref{lg:sec:cocycle-correction} carries out this finite
iteration in \eqref{lg:eq:Q-recursion}, retaining the free phases
for energies near zero. At $E=0$ the state-independent term
$(8Wg^2)^{-1}$ from \eqref{lg:eq:c2} remains, and the resulting
bounded correction $Q_g$ gives
$\E_iY_{i+1}=(8Wg^2)^{-1}+o((Wg^2)^{-1})$ by \eqref{lg:eq:smalldefect}.

Returning to the two sums split at $n_*$, the second contributes
corrected drift of order $C_1\log W$ at $n_0=\lceil C_1Wg^2\log W\rceil$.
By Lemma~\ref{lg:lem:fluctuations}, the martingale costs at most
$C_D\sqrt{C_1}\log W$ and the initial interval at most $C_D\log W$
with high probability. Including the bounded endpoint difference,
\eqref{lg:eq:integrated} yields
\[
 \ell_{n_0}(C_0,w)\ge
 \bigl(cC_1-C_D\sqrt{C_1}-C_D\bigr)\log W
\]
for each fixed initial normal $w$. Taking $C_1$ sufficiently large
and applying \eqref{lg:eq:scalar-reduction} to a fixed orthonormal
basis with one extra power of $W$ in the failure bound gives the
claimed estimate for $\|D_{n_0}(\Phi;0)^{-1}\|$.

\begin{figure}[H]
\centering
\hypersetup{allcolors=black}
\begin{tikzpicture}[
 method/.style={draw=black,rounded corners=2pt,align=center,
   font=\fontsize{8}{9.5}\selectfont,inner xsep=4pt,inner ysep=5pt,
   minimum height=0.9cm},
 framework/.style={method,text width=11.5cm},
 result/.style={framework,fill=black!3},
 route/.style={->,line width=0.5pt,draw=black}]
 \node[framework] (fm) at (7.6,0)
   {\textbf{Fractional moments of the Green function}\\
    $g<W$: fixed energy \qquad $g\ge W$: energy averaging, blocks grouped on scale $g^2$\\
    Lemma~\ref{uva:thm:upper-audit}; fractional-moment method of \cite{AM,ASFH}};
 \node[method,text width=3.1cm] (strong) at (1.7,-2.55)
   {\textbf{Extremely small coupling}\\
    Direct fractional-moment decay from \cite{PSSS}\\
    Theorem~\ref{thm:strong-upper}};
 \node[method,text width=8.0cm] (transfer) at (9.5,-2.55)
   {\textbf{Growth in the weakest boundary direction}\\
    $g<W$: growth at one length\\
    $g\ge W$: negative-moment decay at every length\\
    Proposition~\ref{up:prop:conversion}; boundary comparison of \cite{MS}};
 \draw[route] ([xshift=-4.4cm]fm.south) -- ++(0,-0.4) -| (strong.north);
 \draw[route] ([xshift=1.9cm]fm.south) -- (transfer.north);
 \node[method,text width=6.4cm] (operator) at (3.7,-4.95)
   {\textbf{Contraction of $P_{1/2}$}\\
    Perturb a fixed number of blocks and bound the Gaussian divergence\\
    Lemmas~\ref{up:lem:operator-reduction}, \ref{up:lem:operator-iteration},
    \ref{up:lem:angular}; criterion adapted from \cite{Drogin}};
 \node[method,text width=6.4cm] (direct) at (11.5,-4.95)
   {\textbf{Direct transfer growth}\\
    Apply the corrector method of \cite{SBanomaly,SSB} to the logarithmic volume increment};
 \draw[route] (transfer.south) -- ++(0,-0.45) -| (operator.north);
 \draw[route] (transfer.south) -- ++(0,-0.45) -| (direct.north);
 \node[method,text width=3.15cm] (small) at (1.75,-7.7)
   {$\Delta\Delta$: \textbf{Small coupling}\\
    Match the strongly amplified directions\\
    Lemma~\ref{up:lem:small-contraction}};
 \node[method,text width=3.15cm] (fixed) at (5.65,-7.7)
   {\textbf{Stretch a prescribed direction}\\
    General perturbation, sharper cost away from zero\\
    Lemma~\ref{up:lem:density-contraction}};
 \node[method,text width=3.15cm] (center) at (9.55,-7.7)
   {\textbf{Near zero for $W^\epsilon\le g<W$}\\
    $\Delta\Delta$: growth at length $n_0\asymp Wg^2\log W$\\
    Lemma~\ref{up:lem:center-transfer}};
 \node[method,text width=3.15cm] (uniform) at (13.45,-7.7)
   {\textbf{Full bulk for $g\ge W$}\\
    $\Delta$: negative-moment decay at every length\\
    Theorem~\ref{uvg:thm:transfer}};
 \draw[route] (operator.south) -- ++(0,-0.4) -| (small.north);
 \draw[route] (operator.south) -- ++(0,-0.4) -| (fixed.north);
 \draw[route] (direct.south) -- ++(0,-0.4) -| (uniform.north);
 \draw[route] (direct.south) -- ++(0,-0.4) -| (center.north);
 \node[result] (upper) at (7.6,-10.6)
   {\textbf{Concentration upper bound}\\
    Half the mass in at most $W^\tau\ellstar$ coordinates outside the exceptional set};
 \draw[route] (strong.west) -- (-0.35,-2.55) |- (upper.west);
 \draw[route] (small.south) -- ++(0,-0.45) -| ([xshift=-4.5cm]upper.north);
 \draw[route] (fixed.south) -- ([xshift=-1.95cm]upper.north);
 \draw[route] (center.south) -- ([xshift=1.95cm]upper.north);
 \draw[route] (uniform.south) -- ++(0,-0.45) -| ([xshift=4.5cm]upper.north);
\end{tikzpicture}
\caption{Upper-bound methods. $\Delta\Delta$ marks the small-coupling construction and the
near-zero growth estimate. $\Delta$ marks the full-bulk estimate for $g\ge W$.
Both direct transfer estimates use a corrected drift on the scale $(Wg^2)^{-1}$.
Figure~\ref{fig:upper-ranges} gives the coupling and energy ranges.}
\label{fig:upper-methods}
\end{figure}
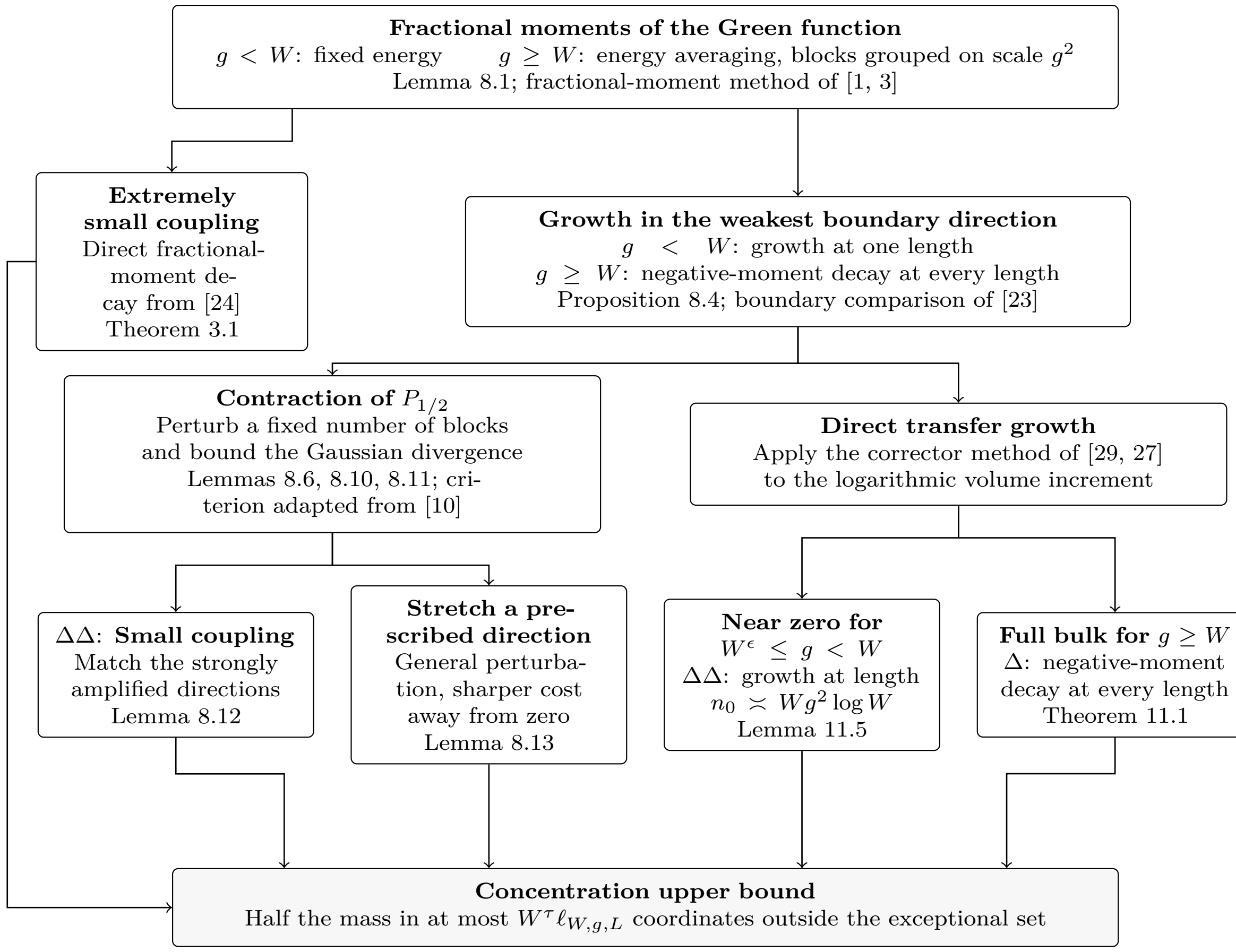

\subsection{Component estimates and proof of the upper bound}
\label{sec:upper-regions}

The deterministic bound $\ell(\psi_k)\le WL$ suffices whenever
$L\le W^\tau\max\{1,Wg^2\}$. For longer cycles we use the four
concentration estimates below. The first uses \cite{PSSS} directly.
The other three use transfer growth, with fixed-energy fractional
moments for $g<W$ and energy-averaged fractional moments for $g\ge W$.
Figures~\ref{fig:upper-methods} and~\ref{fig:upper-ranges} give the
methods and their coupling ranges, respectively.

\begin{figure}[H]
\centering
\renewcommand{\arraystretch}{1.35}
\fontsize{8.5}{10.5}\selectfont
\begin{tabularx}{\textwidth}{@{}p{0.30\textwidth}p{0.40\textwidth}X@{}}
\toprule
\textbf{Coupling} & \textbf{Main idea} & \textbf{Corresponding results}\\
\midrule
\textbf{Extremely small}\newline $g\le c_* W^{-1/2}$
& The fractional-moment estimate of \cite{PSSS} gives concentration on the single-block scale.
& Theorem~\ref{thm:strong-upper}\\[3pt]
\textbf{Small}\newline $c_*W^{-1/2}\le g\le W^{-\epsilon}$
& Match the amplified directions to obtain contraction at the block scale $\max\{1,Wg^2\}$ up to a small power of $W$.
& Lemma~\ref{up:lem:small-contraction}; Theorem~\ref{thm:sharp-small-upper}\\[3pt]
\textbf{Intermediate}\newline $W^{-\epsilon}\le g\le W^\epsilon$
& A fixed matrix direction gives the block scale $W\max\{g^{-2},g^4\}$, sufficient for small $\epsilon$.
& Lemma~\ref{up:lem:density-contraction}(i); Theorem~\ref{thm:compact-upper}\\[3pt]
\textbf{Large}\newline $W^\epsilon\le g<W$
& Away from zero, sharpen the contraction of $P_{1/2}$. Near zero, estimate transfer growth directly. Both give growth at one length.
& Lemmas~\ref{up:lem:density-contraction}(ii), \ref{up:lem:center-transfer}; Theorem~\ref{thm:large-upper}\\[3pt]
\textbf{Largest couplings}\newline $g\ge W$
& Direct negative-moment decay at every length gives the energy-averaged bound throughout the bulk.
& Theorem~\ref{uvg:thm:transfer}; Lemma~\ref{uva:thm:upper-audit}(ii); Theorem~\ref{thm:large-upper}\\
\bottomrule
\end{tabularx}
\normalsize
\caption{Upper-bound regimes for fixed $0<\epsilon<1/4$. For $W^\epsilon\le g<W$ the two energy windows overlap. The component theorems give the precise bounds.}
\label{fig:upper-ranges}
\end{figure}

The four theorems specify an energy window $J$ and a block radius $R$
for which, given $D>0$,
\begin{equation}\label{eq:averaged-tail}
 \frac1{WL}\E\#\left\{k:\lambda_k\in J,\quad
 \min_{a\in\Z_L}\sum_{\|b-a\|_L>R}\|\psi_{k,b}\|_2^2>W^{-D}\right\}
 \le W^{-D}.
\end{equation}
For every remaining eigenvector with $\lambda_k\in J$, a block ball of radius
$R$ contains at least half the mass, so $\ell(\psi_k)\le\min\{N,W(2\lceil R\rceil+1)\}$
for large $W$.

\begin{theorem}[Strong-disorder concentration; finite-volume corollary of
{\cite[Theorem~1]{PSSS}}]\label{thm:strong-upper}
There exists a sufficiently small constant $c_*\in(0,1)$ such that, for every fixed
$D>0$, all sufficiently large $W$, every deterministic
$g\le c_*W^{-1/2}$, and every $L$,
\eqref{eq:averaged-tail} holds with
\begin{equation}\label{eq:strong-radius}
 R\coloneqq\lceil(\log W)^2\rceil,\qquad J=\R.
\end{equation}
The width threshold is independent of $L,g$.
\end{theorem}

The proof is given at the end of Section~\ref{sec:fr-concentration}.

\begin{theorem}[Sharp concentration at small coupling]
\label{thm:sharp-small-upper}
Fix $0<\alpha<1/2$, $c_*,D,\varepsilon>0$, and
$\kappa\in(0,2)$. There is $W_0$ such that, for every
$W\ge W_0$, every $L$, and every deterministic
\[
 c_*W^{-1/2}\le g\le W^{-\alpha},
\]
\eqref{eq:averaged-tail} holds with
\begin{equation}\label{eq:sharp-small-radius}
 R\le W^\varepsilon\max\{1,Wg^2\},\qquad J=[-2+\kappa,2-\kappa].
\end{equation}
The width threshold is independent of $L,g$ in the stated range.
\end{theorem}

The proof in Section~\ref{sec:upper-center-completion} uses
Lemma~\ref{up:lem:small-contraction} from Section~\ref{sec:small-upper-proof}.

\begin{theorem}[Concentration from the general lift]
\label{thm:compact-upper}
Fix $B,E_0,D>0$. For all sufficiently large $W$, every deterministic
$W^{-B}\le g\le W^B$, and every $L$,
\eqref{eq:averaged-tail} holds with
\begin{equation}\label{eq:coarse-radius}
 R\le W\max\{g^{-2},g^4\}(\log W)^2,\qquad J=[-E_0,E_0].
\end{equation}
The width threshold depends only on $D,B,E_0$ and is uniform in $g,L$.
\end{theorem}

The proof in Section~\ref{sec:upper-center-completion} uses the general
lift in Section~\ref{sec:compact-proof}.

\begin{theorem}[Concentration at large coupling]
\label{thm:large-upper}
Fix $\kappa\in(0,2)$ and $D,\epsilon>0$.
There is $W_0$ such that, for every $W\ge W_0$, every deterministic
finite $g\ge W^\epsilon$, and every $L$,
\eqref{eq:averaged-tail} holds with
\begin{equation}\label{eq:large-radius}
 R\le Wg^2(\log W)^2,\qquad J=[-2+\kappa,2-\kappa].
\end{equation}
The width threshold depends only on
$\kappa,D,\epsilon$, and not on $L,g$.
\end{theorem}

The proof concludes Section~\ref{sec:large-upper-proof}.

\begin{proof}[Proof of the upper bound and completion of Theorem~\ref{thm:main}]
Fix $0<\tau<1$ and define
\[
 B_+=\#\{k:\lambda_k\in J_{g,\kappa,L},\ \ell(\psi_k)>W^\tau\ellstar\}.
\]
If $L\le W^\tau\max\{1,Wg^2\}$, then $B_+=0$ by the cap $\ell(\psi_k)\le WL$.
Otherwise $\ellstar=\ellscale$ and $L>W^\tau\ge L_{\kappa/2}$ for large $W$,
so $F_g(-|E|)\ge\kappa/2$ on $J_{g,\kappa,L}$ by \eqref{eq:bulk-sampling}.
Split the remaining couplings at $c_*W^{-1/2}$, $W^{-\tau/8}$ and
$W^{\tau/4}$. In increasing order use Theorem~\ref{thm:strong-upper},
Theorem~\ref{thm:sharp-small-upper} with $(\alpha,\varepsilon)=(\tau/8,\tau/2)$,
Theorem~\ref{thm:compact-upper} with $(B,E_0)=(\tau/4,4)$, and
Theorem~\ref{thm:large-upper} with $\epsilon=\tau/4$.
The first theorem holds on $\R$. For the small- and large-coupling
theorems, \eqref{eq:bulk-bridge} at $\kappa/2$ gives a common bulk margin
$\kappa_0>0$. For the intermediate range, \eqref{eq:bounded-physical-bulk}
gives $J_{g,\kappa,L}\subseteq[-4,4]$.
For Theorem~\ref{thm:compact-upper}, the radius ratio to $\max\{1,Wg^2\}$ costs
at most $W^{\tau/2}(\log W)^2$, since the extra factors are
$g^{-4}$ for $g\le1$ and $g^2$ for $g\ge1$.
The radii in all four theorems therefore satisfy
$W(2\lceil R\rceil+1)\le W^\tau\ellscale$ for large $W$.
Thus \eqref{eq:averaged-tail} in the applicable range gives
\begin{equation}\label{uv:eq:upper-count}
 (WL)^{-1}\E B_+\le W^{-D}.
\end{equation}
By \eqref{eq:bulk-sampling} and Lemma~\ref{lem:rough-density}, the
expected fraction in \eqref{eq:main}, with $0/0=0$, is at least
\begin{equation}\label{uv:eq:sampling}
 1-(1-2\kappa)^{-1}(WL)^{-1}\E(B_-+B_+)
   -C_\kappa e^{-c_\kappa W^2L}\ge1-W^{-D},
\end{equation}
using \eqref{uv:eq:lower-regions-output} and \eqref{uv:eq:upper-count}
with exponent $D+2$.
\end{proof}

\section{Spreading from local resolvent moments}
\label{sec:unrestricted}\label{uv:sec:local-resolvents}
We prove Theorem~\ref{uv:lem:lower-volume} by counting large spectral
atoms. For large volumes, a comparison with a local periodic resolvent
first transfers the moment bound to the full cycle.

\begin{theorem}[Local resolvent comparison]
\label{uv:thm:smooth-cutoff}
Fix a central block $a$ and an integer $R\ge1$ with $2R+1<L$.
Let $I=[a-R,a+R]$, and let $H_I^{\mathrm{per}}$ be the restriction of
$H$ to $I$ with periodic boundary conditions and the same onsite blocks and hopping.
For every $z=E+i\eta$ with $\eta R\ge2$, we have, in the ordering
of positive semidefinite $W\times W$ matrices,
\begin{equation}\label{uv:eq:smooth-comparison}
 C^{-1}\Im G_{I,aa}^{\mathrm{per}}(z)
 \ \le\ \Im G_{L,aa}(z)\
 \le C\Im G_{I,aa}^{\mathrm{per}}(z),
\end{equation}
where $C$ is an absolute constant. The estimate is deterministic.
\end{theorem}
\begin{proof}

Define a rectangular multiplication operator from the full cycle of $L$ blocks
to $I$ by
$\chi(b)=\max\{0,1-|b-a|/R\}$.
Onsite (i.e.\ within the same block) terms commute with $\chi$.
The geometric resolvent identity reads
\[
 G_I^{\mathrm{per}}K G_L
       =\chi G_L-G_I^{\mathrm{per}}\chi, \quad K:=H_I^{\mathrm{per}}\chi-\chi H_L 
\]
For $u\in\C^W$, let $u_L,u_I$ denote its embeddings at block $a$ in the
full cycle and in $I$, respectively. Since $\chi(a)=1$,
$(\chi G_L)_{aa}=(G_L)_{aa}$ and
$(G_I^{\mathrm{per}}\chi)_{aa}=(G_I^{\mathrm{per}})_{aa}$.
Cauchy--Schwarz and Ward's identity give
\[
 \begin{aligned}
 \left|u^*(G_{L,aa}-G_{I,aa}^{\mathrm{per}})u\right|
 &=\left|u_I^*G_I^{\mathrm{per}}KG_Lu_L\right|\\
 &\le\norm K\,\norm{(G_I^{\mathrm{per}})^*u_I}\,\norm{G_Lu_L}
 \le \norm K\, \frac{1}\eta  \left(u^*\Im G_{L,aa}u \right)^{1/2}\left(  u^*\Im G_{I,aa}^{\mathrm{per}}u\right)^{1/2}
 \\
 &\le \norm K\, \frac{1}{2\eta}  \left(u^*\Im G_{L,aa}u  +    u^*\Im G_{I,aa}^{\mathrm{per}}u\right) 
 \end{aligned}
\]
The only potentially nonzero blocks of the $(2R+1)\times L$ block matrix
$K$ are
\[
 K_{bc}=\frac{g}{\max\{1,g\}}\bigl(\chi(c)-\chi(b)\bigr)I_W,
 \qquad b,c\in I,\quad |b-c|=1.
\]
The Schur bound gives $\norm K\le 2g/(\max\{1,g\}R)\le\eta$.
Thus \eqref{uv:eq:smooth-comparison} holds with $C=3$.
\end{proof}

In the application to Theorem~\ref{uv:lem:lower-volume}, it suffices
to bound the resolvent moments on the local cycle of $L_c$ blocks.

The comparison \eqref{uv:eq:smooth-comparison} transfers local moments to the full cycle. We now convert these moments into a bound on the number of eigenvectors with short half-mass length.
\begin{lemma}[Poisson smoothing of heavy spectral atoms]
\label{uv:lem:poisson}
Let $\mu=\sum_k p_k\delta_{\lambda_k}$ be a finite positive measure,
$\eta,h>0$, and
\[
 f(E)=\int\frac{\eta}{(\lambda-E)^2+\eta^2}\,d\mu(\lambda),
 \qquad J_\eta=J+[-\eta,\eta].
\]
Then
\begin{equation}\label{uv:eq:poisson}
 \sum_{\lambda_k\in J}p_k\one_{\{p_k>h\}}
 \le\frac2\pi\int_{J_\eta}
       f(E)\one_{\{f(E)>h/(2\eta)\}}\,dE.
\end{equation}
There is also a matrix version. If
$\boldsymbol\mu=\sum_k v_kv_k^*\delta_{\lambda_k}$,
$m_k=\norm{v_k}^2$, and
$F(E)=\int\eta((\lambda-E)^2+\eta^2)^{-1}\,d\boldsymbol\mu(\lambda)$,
then
\begin{equation}\label{uv:eq:matrix-poisson}
 \sum_{\lambda_k\in J}m_k\one_{\{m_k>h\}}
 \le\frac2\pi\int_{J_\eta}
       \Tr F(E)\one_{\{\norm{F(E)}>h/(2\eta)\}}\,dE.
\end{equation}
\end{lemma}
\begin{proof}
If $p_k>h$ and $|E-\lambda_k|\le\eta$, positivity gives
$f(E)\ge p_k/(2\eta)>h/(2\eta)$.
Also
$\int_{\lambda_k-\eta}^{\lambda_k+\eta}
\eta((\lambda_k-E)^2+\eta^2)^{-1}\,dE=\pi/2$.
On the union of these intervals, $f$ is at least the sum of the
corresponding Poisson contributions, each restricted to its own
interval. Tonelli's theorem proves \eqref{uv:eq:poisson}, even when
the intervals overlap. In the matrix case, the rank-one atom
forces $\norm{F(E)}\ge m_k/(2\eta)$ on the same interval, and the
trace contains its scalar contribution. The identical argument
proves \eqref{uv:eq:matrix-poisson}.
\end{proof}

\begin{proof}[Proof of Theorem~\ref{uv:lem:lower-volume}]
Put $\eta_L=W^\nu/\ellstar$. If $L_c<L$, then $WL_c>\ellscale$,
$\eta_c=\eta_L=\eta$, and $\eta R\ge2$. Theorem~\ref{uv:thm:smooth-cutoff}
and the law of the local periodic restriction transfer
\eqref{uv:eq:lower-volume-moment} to the original cycle. For $L_c=L$,
apply the assumption directly. Thus, with $C_p$ enlarged as needed,
\[
 \sup_{E\in J_{\eta_L}}\max_{1\le x\le WL}
 \E\bigl(\Im G_{xx}(E+i\eta_L)\bigr)^p\le C_p.
\]
Set $r=W^{-\tau}\ellstar$, and let $B_-(r)$ count eigenvectors
with eigenvalue in $J$ and half-mass length at most $r$.
Using \eqref{uv:eq:poisson} at each coordinate, we obtain
\[
 \begin{aligned}
 B_-(r)&\le4\sum_{k:\lambda_k\in J}\sum_x
 |\psi_k(x)|^2\one_{\{|\psi_k(x)|^2>1/(4r)\}}\\
 &\le\frac8\pi\sum_x\int_{J_{\eta_L}}
 \Im G_{xx}(E+i\eta_L)
 \one_{\{\Im G_{xx}(E+i\eta_L)>1/(8r\eta_L)\}}\,dE.
 \end{aligned}
\]
The moment bound above and $|J_{\eta_L}|\le2E_0+2$ then give
\[
 \frac{\E B_-(r)}{WL}\le C_p W^{-(p-1)(\tau-\nu)}.
\]
Since $(p-1)(\tau-\nu)>D+1$,
this is at most $W^{-D}$ for sufficiently large $W$, proving
\eqref{uv:eq:lower-volume-output}, uniformly in $g,L$ under the
stated assumptions.
\end{proof}

\section{The loop hierarchy and propagation}
\label{sec:flow-overview}
{\bf Basic strategy:} Following \cite{YY,TYY}, we obtain resolvent bounds through the resolvent
loops defined in \eqref{fo:eq:loops}. In particular, the mixed two-loop
controls the Hilbert--Schmidt norm of a Green block,
\[
 \cL_2(t;a,b;+,-)=W^{-2}\|(G_t)_{ba}\|_{\HS}^2,
\]
which, together with unitary invariance and fluctuation averaging, yields
entrywise and averaged resolvent bounds. Under the characteristic stochastic
flow \eqref{fo:eq:flow}--\eqref{fo:eq:characteristic}, these loops satisfy the
hierarchy \eqref{fo:eq:hierarchy}: both the drift and the quadratic variations
of the martingales are expressed in terms of loops. The main task is to show
that the random loops $\cL_k(t)$ remain close to the deterministic primitive
loops $\cK_k(t)$, defined by \eqref{fo:eq:primitive}, as the spectral height
decreases to $\eta_L$. For bounded coupling, we follow the loop
hierarchy method and prove its inputs here, adapting the stability analysis of
\cite[Section~8 and Appendix~B.1]{DYYY}.

The main new difficulty is to make this argument uniform for arbitrarily
large $g$. Here we use transfer matrices, a tool more commonly associated
with localization proofs. Their role is natural: after normalization, the
onsite randomness has size $g^{-1}$, while the free transfer matrix in the
bulk is conjugate to a unitary matrix. We exploit this deterministic
structure and treat the onsite randomness as a small centered perturbation,
a viewpoint that also underlies our transfer-matrix analysis of the upper
bound.

Transfer estimates provide three inputs to the loop hierarchy. Combined
with Ward's identity, they give initial resolvent moments at
$\eta\asymp g^{-2}$, as well as coarse moment bounds throughout the subsequent
propagation. They also yield spatial interpolation of Green functions within
cells much shorter than $g^2$. This interpolation, extended to both random
and primitive loops, reduces the spatial estimates to $W^{O(1)}$ sampled
block tuples, even when the cycle length is superpolynomial in $W$. These
inputs allow the loop hierarchy to propagate the sharp estimates down to
$\eta_L$, with constants uniform in $g$.

\subsection{The characteristic stochastic flow}
\label{sec:fo-characteristic}
 For each coupling $g$, define
\begin{equation}\label{fo:eq:flow}
 \begin{gathered}
 A_g\coloneqq\sqrt{\gamma_g}\,gA_L,\quad \gamma_g\coloneqq\max\{1,g\}^{-2}.\\
 H_t\coloneqq H_t(g)=A_g\otimes I_W+\mathcal B_t^{(g)},\qquad 0\le t\le1.
 \end{gathered}
\end{equation}
The process $\mathcal B_t^{(g)}$ is block diagonal, with independent
Hermitian Brownian blocks started at zero and covariance
\[
 \E(d\mathcal B_t^{(g)})_{ai,aj}(d\mathcal B_t^{(g)})_{bk,b\ell}
   =W^{-1}\gamma_g\delta_{ab}\delta_{i\ell}\delta_{jk}\,dt.
\]
At $t=1$, $H_1(g)$ has the same law as the normalized Hamiltonian $H$
in \eqref{eq:model} at coupling $g$.

For $z\in\C_+$, define the finite-cycle Dyson solution by
\begin{equation}\label{fo:eq:dyson-family}
 m_{g,L}(z)=L^{-1}\sum_{p\in(2\pi/L)\Z_L}
       \frac1{2\sqrt{\gamma_g}\,g\cos p-z-\gamma_g m_{g,L}(z)}.
\end{equation}
We take the Herglotz solution with asymptotic $-1/z$ and use its
boundary values at real energies. Equation~\eqref{fo:eq:dyson-family}
specializes \cite[(3.2)--(3.4)]{TYY} to our hopping and covariance.

Fix a real energy $E$ with $\Im m_{g,L}(E)>0$, and define
\begin{equation}\label{fo:eq:characteristic}
 \begin{aligned}
 m&\coloneqq m(g,E)=m_{g,L}(E),\\
 M_B&\coloneqq M_B(g,E)=(A_g-E-\gamma_g m(g,E))^{-1},\\
 M&\coloneqq M(g,E)=M_B(g,E)\otimes I_W,\\
 z_t&\coloneqq z_t(g,E)=E+(1-t)\gamma_g m(g,E),\\
 \eta_t&\coloneqq\eta_t(g,E)=\Im z_t(g,E)
       =(1-t)\gamma_g\Im m(g,E),\\
 G_t&\coloneqq G_t(g,E)=(H_t(g)-z_t(g,E))^{-1},\qquad 0\le t<1.\\
 \wt G_t &\coloneqq\wt G_t(g,E)=G_t(g,E)-M(g,E)
 \end{aligned}
\end{equation}
By definition, $\wt G_0=0$. We will propagate the approximation $G_t(g,E)\approx M(g,E)$ along the
characteristic, with the deterministic matrix $M(g,E)$ independent of $t$. These characteristics are those of \cite[Definitions~3.1--3.2, (3.5)--(3.8)]{TYY}
with variance rate $\gamma_g$; see also \cite[Definition~2.7]{YY}.
We suppress $(g, E)$ within a fixed flow, except when comparing flows.

For an $N\times N$ block matrix $A=(A_{ab})_{a,b\in\Z_L}$ with
$W\times W$ blocks, define
\[
 \begin{aligned}
 \langle A\rangle_{ab}&\coloneqq W^{-1}\Tr A_{ab}
       =W^{-1}\sum_{1\le i\le W} A_{ai,bi}.
 \end{aligned}
\] 
Since $M_{aa}=mI_W$, we have
$\langle\wt G_t\rangle_{aa}=\langle G_t\rangle_{aa}-m$.
It\^o's formula gives
\begin{equation}\label{fo:eq:ito}
 dG_t=-G_t\,d\mathcal B_t^{(g)}\,G_t
       +\gamma_g G_t\diag(\langle\wt G_t\rangle_{aa}I_W:a\in\Z_L)G_t\,dt.
\end{equation}
Equation~\eqref{fo:eq:ito} is \cite[(3.9)--(3.10)]{TYY} with covariance
operator $X\mapsto\gamma_g\diag(\langle X\rangle_{aa}I_W:a\in\Z_L)$;
see also the scalar-matrix equation following \cite[Definition~2.7]{YY}.

\subsection{The loop hierarchy and primitive loops}
\label{sec:fo-loops}\label{sec:fo-primitives}

Recall \(E_a=P_a/W\), where \(P_a\) is the projection onto block
\(a\). Write \(G_t^+\coloneqq G_t\), \(G_t^-\coloneqq G_t^*\), and use the same
convention for $\wt G_t$, $M$, and scalars. For signs
\(\boldsymbol\sigma=(\sigma_1,\ldots,\sigma_k)\) and block labels
\(\boldsymbol a=(a_1,\ldots,a_k)\), define
\begin{equation}\label{fo:eq:loops}
 \cL_k(t;\boldsymbol a,\boldsymbol\sigma)
 \coloneqq\cL_k(t;g,E;\boldsymbol a,\boldsymbol\sigma)
 \coloneqq \Tr\bigl(G_t^{\sigma_1}E_{a_1}\cdots
                 G_t^{\sigma_k}E_{a_k}\bigr).
\end{equation}
This is \cite[Definition~3.4, (3.16)--(3.19)]{TYY},
equivalently \cite[Definition~2.9]{YY}, with block size $W$.
The trace is unnormalized. In particular,
\[
 \cL_1(t;a,+)=\langle G_t\rangle_{aa},\qquad
 \cL_2(t;a,b;+,-)=W^{-2}\|(G_t)_{ba}\|_{\HS}^2.
\]
The mixed two-loop measures the resolvent mass between two blocks.
Longer loops arise when this quantity is differentiated or its
martingale fluctuations are estimated.

We describe the cut operations explicitly, following
\cite[Definition~3.5, (3.20)--(3.22)]{TYY} and
\cite[Definition~2.10]{YY}. Inserting \(E_c\) into the \(i\)-th
resolvent factor means replacing
\(G_t^{\sigma_i}\) by \(G_t^{\sigma_i}E_cG_t^{\sigma_i}\).
Denote the resulting \((k+1)\)-loop by
\(\operatorname{Cut}_i^c\cL_k(t)\).
For \(i<j\), the two cut words have indices
\begin{equation}\label{fo:eq:cuts}
\begin{aligned}
 \operatorname{Cut}_{L,i,j}^c:\quad
 &\text{signs }(\sigma_1,\ldots,\sigma_i,\sigma_j,\ldots,\sigma_k),\\
 &\text{labels }(a_1,\ldots,a_{i-1},c,a_j,\ldots,a_k),\\
 \operatorname{Cut}_{R,i,j}^c:\quad
 &\text{signs }(\sigma_i,\ldots,\sigma_j),\\
 &\text{labels }(a_i,\ldots,a_{j-1},c).
 \end{aligned}
\end{equation}
Thus the two loops share the label \(c\), and their lengths are
\(k+i-j+1\) and \(j-i+1\). Their sum is \(k+2\), and each length
lies between two and \(k\). These operations act on the indices of
any family of loops, including the primitive loops defined by
\eqref{fo:eq:primitive}.

Let $\mathcal M_k=\mathcal M_k(t;g,E;\boldsymbol a,\boldsymbol\sigma)$
be the martingale part of $\cL_k(t)$, with
\(\mathcal M_k(0)=0\), and define the centered insertion term by
\begin{equation}\label{fo:eq:light}
 \mathcal W_k\coloneqq\mathcal W_k(t;g,E;\boldsymbol a,\boldsymbol\sigma)
 \coloneqq \gamma_g W\sum_{i=1}^k\sum_c
           \langle\wt G_t^{\sigma_i}\rangle_{cc}\operatorname{Cut}_i^c\cL_k(t).
\end{equation}
This is the centered insertion in \cite[(3.26)]{TYY}, or
\cite[(2.47)]{YY}, with the projected covariance $\gamma_g I_L$.
Equation~\eqref{fo:eq:ito} and the contraction of two differentiated
resolvent factors yield the \emph{loop hierarchy}:
\begin{equation}\label{fo:eq:hierarchy}
 d\cL_k(t)=d\mathcal M_k+\mathcal W_k\,dt
  +\gamma_g W\sum_{i<j}\sum_c
       (\operatorname{Cut}_{L,i,j}^c\cL_k(t))
       (\operatorname{Cut}_{R,i,j}^c\cL_k(t))\,dt.
\end{equation}
This is \cite[Lemma~3.6, (3.24)]{TYY}, corresponding to
\cite[Lemma~2.11, (2.45)]{YY}, with our covariance and block
normalization. The coefficient is \(\gamma_g W\): Gaussian contraction
contributes \(\gamma_g/W\), and each of the two new normalized block
projectors contributes a factor \(W\).
For two loops $\cL_k,\cL_k'$ with martingale parts
$\mathcal M_k,\mathcal M_k'$, put
\[
 B_i=G_t^{\sigma_i}E_{a_i}\cdots
 G_t^{\sigma_{i+k-1}}E_{a_{i+k-1}}G_t^{\sigma_i},
 \qquad 1\le i\le k,
\]
with indices read modulo $k$. Define $B_j'$ in the same way for
the second loop $\cL_k'$. Their exact martingale contraction is
\begin{equation}\label{gc:eq:bracket}
 \frac{d}{dt}\langle\mathcal M_k,\overline{\mathcal M_k'}\rangle
 =W^{-1} \gamma_g\sum_c\sum_{i,j=1}^k\Tr(P_cB_iP_cB_j'^*).
\end{equation}
The proof is in Appendix~\ref{sec:gc-objects}. Replacing $P_c$ by
$WE_c$ makes the right side $\gamma_gW$ times joined
$(2k+2)$-loops with a summed block label.
Consequently \eqref{fo:eq:hierarchy} cannot be closed at a single
length: the drift contains a \((k+1)\)-loop and the bracket contains
still longer loops.

The primitive loops keep the quadratic part of
\eqref{fo:eq:hierarchy} and its exact initial data. With \(a_0\coloneqq a_k\), set
$\cK_1(t;a,\sigma)\coloneqq\cK_1(t;g,E;a,\sigma)=m(g,E)^\sigma$.
For $k\ge2$, write
$\cK_k(t;\boldsymbol a,\boldsymbol\sigma)
=\cK_k(t;g,E;\boldsymbol a,\boldsymbol\sigma)$ and define
\begin{equation}\label{fo:eq:primitive}
 \begin{split}
 \partial_t\cK_k(t)
 &=\gamma_g W\sum_{i<j}\sum_c
       (\operatorname{Cut}_{L,i,j}^c\cK_k(t))
       (\operatorname{Cut}_{R,i,j}^c\cK_k(t)),\\
 \cK_k(0;\boldsymbol a,\boldsymbol\sigma)
 &=\Tr\prod_{j=1}^k(M^{\sigma_j}E_{a_j})
   =W^{1-k}\prod_{j=1}^k(M_B^{\sigma_j})_{a_{j-1},a_j}.
 \end{split}
\end{equation}
The initial value is an \(M\)-loop. Unlike \(G\)-loops, primitive
loops form a triangular deterministic system once the two-loops have
been solved. Equation~\eqref{fo:eq:primitive} is
\cite[Definition~3.7, (3.27)--(3.29)]{TYY} with the cut coefficient
$\gamma_g W$. The corresponding band-matrix definition is
\cite[Definition~2.12, (2.48)]{YY}, with a scalar initial comparison matrix.

\begin{lemma}[Mixed Ward identity]\label{fo:lem:mixed-ward}
Let $\mathcal X\in\{\cL,\cK\}$ be defined by \eqref{fo:eq:loops} or
\eqref{fo:eq:primitive}. For $k\ge2$, $0\le t<1$ with $\eta_t>0$,
and every choice of the remaining signs and labels,
\begin{equation}\label{gc:eq:ward}
 \begin{split}
 &\sum_{c\in\Z_L}\mathcal X_k
    (t;a_1,\ldots,a_{k-1},c;+,\sigma_2,\ldots,\sigma_{k-1},-)\\
 &\quad=(2i)^{-1}W^{-1}\eta_t^{-1}\bigl[
   \mathcal X_{k-1}(t;a_1,\ldots,a_{k-1};+,\sigma_2,\ldots,\sigma_{k-1})\\
 &\hspace{43mm}-
   \mathcal X_{k-1}(t;a_1,\ldots,a_{k-1};-,\sigma_2,\ldots,\sigma_{k-1})
   \bigr].
 \end{split}
\end{equation}
Cyclic permutation gives the identity at any vertex between opposite
signs.
\end{lemma}
The proof, including the block normalization, is in
Appendix~\ref{sec:gc-objects}.

For the two-loop equation, define matrices on block labels by
\begin{equation}\label{fo:eq:PQ}
 P_{ab}\coloneqq P_{ab}(g,E)=\gamma_g|(M_B(g,E))_{ab}|^2,\qquad
 Q_{ab}\coloneqq Q_{ab}(g,E)=\gamma_g(M_B(g,E))_{ab}^2.
\end{equation}
The matrices $P/\gamma_g$ and $Q/\gamma_g$ are the mixed-sign and
same-sign initial pair matrices in \cite[Definition~3.8, (3.31)]{TYY}.
The symmetry of the hopping matrix $A_g$ in \eqref{fo:eq:flow} gives the two-loop formula
\cite[Example~3.12, (3.47)]{TYY} in the form
\begin{equation}\label{fo:eq:pairs}
 \begin{aligned}
 \cK_2(t;a,b;+,-)&=\gamma_g^{-1}W^{-1} [P(I-tP)^{-1}]_{ab},\\
 \cK_2(t;a,b;+,+)&=\gamma_g^{-1}W^{-1} [Q(I-tQ)^{-1}]_{ab}.
 \end{aligned}
\end{equation}
For example, \(\partial_t\cK_2^{+-}(t)=\gamma_g W(\cK_2^{+-}(t))^2\)
and \(\cK_2^{+-}(0)=P/(\gamma_g W)\), which proves the first
formula. The other signs are obtained by conjugation.
The Dyson equation implies \(\sum_bP_{ab}=1\), hence
\(\sum_b[P(I-tP)^{-1}]_{ab}=(1-t)^{-1}\).
The same-sign inverse stays stable at the regular energies under
consideration. Thus only mixed signs produce the growing propagation
length.

For \(k\ge3\), adjacent cuts, including the cyclic pair \((1,k)\),
produce one two-loop and one \(k\)-loop. The other cuts produce two
loops of lengths strictly less than \(k\). Following
\cite[Definition~5.2]{YY}, denote the resulting linear operator on
block tensors by $\Theta_{t,\boldsymbol\sigma}
=\Theta_{t,\boldsymbol\sigma}(g,E)$, with
$\sigma_{k+1}\coloneqq\sigma_1$. This is the generator in
\cite[Definition~7.3, (7.18)]{TYY} with our covariance:
\begin{equation}\label{fo:eq:generator}
\begin{split}
 (\Theta_{t,\boldsymbol\sigma}F)(\boldsymbol a)
 \coloneqq {}&\gamma_g W\sum_{j=1}^k\sum_c
   \cK_2(t;a_j,c;\sigma_j,\sigma_{j+1})\\
 &\hspace{12mm}\cdot
   F(a_1,\ldots,a_{j-1},c,a_{j+1},\ldots,a_k).
 \end{split}
\end{equation}
Then
\[
 \partial_t\cK_k(t)=\Theta_{t,\boldsymbol\sigma}\cK_k(t)+
 \gamma_g W\sum_{\substack{i<j\text{ nonadjacent}}}\sum_c
       (\operatorname{Cut}_{L,i,j}^c\cK_k(t))
       (\operatorname{Cut}_{R,i,j}^c\cK_k(t)).
\]
The source on the right only uses shorter primitive loops.
The linear evolution acts separately on each block label. Its
one-label kernels $U_{s,t}^{\sigma\sigma'}=U_{s,t}^{\sigma\sigma'}(g,E)$,
as in \cite[(7.19)]{TYY} and \cite[(5.17)]{YY}, are
\begin{equation}\label{fo:eq:propagators}
 \begin{gathered}
 U_{s,t}^{+-}=U_{s,t}^{-+}\coloneqq (I-sP)(I-tP)^{-1},\\
 U_{s,t}^{++}\coloneqq (I-sQ)(I-tQ)^{-1},\qquad
 U_{s,t}^{--}\coloneqq \overline{U_{s,t}^{++}}.
 \end{gathered}
\end{equation}
Here $P_{ab}=\gamma_g|(M_B)_{ab}|^2$ and
$Q_{ab}=\gamma_g(M_B)_{ab}^2$, as in \eqref{fo:eq:PQ}.
For instance,
\(U_{s,t}^{+-}=I+(t-s)P(I-tP)^{-1}\) and
\(\partial_tU_{s,t}^{+-}=P(I-tP)^{-1}U_{s,t}^{+-}\).
With
\(\mathcal U_{s,t}^{\boldsymbol\sigma}
 =\mathcal U_{s,t}^{\boldsymbol\sigma}(g,E)\coloneqq \bigotimes_{j=1}^kU_{s,t}^{\sigma_j\sigma_{j+1}}\),
Duhamel's principle gives the primitive-loop counterpart of
\cite[(7.21)]{TYY}, with the deterministic forcing formed by the
nonadjacent cuts; see \cite[Lemma~5.3]{YY} for the general
variation-of-constants formula:
\begin{equation}\label{fo:eq:duhamel}
 \begin{split}
 \cK_k(t)={}&\mathcal U_{s,t}^{\boldsymbol\sigma}\cK_k(s)\\
 &+\gamma_g W\int_s^t\mathcal U_{u,t}^{\boldsymbol\sigma}
  \sum_{\substack{i<j\text{ nonadjacent}}}\sum_c
       (\operatorname{Cut}_{L,i,j}^c\cK_k(u))
       (\operatorname{Cut}_{R,i,j}^c\cK_k(u))\,du.
 \end{split}
\end{equation}
For \(k=3\) the source is empty. Iterating
\eqref{fo:eq:duhamel} expresses every fixed-order primitive in terms
of its initial \(M\)-loops and the two-loop propagators. The tree
representation in \cite[Lemma~4.16]{TYY}, extending
\cite[Lemma~3.4]{YY} to initial $M$-loops with $M$ not proportional to $I_N$, organizes this expansion.
Its size cannot in general be estimated by taking absolute values
before the cancellations between its terms.

Write $D_k=D_k(t;g,E;\boldsymbol a,\boldsymbol\sigma)
\coloneqq\cL_k(t)-\cK_k(t)$. Subtraction gives the exact error
equation \cite[(7.14)--(7.15)]{TYY}, corresponding to
\cite[(5.12)--(5.13)]{YY}, with cut coefficient $\gamma_g W$:
\begin{equation}\label{fo:eq:error-hierarchy}
\begin{split}
 dD_k={}&d\mathcal M_k+\mathcal W_k\,dt\\
 &+\gamma_g W\sum_{i<j}\sum_c\bigl[
   (\operatorname{Cut}_{L,i,j}^c\cK_k(t))
       (\operatorname{Cut}_{R,i,j}^cD_k)\\
 &\hspace{27mm}+(\operatorname{Cut}_{L,i,j}^cD_k)
       (\operatorname{Cut}_{R,i,j}^c\cK_k(t))\\
 &\hspace{27mm}+(\operatorname{Cut}_{L,i,j}^cD_k)
       (\operatorname{Cut}_{R,i,j}^cD_k)\bigr]dt.
 \end{split}
\end{equation}
For \(k\ge3\), its adjacent linear terms are again generated by
\(\Theta_{t,\boldsymbol\sigma}\).
The remaining terms consist of shorter linear errors, quadratic
errors, centered insertions, and the martingale. This decomposition
is the basis of the time induction in Appendix~\ref{sec:gc-schedule}.

The spatial scale of the mixed kernel depends on the coupling.
The bounded-coupling scale has the form of \cite[(3.42)]{TYY};
the large-coupling scale follows from the kernel estimates in
Lemma~\ref{uvp:lem:actual-kernels}. For the two extensions considered here we use
\begin{equation}\label{fo:eq:scale}
 \ell_t\coloneqq\ell_t(g)
 =\min\left\{L,1+\frac{g}{\sqrt{\gamma_g(1-t)}}\right\},\qquad
 \zeta_t\coloneqq\zeta_t(g)=\frac1{ W\ell_t(g)(1-t)\gamma_g}.
\end{equation}
The length \(\ell_t\) counts blocks, whereas \(\ellstar\) counts
coordinates. In the parameter regimes of Lemma~\ref{fo:lem:deterministic}
below, \(\Im m\asymp1\), so
\(\zeta_t\asymp(W\ell_t\eta_t)^{-1}\).
Directly from the definition, $\zeta_t\ge c/W$ and
$0\le\partial_t\log\zeta_t\le(1-t)^{-1}$ almost everywhere,
including the capped scale; hence
$\zeta_s\le\zeta_t\le(1-s)(1-t)^{-1}\zeta_s$ for $0\le s\le t<1$.
Under the volume assumptions of Theorem~\ref{gc:thm:common},
$\ell_t\asymp\max\{1,g^2\}$ and $\zeta_t\asymp W^{-1}$ for $0\le t\le1/2$.

\begin{lemma}[Size of the primitive loops]\label{fo:lem:deterministic}
Fix an integer $k\ge2$. Consider either of the following parameter regimes for deterministic
$g=g(W)$ and cycle length $L=L(W)$.

\textup{(a) Bounded coupling.} Fix $C,c>0$ and assume $|E|\le C$ and
$\Im m_{g,L}(E)\ge c$ and 
\[
 0<g\le C,\qquad \log W\ll L\le W^2.
\]

\textup{(b) Large coupling.} Fix $\kappa\in(0,2)$.
There exist $C_\kappa>1$, independent of $k$, and $C_*=C_*(\kappa,k)>0$
such that the size estimates below hold when $|E|\le2-\kappa$ and
\begin{equation}\label{fo:eq:polynomial-large-range}
 g\ge C_\kappa,\qquad L\ge C_*g^2\log W.
\end{equation}
For sufficiently large $W$, uniformly on
$0\le t<1$, $1-t\ge W^{-3}$, the primitive loops of length $k\ge3$ satisfy
\begin{equation}\label{fo:eq:primitivebound}
 |\cK_k(t)|\le C'\zeta_t^{k-1},
\end{equation}
for all block labels and sign words. In case \textup{(a)}, $C'$ depends
only on $C,c,k$; in case \textup{(b)}, only on $\kappa,k$, with no
upper bound on $g$ or $L$. At length two,
\[
 \max_{a,b,\sigma_1,\sigma_2}
       |\cK_2(t;a,b;\sigma_1,\sigma_2)|\le C'\zeta_t.
\]
For every fixed $\epsilon,Q>0$, under case \textup{(a)}, or under
case \textup{(b)} with the stronger volume condition $L\gg g^2\log W$,
\begin{equation}\label{fo:eq:primitive-tail}
 \sup_{0\le t\le1-W^{-3}}\max_{\boldsymbol\sigma\in\{+,-\}^k}
 \sup_{\diam(\boldsymbol a)>W^\epsilon\ell_t}
       |\cK_k(t;\boldsymbol a,\boldsymbol\sigma)|\le W^{-Q}.
\end{equation}
In case \textup{(b)}, the constants are independent of $g,L$; the width
threshold for \eqref{fo:eq:primitive-tail} may depend on the divergence
rate of $L/(g^2\log W)$, uniformly over families with a common rate.
\end{lemma}

The lower bound $L\ge C_*g^2\log W$ controls the finite-cycle comparison
at large coupling. Lemma~\ref{uv:lem:short-spreading} covers the
complementary short cycles in the proof of the lower bound in
Theorem~\ref{thm:main}.
Primitive-loop estimates are a central ingredient in
Yau--Yin~\cite[Section~3 and Lemma~2.17]{YY} and
Truong--Yang--Yin~\cite[Section~4 and Lemma~3.15]{TYY}.
Their method combines a tree representation, the primitive Ward identity
and sum-zero cancellation; estimates for individual molecule cores then close under tree
gluing. Our main task here is to extend these estimates uniformly to
arbitrarily large $g$. This requires additional ideas to handle the long
decay scale of $M$ without a $\log g$ loss, but the proof follows the
same underlying method of YY and TYY.
The size estimates in cases \textup{(a)} and \textup{(b)} are proved in
Appendices~\ref{app:bounded-primitives} and~\ref{app:uniform-profile},
respectively; \eqref{fo:eq:primitive-tail} is proved in
Appendix~\ref{app:primitive-tails}.

\noindent\textbf{Definition (Stochastic domination).}
For random variables $X$ and $Y\ge0$ depending on $W$, we write
$X\prec Y$ if for every $\epsilon,D>0$,
\[
 \Prob\bigl(|X|>W^\epsilon Y\bigr)\le W^{-D}
\]
for sufficiently large $W$.
Unless otherwise stated, the estimate is uniform over the deterministic
parameters in the statement.

The basic estimates to be propagated along the flow are the following,
written at a deterministic time $0\le s<1$, with fixed
$\theta\in(0,1/2]$ and $c_{\rm p}>0$:
\begin{equation}\label{gc:eq:one-step-input}
 \begin{gathered}
 \max_{\boldsymbol a,\boldsymbol\sigma}\frac{|\cL_k(s)|}{\zeta_s^{k-1}}\prec1,
 \qquad \max_{\boldsymbol a,\boldsymbol\sigma}\frac{|D_k(s)|}{\zeta_s^k}\prec1,\qquad
 \max_{a,b}\frac{|\langle\wt G_s\rangle_{ab}|}{\zeta_s}\prec1,\qquad
 \frac{\|\wt G_s\|_{\max}}{\sqrt{\zeta_s}}\prec1,\\
 \max_{a,b}\frac{|D_2(s;a,b;+,-)|}
 {\zeta_s^2(e^{-c_{\rm p}(\|a-b\|_L/\ell_s)^\theta}+W^{-F})}\prec1,\\
 \max_{\boldsymbol\sigma\in\{+,-\}^k}
 \max_{\boldsymbol a:\,\diam(\boldsymbol a)>W^\epsilon\ell_s}
       \bigl(|\cL_k(s)|+|D_k(s)|\bigr)\prec W^{-A},\qquad k\ge2.
 \end{gathered}
\end{equation}
They are required for every fixed loop length $k$ and
$\epsilon,A,F>0$, uniformly over the stated deterministic parameters.

At bounded coupling, $\eta_u\asymp1$ uniformly for $0\le u\le1/2$.
We prove the following initial input by adapting the resolvent and stability
method of \cite[Section~8 and Appendix~B.1]{DYYY}.

\begin{lemma}[Initial loop estimates at bounded coupling]
\label{fo:lem:bounded-initial}
Assume the bounded-coupling conditions of
Lemma~\ref{fo:lem:deterministic}(a). There exists $c_{\rm p}>0$ such that
\begin{equation}\label{gc:eq:iteration-input}
 \left.\eqref{gc:eq:one-step-input}\right|_{s=u,\ \theta=1/2},
 \quad 0\le u\le1/2,
 \qquad \sup_{0\le u\le1/2}\max_{a,b}\|(G_u)_{ab}\|_{\op}\prec1.
\end{equation}
The time suprema are inside the probability, and all estimates are
uniform over the stated deterministic parameters.
\end{lemma}
The proof is in Appendix~\ref{sec:bounded-proof}.
Starting from this input at $s=1/2$, the flow propagates
\eqref{gc:eq:one-step-input} to the terminal spectral height $\eta_L$.
For large $g$, already the initial height is of order $g^{-2}$,
and transfer matrices are needed. The main new inputs for this regime
are presented in Section~\ref{sec:fo-transfer-hierarchy}.

For spatial interpolation of primitives, we use an auxiliary flow with
orbital dimension $\wt W$; the parameters of the original flow remain fixed.

\begin{lemma}[Auxiliary orbital limit]\label{uvc:lem:orbital-limit}
Fix $g>0$, $L,W\in\mathbb N$, $E\in\mathbb R$ and $0\le T<1$.
Assume that the Dyson solution $m=m_{g,L}(E)$ in
\eqref{fo:eq:characteristic} satisfies $\Im m>0$.
For each $\wt W\in\mathbb N$, let $\mathcal B_t^{(g,\wt W)}$ be block diagonal,
with independent Hermitian Brownian blocks starting at zero and covariance
\[
 \E\,d(\mathcal B_t^{(g,\wt W)})_{ai,aj}\,
       d(\mathcal B_t^{(g,\wt W)})_{bk,b\ell}
 =\frac{\gamma_g}{\wt W}\delta_{ab}\delta_{i\ell}\delta_{jk}\,dt,
 \qquad a,b\in\Z_L,\quad i,j,k,\ell\in\{1,\ldots,\wt W\}.
\]
Using $A_g$ and $\gamma_g$ from \eqref{fo:eq:flow}, define
\[
 H_t^{(\wt W)}=A_g\otimes I_{\wt W}+\mathcal B_t^{(g,\wt W)},\qquad
 G_t^{(\wt W)}=(H_t^{(\wt W)}-z_tI_{L\wt W})^{-1},\qquad
 z_t=E+(1-t)\gamma_g m.
\]
Let $P_a^{(\wt W)}$ project onto block $a$ and set
$E_a^{(\wt W)}=P_a^{(\wt W)}/\wt W$,
$G_t^{(\wt W),+}=G_t^{(\wt W)}$ and
$G_t^{(\wt W),-}=(G_t^{(\wt W)})^*$.
With the unnormalized trace, put
\[
 \cL_k^{(\wt W)}(t;\boldsymbol a,\boldsymbol\sigma)
 =\Tr\bigl(G_t^{(\wt W),\sigma_1}E_{a_1}^{(\wt W)}\cdots
           G_t^{(\wt W),\sigma_k}E_{a_k}^{(\wt W)}\bigr).
\]
For every fixed $k\ge1$, with $\cK_k$ from \eqref{fo:eq:primitive},
\begin{equation}\label{uvc:eq:orbital-limit}
 \lim_{\wt W\to\infty}\sup_{0\le t\le T}
 \max_{\boldsymbol a\in\Z_L^k,\,\boldsymbol\sigma\in\{+,-\}^k}
 \left|\wt W^{k-1}\E\cL_k^{(\wt W)}(t;\boldsymbol a,\boldsymbol\sigma)
       -W^{k-1}\cK_k(t;\boldsymbol a,\boldsymbol\sigma)\right|=0.
\end{equation}
The parameters $W,g,L,E,T$, and hence $m$ and $z_t$, are fixed in this
limit; no convergence rate uniform in $g,L$ is asserted.
\end{lemma}
The proof is in Appendix~\ref{app:loop-identities}.

\subsection{Applications of transfer matrices to the loop hierarchy}
\label{sec:fo-transfer-hierarchy}
This subsection presents the main new inputs to the lower-bound argument:
transfer-matrix resolvent moments and initialization
(Lemmas~\ref{fo:lem:resolvent-moments} and~\ref{fo:prop:inputs}),
and spatial sampling of random and primitive loops
(Lemma~\ref{uvc:thm:ordinary-sampling}). We also state the residual
and continuity estimates needed to use these inputs in the loop hierarchy.
Their large-coupling proofs are given in Section~\ref{sec:large-proof},
independently of the propagation theorem below.

\begin{lemma}[Resolvent moments at large coupling]
\label{fo:lem:resolvent-moments}
For every fixed $0<\kappa<2$ and $p,C_0>0$, if
$g\ge C_\kappa$, $L\ge4g^2$, $|E|\le2-\kappa$, and
$0<\eta\le C_0/g^2$, then
\begin{equation}\label{fo:eq:coarse-resolvent-moments}
 \sup_{0\le t\le1}\max_{a,b}
 \E\left\|\bigl[(H_t(g)-E-i\eta)^{-1}\bigr]_{ab}\right\|_{\op}^p
 \le C_p(g^2\eta)^{-p}.
\end{equation}
The constants and width threshold depend only on the fixed parameters;
there is no upper bound on $g$.
\end{lemma}
The proof in Section~\ref{sec:uniform-resolvent-seed} combines transfer
estimates with Ward's identity. These moments supply both the
initialization below and the coarse control used during propagation.

\begin{lemma}[Initial loop estimates at large coupling]
\label{fo:prop:inputs}\label{uvu:fo:prop:inputs}
Fix $0<\kappa<2$ and $C\ge2$.
For a sufficiently large $C_\kappa$, assume
\[
 g\ge C_\kappa,\qquad L\gg g^2\log W,\qquad
 L\le Cg^2W^2,\qquad |E|\le2-\kappa.
\]
Uniformity is over deterministic parameter families for which
$L/(g^2\log W)$ diverges uniformly; the width threshold may depend on its rate.
There exists $c_{\rm p}>0$ such that
\eqref{gc:eq:iteration-input} holds.
The time suprema are inside the probability, and there is no upper bound on $g$.
\end{lemma}
The proof is in Section~\ref{sec:uniform-initialization}.

\begin{definition}[Cell partitions $\mathcal P_B$]
\label{def:cells-anchors}
A deterministic partition of the block cycle $\Z_L$ into nonempty
consecutive intervals is a \emph{cell partition}, and its intervals are
called \emph{cells}.
For fixed $g,L,W$ and $B>0$, set $n=\max\{1,\lfloor g^2W^{-B}\rfloor\}$
and fix a deterministic partition $\mathcal P_B(g,L,W)$ into
$\lceil L/n\rceil$ cells whose lengths differ by at most one.
Each cell then has length at most $n$.
\end{definition}

\begin{lemma}[Spatial sampling of loop errors]
\label{uvc:thm:ordinary-sampling}
Fix $0<\kappa<2$, $C>1$, $D>0$, and $K\in\mathbb N$.
Choose a fixed $B>2$ sufficiently large depending on $D$ and the other
fixed parameters. The following holds for all sufficiently large $W$.
Assume that $g,L,E$ are deterministic and
\[
 \begin{gathered}
 4g^2\le L\le g^2W^C,\\
 |E|\le2-\kappa,\qquad
 |m_{g,L}(E)|\le C,\quad \Im m_{g,L}(E)\ge C^{-1}.
 \end{gathered}
\]
Use the partition $\mathcal P_B$ from Definition~\ref{def:cells-anchors};
it satisfies $|\mathcal P_B|\le3W^{B+C}$.
For each $1\le k\le K$, choose
at most $d_k\coloneqq6^k$ deterministic points from each product $I_1\times\cdots\times I_k$
of cells $I_1,\ldots,I_k\in\mathcal P_B$,
and let $\mathcal S_k$ be their union.
The points can be chosen independently of time and of the sign word
so that, with probability at least $1-W^{-D}$,
\begin{equation}\label{fo:eq:spatial-sampling}
 \max_{\boldsymbol a\in\Z_L^k}
       |D_k(t;\boldsymbol a,\boldsymbol\sigma)|
 \le C_k\max_{\boldsymbol a\in\mathcal S_k}
       |D_k(t;\boldsymbol a,\boldsymbol\sigma)|+W^{-D},
 \qquad D_k=\cL_k-\cK_k,
\end{equation}
simultaneously for $1\le k\le K$, all sign words $\boldsymbol\sigma$,
and $0\le t\le1-W^{-3}$. The same assertion holds with $D_k$ replaced
by $\cL_k$. Here $C_k$ depends only on $k$; all other constants and the
width threshold depend only on the fixed parameters, not on $g,L,E$.
\end{lemma}

If $g^2<W^B$, all cells are singletons and $L\le W^{C+B}$, so taking
$\mathcal S_k=\Z_L^k$ gives exact sampling. The proof in
Section~\ref{app:uniform-compression} therefore only needs interpolation
when $g^2\ge W^B$. Green-block interpolation, together with an auxiliary
orbital limit, places random and primitive loops in a common approximate
finite-dimensional space on each product of cells. Sampling points for
this space are chosen within that product. The argument does not use
the sharp local law being proved.

For the flow \eqref{fo:eq:flow}, the Dyson quantities at time
$0\le u<1$ and spectral point $z\in\C_+$, together with their values
on the characteristic \eqref{fo:eq:characteristic}, are
\begin{equation}\label{gc:eq:positive-dyson}
 \begin{gathered}
 M_{B,u}(z;g)\coloneqq(A_g-z-u\gamma_gm_u(z;g))^{-1},\quad
 m_u(z;g)\coloneqq L^{-1}\Tr M_{B,u}(z;g),\\
 M_u(z;g)\coloneqq M_{B,u}(z;g)\otimes I_W,\qquad
 \Im m_u(z;g)>0,\\
 m_u(z_u;g)=m(g,E),\qquad M_{B,u}(z_u;g)=M_B(g,E).
 \end{gathered}
\end{equation}
For the next two lemmas and the accompanying discussion, suppress $g$
in these Dyson quantities and write
$G=(H_u(g)-z)^{-1}$ and $\wt G=G-M_u(z)$.
Following \cite[definition preceding (6.37)]{TYY}, define the residual by
\begin{equation}\label{gc:eq:residual-input}
 \begin{aligned}
 \Pi(G)&\coloneqq I_N+(z-A_g\otimes I_W)G
              +u\gamma_g\diag(\langle G\rangle_{aa}I_W)G\\
 &=\mathcal B_u^{(g)}G+u\gamma_g\diag(\langle G\rangle_{aa}I_W)G,
 \qquad R_{ab}\coloneqq\langle M_u(z)\Pi(G)\rangle_{ab}.
 \end{aligned}
\end{equation}
The first expression gives $\Pi(M_u(z))=0$; the second follows from
$(H_u-z)G=I_N$. For $a=b$, $R_{aa}$ is the averaged residual in
\cite[(6.50)]{TYY}.

The next two lemmas organize the uniform coarse moments of
Lemma~\ref{fo:lem:resolvent-moments} and the spatial interpolation
underlying Lemma~\ref{uvc:thm:ordinary-sampling} into propagation inputs.
For arbitrarily large $g$, the coarse moments control exceptional-event
contributions in fluctuation averaging, while the interpolation reduces
residual maxima to polynomially many sampled pairs.

\begin{lemma}[Fluctuation averaging]\label{gc:lem:residual}
Use the parameter ranges of Lemma~\ref{fo:lem:bounded-initial} or
Lemma~\ref{fo:prop:inputs}, with fixed energy and density margins.
Fix deterministic $0\le u\le1-W^{-3}$ and take $z=z_u(g,E)$.
For deterministic $1/2\le s\le u$, also allow any deterministic point
\[
 z=\Re z_u+iy,\qquad
 \eta_u\le y\le\xi^{-1}\eta_s(g',E'),\qquad
 \xi=\sqrt{\frac{s\gamma_{g'}}{u\gamma_g}},
\]
where $(g',E')$ is supplied by Lemma~\ref{fo:lem:families}.

\noindent\textup{(i) Averaged residual.}
For deterministic polynomially bounded $K\ge1$ and $\chi\ge W^{-1/2}$,
\begin{equation}\label{gc:eq:averaging-input}
 \max_{a,b}|\langle G\rangle_{ab}|\prec K,\quad
 \max_{ai,bj}|G_{ai,bj}-\delta_{ij}\langle G\rangle_{ab}|\prec K\chi
 \quad\Longrightarrow\quad
 \max_{a,b}|R_{ab}|\prec K^2\chi^2.
\end{equation}
\noindent\textup{(ii) Entrywise estimate.}
For each pair $a,b$ and deterministic $\Psi>0$,
\begin{equation}\label{gc:eq:haar}
 W^{-2}\|G_{ab}\|_{\rm HS}^2\prec\Psi^2
 \quad\Longrightarrow\quad
 \max_{i,j}|G_{ai,bj}-\delta_{ij}\langle G\rangle_{ab}|\prec\Psi.
\end{equation}
\noindent\textup{(iii) Residual moments.}
For each fixed even $p\ge2$ and $C>0$, there is $C_p$ such that,
for deterministic $\chi\ge W^{-1/2}$,
\begin{equation}\label{fo:eq:residual-moments}
 \max_{ai,bj}\|\wt G_{ai,bj}\|_{L^{2p}}\le C\chi
 \quad\Longrightarrow\quad
 \max_{a,b}\E|R_{ab}|^p\le C_p\chi^{2p}.
\end{equation}
All estimates are uniform in the allowed deterministic parameters.
Probabilities and moments use the original Gaussian measure, without
conditioning on a stopping event.
\end{lemma}

On the characteristic $z=z_u$, resolvent subtraction gives, with
$Q$ from \eqref{fo:eq:PQ},
\begin{equation}\label{gc:eq:scalar}
 \begin{gathered}
 \langle\wt G\rangle_{ab}
 =-R_{ab}+u\gamma_g\sum_cM_{B,ac}
       \langle\wt G\rangle_{cc}\langle G\rangle_{cb},\\
 (I-uQ)(\langle\wt G\rangle_{aa})_a=-(R_{aa})_a
       +u\gamma_g\left(\sum_cM_{B,ac}
         \langle\wt G\rangle_{cc}\langle\wt G\rangle_{ca}\right)_a.
 \end{gathered}
\end{equation}
These identities convert residual estimates into the averaged local law.

The proof in Section~\ref{sec:input-averaging} uses the Gaussian
calculations in Appendix~\ref{sec:classical-random-matrix}, following
\cite[(6.40), (6.50)]{TYY} and the averaged improvement in
\cite[proof of (4.5), (4.12)]{YY}.
For the random centering $\langle G\rangle_{ab}I_W$ in part (i),
that appendix gives the block and orbital index counts, pairing the
covariance factors $u\gamma_g$ with absolute or squared Dyson row sums.
Part (ii) uses common Haar conjugation.

For the next lemma, the same coarse moments yield time-increment bounds
through It\^o's formula and BDG. Unweighted sampling of loops and loop
errors follows directly from Lemma~\ref{uvc:thm:ordinary-sampling}.
Weighted loop sampling and sampling of centered block traces and
traceless entries use the same interpolation mechanism; residual sampling
requires the additional proof in
Proposition~\ref{uvc:prop:residual-compression}.

\begin{lemma}[Sampling and continuity for propagation]
\label{fo:lem:sampling-continuity}
Use the coupling, volume and regularity assumptions of
Lemma~\ref{gc:lem:residual}. Fix a finite set of loop lengths and $D,D'>0$.
The following estimates are uniform in the deterministic parameters,
for $0\le u,v\le1-W^{-3}$.

\textup{(a) Coarse moments and time increments.}
Along the fixed characteristic $(g,E)$, let $X(u)$ be any of
\[
 \begin{gathered}
 (G_u)_{ab},\qquad \langle\wt G_u\rangle_{ab},\qquad
 (G_u)_{ai,bj}-\delta_{ij}\langle G_u\rangle_{ab},\\
 \cL_k(u;\boldsymbol a,\boldsymbol\sigma),\qquad
 D_k(u;\boldsymbol a,\boldsymbol\sigma),\qquad R_{ab}(u).
 \end{gathered}
\]
The block indices $a,b,\boldsymbol a$, orbital indices $i,j$, and signs
$\boldsymbol\sigma$ are fixed when comparing $u$ and $v$.
Here $R_{ab}(u)$ is the residual in \eqref{gc:eq:residual-input} evaluated
at $z=z_u(g,E)$. Use operator norm for $(G_u)_{ab}$ and absolute value
for the other five quantities.
There is a fixed exponent $C$ such that, for every fixed $p\ge2$,
\begin{equation}\label{fo:eq:coarse-time-input}
 \max_{a,b}\E\|(G_u)_{ab}\|_{\op}^p\le C_pW^{Cp},\qquad
 \E\|X(v)-X(u)\|^p\le C_pW^{Cp}|v-u|^{p/2}.
\end{equation}
The exponent may depend on the fixed set of loop lengths and regularity
margins, but not on $p,g,L$ in the arbitrary-coupling range.

\textup{(b) Spatial suprema.}
For fixed signs and orbital indices, choose one of the scalar fields
\[
 \begin{aligned}
 X(u;\boldsymbol a)&\in\{\cL_k(u;\boldsymbol a,\boldsymbol\sigma),
                   D_k(u;\boldsymbol a,\boldsymbol\sigma)\},
                   &&\boldsymbol a\in\Z_L^k,\\
 X(u;\boldsymbol a)&\in\{\langle\wt G_u\rangle_{ab},
   (G_u)_{ai,bj}-\delta_{ij}\langle G_u\rangle_{ab},R_{ab}(u)\},
                   &&\boldsymbol a=(a,b)\in\Z_L^2.
 \end{aligned}
\]
Take $w_u=1$; for loops and loop errors one may instead fix $c,F>0$,
$0<\theta\le1/2$ and take
$w_u(\boldsymbol a)=e^{-c(\diam(\boldsymbol a)/\ell_u)^\theta}+W^{-F}$.
There is a deterministic set $\mathcal S_X$ of at most $W^C$ block tuples,
independent of time and allowed to depend on the fixed weight parameters,
$D,D'$, such that, with probability at least $1-W^{-D'}$, simultaneously
for all $0\le u\le1-W^{-3}$,
\begin{equation}\label{fo:eq:weighted-sampling-input}
 \max_{\boldsymbol a}\frac{|X(u;\boldsymbol a)|}{w_u(\boldsymbol a)}
 \le C\max_{\boldsymbol a\in\mathcal S_X}
             \frac{|X(u;\boldsymbol a)|}{w_u(\boldsymbol a)}+W^{-D}.
\end{equation}
For joint-tail restrictions, sampling is cellwise and enlarges the
diameter boundary by at most twice the cell length, which is
$o(W^\epsilon\ell_u)$ for every fixed $\epsilon>0$.
For entries there is an additional union over the $W^2$ orbital indices.
In polynomial volume one may take all block tuples. In the
arbitrary-coupling range, the constants are independent of $g,L$;
uniformity of the width threshold is as in Lemma~\ref{fo:prop:inputs}.

At each fixed deterministic $(u,z)$ on the comparison segments in
Lemma~\ref{gc:lem:residual}, the first bound in
\eqref{fo:eq:coarse-time-input} holds for $G=(H_u(g)-z)^{-1}$.
At this point, \eqref{fo:eq:weighted-sampling-input} with $w=1$ holds
for $\langle G-M_u(z)\rangle_{ab}$,
$G_{ai,bj}-\delta_{ij}\langle G\rangle_{ab}$ and $R_{ab}$;
the deterministic sample sets may depend on $(u,z)$.
\end{lemma}
The proof is in Section~\ref{sec:sampling-propagation-extensions}.
At each fixed time, the derivative bound \eqref{fo:eq:vertical-modulus}
is used to extend estimates from a discrete height mesh to the full
vertical comparison segment.
In propagation, spatial sampling is applied to the actual fields
after adapted stopping and moment estimates, not to their martingale
differentials.

\subsection{Propagation of the loop estimates}
\label{sec:fo-inputs}\label{sec:fo-continuation}

We propagate the estimates along a finite sequence of times
$1/2=s_0<s_1<\cdots<1$, using the estimates obtained at each step
as the input for the next. The rough bounds at the start of a step for
$G_{s_{i+1}}(g,E)$ come from the earlier resolvent
$G_{s_i}(g',E')$, with adjusted coupling and energy.
Lemma~\ref{fo:lem:families} gives the relation between $g'$ and $g$
and controls $E'$, so that the new characteristic retains the required
bulk regularity. In particular, $E'$ remains a bulk anchor for coupling
$g'$, with $\Im m_{g',L}(E')$ bounded below by a positive constant.

%
\begin{lemma}[Resolvent comparison at two times]\label{fo:lem:families}
Fix deterministic $1/2\le s\le u\le1-W^{-3}$. Use \eqref{fo:eq:flow},
\eqref{fo:eq:characteristic} and \eqref{fo:eq:scale}, and set
$$g'\coloneqq g\sqrt{s/u},\quad R\coloneqq(1-s)/(1-u).$$
\begin{enumerate}
\item[(a)] Assume the coupling and volume bounds of
Lemma~\ref{fo:lem:deterministic}(a), and suppose
$|E|\le C'$ and $\Im m_{g,L}(E)\ge c'$ for fixed $C',c'>0$.
There is a constant $C''>0$, depending only on the fixed coupling and
regularity bounds, and a deterministic $E'\in\R$ such that
\begin{equation}\label{fo:eq:backward-matching}
 \begin{gathered}
 \Im G_u(g,E)\overset{\mathrm d}{\preceq} C''R\,\Im G_s(g',E'),\\
 \max_{a,b}\|(G_u(g,E))_{ab}\|_{\op}
 \overset{\mathrm d}{\le} C''R\max_{a,b}\|(G_s(g',E'))_{ab}\|_{\op},\\
 \max_{a,b}|\langle G_u(g,E)\rangle_{ab}|
 \overset{\mathrm d}{\le} C''R\max_{a,b}|\langle G_s(g',E')\rangle_{ab}|,\\
 \zeta_s(g')\asymp\zeta_s(g),\qquad
 \zeta_s(g')\le\zeta_u(g), \qquad
 |E'|\le C'',\qquad \Im m_{g',L}(E')\ge(C'')^{-1}.
 \end{gathered}
\end{equation}
Here $\mathrm d$ means that a coupling exists under which the three
comparisons hold simultaneously, and $\preceq$ denotes the
positive-semidefinite order.

\item[(b)] Assume the coupling and volume bounds of
Lemma~\ref{fo:lem:deterministic}(b), with $|E|\le2-\kappa$,
where $C_\kappa$ is sufficiently large and $C_*$ is chosen after the
prescribed finite orders and accuracy exponents.
There are a constant $C''>0$ and a deterministic $E'\in\R$
such that the three resolvent comparisons and the scale bounds in
\eqref{fo:eq:backward-matching} hold, together with
\begin{equation}\label{fo:eq:large-backward-domain}
 |E'|\le2-\kappa/2,\qquad
 |E'-E|\le C''\bigl((g')^{-2}-g^{-2}\bigr).
\end{equation}
\end{enumerate}
\end{lemma}
The proof is given in Appendix~\ref{app:matching-parameters}.
Following the two-time comparison of \cite[Lemma~2.8]{YY} and
\cite[Lemma~3.3]{TYY}, we first use Gaussian scaling to identify the two
matrix laws up to a bounded factor. We then choose $E'$ so that the corresponding spectral points
have the same real part and comparable imaginary parts, with loss at most
$CR$. Spectral calculus and Ward's identity yield the three resolvent
comparisons.

On a time step $[s,t]$, Lemma~\ref{fo:lem:families} supplies rough
bounds with losses in powers of $R=(1-s)/(1-t)$.
Theorem~\ref{gc:thm:common}(i) uses the loop hierarchy to obtain sharp
bounds that serve as the input for the next step; part (ii) iterates
these bounds to the target time.

\begin{theorem}[Propagation of the loop estimates]\label{gc:thm:common}
Use the flow, loops and scales in
\eqref{fo:eq:flow}--\eqref{fo:eq:scale}, in the parameter ranges of
Lemma~\ref{fo:lem:bounded-initial} or Lemma~\ref{fo:prop:inputs}.
Below, case \textup{(a)} denotes
bounded coupling and case \textup{(b)} denotes large coupling, including
arbitrarily large $g$.
Assume the primitive-loop estimates of Lemma~\ref{fo:lem:deterministic}
and the following inputs from Section~\ref{sec:fo-transfer-hierarchy},
in their respective parameter ranges:
\begin{itemize}
\item the transfer-matrix resolvent moment bounds of
Lemma~\ref{fo:lem:resolvent-moments}, at large coupling;
\item the spatial sampling estimates for random loops and loop errors in
Lemma~\ref{uvc:thm:ordinary-sampling}, at large coupling;
\item the fluctuation-averaging and entrywise estimates of
Lemma~\ref{gc:lem:residual};
\item the coarse moments, weighted spatial sampling and time-continuity
estimates of Lemma~\ref{fo:lem:sampling-continuity}, including its
estimates at the comparison spectral points.
\end{itemize}
All inputs are required uniformly over the deterministic parameters and
comparison points used below, with fixed regularity margins.
There is no separate polynomial upper bound on $g$ or $L$.

\noindent\textup{(i) One time step.}
Fix $c,c_{\rm p}>0$, $\theta\in(0,1/2]$ and deterministic times
$1/2\le s<t<1$. Take a deterministic pair $(g,E)$ in the range
\begin{equation}\label{gc:eq:step-parameters}
 \textup{(a)}\quad |E|\le C',\quad \Im m_{g,L}(E)\ge c';
 \qquad \textup{(b)}\quad |E|\le2-\kappa/2,
\end{equation}
where $C'>0$ and $c'\in(0,1)$ are fixed in case \textup{(a)}.
Assume \eqref{gc:eq:one-step-input} at $(g,E)$ for every fixed $k$ and
$\epsilon,A,F>0$. The times satisfy
\begin{equation}\label{gc:eq:one-step}
 1-t\ge W^{-3},\qquad
 \left(\frac{1-s}{1-t}\right)^{20}\zeta_t\le W^{-c}.
\end{equation}
We also assume
\[
 \left.\eqref{gc:eq:one-step-input}\right|_{g\to g',\,E\to E'}
 \quad\text{uniformly for }s\le u\le t,
\]
where Lemma~\ref{fo:lem:families} supplies $(g',E')$ for the given
$g,E,s,u$.
Then there exists a constant $c_{\rm p}'>0$ such that
\begin{equation}\label{gc:eq:one-step-output}
 \left.\eqref{gc:eq:one-step-input}\right|_{
 s\to u,\,\theta\to\theta/2,\,c_{\rm p}\to c_{\rm p}'},
 \qquad u\in[s,t].
\end{equation}
The supremum over $u\in[s,t]$ of each
normalized left-hand side is taken inside the probability; the supremum
over the allowed deterministic $(g,E)$ is outside.

\medskip

\noindent\textup{(ii) Iteration.}
Assume in addition the initial loop estimates of
Lemma~\ref{fo:lem:bounded-initial} or Lemma~\ref{fo:prop:inputs},
uniformly at the required comparison pairs.
The initial pair $(g,E)$ satisfies
\begin{equation}\label{gc:eq:initial-parameters}
 \textup{(a)}\quad |E|\le C',\quad \Im m_{g,L}(E)\ge c';
 \qquad \textup{(b)}\quad |E|\le2-3\kappa/4.
\end{equation}
Fix $c_0>0$. There exist constants $\theta_*=\theta_*(c_0)>0$
and $c_{\rm p}^{*}>0$ such that for every deterministic $t_*$ satisfying
\begin{equation}\label{gc:eq:iteration-domain}
 \frac12\le t_*\le1-W^{-3},\qquad \zeta_{t_*}\le W^{-c_0},
\end{equation}
the conclusion is
\begin{equation}\label{gc:eq:iteration-output}
 \left.\eqref{gc:eq:one-step-input}\right|_{
 s\to u,\,\theta\to\theta_*,\,c_{\rm p}\to c_{\rm p}^{*}},
 \qquad u\in[1/2,t_*],
\end{equation}
in the same probability sense as in part (i).
\end{theorem}

The deterministic inputs are the propagator and Dyson bounds of
Lemmas~\ref{da:lem:scales-propagators} and~\ref{da:lem:positive-height},
and the spatial budgets of Lemma~\ref{da:lem:scale-budgets}.
Lemma~\ref{da:lem:comparison-parameters}(b) controls the finite comparison
chains in part (ii): it preserves bounded energies and a fixed positive
density margin at bounded coupling, and carries the initial margin
$|E|\le2-3\kappa/4$ into $|E'|\le2-\kappa/2$ at large coupling,
with at most a factor two in the upper-volume constant.

The proof in Appendix~\ref{app:gaussian-continuation} follows the
five-step propagation and finite induction of
\cite[Theorem~2.21, (2.68)--(2.72), Section~5]{YY} and
\cite[Theorem~5.4, (5.6)--(5.10), Section~7]{TYY}.
The main new work is the transfer-matrix and sampling inputs proved in
Section~\ref{sec:large-proof}; Appendix~\ref{sec:gc-one-step} checks
their compatibility with spatial sums, stopping and time suprema.

\subsection{Completion of the resolvent moment bounds}
\label{sec:fo-completion}
We complete Theorems~\ref{thm:bounded-lower} and~\ref{uvu:thm:terminal}.

\begin{lemma}[Deterministic bounds for the terminal comparison]
\label{fo:lem:terminal-deterministic}
\textup{(a)} Fix $C>0$ and $0<\kappa<1/2$. For $g\le C$ and all
sufficiently large $L$, the condition $F_g(-|E|)\ge\kappa$ implies
\begin{equation}\label{fo:eq:terminal-bounded-regularity}
 |E|\le C',\qquad \Im m_{g,L}(E)\ge c'>0.
\end{equation}
Also $|m_{g,L}(e)|\le\max\{1,g\}$ for every real $e$.
The constants and the volume threshold depend only on $C,\kappa$.

\textup{(b)} Fix $0<\kappa<2$. There is $C_\kappa>1$ such that, if
$g\ge C_\kappa$ and $L\ge C_\kappa g^2$, then for every
$|E|\le2-\kappa$ and $0\le\delta\le1$ there exists a deterministic
$E_0\in[-2+3\kappa/4,2-3\kappa/4]$ satisfying
\begin{equation}\label{fo:eq:terminal-large-anchor}
 E=E_0+\delta g^{-2}\Re m_{g,L}(E_0),\qquad
 c\le\Im m_{g,L}(E_0)\le |m_{g,L}(E_0)|\le C.
\end{equation}
The constants are independent of $g,L,E,\delta$.
\end{lemma}
Parts \textup{(a)} and \textup{(b)} are proved in
Appendices~\ref{sec:bl-kernels} and~\ref{app:terminal-deterministic-large}, respectively.

\begin{lemma}[Terminal matching at bounded coupling]
\label{fo:lem:bounded-terminal-matching}
Fix $C>0$ and $0<\kappa<1/2$. Let $g,L,E,\eta$ be deterministic
parameters, with $L$ an integer, satisfying
\[
 0<g\le C,\qquad \log W\ll L\le W^2,\qquad
 F_g(-|E|)\ge\kappa,\qquad 0<\eta=o(1).
\]
Set $m\coloneqq m_{g,L}(E+i\eta)$ and define
\begin{equation}\label{bl:eq:terminal}
 \begin{gathered}
 t_*\coloneqq\frac{\gamma_g\Im m}{\gamma_g\Im m+\eta},\qquad
 g_*\coloneqq\sqrt{t_*}\,g,\qquad
 c_*\coloneqq\sqrt{\frac{t_*\gamma_{g_*}}{\gamma_g}},\qquad
 E_*\coloneqq c_*\left(E-\eta\frac{\Re m}{\Im m}\right).
 \end{gathered}
\end{equation}
For sufficiently large $W$, these parameters satisfy
\begin{equation}\label{fo:eq:bounded-terminal-parameters}
 \tfrac12\le t_*<1,\qquad 1-t_*\asymp\eta,\qquad
 |E_*|\le C',\qquad \Im m_{g_*,L}(E_*)\ge c',
\end{equation}
and
\begin{equation}\label{fo:eq:bounded-terminal-resolvent}
 G(E+i\eta)\stackrel d=c_*G_{t_*}(g_*,E_*).
\end{equation}
Here $C',c'>0$ and the comparison constants depend only on $C,\kappa$.
The threshold in $W$ may depend on the rates in the asymptotic assumptions
and is uniform over families with uniform rates.
\end{lemma}
The proof is given under \emph{Terminal matching} in
Appendix~\ref{app:deterministic}.

\begin{proof}[Proof of Theorem~\ref{thm:bounded-lower}]
By \eqref{eq:scales}, $\ellstar\ge W$. Set
\begin{equation}\label{fo:eq:proof-parameters}
 \eta\coloneqq W^\nu/\ellstar\le W^{\nu-1}=o(1),\qquad
 c_0\coloneqq\nu/3.
\end{equation}
For each deterministic $E$ with $F_g(-|E|)\ge\kappa$,
define $t_*,g_*,c_*,E_*$ by \eqref{bl:eq:terminal}.
Lemma~\ref{fo:lem:bounded-terminal-matching} gives
\eqref{fo:eq:bounded-terminal-parameters} and
\eqref{fo:eq:bounded-terminal-resolvent}.
In particular, $0<g_*\le g\le C$ and $\gamma_{g_*}\asymp1$.
Together with the energy and density bounds, this gives the bounded-coupling
input range \eqref{gc:eq:initial-parameters}(a); it allows the comparison
coupling $g_*$ to lie below $W^{-1/2}$.
Here $\ellstar=\min\{WL,W^2g^2\}$, so
$\ellstar\le WL$ and $\ellstar\le W^2g^2\le C^2W^2$.
Using \eqref{fo:eq:scale} at coupling $g_*$ and $1-t_*\asymp\eta$, we obtain
\begin{equation}\label{bl:eq:endbudget}
 \begin{aligned}
 \zeta_{t_*}(g_*)&\le C'\left(\frac1{WL\eta}+\frac1{Wg\sqrt\eta}\right)\quad
  \le C'(W^{-\nu}+W^{-\nu/2})\le W^{-c_0},\\
 1-t_*&\asymp\eta\ge cW^{\nu-2}\gg W^{-3}.
 \end{aligned}
\end{equation}
Lemmas~\ref{fo:lem:deterministic} and~\ref{fo:lem:bounded-initial}
supply the primitive-loop and initial estimates, respectively, while
Lemmas~\ref{gc:lem:residual} and~\ref{fo:lem:sampling-continuity}
supply the residual and continuity inputs. These inputs are uniform
for fixed energy and density margins, including the comparison pairs
controlled by Lemma~\ref{da:lem:comparison-parameters}(b).
Theorem~\ref{gc:thm:common}(ii), with \eqref{bl:eq:endbudget}, therefore gives
\[
 \|\wt G_{t_*}(g_*,E_*)\|_{\max}
 \prec\sqrt{\zeta_{t_*}(g_*)}\le W^{-c_0/2}.
\]
Since $M_{xx}(g_*,E_*)=m_{g_*,L}(E_*)$ and
Lemma~\ref{fo:lem:terminal-deterministic}(a) gives
$|m_{g_*,L}(E_*)|\le\max\{1,g_*\}\le\max\{1,C\}$,
the decomposition $G_{t_*}=M+\wt G_{t_*}$ bounds the diagonal imaginary
parts of $G_{t_*}(g_*,E_*)$ by a fixed constant except on an event of
probability at most $W^{-A}$, for every fixed $A>0$.
Combining \eqref{fo:eq:bounded-terminal-resolvent} with $0<c_*\le1$,
which follows from \eqref{bl:eq:terminal} and
\eqref{fo:eq:bounded-terminal-parameters}, yields
\begin{equation}\label{fo:eq:bounded-terminal-bound}
 \max_x\Im G_{xx}(E+i\eta)\le C'
\end{equation}
with probability at least $1-W^{-A}$, uniformly over the stated
deterministic parameters.
Since $\Im G_{xx}\le\eta^{-1}\le C'W^2$, choosing $A>2p+1$
bounds the exceptional-event contribution to the $p$th moment by
$C_pW^{2p-A}$ and proves \eqref{eq:bounded-lower-moment}.
\end{proof}

\begin{proof}[Proof of Theorem~\ref{uvu:thm:terminal}]
Fix a deterministic $E\in[-2+\kappa,2-\kappa]$. Put $q=L/g^2$, so
$\log W\ll q\le W^2$, and set
$\delta_*=g^2\eta_L=W^\nu/[W\min(q,W)]\le W^{\nu-1}=o(1)$,
$t_*=1-\delta_*$ and $\widetilde g=\sqrt{t_*}\,g$.
For sufficiently large $W$, $t_*\ge1/2$ and $L\ge C_*g^2\log W$.
Lemma~\ref{fo:lem:matching}, with $g_*=\widetilde g$, therefore gives
$E_0\in[-2+3\kappa/4,2-3\kappa/4]$ such that, with
\[
 \eta(E):=\Im z_{t_*}(\widetilde g,E_0)
 =\eta_Lt_*^{-1}\Im m_{\widetilde g,L}(E_0)\asymp\eta_L,
\]
we have $z_{t_*}(\widetilde g,E_0)=E+i\eta(E)$ and
$G_{t_*}(\widetilde g,E_0)\stackrel d=G(E+i\eta(E))$
by \eqref{fo:eq:large-match-domain}.

We have $1-t_*=\delta_*\ge W^{\nu-2}\gg W^{-3}$, and
\eqref{fo:eq:scale} gives
\[
 \zeta_{t_*}(\widetilde g)
 \le C\left(\frac1{Wq\delta_*}+\frac1{W\sqrt{\delta_*}}\right)
 \le C(W^{-\nu}+W^{-\nu/2})\le W^{-\nu/3}.
\]
We apply Theorem~\ref{gc:thm:common}(ii) with bulk parameter $\kappa/2$
and $c_0=\nu/3$.
The comparison-parameter calculation \eqref{fo:eq:large-input-domain}
in its proof gives $g/2\le g'\le g$ and $|E'|\le2-\kappa/4$.
Consequently,
\[
 \frac{L}{(g')^2\log W}\ge\frac{L}{g^2\log W}\longrightarrow\infty,
 \qquad L\le4(g')^2W^2.
\]
Lemma~\ref{fo:lem:terminal-deterministic}(b), with $\delta=0$ and
bulk margin $\kappa/4$, bounds $|m_{g',L}(E')|$ above and
$\Im m_{g',L}(E')$ below uniformly.
After enlarging the coupling threshold, Lemmas~\ref{fo:lem:deterministic}(b),
\ref{fo:lem:resolvent-moments} and~\ref{fo:prop:inputs} supply the primitive,
resolvent-moment and initial estimates at every required pair; use
upper-volume constant $4$ in the initial lemma.
The displayed volume bounds and the verified regularity bounds also
give the hypotheses of Lemma~\ref{uvc:thm:ordinary-sampling} with a fixed
$C\ge3$ chosen also to cover the regularity bounds,
which supplies spatial sampling without any additional coupling threshold.
Lemmas~\ref{gc:lem:residual} and~\ref{fo:lem:sampling-continuity}
supply the remaining stochastic inputs.
All conditions of Theorem~\ref{gc:thm:common}(ii) are therefore satisfied,
and \eqref{gc:eq:iteration-output} gives
\begin{equation}\label{uvu:eq:terminal-law}
 \|\wt G_{t_*}(\widetilde g,E_0)\|_{\max}
 \prec\sqrt{\zeta_{t_*}(\widetilde g)}\le W^{-\nu/6}.
\end{equation}
Since $M_{ii}(\widetilde g,E_0)=m_{\widetilde g,L}(E_0)$ is uniformly
bounded by \eqref{fo:eq:terminal-large-anchor}, the matrix identity in law
above yields $\max_i\Im G_{ii}(E+i\eta(E))\le C$ outside an event
of probability at most $W^{-D}$, for every fixed $D>0$.

To obtain unconditional moments, take $\eta_0=g^{-2}$.
Since $\eta(E)\asymp\delta_*\eta_0$ and $\delta_*=o(1)$,
comparison of Poisson kernels gives
$\Im G_{ii}(E+i\eta(E))\le C\delta_*^{-1}\Im G_{ii}(E+i\eta_0)$.
Lemma~\ref{fo:lem:resolvent-moments}, with $t=1$, $C_0=1$
and moment order $2p$, bounds the seed moment uniformly.
Cauchy--Schwarz therefore bounds the exceptional-event contribution
to the $p$th moment by
$C\delta_*^{-p}W^{-D/2}\le CW^{(2-\nu)p-D/2}=o(1)$ when $D>4p$.
Finally, $\eta(E)\asymp\eta_L$ gives pointwise comparable Poisson
kernels, so the same moment bound holds at $\eta_L$, uniformly over
the stated $E,g,L$ and all coordinates, proving \eqref{uvu:eq:local-moments}.
\end{proof}

\section{Dyson solutions, pair kernels, and propagation scales}
\label{sec:lower-deterministic}\label{sec:scalar-dyson-regularity}\label{sec:deterministic-inputs}
We state the deterministic inputs used by the main-text proofs, in
proof-dependency order. Complete proofs and appendix-only refinements
are in Appendix~\ref{app:dyson-foundations}.
Use the characteristic, Dyson matrices, pair kernels and scales defined in
\eqref{fo:eq:characteristic}, \eqref{gc:eq:positive-dyson},
\eqref{fo:eq:PQ} and \eqref{fo:eq:scale}.

\begin{lemma}[Scales and deterministic propagators]\label{da:lem:scales-propagators}
Use $\gamma=\gamma_g$ and the flow, scales, pair kernels and propagators
in \eqref{fo:eq:characteristic}, \eqref{fo:eq:scale}, \eqref{fo:eq:PQ}
and \eqref{fo:eq:propagators}. Assume either
\begin{equation}\label{da:eq:propagation-domain}
 \begin{aligned}
 &g\le C,\quad |E|\le B_0,\quad \Im m_{g,L}(E)\ge c',\quad
       \log W\ll L\le W^{B_L};\\
 &g\ge C_\kappa,\quad |E|\le2-\kappa,\quad L\gg g^2\log W.
 \end{aligned}
\end{equation}
The constants $C,B_0,c',B_L,\kappa$ are fixed.
Uniformly for $1/2\le s\le u\le t<1$, $1-t\ge W^{-3}$, we have
\begin{equation}\label{fo:eq:propagation-dyson}
 \begin{gathered}
 c\le\Im m\le C',\qquad \max_{a,b}|M_{B,ab}|\le C',\\
 \gamma\max_a\sum_b\bigl(|M_{B,ab}|+|M_{B,ab}|^2\bigr)\le C',
 \qquad \max_a\sum_b|[(I-uQ)^{-1}]_{ab}|\le C'.
 \end{gathered}
\end{equation}

For every prescribed $A>0$, the following bounds hold with $c,C'$
independent of $A$; $C_A$ and the width threshold may depend on $A$:
\begin{equation}\label{fo:eq:kernel-bounds}
 0\le[P(I-uP)^{-1}]_{ab}
 \le\frac {C'}{(1-u)\ell_u}e^{-c\|a-b\|_L/\ell_u}+C_A\ell_u^{-1}W^{-A}.
\end{equation}

Set $\theta_u=(1-u)P(I-uP)^{-1}$. Since $P\mathbf1=\mathbf1$,
this is a probability profile, with $\|\theta_u\|_\infty\le C/\ell_u$
by \eqref{fo:eq:kernel-bounds}.

For a scalar function $h$ on the cycle and $n\ge2$, define the
row-restoration commutator by
\begin{equation}\label{gc:eq:commutator}
 \begin{split}
 C_v(h)(\boldsymbol a)\coloneqq {}&\sum_c[P(I-vP)^{-1}]_{a_1c}h(c) \cdot\left[\prod_{j=2}^n\theta_v(c-a_j)
             -\prod_{j=2}^n\theta_v(a_1-a_j)\right].
 \end{split}
\end{equation}
Then, with $\operatorname{row}A(a_1)=\sum_{a_2,\ldots,a_n}A(\boldsymbol a)$,
\begin{equation}\label{da:eq:commutator-bound}
 \operatorname{row}C_v(h)=0,\qquad
 \|C_v(h)\|_\infty\le\frac{C_n\|h\|_\infty}{(1-v)\ell_v^{n-1}}.
\end{equation}

Constants are independent of $g,L$; the width threshold may depend on
the rate of divergence in the volume lower bound.
\end{lemma}
The proof is in Appendix~\ref{app:main-propagators}.

\begin{lemma}[Spatial summation and weights]\label{da:lem:scale-budgets}
Use the parameter ranges of Theorem~\ref{gc:thm:common}, with fixed
regularity margins, and $0\le u\le1-W^{-3}$.
For a fixed $C$, the scales in \eqref{fo:eq:scale} satisfy
\begin{equation}\label{fo:eq:spatial-sum-budget}
 \gamma_g L+L/\ell_u\le W^C,\qquad
 \gamma_gW\ell_u\zeta_u=(1-u)^{-1}.
\end{equation}
For fixed $c_1,F>0$ and $0<\theta\le1/2$, define
$w_u(r)=e^{-c_1(r/\ell_u)^\theta}+W^{-F}$. Then
\begin{equation}\label{da:eq:weight-comparison}
 e^{-C(n/\ell_u)^\theta}\le w_u(r')/w_u(r)
       \le e^{C(n/\ell_u)^\theta}\quad(|r-r'|\le2n,\ n\ge0).
\end{equation}
In particular the weight ratio is at most two on cells of length
$n\le g^2W^{-B}$ in the large-coupling range. If
$LW^{-F/2}\le\ell_u$, then
\begin{equation}\label{da:eq:weight-convolution}
 \sum_b w_u(d(a,b))w_u(d(b,c))\le C_\theta\ell_u w_u(d(a,c)),
 \qquad
 \sum_b\sqrt{w_u(d(a,b))w_u(d(b,c))}
       \le C_\theta\ell_u\sqrt{w_u(d(a,c))}.
\end{equation}
All constants and exponents are independent of $g,L$ in the stated ranges.
\end{lemma}
The proof is in Appendix~\ref{app:scale-budgets}.

\begin{lemma}[Comparison parameters through finite iteration]\label{da:lem:comparison-parameters}
\textup{(a)} Use the parameter hypotheses of Lemma~\ref{fo:lem:families}, and fix
$1/2\le s\le u<1$. Set $g'=g\sqrt{s/u}$ and
$\xi=(s\gamma_{g'}/(u\gamma_g))^{1/2}$. There are deterministic $E',w$,
with the energy, density and scale bounds of Lemma~\ref{fo:lem:families}, such that
\begin{equation}\label{da:eq:comparison-parameters}
 \Re w=\Re z_u(g,E),\qquad
 1\le\frac{\Im w}{\eta_u(g,E)}\le C\frac{1-s}{1-u},\qquad
 w=\xi^{-1}z_s(g',E').
\end{equation}
At large coupling, $\xi=1$. The corresponding density and energy bounds
hold along the vertical segment from $z_u(g,E)$ to $w$.

\textup{(b)} Start from $(g_0,E_0)$ in the initial parameter range of
Theorem~\ref{gc:thm:common}(ii). Consider a chain $(g_j,E_j)_{j=0}^m$ of $m\le J$ comparisons from
part \textup{(a)}, where $J$ is fixed, interspersed with movement to
earlier times at fixed coupling, all times in $[1/2,1)$.
At bounded coupling, all resulting energies are bounded and all
$\Im m_{g,L}(E)$ have a positive lower bound depending only on the
initial bounds and $J$; also $c_Jg_0\le g\le C$.
At large coupling the domain
\begin{equation}\label{fo:eq:large-time-domain}
 u\gamma_g\le g_0^{-2},\qquad s\gamma_{g'}=u\gamma_g
\end{equation}
is preserved. For sufficiently large $C_\kappa$, the chain satisfies
\begin{equation}\label{fo:eq:large-input-domain}
 g_0/\sqrt2\le g\le g_0,\qquad |E|\le2-\kappa/2.
\end{equation}
The volume bounds remain valid, with at most a factor two in the
large-coupling upper-volume constant. Earlier comparisons preserve both
$1-u\ge W^{-3}$ and an upper bound $\zeta_u(g)\le W^{-c_0}$.
\end{lemma}
Parts \textup{(a)} and \textup{(b)} are proved in Appendices~\ref{app:matching-parameters} and \ref{app:deterministic-iteration}, respectively.

\begin{lemma}[Dyson bounds on comparison segments]\label{da:lem:positive-height}
The scalar equation in \eqref{gc:eq:positive-dyson} has a unique
Herglotz solution; at the characteristic it has the values given there.
In either parameter range of Lemma~\ref{fo:lem:families}, fix
$1/2\le s\le u<1$. Let the deterministic comparison point be
$w=\xi^{-1}z_s(g',E')$, where $g'=g\sqrt{s/u}$,
$\xi=(s\gamma_{g'}/(u\gamma_g))^{1/2}$ and $E'$ is the comparison
energy in Lemma~\ref{fo:lem:families}. Let
$z=\Re z_u(g,E)+iy$, $\eta_u(g,E)\le y\le\Im w$.
Use the Dyson solution and matrix in \eqref{gc:eq:positive-dyson},
and put $\gamma=\gamma_g$.
Uniformly on this segment,
\begin{equation}\label{gc:eq:positive-stability}
 \begin{gathered}
 c\le\Im m_u(z)\le C',\qquad
 \max_{a,b}|M_{B,u}(z)_{ab}|\le C',\\
 \gamma\max_a\sum_b
       \bigl(|M_{B,u}(z)_{ab}|+|M_{B,u}(z)_{ab}|^2\bigr)\le C',\\
 \max_a\sum_b
 \left|\left[(I-u\gamma(M_{B,u}(z)_{ab}^{\,2})_{a,b})^{-1}
       \right]_{ab}\right|\le C'.
 \end{gathered}
\end{equation}
At large coupling, these bounds also hold directly for
$|x|\le2-\kappa/2$, $z=x+iy$, $0<y\le C/g^2$ and $1/2\le u\le1$,
provided $g\ge C_\kappa$ and $L\gg g^2\log W$.
At the characteristic for every $0\le u<1$, and on
the comparison segments above, one also has
\begin{equation}\label{gc:eq:residual-row-bounds}
 \max_{a,c}|M_{B,ac}|\le C,\qquad
 u\gamma_g\max_a\sum_c\bigl(|M_{B,ac}|+|M_{B,ac}|^2\bigr)\le C.
\end{equation}
Here $M_B=M_{B,u}(z)$. At the upper comparison point,
$M_{B,u}(w)=\xi M_B(g',E')$. For $r=\Im z/\gamma_g$, with
$u,g,\Re z$ held fixed, one has
\begin{equation}\label{da:eq:height-derivative}
 \max_{a,b}|\partial_rM_{B,u}(z)_{ab}|
 +\gamma_g\max_a\sum_b|\partial_rM_{B,u}(z)_{ab}|\le C.
\end{equation}
Constants are independent of $g,L$.
\end{lemma}
The proof is in Appendix~\ref{sec:bl-matching}.

\begin{lemma}[Free rows and interval interpolation]\label{da:lem:free-interpolation}
Fix $|E|\le2-\kappa$, $g\ge C_\kappa$, $L\ge4g^2$, and
$m=m_{g,L}(E)$ with $c\le\Im m\le|m|\le C$.
For $M_B=(A_L-E-g^{-2}m)^{-1}$,
\begin{equation}\label{uvc:eq:free-row-bounds}
 \max_{b,c}|M_{B,bc}|\le C,\qquad
 g^{-2}\max_b\sum_c|M_{B,bc}|\le C.
\end{equation}
Let $I$ be a deterministic consecutive interval with $|I|\le cg^2$
and first block $b_0$; block indices are taken modulo $L$.
Write $E=2\sin\theta$, $-\pi/2<\theta<\pi/2$, and lift $I$ to an
integer interval starting at $b_0$. For $a\in I$ and any block $b$, define
\begin{equation}\label{da:eq:cell-interpolation-coefficients}
 \begin{aligned}
 c_{a,b_0}&=\frac{\sin((a-b_0+1)(\pi/2-\theta))}{\cos\theta},\\
 c_{a,b_0-1}&=-\frac{\sin((a-b_0)(\pi/2-\theta))}{\cos\theta},\\
 q_{ab}&=\one_{\{b_0\le b<a\}}
          \frac{\sin((a-b)(\pi/2-\theta))}{\cos\theta}.
 \end{aligned}
\end{equation}
Set $c_{a\alpha}=0$ for $\alpha\notin\{b_0,b_0-1\}$.
In the indicator, lift $b$ to $[b_0,b_0+L)$. The coefficients $c,q$
are real and bounded by $C$, independent of $g$, time and the random matrices
once $E$ and $I$ are fixed. Then
\begin{equation}\label{da:eq:dyson-interpolation}
 M_{B,ab}=\sum_\alpha c_{a\alpha}M_{B,\alpha b}+q_{ab}+e^M_{ab},\qquad
 |e^M_{ab}|\le C|I|/g^2.
\end{equation}
These are the coefficients of the real free forced recurrence, so the
same $c,q$ apply when that recurrence is perturbed by onsite matrices.
The same interpolation bound holds for the large-coupling
comparison matrices of Lemma~\ref{da:lem:positive-height}.
Interpolation in both coordinates follows by applying
\eqref{da:eq:dyson-interpolation} separately in each coordinate on its
respective interval, using symmetry.
\end{lemma}
The proof is in Appendix~\ref{app:free-interpolation}.

For the upper-bound argument in Section~\ref{sec:small-upper-proof},
we also need the following one-block Dyson profile.

\begin{lemma}[Dyson profile for a regular inverse]\label{da:lem:regular-dyson-profile}
Fix $\kappa_0\in(0,2)$ and $C_*>0$. Let $g=o(1)$ and let $R$ be
Hermitian with $\rank\one_{\{|g^2R|>g\}}\le Wg^{1/2}$.
For $|E|\le2-\kappa_0$, define the Herglotz solution
$M(i\eta)=(-EI-g^2R-i\eta-m(i\eta))^{-1}$ and
$m(i\eta)=W^{-1}\Tr M(i\eta)$. Then, for sufficiently large $W$,
\begin{equation}\label{us:eq:mde}
 \Im m(i\eta)\asymp1,\qquad
 \Im M(i\eta)\asymp(I+g^4R^2)^{-1},\qquad 0<\eta\le C_*.
\end{equation}
The matrix comparison is in the sense of quadratic forms, with constants
depending only on $\kappa_0,C_*$.
\end{lemma}
The proof is in Appendix~\ref{app:deformed-dyson}.

\section{Transfer estimates, spatial sampling, and loop initialization}
\label{sec:large-proof}\label{sec:input-proofs}
We prove the large-coupling inputs stated in
Section~\ref{sec:fo-transfer-hierarchy}, which allow the common
Theorem~\ref{gc:thm:common} to be used without an upper bound on $g$. The transfer estimate of Lemma~\ref{uvs:lem:transfer-moments} gives the
resolvent seed through Ward's identity and also supplies the coarse moments
used during propagation. The same estimates prove the short-cycle spreading
bound of Lemma~\ref{uv:lem:short-spreading} in
Section~\ref{sec:uniform-resolvent-seed}. The weighted Ward estimate in
Section~\ref{sec:uniform-resolvent-seed} adds exponential spatial decay.
Section~\ref{app:uniform-compression} then proves the spatial sampling
of Lemma~\ref{uvc:thm:ordinary-sampling}:
short-cell transfer control gives Green-function interpolation; the auxiliary
orbital limit transfers this interpolation to primitive loops; and a common
finite-dimensional space reduces the required spatial maxima to polynomially
many ordinary tuples. Section~\ref{sec:sampling-propagation-extensions}
establishes the continuity, weighted sampling and residual interpolation
needed for propagation. Sections~\ref{sec:input-averaging}
and~\ref{sec:uniform-initialization} prove fluctuation averaging and the
complete initial estimates, respectively, uniformly in $g$.
The five-step propagation is proved in
Appendix~\ref{sec:gc-one-step}.

\begin{definition}[Transfer matrices and free coordinates]
\label{def:transfer-matrices}
For $g>0$, $z\in\C$, and Hermitian onsite blocks $V_j$, define
\[
 T_j(z):=\begin{pmatrix}zI_W-V_j/g&-I_W\\I_W&0\end{pmatrix},
 \qquad T(z):=\begin{pmatrix}z&-1\\1&0\end{pmatrix}\otimes I_W.
\]
For time-dependent Hermitian onsite blocks $V_j(t)$, define
\[
 T_j(z,t):=\begin{pmatrix}zI_W-V_j(t)/g&-I_W\\I_W&0\end{pmatrix}.
\]
The free transfer matrix $T(z)$ corresponds to zero onsite blocks;
for fixed $z$, it is independent of $j$ and $t$.
For $E\in(-2,2)$, write $E=2\sin\theta$ with $-\pi/2<\theta<\pi/2$ and set
\[
 \begin{aligned}
 S_E&:=(2\cos\theta)^{-1/2}
       \begin{pmatrix}1&1\\ie^{-i\theta}&-ie^{i\theta}\end{pmatrix}
       \otimes I_W,\\
 R_E&:=S_E^{-1}T(E)S_E
       =\diag(-ie^{i\theta},ie^{-i\theta})\otimes I_W.
 \end{aligned}
\]
The matrix $R_E$ is unitary, and
$\|S_E\|+\|S_E^{-1}\|\le C_\kappa$ for $|E|\le2-\kappa$.
Write $S:=S_E|_{E=0}$, which is also unitary.
\end{definition}

\subsection{Transfer moments on the \texorpdfstring{$g^2$}{g squared} scale}
\label{sec:uniform-seed}
The products in Lemma~\ref{uvs:lem:transfer-moments} are controlled at each
deterministic energy, with constants independent of the coupling and the
number of blocks in a window of length at most $C_0g^2$.
\begin{lemma}[Maximal transfer moments on scattering-length windows]
\label{uvs:lem:transfer-moments}
Fix $\kappa>0$ and $C_0,C_1<\infty$.
There is $C_\kappa<\infty$ such that, for $g\ge C_\kappa$,
$|E|\le2-\kappa$, $n\le C_0g^2$, and
$0\le\eta\le1$ with $n\eta\le C_1$, the matrices $T_j(z)$ from
Definition~\ref{def:transfer-matrices}, with $z=E+i\eta$, satisfy,
for every fixed $p>0$ and all sufficiently large $W$,
\begin{equation}\label{uvs:eq:transfer-moments}
 \E\max_{0\le j\le n}
 \left(\norm{T_j(z)\cdots T_1(z)}^p+
       \norm{(T_j(z)\cdots T_1(z))^{-1}}^p\right)
 \le C.
\end{equation}
Products at $j=0$ are identities. The constants and width threshold
are independent of $g,n,E,\eta$ in the displayed range. The same
bound holds when the onsite variance is any number in $[0,1]$.
\end{lemma}
\begin{proof}
We use the real-energy normal form of Sadel--Schulz-Baldes
\cite[Section~5.1, (22)]{SSB}. In the coordinates $S_E$ of
Definition~\ref{def:transfer-matrices}, substituting $z=E+i\eta$ gives
\[
 M_j=S_E^{-1}T_j(z)S_E=R_E+i\eta B_E+g^{-1}X_j,
\]
where $\norm{B_E}\le C$ and $X_j$ is a
centered real-linear function of the GUE block.
For every deterministic unit vector $v$, the real Gaussian vector
$\xi=X_jv\in\R^{4W}$ has covariance $\Sigma_v$ satisfying
\[
 \Tr\Sigma_v\le C,\qquad
 \norm{\Sigma_v}\le C/W.
\]
Indeed, for normalized GUE,
$\E\norm{Vw}^2=\norm{w}^2$ and the variance of every real linear
functional of $Vw$ is at most $C\norm{w}^2/W$.
The bounded deterministic conjugations preserve these estimates.

We record the needed Gaussian calculation. Put
$w=R_Ev+i\eta B_Ev$, so
$\norm{w}^2\le1+C\eta$. For $q\ge2$,
$x^{q/2}\le e^{q(x-1)/2}$ gives
\[
 \E\norm{M_jv}^q
 \le e^{Cq\eta}
 \E\exp\!\left\{\frac qg\Re\langle w,\xi\rangle
                   +\frac{q}{2g^2}\norm{\xi}^2\right\}.
\]
The Gaussian quadratic-moment formula, diagonalizing $\Sigma_v$,
bounds the Gaussian expectation by
\[
 \exp\!\left\{Cq/g^2+Cq^2/(Wg^2)\right\},
 \qquad q/(Wg^2)\le c.
\]
To check the constants directly, the logarithmic determinant term
is at most $C(q/g^2)\Tr\Sigma_v$, and the linear term is at most
$C(q^2/g^2)\norm{\Sigma_v}\norm{w}^2$.
Choose $q=W$ and increase $C_\kappa$ if necessary. Homogeneity
then gives, for every past-measurable vector $v$,
\begin{equation}\label{uvs:eq:one-step-moment}
 \E[\norm{M_jv}^q\mid V_1,\ldots,V_{j-1}]
 \le e^{a}\norm{v}^q,\qquad
 a=Cq(g^{-2}+\eta).
\end{equation}

For deterministic unit $v$,
$e^{-aj}\norm{M_j\cdots M_1v}^q$ is a nonnegative
supermartingale. The maximal inequality therefore yields
\[
 \Prob\{\max_{j\le n}\norm{M_j\cdots M_1v}>u\}
 \le e^{an}u^{-q}.
\]
A $1/2$-net of the unit sphere of $\C^{2W}$ has at most
$5^{4W}$ elements, and the operator norm is bounded by twice
the maximum over that net. Since $an\le Cq(C_0+C_1)$ and $q=W$,
\begin{equation}\label{uvs:eq:transfer-tail}
 \Prob\{\max_{j\le n}\norm{M_j\cdots M_1}>u\}
 \le\min\{1,(C/u)^W\}.
\end{equation}
For $W\ge2p$, integrating this tail gives a bound $C$
independent of $g,n$. The conjugation by $S_E$ changes only constants.

For the inverse prefixes, apply the identical argument to the
fresh left-multiplying factors $M_j^{-*}$:
\[
 (M_j\cdots M_1)^{-*}=M_j^{-*}\cdots M_1^{-*}.
\]
They have the same structure as above. In fact
\[
 T_j(z)^{-1}=
 \begin{pmatrix}0&I_W\\-I_W&zI_W-V_j/g\end{pmatrix}
\]
is affine in the same fresh Gaussian block, and the free
inverse-adjoint conjugated matrix is $R_E^{-*}=R_E$.
Thus the same covariance and deterministic-perturbation bounds
hold. This proves \eqref{uvs:eq:transfer-moments}.
\end{proof}

\subsection{Resolvent moments and exponential spatial decay}
\label{sec:uniform-resolvent-seed}
\begin{proof}[\textbf{Proof of Lemma~\ref{fo:lem:resolvent-moments}: \eqref{fo:eq:coarse-resolvent-moments}}]
We first prove a pathwise bound and coarse moments for general $L$.
Fix $g\ge C_\kappa$, $|E|\le2-\kappa$, and $t\in[0,1]$, and set
\begin{equation}\label{uvs:eq:window-length}
 n=\min\!\left\{\left\lfloor\frac{L-1}{2}\right\rfloor,
                    \lfloor g^2\rfloor\right\}.
\end{equation}
Write $G=(H_t(g)-E-i\eta)^{-1}$ and let $f_x=G_{xb}u$ for a
deterministic vector $u$. At every block except $b$, $f$ obeys the
homogeneous transfer recurrence. For each $a$, select a deterministic
direction so that the arc from $a$ through $a+n$ has no $b$ in its
interior; such a direction exists because $2n<L$.
Let $K_a$ be the maximum inverse-transfer norm along this arc,
including the identity. The transfer recurrence and summation over
$0\le j<n$ give
\[
 \norm{f_a}^2
 \le K_a^2(\norm{f_{a+j}}^2+\norm{f_{a+j+1}}^2),\qquad
 \norm{f_a}^2\le\frac{2K_a^2}{n}\norm f^2,
\]
since each block occurs at most twice in the sum. This applies also at
$a=b$, because the source is an endpoint and its inhomogeneous equation
is not used. Ward's identity gives
$\eta\norm f^2=\Im\langle u,f_b\rangle\le\norm u\norm{f_b}$.
Combining these inequalities first at $b$ and then at $a$ yields
\begin{equation}\label{uvs:eq:pathwise-green-block}
 \norm{\bigl[(H_t(g)-E-i\eta)^{-1}\bigr]_{ab}}_{\op}
 \le\frac{2K_aK_b}{n\eta}.
\end{equation}
For $0<\eta\le\min\{1,C_0/n\}$, we have $n\le g^2$ and
$n\eta\le C_0$, so Lemma~\ref{uvs:lem:transfer-moments} bounds
$\E K_a^{2p}$ and $\E K_b^{2p}$ uniformly in the onsite variance
$t\in[0,1]$. Cauchy--Schwarz, without any independence requirement, gives
\begin{equation}\label{uvs:eq:coarse-moments}
 \sup_{0\le t\le1}\max_{a,b}
 \E\norm{\bigl[(H_t(g)-E-i\eta)^{-1}\bigr]_{ab}}_{\op}^p
 \le C_p(n\eta)^{-p},\qquad 0<\eta\le\min\{1,C_0/n\}.
\end{equation}
At the seed height $\eta=c_0/n$, with fixed
$0<c_0\le\min\{1,C_0\}$, this gives the bounded resolvent moments
needed for initialization, without an upper bound on $g$ or $L$.

Finally, if $L\ge4g^2$, then $n=\lfloor g^2\rfloor\asymp g^2$;
thus $0<\eta\le\min\{1,C_0/g^2\}$ implies $n\eta\le C_0$, and
\eqref{uvs:eq:coarse-moments} gives
\eqref{fo:eq:coarse-resolvent-moments}.
If $1<\eta\le C_0/g^2$, then $g^2<C_0$ and the deterministic bound
$\|(H_t(g)-E-i\eta)^{-1}\|\le\eta^{-1}\le C_0(g^2\eta)^{-1}$
proves the remaining case.
\end{proof}

\begin{proof}[Proof of Lemma~\ref{uv:lem:short-spreading}]
For $L<L_{\kappa/2}$, Lemma~\ref{thm:orbital} proves the claim after
reducing $c_\kappa$, since $\min\{L,g^2\}\le L_{\kappa/2}$.
Assume henceforth $L\ge L_{\kappa/2}$ and write $J=J_{g,\kappa,L}$.
By Lemma~\ref{da:lem:distribution}\textup{(iii)} with margin $\kappa/2$, choosing $C_\kappa$
large enough puts $J$ in a fixed free-bulk interval. Take $n$ from
\eqref{uvs:eq:window-length} and $\eta=c_0/n$, decreasing $c_0$ so that
$J_\eta=J+[-\eta,\eta]$ remains in the free bulk. All constants below
depend only on $\kappa$.

First we prove the seed tail estimates \eqref{uv:eq:seed-exponential-tail}.
At $t=1$, the inverse-product version of \eqref{uvs:eq:transfer-tail} and
\eqref{uvs:eq:pathwise-green-block}, with $a=b$, give
$\Prob(K_a>v)\le\min\{1,(C/v)^W\}$ and
$X_a(E)\coloneqq\|G_{aa}(E+i\eta)\|\le2K_a^2/c_0$.
Consequently, uniformly for $E\in J_\eta$ and every block $a$,
\begin{equation}\label{uv:eq:seed-exponential-tail}
 \begin{gathered}
 \Prob(X_a(E)>v)\le\min\{1,(C/v)^{W/2}\},\qquad v>0,\\
 \E[X_a(E)\one_{\{X_a(E)>A\}}]\le C A(C/A)^{W/2},
 \qquad A>C,\quad W\ge4.
 \end{gathered}
\end{equation}
The second inequality follows from
$\E[X\one_{\{X>A\}}]=A\Prob(X>A)+\int_A^\infty\Prob(X>v)\,dv$.

Next we obtain the heavy-block mass estimate
\eqref{uv:eq:constant-heavy-block}.
Since $n\ge\frac14\min\{L,g^2\}$,
$\frac{C'_\kappa}{2\eta\min\{L,g^2\}}\ge\frac{C'_\kappa}{8c_0}$;
choose $C'_\kappa$ so that the last constant exceeds $2C$ in
\eqref{uv:eq:seed-exponential-tail}.
Apply \eqref{uv:eq:matrix-poisson} to each block spectral measure,
using $\|\Im G_{aa}\|\le X_a(E)$ and $\Tr\Im G_{aa}\le WX_a(E)$.
Summing over the blocks gives
\begin{equation}\label{uv:eq:constant-heavy-block}
 \begin{aligned}
 &\frac1{WL}\E\sum_{k:\lambda_k\in J}\sum_a
 \|\psi_{k,a}\|^2\one_{\{\|\psi_{k,a}\|^2>C'_\kappa/\min\{L,g^2\}\}}\\
 &\quad\le\frac C L\sum_a\int_{J_\eta}
 \E[X_a(E)\one_{\{X_a(E)>C'_\kappa/(8c_0)\}}]\,dE
 \le Ce^{-cW},
 \end{aligned}
\end{equation}
where the last inequality uses \eqref{uv:eq:seed-exponential-tail}
and the bounded length of $J_\eta$.

Finally we deduce \eqref{uv:eq:constant-scattering-outline}.
By \eqref{uv:eq:constant-heavy-block}, the expected number of
eigenvectors with $\lambda_k\in J$ whose total mass on these heavy blocks
exceeds $1/8$ is at most $8CWL e^{-cW}$. For every remaining eigenvector,
a cyclic coordinate interval $I$ meets at most $|I|/W+2$ blocks, so
\[
 \sum_{x\in I}|\psi_k(x)|^2
 \le\frac18+\left(\frac{|I|}{W}+2\right)
                   \frac{C'_\kappa}{\min\{L,g^2\}}<\frac12
\]
whenever $|I|\le c_\kappa W\min\{L,g^2\}$, $c_\kappa$ is sufficiently
small, and $\min\{L,g^2\}$ is sufficiently large.
For bounded $\min\{L,g^2\}$, reduce $c_\kappa$ so that the target
length is at most the single-block scale in \eqref{eq:haar-count};
that estimate covers the remaining case. Absorbing its factor $W$
into the exponential proves \eqref{uv:eq:constant-scattering-outline},
uniformly in $L,g$.
\end{proof}

\begin{lemma}[Weighted Ward estimate]\label{uvc:lem:weighted-ward}
For $g\ge1$ and $0\le t\le1$, let $G(z)=(H_t(g)-z)^{-1}$, and let
$\iota_b:\C^W\to\C^{WL}$ insert a vector into block $b$.
For $z=E+i\eta$, $0<\eta\le1$, and $u\in\C^W$, put
\[
 f=G(z)\iota_bu,\qquad f_c:=G_{cb}(z)u\in\C^W,\qquad
 w_c=\exp(\alpha d(c,b)),\qquad c\in\Z_L,
\]
where $0\le\alpha\le\eta/8$.
Then
\[
 \eta\sum_c w_c^2\|f_c\|^2
       \le 2\Im\langle u,G_{bb}(z)u\rangle.
\]
The estimate holds pathwise, with the same constant in the auxiliary
orbital model.
\end{lemma}
\begin{proof}
We use a Combes--Thomas-type exponential-weight argument; see
\cite[Appendix~A.2, p.~827]{GKT} for the classical exponential-conjugation
method. The calculation below retains the local imaginary part on the
right-hand side.

Fix a realization and use $\langle x,y\rangle=x^*y$.
Since $g\ge1$, the normalization in \eqref{fo:eq:flow} gives
$\sqrt{\gamma_g}\,g=1$, so the hopping blocks are $I_W$.
Writing $B_c=B_c^*$ for the onsite blocks, the resolvent equation is
\[
 (B_c-E)f_c+f_{c-1}+f_{c+1}-i\eta f_c
       =\one_{\{c=b\}}u,
\]
with all block indices taken modulo $L$. Set
$S:=\sum_cw_c^2\|f_c\|^2$. Multiply the equation by $w_c^2f_c^*$
from the left and sum over $c$. The onsite contribution
$\sum_cw_c^2\langle f_c,(B_c-E)f_c\rangle$ is real.
Reindexing the term involving $f_{c-1}$ writes the hopping contribution as
\[
 \sum_c\left[
 w_c^2\langle f_c,f_{c+1}\rangle
 +w_{c+1}^2\langle f_{c+1},f_c\rangle\right].
\]
Since $\langle f_{c+1},f_c\rangle
=\overline{\langle f_c,f_{c+1}\rangle}$, its imaginary part is
$\sum_c(w_c^2-w_{c+1}^2)\Im\langle f_c,f_{c+1}\rangle$.
On the right, $w_b=1$ and $f_b=G_{bb}(z)u$ give
\[
 \Im\langle f_b,u\rangle
       =-\Im\langle u,G_{bb}(z)u\rangle.
\]
Taking imaginary parts and rearranging therefore yields the exact identity
\[
 \eta S=\Im\langle u,G_{bb}(z)u\rangle+\mathcal E,\qquad
 \mathcal E:=\sum_c(w_c^2-w_{c+1}^2)
                    \Im\langle f_c,f_{c+1}\rangle.
\]

The distance satisfies $|d(c+1,b)-d(c,b)|\le1$. Hence
 \[ \begin{aligned}
 |\mathcal E|
 &\le2\sinh\alpha\sum_cw_cw_{c+1}\|f_c\|\|f_{c+1}\|\,\le\sinh\alpha\sum_c
       \bigl(w_c^2\|f_c\|^2+w_{c+1}^2\|f_{c+1}\|^2\bigr)
 =2\sinh\alpha\,S.
 \end{aligned}
\]
Since $\alpha\le\eta/8\le1$, we have
$2\sinh\alpha\le4\alpha\le\eta/2$. Consequently
\[
 \eta S\le\Im\langle u,G_{bb}(z)u\rangle+\frac\eta2S,
 \qquad
 \eta\sum_cw_c^2\|f_c\|^2
       \le2\Im\langle u,G_{bb}(z)u\rangle.
\]
Only Hermiticity of the onsite blocks and the normalized nearest-neighbor
hopping were used. The calculation is pathwise and independent of the
orbital dimension, so the same constant applies to the auxiliary model.
\end{proof}

\paragraph{Exponential spatial decay.}
For $L\ge4g^2$, $g\ge C_\kappa$,
$|E|\le2-\kappa$ and fixed $0<c_0\le C_0$, the normalized model satisfies
\begin{equation}\label{uvc:eq:seed-tails}
 \E\|G_{ab}(E+i\eta)\|_{\op}^p
 \le C_p e^{-pc\,d(a,b)/g^2},\qquad
 c_0/g^2\le\eta\le C_0/g^2,
\end{equation}
for every fixed $p>0$, with $G(z)=(H-z)^{-1}$. The same bound holds
with $G(z)=(H_t(g)-z)^{-1}$ at each fixed $0\le t\le1/2$.
Constants are uniform in $g,L,t$, and $c$ is independent of $p$;
no polynomial volume upper bound is needed here.
\begin{proof}[\textbf{Proof of \eqref{uvc:eq:seed-tails}}]
For $\eta>1$, the weighted identity in the proof of
Lemma~\ref{uvc:lem:weighted-ward} with fixed
small $\alpha>0$ and $\|G_{bb}\|\le\eta^{-1}$ gives
$\|G_{ab}\|\le\sqrt2\,\eta^{-1}e^{-\alpha d(a,b)}$, which suffices since $g\ge1$.
Assume henceforth $\eta\le1$. For $u\in\C^W$, put $f=G(E+i\eta)\iota_bu$.
Choose from $a$ a source-avoiding arc of $n\asymp g^2$ sites,
as in the derivation of \eqref{uvs:eq:pathwise-green-block}, and let $K_a$ be the maximal
inverse-transfer norm there.  The homogeneous recurrence implies
\[
 \|f_a\|^2\le\frac{2K_a^2}{n}\sum_{c\text{ on arc}}\|f_c\|^2.
\]
On this arc $d(c,b)\ge d(a,b)-n-1$.  Apply Lemma~\ref{uvc:lem:weighted-ward}
with $\alpha=c\eta$, observing $\alpha n=O(1)$:
\[
 \|G_{ab}u\|^2
 \le \frac{C K_a^2}{n\eta}
       e^{-2\alpha d(a,b)}\|G_{bb}\|\,\|u\|^2.
\]
Here $n\eta\asymp1$. Lemmas~\ref{uvs:lem:transfer-moments}
and~\ref{fo:lem:resolvent-moments} bound every fixed moment of
$K_a$ and $G_{bb}$, since $g^2\eta\ge c_0$;
no independence between them is required.
Taking the supremum over unit vectors $u$ and applying H\"older proves
\eqref{uvc:eq:seed-tails}, since $\alpha\ge c/g^2$.
The Brownian-time version uses the same one-step covariance bounds
with variance no larger than one.
\end{proof}

\subsection{Spatial sampling}
\label{app:uniform-compression}
Sampling controls maxima over all block tuples using only polynomially
many sampling points in $W$, avoiding a union bound over all $L^k$ tuples.
We deduce Lemma~\ref{uvc:thm:ordinary-sampling} from scalar Green-function
interpolation with a uniform coarse envelope, a common finite-dimensional
space for random and primitive loops, and a classical sampling inequality.

Throughout this subsection, unless otherwise stated, the volume and
regularity assumptions of Lemma~\ref{uvc:thm:ordinary-sampling}
are in force, with $g\ge C_\kappa$ sufficiently large for the auxiliary
transfer estimates. In the sampling proof these estimates are used only
when $g^2\ge W^B$, which ensures this fixed lower bound for sufficiently
large $W$.
Recall the loop definition \eqref{fo:eq:loops}.
Estimate constants are denoted by $C''$ to distinguish them from the
fixed input bound $C$ and may change between estimates; they may depend
on the fixed parameters of Lemma~\ref{uvc:thm:ordinary-sampling}.

\paragraph{Green-function interpolation and its uniform envelope.}
\begin{proposition}[Scalar Green-block interpolation]\label{uvc:prop:green}
Fix $0<\kappa<2$, $C>1$, $D>0$ and $B>2$.
For all sufficiently large $W$, assume that $g,L,E$ are deterministic and
\[
 \begin{gathered}
 g^2\ge W^B,\qquad 4g^2\le L\le g^2W^C,\\
 |E|\le2-\kappa,\qquad
 |m_{g,L}(E)|\le C,\quad \Im m_{g,L}(E)\ge C^{-1}.
 \end{gathered}
\]
For the flow \eqref{fo:eq:flow}, set $G_t=(H_t(g)-z_t)^{-1}$ along
the characteristic \eqref{fo:eq:characteristic}, $0\le t\le1-W^{-3}$.
All constants and the width threshold depend only on the fixed parameters
$\kappa,C,B,D$, not on $g,L,E$.
For $\mathcal P=\mathcal P_B(g,L,W)=\{I_j\}$ from
Definition~\ref{def:cells-anchors}, let
$b_j$ be the first block of $I_j$, with block indices taken modulo $L$.
Write
\[
 \mathcal G_t(\mathcal P):=
 \max_{\substack{I_i,I_j\in\mathcal P\\
                  \alpha\in\{b_i,b_i-1\},\ \beta\in\{b_j,b_j-1\}}}
 \|(G_t)_{\alpha\beta}\|_{\op}.
\]

For $a\in I_j$ and any block $b$, define $c_{a\alpha}$ and $q_{ab}$ by
\eqref{da:eq:cell-interpolation-coefficients} with $I=I_j$ and $b_0=b_j$,
and set $c_{a\alpha}=0$ for $\alpha\notin\{b_j,b_j-1\}$.
Thus the first index of $c$ or $q$ determines the cell used in its definition.
Set
\[
 d_{ab}=q_{ab}+\sum_{\alpha\in\{b_j,b_j-1\}}c_{a\alpha}q_{b\alpha},
\]
where $q_{b\alpha}$ is defined using the cell containing $b$.
All coefficients are real and bounded by $C''$.
Once the partition is fixed, they
depend only on $E$.
With probability at least $1-W^{-D}$, simultaneously for all
$a,b$ and $0\le t\le1-W^{-3}$,
\begin{equation}\label{uvc:eq:interpolate}
 (G_t)_{ab}
 =\sum_{\alpha,\beta}c_{a\alpha}c_{b\beta}
                       (G_t)_{\alpha\beta}
       +d_{ab}I_W+\mathcal E_{ab}(t),\qquad
 \norm{\mathcal E_{ab}(t)}_{\op}
       \le C''W^{-B/4}(1+\mathcal G_t(\mathcal P)).
\end{equation}
On the same event,
\begin{equation}\label{uvc:eq:anchor-time-envelope}
 \sup_{0\le t\le1-W^{-3}}\mathcal G_t(\mathcal P)\le W^{C+10}.
\end{equation}
The envelope exponent is independent of $B,D,g,L$.

For a fixed deterministic comparison $1/2\le s\le u\le1-W^{-3}$ as in
Lemma~\ref{da:lem:comparison-parameters}(a), with upper comparison point $w$,
the interpolation and error bound in \eqref{uvc:eq:interpolate} hold
on the same event and with the same coefficients for $G=(H_u(g)-z)^{-1}$,
uniformly over $z=\Re z_u+iy$, $\eta_u\le y\le\Im w$, with $G$ replacing
$G_t$ also in the maximum defining $\mathcal G_t(\mathcal P)$.
For the adjoint resolvents, the same bound holds on the same event,
with the same homogeneous coefficients and with $d_{ab}$ replaced by $d_{ba}$.
\end{proposition}

\paragraph{A common interpolation space for loops and primitives.}
\begin{lemma}[A common scalar interpolation space]\label{uvc:lem:space}
Under the assumptions of Proposition~\ref{uvc:prop:green}, use its
partition and envelope. Fix $1\le k\le K$.
Fix cells $I_1,\ldots,I_k\in\mathcal P$. There is a deterministic
space $\mathcal V$ of scalar functions on $X=I_1\times\cdots\times I_k$, with
$1\le\dim\mathcal V\le d_k=6^k$.
For $0\le t\le1-W^{-3}$, the second bound below holds deterministically,
and the first holds on the event of Proposition~\ref{uvc:prop:green}:
\begin{equation}\label{uvc:eq:common-space-errors}
 \begin{aligned}
 \operatorname{dist}_{\infty}(W^{k-1}\cL_k(t),\mathcal V)
 &\le W^{-B/8}(1+\mathcal G_t(\mathcal P))^k,\\
 \operatorname{dist}_{\infty}(W^{k-1}\cK_k(t),\mathcal V)
 &\le W^{-B/4}(1-t)^{-C''}.
 \end{aligned}
\end{equation}
\[
\operatorname{dist}_\infty(f,\mathcal V)
:=\inf_{v\in\mathcal V}\max_{\boldsymbol x\in X}
|f(\boldsymbol x)-v(\boldsymbol x)|.
\]
The exponent $C''$ depends only on $k$ and the fixed volume and
regularity parameters, not on $B,g,L$; the width threshold may also
depend on $B,K$. The same
space works for every time and every sign tuple of length $k$.
\end{lemma}

\paragraph{Finite-dimensional sampling.}
\begin{lemma}[Finite maximum-volume sampling]\label{uvc:lem:fekete}
If $\mathcal V$ is a $d_0$-dimensional space, $d_0\ge1$, of complex functions on
a finite set $X$, there exist $d_0$ points $\xi_1,\ldots,\xi_{d_0}$ such
that, for all $Z,F:X\to\C$,
\begin{equation}\label{uvc:eq:finite-norming}
 \begin{aligned}
 \sup_X|Z-F|&\le d_0\max_{1\le j\le d_0}|Z(\xi_j)-F(\xi_j)|\\
 &\quad +(d_0+1)\bigl(\operatorname{dist}_\infty(Z,\mathcal V)
              +\operatorname{dist}_\infty(F,\mathcal V)\bigr).
 \end{aligned}
\end{equation}
\end{lemma}

\begin{proof}
This is the classical maximum-volume (Fekete-point) argument.
Choose a basis and evaluation points maximizing the absolute determinant
of the evaluation matrix, which is nonzero by linear independence.
Cramer's rule expresses any evaluation as a combination of the selected
ones with coefficients of modulus at most one: each numerator replaces
one selected row and cannot exceed the maximal determinant.
Apply this representation to the difference of best approximations in
$\mathcal V$ and add their errors to obtain \eqref{uvc:eq:finite-norming}.
\end{proof}

\begin{proof}[\textbf{Proof of Lemma~\ref{uvc:thm:ordinary-sampling}}]
If $g^2<W^B$, take $\mathcal S_k=\Z_L^k$ and sampling is exact.
Otherwise, work on the event of Proposition~\ref{uvc:prop:green}, of
probability at least $1-W^{-D}$. By
\eqref{uvc:eq:anchor-time-envelope}--\eqref{uvc:eq:common-space-errors},
the sum of the two approximation errors in Lemma~\ref{uvc:lem:space}
is at most
\[
 2^K W^{-B/8+K(C+10)}+W^{-B/4+3C''}\le W^{-D-1},
\]
uniformly over all cell products, $k\le K$, signs, and times, after
choosing $B$ sufficiently large; both exponents $C+10$ and $C''$ are
independent of $B$.
On each product, take the at most $d_k$ sampling points supplied by
Lemma~\ref{uvc:lem:fekete} for its common space, and let $\mathcal S_k$
be their union. These sets are deterministic and independent of time
and signs. Applying \eqref{uvc:eq:finite-norming} and dividing by
$W^{k-1}$ proves \eqref{fo:eq:spatial-sampling} with $C_k=d_k$,
since the error is at most $(d_k+1)W^{-D-k}\le W^{-D}$ for large $W$.
Taking $F=0$ in the same inequality proves the assertion for $\cL_k$.
\end{proof}

\subsubsection{Green-block interpolation and uniform envelope}
\label{sec:green-interpolation-proof}
We prove Proposition~\ref{uvc:prop:green}. The short-cell transfer
estimates below give \eqref{uvc:eq:interpolate}; the subsequent
maximal-increment estimate gives \eqref{uvc:eq:anchor-time-envelope}.
Use the normalized Brownian blocks
$B_j(t):=g(\mathcal B_t^{(g)})_{jj}$ from \eqref{fo:eq:flow}.
Fix $|E|\le2-\kappa$ and use $T_j(z,t)$ from
Definition~\ref{def:transfer-matrices}
with $V_j(t)=B_j(t)$. For prefix products on a cell or interval,
number its blocks from $1$ in consecutive order, starting at its
first block, and put
\[
 P_j(z,t):=T_j(z,t)\cdots T_1(z,t).
\]
Here $1\le j\le |I|$ on the chosen interval $I$, and
$P_0(z,t)=I_{2W}$.

\begin{lemma}[Polynomially many cells]\label{uvc:lem:cells}
Fix $B>2$ and $D>0$, and assume $g^2\ge W^B$.
With probability at least $1-W^{-D}$, simultaneously for every cell
$I\in\mathcal P_B$, every $1\le j\le |I|\le\lfloor g^2W^{-B}\rfloor$,
all $0\le t\le1$, and all $z=E+i\eta$ with $0\le\eta\le g^{-2}$,
\[
 T_j(z,t)\cdots T_1(z,t)=T(E)^j+O(W^{-B/4}).
\]
The error is in operator norm. The assertion holds for sufficiently
large $W$ depending on the
fixed parameters, uniformly in $g,L$.
\end{lemma}
\begin{proof}
Set $Q_j(t):=R_E^{-j}S_E^{-1}P_j(E,t)S_E$, with $Q_0(t)=I_{2W}$.
For every fixed $p\ge2$, in any orbital dimension $W$ and on any
interval of length $n$ with $n/g^2\le c_p$, we first prove
\begin{equation}\label{uvc:eq:prefix-moments}
 \E\sup_{0\le t\le1}\max_{1\le j\le n}
        \norm{Q_j(t)-I}_{\mathrm{HS}}^p
 \le C_pW^{p/2}(n/g^2)^{p/2},
\end{equation}
where $c_p,C_p>0$ are independent of $g,n,W$.
At time $1$, the sequence $Q_j(1)$ is a spatial martingale:
conditional expectation over the new block replaces $T_j(E,1)$ by
$T(E)$. The conditional Gaussian moment bound gives
\[
 \E[\norm{Q_j-Q_{j-1}}_{\mathrm{HS}}^p\mid\mathcal F_{j-1}]
 \le C_p g^{-p}\norm{Q_{j-1}}_{\mathrm{HS}}^p.
\]
The Hilbert-space martingale maximal inequality and H\"older imply
\[
 X:=\E\max_{1\le j\le n}\norm{Q_j-I}_{\mathrm{HS}}^p
 \le C_p(n/g^2)^{p/2}
       \E\max_{0\le j\le n}\norm{Q_j}_{\mathrm{HS}}^p
 \le C_p(n/g^2)^{p/2}(W^{p/2}+X).
\]
For $n/g^2$ sufficiently small, absorb the last term.

For every fixed prefix $j$, $Q_j(t)$ is also a matrix martingale
with respect to the joint Brownian filtration of $B_1,\ldots,B_n$.
Indeed each matrix product contains each independent
block Brownian motion at most once, and cross quadratic variations
between different blocks vanish. The finite maximum
$\max_{1\le j\le n}\norm{Q_j(t)-I}_{\mathrm{HS}}$ is a nonnegative
submartingale by convexity. Doob's $L^p$ inequality in time bounds its
supremum moment by the terminal spatial maximum just estimated,
proving \eqref{uvc:eq:prefix-moments}.

Now each cell of $\mathcal P_B$ has length at most
$\lfloor g^2W^{-B}\rfloor$, and there are at most $3W^{C+B}$ cells.
Apply \eqref{uvc:eq:prefix-moments}, Markov's inequality, and a union over
these cells. The failure bound is
$C_pW^{C+B+p/2-Bp/4}$; since $B>2$, choose the fixed moment order $p$
sufficiently large after $C+B,D$. The identity
$P_j(E,t)-T(E)^j=S_ER_E^j(Q_j(t)-I)S_E^{-1}$, with $R_E$ unitary
and $S_E,S_E^{-1}$ uniformly bounded in the bulk, gives
$\|P_j(E,t)-T(E)^j\|=O(W^{-B/4})$ on this event.

Fix a cell $I$ of length $n=|I|$ and a Brownian time $u\in[0,1]$,
and suppress this time in
$T_j(z,u)$ and $P_j(z,u)$.
The real-energy prefixes and inverse prefixes have bounded norm.
Discrete variation of constants gives
\begin{equation}\label{uvc:eq:transfer-duhamel}
 P_j(z)=P_j(E)+\sum_{r=1}^j
 P_j(E)P_r(E)^{-1}[T_r(z)-T_r(E)]P_{r-1}(z).
\end{equation}
The inserted difference has norm $|z-E|$. Since $n\le g^2W^{-B}$,
discrete Gronwall and the real-energy estimate give
\[
 \|P_j(z)-T(E)^j\|
 \le C''W^{-B/4}+C''n|z-E|e^{C''n|z-E|}
 =O(W^{-B/4}),
\]
uniformly over all prefixes, cells and times whenever
$|z-E|\le C''/g^2$, for any fixed $C''$. This proves the assertion
on the same event.
\end{proof}

\begin{proof}[\textbf{Proof of Proposition~\ref{uvc:prop:green}}]
Work on the event of Lemma~\ref{uvc:lem:cells}, fix a cell $I$ of length
$n=|I|$ and a Brownian time $u\in[0,1]$, and suppress this time in
$T_j(z,u)$ and $P_j(z,u)$. The lemma's proof gives
$P_j(z)=T(E)^j+O(W^{-B/4})$ throughout $|z-E|\le C''/g^2$.
Since the real free powers and their inverses are uniformly bounded,
the inverse identity and multiplication show that
$P_j(z)$, $P_j(z)^{-1}$ and $P_j(z)P_r(z)^{-1}$ are within
$O(W^{-B/4})$ in norm of the free products
$T(E)^j$, $T(E)^{-j}$ and $T(E)^{j-r}$, respectively, uniformly over
$0\le r\le j\le n$, all cells, $u\in[0,1]$, and $|z-E|\le C''/g^2$.

Along the characteristic \eqref{fo:eq:characteristic}, $|m|\le C$ gives
$|z_t-E|\le C/g^2$. We now prove \eqref{uvc:eq:interpolate} on this
event for $0\le t\le1-W^{-3}$, writing $G=G_t$.
In the original block indices, the Green column obeys the forced recurrence
\[
 \binom{G_{r+1,b}}{G_{r,b}}
 =T_r(z_t,t)\binom{G_{r,b}}{G_{r-1,b}}
                 +\binom{\one_{r=b}I_W}{0}.
\]
For $a\in I_j$, iterate in the first coordinate from
$(G_{b_j,b},G_{b_j-1,b})$ to $G_{ab}$.
There is at most one source term, at $b$, and the homogeneous
and source coefficients differ from their scalar free counterparts
by $O(W^{-B/4})$.
Thus
\[
 \begin{aligned}
 G_{ab}&=\sum_\alpha c_{a\alpha}G_{\alpha b}
                  +q_{ab}I_W+\mathcal R_{ab},\qquad
 \|\mathcal R_{ab}\|&\le C''W^{-B/4}
       (1+\|G_{b_j,b}\|+\|G_{b_j-1,b}\|).
 \end{aligned}
\]
For $b\in I_l$, apply the same first-coordinate recurrence on $I_l$ to
$(G^*)_{b\alpha}$ at $\overline{z_t}$ and take its adjoint. The real free coefficients
give second-coordinate coefficients $c_{b\beta}$,
$\beta\in\{b_l,b_l-1\}$, and source term $q_{b\alpha}I_W$. In particular,
\[
 \max_{\substack{I_i\in\mathcal P\\\alpha\in\{b_i,b_i-1\}}}
 \|G_{\alpha b}\|\le C''(1+\mathcal G_t(\mathcal P)).
\]
This bounds $\mathcal R_{ab}$ by $C''W^{-B/4}(1+\mathcal G_t(\mathcal P))$;
substituting the second-coordinate interpolation gives
\eqref{uvc:eq:interpolate} with the stated $d_{ab}$.

For $G=(H_u(g)-z)^{-1}$ on the comparison segment in
Proposition~\ref{uvc:prop:green}, \eqref{da:eq:comparison-parameters}
and its parameter bounds give
$|z-E|\le|\Re z_u-E|+\Im w\le C''/g^2$ throughout the comparison segment.
The transfer bound obtained from \eqref{uvc:eq:transfer-duhamel} and
the Green recurrence give \eqref{uvc:eq:interpolate} for these resolvents
and, at conjugate spectral parameters, for their adjoints.

To prove \eqref{uvc:eq:anchor-time-envelope} and the increments used in
Lemma~\ref{fo:lem:sampling-continuity}, we establish
\eqref{uvc:eq:block-time-increments}.
For this increment estimate, assume only $L\gamma_g\le CW^{C_L}$
for fixed $C_L>0$ and $|m(g,E)|\le C$.
Take $0<T<1$ and use $G_s$ on the characteristic;
then $|z_s'|=\gamma_g|m(g,E)|\le C\gamma_g$. Assume that, for every
fixed $q>0$,
\begin{equation}\label{uvc:eq:momentinput}
 \sup_{a,b,\,0\le s\le T}\E\|(G_s)_{ab}\|^q\le C_qW^{3q},
\end{equation}
where $C_q$ is independent of $W,g,L$. For every fixed $p\ge2$
and block pair $a,b$, we claim that
\begin{equation}\label{uvc:eq:block-time-increments}
 \begin{aligned}
 \E\sup_{0\le s\le T}\|(G_s)_{ab}\|^p
     &\le C_pW^{(C_L+9)p},\\
 \E\|(G_v-G_u)_{ab}\|^p
     &\le C_pW^{(C_L+9)p}(v-u)^{p/2},\qquad 0\le u\le v\le T.
 \end{aligned}
\end{equation}

It\^o's formula for the flow \eqref{fo:eq:flow} gives
$d(G_s)_{ab}=A_{ab}(s)\,ds+d\mathcal M_{ab}(s)$, where
\[
 \begin{aligned}
 A_{ab}(s)&=z_s'\sum_c(G_s)_{ac}(G_s)_{cb}
       +\gamma_g\sum_c(G_s)_{ac}\langle G_s\rangle_{cc}(G_s)_{cb},\\
 d\mathcal M_{ab}(s)&=-\sum_c(G_s)_{ac}\,d(\mathcal B_s^{(g)})_{cc}\,(G_s)_{cb}.
 \end{aligned}
\]
Define the scalar quadratic variation
$[\mathcal M_{ab}]_{\HS}(s):=\sum_{\mu,\nu=1}^W
[(\mathcal M_{ab})_{\mu\nu},\overline{(\mathcal M_{ab})_{\mu\nu}}]_s$.
Independence of the Brownian blocks and the GUE covariance give
\[
 \frac{d}{ds}[\mathcal M_{ab}]_{\HS}(s)
 =\frac{\gamma_g}{W}\sum_c
       \|(G_s)_{ac}\|_{\HS}^2\|(G_s)_{cb}\|_{\HS}^2.
\]
Using \eqref{uvc:eq:momentinput}, $\|X\|_{\HS}\le\sqrt W\|X\|$
and $|z_s'|\le C\gamma_g$, Minkowski and H\"older give
\[
 \begin{aligned}
 \bigl\|\|A_{ab}(s)\|\bigr\|_{L^p}
 &\le C_p(L\gamma_g)W^9,\qquad
 \left\|\frac{d}{ds}[\mathcal M_{ab}]_{\HS}(s)\right\|_{L^{p/2}}
 &\le C_p(L\gamma_g)W^{13}.
 \end{aligned}
\]
Hilbert-space BDG and H\"older in time, with $L\gamma_g\le CW^{C_L}$ and
$T<1$, prove \eqref{uvc:eq:block-time-increments}, using
\eqref{uvc:eq:momentinput} at $s=0$ for the first bound.

We return to the large-coupling assumptions and take $T=1-W^{-3}$.
Here $|\Re z_t|\le2-\kappa/2$ and $0<\eta_t\le C/g^2$ for large $W$.
Since $g^2\eta_t=(1-t)\Im m\asymp1-t$ and $1-t\ge W^{-3}$,
\eqref{fo:eq:coarse-resolvent-moments} with margin $\kappa/2$ verifies
\eqref{uvc:eq:momentinput} uniformly in $g,L$;
$L\gamma_g=L/g^2\le W^C$ allows $C_L=C$ in the increment estimate.

Apply the first estimate in \eqref{uvc:eq:block-time-increments} and
Markov's inequality to the $O(W^{2(C+B)})$ endpoint pairs
in $\mathcal G_t(\mathcal P)$. A union bound gives
\[
 \Prob\!\left(\sup_{t\le1-W^{-3}}\mathcal G_t(\mathcal P)>W^{C+10}\right)
 \le C_pW^{2(C+B)-p}.
\]
Choose $p>2(C+B)+D+1$ and use failure exponent $D+1$ in
Lemma~\ref{uvc:lem:cells}. For sufficiently large $W$, the combined
failure probability is at most $W^{-D}$, proving
\eqref{uvc:eq:anchor-time-envelope} together with
\eqref{uvc:eq:interpolate}.
\end{proof}

\subsubsection{A common interpolation space}
\label{sec:common-interpolation-proof}

\begin{proof}[\textbf{Proof of Lemma~\ref{uvc:lem:space}}]
Set $n=\lfloor g^2W^{-B}\rfloor$ and let $b_j$ be the first block of
$I_j$, $1\le j\le k$, with cyclic subscripts.
On $I_j\times I_{j+1}$, let $\mathcal F_j$ consist of the four functions
$c_{x\alpha}c_{y\beta}$, with
$\alpha\in\{b_j,b_j-1\}$ and $\beta\in\{b_{j+1},b_{j+1}-1\}$,
and the two functions $d_{xy},d_{yx}$ from Proposition~\ref{uvc:prop:green}.
For $x_{k+1}=x_1$, set
\[
 \mathcal V=\operatorname{span}\left\{
   \prod_{j=1}^k f_j(x_j,x_{j+1}):f_j\in\mathcal F_j\right\}.
\]
One homogeneous product takes value one at $(b_1,\ldots,b_k)$,
so $1\le\dim\mathcal V\le6^k$.
The coefficients are deterministic and, for fixed cells and $E$,
independent of time, signs and orbital dimension.
Expanding \eqref{uvc:eq:interpolate} in the cyclic trace word, using
$d_{yx}$ for an adjoint block, places its interpolant in $\mathcal V$.
On the event of
Proposition~\ref{uvc:prop:green}, telescoping the product in
\eqref{uvc:eq:interpolate} and using $W^{-1}|\Tr X|\le\|X\|$
prove the random bound in \eqref{uvc:eq:common-space-errors} after
absorbing the fixed prefactor into $W^{B/8}$ for sufficiently large $W$.

\smallskip\noindent\emph{Uniform interpolation error before taking the limit.}
For the primitive bound in \eqref{uvc:eq:common-space-errors}, use the
auxiliary model of Lemma~\ref{uvc:lem:orbital-limit}, keeping the original
cells fixed. At fixed $0\le t\le1-W^{-3}$, use the transfer matrices of
Definition~\ref{def:transfer-matrices} and the cellwise prefix notation of
Section~\ref{sec:green-interpolation-proof} at width $\wt W$, with
$1\le j\le |I|\le n$. Put
$\|Y\|_{2,\wt W}=\wt W^{-1/2}\|Y\|_{\mathrm{HS}}$.
Estimate~\eqref{uvc:eq:prefix-moments} bounds the normalized Hilbert--Schmidt
$L^p$ norm of $P_j(E,t)-T(E)^j$ by $C_p\sqrt{n/g^2}$,
uniformly in $\wt W,g,L$.
Since $\Re z_t=E+O(g^{-2})$ remains in the buffered free bulk after
increasing $C_\kappa$, \eqref{uvc:eq:transfer-duhamel} and the prefix moments
of Lemma~\ref{uvs:lem:transfer-moments}, with H\"older and Minkowski, bound
the same norm of $P_j(z_t,t)-P_j(E,t)$ by
$C_p n|z_t-E|\le C_p n/g^2$.
The inverse-prefix moments and H\"older then give a
$C_p\sqrt{n/g^2}$ bound in this norm for
$P_j(z_t,t)P_r(z_t,t)^{-1}-T(E)^{j-r}$, $0\le r\le j\le |I|$.
Applying \eqref{fo:eq:coarse-resolvent-moments} at width $\wt W$ with
bulk margin $\kappa/2$ and $g^2\eta_t=(1-t)\Im m\asymp1-t$ gives
\begin{equation}\label{uvc:eq:auxiliary-block-moments}
 \sup_{b,c}\E\|(G_t^{(\wt W)})_{bc}\|^p\le C_p(1-t)^{-p}
\end{equation}
for every fixed $p$ and all $\wt W\ge \wt W_0(p,\kappa,C)$, uniformly in
$g,L,E$ in the present parameter range.
Let $\mathcal E_{ab}^{(\wt W)}(t)$ be the algebraic remainder in
\eqref{uvc:eq:interpolate} for the auxiliary model, with the same
coefficients and summation ranges. Iterate the forced Green recurrence
from the proof of Proposition~\ref{uvc:prop:green} in each coordinate.
The preceding prefix bounds, \eqref{uvc:eq:auxiliary-block-moments}, and
H\"older give
\[
 \sup_{a,b}\bigl\|\|\mathcal E_{ab}^{(\wt W)}(t)\|_{2,\wt W}\bigr\|_{L^p}
 \le C_p\sqrt{n/g^2}(1-t)^{-1}.
\]
Replacing each Green block by its
scalar interpolation in the cyclic word defines
$\wt W^{k-1}\widetilde{\cL}_k^{(\wt W)}\in\mathcal V$.
Telescoping this product, using H\"older and
$|\wt W^{-1}\Tr AB|\le\|A\|_{2,\wt W}\|B\|_{2,\wt W}$, gives
\[
 \sup_{\boldsymbol x\in X}
 \E|\wt W^{k-1}\cL_k^{(\wt W)}(t;\boldsymbol x)
       -\wt W^{k-1}\widetilde{\cL}_k^{(\wt W)}(t;\boldsymbol x)|
 \le C''\sqrt{n/g^2}(1-t)^{-C''}.
\]
The constants and width threshold remain independent of $g,L$.
This is a supremum of pointwise expected errors, not an expectation
of a spatial supremum; no microscopic union bound is used.
Since $\wt W^{k-1}\E\widetilde{\cL}_k^{(\wt W)}\in\mathcal V$,
the same bound holds for
$\operatorname{dist}_\infty(\wt W^{k-1}\E\cL_k^{(\wt W)},\mathcal V)$.

\smallskip\noindent\emph{Passage to the primitive at fixed parameters.}
At fixed $W,g,L,E$, cells and $t<1$, Lemma~\ref{uvc:lem:orbital-limit}
gives convergence in $\ell^\infty(X)$ because $X$ is finite;
since $\mathcal V$ is closed and independent of $\wt W$,
\[
 \operatorname{dist}_{\ell^\infty(X)}(W^{k-1}\cK_k,\mathcal V)
 =\lim_{\wt W\to\infty}
   \operatorname{dist}_{\ell^\infty(X)}(\wt W^{k-1}\E\cL_k^{(\wt W)},\mathcal V)
 \le C''\sqrt{n/g^2}(1-t)^{-C''}.
\]
Since $n/g^2\le W^{-B}$, absorbing the fixed prefactor into $W^{B/4}$
for sufficiently large $W$ proves the primitive bound in
\eqref{uvc:eq:common-space-errors}.
Uniformity in $g,L$ comes from the prelimit bound, without requiring
a uniform convergence rate.
\end{proof}

\subsection{Sampling and continuity for propagation}
\label{sec:sampling-propagation-extensions}
We prove Lemma~\ref{fo:lem:sampling-continuity}. We first extend the
Green-block increments proved in Section~\ref{app:uniform-compression}
to the fields used in propagation. We then treat spatial weights,
centered traces and entries, and the averaged residual, before proving
the extension from height meshes to full vertical comparison segments
at fixed times.

\subsubsection{Time increments and weighted sampling}
We prove the time-increment bounds \eqref{fo:eq:coarse-time-input} in
Lemma~\ref{fo:lem:sampling-continuity}(a) and the weighted sampling bound
\eqref{uvc:eq:weighted-pair-sampling} in Corollary~\ref{uvc:cor:weighted-sampling}.
Together, they reduce the weighted suprema used in spatial stopping,
to any prescribed polynomial accuracy, to ordinary spatial samples and
a polynomial-size time net.

\begin{proof}[\textbf{Proof of Lemma~\ref{fo:lem:sampling-continuity}(a): \eqref{fo:eq:coarse-time-input}}]
Under the volume and regularity assumptions of
Lemma~\ref{uvc:thm:ordinary-sampling}, with $g\ge C_\kappa$ sufficiently large,
the block moments verified in the proof of
Proposition~\ref{uvc:prop:green} give
$\sup_{a,b,\,0\le s\le1-W^{-3}}\E\|(G_s)_{ab}\|_{\op}^q\le C_qW^{3q}$.
The second estimate in \eqref{uvc:eq:block-time-increments}, with $C_L=C$,
therefore supplies the Green-block increments in this range.
We next prove the normalized-loop bounds
\begin{equation}\label{uvc:eq:loop-time-increments}
 \E|X(v)-X(u)|^p\le C_pW^{C_*p}|v-u|^{p/2},\qquad
 X\in\{W^{k-1}\cL_k,W^{k-1}(\cL_k-\cK_k)\},
\end{equation}
for $p\ge2$, fixed block and sign labels, and
$0\le u\le v\le1-W^{-3}$. The constant $C_*$ depends only on $C$, the fixed loop
lengths and regularity margins, not on $p,g,L$.

For $W^{k-1}\cL_k$, telescope the difference of its $k$ Green-block factors.
H\"older, \eqref{uvc:eq:block-time-increments}, and the coarse block
moments bound each term: its $L^p$ norm is at most
$C_pW^{C+9+3(k-1)}|v-u|^{1/2}$. This proves
\eqref{uvc:eq:loop-time-increments} for $X=W^{k-1}\cL_k$.
For the primitive, the following deterministic argument uses only
$\max_{a,b}|M_{B,ab}|\le C_M$, the exact Ward identity and
\eqref{fo:eq:primitive}; its constants depend only on $k,C_M$.
Set $F_k=W^{k-1}\cK_k$ and
\[
 R_k(t):=\max_{1\le j\le k}\max_{\boldsymbol a_{\ne j},\,\boldsymbol\sigma}
                  \sum_{a_j\in\Z_L}|F_k(t;\boldsymbol a,\boldsymbol\sigma)|.
\]
Here $\|F_k\|_\infty$ denotes the maximum over all labels and signs.
The entry bound and
$\gamma_g\sum_b|M_{B,ab}|^2=1$ give
$\|F_k(0)\|_\infty+\gamma_gR_k(0)\le C_k$ for $k\ge2$:
for the row, apply Cauchy--Schwarz to the two edges incident to $a_j$.
The pair formula \eqref{fo:eq:pairs}, $|Q|\le P$ and $P\mathbf1=\mathbf1$
give $\|F_2(t)\|_\infty+\gamma_gR_2(t)\le C(1-t)^{-1}$.
By \eqref{fo:eq:propagators}, the absolute row and column sums of each
one-label propagator are at most $(1-u)/(1-t)$.
Thus the $k$-label propagator has norm at most
$((1-u)/(1-t))^k$ for both $\|\cdot\|_\infty$ and $R_k$.

For a nonadjacent cut in \eqref{fo:eq:cuts}, let $k_L,k_R\ge3$ be the
two child lengths, with $k_L+k_R=k+2$. The normalized source is
$\gamma_g\sum_cF_{k_L}F_{k_R}$. Its supremum is bounded by
$\gamma_g\|F_{k_L}\|_\infty R_{k_R}$, and its one-label row by
$\gamma_gR_{k_L}R_{k_R}$: the free external label belongs to one child,
whereas $c$ belongs to both. Induction in \eqref{fo:eq:duhamel} therefore gives
\begin{equation}\label{uvc:eq:primitive-coarse-row}
 \|F_k(t)\|_\infty+\gamma_gR_k(t)
 \le C_k(1-t)^{-(2k-3)},\qquad k\ge2,\quad 0\le t<1.
\end{equation}
Indeed, $k=3$ has no source; for $k\ge4$ the Duhamel bound is
$C_k(1-t)^{-k}[1+\int_0^t(1-u)^{2-k}\,du]
\le C_k(1-t)^{-(2k-3)}$.
For $k=2$, the pair Riccati equation gives the derivative bound below.
For $k\ge3$, the generator $\Theta_{t,\boldsymbol\sigma}$ in
\eqref{fo:eq:generator} satisfies
$\|\Theta_{t,\boldsymbol\sigma}\|_{\ell^\infty\to\ell^\infty}\le k/(1-t)$.
Together with the nonadjacent-source supremum estimate and
\eqref{uvc:eq:primitive-coarse-row}, this proves
\begin{equation}\label{uvc:eq:primitive-coarse-derivative}
 \|\partial_t(W^{k-1}\cK_k(t))\|_\infty
 \le C_k(1-t)^{-(2k-2)},\qquad k\ge2.
\end{equation}
Consequently the normalized primitive increments are bounded by
$C_kW^{6k-6}|v-u|$ on $[0,1-W^{-3}]$; $\cK_1=m$ is constant.
The triangle inequality and the already proved case $X=W^{k-1}\cL_k$
of \eqref{uvc:eq:loop-time-increments} give the case
$X=W^{k-1}(\cL_k-\cK_k)$, after increasing $C_*$ if necessary.

For Lemma~\ref{fo:lem:sampling-continuity}(a), both parameter ranges
satisfy $L\gamma_g\le CW^2$. At bounded coupling,
$\eta_s^{-1}\le C''W^3$ gives the same block moment input, while
$|z_s'|=\gamma_g|m|\le C\gamma_g$ supplies the characteristic derivative bound.
Thus \eqref{uvc:eq:block-time-increments}, with $C_L=2$, gives
$\E\|(G_v-G_u)_{ab}\|_{\op}^p\le C_pW^{11p}|v-u|^{p/2}$
in both cases.
Normalized traces and entries inherit the block increment bound;
centering does not change the trace increment because $M(g,E)$ is
independent of time. Dividing \eqref{uvc:eq:loop-time-increments} by
$W^{(k-1)p}$ gives the bounds for $\cL_k$ and $D_k$ at large coupling. At bounded
coupling, $\max_{a,b}|M_{B,ab}|\le(\gamma_g\Im m)^{-1}\le C'$,
so the same product telescoping and the primitive Lipschitz
bound \eqref{uvc:eq:primitive-coarse-derivative} give these two bounds, including
the constant length-one primitive.

For the residual, solve the first identity in \eqref{gc:eq:scalar} for
$R_{ab}$ and subtract it at $v,u$. In the quadratic term, separate the
increment of the time prefactor and telescope the two trace factors.
Using $\gamma_g\sum_c|M_{B,ac}|\le C$, Minkowski and H\"older with
\eqref{uvc:eq:momentinput} and the second estimate in
\eqref{uvc:eq:block-time-increments}, with $C_L=2$, give
$\|R_{ab}(v)-R_{ab}(u)\|_{L^p}\le C_{p}W^{C}|v-u|^{1/2}$,
with $C$ independent of $p,g,L$ in the arbitrary-coupling range.
This proves the residual case of \eqref{fo:eq:coarse-time-input}.

Finally, fix $D>0$ and a deterministic collection $\mathcal S$ of at most
$W^{C_{\rm S}}$ tuples. Choose $p>4$ sufficiently large after $D,C_{\rm S}$.
Markov's inequality on successive dyadic increments, followed by summation
of their geometric probabilities and a union over $\mathcal S$, gives
a polynomial $1/4$-H\"older seminorm: for a fixed $C_{\rm H}$,
\begin{equation}\label{uvc:eq:time-holder-envelope}
 \Prob\!\left(
 \max_{\boldsymbol a\in\mathcal S}
 \sup_{0\le u<v\le1-W^{-3}}
 \frac{\|X(v;\boldsymbol a)-X(u;\boldsymbol a)\|}{|v-u|^{1/4}}
 >W^{C_{\rm H}}\right)\le W^{-D}.
\end{equation}
This holds for each field $X$ in Lemma~\ref{fo:lem:sampling-continuity}(a);
for Green blocks, normalized block traces, $W^{k-1}\cL_k$ and $W^{k-1}(\cL_k-\cK_k)$, it
holds under the assumptions of Lemma~\ref{uvc:thm:ordinary-sampling}
with $g\ge C_\kappa$ sufficiently large.
The exponent is independent of $g,L$. No energy or coupling union is used.
\end{proof}

\begin{corollary}[Weighted spatial sampling]
\label{uvc:cor:weighted-sampling}
Fix $D,D',F,c_1>0$ and $0<\theta\le1/2$.
Under the assumptions of Lemma~\ref{uvc:thm:ordinary-sampling}, with
$K\ge2$ and $g\ge C_\kappa$ sufficiently large, take $\ell_t$ from
\eqref{fo:eq:scale} and set
\[
 w_t(r)=\exp[-c_1(r/\ell_t)^\theta]+W^{-F}.
\]
For $B$ sufficiently large depending on these fixed parameters,
the deterministic sample sets $\mathcal S_k$ from that lemma have
polynomial cardinality and are independent of time and signs.
With probability at least $1-W^{-D'}$,
\begin{equation}\label{uvc:eq:weighted-pair-sampling}
 \max_{a,b\in\Z_L}\frac{|D_2(t;a,b,\boldsymbol\sigma)|}{w_t(d(a,b))}
 \le 2d_2\max_{(a,b)\in\mathcal S_2}
       \frac{|D_2(t;a,b,\boldsymbol\sigma)|}{w_t(d(a,b))}+W^{-D},
\end{equation}
simultaneously for $0\le t\le1-W^{-3}$ and all sign pairs
$\boldsymbol\sigma$, with $D_2$ as in \eqref{fo:eq:spatial-sampling}.
For $1\le k\le K$, the same bound holds for $D_k$ or $\cL_k$ on
order-$k$ tuples, with $d_k$, $\mathcal S_k$ and the weight
$w_t(\diam(\boldsymbol a))$, simultaneously for all sign words.
Sample tuples move any diameter by at most $2n$, so joint-tail
conditions transfer with this explicit enlargement of their boundary.
Together with \eqref{uvc:eq:time-holder-envelope}, any fixed polynomial
accuracy in such time-dependent suprema is determined by a
polynomial-size time net and ordinary sample set.
\end{corollary}
\begin{proof}
If $g^2<W^B$, all cells are singletons and spatial sampling is exact.
Otherwise, within each product of cells, \eqref{da:eq:weight-comparison} gives a
weight ratio at most $\exp(C(n/\ell_t)^\theta)\le2$ for large $W$,
since each distance or diameter changes by at most $2n$.
Apply Lemma~\ref{uvc:thm:ordinary-sampling} with its exponent at least
$\max\{D+F,D'\}+1$, using Lemma~\ref{uvc:lem:fekete} on each product of cells.
Dividing the remainder by the floor $W^{-F}$ leaves an error at most
$W^{-D}$, proving \eqref{uvc:eq:weighted-pair-sampling} and its
higher-order version with failure probability at most $W^{-D'-1}$.
Direct differentiation of $w_u$, using $w_u\ge W^{-F}$ and
$1-u\ge W^{-3}$, gives polynomial time moduli for $w_u$ and $w_u^{-1}$,
uniformly in $g,L$. Use \eqref{uvc:eq:time-holder-envelope} with failure
exponent $D'+1$ and combine the two events. This proves the last
assertion on a mesh $W^{-M}$ with $M/4$ larger than
the fixed accuracy and modulus exponents.
\end{proof}

\begin{remark}[Trace and entry envelopes]
The scalar traces $\langle\wt G\rangle_{ab}$ obey the Green-block
interpolation of Proposition~\ref{uvc:prop:green} directly: the particular
scalar term cancels against the same particular term in the deterministic
$M_B$ recurrence.
Thus their supremum reduces to polynomially many pairs formed from the
first blocks of the cells and their predecessors, plus an arbitrarily
small remainder. The same statement holds for
individual orbital entries after subtracting the scalar diagonal
part, with an additional polynomial union over orbital indices.
This observation concerns norming and continuity only.  It does not
replace the residual or fluctuation-averaging estimate at a fixed
pair.
\end{remark}

\subsubsection{Residual sampling and continuity in spectral height}
We prove residual sampling in Proposition~\ref{uvc:prop:residual-compression}
and the spectral-height derivative bound \eqref{fo:eq:comparison-height-interpolation}.
We combine these results with Corollary~\ref{uvc:cor:weighted-sampling}
and the preceding trace and entry envelopes to prove
Lemma~\ref{fo:lem:sampling-continuity}(b).

\begin{proposition}[Compression of the averaged Dyson residual]
\label{uvc:prop:residual-compression}
Define the exact averaged residual
\[
 R_{ab}=-\langle\wt G\rangle_{ab}
       +\frac t{g^2}\sum_cM_{B,ac}\langle\wt G\rangle_{cc}\langle G\rangle_{cb}.
\]
Under the assumptions of Lemma~\ref{uvc:thm:ordinary-sampling},
with $g\ge C_\kappa$ sufficiently large, for every
prescribed $D,D'>0$ the cell exponent $B$ can be chosen so that, with
probability at least $1-W^{-D'}$, simultaneously for every microscopic pair
and $0\le t\le1-W^{-3}$,
\[
 R_{ab}=\sum_{\alpha,\beta}c_{a\alpha}c_{b\beta}R_{\alpha\beta}
                   +O(W^{-D}).
\]
The coefficients are those of Proposition~\ref{uvc:prop:green}: for
$a\in I_j$ and $b\in I_l$, the sums use $\alpha\in\{b_j,b_j-1\}$ and
$\beta\in\{b_l,b_l-1\}$, where $b_j$ and $b_l$ are the first blocks of
the respective cells.
Consequently a pointwise residual estimate on polynomially many
pairs of these blocks yields its full microscopic supremum. The assertion
requires only a crude polynomial Green envelope and does not assume
the small-branch local law which the residual is later used to prove.
\end{proposition}
\begin{proof}
If $g^2<W^B$, retain every block and take
$c_{a\alpha}=\one_{a=\alpha}$; the assertion is then exact and
$L\le W^{C+B}$. Otherwise, on the prefix event,
Proposition~\ref{uvc:prop:green} and the coarse time envelope give
$\max_{a,b,t}\|G_{ab}(t)\|\le W^{C''}$ for a fixed exponent $C''$
independent of $B,g,L$.
Thus $|\langle\wt G\rangle_{cc}|+|\langle G\rangle_{cb}|\le C''W^{C''}$.
Use \eqref{uvc:eq:free-row-bounds} and apply
\eqref{da:eq:dyson-interpolation} on the cell containing $a$.
The hypotheses of Lemma~\ref{da:lem:free-interpolation} follow from those of
Lemma~\ref{uvc:thm:ordinary-sampling}.
Use second-coordinate interpolation on the scalar trace, with
$c_{b\beta}$ and $q_{bc}$ defined by the cell containing $b$:
\[
 \langle G\rangle_{cb}=\sum_\beta c_{b\beta}\langle G\rangle_{c\beta}
                 +q_{bc}+e_{cb}^{\mathrm{tr}},\qquad
 |e_{cb}^{\mathrm{tr}}|\le C''W^{-B/4}(1+W^{C''}).
\]
Again $q_{bc}$ is bounded and supported on at most $n$ values of
$c$.  Inserting these two formulas in
$K_{ab}=t g^{-2}\sum_cM_{B,ac}\langle\wt G\rangle_{cc}\langle G\rangle_{cb}$ gives its homogeneous
part $\sum_{\alpha,\beta}c_{a\alpha}c_{b\beta}K_{\alpha\beta}$.
Terms containing a particular solution cost at most
$C''(n/g^2)W^{2C''}$.  Terms containing $e^{\mathrm{tr}}$ are bounded with the
absolute row of $M_B$, and cost at most $C''W^{-B/4}W^{2C''}$.
Terms containing $e^M$ may be bounded crudely with
$L/g^2\le W^C$, and cost at most $C''(n/g^2)W^{C+2C''}$.
Products of errors or particular terms satisfy the same larger
bound.  Hence
\[
 \left|K_{ab}-\sum_{\alpha,\beta}
                  c_{a\alpha}c_{b\beta}K_{\alpha\beta}\right|
 \le C''W^{C+2C''}(W^{-B/4}+n/g^2).
\]
Applying \eqref{da:eq:dyson-interpolation} in both coordinates gives the
same particular term $d_{ab}$ as in \eqref{uvc:eq:interpolate}, so these
terms cancel in $M_{B,ab}-\langle G\rangle_{ab}$. Their remaining
error is at most $C''W^{-B/4}(1+W^{C''})+C''n/g^2$.
Choosing $B$ after $C,C'',D,D'$ proves the proposition.  As in the proof of
Proposition~\ref{uvc:prop:green}, the exponent $C''$ of the crude
envelope can be fixed before $B$; the moment order is chosen after
the polynomial number of cells.
\end{proof}

\begin{proof}[\textbf{Proof of Lemma~\ref{fo:lem:sampling-continuity}(b) and the comparison-height estimates}]
\noindent\textbf{Proof of \eqref{fo:eq:weighted-sampling-input}.}
At bounded coupling, or at large coupling with $g^2<W^B$, the volume
is at most $CW^{B+2}$. Retain every block tuple, so spatial norming is
an identity. For the remaining range $g^2\ge W^B$,
Lemma~\ref{uvc:thm:ordinary-sampling} treats
$\cL_k,D_k$; their weighted versions follow from
\eqref{uvc:eq:weighted-pair-sampling} and its higher-order form in
Corollary~\ref{uvc:cor:weighted-sampling}, which also gives the diameter
enlargement. The Green-block
interpolation of Proposition~\ref{uvc:prop:green} also applies to normalized
traces and orbital entries. The particular scalar terms in the random and deterministic
recurrences agree and cancel in $\langle\wt G\rangle$, and those in
an orbital entry cancel against its block trace. Thus the same
polynomial collection of pairs formed from cell first blocks and their
predecessors controls both centered fields.
Proposition~\ref{uvc:prop:residual-compression} proves the residual case.
For a weight with floor $W^{-F}$, choose the unified exponent in
Lemma~\ref{uvc:thm:ordinary-sampling} at least $\max\{D+F,D'\}$,
and increase it further by the fixed normalization and mesh costs.
This gives precisely \eqref{fo:eq:weighted-sampling-input}.

For vertical comparison in the bounded-coupling range,
$y\ge\eta_u\ge cW^{-3}$ gives $\|(H_u-z)^{-1}\|\le CW^3$,
and retaining all block tuples gives exact norming in polynomial volume.
In the arbitrary-coupling range, fix an allowed deterministic $u,z$.
The real part remains in the buffered bulk, and
$W^{-3}\lesssim\Im z/\gamma_g\lesssim1$. Moreover,
$|z-E|\le Cg^{-2}$ and
$|z+u\gamma_gm_u(z)-E|\le Cg^{-2}$, by
\eqref{da:eq:comparison-parameters} and \eqref{gc:eq:positive-stability}.
The comparison-point clause of Proposition~\ref{uvc:prop:green} supplies
the random Green interpolation at $(u,z)$. The comparison-matrix clause
of Lemma~\ref{da:lem:free-interpolation} uses the same real free energy
$E$ and the same scalar source. Equation~\eqref{gc:eq:positive-stability}
bounds the $\gamma_g$-weighted row sums of $M_{B,u}(z)$ needed for
interpolation. The coefficient in the residual identity is $u\gamma_g$.
The proof of Proposition~\ref{uvc:prop:residual-compression} consequently
applies with $M_{B,u}(z)$. For a coarse Green envelope $W^{C_G}$,
terms retaining $M_{B,u}(z)$ use its bounded $\gamma_g$-weighted
absolute row. Particular-source terms have support on at most $n$
labels and cost $C(n/g^2)W^{2C_G}$. The comparison-matrix clause of
Lemma~\ref{da:lem:free-interpolation}, applied on the cell containing $a$,
gives $|e^M_{ac}|\le Cn/g^2$,
so terms where this error replaces $M$ cost at most
\[
 C(n/g^2)(L/g^2)W^{2C_G}\le C(n/g^2)W^{2+2C_G}.
\]
Choose the cell exponent $B$ after the fixed coarse-envelope and
accuracy exponents. If $g^2<W^B$, retain every block: norming is
exact and $L\le CW^{B+2}$. Otherwise choose $n\le g^2W^{-B}$ and
apply the displayed error budgets. The scalar particular terms in
$M_{B,u}(z)_{ab}-\langle G\rangle_{ab}$ cancel as before.
Lemma~\ref{fo:lem:resolvent-moments} supplies the fixed-height block
moments. Minkowski and H\"older in the exact residual identity give,
with a fixed $C_0$ independent of $p,g,L$,
\[
 \begin{aligned}
 \|R_{ab}(u,z)\|_{L^p}
 &\le\|\langle G-M_u(z)\rangle_{ab}\|_{L^p}
  +u\gamma_g\sum_c|M_{B,u}(z)_{ac}|
   \|\langle G-M_u(z)\rangle_{cc}\|_{L^{2p}}
   \|\langle G\rangle_{cb}\|_{L^{2p}}\\
 &\le C_pW^{2C_0}.
 \end{aligned}
\]
This proves fixed-height norming and coarse moments without a small
residual estimate.

\noindent\textbf{Interpolation on comparison-height segments.}
We prove \eqref{fo:eq:comparison-height-interpolation} to extend residual
estimates from discrete height meshes to full vertical comparison segments
at fixed times, before applying Dyson stability.
Fix any polynomial collection of allowed deterministic times and their
comparison segments. For each such time $u$, write
$z=x+i\gamma_g r$ and $G_u(z)=(H_u-z)^{-1}$, with $u,x,g$ fixed and
$cW^{-3}\le r\le C$. Thus $r$ is the normalized spectral height.
Since $\partial_rG_u(z)=i\gamma_gG_u(z)^2$, the two block Ward identities and
Cauchy--Schwarz give
\begin{equation}\label{fo:eq:vertical-modulus}
 \|\partial_rG_u(z)_{ab}\|_{\op}
 \le r^{-1}\sqrt{\|\Im G_u(z)_{aa}\|_{\op}
                         \|\Im G_u(z)_{bb}\|_{\op}}.
\end{equation}
Choose a fixed coarse envelope exponent $C_G>C_0+2$ before choosing
any sampling accuracy or moment order. At deterministic height mesh
points the preceding fixed-height moments, trace and traceless-entry
samples, and the union over $W^2$ orbital indices give
$\max_{a,b}\|G_u(z)_{ab}\|_{\op}\le W^{C_G}$ with arbitrarily high
probability. Here entries are reconstructed by adding the centered
trace and the bounded Dyson entry, and
$\|A\|_{\op}\le W\max_{i,j}|A_{ij}|$ accounts for the extra factor $W$.
The moment order is chosen after the polynomial number of samples and
mesh points, so $C_G$ is independent of their number.
By \eqref{fo:eq:vertical-modulus}, while the envelope is at most
$2W^{C_G}$ every block derivative is at most $CW^{C_G+3}$.
A mesh of spacing $o(W^{-3})$ therefore excludes a first exit from
$2W^{C_G}$ and gives this derivative bound throughout each segment.

In the notation of Lemma~\ref{gc:lem:residual}, the residual identity is
\[
 R_{ab}(u,z)=-\langle\wt G_u(z)\rangle_{ab}
 +u\gamma_g\sum_cM_{B,u}(z)_{ac}
       \langle\wt G_u(z)\rangle_{cc}\langle G_u(z)\rangle_{cb}.
\]
Differentiating at fixed $u$ and applying the covariance-weighted row bounds
\eqref{gc:eq:positive-stability} and \eqref{da:eq:height-derivative} for
$M_{B,u}(z)$ and its derivative, together with the preceding Green-block
and derivative bounds, gives $|\partial_rR_{ab}(u,z)|\le CW^{2C_G+3}$,
with no extra factor $L$ or $g$.
Consequently, for a fixed $C_H>2C_G+4$ and arbitrary $D'>0$,
\begin{equation}\label{fo:eq:comparison-height-interpolation}
 \Prob\!\left(\max_{u}\sup_{r}\max_{a,b}
 \left[\|G_u(z)_{ab}\|_{\op}+\|\partial_rG_u(z)_{ab}\|_{\op}
       +|\partial_rM_{B,u}(z)_{ab}|+|\partial_rR_{ab}(u,z)|\right]
       >W^{C_H}\right)\le W^{-D'}.
\end{equation}
The maximum over $u$ is over the prescribed deterministic collection;
each supremum is over its comparison segment. Given $D>0$, choose a
height mesh of spacing $W^{-q_h}$ with $q_h>C_H+D+1$, and then choose
the sampling and pointwise moment failure exponents after the mesh
cardinality. Outside the exceptional event in \eqref{fo:eq:comparison-height-interpolation},
$G_u(z),M_u(z),R(u,z)$ differ from their nearest mesh values by at most $W^{-D}$.
This is a height assertion at fixed times; it does not use or assert
an off-characteristic time increment estimate.

For the time suprema, use \eqref{uvc:eq:time-holder-envelope} and a
mesh of spacing $W^{-M}$, with $M/4$ exceeding the fixed modulus,
weight and accuracy exponents. All assertions are uniform in fixed
deterministic parameters; no energy or coupling union and no stopped
Gaussian integration by parts is used.
\end{proof}

\subsection{Fluctuation averaging and uniform entrywise estimates}
\label{app:uniform-continuation}\label{sec:input-averaging}

\begin{proof}[Proof of Lemma~\ref{gc:lem:residual}]
Fix deterministic $u,z$ in either parameter range of the lemma,
including its permitted comparison points. All moments below use the
original Gaussian law.

\medskip
\noindent\textbf{(i) and (iii): residual estimates.}
Set $Z_{ab}\coloneqq G_{ab}-\langle G\rangle_{ab}I_W$.
For each fixed even integer $p\ge2$, fixed $C>0$, and deterministic
$K,\chi$ as in part~(i), we use the moment implication
\begin{equation}\label{gc:eq:residual-fixed-moment}
 \left.
 \begin{aligned}
 \max_{a,b}\|\langle G\rangle_{ab}\|_{L^{2p}}&\le CK,\\
 \max_{a,b,i,j}\|(Z_{ab})_{ij}\|_{L^{2p}}&\le CK\chi
 \end{aligned}\right\}
 \quad\Longrightarrow\quad
 \max_{a,b}\E|R_{ab}|^p\le C_p(K^2\chi^2)^p.
\end{equation}
The constant $C_p$ depends only on $p$, $C$ and the fixed regularity
margins, uniformly in the allowed parameters.
Appendix~\ref{sec:classical-random-matrix} proves
\eqref{gc:eq:residual-fixed-moment} and \eqref{fo:eq:residual-moments}
by the same Gaussian expansion, using the covariance-weighted Dyson
rows \eqref{gc:eq:residual-row-bounds}. This establishes part~(iii).

To deduce part~(i), fix $p$ and $\varepsilon>0$.
The stochastic-domination hypotheses, together with H\"older and the
higher coarse moments \eqref{fo:eq:coarse-time-input} on their exceptional
events, give the unconditional moment assumptions in
\eqref{gc:eq:residual-fixed-moment} with $K$ replaced by
$W^\varepsilon K$. The failure exponent is chosen after the required
moment order and accuracy. Apply \eqref{gc:eq:residual-fixed-moment}
and Markov with sufficiently large fixed $p$ at the polynomial residual
sample set from \eqref{fo:eq:weighted-sampling-input}; its norming
inequality gives \eqref{gc:eq:averaging-input}. In polynomial volume
this set can be all pairs. No sampling or stopping event is inserted
into Gaussian integration by parts.

\medskip
\noindent\textbf{(ii): entrywise estimate.}
For an independent Haar unitary $U\in U(W)$, common conjugation preserves
the deterministic hopping and the joint law of the GUE blocks. Thus
\[
 (G_{ab})_{a,b}\stackrel d=(UG_{ab}U^*)_{a,b},
 \qquad \Tr(UG_{ab}U^*)=\Tr G_{ab}.
\]
For every fixed integer $p\ge1$, deterministic
$A\in\C^{W\times W}$ and $1\le i,j\le W$, Haar integration gives
\begin{equation}\label{gc:eq:haar-fixed-moment}
 \E_U\left|(UAU^*)_{ij}-\delta_{ij}W^{-1} \Tr A\right|^{2p}
       \le CW^{-2p}\|A\|_{\HS}^{2p}.
\end{equation}
The constant depends only on $p$; the proof is in
Appendix~\ref{sec:classical-random-matrix}.
Apply \eqref{gc:eq:haar-fixed-moment} conditionally on $G_{ab}$,
then use Markov and a union bound over $i,j$ to obtain
\eqref{gc:eq:haar}. The same auxiliary $U$ acts on all blocks;
independence between resolvent blocks is not required.
\end{proof}

We next combine \eqref{gc:eq:haar-fixed-moment} with the sampling and
time increments of Section~\ref{sec:sampling-propagation-extensions}.
Fix a characteristic and a deterministic interval $I\subset[0,1-W^{-3}]$.
Let $\Psi_u$ be a positive deterministic scale satisfying
$W^{-C_*}\le\Psi_u\le W^{C_*}$ and
$|\partial_u\log\Psi_u|\le W^{C_*}$ almost everywhere, with
$\Psi$ absolutely continuous. If the fixed-pair Hilbert--Schmidt
hypothesis holds uniformly in deterministic $u,a,b$, then
\begin{equation}\label{gc:eq:haar-supremum}
 \left[W^{-2}\|(G_u)_{ab}\|_{\HS}^2\prec\Psi_u^2
       \text{ uniformly in }u,a,b\right]
 \quad\Longrightarrow\quad
 \sup_{u\in I}\max_{ai,bj}
 \frac{|(G_u)_{ai,bj}-\delta_{ij}\langle G_u\rangle_{ab}|}
      {\Psi_u}\prec1.
\end{equation}
Here and in this implication the deterministic characteristic parameters
are fixed, and the constants are uniform over their allowed range.

\begin{proof}[Proof of \eqref{gc:eq:haar-supremum}]
Fix $\epsilon,D,Q>0$ and put
$Z_{ai,bj}(u)=(G_u)_{ai,bj}-\delta_{ij}\langle G_u\rangle_{ab}$.
The conditional Haar bound \eqref{gc:eq:haar-fixed-moment} and removal
of exceptional events with \eqref{fo:eq:coarse-time-input} give
\[
 \E|Z_{ai,bj}(u)/\Psi_u|^{2p}\le C_{p,\delta}W^{2p\delta}
\]
for every fixed $p$ and $\delta>0$.
Choose deterministic traceless-entry samples from
\eqref{fo:eq:weighted-sampling-input} with failure exponent $D+3$
and accuracy $Q'>C_*+Q+2$. Their cardinality is at most
$W^{C_{\rm S}}$ for each orbital pair, with $C_{\rm S}$ now fixed.
The coarse time increments and the scale assumptions give a fixed
$C_Z$, independent of the moment order, such that
\[
 \|Z_{ai,bj}(v)/\Psi_v-Z_{ai,bj}(u)/\Psi_u\|_{L^p}
 \le C_pW^{C_Z}|v-u|^{1/2}.
\]
Indeed, split the difference into increments of $Z$ and $\Psi^{-1}$,
and use $|\partial_u\Psi_u^{-1}|\le W^{2C_*}$.
For a mesh of spacing $h=W^{-q}$, Markov at the mesh points and
within-interval dyadic subdivision (with $p>2$) bound the probabilities
of a sampled normalized value exceeding $W^\epsilon$ and of a normalized
increment exceeding $W^{-Q}$, respectively, by
\[
 C_{p,\delta}W^{C_{\rm S}+q+2-2p(\epsilon-\delta)},
 \qquad C_pW^{C_{\rm S}+q+2+p(C_Z+Q-q/2)}.
\]
Here the exponent $2$ accounts for the orbital pairs, and the dyadic
maximum has $L^p$ norm at most $C_pW^{C_Z}h^{1/2}$.
First choose $q>2(C_Z+Q+1)$, then $\delta<\epsilon$, and then the
moment orders so that both probabilities are at most $W^{-D-1}$.
The sample failure probability is at most $W^{-D-1}$ after the
orbital union; the additive norming error divided by $\Psi_u$ is
$O(W^{-Q-2})$. Spatial norming proves
\eqref{gc:eq:haar-supremum}. All moments use the original Gaussian law,
and interpolation applies only to the actual entry field.
\end{proof}

\subsection{Initial loop estimates at large coupling}
\label{sec:uniform-initialization}
We prove Lemma~\ref{fo:prop:inputs} on $0\le t\le1/2$, using
the spatial decay \eqref{uvc:eq:seed-tails} and the residual and Haar
estimates of Lemma~\ref{gc:lem:residual}.
The spatial sums below use \eqref{fo:eq:spatial-sum-budget}, uniformly
even when $g,L$ are superpolynomial. Choose each tail-floor exponent
after the fixed polynomial costs and required accuracy.

\begin{proof}[\textbf{Proof of Lemma~\ref{fo:prop:inputs}}]
Fix $g,L,E$ as in Lemma~\ref{fo:prop:inputs}. Use the flow
$H_t(g)$ and characteristic $z_t=z_t(g,E)$ from \eqref{fo:eq:flow} and
\eqref{fo:eq:characteristic}, with all Dyson quantities and scales
evaluated at $(g,E)$. The estimates are uniform in these
deterministic parameters, with spatial maxima and time suprema inside
probability. The scale definition \eqref{fo:eq:scale} gives
$\ell_t\asymp g^2$ and $\zeta_t\asymp W^{-1}$ on $0\le t\le1/2$; hence
\eqref{gc:eq:iteration-input} requires actual loops of size
$W^{1-k}$, errors of size $W^{-k}$, and the pair profile and joint
tails stated in \eqref{gc:eq:one-step-input}.
For the block operator bound, at every deterministic time
$t\le1/2$, \eqref{fo:eq:coarse-resolvent-moments} applies at $z_t$,
since $g^2\Im z_t\asymp1$ and its real part stays in the
buffered bulk. Thus $\E\|G_{ab}(t)\|^p\le C_p$ for every fixed $p$.
For sufficiently large $W$, we have
$4\le L/g^2\le CW^2\le W^3$, so the spatial sampling
estimates apply with volume exponent three on $0\le t\le1/2$.
Enlarge the fixed constant $C\ge3$ to cover the volume and regularity bounds, and take
$\mathcal P=\mathcal P_B(g,L,W)$ from Definition~\ref{def:cells-anchors}.
Let $b_j$ be the first block of each cell $I_j\in\mathcal P$.
If $g^2<W^B$, the cells are singletons and the $b_j$ run
over all blocks; otherwise Proposition~\ref{uvc:prop:green} applies.
For the polynomial collection of pairs with
$\alpha\in\{b_i,b_i-1\}$ and $\beta\in\{b_j,b_j-1\}$,
\eqref{uvc:eq:time-holder-envelope} gives a simultaneous
$1/4$-H\"older seminorm at most $W^{C''}$, with arbitrary polynomial
probability. Choose a time mesh of spacing $W^{-M}$ with
$M>4(C''+1)$. For each $\epsilon,D>0$, Markov's inequality at the
mesh points gives
\[
 \Prob\!\left(\max_{t\ {\rm on\ the\ mesh}}
       \max_{I_i,I_j\in\mathcal P}
       \max_{\alpha\in\{b_i,b_i-1\}}
       \max_{\beta\in\{b_j,b_j-1\}}
       \|G_{\alpha\beta}(t)\|>W^\epsilon\right)
 \le C_p W^{2(C+B)+M-p\epsilon}\le W^{-D},
\]
by choosing $p$ after the cell and mesh exponents. H\"older continuity
then gives $\sup_{t\le1/2}\mathcal G_t(\mathcal P)\prec1$.
Finally \eqref{uvc:eq:interpolate} has four homogeneous terms and a
bounded particular term (and no interpolation is needed for singleton cells), so
\[
 \sup_{t\le1/2}\max_{a,b}\|G_{ab}(t)\|
 \le C''\bigl(1+\sup_{t\le1/2}\mathcal G_t(\mathcal P)\bigr)\prec1.
\]
If the cells are singletons, the same mesh argument already includes
all blocks. No interpolation of a matrix norm by ordinary loops is used.
The operator bound gives
$W^{-2}\|G_{ab}(t)\|_{\HS}^2\le W^{-1}\|G_{ab}(t)\|_{\op}^2
\prec W^{-1}$. Apply \eqref{gc:eq:haar-supremum} with the constant
scale $\Psi_t=W^{-1/2}$ on $[0,1/2]$ to obtain the simultaneous
traceless-entry bound. The residual estimate \eqref{gc:eq:averaging-input}
gives $R\prec W^{-1}$. The scalar equations and
stable inverse imply, for every fixed sufficiently small $\nu>0$,
$U\le CW^{-1+\nu}+CW^\nu U^2$ with arbitrarily high probability,
where $U=\max_a|\langle\wt G\rangle_{aa}|$. At time zero $G=M$;
continuity selects the small branch, giving $\max_{a,b}|\langle\wt G\rangle_{ab}|\prec W^{-1}$.
Residual compression and the polynomial time modulus justify this
simultaneous continuity statement without a microscopic union.

For a word with one freely summed label $c$, the two incident Green
edges, Cauchy--Schwarz, and exact Ward give absolute row
$\prec g^2W^{1-k}$.  In fact the highest-length linear coefficient
uses a pair row bounded by $Cg^2/W$ with a fixed constant, because
$\max_a\Im\langle G\rangle_{aa}$ is now bounded.
For primitives, including $k=2$, \eqref{uvc:eq:primitive-coarse-row} at
coupling $g$ and $t\le1/2$ gives the deterministic row bound
$C_kg^2W^{1-k}$. The assumption $L/(g^2\log W)\to\infty$ ensures that the
primitive size and tail bounds of
Lemma~\ref{fo:lem:deterministic}(b) apply on this time interval.
Subtract the primitive equation as
$\cL_{k_L}D_{k_R}+D_{k_L}\cK_{k_R}$, where $k_L,k_R$ are the child lengths
and $k_L+k_R=k+2$. The interior cuts in \eqref{fo:eq:error-hierarchy},
the light insertion \eqref{fo:eq:light}, and the bracket
\eqref{gc:eq:bracket} have respective budgets
\[
 \begin{aligned}
 \frac W{g^2}(g^2W^{1-k_L})W^{-k_R}&=W^{-k},\qquad
 \frac W{g^2}W^{-1}(g^2W^{-k})&=W^{-k},\qquad
 \frac W{g^2}(g^2W^{-2k-1})&=W^{-2k}.
 \end{aligned}
\]
The first uses the shorter-error induction, the second the centered
one-loop bound, and the third only raw loops of length $2k+2$.
Predictably stop only on already established crude and shorter-word envelopes.
BDG at the polynomial ordinary sample tuples and the norming theorem
control the full supremum of $D_k$.  Pathwise Gronwall has a fixed
coefficient $C$, so induction gives $D_k\prec W^{-k}$.
No interpolation property of the martingale itself is assumed.

The exponential moments \eqref{uvc:eq:seed-tails}, spatial sampling
and \eqref{uvc:eq:time-holder-envelope} give uniform microscopic tails
with arbitrary polynomial floors on $0\le t\le1/2$. In particular,
the raw pair bound $W^{-1}e^{-c d(a,b)/g^2}$ and the sharp
centered pair bound $W^{-2}$ combine into
$W^{-2+o(1)}e^{-c\sqrt{d(a,b)/g^2}}$, with any prescribed
polynomial floor; their crossover is at $d(a,b)/g^2=O(\log W)$.
The tail bound has no prefactor
$\eta^{-1}\asymp g^2$, while the sharp centered bound comes
from the residual estimate and early-time loop induction above.
A long edge controls every longer word, and
\eqref{fo:eq:primitive-tail} supplies the primitive tails.
Together with Lemma~\ref{fo:lem:deterministic}(b), these estimates prove
\eqref{gc:eq:iteration-input}.
\end{proof}

The residual and initial estimates proved in
Sections~\ref{sec:input-averaging} and~\ref{sec:uniform-initialization}
supply the corresponding inputs to Theorem~\ref{gc:thm:common}(ii),
whose five-step proof and
finite time iteration are in Appendices~\ref{sec:gc-one-step}
and~\ref{sec:gc-schedule}.

\section{Transfer estimates and eigenvector concentration}
\label{sec:upper-method}

We first convert fractional-moment decay of the resolvent to
concentration, then obtain the required resolvent estimates from growth
of $D_n$. The growth estimates are proved through contraction of
$P_{1/2}$ or directly through the corrected scalar increments of
Section~\ref{sec:large-upper-proof}, according to the coupling and energy.

\subsection{Reduction to fractional-moment bounds}
\label{sec:fr-concentration}

Following the fractional-moment method of Aizenman--Molchanov
\cite{AM}, we give two sufficient resolvent inputs for
\eqref{eq:averaged-tail}: fixed-energy decay with polynomially bounded
parameters, and energy-averaged decay with constants uniform in the
coupling. We use the split at $g=W$ explained in
Section~\ref{sec:intro-upper-reductions}.
For the first route, the finite-volume criterion below reduces decay
to contraction at one deterministic length.

For a finite path or cycle $\Lambda$, let $H_\Lambda$ denote the
corresponding Hamiltonian, with Dirichlet boundary conditions on a path,
and let $\iota_x$ denote the inclusion of block $x$. Write $L=|\Lambda|$;
$\|x-y\|_L$ denotes cyclic block distance on a cycle and distance along
the path otherwise. Set
\begin{equation}\label{up:eq:local-resolvent}
 G_\Lambda(E)\coloneqq(H_\Lambda-E)^{-1},\qquad
 G_\Lambda(x,y;E)\coloneqq\iota_x^*G_\Lambda(E)\iota_y.
\end{equation}

The following conversion uses classical spectral averaging.
On bounded energy sets, the argument of
Aizenman--Schenker--Friedrich--Hundertmark already accepts
energy-averaged fractional-moment bounds; see
\cite[Appendix~A, proof of Theorem~A.1, (A.10)--(A.12)]{ASFH}.
For the energy-averaged input, we average over cells of about $g^2$
blocks. Transfer moments remain bounded on this scale, and normalization
by the cell dimension removes the coupling loss in single-block
spectral averaging. The proof retains the transfer-dependent weight
to obtain constants uniform in $g$.

\begin{lemma}[From fractional-moment bounds to concentration]
\label{uva:thm:upper-audit}
Consider the model \eqref{eq:model}. Fix $E_0>0$, $q\in(0,1)$
and $C,c>0$, and set $J=[-E_0,E_0]$. Assume one of the following
bounds.

\noindent\textup{(i) Fixed-energy bounds.} For a fixed $C_0\ge1$, suppose
$0\le g\le W^{C_0}$, $1\le\mathfrak b\le W^{C_0}$, and
\begin{equation}\label{up:eq:fmr-pointwise}
 \sup_{E\in J}\E\|G_\Lambda(x,y;E)\|^q
 \le W^C e^{-c\|x-y\|_L/\mathfrak b}.
\end{equation}
The bound is uniform over finite Dirichlet paths or cycles $\Lambda$
and blocks $x,y\in\Lambda$.

\noindent\textup{(ii) Energy-averaged bounds.} Suppose $E_0<2$,
$g\ge W$, $\mathfrak b=Wg^2$, and
\begin{equation}\label{uva:eq:conversion-input}
 \int_J\E\|G_\Lambda(x,y;E)\|^q\,dE
 \le W^C e^{-c\|x-y\|_L/\mathfrak b}.
\end{equation}
The bound is uniform over finite Dirichlet paths $\Lambda$ and
blocks $x,y\in\Lambda$.

Then, for every fixed $D>0$ and all sufficiently large $W$,
\eqref{eq:averaged-tail} holds on every finite cycle with a block radius
\begin{equation}\label{uva:eq:conversion-radius}
 R\le C_D\mathfrak b\log W.
\end{equation}
In particular, the expected number of eigenvectors in $J$ with
$\ell>C_DW\mathfrak b\log W$, divided by $WL$,
is at most $W^{-D}$. The output constants and width threshold are
independent of $g,L,\mathfrak b$ within the stated ranges.
\end{lemma}

The proof is given in Appendix~\ref{app:fr-concentration-proof}.

\noindent\textbf{Remark.}
For $g\le c_*W^{-1/2}$, the required fractional-moment estimate
is due to Peled--Schenker--Shamis--Sodin \cite[Theorem~1]{PSSS};
its application is given at the end of this subsection. At coupling of
order one, the coordinate scale $W^2$ was obtained by Drogin
\cite[Theorem~1.2 and Corollary~1.4]{Drogin} for GOE blocks and $g=1$,
in the form of Lyapunov gaps and asymptotic Green-function decay.
Our small- and intermediate-coupling proofs use his quantitative
nonconcentration criterion
\cite[Proposition~2.6 and Claims~2.7--2.8]{Drogin}; the general lift
adapts \cite[Proposition~4.4]{Drogin}.
We work with GUE blocks, retain finite-length and coupling-uniform
estimates, and supply the small-coupling construction in
Section~\ref{sec:small-upper-proof}.

For polynomially bounded coupling, the next criterion reduces
fractional-moment decay to contraction at one deterministic length.

For a finite block interval $\Lambda=[a,b]$ and deterministic Hermitian
$W\times W$ endpoint shifts $L_0,R_0$, define
\begin{equation}\label{up:eq:shifted-resolvent}
 H_\Lambda^{L_0,R_0}\coloneqq H_\Lambda+\iota_aL_0\iota_a^*
                       +\iota_bR_0\iota_b^*,\qquad
 G_\Lambda^{L_0,R_0}(E)\coloneqq(H_\Lambda^{L_0,R_0}-E)^{-1}.
\end{equation}
Write $G_\Lambda^{L_0,R_0}(x,y;E)\coloneqq
\iota_x^*G_\Lambda^{L_0,R_0}(E)\iota_y$.

The following criterion is a block version of the finite-volume
iteration in \cite[proof of Theorem~1.1, (2.23)--(2.24)]{ASFH}.

\begin{proposition}[Finite-volume criterion for fractional-moment decay]
\label{up:prop:fr-concentration}
Take the model and $C_0,E_0,g,\mathfrak b$ as in
Lemma~\ref{uva:thm:upper-audit}\textup{(i)}, with $g>0$.
Fix $s\in(0,1/2)$ and $C_1\ge1$, and allow $J$ to be any deterministic subinterval
of $[-E_0,E_0]$. In place of \eqref{up:eq:fmr-pointwise}, assume that
one deterministic integer $1\le n\le C_1\mathfrak b\log W$ satisfies
\begin{equation}\label{uv:eq:upper-volume-input}
 \sup_{E\in J}\ \sup_{L_0=L_0^*,\,R_0=R_0^*}
 \E\|G_{[0,n]}^{L_0,R_0}(0,n;E)\|^s\le W^{-2}.
\end{equation}
The endpoint shifts have no norm restriction; the expectation is over
the independent GUE blocks in $[0,n]$, and the same $n$ works for all
$E,L_0,R_0$.

Then, for all sufficiently large $W$, \eqref{up:eq:fmr-pointwise}
holds with $q=s/2$ on every finite Dirichlet path and cycle, with the
same scale $\mathfrak b$. The output constants and width threshold are
independent of $L,g,\mathfrak b,J$.
When $J=[-E_0,E_0]$, Lemma~\ref{uva:thm:upper-audit}\textup{(i)}
converts this bound to \eqref{eq:averaged-tail}.
\end{proposition}

The proof is given in Appendix~\ref{app:finite-volume-criterion}.

We record the a priori bound used in that proof and in the transfer
comparison below. For each fixed $0<p<1$, every finite block interval
$\Lambda$ and $x,y\in\Lambda$,
\begin{equation}\label{eq:c-apriori}
 \sup_{E\in\R}\sup_{L_0=L_0^*,\,R_0=R_0^*}
 \E\|G_\Lambda^{L_0,R_0}(x,y;E)\|^p
 \le C\max\{1,g\}^pW^{2p}.
\end{equation}
The constant depends only on $p$; the same bound holds for unshifted
resolvents on finite cycles. To obtain it, write
$V_x=V_x^\circ+\xi_xI_W$, where
$\xi_x\coloneqq W^{-1}\Tr V_x\sim N(0,W^{-2})$ is independent of
$V_x^\circ$. For $x=y$, conditioning on the remaining variables gives
the diagonal Schur complement
$((\xi_x/\max\{1,g\})I_W-A)^{-1}$ with $A$ Hermitian.
For $x\ne y$, condition on the traceless parts, all other blocks, and $\xi_x-\xi_y$,
leaving the common scalar $(\xi_x+\xi_y)/2$ in the two-block compression.
After scaling by $\max\{1,g\}^{-1}$, the scalar density is bounded by
$CW\max\{1,g\}$. A union over at most $2W$ eigenvalues therefore gives
the tail $\min\{1,CW^2\max\{1,g\}/h\}$ at level $h$.
Integrating it proves \eqref{eq:c-apriori}, uniformly in the endpoint
shifts and the volume, and also on cycles.

\begin{proof}[\textbf{Proof of Theorem~\ref{thm:strong-upper}}]
For sufficiently small $c_*$, \cite[Theorem~1]{PSSS}, with $W$
orbitals and exponent $1/2$, followed by summation over the $W$
coordinate vectors, gives
\[
 \sup_{E\in\R}\E\|G_\Lambda(a,b;E)\|^{1/2}
 \le W^C e^{-c|a-b|}
\]
on every finite Dirichlet path, uniformly for $g\le c_*W^{-1/2}$.
On cycles, the Schur identity \eqref{eq:c-cycle-Schur},
\eqref{eq:c-apriori} and Cauchy--Schwarz give the same decay with
exponent $q=1/4$; Jensen's inequality gives this exponent on paths
as well. Thus \eqref{up:eq:fmr-pointwise} holds with $\mathfrak b=1$.

Set $J=[-3,3]$. Since the hopping has norm $2g\le1/2$ for large $W$,
Weyl's inequalities and the one-block GUE norm tail give
\[
 \frac1{WL}\E\#\{k:\lambda_k\notin J\}
 \le\frac1L\sum_a\Prob(\|V_a\|>5/2)\le Ce^{-cW},
\]
uniformly in $L,g$.
Applying Lemma~\ref{uva:thm:upper-audit}\textup{(i)} on $J$ with
exponent $D+1$ gives an exceptional proportion at most $W^{-D-1}$,
with mass threshold $W^{-D-1}$ and radius $C_{D+1}\log W$.
For large $W$, this radius is at most $\lceil(\log W)^2\rceil$;
adding the spectral-tail count proves \eqref{eq:averaged-tail} on
$\R$, since $W^{-D-1}+Ce^{-cW}\le W^{-D}$.
\end{proof}

\subsection{From transfer growth to fractional resolvent moments}
\label{sec:upper-conversion}

The basic reduction is explained in Section~\ref{sec:intro-concentration-methods}.
For the transfer-growth route with $g<W$, high-probability growth at
one length gives the local bound \eqref{uv:eq:upper-volume-input}.
Proposition~\ref{up:prop:fr-concentration} then gives the
fixed-energy input of Lemma~\ref{uva:thm:upper-audit}\textup{(i)}.
For $g\ge W$, negative-moment decay gives the energy-averaged input
of part~\textup{(ii)} directly.

\begin{definition}[Boundary frames and transfer amplitudes]\label{def:boundary-frames}
Recall the fixed unitary matrix
\[
 S=\frac1{\sqrt2}\begin{pmatrix}1&1\\i&-i\end{pmatrix}\otimes I_W
\]
from Definition~\ref{def:transfer-matrices}, and define
\begin{equation}\label{up:eq:boundary-frame}
 \mathcal F_W\coloneqq
 \left\{\Phi\in\C^{2W\times W}:\Phi^*\Phi=I_W,\quad
                      \Phi^*\diag(I_W,-I_W)\Phi=0\right\}.
\end{equation}
Take an i.i.d. sequence of GUE blocks
with covariance \eqref{eq:covariance}, without cyclic identification.
For $g>0$ and a real energy $E$ of the normalized model
\eqref{eq:model}, use Definition~\ref{def:transfer-matrices} and abbreviate
\begin{equation}\label{eq:c-transfer}
 T_i:=T_i\!\left(\frac{\max\{1,g\}E}{g}\right).
\end{equation}
For every deterministic $\Phi\in\mathcal F_W$ and $n\ge1$, define
\begin{equation}\label{up:eq:D-definition}
 \mathcal T_n\coloneqq S^*(T_n\cdots T_1)S,\qquad
 D_n(\Phi;E)\coloneqq\sqrt2\begin{pmatrix}I_W&0\end{pmatrix}\mathcal T_n\Phi.
\end{equation}
\end{definition}

\begin{proposition}[From transfer growth to fractional resolvent moments]
\label{up:prop:conversion}
Consider the model \eqref{eq:model} and $D_n$ defined in
\eqref{up:eq:D-definition}.

\noindent\textup{(i) High-probability growth.} Fix $C_0\ge1$, $E_0>0$ and $t\in(0,1/2)$, and let
$g,\mathfrak b$ and the energy interval $J$ be deterministic, with
\[
 0<g\le W^{C_0},\qquad 1\le\mathfrak b\le W^{C_0},\qquad
 J\subset[-E_0,E_0].
\]
For a sufficiently large fixed $D'$, suppose that there are constants
$C,W_0$, independent of $W,g,\mathfrak b,J$, such that for $W\ge W_0$ some
\[
 2\log W\le n_0=n_0(W,g,\mathfrak b,J,D')
       \le C\mathfrak b\log W
\]
satisfies
\begin{equation}\label{up:eq:input}
 \sup_{E\in J}\sup_{\Phi\in\mathcal F_W}
 \Prob\bigl(\|D_{n_0}(\Phi;E)^{-1}\|>W^{-D'}\bigr)
 \le W^{-D'}.
\end{equation}
The same deterministic length
$n_0$ works for all displayed $E,\Phi$, with both suprema outside
the probability. Then \eqref{uv:eq:upper-volume-input} holds with
$s=t$ and $n=n_0$, with $C_1$ in
Proposition~\ref{up:prop:fr-concentration} chosen to cover the length
constant after $D'$ is fixed. This bound is uniform in $W,g,\mathfrak b,J$
and the Hermitian endpoint shifts; $n_0$ is independent of $E,L_0,R_0$.

\noindent\textup{(ii) Negative-moment decay.} Take $E_0,g,\mathfrak b,J$
as in Lemma~\ref{uva:thm:upper-audit}\textup{(ii)}. Suppose that, for
fixed $s\in(0,1)$ and $C,c>0$,
\begin{equation}\label{up:eq:negative-transfer-input}
 \sup_{E\in J}\sup_{\Phi\in\mathcal F_W}
 \E\|D_n(\Phi;E)^{-1}\|^s
 \le CW e^{-cn/(Wg^2)},\qquad n\ge1,
\end{equation}
uniformly in $W,g,n$. Then \eqref{uva:eq:conversion-input} holds with
$q=s/4$, uniformly on finite Dirichlet paths. The output constants
and width threshold are independent of $g$ and the path length.
\end{proposition}

\begin{proof}[Proof of Proposition~\ref{up:prop:conversion}\textup{(i)}]
Fix $0<t<1/2$ and select the input exponent $D'$ below.
For $G=G_{[0,n_0]}^{L_0,R_0}(E)$ and $u\in\C^W$, set
$f_i\coloneqq G(i,0)u$, $0\le i\le n_0$. The right boundary condition is
$f_{n_0+1}=\max\{1,g\}R_0f_{n_0}/g$. Choose a deterministic
$\Psi=\Psi(g,R_0)\in\mathcal F_W$ so that $S\Psi$ is an orthonormal
frame of the boundary subspace
\[
 \left\{\binom{v}{\max\{1,g\}R_0v/g}:v\in\C^W\right\}.
\]
This choice depends only on $g,R_0$, not on the onsite blocks or $L_0$.
Direct multiplication gives
\[
 (S^*T_iS)^*\diag(I_W,-I_W)(S^*T_iS)=\diag(I_W,-I_W).
\]
Hence, for $\Phi\in\mathcal F_W$ and $D_n=D_n(\Phi;E)$,
there is a unitary $U_n$ such that
\begin{equation}\label{eq:c-D}
 \mathcal T_n\Phi
 =2^{-1/2}\binom{D_n}{U_nD_n},\qquad
 \|(T_n\cdots T_1)S\Phi v\|_2=\|D_nv\|_2.
\end{equation}
Thus $\|D_n^{-1}\|$ is the reciprocal of the least expansion on
$S\ran\Phi$. Applying the boundary comparison of
\cite[Claim~2.2(1) and proof of Claim~2.4]{MS} with the fixed right
boundary frame $\Psi$ and the reversed product $T_1\cdots T_{n_0}$ gives
\begin{equation}\label{eq:c-boundary-comparison}
 \|G(n_0,0)\|
 \le\|D_{n_0}(\Psi;E)^{-1}\|
       \bigl(\|G(0,0)\|+\|G(1,0)\|\bigr).
\end{equation}
Here $D_{n_0}(\Psi;E)$ is formed with the reversed product,
which has the same law as $T_{n_0}\cdots T_1$ for this deterministic frame.
The boundary-pair norm dominates $\|f_{n_0}\|_2$, so no $\|R_0\|$
factor occurs; $L_0$ affects only the source equation.
On $\{\|D_{n_0}(\Psi;E)^{-1}\|\le W^{-D'}\}$ use this comparison.
The complement has probability at most $W^{-D'}$ by \eqref{up:eq:input};
there, use Cauchy--Schwarz and \eqref{eq:c-apriori} with exponent
$2t<1$. Thus
\[
 \E\|G(n_0,0)\|^t
 \le C\max\{1,g\}^tW^{2t}
          \bigl(W^{-D't}+W^{-D'/2}\bigr).
\]
Taking $D'$ sufficiently large in terms of $t,C_0$ makes the
right-hand side at most $W^{-2}$. Adjunction gives
\eqref{uv:eq:upper-volume-input} with $s=t$ and $n=n_0$.
The same deterministic $n_0$ works for every energy and both endpoint shifts.
\end{proof}

\begin{proof}[Proof of Proposition~\ref{up:prop:conversion}\textup{(ii)}]
On a finite Dirichlet path $\Lambda$, take $x<y$ and $n=y-x$.
Applying \eqref{eq:c-boundary-comparison} to the Schur complement
onto $[x,y]$ gives
\begin{equation}\label{uvg:eq:boundary}
 \norm{G_\Lambda(y,x;E)}
 \le Z(E)\bigl(\norm{G_\Lambda(x,x;E)}
                    +\norm{G_\Lambda(x+1,x;E)}\bigr).
\end{equation}
Here $Z(E)\coloneqq\norm{D_n(\Phi;E)^{-1}}$ is formed from the reversed
segment product $T_{x+1}\cdots T_y$ and the boundary frame
$\Phi\in\mathcal F_W$ determined by the right exterior.
Conditioning on the right exterior at each fixed energy,
\eqref{up:eq:negative-transfer-input} and reversal of the independent
blocks give \(\E Z(E)^s\le CW e^{-cn/\mathfrak b}\).

To control the resolvent factor after Cauchy--Schwarz in energy,
we use Boole's identity \cite[(1.5)]{PSZ}. For each realization and
vector $v\in\C^{W|\Lambda|}$, its spectral measure for $H_\Lambda$
has mass $\|v\|^2$, so
\[
 \operatorname{Leb}\{E\in\R:|\langle v,G_\Lambda(E)v\rangle|>h\}
 =\frac{2\|v\|^2}{h},\qquad h>0.
\]
Integrating the tail truncated by $|J|$ gives
$\int_J|\langle v,G_\Lambda(E)v\rangle|^s\,dE
\le C_s|J|^{1-s}\|v\|^{2s}$.
Polarization, subadditivity for $s<1$, and summation over the $W^2$
entries of a block therefore give
\begin{equation}\label{uvg:eq:energyapriori}
 \int_J\|G_\Lambda(a,b;E)\|^s\,dE
 \le C_sW^2|J|^{1-s},\qquad a,b\in\Lambda,
\end{equation}
uniformly in the volume.

Take \(p=s/2\). Cauchy--Schwarz on \(dE\,d\Prob\), using
\eqref{uvg:eq:boundary} and \eqref{uvg:eq:energyapriori}, now gives
\[
 \int_J\E\norm{G_\Lambda(y,x;E)}^p\,dE
 \le CW^{3/2}e^{-cn/(2\mathfrak b)}.
\]
Use adjunction for the other direction and \eqref{uvg:eq:energyapriori}
for the diagonal.
H\"older on the finite measure \(J\times\Omega\) lowers the exponent to
\(q=p/2=s/4\), proving \eqref{uva:eq:conversion-input} on paths.
\end{proof}

\subsection{From operator \texorpdfstring{$P_{1/2}$}{P(1/2)} contraction to transfer growth}
\label{sec:upper-inputs}

The reduction is explained in Section~\ref{sec:intro-operator-contraction}.
We give the definitions, results and proofs needed for \eqref{up:eq:input}.

\begin{definition}[Transfer batches and the operator $P_{1/2}$]
\label{up:def:transfer-operator}
Let $\Phi\in\mathcal F_W$ from \eqref{up:eq:boundary-frame},
$Y=\ran\Phi$, and $U\subset Y$ with $\dim U=W-1$.
For an invertible matrix $\mathcal T$ and either subspace $Z=U,Y$,
choose an orthonormal frame $\Phi_Z$
whose columns span $Z$ and define
\[
 \operatorname{vol}(\mathcal T|_Z)\coloneqq
       \det(\Phi_Z^*\mathcal T^*\mathcal T\Phi_Z)^{1/2},\qquad
 \operatorname{vol}(\mathcal T|_{\{0\}})\coloneqq1.
\]
Define the logarithmic volume ratio by
\begin{equation}\label{eq:c-volume-cocycle}
 c(\mathcal T;U,Y)\coloneqq
 \log\frac{\operatorname{vol}(\mathcal T|_Y)}{\operatorname{vol}(\mathcal T|_U)}.
\end{equation}
Let $q$ be a fixed positive integer, chosen separately in the lemmas
below, and set
\begin{equation}\label{up:eq:transfer-batches}
 F_\ell\coloneqq S^*T_{\ell q}\cdots T_{(\ell-1)q+1}S.
\end{equation}
These batches are independent and identically distributed, with
$F_1=\mathcal T_q$ and $F_k\cdots F_1=\mathcal T_{qk}$.
For a bounded measurable function $h$ of the pair $(U,Y)$, define
\begin{equation}\label{cc:eq:Ps}
 (P_{1/2}h)(U,Y)\coloneqq
 \E\bigl[e^{-c(\mathcal T_q;U,Y)/2}
       h(\mathcal T_q U,\mathcal T_q Y)\bigr].
\end{equation}
Write $\|h\|_\infty\coloneqq\sup_{U,Y}|h(U,Y)|$.
\end{definition}

\begin{lemma}[From $\|P_{1/2}^k1\|_\infty$ to the transfer input]
\label{up:lem:operator-reduction}
Use Definition~\ref{up:def:transfer-operator} with a fixed batch length
$q$, under the model and parameter assumptions of
Proposition~\ref{up:prop:conversion}\textup{(i)}, with the additional lower bound
$g\ge W^{-C_0}$. For each deterministic $\Phi\in\mathcal F_W$,
let $Y=\ran\Phi$ and let
$U_j$ be the span of its columns other than the $j$th.
Use identity matrices for empty products; then
$D_0=\sqrt2[I_W\ 0]\Phi$ is unitary by \eqref{up:eq:boundary-frame}.
Thus
\begin{equation}\label{cc:eq:iterate}
 \begin{gathered}
 (P_{1/2}^k1)(U,Y)=\E e^{-c(F_k\cdots F_1;U,Y)/2},\\
 \E\|D_{qk}(\Phi;E)^{-1}\|^{1/2}
 \le\sum_{j=1}^W(P_{1/2}^k1)(U_j,Y),\qquad k\ge0.
 \end{gathered}
\end{equation}
Suppose that constants $c,C>0$, uniform in the stated deterministic
$g,\mathfrak b,E$, satisfy
\begin{equation}\label{up:eq:power-input}
 \|P_{1/2}^k1\|_\infty\le W^C e^{-ck/\mathfrak b}\quad(k\ge0).
\end{equation}
For every fixed $D'>0$, a sufficiently large fixed $C_{D'}$ gives
\begin{equation}\label{up:eq:transfer-length}
 n_0\coloneqq q\left\lceil C_{D'}\mathfrak b\log W/q\right\rceil,
 \qquad 2\log W\le n_0\le 2C_{D'}\mathfrak b\log W,
\end{equation}
at which \eqref{up:eq:input} holds for all sufficiently large $W$.
The same deterministic $n_0$ works for all displayed energies and
frames, with the suprema outside the probability.
\end{lemma}

\begin{proof}[Proof of Lemma~\ref{up:lem:operator-reduction}]
The Gram matrix of the image frame in \eqref{eq:c-D} is $D_n^*D_n$. Its determinant
and the $j$th principal minor are the squared volumes on $Y$ and $U_j$.
Their ratio, by the diagonal cofactor formula for $(D_n^*D_n)^{-1}$,
gives (see the proof in Appendix~\ref{app:compact-circle})
\begin{equation}\label{eq:c-inverse-cocycle}
 \begin{aligned}
 \|e_j^*D_n^{-1}\|&=e^{-c(\mathcal T_n;U_j,Y)},\quad
 \|D_n^{-1}\|^{1/2}&\le\sum_{j=1}^W
     e^{-c(\mathcal T_n;U_j,Y)/2}.
 \end{aligned}
\end{equation}
For transfer products $\mathcal T,\mathcal T'$, multiplicativity of restricted volumes gives
\[
 c(\mathcal T'\mathcal T;U,Y)
 =c(\mathcal T;U,Y)+c(\mathcal T';\mathcal T U,\mathcal T Y).
\]
This identity and successive conditioning on the independent
batches give the first line of \eqref{cc:eq:iterate};
\eqref{eq:c-inverse-cocycle} then gives the second line. The product order is
$F_k\cdots F_1$: $F_1$ acts first.

At the multiple of $q$ in \eqref{up:eq:transfer-length},
\eqref{cc:eq:iterate} and \eqref{up:eq:power-input} give
\[
 \sup_{E\in J}\sup_{\Phi\in\mathcal F_W}
 \E\|D_{n_0}(\Phi;E)^{-1}\|^{1/2}
 \le W^{C+1-cC_{D'}/q}.
\]
Markov's inequality therefore gives
\[
 \sup_{E\in J}\sup_{\Phi\in\mathcal F_W}
 \Prob\bigl(\|D_{n_0}(\Phi;E)^{-1}\|>W^{-D'}\bigr)
 \le W^{C+1+D'/2-cC_{D'}/q}\le W^{-D'}
\]
for a sufficiently large fixed $C_{D'}$ and then large $W$.
This is exactly \eqref{up:eq:input}.
\end{proof}

\begin{definition}[Isotropic flags and the transfer group]
\label{up:def:isotropic-group}
\begin{enumerate}
\item[\textup{(i)}] Set $J_{\rm sig}\coloneqq\diag(I_W,-I_W)$.
\item[\textup{(ii)}] A subspace $Y$ is
\emph{$J_{\rm sig}$-isotropic} if $x^*J_{\rm sig}y=0$ for all $x,y\in Y$.
\item[\textup{(iii)}] A $J_{\rm sig}$-isotropic subspace $Y$ is
\emph{maximal isotropic} if $\dim Y=W$.
Thus the maximal isotropic subspaces are precisely $Y=\ran\Phi$
with $\Phi\in\mathcal F_W$.
\item[\textup{(iv)}] A pair $(U,Y)$ from
Definition~\ref{up:def:transfer-operator} is called an \emph{isotropic flag}.
\item[\textup{(v)}] Set
\begin{equation}\label{eq:c-group}
 \mathcal G\coloneqq
 \{\mathcal T\in\C^{2W\times2W}:\mathcal T^*J_{\rm sig}\mathcal T=J_{\rm sig},\ \det \mathcal T=1\}.
\end{equation}
\end{enumerate}
\end{definition}

The conjugated transfers $S^*T_iS$ belong to
$\mathcal G$: they preserve $J_{\rm sig}$, as shown in the proof of
Proposition~\ref{up:prop:conversion}\textup{(i)}, and $\det T_i=1$.

\begin{definition}[Angular $L^4$ norm]\label{up:def:angular-norm}
For each $U$ as in Definition~\ref{up:def:transfer-operator}, let
$\mathcal C_U$ be the circle of its maximal isotropic extensions $Y$.
If $e_+,e_-$ are unit eigenvectors of $J_{\rm sig}$ with eigenvalues
$+1,-1$, respectively, in $(U+J_{\rm sig}U)^\perp$, its normalized angular measure is
\[
 Y_\theta=U\oplus\C(e_++e^{i\theta}e_-)/\sqrt2,
 \qquad dm_U=d\theta/(2\pi).
\]
For a bounded measurable function $h$ of the pair $(U,Y)$, define
\[
 \|h\|_{\mathrm{ang},4}
 \coloneqq\sup_U\|h(U,\cdot)\|_{L^4(\mathcal C_U,dm_U)}
 =\sup_U\left(\int_{\mathcal C_U}|h(U,Y)|^4\,dm_U(Y)\right)^{1/4}.
\]
\end{definition}

\noindent\textbf{Definition (Simultaneous unitary invariance).}
A scalar function $h$ on these pairs $(U,Y)$ is called
\emph{simultaneously unitarily invariant} if
\begin{equation}\label{up:eq:invariant-function}
 h(\diag(\widetilde V,\widetilde V)U,
   \diag(\widetilde V,\widetilde V)Y)=h(U,Y)
\end{equation}
for every such flag $(U,Y)$ and every $W\times W$ unitary matrix
$\widetilde V$.

For the density decomposition below, let $d\mathcal T$ denote
Haar measure on $\mathcal G$.

\begin{lemma}[Good and bad components and their operator bounds]
\label{up:lem:good-bad}
For the model \eqref{eq:model}, fix $B,E_0>0$ and an integer
$q\ge1000$ divisible by 1000. Let $\mu$ be the law of
$F_1=\mathcal T_q$ from Definition~\ref{up:def:transfer-operator}.
Uniformly for deterministic $W^{-B}\le g\le W^B$ and $|E|\le E_0$, $\mu$ has,
for all sufficiently large $W$, a decomposition with a weight
$\widetilde{\pi}\in[1/2,1]$:
\begin{equation}\label{cc:eq:decomp}
 \begin{gathered}
 \mu=\widetilde{\pi} f\,d\mathcal T+\mu_{\rm bad},\qquad
 \mu_{\rm bad}\ge0,\qquad \mu_{\rm bad}(\mathcal G)=1-\widetilde{\pi},\\
 f\in C_c^\infty(\mathcal G),\qquad f\ge0,\qquad\int f\,d\mathcal T=1.
 \end{gathered}
\end{equation}
The two measures on the right of \eqref{cc:eq:decomp} are the pushforwards
of the input law weighted by $\chi_P$ and $1-\chi_P$, respectively,
under the transfer-product map. Here $\chi_P$ is a smooth cutoff.
The cutoff is independent of the incoming flag, and
\begin{equation}\label{up:eq:coarse-density}
 \begin{gathered}
 \widetilde{\pi}\ge1-Ce^{-cW}\ge\tfrac12,\qquad
 \sup_{\mathcal T\in\supp f}\max\{\|\mathcal T\|,\|\mathcal T^{-1}\|\}\le W^C,\\
 \int\|\mathcal T\|^{1/2}\,d\mu_{\rm bad}(\mathcal T)\le\tfrac12.
 \end{gathered}
\end{equation}
With
\begin{equation}\label{eq:c-good-bad-operators}
 \begin{gathered}
 (G_{1/2}h)(U,Y)\coloneqq
       \widetilde{\pi}\int e^{-c(\mathcal T;U,Y)/2}h(\mathcal T U,\mathcal T Y)f(\mathcal T)\,d\mathcal T,\\
 (B_{1/2}h)(U,Y)\coloneqq
       \int e^{-c(\mathcal T;U,Y)/2}h(\mathcal T U,\mathcal T Y)\,d\mu_{\rm bad}(\mathcal T),\qquad
 P_{1/2}=G_{1/2}+B_{1/2},
 \end{gathered}
\end{equation}
every bounded measurable $h$ satisfies
\begin{equation}\label{up:eq:GB-bounds}
 \|G_{1/2}h\|_\infty
 \le W^C\sup_U\|h(U,\cdot)\|_{L^4(\mathcal C_U,dm_U)},\qquad
 \|B_{1/2}h\|_\infty\le\tfrac12\|h\|_\infty.
\end{equation}
All constants and width thresholds depend only on $B,E_0,q$.
\end{lemma}

The proof is given in Appendix~\ref{sec:c-component-proof}. It follows
the product-differential and Gaussian integration-by-parts arguments of
\cite[Proposition~4.1, Lemma~4.2 and Section~4.3]{Drogin}, with the
GUE coercivity estimate of Appendix~\ref{app:compact-gaussian} and
the $G_{1/2},B_{1/2}$ estimates of
Appendix~\ref{app:c-operator-proofs}.

\begin{lemma}[From $L^4$ contraction to decay of $\|P_{1/2}^k1\|_\infty$]
\label{up:lem:operator-iteration}
Use the model and deterministic parameter range of
Proposition~\ref{up:prop:conversion}\textup{(i)}, with $g\ge W^{-C_0}$, and suppose that
Lemma~\ref{up:lem:good-bad} applies to the chosen fixed $q$.
Assume that, for every bounded nonnegative measurable function $h$
that is simultaneously unitarily invariant in the sense of
\eqref{up:eq:invariant-function},
one has, uniformly in the stated deterministic parameters,
\begin{equation}\label{up:eq:P-contraction}
 \sup_U\|(P_{1/2}h)(U,\cdot)\|_{L^4(\mathcal C_U,dm_U)}
 \le(1-c/\mathfrak b)
       \sup_U\|h(U,\cdot)\|_{L^4(\mathcal C_U,dm_U)},
 \qquad 0<c\le\tfrac14.
\end{equation}
Then \eqref{up:eq:power-input} holds.
\end{lemma}

\begin{proof}[Proof of Lemma~\ref{up:lem:operator-iteration}]
Simultaneous conjugation $V_i\mapsto\widetilde V V_i\widetilde V^*$
preserves the joint GUE law, so \eqref{cc:eq:Ps} implies that
$P_{1/2}^r1$ satisfies \eqref{up:eq:invariant-function} for every $r\ge0$.
Set $\rho\coloneqq1-c/\mathfrak b\ge3/4$.
Applying \eqref{up:eq:P-contraction} repeatedly to $P_{1/2}^r1$ gives
angular norm at most $\rho^r$. In the exact expansion
\begin{equation}\label{cc:eq:firstgood}
 P_{1/2}^k=B_{1/2}^k+
       \sum_{j=0}^{k-1}B_{1/2}^jG_{1/2}P_{1/2}^{k-j-1},
\end{equation}
apply the angular bound to the rightmost factor, the smoothing bound
once, and then the pointwise bound for $B_{1/2}$. After increasing $C$, this yields
\begin{equation}\label{cc:eq:pointwise}
 \|P_{1/2}^k1\|_\infty
 \le2^{-k}+W^C\sum_{j=0}^{k-1}2^{-j}\rho^{k-j-1}
 \le\left(1+\frac{W^C}{\rho-1/2}\right)\rho^k
 \le W^C e^{-ck/\mathfrak b}.
\end{equation}
Only the rightmost powers applied to $1$ need be simultaneously
unitarily invariant in the sense of \eqref{up:eq:invariant-function}.
\end{proof}

\subsection{A criterion for contraction of \texorpdfstring{$P_{1/2}$}{P(1/2)}}
\label{sec:upper-contraction-rates}
\label{sec:upper-angular}

The criterion is explained in Section~\ref{sec:intro-contraction-criterion}.
For a field
$X=(X_1,\ldots,X_q)$ on $\operatorname{Herm}(W)^q$, each $X_i$ is
Hermitian and may depend on all input blocks $V=(V_1,\ldots,V_q)$.
For a smooth scalar function
$\varphi:\operatorname{Herm}(W)^q\to\R$, write
$D_X\varphi=\sum_i\langle\nabla_{V_i}\varphi,X_i\rangle$, with
real Hilbert--Schmidt inner products. Since each real orthonormal
Hermitian coordinate of $V_i$ has variance $W^{-1}$, the Gaussian
divergence is
\begin{equation}\label{up:eq:gaussian-divergence}
 \delta_{\rm G}X\coloneqq
 W\sum_{i=1}^q\langle V_i,X_i\rangle-\operatorname{div}X,
 \qquad
 \E D_X\varphi=\E[\varphi\,\delta_{\rm G}X],
 \quad\varphi\in C_c^\infty(\operatorname{Herm}(W)^q).
\end{equation}
Here $\operatorname{div}X$ is the ordinary Euclidean divergence in
these real coordinates.

\begin{lemma}[From Gaussian divergence to contraction]\label{up:lem:angular}
Work with the model \eqref{eq:model}. Fix $B_g,B_b,E_0>0$ and a
positive integer $q$, and let $g,\mathfrak b,J$ be deterministic with
\[
 W^{-B_g}\le g\le W^{B_g},\qquad
 1\le\mathfrak b\le W^{B_b},\qquad J\subset[-E_0,E_0].
\]
Use Definition~\ref{up:def:transfer-operator} with this $q$, writing
$\mathcal T_q=\mathcal T_q(V;E)$.
Suppose that for every deterministic $E\in J$ and incoming flag
$(U,Y)$ there exist a cutoff and a vector field
\[
 \chi:\operatorname{Herm}(W)^q\longrightarrow[0,1],\qquad
 X=(X_1,\ldots,X_q):\operatorname{Herm}(W)^q
       \longrightarrow\operatorname{Herm}(W)^q.
\]
These are deterministic functions of $V$, allowed to depend on
$W,g,E,U,Y$. Input derivatives and expectations are taken with
these parameters fixed.

Require that, for every $V\in\supp\chi$, there exists a matrix
$K_{\mathrm{out}}=K_{\mathrm{out}}(V;W,g,E,U,Y)\in\C^{W\times W}$ such that
\begin{equation}\label{up:eq:angular-geometry}
 \begin{gathered}
 D_Xc(\mathcal T_q;U,Y)=1,\\
 D_X(\mathcal T_q U,\mathcal T_q Y)=\left.\frac{d}{dt}\right|_{t=0}
 \bigl(\diag(e^{tK_{\mathrm{out}}},e^{tK_{\mathrm{out}}})\mathcal T_q U,\diag(e^{tK_{\mathrm{out}}},e^{tK_{\mathrm{out}}})\mathcal T_q Y\bigr),
 \qquad K_{\mathrm{out}}^*=-K_{\mathrm{out}}.
 \end{gathered}
\end{equation}
In the right-hand derivative, $\mathcal T_q,U,Y$ and the value of
$K_{\mathrm{out}}$ at the current input are held fixed.
Assume that $\chi X$ is smooth in its active input coordinates,
extends by zero across the cutoff boundary, and is integrable
together with the derivatives needed for Gaussian integration by
parts. These regularity conditions may be checked after conditioning on
coordinates whose field components vanish identically. With the original,
unconditional Gaussian law, suppose that
\begin{equation}\label{up:eq:field}
 \E\chi\ge\tfrac12,\qquad
 \E|\delta_{\rm G}(\chi X)|^2\le C\mathfrak b,
\end{equation}
with $\delta_{\rm G}$ as in \eqref{up:eq:gaussian-divergence}.

Then the contraction estimate \eqref{up:eq:P-contraction} holds
for every bounded nonnegative measurable function $h$ that is
simultaneously unitarily invariant in the sense of
\eqref{up:eq:invariant-function}. If the fields can be chosen with
$K_{\mathrm{out}}=0$ throughout $\supp\chi$, the estimate holds for every bounded
measurable $h$, without the invariance assumption. The contraction
constant and width threshold are uniform over the deterministic
parameters whenever the field bounds are uniform.
\end{lemma}

\begin{proof}[Proof of Lemma~\ref{up:lem:angular}]
We adapt the circle nonconcentration argument of
\cite[Proposition~2.6 and Claims~2.7--2.8]{Drogin}, using Gaussian
integration by parts in the input variables. By positivity and
homogeneity, take a simultaneously unitarily invariant $h\ge0$
in the sense of \eqref{up:eq:invariant-function}, with angular
$L^4$ norm at most one. Initially take $h$ strictly positive.
Fix $U$ and a positive probability density
$\rho:\mathcal C_U\to(0,\infty)$. We test the angular $L^4$ norm
against $\rho^{3/4}$, whose $L^{4/3}(dm_U)$ norm is one, and define
\[
 R_{\mathcal T_q}(Y)\coloneqq
 \rho(Y)^{-1} e^{-2c(\mathcal T_q;U,Y)}h(\mathcal T_q U,\mathcal T_q Y)^4.
\]
Keep $U,Y,\rho(Y)$ fixed and vary only the GUE blocks. Along the local
flow of the smooth cut field $\chi X$, every smooth invariant
function of the terminal flag has zero derivative by
\eqref{up:eq:angular-geometry}. Such functions separate the orbits
of the compact unitary group, so the terminal flag stays in one orbit;
no regularity of $K_{\mathrm{out}}$ is needed. Thus even a measurable
invariant $h$ is constant along this flow. Together with
the unit cocycle derivative, this gives, in the flow sense,
\[
 D_{\chi X}\log R_{\mathcal T_q}(Y)=-2\chi.
\]
Let $\varphi_t$ be a primitive of $\mathbf1_{[-t,t]}$ with
$\|\varphi_t\|_\infty\le t$, first smoothed if necessary.
Changing variables along the local flow justifies Gaussian integration
by parts for the bounded function $\varphi_t(\log R_{\mathcal T_q}(Y))$ using its
flow derivative. Equation~\eqref{up:eq:gaussian-divergence} gives
\[
 2\E[\chi\mathbf1_{\{|\log R_{\mathcal T_q}(Y)|\le t\}}]
 =-\E[\varphi_t(\log R_{\mathcal T_q}(Y))\,\delta_{\rm G}(\chi X)]
 \le t\,\E|\delta_{\rm G}(\chi X)|
 \le Ct\sqrt{\mathfrak b}.
\]
Taking $t=c_0/\sqrt{\mathfrak b}\le1$ with sufficiently small
$c_0>0$ and using $\E\chi\ge1/2$ yields, for every fixed $Y$,
\[
 \Prob\{|\log R_{\mathcal T_q}(Y)|>t\}\ge\tfrac14.
\]
The scalar deficit satisfies
\[
 \frac34+\frac r4-r^{1/4}\ge0,\qquad
 \inf_{|\log r|>t}\left(\frac34+\frac r4-r^{1/4}\right)
 \ge ct^2\qquad(0<t\le1).
\]
For $\mathcal T\in\mathcal G$, the map $Y\mapsto \mathcal T Y$ has the exact angular Jacobian
\begin{equation}\label{cc:eq:jac}
 \frac{dm_{\mathcal T U}(\mathcal T Y)}{dm_U(Y)}=e^{-2c(\mathcal T;U,Y)}.
\end{equation}
Its proof is given under \textbf{Proof of \eqref{cc:eq:jac}} in
Appendix~\ref{app:compact-circle}.
Applying \eqref{cc:eq:jac} with $\mathcal T=\mathcal T_q$ gives
\[
 \int_{\mathcal C_U}\rho(Y)R_{\mathcal T_q}(Y)\,dm_U(Y)
 =\int_{\mathcal C_{\mathcal T_q U}}h(\mathcal T_q U,Y')^4\,dm_{\mathcal T_q U}(Y')\le1.
\]
Thus the deficit and the pointwise nonconcentration bound imply
\[
 \begin{aligned}
 \int\rho(Y)^{3/4}(P_{1/2}h)(U,Y)\,dm_U(Y)
 &=\int\rho(Y)\E R_{\mathcal T_q}(Y)^{1/4}\,dm_U(Y)\\
 &\le1-c/\mathfrak b.
 \end{aligned}
\]
Adding a positive constant and normalizing, then using dominated
convergence, removes the strict positivity assumptions on $h$ and $\rho$.
Taking the supremum over the dual tests $\rho^{3/4}$ and then over $U$
proves \eqref{up:eq:P-contraction}.
If $K_{\mathrm{out}}=0$, the endpoint flag is constant along the field,
so the same argument requires no invariance of $h$. In this case,
$|P_{1/2}h|\le P_{1/2}|h|$ extends the estimate to arbitrary bounded
measurable $h$.
\end{proof}

\begin{lemma}[Gaussian divergence bounds at small coupling]
\label{up:lem:small-contraction}
Use the model \eqref{eq:model}, fix $\kappa_0\in(0,2)$, and let
$J=[-2+\kappa_0,2-\kappa_0]$. Consider the following two ranges.

\noindent\textup{(i)} Fix $0<\alpha<\beta<1/2$ and $\varepsilon>0$.
For each deterministic $a\in[\alpha,\beta]$, set $g=W^{-a}$.

\noindent\textup{(ii)} Fix $c_*>0$ and $0<\theta\le1/100$, and let
$g$ be any deterministic coupling in the critical window
\begin{equation}\label{eq:critical-window}
 c_*W^{-1/2}\le g\le W^{-1/2+\theta/2}.
\end{equation}

In each range there is a fixed local length $q$ such that, for all
sufficiently large $W$, every deterministic $E\in J$ and incoming
flag admit $\chi,X$ satisfying the geometric and regularity assumptions
of Lemma~\ref{up:lem:angular}, with $0\le\chi\le1$, $\E\chi\ge1/2$, and
\begin{equation}\label{up:eq:small-transfer-scales}
 \E|\delta_{\rm G}(\chi X)|^2\le C\mathfrak b,
 \qquad \mathfrak b\coloneqq
 \begin{cases}
 W^{\varepsilon/2}Wg^2,&\text{in (i)},\\
 W^{18\theta},&\text{in (ii)}.
 \end{cases}
\end{equation}
The local length, constant $C$ and width threshold depend only on
$\alpha,\beta,\varepsilon,\kappa_0$ in (i), and on
$c_*,\theta,\kappa_0$ in (ii). They are uniform over the displayed
deterministic couplings, energies and flags.
\end{lemma}

The proof of Lemma~\ref{up:lem:small-contraction} is given in
Section~\ref{su:sec:generator-bounds}.

\begin{lemma}[Gaussian divergence bounds from the product differential]
\label{up:lem:density-contraction}
Use the model \eqref{eq:model}, fix $B,E_0>0$, and take $q=1000$.
For all sufficiently large $W$ and every deterministic
$W^{-B}\le g\le W^B$, $|E|\le E_0$ and incoming flag, there are
$\chi_P,X$ satisfying the geometric and regularity assumptions of
Lemma~\ref{up:lem:angular}, with $K_{\mathrm{out}}=0$, $0\le\chi_P\le1$ and
$\E\chi_P\ge1/2$, and the following estimates.

\noindent\textup{(i)} Throughout this range,
\begin{equation}\label{eq:c-common-scale}
 \E|\delta_{\rm G}(\chi_P X)|^2\le C\mathfrak b,
 \qquad \mathfrak b\coloneqq W\max\{g^{-2},g^4\}.
\end{equation}
The constant and width threshold depend only on $B,E_0$.

\noindent\textup{(ii)} If $g\ge1$ and
$\kappa_1\le|E|\le E_0$, where $0<\kappa_1<E_0$ is fixed, one can
choose a field with
\[
 \E|\delta_{\rm G}(\chi_P X)|^2\le C_{\kappa_1}Wg^2.
\]
The constant and width threshold may additionally depend on $\kappa_1$.
The same cutoff is used in both estimates and is independent of the
incoming flag; the fields may differ.
\end{lemma}

The proof is given in Section~\ref{sec:compact-proof}.

\subsection{Concentration at small and intermediate coupling}
\label{sec:upper-center-completion}

The constructions in Sections~\ref{sec:small-upper-proof}
and~\ref{sec:compact-proof} provide the transfer input
\eqref{up:eq:input}. Proposition~\ref{up:prop:conversion}\textup{(i)}
then gives the local resolvent bound
\eqref{uv:eq:upper-volume-input}; the finite-volume criterion and
Lemma~\ref{uva:thm:upper-audit}\textup{(i)} give concentration.

\begin{proof}[Proof of Theorems~\ref{thm:sharp-small-upper}
and~\ref{thm:compact-upper}]
For either construction below, apply Lemma~\ref{up:lem:angular}
to a field of local length $q_0$.
For $Q=\operatorname{lcm}(q_0,1000)$, independence and the cocycle
identity imply that the operator of a $Q$-block batch is
$P_{1/2}^{Q/q_0}$, where $P_{1/2}$ denotes the local $q_0$-block
operator. Since $P_{1/2}$ preserves simultaneous unitary invariance
\eqref{up:eq:invariant-function},
the contraction factor of $P_{1/2}^{Q/q_0}$ is at most
$(1-c/\mathfrak b)^{Q/q_0}\le1-c/\mathfrak b$.
At this fixed length $Q$, Lemmas~\ref{up:lem:good-bad}
and~\ref{up:lem:operator-iteration} therefore apply.
This uses operator iteration and requires no extension of the local
field to the longer product.

For Theorem~\ref{thm:sharp-small-upper}, set
\[
 \theta=\min\left\{\frac{\varepsilon}{200},
                \frac{1-2\alpha}{4},\frac1{100}\right\},
 \qquad \beta=\frac12-\frac\theta2>\alpha.
\]
If $c_*W^{-1/2}\le g\le W^{-\beta}$, use
Lemma~\ref{up:lem:small-contraction}\textup{(ii)} with
$\kappa_0=\kappa$; its scale is $W^{18\theta}\le W^{\varepsilon/2}$.
If $W^{-\beta}\le g\le W^{-\alpha}$, use part~\textup{(i)} of
that lemma with loss parameter $\varepsilon$; its scale is
$W^{\varepsilon/2}Wg^2$. In both cases,
Lemmas~\ref{up:lem:angular}, \ref{up:lem:operator-iteration}
and~\ref{up:lem:operator-reduction} give \eqref{up:eq:input}
on $J=[-2+\kappa,2-\kappa]$ with a scale at most
$\mathfrak b=W^{\varepsilon/2}\max\{1,Wg^2\}$.
Fix $0<q<1/4$. Proposition~\ref{up:prop:conversion}\textup{(i)},
with $t=2q$, gives \eqref{uv:eq:upper-volume-input}.
Proposition~\ref{up:prop:fr-concentration} now proves
\eqref{up:eq:fmr-pointwise} on paths and cycles.
In each application, first choose $C_0$ to cover the polynomial ranges
of $g,\mathfrak b$, then $D'$ as required by
Proposition~\ref{up:prop:conversion}\textup{(i)}, and finally $C_1$ in
Proposition~\ref{up:prop:fr-concentration} to cover the resulting length constant.
All parameters are fixed before $W,g,L$; the two ranges meet at
$g=W^{-\beta}$, so taking the larger constants gives uniformity
on their union.
Lemma~\ref{uva:thm:upper-audit}\textup{(i)} gives
\eqref{eq:averaged-tail} with
\[
 R\le C_DW^{\varepsilon/2}\max\{1,Wg^2\}\log W.
\]
For large $W$, absorbing $C_D\log W$ into the remaining power
proves \eqref{eq:sharp-small-radius} and the coordinate bound
$\min\{WL,W^\varepsilon\max\{W,W^2g^2\}\}$.

For Theorem~\ref{thm:compact-upper}, use
Lemma~\ref{up:lem:density-contraction}\textup{(i)} with batch length
$1000$ and $J=[-E_0,E_0]$. The same chain gives
\eqref{up:eq:input}, then \eqref{uv:eq:upper-volume-input}
and \eqref{up:eq:fmr-pointwise}, with
$\mathfrak b=W\max\{g^{-2},g^4\}\le W^{1+4B}$.
Lemma~\ref{uva:thm:upper-audit}\textup{(i)} gives
$R\le C_DW\max\{g^{-2},g^4\}\log W$; absorbing $C_D$ into
$\log W$ proves \eqref{eq:coarse-radius}.
Both bounds hold for every deterministic $g,L$ in the stated
ranges, with uniform constants and the finite-volume cap $WL$.
\end{proof}

\section{Vector-field construction and estimates at small coupling}
\label{sec:small-upper-proof}

The construction is explained in Section~\ref{sec:intro-small-fields}.
Here we prove Lemma~\ref{up:lem:small-contraction}; the passage to
contraction is given by Lemma~\ref{up:lem:angular}, following
\cite[Proposition~2.6 and Claims~2.7--2.8]{Drogin}.

\subsection{Constructing \texorpdfstring{$X\in(\C^{W\times W})^q$}{X in (C\string^(W x W))\string^q} from \texorpdfstring{$K\in\C^{W\times W}$}{K in C\string^(W x W)}}
\label{su:sec:construction}

Fix $0<g\le1$, $E\in\R$, an isotropic flag $(U,Y)$, and an integer
$q\ge7$. For the original input variables and transfer product, write
\[
 V=(V_1,\ldots,V_q)\in\operatorname{Herm}(W)^q,\qquad
 \mathcal T_q=\mathcal T_q(V;E)=S^*T_q\cdots T_1S.
\]
Here $T_i$ is defined in \eqref{eq:c-transfer}.
\paragraph{Data used to define the field.}
We apply the original transfers $T_j$ to $SU,SY$, using the fixed
unitary matrix $S$ from Definition~\ref{def:transfer-matrices}.
For each fixed incoming flag, the first-block projections of
$T_j\cdots T_1SY$, $1\le j\le q$, are almost surely invertible.
Define the Hermitian matrices $R_j$ and the normal lines $\C v_j$ by
\begin{equation}\label{su:eq:graph-input}
 T_j\cdots T_1SY
   =\operatorname{ran}\begin{pmatrix}I_W\\-gR_j\end{pmatrix}, \quad T_j\cdots T_1SU
   =\left\{\begin{pmatrix}x\\-gR_jx\end{pmatrix}:
                     x\in\C^W,\ v_j^*x=0\right\}.
\end{equation}
Choose a nonzero representative $v_1$. The subsequent representatives
can be chosen so that
\begin{equation}\label{su:eq:transported-normal}
 R_j=(V_j-EI_W-g^2R_{j-1})^{-1},\qquad
 v_j=R_jv_{j-1},\qquad 2\le j\le q.
\end{equation}
Call an inverse state $R$ regular if
\begin{equation}
 \norm R\le W^3,\qquad \rank\one_{\{|g^2 R|>g\}}\le Wg^{1/2}.
 \label{us:eq:regular}
\end{equation}
For the two regimes of Lemma~\ref{up:lem:small-contraction} set
\begin{equation}\label{su:eq:source-projection}
 P\coloneqq
 \begin{cases}
 \mathbf1_{\{|g^2R_4|\le W^{\varepsilon/800}\}},&\text{in regime (i)},\\
 I_W,&\text{in regime (ii)}.
 \end{cases}
\end{equation}

\begin{lemma}[Estimates after the first four blocks]\label{su:lem:preparation}
Fix a deterministic incoming flag. In either regime, with probability
$1-o(1)$ uniformly over the flag and the stated energies and couplings,
$R_4$ is regular and the following bounds hold.
In regime (i),
\begin{equation}\label{su:eq:prep}
 \norm{R_4}\le W^3,\qquad \|Pv_4\|\ge W^{-10}\|v_4\|,\qquad
 b_4=v_4^*(I+g^4R_4^2)^{-1}v_4\le C\|Pv_4\|^2,
\end{equation}
and
\begin{equation}\label{su:eq:prepared-large-rank}
 \rank\one_{\{|g^2R_4|>W^{-\varepsilon/800}\}}
 \le Wg^2W^{\varepsilon/400}.
\end{equation}
In regime (ii),
\begin{equation}\label{nc:eq:prepared}
 \norm{R_4}\le W^{1+\theta},\qquad
 \rank\one_{\{|g^2R_4|>W^{-\theta}\}}\le W^{3\theta},\qquad
 g^2\norm{R_4}\le W^{2\theta}.
\end{equation}
\end{lemma}
The proof, including first-step invertibility, is in
Appendix~\ref{su:sec:preparation}.

We hold the chosen representative $v_4$ fixed when varying
$V_5,\ldots,V_q$.
For fixed $E,g,U,Y$, let $\chi\in[0,1]$ and $K\in\C^{W\times W}$
depend only on $V_1,\ldots,V_{q-1}$.
Write $R_{[q-1,5]}\coloneqq R_{q-1}\cdots R_5$ and assume
\begin{equation}\label{su:eq:generator-normalization}
 K^*=-K,\qquad
 \operatorname{Re}\bigl(v_{q-1}^*R_{[q-1,5]}Kv_4\bigr)=\|v_{q-1}\|^2
       \quad\text{on }\supp\chi.
\end{equation}

\paragraph{The formula for $X$.}
For $V\in\supp\chi$, first form the correction vector and matrix
\begin{equation}\label{us:eq:L}
 \begin{aligned}
 v_{\rm corr}&\coloneqq R_{[q-1,5]}Kv_4-v_{q-1},\\
 L&\coloneqq
 \|v_{q-1}\|^{-2}\bigl(v_{\rm corr}v_{q-1}^*-v_{q-1}v_{\rm corr}^*\bigr)
 -\|v_{q-1}\|^{-4} v_{q-1}^*v_{\rm corr}\,v_{q-1}v_{q-1}^*.
 \end{aligned}
\end{equation}
With $[A,B]=AB-BA$, define $X$ on $\supp\chi$ by the following
formula and set $X=0$ outside $\supp\chi$:
\begin{equation}\label{us:eq:physical}
 \boxed{X_i(V)\coloneqq
 \begin{cases}
 0,&1\le i\le4,\\
 [K,V_5]-g^2[K,R_4],&i=5,\\
 [K,V_i],&6\le i\le q-1,\\
 [K,V_q]-[L,V_q-g^2R_{q-1}],&i=q.
 \end{cases}}
\end{equation}

\begin{lemma}[The constructed field satisfies the geometric conditions]
\label{su:lem:field-geometry}
Let $\chi,K$ be as above, satisfying \eqref{su:eq:generator-normalization},
and define $L,X$ by \eqref{us:eq:L} and \eqref{us:eq:physical}.
Then every $X_i$ is Hermitian and, for each $V\in\supp\chi$, we have
\begin{equation}\label{su:eq:field-geometry}
 \begin{gathered}
 D_Xc(\mathcal T_q;U,Y)=1,\\
 D_X(\mathcal T_q U,\mathcal T_q Y)=\left.\frac{d}{dt}\right|_{t=0}
 \bigl(\diag(e^{t(K-L)},e^{t(K-L)})\mathcal T_q U,
       \diag(e^{t(K-L)},e^{t(K-L)})\mathcal T_q Y\bigr).
 \end{gathered}
\end{equation}
\end{lemma}
\begin{proof}
Equations~\eqref{su:eq:generator-normalization} and \eqref{us:eq:L} give
\begin{equation}\label{su:eq:correction-identities}
 \Re(v_{q-1}^*v_{\rm corr})=0,\qquad
 L^*=-L,\qquad Lv_{q-1}=v_{\rm corr}.
\end{equation}
Hermiticity follows from \eqref{us:eq:physical} and the skew-Hermitian
properties in \eqref{su:eq:generator-normalization} and
\eqref{su:eq:correction-identities}. Since the first four components vanish,
$D_XR_4=D_Xv_4=0$. Put $M_q\coloneqq V_q-EI_W-g^2R_{q-1}$.
Differentiating \eqref{su:eq:transported-normal},
using \eqref{us:eq:physical} and \eqref{su:eq:correction-identities}, gives
\begin{equation}\label{us:eq:orbit}
 \begin{aligned}
 D_XR_j&=-R_j\bigl(X_j-g^2D_XR_{j-1}\bigr)R_j
          =[K,R_j],\qquad 5\le j\le q-1,\\
 D_Xv_{q-1}&=[K,R_{[q-1,5]}]v_4=(K-L)v_{q-1}-v_{q-1},\\
 D_XM_q&=[K-L,M_q].
 \end{aligned}
\end{equation}
The terminal spaces have the parametrizations
\[
 \begin{aligned}
 S\mathcal T_qY&=\{(-g^{-1}M_qx,x):x\in\C^W\},\\
 S\mathcal T_qU&=\{(-g^{-1}M_qx,x):v_{q-1}^*x=0\},
 \end{aligned}
\]
which do not require $M_q$ to be invertible. The base-times-height
formula and cocycle additivity give
\begin{equation}\label{su:eq:graph-volume-ratio}
 \begin{aligned}
 e^{-2c(\mathcal T_q;U,Y)}
 &=e^{-2c(\mathcal T_4;U,Y)} \cdot
 g^{2(q-4)}\cdot \frac{{v_{q-1}^*(g^2I_W+M_q^2)^{-1}v_{q-1}}}
        {v_4^*(I_W+g^2R_4^2)^{-1}v_4}
 \end{aligned}
\end{equation}

Where $R_q$ exists, the numerator equals
$v_q^*(I_W+g^2R_q^2)^{-1}v_q$; the displayed form is defined for every $V_q$.
Only this terminal quadratic form varies in \eqref{su:eq:graph-volume-ratio}.
The terminal identities in \eqref{us:eq:orbit} cancel the commutator terms and give
\[
 D_X\bigl(v_{q-1}^*(g^2I_W+M_q^2)^{-1}v_{q-1}\bigr)
 =-2v_{q-1}^*(g^2I_W+M_q^2)^{-1}v_{q-1}.
\]
Taking a logarithmic derivative of \eqref{su:eq:graph-volume-ratio}
yields $D_Xc(\mathcal T_q;U,Y)=1$.

At the current input $V$, freeze the values $K(V),L(V)$.
The last two equations in \eqref{us:eq:orbit} are the infinitesimal variation of
\[
 (M_q,v_{q-1})\longmapsto
 \bigl(e^{t(K-L)}M_qe^{-t(K-L)},\ e^{-t}e^{t(K-L)}v_{q-1}\bigr).
\]
In the terminal parametrizations above, this sends
$x$ to $e^{t(K-L)}x$; the scalar factor $e^{-t}$ does not change the
normal's kernel. Thus $S\mathcal T_q U,S\mathcal T_q Y$ undergo the same block-diagonal
unitary rotation. Since this rotation commutes with $S$, the second
conclusion in \eqref{su:eq:field-geometry} follows.
\end{proof}

\subsection{Construction of \texorpdfstring{$K$}{K}}
\label{su:sec:generator-choice}
\label{su:sec:explicit-generator}

Fix either regime of Lemma~\ref{up:lem:small-contraction} and condition on
the first four blocks. Use the factorization
$R_{[q-1,5]}=A\widetilde R$ from Section~\ref{sec:intro-small-generator}.

\begin{proof}[Proof of \eqref{us:eq:cavity}]
Consider the block tridiagonal matrix on $\{5,\ldots,q-1\}$
with first diagonal block $V_5-EI_W-g^2R_4$, remaining diagonal
blocks $V_j-EI_W$, and adjacent off-diagonal blocks $gI_W$.
Whenever the inverses exist, successive Schur complements from
the left and from the right express its $(q-1,5)$ inverse block
as $(-g)^{q-6}R_{[q-1,5]}$ and
$(-g)^{q-6}A\widetilde R$, respectively.
Cancelling the common factor proves \eqref{us:eq:cavity}.
\end{proof}
At the deterministic threshold $t_{\mathrm{cut}}>0$ in
\eqref{su:eq:construction-threshold}, set
\[
 A_>\coloneqq A\cdot\mathbf1_{(t_{\mathrm{cut}}^2,\infty)}(A^*A),\qquad
 A_<\coloneqq A-A_>.
\]
Let $Q$ be the orthogonal projection with
\begin{equation}\label{su:eq:retained-projection}
 \ran Q=\widetilde R\bigl(\ran(A_>^*)\bigr).
\end{equation}
Thus $A_>\widetilde R(I-Q)=0$ since $\widetilde R$ is Hermitian,
and $\|A_<\|\le t_{\mathrm{cut}}$ by the singular-value cutoff.
Use $P$ from \eqref{su:eq:source-projection}.
In (i) it suffices to take $\varepsilon>0$ small in terms of
$\alpha,\beta$, since decreasing $\varepsilon$ strengthens the bounds.
For $Pv_4\ne0$ let $\Pi$ project onto $\ran P\cap v_4^\perp$ and set
\[
 \Pi\coloneqq P-\|Pv_4\|^{-2} (Pv_4)\cdot(Pv_4)^*,\qquad
 B_Q\coloneqq Q\Pi Q\big|_{\ran Q}.
\]
If $B_Q>0$, define
\begin{equation}
 y\coloneqq\Pi QB_Q^{-1}Qv_4,\qquad
 K_0\coloneqq\|Pv_4\|^{-2}\bigl(y\cdot(Pv_4)^*-Pv_4\cdot y^*\bigr).
 \label{us:eq:H0}
\end{equation}
Then $Py=y$, $v_4^*y=0$ and $Qy=Qv_4$, since $Q\Pi Q=B_Q$ on
$\ran Q$. Hence
\[
 K_0^*=-K_0,\qquad K_0=PK_0P,\qquad \rank K_0\le2,\qquad K_0v_4=y.
\]
Since $A_>\widetilde R(y-v_4)=0$, the matching error is
\begin{equation}\label{us:eq:error}
 e_{q-1}\coloneqq A_<\widetilde R(y-v_4)=R_{[q-1,5]}(y-v_4).
\end{equation}
For the normalization \eqref{su:eq:generator-normalization}, set
\begin{equation}\label{us:eq:normalize}
 \xi\coloneqq\|v_{q-1}\|^{-2} \Re(v_{q-1}^*R_{[q-1,5]}y),\qquad
 K\coloneqq K_0/\xi.
\end{equation}
These formulas are used on the cutoff support, where
\begin{equation}\label{su:eq:generator-admissibility}
 \|Pv_4\|>0,\qquad B_Q\ge\tfrac12 I_{\ran Q},\qquad
 \|y\|\le\|v_4\|,\qquad \|A_<\widetilde R\|\,\|v_4\|\le\tfrac14\|v_{q-1}\|.
\end{equation}
Under these conditions, $\|B_Q^{-1}\|\le2$, and the two orthogonal
rank-one terms in \eqref{us:eq:H0} give
\begin{equation}\label{us:eq:sourceenergy}
 \|y\|\le2\|Qv_4\|,\qquad
 \|K_0\|_{\HS}^2=2\|Pv_4\|^{-2} \|y\|^2\le8\|Pv_4\|^{-2} \|Qv_4\|^2.
\end{equation}
By \eqref{us:eq:error} and $\|y\|\le\|v_4\|$,
\[
 \xi=1+\|v_{q-1}\|^{-2} \Re(v_{q-1}^*e_{q-1}),\qquad
 \|e_{q-1}\|\le2\|A_<\widetilde R\|\,\|v_4\|,\qquad
 \tfrac12\le\xi\le\tfrac32.
\]
Combining these bounds with \eqref{us:eq:normalize} and the definition
of $L$ in \eqref{us:eq:L} yields
\begin{equation}\label{su:eq:normalized-generator-costs}
 \begin{gathered}
 |\xi-1|\le2\|v_{q-1}\|^{-1} \|A_<\widetilde R\|\,\|v_4\|,\qquad \|K\|_{\HS}^2\le32\|Pv_4\|^{-2} \|Qv_4\|^2,\\
 v_{\rm corr}=\xi^{-1}\bigl(e_{q-1}-(\xi-1)v_{q-1}\bigr),\qquad
 \|v_{\rm corr}\|\le8\|A_<\widetilde R\|\,\|v_4\|,\qquad
 \|L\|_{\HS}\le24\|v_{q-1}\|^{-1} \|A_<\widetilde R\|\,\|v_4\|.
 \end{gathered}
\end{equation}

\begin{lemma}[Bounds and equivariance of $K$]
\label{su:lem:chosen-generator}
In each regime of Lemma~\ref{up:lem:small-contraction} one can choose a
fixed length $q$ and a cutoff $\chi$ so that \eqref{su:eq:generator-admissibility}
holds and $K$ in \eqref{us:eq:normalize} has the following properties.
Pointwise bounds hold on $\supp\chi$. Expectations include the preparation blocks.

\smallskip\noindent
\textup{(i) Support and normalization.}
\begin{equation}\label{su:eq:chosen-generator-geometry}
 K^*=-K,\qquad K=PKP,\qquad
 \Re\bigl(v_{q-1}^*R_{[q-1,5]}Kv_4\bigr)=\|v_{q-1}\|^2.
\end{equation}
\smallskip\noindent
\textup{(ii) Source dependence and unitary symmetry.}
In the construction replace $v_4$ by a variable $\widehat v_4\ne0$,
keeping $E,g,R_4,P$ and the preparation indicator fixed. This gives a
family smooth near the support of the source cutoff $\widehat\chi$
constructed below, with
\begin{equation}\label{su:eq:generator-source-family}
 K=\widehat K(V_5,\ldots,V_{q-1};v_4;R_4).
\end{equation}
The family is invariant under $\widehat v_4\mapsto a\widehat v_4$ for
$a\in\C\setminus\{0\}$. For a unitary $\mathcal U$ acting identically
on $\ker P$, the transformation
\begin{equation}\label{su:eq:cut-generator-equivariance}
 V_5-g^2R_4\longmapsto\mathcal U\cdot(V_5-g^2R_4)\cdot\mathcal U^*,\quad
 V_i\longmapsto\mathcal U\cdot V_i\cdot\mathcal U^*,
       \quad6\le i\le q-1,\quad
 \widehat v_4\longmapsto\mathcal U\widehat v_4
\end{equation}
sends $\widehat K$ to $\mathcal U\widehat K\mathcal U^*$.

\smallskip\noindent
\textup{(iii) Size and matching error.}
We have
\begin{equation}\label{su:eq:pair-generator-size}
 \|K\|_{\HS}^2\le32,\qquad
 \E\|\chi K\|_{\HS}^2\le
 \begin{cases}
 Cg^2W^{\varepsilon/4},&\text{in regime (i)},\\
 CW^{-1+14\theta},&\text{in regime (ii)},
 \end{cases}
\end{equation}
and
\begin{equation}\label{su:eq:pair-matching-error}
 \|R_{[q-1,5]}Kv_4-v_{q-1}\|\le CW^{-4}\|v_{q-1}\|.
\end{equation}

\smallskip\noindent
\textup{(iv) Derivative bounds in $v_4$.}
Differentiate the vector argument with all matrices fixed, in the real
tangent space
\[
 \mathcal H_{v_4}\coloneqq\{w\in\ran P:\Re(v_4^*w)=0\}
 =\{Av_4:A^*=-A,\ A=PAP\}.
\]
There is a real subspace $\mathcal E_{v_4}\subseteq\mathcal H_{v_4}$ such that
\begin{equation}\label{su:eq:pair-source-regularity}
 \begin{aligned}
 \|Pv_4\|\cdot\|(w\cdot\nabla_{v_4})\widehat K\|_{\HS}
   &\le C\|w\|,&&w\in\mathcal H_{v_4},\\
 \|Pv_4\|\cdot\|(w\cdot\nabla_{v_4})\widehat K\|_{\HS}
   &\le C\|K\|_{\HS}\|w\|,
       &&w\in\mathcal H_{v_4}\cap\mathcal E_{v_4}^\perp,
 \end{aligned}
\end{equation}
where orthogonality is for the real Euclidean inner product and
\begin{equation}\label{su:eq:pair-exceptional-dimension}
 \dim_{\R}\mathcal E_{v_4}\le
 \begin{cases}
 CWg^2W^{\varepsilon/4},&\text{in regime (i)},\\
 CW^{14\theta},&\text{in regime (ii)}.
 \end{cases}
\end{equation}
Summing over an adapted real orthonormal basis $(e_a)$ of
$\mathcal H_{v_4}$ and using $\dim_{\R}\mathcal H_{v_4}\le2W$ gives
\begin{equation}\label{su:eq:pair-source-derivative-sum}
 \|Pv_4\|^2\cdot\sum_a
 \|(e_a\cdot\nabla_{v_4})\widehat K\|_{\HS}^2
 \le C\dim_{\R}\mathcal E_{v_4}+CW\|K\|_{\HS}^2.
\end{equation}
\smallskip\noindent
\textup{(v) Cutoff conditions.}
The cutoff family $\widehat\chi$ has the parameter dependence and
invariances in part (ii), with
\[
 0\le\widehat\chi\le1,\qquad
 \chi=\left.\widehat\chi\right|_{\widehat v_4=v_4},\qquad
 \E\chi\ge\tfrac12.
\]
On its support,
\[
 |R_{q-1}|_{\mathrm{op}}\le CW^3,\qquad
 g^2|R_4|_{\mathrm{op}}\le CW^{2\theta}\quad\text{in regime (ii)}.
\]
Its derivative satisfies
\begin{equation}\label{su:eq:pair-cutoff-derivative}
 \E\left[\|Pv_4\|^2\|K\|_{\HS}^2
       \|\nabla_{v_4}\widehat\chi\|^2\right]\le C.
\end{equation}
The gradient is taken along $\mathcal H_{v_4}$.
The local length, constants, and width threshold have the dependence
and uniformity specified in Lemma~\ref{up:lem:small-contraction}.
\end{lemma}

The proof is given in Section~\ref{su:sec:construction-estimates}.

\subsection{The Gaussian divergence bound}
\label{su:sec:generator-bounds}

\begin{lemma}[Gaussian divergence bound for the constructed field]
\label{us:lem:score}
Take $P,K,\chi$ from Lemma~\ref{su:lem:chosen-generator} and define
$L,X$ by \eqref{us:eq:L} and \eqref{us:eq:physical}. Conditional on
$V_1,\ldots,V_4$, write $\E$ for expectation over $V_5,\ldots,V_q$.
If $g^2\|PR_4P\|\le M_0$ on $\supp\chi$ for a deterministic $M_0>0$, then
\begin{equation}\label{us:eq:score}
 \begin{aligned}
  \E|\delta_{\rm G}(\chi X)|^2
 &\le C\E\left[\|Pv_4\|^2\sum_a
       \left\|(e_a\cdot\nabla_{v_4})(\widehat\chi\,\widehat K)\right\|_{\HS}^2\right]\\
 &\quad+CW(1+M_0^2)\E\|\chi K\|_{\HS}^2
       +Wg^4\E\|[\chi L,R_{q-1}]\|_{\HS}^2
       +CWM_0e^{-cW}.
 \end{aligned}
\end{equation}
Here $C,c>0$ are absolute, $(e_a)$ is a real orthonormal basis of
$\mathcal H_{v_4}$, and the derivative is evaluated at the prepared $v_4$.
\end{lemma}

\begin{proof}
To bound the divergence from blocks $5,\ldots,q-1$ in
\eqref{su:eq:active-matrix-derivative-bound}, expand
$\chi K=\sum_l a_l\mathsf K_l$ in a real Hilbert--Schmidt orthonormal
basis of $\{H:H^*=-H,\ H=PHP\}$. The matrix field multiplying $a_l$
in \eqref{us:eq:physical} differentiates along
\[
 Z_l\coloneqq[\mathsf K_l,V_5-g^2R_4]\cdot\nabla_{V_5}
       +\sum_{i=6}^{q-1}[\mathsf K_l,V_i]\cdot\nabla_{V_i}.
\]
The commutators have zero Euclidean divergence, and
$\langle V_i,[\mathsf K_l,V_i]\rangle=0$ for
$\langle A,B\rangle=\Re\Tr(A^*B)$. Thus only the shift at block $5$
contributes to the Gaussian adjoint,
\[
 Z_l^*=b_l-Z_l,\qquad
 b_l\coloneqq Wg^2\langle R_4,[\mathsf K_l,V_5-g^2R_4]\rangle.
\]
Writing $\delta_{\rm act}=\sum_l Z_l^*a_l$, two integrations by parts give
\begin{equation}\label{su:eq:active-divergence-decomposition}
 \E\delta_{\rm act}^2
 =\E\sum_{l,j}(Z_ja_l)(Z_la_j)
  +\E\sum_{l,j}a_la_jZ_lb_j
  +\E\sum_{l,j}a_l[Z_j,Z_l]a_j.
\end{equation}
By Cauchy--Schwarz and orthonormality, the first term is bounded in
absolute value by $\E\sum_l\|Z_l(\chi K)\|_{\HS}^2$. The second equals
\[
 Wg^2\E\Re\Tr\!\left(R_4[\chi K,[\chi K,V_5-g^2R_4]]\right).
\]
Using $PR_4=R_4P$, $\chi K=P(\chi K)P$, and
$g^2\|PR_4P\|\le M_0$, we obtain
\[
 \begin{aligned}
 \left|\E\sum_{l,j}a_la_jZ_lb_j\right|
 &\le CWM_0\E\bigl[(\|V_5\|+M_0)\|\chi K\|_{\HS}^2\bigr]\\
 &\le CW(1+M_0^2)\E\|\chi K\|_{\HS}^2+CWM_0e^{-cW}.
 \end{aligned}
\]
The last inequality uses the GUE norm tail and the uniform bound
on $\|\chi K\|_{\HS}$ in \eqref{su:eq:pair-generator-size}.
For the third term, summing in a matrix-unit basis gives
\begin{equation}\label{su:eq:commutator-basis-bound}
 \sum_l\|[\mathsf K_l,\chi K]\|_{\HS}^2
 \le2W\|\chi K\|_{\HS}^2.
\end{equation}
Since $[Z_j,Z_l]=Z_{[\mathsf K_l,\mathsf K_j]}$, with $Z_H$ defined
by replacing $\mathsf K_l$ by $H$, the invariance of the inner product gives
\[
 \begin{aligned}
 \left|\E\sum_{l,j}a_l[Z_j,Z_l]a_j\right|
 &=\left|\E\sum_l\langle[\mathsf K_l,\chi K],Z_l(\chi K)\rangle\right|\\
 &\le\tfrac12\E\sum_l\|Z_l(\chi K)\|_{\HS}^2
       +W\E\|\chi K\|_{\HS}^2.
 \end{aligned}
\]
Combining the three estimates yields
\begin{equation}\label{su:eq:active-matrix-derivative-bound}
 \E\delta_{\rm act}^2
 \le C\E\sum_l\|Z_l(\chi K)\|_{\HS}^2
       +CW(1+M_0^2)\E\|\chi K\|_{\HS}^2+CWM_0e^{-cW}.
\end{equation}

To bound the derivatives in \eqref{su:eq:active-matrix-derivative-bound},
differentiate \eqref{su:eq:cut-generator-equivariance} with
$\mathcal U=e^{t\mathsf K_l}$. By Lemma~\ref{su:lem:chosen-generator}(ii),(v),
\[
 Z_l(\chi K)=[\mathsf K_l,\chi K]
       -\bigl((\mathsf K_lv_4)\cdot\nabla_{v_4}\bigr)
          (\widehat\chi\,\widehat K).
\]

In coordinates on $\ran P$
with $Pv_4$ in the first coordinate direction, the real covariance of
$(\mathsf K_lv_4)_l$ has eigenvalue $\|Pv_4\|^2$ in the phase direction
$iPv_4$ and $\|Pv_4\|^2/2$ in the remaining tangent directions.
Applying this to each real coordinate of the derivative yields
\[
 \sum_l\left\|\bigl((\mathsf K_lv_4)\cdot\nabla_{v_4}\bigr)
            (\widehat\chi\,\widehat K)\right\|_{\HS}^2
 \le\|Pv_4\|^2\sum_a
       \left\|(e_a\cdot\nabla_{v_4})(\widehat\chi\,\widehat K)\right\|_{\HS}^2.
\]
Thus the chain rule and \eqref{su:eq:commutator-basis-bound} give
\begin{equation}\label{su:eq:matrix-source-derivative-comparison}
 \sum_l\|Z_l(\chi K)\|_{\HS}^2
 \le C\|Pv_4\|^2\sum_a
       \left\|(e_a\cdot\nabla_{v_4})(\widehat\chi\,\widehat K)\right\|_{\HS}^2
       +CW\|\chi K\|_{\HS}^2.
\end{equation}

Finally condition on $V_1,\ldots,V_{q-1}$, which fixes $\chi,K,L,R_{q-1}$.
The last component of \eqref{us:eq:physical} is
\[
 \chi X_q=[\chi K-\chi L,V_q]+g^2[\chi L,R_{q-1}].
\]
The commutator field has zero Gaussian divergence, while the constant
field contributes
\[
 \delta_{\rm last}\coloneqq Wg^2\langle V_q,[\chi L,R_{q-1}]\rangle,\qquad
 \delta_{\rm G}(\chi X)=\delta_{\rm act}+\delta_{\rm last}.
\]
The conditional GUE mean and variance give
\begin{equation}\label{su:eq:last-divergence-moments}
 \E[\delta_{\rm last}\mid V_1,\ldots,V_{q-1}]=0,\qquad
 \E[\delta_{\rm last}^2\mid V_1,\ldots,V_{q-1}]
       =Wg^4\|[\chi L,R_{q-1}]\|_{\HS}^2.
\end{equation}
Since $\delta_{\rm act}$ is fixed under this last conditional expectation,
the cross term vanishes. Combining
\eqref{su:eq:matrix-source-derivative-comparison},
\eqref{su:eq:active-matrix-derivative-bound}, and
\eqref{su:eq:last-divergence-moments} proves \eqref{us:eq:score}.
\end{proof}

\begin{proof}[Proof of Lemma~\ref{up:lem:small-contraction}]
For the field in \eqref{us:eq:physical}, \eqref{su:eq:chosen-generator-geometry}
and Lemma~\ref{su:lem:field-geometry} give \eqref{up:eq:angular-geometry}.
The first four components vanish, so integration by parts is conditional
on preparation. The smooth cutoff gives the required regularity in the
active blocks, and truncation is valid in $V_q$ since the field is affine
in $V_q$ with the other inputs fixed.

We now average \eqref{us:eq:score} over preparation. For its derivative
term the product rule and \eqref{su:eq:pair-source-derivative-sum} give
\begin{equation}\label{su:eq:chosen-generator-cut-costs}
 \begin{aligned}
 \|Pv_4\|^2\sum_a
 \left\|(e_a\cdot\nabla_{v_4})(\widehat\chi\,\widehat K)\right\|_{\HS}^2
 &\le C\chi^2\bigl(\dim_{\R}\mathcal E_{v_4}+W\|K\|_{\HS}^2\bigr)\\
 &\quad+C\|Pv_4\|^2\|K\|_{\HS}^2
                    \|\nabla_{v_4}\widehat\chi\|^2.
 \end{aligned}
\end{equation}
Taking expectations and using the size, dimension, and cutoff estimates
in \eqref{su:eq:pair-generator-size}, \eqref{su:eq:pair-exceptional-dimension},
and \eqref{su:eq:pair-cutoff-derivative} gives
\begin{equation}\label{su:eq:pair-source-energy}
 \E\left[\|Pv_4\|^2\sum_a
       \left\|(e_a\cdot\nabla_{v_4})(\widehat\chi\,\widehat K)\right\|_{\HS}^2\right]\le
 \begin{cases}
 CWg^2W^{\varepsilon/4},&\text{in regime (i)},\\
 CW^{14\theta},&\text{in regime (ii)}.
 \end{cases}
\end{equation}
Here $Wg^2\ge W^{1-2\beta}$ in (i) and $W^{14\theta}\ge1$ in (ii)
absorb the constant cutoff cost. For the last block, \eqref{us:eq:L}
and \eqref{su:eq:pair-matching-error} give $\|L\|_{\HS}\le CW^{-4}$
on $\supp\chi$. There $\|R_{q-1}\|\le CW^3$, so
\[
 Wg^4\E\|[\chi L,R_{q-1}]\|_{\HS}^2
 \le CWg^4W^6W^{-8}\le CW^{-1},
\]
since $g\le1$ in both parameter regimes for large $W$.

In (i), \eqref{su:eq:source-projection} permits $M_0=W^{\varepsilon/8}$.
Substituting \eqref{su:eq:pair-generator-size} and
\eqref{su:eq:pair-source-energy} into \eqref{us:eq:score} gives
\begin{equation}\label{su:eq:local}
 \E|\delta_{\rm G}(\chi X)|^2
 \le C\left(Wg^2W^{\varepsilon/4}
           +W^{1+\varepsilon/4}g^2W^{\varepsilon/4}+W^{-1}\right)
 \le CWg^2W^{\varepsilon/2}.
\end{equation}
In (ii), Lemma~\ref{su:lem:chosen-generator}(v) permits $M_0=CW^{2\theta}$,
and the same substitution gives
\begin{equation}\label{nc:eq:finalscore}
 \E|\delta_{\rm G}(\chi X)|^2
 \le C\left(W^{14\theta}+W^{1+4\theta}W^{-1+14\theta}+W^{-1}\right)
 \le CW^{18\theta}.
\end{equation}
In both cases $WM_0e^{-cW}=O(W^{-1})$. The cutoff mass
$\E\chi\ge1/2$, the fixed length $q$ and the uniformity of the bounds
follow from Lemma~\ref{su:lem:chosen-generator}.
\end{proof}

\subsection{Proof of Lemma~\ref*{su:lem:chosen-generator}}
\label{su:sec:construction-estimates}
\label{sec:critical-upper-proof}

\subsubsection{Step 1. Retained directions and overlap}
\label{su:sec:retained-geometry}
\label{su:sec:range-proof}
\label{su:sec:packet}
\label{su:sec:range}
\label{su:sec:field}

In the two regimes of Lemma~\ref{up:lem:small-contraction}, choose the
fixed length $q$ by \eqref{su:eq:packet2} and \eqref{nc:eq:parameters},
respectively, and set
\begin{equation}\label{su:eq:construction-threshold}
 t_{\mathrm{cut}}\coloneqq
 \begin{cases}
 W^{-40}g^{-q},&\text{in regime (i)},\\
 W^{-30+(1-2\theta)(q-6)/2},&\text{in regime (ii)}.
 \end{cases}
\end{equation}
The goal is to bound the rank of $Q$ in \eqref{su:eq:retained-projection}
and its overlap with the prepared vector $v_4$, as in
\eqref{su:eq:construction-overlap-goals}. We first establish the two
inverse-image estimates used in regimes (i) and (ii).

\begin{lemma}[Relative range]\label{us:lem:range}
Fix $\kappa_0\in(0,2)$ and $\delta>0$. Let $|E|\le2-\kappa_0$ and
let $V$ be normalized complex GUE of size $W$. The data
$E,B_0,B_1,Q_0,U,v$ below are deterministic.
Let $B_0=B_{\rm bad}+B_{\rm small}$ be Hermitian, with a spectral projection $Q_0$ of $B_0$
and $B_{\rm bad}\coloneqq Q_0B_0Q_0$, with
$\norm{B_{\rm small}}=o(1)$. Let $B_1$ be Hermitian of norm $o(1)$,
$U$ a fixed subspace, and $v$ a fixed vector. Set
\[
 R\coloneqq (V-EI-B_0-B_1)^{-1}.
\]
If $W^\delta\le d=o(W)$, $d\ge\dim U+2\rank Q_0+1$, and
$\norm{Q_0P_U}\le C\sqrt{d/W}$, then, with failure
$W^{-D}+Ce^{-cd}$ for each prescribed $D$,
\begin{equation}
 \norm{P_{RU}v}^2\le C(d/W)(\log W)^2 v^*(I+B_0^2)^{-1}v.
 \label{us:eq:relative}
\end{equation}
Here $P_{RU}$ is the orthogonal projection onto $RU$, and the
probability event may depend on the fixed vector $v$.
The constants do not depend on $\norm{B_{\rm bad}}$.
\end{lemma}

In Lemma~\ref{us:lem:range}, the fresh cross block is used both to
form the inverse-image subspace and to test its overlap. We first
prove the Gaussian estimate that controls this dependence.

\begin{lemma}[Dependent Gaussian projection]\label{us:lem:self}
Fix $\delta>0$ and constants $C,C_1>0$. Let $d\ge W^\delta$, $d=o(W)$, $\eta\coloneqq Ad/W$ with fixed sufficiently
large $A$, and let $R$ be a fixed Hermitian inverse on a space of dimension
at most $W$. Suppose
\begin{equation}
 C_0d\le\rank\one_{|R|\ge\eta^{-1}}\le C_1W\eta,\quad
 \Tr|R_<|\le CW\log W,\quad
 \Tr R_<^2\le CW/\eta,\quad \norm{R_<}\le\eta^{-1},
 \label{us:eq:spectrum}
\end{equation}
where $R_<\coloneqq R\one_{|R|<\eta^{-1}}$ and $C_0$ is sufficiently large.
Let $F$ be a deterministic subspace of the space on which $R$ acts,
with dimension at most $d$, and let $Y$ have independent centered complex
Gaussian entries of variance $1/W$, one row per dimension of this space,
and at most $d$ columns. All projections in this lemma are orthogonal.
Then, with failure $Ce^{-cd}$,
\begin{equation}
 \norm{Y^*P_{F+\ran(RY)}}\le C\sqrt\eta\log W.
 \label{us:eq:self}
\end{equation}
\end{lemma}
\begin{proof}
The dependence appears in the quadratic form $Y^*(I-P_F)RY$. We split the
spectrum so that the large inverse directions provide a lower bound
for the generated vectors, while the remaining quadratic form has only
a logarithmic cost. Set
\[
 Y_\perp\coloneqq(I-P_F)RY,\qquad
 E_>\coloneqq\ran\one_{\{|R|\ge\eta^{-1}\}}.
\]
Set $H=E_>\cap F^\perp$, so $\dim H\ge(C_0-1)d$ and
$P_HY_\perp=P_HRY$. Each column has covariance
$W^{-1}P_HR^2P_H\ge(W\eta^2)^{-1}P_H$.
The rectangular Gaussian lower bound therefore gives
\begin{equation}\label{us:eq:lowerK}
 \norm{Y_\perp b}\ge c\eta^{-1/2}\norm b.
\end{equation}
The Gaussian norm and quadratic-form bounds also give
\begin{equation}\label{us:eq:Gaussnorm}
 \norm{Y^*P_F}+\norm{Y^*P_{(I-P_F)E_>}}\le C\sqrt\eta,
 \qquad \norm{(I-P_F)R_<Y}\le C\eta^{-1/2},
\end{equation}
and
\begin{equation}\label{us:eq:quadratic}
 \norm{Y^*(I-P_F)R_<Y}\le C\log W.
\end{equation}
These three estimates hold simultaneously with failure $Ce^{-cd}$;
their elementary verifications are in
Appendix~\ref{app:gaussian-projection-inputs}.

Write $Y_\perp=Y_{\perp,>}+Y_{\perp,<}$ according to the split of $R$. Equations
\eqref{us:eq:lowerK}--\eqref{us:eq:quadratic} imply
\[
 \norm{Y^*Y_\perp b}
 \le C\sqrt\eta\norm{Y_{\perp,>}b}+C\log W\norm b
 \le C\sqrt\eta\log W\norm{Y_\perp b}.
\]
Here $\norm{Y_{\perp,<}b}\le C\eta^{-1/2}\norm b$ is absorbed using
\eqref{us:eq:lowerK}. Since $F$ and $\ran Y_\perp$ are orthogonal and span
$F+\ran(RY)$, \eqref{us:eq:self} follows.
\end{proof}

\begin{proof}[Proof of Lemma~\ref{us:lem:range}]
Write $\Gamma=R^{-1}=V-EI-B_0-B_1$.
Choose a $d$-dimensional space $Q$ containing $\ran Q_0,U,\ran(B_1Q_0)$ and
$(I-Q_0)v$. Decompose $Q=\ran Q_0\oplus Q_1$ and then $Q\oplus Q^\perp$.
The complementary matrix $\Gamma_c$ has only a small-norm deformation.
For every prescribed $D>0$, its spectral input is
\begin{equation}\label{us:eq:complement-spectrum}
 \left.\eqref{us:eq:spectrum}\right|_{R=\Gamma_c^{-1},\,\eta=Ad/W}
 \quad\text{with probability at least }1-W^{-D}.
\end{equation}
Here $A$ is a sufficiently large fixed constant. The counting-law and
dyadic-summation proof is given under \textbf{Proof of
\eqref{us:eq:complement-spectrum}} in
Appendix~\ref{app:inverse-image-spectrum}. Write the cross block as
$(Y_0,Y_1)$ according to $\ran Q_0,Q_1$. Because $\ran(B_1Q_0)\subset Q$,
$Y_0$ has zero mean and is independent of both $\Gamma_c$ and $Y_1$.
The latter block can have a deterministic mean. Condition on these two
objects and put $F\coloneqq \ran(\Gamma_c^{-1}Y_1)$. Lemma \ref{us:lem:self} gives
\begin{equation}
 \norm{Y_0^*P_{F+\ran(\Gamma_c^{-1}Y_0)}}
 \le C\sqrt\eta\log W,\qquad \eta=Ad/W.
 \label{us:eq:Y0}
\end{equation}

There is also the following uniform inverse-image compression:
\begin{equation}
 \norm{x_Q}\le C\sqrt\eta\norm x
 \quad\text{whenever }\Gamma x\in Q.
 \label{us:eq:absolute}
\end{equation}
Indeed $x_c=-\Gamma_c^{-1}Yx_Q$. Let $\mathcal E$ be the spectral subspace
of $\Gamma_c$ for $[-\eta,\eta]$, of dimension at least $C_0d$.
Conditional on $\Gamma_c$, write $\overline Y\coloneqq \E[Y\mid\Gamma_c]$ and let
$\mathcal N\subset\mathcal E$ be the kernel of
$(P_{\mathcal E}\overline Y)^*$ restricted to $\mathcal E$.
Then $\dim\mathcal N\ge\dim\mathcal E-d$. For every column vector $u$,
the inverse spectral bound is applied before the row projection:
\[
 \norm{\Gamma_c^{-1}Yu}
 \ge\eta^{-1}\norm{P_{\mathcal E}Yu}
 \ge\eta^{-1}\norm{P_{\mathcal N}P_{\mathcal E}Yu}.
\]
The remaining centered Gaussian block $P_{\mathcal N}P_{\mathcal E}Y$
has least singular value at least $c\sqrt{d/W}$, in orthonormal coordinates
on $\mathcal N$. Taking $u=x_Q$ gives
$\norm{x_c}\ge c\eta^{-1}\sqrt{d/W}\norm{x_Q}$.
This event and \eqref{us:eq:Y0} can be intersected without any independence
assertion between them.

For $x\in RU$, let $u\coloneqq \Gamma x\in U$. Transversality and the usual GUE
norm bound imply
\[
 \norm u\le2\norm{(I-Q_0)u}
 \le2\norm{V-EI-B_{\rm small}-B_1}\norm x\le C\norm x,
 \qquad \norm{Q_0u}\le C\sqrt\eta\norm x.
\]
Writing $x_0=Q_0x$, the $Q_0$ component of $\Gamma x=u$ is
\[
 B_{\rm bad}x_0
 =Q_0(V-EI-B_{\rm small}-B_1)x_Q+Y_0^*x_c-Q_0u.
\]
The deterministic $Q_0$-to-complement block vanishes by the choice of $Q$;
its Gaussian part is exactly $Y_0^*$.
Since $x_c=-\Gamma_c^{-1}Yx_Q\in F+\ran(\Gamma_c^{-1}Y_0)$,
\eqref{us:eq:Y0} bounds $\|Y_0^*x_c\|$ by
$C\sqrt\eta\log W\|x_c\|$. The bounded norm of
$V-EI-B_{\rm small}-B_1$, \eqref{us:eq:absolute}, and the bound on $Q_0u$
above now give
\begin{equation}
 \norm{B_{\rm bad}x_0}+\norm{x_Q}
 \le C\sqrt\eta\log W\norm x.
 \label{us:eq:weightedcompression}
\end{equation}
Finally split $v$ into its good
and bad components and use weighted Cauchy--Schwarz:
\[
 \begin{aligned}
 |v^*x|&\le \norm{v_{\rm good}}\norm{x_Q}
 +\norm{(I+B_{\rm bad}^2)^{-1/2}v_{\rm bad}}
  \norm{(I+B_{\rm bad}^2)^{1/2}x_0}\\
 &\le C\sqrt\eta\log W\sqrt{v^*(I+B_0^2)^{-1}v}\norm x.
 \end{aligned}
\]
Take the supremum over unit vectors in $RU$.
\end{proof}

When the deformation profile is not needed, the compression step in
Lemma~\ref{us:lem:range} gives the following consequence.

\begin{corollary}[Fresh-matrix compression]\label{nc:lem:compression}
Fix $\kappa_0\in(0,2)$ and $\delta>0$. All matrices and vectors below
act on $\C^W$.
Let $B=B_{\rm small}+B_{\rm bad}$ be deterministic Hermitian,
$\norm{B_{\rm small}}=o(1)$ and $\rank B_{\rm bad}\le r_0$.
Let $U$ be a deterministic subspace of dimension at most $m$, and let
$v$ be a deterministic unit vector. Put $R\coloneqq (V-EI-B)^{-1}$, where
$V$ is fresh normalized complex GUE and $|E|\le2-\kappa_0$.
If $W^\delta\le d=o(W)$ and $d\ge m+r_0+1$, then, for each fixed $D'>0$,
\begin{equation}
 \norm{P_{RU}v}^2\le Cd/W
 \label{nc:eq:compression}
\end{equation}
with failure at most $W^{-D'}+Ce^{-cd}$, where $P_{RU}$ is the
orthogonal projection onto $RU$. The constants are independent
of $\norm{B_{\rm bad}}$ and of any relative orientation of $U,v$ and
$\ran B_{\rm bad}$.
\end{corollary}
\begin{proof}
Choose a deterministic \(d\)-dimensional space \(\mathcal Q\)
containing \(U,v\), and \(\operatorname{ran}B_{\rm bad}\).
The complementary block of \(V-EI-B\) has only a small-norm deformation.
The argument proving \eqref{us:eq:absolute} applies without the
transversality hypothesis: at \(\eta=Ad/W\), its near-zero spectral
space has dimension at least a large constant times \(d\).
After removing the deterministic cross-block mean inside that space,
the remaining Gaussian matrix has least singular value at least
\(c\sqrt{d/W}\). The complementary row equation consequently gives
\(\|x_{\mathcal Q}\|\le C\sqrt{d/W}\|x\|\) for every \(x\in RU\).
Since \(v\in\mathcal Q\), take the supremum of \(|v^*x|\) over unit
vectors in \(RU\). The local-law and Gaussian failures are respectively
\(W^{-D'}\) and \(Ce^{-cd}\), as in the proof of \eqref{us:eq:absolute}.
\end{proof}

\paragraph{Rank and overlap for the retained space.}
On the event in Lemma~\ref{su:lem:preparation}, fix $V_1,\ldots,V_4$.
A singular-value moment bounds the retained dimension; the fresh matrix
$V_5$ then controls its overlap with $v_4$. The cutoff in
Section~\ref{su:sec:cutoff-estimates} will enforce the resulting good events.

The right product $A$ is the adjoint of a forward product in reversed
order with initial state zero. The product moment bound \eqref{us:eq:product} in
Appendix~\ref{app:product-moment-proof}, with $h=q-6$, gives
\begin{equation}\label{su:eq:retained-rank-moment}
 \rank Q=\rank A_>\le t_{\mathrm{cut}}^{-p}\Tr|A|^p
 \le t_{\mathrm{cut}}^{-p}W^{1+e},\qquad 0<p<\frac2{q-5},
\end{equation}
on an invariant event of probability $1-o(1)$, for each fixed $e>0$.

\paragraph{Regime (i).}
For $g=W^{-a}$ with $a\in[\alpha,\beta]\subset(0,1/2)$, choose
\begin{equation}\label{su:eq:packet2}
 q=5+\left\lceil25600\varepsilon^{-1} (40+\beta)\right\rceil,\qquad
 p=\frac2{q-5}(1-\varepsilon/6400).
\end{equation}
Using the threshold \eqref{su:eq:construction-threshold} and
$e=\varepsilon/3200$ in \eqref{su:eq:retained-rank-moment} yields
\begin{equation}\label{su:eq:rank}
 \rank Q\le W^{1+\varepsilon/3200+40p}g^{qp}
 <Wg^2W^{\varepsilon/800}.
\end{equation}
The last inequality follows by substitution of \eqref{su:eq:packet2},
uniformly in $a$; in particular $q,p$ do not depend on $a,W$.

\smallskip\noindent
\emph{Small overlap and compatibility with $P$.}
To control $\|Qv_4\|/\|Pv_4\|$, apply Lemma~\ref{us:lem:range} with
incoming deformation $g^2R_4$ and spectral projection
\[
 Q_0\coloneqq\mathbf1_{\{|g^2R_4|>W^{-\varepsilon/800}\}}.
\]
The rank of $Q_0$ is bounded by \eqref{su:eq:prepared-large-rank}. We first
check the required transversality and remove the large right
self-energy from the deformation.

Use the retained space $\ran(A_>^*)$ from \eqref{su:eq:rank}, and split $g^2 \overleftarrow R_6=B_{1,>}+B_{1,<}$ at $W^{-\varepsilon/800}$.
Define
\[
 U\coloneqq \ran(A_>^*)+\ran B_{1,>}.
\]
The joint right-cavity law is invariant under
common unitary conjugation, independently of the incoming data. Therefore
$U$ is Haar conditional on its dimension and invariant events.
The Gaussian representation of a Haar subspace gives
\begin{equation}
 \norm{Q_0P_U}\le C\sqrt{(\dim U+\rank Q_0+W^{(1-2\beta)/4})/W}
 \label{us:eq:Haar}
\end{equation}
with exponential failure in the numerator's dimension scale.
This uses only common orbital invariance, not independent cavity matrices.
Since $\ran Q=\widetilde R\ran(A_>^*)$, the enlarged space satisfies
\begin{equation}
 \ran Q\subset\widetilde RU
  =(V_5-EI-g^2R_4-B_{1,<})^{-1}U.
 \label{us:eq:erase}
\end{equation}
Both inverse images are the set of $x$ for which the corresponding
matrix times $x$ belongs to $U$; their difference has range in $U$.
We apply the range estimates to this enlarged space and then restrict
them to $\ran Q$.
For each new block, condition on its Hermitian shift and write
$\Gamma=V-EI_W-g^2R$ with $V$ normalized complex GUE.
For an interval $J$ and $\delta>0$, the single-block Gaussian regularization
estimates of \cite[Theorem~1]{APSSS} (see also \cite[Theorem~2, (1.10)]{PSSS}) give
\begin{equation}
 \E N(\Gamma,J)\le CW|J|,\qquad
 \Prob(\norm{\Gamma^{-1}}>W^{1+\delta})\le CW^{-\delta}.
 \label{us:eq:wegner}
\end{equation}
Applying the count bound to the right recursion gives
$\rank B_{1,>}\le Wg^2 W^{\varepsilon/400}$ with failure $O(W^{-\varepsilon/800})$.
For the dimension parameter in Lemma~\ref{us:lem:range}, with
$\delta=(1-2\beta)/4$, choose
\begin{equation}\label{su:eq:q}
 d\coloneqq \lceil8Wg^2 W^{\varepsilon/200}\rceil.
\end{equation}
Equations \eqref{su:eq:rank} and \eqref{su:eq:prepared-large-rank} give
$d\ge\dim U+2\rank Q_0+1$, while
\[
 d\ge W^{1-2\beta},\qquad d/W\le CW^{-2\alpha+\varepsilon/200}=o(1).
\]
Thus $d$ dominates the numerator in \eqref{us:eq:Haar}, and the spectral
scale $d/W$ is mesoscopic. Also $\ran(I-P)\subset\ran Q_0$.
Use the compression bound \eqref{us:eq:absolute} in the proof of
Lemma~\ref{us:lem:range}, and apply \eqref{us:eq:relative} separately
to the two fixed vectors $v_4,Pv_4$. Put
$b_4=v_4^*(I+g^4R_4^2)^{-1}v_4$; preparation in \eqref{su:eq:prep}
gives $b_4\le C\|Pv_4\|^2$. For the projection $Q$ in
\eqref{su:eq:retained-projection}, this gives
\begin{equation}\label{su:eq:geometry}
 \norm{(I-P)Q}^2\le Cd/W=o(1),\qquad
 \norm{Qv_4}^2\le g^2 W^{\varepsilon/100}b_4,\qquad
 \norm{QPv_4}^2\le g^2 W^{\varepsilon/100}\|Pv_4\|^2.
\end{equation}
The fixed constants and $(\log W)^2$ are absorbed into $W^{\varepsilon/200}$
with a common threshold.

\paragraph{Regime (ii).}
Here $P=I_W$, and preparation gives \eqref{nc:eq:prepared}.
To control $\|R_{[q-1,5]}v_4\|$, the physical scale $g^2$ may be too small for
the local law. Introduce the larger spectral scale
\begin{equation}\label{nc:eq:buffer-scale}
 \upsilon=W^{-1+2\theta},\qquad
 c_*W^{-1/2}\le g\le W^{-1/2+\theta/2}.
\end{equation}
Then $g^2\le\upsilon$ and $W\upsilon=W^{2\theta}$, so
Lemma~\ref{us:lem:local-law} applies at $\eta=\upsilon$.
We use this scale in Step~2 to bound $\|v_{q-1}\|$ from below.

Choose
\begin{equation}\label{nc:eq:parameters}
 q\coloneqq5+\lceil10/\theta\rceil,\qquad
 p\coloneqq\frac2{q-5}(1-\theta/4).
\end{equation}
Using the threshold \eqref{su:eq:construction-threshold} and
$e=\theta/4$ in \eqref{su:eq:retained-rank-moment} gives
\begin{equation}\label{nc:eq:dimension}
 \rank Q\le W^{1+\theta/4+30p}\upsilon^{(q-6)p/2}\le W^{9\theta},
\end{equation}
by \eqref{nc:eq:parameters}. Again $q$ is fixed independently of $g,W$.

\smallskip\noindent
\emph{Small overlap.}
Since $P=I_W$, the absolute compression bound in
Corollary~\ref{nc:lem:compression} suffices for the overlap.
Split each of $g^2 R_4$ and $g^2 \overleftarrow R_6$ spectrally at $W^{-\theta}$.
The two large pieces have ranks at most $W^{3\theta}$ with failure
$O(W^{-\theta})$: their expected counts are at most
$CWg^2 W^\theta\le CW^{2\theta}$. The small pieces have total norm
at most $2W^{-\theta}$. Conditional on all right-cavity data, apply
Corollary~\ref{nc:lem:compression} to
\[
 B=g^2 R_4+g^2 \overleftarrow R_6,\qquad U=\ran A_>^*,\qquad
 d\coloneqq \lceil W^{12\theta}\rceil.
\]
Indeed $d\ge W^{9\theta}+2W^{3\theta}+1$ for large $W$ and $d=o(W)$.
The source vector $v_4$ and all these data are fixed before $V_5$.
Applying Corollary~\ref{nc:lem:compression} to $v_4/\|v_4\|$ and using
$\ran Q=\widetilde R U$ yields
\begin{equation}
 \|v_4\|^{-2} \|Qv_4\|^2\le CW^{-1+12\theta}
                 \le W^{-1+13\theta}
 \label{nc:eq:geometry}
\end{equation}
with probability $1-o(1)$, uniformly in \eqref{eq:critical-window}.

\paragraph{Common conclusion.}
In both regimes the overlap estimates give, on $\ran Q$,
\[
 B_Q\ge\left(1-\|(I-P)Q\|^2-
       \|Pv_4\|^{-2} \|QPv_4\|^2\right)I_{\ran Q}
 \ge\tfrac12 I_{\ran Q},\qquad
 \|Pv_4\|^{-2} \|Qv_4\|^2=o(1).
\]
Here the first loss vanishes when $P=I_W$. By \eqref{us:eq:sourceenergy},
$\|y\|=o(\|v_4\|)$, so the first three conditions in
\eqref{su:eq:generator-admissibility} hold with room for the cutoff
margins. Combining \eqref{su:eq:rank}, \eqref{su:eq:geometry},
\eqref{nc:eq:dimension} and \eqref{nc:eq:geometry} gives, for any
$0\le\chi\le1$ supported on these events,
\begin{equation}\label{su:eq:construction-overlap-goals}
 \begin{aligned}
 \rank Q&\le
 \begin{cases}CWg^2W^{\varepsilon/4},&\text{in (i)},\\
 CW^{14\theta},&\text{in (ii)},\end{cases}
 &\E\!\left[\chi^2\|Pv_4\|^{-2} \|Qv_4\|^2\right]&\le
 \begin{cases}Cg^2W^{\varepsilon/4},&\text{in (i)},\\
 CW^{-1+14\theta},&\text{in (ii)}.\end{cases}
 \end{aligned}
\end{equation}
The rank bound is pointwise on the support, where $\|Qv_4\|\le\|Pv_4\|$.
The expectation includes the preparation blocks.

\subsubsection{Step 2. The matching error}
\label{su:sec:error-derivative-estimates}
\label{su:sec:discarded-output}

We prove the relative error bound \eqref{su:eq:construction-error-goal}.
Since $\|A_<\|\le t_{\mathrm{cut}}$, the numerator is at most
$t_{\mathrm{cut}}\|\widetilde R\|\,\|v_4\|$. We therefore need an upper bound on
$\|\widetilde R\|$ and a lower bound on the transported vector
$\|v_{q-1}\|$. The inverse bound is common to both regimes.

The active-chain matrix has diagonal blocks
$V_5-EI-g^2R_4,V_6-EI,\ldots,V_{q-1}-EI$ and adjacent off-diagonal
blocks $gI$. Impose the bound $W^3$ on its inverse. It has
failure $C(q-5)^{3/2}/W$, uniformly over the prepared boundary: the common
scalar Gaussian shift has variance $1/((q-5)W^2)$, independently of the
remaining coordinates, and a union bound covers the $(q-5)W$ eigenvalues.
The exact Schur identities give
\begin{equation}\label{su:eq:active-chain-inverse}
 \norm{\widetilde R},\norm{R_{q-1}}\le W^3,\qquad
 \norm{R_{[q-1,5]}}\le W^3g^{-(q-6)}.
\end{equation}
For the output cutoff, set
\begin{equation}\label{su:eq:output-cutoff-scale}
 v_{\min}\coloneqq
 \begin{cases}
 g^{-(q-5)}\|Pv_4\|W^{-3\varepsilon/800},&\text{in (i)},\\
 \upsilon^{-(q-5)/2}\|v_4\|W^{-4\theta},&\text{in (ii)}.
 \end{cases}
\end{equation}
For every fixed $D>0$ and sufficiently large $W$, the probability that the
chain is regular and $\|v_{q-1}\|<2v_{\min}$ is at most $W^{-D}$.

\paragraph{Regime (i).}
To bound the denominator $\|v_{q-1}\|$ from below, define
\[
 F_j\coloneqq(I+g^4R_j^2)^{-1},\qquad4\le j\le q-1.
\]
Write $R_{[j,5]}=R_j\cdots R_5$ for $5\le j\le q-1$.
For a regular past, $u$ fixed before $V_j$, and $\eta=g^2$, the
relative local law \eqref{us:eq:relative-local-law} in
Appendix~\ref{app:relative-local-law-proof} gives
\[
 \Im\,u^*(R_j^{-1}-ig^2 I_W)^{-1}u
 \asymp u^*(I_W+g^4R_{j-1}^2)^{-1}u
\]
with failure probability at most $W^{-D}$ for every fixed $D>0$ and all
sufficiently large $W$. The probability that the chain is regular and the
following comparison fails for some $5\le j\le q-1$ is also at most $W^{-D}$:
\begin{equation}
 v_4^*R_{[j,5]}^*F_jR_{[j,5]}v_4\asymp_{j-4}g^{-2(j-4)}v_4^*F_4v_4,
 \qquad 5\le j\le q-1.
 \label{us:eq:profile}
\end{equation}
Indeed, $R_jF_jR_j=g^{-2}\Im(R_j^{-1}-ig^2)^{-1}$, and the Dyson
profile at the next step is comparable to $F_{j-1}$. Normalize the
conditional vector at each step and sum the failure probabilities over
the fixed number of steps, taking the one-step exponent larger than $D$.
Regularity itself has probability $1-o(1)$.

Since $b_4\ge(1+W^{\varepsilon/400})^{-1}\|Pv_4\|^2$ by the definition of $P$,
\eqref{us:eq:profile} gives
\[
 \norm{v_{q-1}}\ge c_{q-5}g^{-(q-5)}\|Pv_4\|W^{-\varepsilon/800}.
\]
The comparison constants are uniform in $a$ and $q$ is fixed, so this
exceeds $2v_{\min}$ for large $W$.

\paragraph{Regime (ii).}
Condition on the prepared $R_4,v_4$. By \eqref{nc:eq:prepared},
$\|R_4\|\le W^{1+\theta}$ and $R_4$ is regular on the preparation event.
For this fixed $v_4$, the lower bound proved in
Appendix~\ref{app:buffered-output} is
\begin{equation}
 \norm{R_{[q-1,5]}v_4}\ge c_{q-5}\upsilon^{-(q-5)/2}\|v_4\|W^{-3\theta}.
 \label{nc:eq:profile}
\end{equation}
The probability that this bound fails and the entire chain is regular
is at most $W^{-D'}$ for every fixed $D'>0$. It gives
$\|v_{q-1}\|\ge2v_{\min}$ since $c_{q-5}W^\theta>2$ for large $W$.
The estimate uses forward conditioning, whereas
\eqref{nc:eq:geometry} uses right-cavity conditioning; we combine their
events by a union bound, without conditioning either estimate on the
other's data.

\paragraph{Common error bounds.}
On $\|v_{q-1}\|\ge v_{\min}$, substitute
\eqref{su:eq:construction-threshold} into
$t_{\mathrm{cut}}W^3\|v_4\|/v_{\min}$ and, in regime (i), use
$\|v_4\|\le W^{10}\|Pv_4\|$ from \eqref{su:eq:prep}. This gives
\begin{equation}\label{su:eq:error}
 \|v_{q-1}\|^{-1} \|A_<\widetilde R\|\,\|v_4\|\le
 \begin{cases}
 g^{-5}W^{-27+3\varepsilon/800},&\text{in (i)},\\
 \upsilon^{1/2}W^{-27+4\theta},&\text{in (ii)}.
 \end{cases}
\end{equation}
Using $g^{-5}\le W^{5\beta}$ and $\beta<1/2$ in (i), we obtain in both regimes
\begin{equation}\label{su:eq:construction-error-goal}
 \|v_{q-1}\|^{-1} \|A_<\widetilde R\|\,\|v_4\|\le CW^{-4}.
\end{equation}
This verifies the last condition in \eqref{su:eq:generator-admissibility}
and controls $|\xi-1|$ by \eqref{su:eq:normalized-generator-costs}.
For the normalization and cutoff derivatives, \eqref{su:eq:active-chain-inverse}
gives the common ratio bound
\begin{equation}\label{su:eq:output-derivative-ratio}
 v_{\min}^{-2} \|Pv_4\|^2\|R_{[q-1,5]}\|^2\le
 \begin{cases}
 g^2W^{6+3\varepsilon/400},&\text{in (i)},\\
 C_{c_*,q}W^{5+8\theta+2\theta(q-5)}\le C_{c_*,q}W^{25+10\theta},
       &\text{in (ii)}.
 \end{cases}
\end{equation}
In the second line we used $g\ge c_*W^{-1/2}$ and
$2\theta(q-5)\le20+2\theta$. All restrictions above have common
probability $1-o(1)$ and will be imposed by the cutoff.
In particular, \eqref{su:eq:normalized-generator-costs} gives the matching
estimate \eqref{su:eq:pair-matching-error} and, together with Step~1,
both size estimates in \eqref{su:eq:pair-generator-size}. The real
normalization \eqref{us:eq:normalize} gives \eqref{su:eq:chosen-generator-geometry}.

\subsubsection{Step 3. Derivatives in \texorpdfstring{$v_4$}{v4}}
\label{su:sec:generator-estimates}

Keep the matrices fixed. Both thresholds in
\eqref{su:eq:construction-threshold} are independent of $v_4$, so
$P,Q$ are fixed as $v_4$ varies. We prove a common directional estimate,
then use the error bounds from Section~\ref{su:sec:discarded-output}
to control the scalar normalization.

We first isolate the directions in which the derivative can be large. Define
\[
 u\coloneqq QB_Q^{-1}Qv_4.
\]
Since $u\in\ran Q$ and $y=\Pi u$, we have
\[
 \|y\|^2=u^*B_Qu\ge\tfrac12\|u\|^2.
\]
For $w\in\mathcal H_{v_4}$, the directional derivatives of $P,Q$ and
$\|Pv_4\|$ vanish. Differentiation therefore yields
\[
 \begin{aligned}
  (w\cdot\nabla_{v_4})\Pi
   &=-\|Pv_4\|^{-2}\bigl(w\cdot(Pv_4)^*+Pv_4\cdot w^*\bigr),\\
  (w\cdot\nabla_{v_4})y
   &=(I-\Pi QB_Q^{-1}Q)\cdot\bigl((w\cdot\nabla_{v_4})\Pi\bigr)\cdot u
      +\Pi QB_Q^{-1}Qw.
 \end{aligned}
\]
Using $\|B_Q^{-1}\|\le2$ and $\|u\|\le\sqrt2\|y\|$ gives
\begin{equation}\label{su:eq:matching-source-derivative}
 \|(w\cdot\nabla_{v_4})y\|
 \le C\|Qw\|+C\|Pv_4\|^{-1} \|y\|\|w\|.
\end{equation}
For the source derivative of the unnormalized generator, define
\[
 \widehat K_0(V_5,\ldots,V_{q-1};\widehat v_4;R_4)
 \coloneqq\left.K_0\right|_{v_4=\widehat v_4}.
\]
The product rule in \eqref{us:eq:H0}, together with
\eqref{su:eq:matching-source-derivative} and \eqref{us:eq:sourceenergy}, gives
\begin{equation}\label{su:eq:unnormalized-source-response}
 \|Pv_4\|\,\|(w\cdot\nabla_{v_4})\widehat K_0\|_{\HS}
 \le C\|Qw\|+C\|K_0\|_{\HS}\|w\|.
\end{equation}
In particular, $\|Qv_4\|\le\|Pv_4\|$ implies
\begin{equation}\label{su:eq:matching-error-derivative}
 \|(w\cdot\nabla_{v_4})y-w\|\le C\|w\|.
\end{equation}

It remains to control differentiation of the scalar normalization.
With the matrices fixed,
\[
 (w\cdot\nabla_{v_4})v_{q-1}=R_{[q-1,5]}w,\qquad
 (w\cdot\nabla_{v_4})e_{q-1}=A_<\widetilde R\bigl((w\cdot\nabla_{v_4})y-w\bigr).
\]
Apply the quotient rule to
$\xi=1+\|v_{q-1}\|^{-2}\Re(v_{q-1}^*e_{q-1})$.
Using \eqref{su:eq:matching-error-derivative} and
$\|e_{q-1}\|\le2\|A_<\widetilde R\|\,\|v_4\|$, then taking the supremum
over unit $w\in\mathcal H_{v_4}$, gives
\begin{equation}\label{su:eq:normalization-source-response}
 \|Pv_4\|^2\|\nabla_{v_4}\xi\|^2
 \le C\|v_{q-1}\|^{-2} \|A_<\widetilde R\|^2\|v_4\|^2\left(1+
       \|v_{q-1}\|^{-2} \|Pv_4\|^2\|R_{[q-1,5]}\|^2\right).
\end{equation}
Substituting \eqref{su:eq:error} and
\eqref{su:eq:output-derivative-ratio} into
\eqref{su:eq:normalization-source-response} yields
\begin{equation}\label{su:eq:normalization-derivative-bound}
 \|Pv_4\|^2\|\nabla_{v_4}\xi\|^2\le
 \begin{cases}
 C\bigl(g^{-10}W^{-54+3\varepsilon/400}+g^{-8}W^{-48+3\varepsilon/200}\bigr),
     &\text{in (i)},\\
 CW^{-30+20\theta},&\text{in (ii)}.
 \end{cases}
\end{equation}
These are uniformly bounded: in (i), use $g\ge W^{-\beta}$ and
$\beta<1/2$, and in (ii), use the stipulated small $\theta$.
With the normalization derivative controlled in both regimes, define
the exceptional real subspace by
\[
 \mathcal E_{v_4}\coloneqq
 \bigl\{w\in\mathcal H_{v_4}:Qw=0,\ (w\cdot\nabla_{v_4})\xi=0\bigr\}^{\perp_{\mathcal H_{v_4}}}.
\]
Thus the preceding bound on $\nabla_{v_4}\xi$ and the real codimension
of this common kernel give
\begin{equation}\label{su:eq:construction-derivative-goals}
 \|Pv_4\|^2\|\nabla_{v_4}\xi\|^2\le C,\qquad
 \dim_{\R}\mathcal E_{v_4}\le2\rank Q+1.
\end{equation}
For $w\in\mathcal E_{v_4}^\perp\cap\mathcal H_{v_4}$, both $Qw$ and
$(w\cdot\nabla_{v_4})\xi$ vanish. The identity
\[
 (w\cdot\nabla_{v_4})\widehat K
 =\xi^{-1} (w\cdot\nabla_{v_4})\widehat K_0
      -\xi^{-2} K_0(w\cdot\nabla_{v_4})\xi
\]
then proves the second estimate in \eqref{su:eq:pair-source-regularity}
by \eqref{su:eq:unnormalized-source-response}.
For an arbitrary $w\in\mathcal H_{v_4}$, the same identity,
$1/2\le \xi\le3/2$, $\|Qv_4\|\le\|Pv_4\|$, and
\eqref{su:eq:normalization-derivative-bound} prove the first estimate in
\eqref{su:eq:pair-source-regularity}.
Step~1 now gives \eqref{su:eq:pair-exceptional-dimension}. Summing
\eqref{su:eq:pair-source-regularity} over an adapted orthonormal basis
gives \eqref{su:eq:pair-source-derivative-sum}.

\subsubsection{Step 4. The cutoff and its cost}
\label{su:sec:cutoff-estimates}
\label{su:par:common-cutoffs}

Use nested smooth cutoffs with support restrictions
\begin{equation}\label{su:eq:cuts}
 \|Pv_4\|^{-2} \|Qv_4\|^2\le Cg^2W^{3\varepsilon/200},\qquad
 \|Pv_4\|^{-2} \|QPv_4\|^2\le\tfrac18,\qquad
 \|v_{q-1}\|\ge v_{\min}\quad\text{in (i)},
\end{equation}
and
\begin{equation}\label{nc:eq:sourcecuts}
 \|v_4\|^{-2} \|Qv_4\|^2\le2W^{-1+14\theta},\qquad
 \|v_{q-1}\|\ge v_{\min}\quad\text{in (ii)}.
\end{equation}
Take fixed relative transition widths and the output plateau
$\|v_{q-1}\|\ge2v_{\min}$. By \eqref{su:eq:geometry} and
\eqref{nc:eq:geometry}, the overlap plateaux have probability $1-o(1)$.
By \eqref{us:eq:profile} and \eqref{nc:eq:profile}, the probability
that the output lies in the transition region and the chain is regular is
at most $W^{-D}$ for every fixed $D>0$ and sufficiently large $W$.
The $QPv_4$ restriction is unnecessary when $P=I_W$.

For the source family, multiply these cutoffs by
\[
 \varphi\!\left(W^{20}\|P\widehat v_4\|^2/\|\widehat v_4\|^2\right),
 \qquad \widehat v_4\ne0,
\]
where $\varphi\in C^\infty([0,\infty);[0,1])$ vanishes on $[0,1/4]$
and equals one on $[1,\infty)$. This factor equals one at every
prepared $v_4$ by \eqref{su:eq:prep}; in regime (ii), $P=I_W$.
On the resulting support, $\|\widehat v_4\|\le2W^{10}\|P\widehat v_4\|$, so
\eqref{su:eq:error} holds for the source family with at most an extra
factor two. Thus $1/2\le\xi\le3/2$ there for large $W$.
Define the source cutoff to be zero where $P\widehat v_4=0$;
the new factor vanishes in a neighborhood of this set.
It is homogeneous of degree zero, invariant under
\eqref{su:eq:cut-generator-equivariance}, and has zero derivative along
$\mathcal H_{\widehat v_4}$, since both squared norms have zero derivative
in these directions. It therefore adds no source derivative cost.

For $w\in\mathcal H_{v_4}$, the directional derivatives of
$P,Q$ and the norms $\|Pv_4\|,\|v_4\|$ vanish, so
\[
 \|Pv_4\|^2\bigl\|\nabla_{v_4}
       \bigl(\|Pv_4\|^{-2}\|Qv_4\|^2\bigr)\bigr\|^2
 \le4\|Pv_4\|^{-2}\|Qv_4\|^2.
\]
For the overlap cutoff $\widehat\chi_Q$, the fixed relative transition
width and chain rule therefore give, on its transition,
\[
 \|Pv_4\|^2\|\nabla_{v_4}\widehat\chi_Q\|^2
 \le C\|Qv_4\|^{-2} \|Pv_4\|^2.
\]
Together with $\|K\|_{\HS}^2\le32\|Qv_4\|^2/\|Pv_4\|^2$, this gives
\begin{equation}\label{su:eq:common-overlap-cutoff-cost}
 \|Pv_4\|^2\|K\|_{\HS}^2
       \|\nabla_{v_4}\widehat\chi_Q\|^2\le C.
\end{equation}
The same calculation with $QPv_4$ in place of $Qv_4$ shows that the
fixed-level cutoff for $\|QPv_4\|^2/\|Pv_4\|^2$ costs at most
$C\|K\|_{\HS}^2$.

For the lower cutoff $\widehat\chi_{q-1}$, take $v_{\min}$ with
$\nabla_{v_4}v_{\min}=0$ and require
\[
 \widehat\chi_{q-1}=0\quad(\|v_{q-1}\|\le v_{\min}),\qquad
 \widehat\chi_{q-1}=1\quad(\|v_{q-1}\|\ge2v_{\min}).
\]
Differentiating a smooth function of $\|v_{q-1}\|^2/v_{\min}^2$ gives
\begin{equation}\label{su:eq:common-norm-cutoff-cost}
 \|Pv_4\|^2\|\nabla_{v_4}\widehat\chi_{q-1}\|^2
 \le Cv_{\min}^{-2} \|Pv_4\|^2\|R_{[q-1,5]}\|^2
       \mathbf1_{\{v_{\min}\le\|v_{q-1}\|\le2v_{\min}\}}.
\end{equation}
The ratio in \eqref{su:eq:common-norm-cutoff-cost} is polynomially
bounded by \eqref{su:eq:output-derivative-ratio}. Since
\eqref{su:eq:normalized-generator-costs} bounds $\|K\|_{\HS}^2$ on the
support, the arbitrarily small polynomial transition probability makes
this expected weighted cost bounded in both regimes.

Impose regularity, \eqref{su:eq:active-chain-inverse}, and the spectral
moment and rank restrictions of Section~\ref{su:sec:retained-geometry},
with nested margins. Also impose the matrix support restriction
\begin{equation}\label{su:eq:compatibility-cutoff}
 \|(I-P)Q\|^2\le\tfrac18.
\end{equation}
Its plateau has probability $1-o(1)$ by \eqref{su:eq:geometry}; it is
automatic when $P=I_W$. Together with the source cutoffs
\eqref{su:eq:cuts}--\eqref{nc:eq:sourcecuts}, it gives
$B_Q\ge\tfrac34 I_{\ran Q}$ throughout the support for large $W$.
This restriction is invariant under the allowed unitaries and
independent of $v_4$, so it adds no source derivative cost.
Before compact exhaustion, remove the inverse-pole
sets and the set where a singular value of $A$ equals $t_{\mathrm{cut}}$.
The threshold is deterministic, and the latter set
is null: clearing inverse denominators gives a nonzero real polynomial,
since sufficiently large scalar right onsite matrices make all singular
values smaller than the threshold. On this invariant open set, the
retained projection and $Q$ are smooth.

Choose a compact subset losing at most $W^{-1}$ mass and saturate it
under the allowed compact unitary group. Average a smooth cutoff equal
to one on this saturation and supported in the invariant open set.
The average is still one on the saturation, so the mass loss remains
at most $W^{-1}$ without any invariance assumption on the Gaussian law.
Choose the exhaustion measurably in $R_4,P$ and independently of $v_4$,
so this additional cutoff has zero derivative in $v_4$.
Multiply by the source cutoffs above and set the result to zero on
failed preparation. This last indicator is evaluated at the actual
prepared vector and held fixed throughout the source family. Only the
active source cutoffs are varied under \eqref{su:eq:cut-generator-equivariance};
they are invariant because the source norms are preserved and $Q$ and
the output transform jointly. No preparation derivatives enter because
$X_1=\cdots=X_4=0$.

All plateaux have common probability $1-o(1)$. The product rule and
\eqref{su:eq:common-overlap-cutoff-cost}--\eqref{su:eq:common-norm-cutoff-cost}
give
\begin{equation}\label{su:eq:construction-cutoff-goals}
 0\le\chi\le1,\qquad \E\chi\ge\tfrac12,\qquad
 \E\left[\|Pv_4\|^2\|K\|_{\HS}^2
       \|\nabla_{v_4}\widehat\chi\|^2\right]\le C.
\end{equation}
This is \eqref{su:eq:pair-cutoff-derivative} together with the required
cutoff mass. Near the support of $\widehat\chi$, the inverses and the
retained spectral projection are smooth, and $\xi$ is bounded away from zero.
Under \eqref{su:eq:cut-generator-equivariance}, $P$ is fixed,
$A,\widetilde R,Q,\Pi$ are conjugated by $\mathcal U$, and $y$ becomes
$\mathcal U y$. Thus $K_0$ is conjugated and $\xi$ is unchanged.
Rescaling $\widehat v_4$ by $a\ne0$ scales $y$ by $a$ and leaves
$K_0,\xi$ unchanged. This gives the smooth family
\eqref{su:eq:generator-source-family} and its equivariance.
The support bounds follow from \eqref{su:eq:active-chain-inverse} and
\eqref{nc:eq:prepared}. The lengths in \eqref{su:eq:packet2} and
\eqref{nc:eq:parameters} are fixed before $W$ grows, and all estimates
are uniform in the energy, incoming flag, and coupling ranges of
Lemma~\ref{up:lem:small-contraction}. \qed

\section{Vector fields \texorpdfstring{$X$}{X} from a fixed matrix direction}
\label{sec:compact-proof}

We prove Lemma~\ref{up:lem:density-contraction} by the construction
explained in Section~\ref{sec:intro-fixed-fields}, adapting Drogin's
construction \cite[Proposition~4.1, Lemma~4.2 and Proposition~4.4]{Drogin}
to the GUE inputs.
Use the transfers \eqref{eq:c-transfer} and set $q=1000$,
$F=T_1\cdots T_q$ and $\mathcal T_q^{\mathrm{rev}}=S^*FS$.
The Lie algebra of the group \eqref{eq:c-group} is
\begin{equation}\label{up:eq:lie-algebra}
 \mathfrak g\coloneqq
 \{M\in\C^{2W\times2W}:M^*J_{\rm sig}+J_{\rm sig}M=0,\ \Tr M=0\}.
\end{equation}
For a fixed incoming flag $(U,Y)$ and any unit $e\in Y\cap U^\perp$, set
\[
 M=\tfrac1{\sqrt2}\,S\bigl(ee^*-(J_{\rm sig}e)(J_{\rm sig}e)^*\bigr)S^*.
\]
Since $J_{\rm sig}e\perp Y$, $M=M^*$, $\|M\|_{\rm HS}=1$ and
$S^*MS\in\mathfrak g$. The flow of $S^*MS$ fixes $U$ pointwise and
multiplies $e$ by $e^{t/\sqrt2}$. As $Y=U\oplus\C e$, this gives
\begin{equation}\label{up:eq:lift-flag-flow}
 e^{tS^*MS}U=U,\qquad e^{tS^*MS}Y=Y,\qquad
 c(e^{tS^*MS};U,Y)=t/\sqrt2.
\end{equation}
A general target $M$ with $S^*MS\in\mathfrak g$ has the block form
\begin{equation}\label{eq:c-Lie}
 M=\begin{pmatrix}A&B_1\\B_2&-A^*\end{pmatrix},\qquad
 B_1=B_1^*,\quad B_2=B_2^*,\quad\Im\Tr A=0.
\end{equation}
We seek a Hermitian input field $X$ with $D_XF=\sqrt2\,FM$.
Use the real Hilbert--Schmidt inner product
$\langle A,B\rangle=\Re\Tr(A^*B)$ for matrix targets and its sum over
components for input arrays. Linear maps carry the induced operator
norm and Jacobians the Hilbert--Schmidt norm. For a smooth compactly
supported field $\chi_P X$, Gaussian integration by parts with covariance
$I/W$ gives
\begin{equation}\label{up:eq:lift-divergence-reduction}
 \E|\delta_{\rm G}(\chi_P X)|^2
 \le W\E\|\chi_P X\|^2
      +\E\|\nabla_V(\chi_P X)\|_{\rm HS}^2.
\end{equation}
See the proof of \eqref{up:eq:lift-divergence-reduction} in
Appendix~\ref{sec:c-component-proof}.

\begin{lemma}[Gaussian divergence of a lift of a Lie direction]
\label{lem:c-gaussian-lift}
In the model \eqref{eq:model}, fix $B,E_0>0$ and take $q,F$ as above.
There is a smooth compactly supported cutoff $\chi_P$ in the GUE inputs with
\begin{equation}\label{eq:c-main-cutoff}
 \begin{gathered}
 0\le\chi_P\le1,\qquad \Prob(\chi_P=1)\ge1-Ce^{-cW},\qquad
 \|\nabla_V\chi_P\|\le C,\\
 \supp\chi_P\subset\{\max_i\|V_i\|_{\rm op}\le C\}.
 \end{gathered}
\end{equation}
Uniformly for deterministic
$W^{-B}\le g\le W^B$, $|E|\le E_0$, and every fixed $M$ of the form
\eqref{eq:c-Lie} with $\|M\|_{\rm HS}=1$, there is a field
$X$ in the original GUE variables such that
\[
 D_XF=\sqrt2\,FM\quad\text{on a neighborhood of }\supp\chi_P,
\]
and, for all sufficiently large $W$,
\begin{equation}\label{eq:c-gaussian-general}
 \E|\delta_{\rm G}(\chi_P X)|^2
 \le CW\max\{g^{-2},g^4\}.
\end{equation}
If $g\ge1$, $\kappa_1\le|E|\le E_0$, and $M=M^*$, where
$0<\kappa_1<E_0$ is fixed, a lift using the same cutoff satisfies
\begin{equation}\label{eq:c-gaussian-away}
 \E|\delta_{\rm G}(\chi_P X)|^2
 \le C_{\kappa_1}Wg^2.
\end{equation}
The cut fields satisfy the integration-by-parts regularity conditions
of Lemma~\ref{up:lem:angular}. Constants and width thresholds depend
only on $B,E_0$ and, for the second estimate, $\kappa_1$.
\end{lemma}

Set $Q_i\coloneqq T_{i+1}\cdots T_q$ for $1\le i<q$ and $Q_q=I$.
For a Hermitian input array $X$, the product rule at fixed $E,g$ gives
\begin{equation}\label{eq:c-product-map}
 \mathcal A(X)\coloneqq F^{-1}D_XF
 =g^{-1}\sum_{i=1}^qQ_i^{-1}
 \begin{pmatrix}0&0\\X_i&0\end{pmatrix}Q_i.
\end{equation}
Thus both fields solve $\mathcal A(X)=\sqrt2\,M$. For the cutoff
\eqref{eq:c-main-cutoff}, the product rule and
\eqref{up:eq:lift-divergence-reduction} give
\begin{equation}\label{eq:c-input-divergence}
 \E|\delta_{\rm G}(\chi_P X)|^2
 \le C\E\!\left[\mathbf1_{\supp\chi_P}
       \left(W\|X\|^2+\|\nabla_VX\|_{\rm HS}^2\right)\right].
\end{equation}

\subsection{The general perturbation and its Gaussian cost}

\begin{proof}[Proof of \eqref{eq:c-gaussian-general} in Lemma~\ref{lem:c-gaussian-lift}]
On a neighborhood of $\supp\chi_P$, the reconstruction
\eqref{eq:c-observe}--\eqref{eq:c-reconstruction} in
Appendix~\ref{sec:c-component-proof} gives a real linear map $\mathcal K_0$
on Hermitian input arrays such that, for every test matrix $Z$ of
the form \eqref{eq:c-Lie},
\begin{equation}\label{eq:c-main-left-inverse}
 \mathcal K_0(g\mathcal A^*Z)=Z.
\end{equation}
Taking adjoints gives $\mathcal A(g\mathcal K_0^*)=I$, so
$X=\sqrt2\,g\mathcal K_0^*M$ solves $\mathcal A(X)=\sqrt2 M$.
We use the following block formulas from the appendix to bound
$\mathcal K_0$ on arbitrary Hermitian data arrays $\mathbf R=(R_1,\ldots,R_q)$.
Set
\begin{equation}\label{eq:c-transfer-block}
 C_i\coloneqq (\max\{1,g\}E I_W-V_i)/g.
\end{equation}
Define $E_q=0$ and
$E_{i-1}=-R_iC_i-E_i^*$ for $2\le i\le q$. For odd
$3\le i\le q-1$, set
\begin{equation}\label{eq:c-data-map}
 F_i(\mathbf R)\coloneqq
 E_iC_i+C_iE_i^*+C_iR_iC_i+R_{i+1}-R_{i-1}.
\end{equation}
When $\mathbf R=g\mathcal A^*Z$, let $P$ be the upper-left block of $Z$.
The odd-index identity \eqref{eq:c-pair} gives
\begin{equation}\label{eq:c-block-equations}
 C_iP+P^*C_i=F_i(\mathbf R).
\end{equation}
For the disjoint odd-index pairs $(i_j,k_j)=(4j+1,4j-1)$, set
\[
 \widehat V_j\coloneqq2^{-1/2}g(C_{i_j}-C_{k_j})
 =2^{-1/2}(V_{k_j}-V_{i_j}),\qquad 1\le j\le128.
\]
Define
\[
 \mathcal F\mathbf R\coloneqq
 (F_{i_j}(\mathbf R)-F_{k_j}(\mathbf R))_{j=1}^{128},\qquad
 \mathcal L(P)\coloneqq(\widehat V_jP+P^*\widehat V_j)_{j=1}^{128},
\]
so that the paired equations give $\mathcal LP=2^{-1/2}g\mathcal F\mathbf R$.
The normalized differences $\widehat V_j$ are independent GUE matrices.
Appendix~\ref{app:compact-gaussian}, \eqref{cg:eq:SG}--\eqref{cg:eq:cutoff},
constructs the cutoff \eqref{eq:c-main-cutoff} with the additional bound
\begin{equation}\label{eq:c-main-coercivity}
 s_{\min}(\mathcal L|_{\{\Im\Tr P=0\}})\ge c
 \quad\text{on a neighborhood of }\supp\chi_P.
\end{equation}
Taking the left inverse on $\{\Im\Tr P=0\}$ gives
\[
 \mathcal J\coloneqq(\mathcal L^*\mathcal L)^{-1}\mathcal L^*,\qquad
 \mathcal K_P\coloneqq2^{-1/2}g\mathcal J\mathcal F,\qquad
 P=\mathcal K_P\mathbf R.
\]
Back substitution in \eqref{eq:c-reconstruction} gives
\begin{equation}\label{eq:c-main-reconstruction}
 \mathcal K_0\mathbf R=
 \begin{pmatrix}
 P&PC_q+C_qP^*+C_qR_qC_q-R_{q-1}\\
 R_q&-P^*
 \end{pmatrix},\qquad P=\mathcal K_P\mathbf R.
\end{equation}
This formula defines $\mathcal K_0$ on all data arrays. The left-inverse
identity \eqref{eq:c-main-left-inverse} applies when $\mathbf R=g\mathcal A^*Z$.

\smallskip\noindent\emph{Size of the perturbation.}
On $\supp\chi_P$, $\|C_i\|_{\rm op}\le C\max\{1,g^{-1}\}$. Each term
of $F_i(\mathbf R)$ contains one data block and at most two $C$ factors,
by the recursion for $E_i$. Since $q$ is fixed,
\begin{equation}\label{eq:c-Fi-size}
 \|F_i(\mathbf R)\|_{\rm HS}\le C\max\{1,g^{-2}\}\|\mathbf R\|.
\end{equation}
Taking the fixed number of differences defining $\mathcal F$ gives
\begin{equation}\label{eq:c-F-size}
 |\mathcal F|_{\mathrm{op}}\le C\max\{1,g^{-2}\}.
\end{equation}

The lower bound on $s_{\min}(\mathcal L)$ and \eqref{eq:c-F-size} give
\[
 \|\mathcal J\|,\ \|(\mathcal L^*\mathcal L)^{-1}\|\le C,\qquad
 \|\mathcal K_P\|\le Cg\|\mathcal F\|
 \le C\max\{g^{-1},g\}.
\]
Bounding the four blocks in \eqref{eq:c-main-reconstruction} for an
arbitrary input array gives
\[
 \|\mathcal K_0\|
 \le C\big[(1+\|C_q\|_{\rm op})\|\mathcal K_P\|+1+\|C_q\|_{\rm op}^2\big].
\]
Substituting the preceding bounds and using $X=\sqrt2\,g\mathcal K_0^*M$ yields
\begin{equation}\label{eq:c-lift-size}
 \|\mathcal K_0\|\le C\max\{g^{-2},g\},\qquad
 \|X\|\le C\max\{g^{-1},g^2\}\|M\|.
\end{equation}

\smallskip\noindent\emph{Input derivatives and divergence.}
Let $(H_\alpha)_{\alpha=1}^{W^2}$ be a trace-orthonormal basis of
Hermitian matrices and write $D_{i,\alpha}=H_\alpha\cdot\nabla_{V_i}$.
For any $W\times W$ matrix $B$ we have
\begin{equation}\label{eq:c-completeness}
 \sum_{\alpha=1}^{W^2}H_\alpha^2=WI_W,\qquad
 \sum_{\alpha=1}^{W^2}\|H_\alpha B\|_{\rm HS}^2
 =\sum_{\alpha=1}^{W^2}\|BH_\alpha\|_{\rm HS}^2=W\|B\|_{\rm HS}^2.
\end{equation}
Test the adjoints of $\mathcal F$ and $\mathcal L$ on a Hermitian array
$\mathbf Y=(Y_j)_{j=1}^{128}$. For example, differentiating the left
factor of $C_aR_kC_b$ in direction $D_{a,\alpha}$ gives, in the adjoint,
the Hermitian projection of $-g^{-1}H_\alpha Y_jC_b$. Its squared sum
over $\alpha$ is bounded by $Wg^{-2}\|Y_jC_b\|_{\rm HS}^2$ by
\eqref{eq:c-completeness}. The same bound with at most one remaining
$C$ factor applies to each differentiated term of $F_i$, giving
\begin{equation}\label{eq:c-F-bounds}
 \sum_{i=1}^q\sum_{\alpha=1}^{W^2}\|(D_{i,\alpha}\mathcal F)^*\mathbf Y\|^2
 \le CW\max\{g^{-2},g^{-4}\}\|\mathbf Y\|^2.
\end{equation}
For $\mathcal L$, if $i\in\{i_j,k_j\}$, then $V_i$ occurs only in
$\widehat V_j$, with $D_{i,\alpha}\widehat V_j=\pm2^{-1/2}H_\alpha$. Its contribution is
\[
 ((D_{i,\alpha}\mathcal L)P)_j
 =\pm2^{-1/2}(H_\alpha P+P^*H_\alpha),\qquad
 (D_{i,\alpha}\mathcal L)^*\mathbf Y=\pm\sqrt2\,H_\alpha Y_j.
\]
All unlisted components vanish. The adjoint lies in
$\{P:\Im\Tr P=0\}$ since $H_\alpha,Y_j$ are Hermitian, and
\eqref{eq:c-completeness} gives
\begin{equation}\label{eq:c-L-bounds}
 \sum_{i=1}^q\sum_{\alpha=1}^{W^2}\|(D_{i,\alpha}\mathcal L)P\|^2\le CW\|P\|^2,\qquad
 \sum_{i=1}^q\sum_{\alpha=1}^{W^2}\|(D_{i,\alpha}\mathcal L)^*\mathbf Y\|^2\le CW\|\mathbf Y\|^2.
\end{equation}

Differentiating $\mathcal J^*=\mathcal L(\mathcal L^*\mathcal L)^{-1}$ gives
\begin{equation}\label{eq:c-lift-derivative}
 D_{i,\alpha}\mathcal J^*
 =(I-\mathcal L\mathcal J)(D_{i,\alpha}\mathcal L)(\mathcal L^*\mathcal L)^{-1}
       -\mathcal J^*(D_{i,\alpha}\mathcal L)^*\mathcal J^*.
\end{equation}
Here $I-\mathcal L\mathcal J$ is an orthogonal projection. Applying
\eqref{eq:c-L-bounds} to each of the two terms, with its argument fixed
at the current input, yields
\begin{equation}\label{eq:c-J-bound}
 \sum_{i=1}^q\sum_{\alpha=1}^{W^2}\|(D_{i,\alpha}\mathcal J^*)P\|^2
 \le CW\bigl(\|(\mathcal L^*\mathcal L)^{-1}P\|^2
                 +\|\mathcal J^*P\|^2\bigr)
 \le CW\|P\|^2.
\end{equation}

For $\mathcal K_P^*=2^{-1/2}g\mathcal F^*\mathcal J^*$, the product
rule gives one term with $D_{i,\alpha}\mathcal F^*$ and one with
$D_{i,\alpha}\mathcal J^*$. Applying \eqref{eq:c-F-bounds} to
$\mathcal J^*P$ in the first, and \eqref{eq:c-F-size} and
\eqref{eq:c-J-bound} in the second, gives
\begin{equation}\label{eq:c-KP-bound}
 \sum_{i,\alpha}\|(D_{i,\alpha}\mathcal K_P)^*P\|^2
 \le CWg^2\big(\max\{g^{-2},g^{-4}\}+\max\{1,g^{-4}\}\big)\|P\|^2
 \le CW\max\{g^{-2},g^2\}\|P\|^2 .
\end{equation}

To differentiate $\mathcal K_0^*M$ with $M$ fixed, take the adjoint of
\eqref{eq:c-main-reconstruction}. With $\iota_i$ denoting insertion in
data slot $i$, cyclicity of the trace gives
\begin{equation}\label{eq:c-K0-adjoint}
 \mathcal K_0^*M
 =\mathcal K_P^*\bigl(2(A+B_1C_q)\bigr)
       +\iota_q(C_qB_1C_q+B_2)-\iota_{q-1}B_1 .
\end{equation}
Here $\Im\Tr(A+B_1C_q)=0$. When the derivative hits $\mathcal K_P^*$,
apply \eqref{eq:c-KP-bound} to $2(A+B_1C_q)$, whose norm is at most
$C\max\{1,g^{-1}\}\|M\|$. When it hits $C_q$, the identity
$D_{i,\alpha}C_q=-g^{-1}\delta_{iq}H_\alpha$ and \eqref{eq:c-completeness}
bound their contribution to the squared sum by
$CWg^{-2}(\|\mathcal K_P\|^2+\|C_q\|_{\rm op}^2)\|B_1\|^2$. Thus
\[
 \sum_{i,\alpha}\|D_{i,\alpha}(\mathcal K_0^*M)\|^2
 \le CW\max\{g^{-4},g^2\}\|M\|^2 .
\]
Since $X=\sqrt2\,g\mathcal K_0^*M$, it follows that
\begin{equation}\label{eq:c-lift-cost}
 \|\nabla_VX\|_{\rm HS}\le C\sqrt W\max\{g^{-1},g^2\}\|M\|.
\end{equation}
For $\|M\|_{\rm HS}=1$, inserting \eqref{eq:c-lift-size} and
\eqref{eq:c-lift-cost} into \eqref{eq:c-input-divergence} proves
\eqref{eq:c-gaussian-general}.
\end{proof}

\subsection{The Hermitian perturbation away from zero energy}

\begin{proof}[Proof of \eqref{eq:c-gaussian-away} in Lemma~\ref{lem:c-gaussian-lift}]
Assume $g\ge1$, $\kappa_1\le|E|\le E_0$, and $M=M^*$.
For a test matrix $Z$ write $(P,U,R_q)$ for the blocks $(A,B_1,B_2)$
in \eqref{eq:c-Lie}, and let $\mathbf R=g\mathcal A^*Z$.
Since $C_i=EI_W-g^{-1}V_i$, at an interior odd index
\eqref{eq:c-block-equations} reads
\[
 E(P+P^*)=F_i(\mathbf R)+g^{-1}(V_iP+P^*V_i).
\]
By \eqref{eq:c-lift-size},
$\|P\|_{\rm HS}\le\|Z\|_{\rm HS}=\|\mathcal K_0\mathbf R\|
\le Cg\|\mathbf R\|$. Thus \eqref{eq:c-Fi-size} and
$\|V_i\|_{\rm op}\le C$ give
\[
 \kappa_1\|P+P^*\|_{\rm HS}
 \le C\|\mathbf R\|+Cg^{-1}\|P\|_{\rm HS}
 \le C\|\mathbf R\|.
\]
By \eqref{eq:c-main-reconstruction},
\[
 U=E(P+P^*)-g^{-1}(PV_q+V_qP^*)+C_qR_qC_q-R_{q-1},
\]
so $\|U\|_{\rm HS}\le E_0\|P+P^*\|_{\rm HS}
+Cg^{-1}\|P\|_{\rm HS}+C\|\mathbf R\|$. Together with
$\|R_q\|_{\rm HS}\le\|\mathbf R\|$ this gives
\begin{equation}\label{lg:eq:hermitian-blocks}
 \|P+P^*\|_{\rm HS}+\|U\|_{\rm HS}+\|R_q\|_{\rm HS}
 \le C_{\kappa_1}g\|\mathcal A^*Z\|.
\end{equation}
For the Hermitian target $M$ in \eqref{eq:c-Lie}, we have $A=A^*$
and $B_1=B_2\coloneqq B$. The pairing is
\[
 \langle M,Z\rangle=\Re\Tr[A(P+P^*)+B(U+R_q)].
\]
Cauchy--Schwarz, with $\|M\|_{\rm HS}=1$, gives the dual estimate
\begin{equation}\label{lg:eq:hermitian-dual}
 |\langle M,Z\rangle|\le C_{\kappa_1}g\|\mathcal A^*Z\|.
\end{equation}
Set $\Gamma=\mathcal A\mathcal A^*$. Since
$\|Z\|\le Cg\|\mathbf R\|=Cg^2\|\mathcal A^*Z\|$, we have
\begin{equation}\label{lg:eq:coercivity}
 \Gamma\ge cg^{-4}I .
\end{equation}
Take the minimum-norm solution
\begin{equation}\label{lg:eq:minimal-lift}
 X=\sqrt2\,\mathcal A^*\Gamma^{-1}M,\qquad
 \mathcal A(X)=\sqrt2 M,\qquad X\perp\ker\mathcal A .
\end{equation}
Applying \eqref{lg:eq:hermitian-dual} to $Z=\Gamma^{-1}M$ yields
\[
 \|\mathcal A^*\Gamma^{-1}M\|^2
 =\langle M,\Gamma^{-1}M\rangle
 \le C_{\kappa_1}g\|\mathcal A^*\Gamma^{-1}M\|
\]
and hence $\langle M,\Gamma^{-1}M\rangle\le C_{\kappa_1}g^2$.
Using $\|\Gamma^{-1}M\|^2\le\|\Gamma^{-1}\|\langle M,\Gamma^{-1}M\rangle$
and \eqref{lg:eq:coercivity}, we obtain
\begin{equation}\label{lg:eq:hermlift}
 \|X\|^2=2\langle M,\Gamma^{-1}M\rangle\le C_{\kappa_1}g^2,\qquad
 \|\Gamma^{-1}M\|\le C_{\kappa_1}g^3.
\end{equation}

\smallskip\noindent\emph{Input derivatives and divergence.}
By \eqref{eq:c-transfer},
differentiating a factor $T_\ell$ or $T_\ell^{-1}$ inserts $H_\alpha$
with a factor $g^{-1}$. The map in \eqref{eq:c-product-map} and its
adjoint already have a prefactor $g^{-1}$, while all
undifferentiated factors are bounded on the cutoff neighborhood.
Thus a typical term in $(D_{i,\alpha}\mathcal A^*)Z$ has the form
$g^{-2}B_0H_\alpha B_1Z_{ab}B_2$, with bounded $B_j$, or the same
form with $H_\alpha$ to the right of $Z_{ab}$. By
\eqref{eq:c-completeness}, its square sum over $\alpha$ is at most
$CWg^{-4}\|Z_{ab}\|_{\rm HS}^2$.
For a Hermitian input array $H'$, the same calculation applies to
$(D_{i,\alpha}\mathcal A)H'$ with $H'_k$ in place of $Z_{ab}$.
Summing over the fixed number of
terms, blocks, and input slots gives
\begin{equation}\label{lg:eq:A-derivative-bounds}
 \sum_{i=1}^q\sum_{\alpha=1}^{W^2}\|(D_{i,\alpha}\mathcal A^*)Z\|^2\le CWg^{-4}\|Z\|^2,\qquad
 \sum_{i=1}^q\sum_{\alpha=1}^{W^2}\|(D_{i,\alpha}\mathcal A)H'\|^2\le CWg^{-4}\|H'\|^2.
\end{equation}

The inverse differentiation identity \eqref{eq:c-lift-derivative},
with $\mathcal L=\mathcal A^*$, applied to \eqref{lg:eq:minimal-lift} gives
\[
 D_{i,\alpha}X=\sqrt2\,(I-\mathcal A^*\Gamma^{-1}\mathcal A)
       (D_{i,\alpha}\mathcal A^*)\Gamma^{-1}M
       -\mathcal A^*\Gamma^{-1}(D_{i,\alpha}\mathcal A)X.
\]
The first prefactor in parentheses is an orthogonal projection, and
$\|\mathcal A^*\Gamma^{-1}\|^2=\|\Gamma^{-1}\|\le Cg^4$ by
\eqref{lg:eq:coercivity}. Applying \eqref{lg:eq:A-derivative-bounds}
to the two terms and then \eqref{lg:eq:hermlift} gives
\begin{equation}\label{lg:eq:DV}
 \begin{aligned}
 \|\nabla_VX\|_{\rm HS}^2
 &\le CWg^{-4}\|\Gamma^{-1}M\|^2
       +CWg^{-4}\|\mathcal A^*\Gamma^{-1}\|^2\|X\|^2\\
 &\le CW\bigl(g^{-4}\|\Gamma^{-1}M\|^2+\|X\|^2\bigr)
 \le C_{\kappa_1}Wg^2.
 \end{aligned}
\end{equation}
Inserting \eqref{lg:eq:hermlift} and~\eqref{lg:eq:DV} into
\eqref{eq:c-input-divergence} proves \eqref{eq:c-gaussian-away}.
For both fields the inverses in their construction are smooth on a
neighborhood of $\supp\chi_P$ by \eqref{eq:c-main-coercivity} and
\eqref{lg:eq:coercivity}, respectively. Multiplication by $\chi_P$
and extension by zero give the required smooth compactly supported fields.
\end{proof}

\begin{proof}[Proof of Lemma~\ref{up:lem:density-contraction}]
Take the unit Hermitian target $M$ of \eqref{up:eq:lift-flag-flow}.
For either field of Lemma~\ref{lem:c-gaussian-lift},
$D_X\mathcal T_q^{\mathrm{rev}}=\sqrt2\,\mathcal T_q^{\mathrm{rev}}S^*MS$
on $\supp\chi_P$. Equation~\eqref{up:eq:lift-flag-flow} and the cocycle
identity give $D_Xc=1$ and zero variation of both terminal subspaces,
as required in \eqref{up:eq:angular-geometry} with $K_{\mathrm{out}}=0$.
The cutoff bound \eqref{eq:c-main-cutoff} gives $\E\chi_P\ge1/2$ for
large $W$, and \eqref{eq:c-gaussian-general} and \eqref{eq:c-gaussian-away}
give the two divergence bounds. Reversing the inputs, field and cutoff
restores the order of Definition~\ref{up:def:transfer-operator} and
preserves their bounds by the iid law.
\end{proof}

\section{Large-coupling transfer growth}
\label{sec:large-upper-proof}

The direct transfer argument is explained in Section~\ref{sec:intro-direct-transfer}.

\subsection{One-step logarithmic growth and its conditional mean}
\label{lg:sec:scalar-growth}

Fix $g\ge1$, $|E|\le E_0<2$ and $\Phi\in\mathcal F_W$. To estimate
$D_n(\Phi;E)^{-1}$ from \eqref{up:eq:D-definition}, use the free
coordinates $S_E$, rotation $R_E$, and angle $E=2\sin\theta$ of
Definition~\ref{def:transfer-matrices}. Set
\begin{equation}\label{lg:eq:normal}
 h_\theta:=(2\cos\theta)^{-1}.
\end{equation}
Set
\[
 J_{\rm ph}\coloneqq\begin{pmatrix}0&-iI_W\\iI_W&0\end{pmatrix},
 \qquad S^*J_{\rm ph}S=S_E^*J_{\rm ph}S_E=J_{\rm sig}.
\]
These identities show that $S_E^{-1}S$ preserves $J_{\rm sig}$-isotropy.
For each transfer $T_i$ in \eqref{eq:c-transfer}, multiplication gives
\[
 S_E^{-1}T_iS_E=R_E
 \left[I-ih_\theta g^{-1}\begin{pmatrix}V_i&V_i\\-V_i&-V_i\end{pmatrix}\right].
\]
Every maximal $J_{\rm sig}$-isotropic subspace has the graph representation
\begin{equation}\label{up:eq:graph-subspace}
 Y_C\coloneqq\left\{\begin{pmatrix}u\\Cu\end{pmatrix}:u\in\C^W\right\},
 \qquad C\in U(W).
\end{equation}
Since $\Phi\in\mathcal F_W$, the transformed
initial frame $S_E^{-1}S\Phi$ is also isotropic. Write
\[
 S_E^{-1}S\Phi=\begin{pmatrix}I_W\\C_0\end{pmatrix}B, \quad C_0\in U(W).
\]
Thus the initial subspace in these coordinates is $Y_{C_0}$.
Conditional on $V_1,\ldots,V_{i-1}$, the graph $Y_{C_{i-1}}$ is fixed.
The independent GUE block $V_i$ gives the next graph $Y_{C_i}$ by
\begin{equation}\label{lg:eq:matrix-update}
 \begin{gathered}
 F(C,V)\coloneqq I-ih_\theta g^{-1}V(I+C),\qquad
 C_i\coloneqq -e^{-2i\theta}[C_{i-1}+ih_\theta g^{-1}V_i(I+C_{i-1})]F(C_{i-1},V_i)^{-1},\\
 S_E^{-1}T_iS_E\begin{pmatrix}I_W\\C_{i-1}\end{pmatrix}
 =-ie^{i\theta}\begin{pmatrix}I_W\\C_i\end{pmatrix}F(C_{i-1},V_i).
 \end{gathered}
\end{equation}
In particular $C_i$ is unitary by isotropy.
For a general initial subspace $Y_{C_0}$, define the normal-coordinate
version of $D_n$ by
\begin{equation}\label{lg:eq:normal-transfer}
 D_n( {C_0};E)\coloneqq
 \begin{pmatrix}I_W&0\end{pmatrix}
 S_E^{-1}(T_n\cdots T_1)S_E\begin{pmatrix}I_W\\C_0\end{pmatrix}.
\end{equation}
Iterating \eqref{lg:eq:matrix-update} gives the
factorization
\begin{equation}\label{lg:eq:matrix-product}
 \begin{aligned}
 D_n(C_0;E)&=(-i)^ne^{in\theta}F(C_{n-1},V_n)\cdots F(C_0,V_1).
 \end{aligned}
\end{equation}
Returning to \eqref{up:eq:D-definition} gives
\[
 D_n(\Phi;E)=
 \sqrt2\begin{pmatrix}I_W&0\end{pmatrix}S^*S_E
 \begin{pmatrix}I_W\\C_n\end{pmatrix}D_n(C_0;E)B.
\]
Isotropy and the uniform bounds on the coordinate changes give
\[
 2\left\|\begin{pmatrix}I&0\end{pmatrix}S^*S_E
       \binom{I}{C_n}x\right\|^2
 =\left\|S^*S_E\binom{I}{C_n}x\right\|^2\asymp\|x\|^2.
\]
Similarly, $2\|Bx\|^2=\|S_E^{-1}S\Phi x\|^2\asymp\|x\|^2$
because $\Phi\in\mathcal F_W$. Thus both exterior factors are uniformly
invertible on $|E|\le E_0<2$, and
\begin{equation}\label{lg:eq:coordinate-comparison}
 \|D_n(\Phi;E)^{-1}\|
 \asymp\|D_n(C_0;E)^{-1}\|.
\end{equation}
We use $A^{-*}\coloneqq(A^{-1})^*=(A^*)^{-1}$.
For fixed $C_0$ and a fixed unit vector $v_0$,
\eqref{lg:eq:matrix-product} therefore gives
\[
 D_i(C_0;E)^{-*}v_0
 =-ie^{i\theta}F(C_{i-1},V_i)^{-*}D_{i-1}(C_0;E)^{-*}v_0.
\]
Normalize after each step. The incoming direction $v_{i-1}$ then depends
only on $V_1,\ldots,V_{i-1}$ and is independent of $V_i$. Define
\begin{equation}\label{lg:eq:update}
 \begin{aligned}
 v_i&\coloneqq-ie^{i\theta}
 \frac{F(C_{i-1},V_i)^{-*}v_{i-1}}
      {\|F(C_{i-1},V_i)^{-*}v_{i-1}\|},\quad
 c_i\coloneqq-\log\|F(C_{i-1},V_i)^{-*}v_{i-1}\|.
 \end{aligned}
\end{equation}
Multiplication of the normalization factors gives
\begin{equation}\label{lg:eq:additive-growth}
 \ell_n(C_0,v_0)\coloneqq-\log\|D_n(C_0;E)^{-*}v_0\|
 =\sum_{i=1}^n c_i.
\end{equation}
Geometrically, $v_{i-1}$ specifies the subspace
$\{(u,C_{i-1}u):v_{i-1}^*u=0\}$ inside $Y_{C_{i-1}}$, so $c_i$ is also its one-step
logarithmic volume ratio \eqref{eq:c-volume-cocycle}.

Taking $v_0$ in an orthonormal basis $(e_a)_{a=1}^W$ and using
\eqref{lg:eq:coordinate-comparison} bounds the matrix inverse by
\begin{equation}\label{lg:eq:scalar-reduction}
 \|D_n(\Phi;E)^{-1}\|^2
 \lesssim \sum_{a=1}^W e^{-2\ell_n(C_0,e_a)}
 \le W\max_{1\le a\le W} e^{-2\ell_n(C_0,e_a)} .
\end{equation}
Write $\E_{i-1}$ for expectation conditional on $V_1,\ldots,V_{i-1}$.
The trace and directional moments at the incoming state are
\begin{equation}\label{lg:eq:observables}
 f_k(C_{i-1},v_{i-1})\coloneqq v_{i-1}^*C_{i-1}^kv_{i-1},\qquad
 m_k(C_{i-1})\coloneqq W^{-1}\Tr C_{i-1}^k,\qquad f_0=m_0=1.
\end{equation}
For $k\in\mathbb Z$, unitarity gives $|m_k|,|f_k|\le1$.

For functions of the state, write
\begin{equation}\label{lg:eq:transition}
 (P_gf)(C_{i-1},v_{i-1})\coloneqq\E_{i-1}f(C_i,v_i).
\end{equation}
In the calculations below, $m_k,f_k$ without an argument are
evaluated at the incoming state, with
$m_{-k}=\overline{m_k}$ and $f_{-k}=\overline{f_k}$.

\noindent\textbf{The second-order conditional mean.}
Uniformly in the incoming state, $W$, and $|E|\le E_0<2$, for all
sufficiently large $g$ we have
\begin{equation}\label{lg:eq:c2}
 \E_{i-1}c_i=h_\theta^2g^{-2}\left[\tfrac12 W^{-1}+\Re f_1+\Re(m_1f_1)
                         -\tfrac12 W^{-1}\Re(f_1^2)\right]+O(g^{-4}).
\end{equation}

\begin{proof}[Proof of \eqref{lg:eq:c2}]
By \eqref{lg:eq:update},
\[
 c_i=-\tfrac12\log\left[
 v_{i-1}^*\bigl(F(C_{i-1},V_i)^*F(C_{i-1},V_i)\bigr)^{-1}v_{i-1}
 \right].
\]
Multiplication gives
\[
 F(C_{i-1},V_i)^*F(C_{i-1},V_i)
 =I+h_\theta g^{-1}A+h_\theta^2g^{-2}B,
\]
where the linear and quadratic terms in $V_i$ are
\[
 A\coloneqq-i(V_iC_{i-1}-C_{i-1}^*V_i),\qquad
 B\coloneqq(I+C_{i-1}^*)V_i^2(I+C_{i-1}).
\]
 
Taking the quadratic form against $v_{i-1}$ and then expanding the
scalar logarithm gives
\[
 \begin{aligned}
 \E_{i-1}c_i=h_\theta^2g^{-2}\Big[
 &\tfrac12\E_{i-1}(v_{i-1}^*Bv_{i-1})
 -\tfrac12\E_{i-1}(v_{i-1}^*A^2v_{i-1}) +\tfrac14\E_{i-1}\bigl[(v_{i-1}^*Av_{i-1})^2\bigr]\Big]+O(g^{-4}).
 \end{aligned}
\]
The GUE covariance \eqref{eq:covariance} gives
\[
 \begin{aligned}
 \E_{i-1}(v_{i-1}^*Bv_{i-1})
 &=v_{i-1}^*(2I+C_{i-1}+C_{i-1}^*)v_{i-1}
 =2+2\Re f_1,\\
 \E_{i-1}(v_{i-1}^*A^2v_{i-1})
 &=v_{i-1}^*(2I-m_1C_{i-1}-\overline{m_1}C_{i-1}^*)v_{i-1}
 =2-2\Re(m_1f_1).
 \\
 \E_{i-1}\bigl[(v_{i-1}^*Av_{i-1})^2\bigr]
 &=2\E_{i-1}|v_{i-1}^*V_iC_{i-1}v_{i-1}|^2
   -2\Re\E_{i-1}\bigl[(v_{i-1}^*V_iC_{i-1}v_{i-1})^2\bigr]\\
 &=2W^{-1}(1-\Re(f_1^2)).
 \end{aligned}
\]
Substituting these expressions gives the coefficient in \eqref{lg:eq:c2}.
For the remainder, \eqref{lg:eq:matrix-update} gives the exact identity
\[
 F(C_{i-1},V_i)^{-1}
 =I+ih_\theta g^{-1}V_i(I-e^{2i\theta}C_i).
\]
Since $C_{i-1},C_i$ are unitary, it follows that
\begin{equation}\label{uvg:eq:global}
 \norm{F(C_{i-1},V_i)},\norm{F(C_{i-1},V_i)^{-1}}\le1+Cg^{-1}\norm{V_i},\qquad
 |c_i|\le\log(1+Cg^{-1}\norm{V_i}).
\end{equation}
To justify the expected remainder, expand through cubic order on
\(\{\norm{V_i}\le c_0g\}\), with a sufficiently small fixed \(c_0>0\).
Odd terms vanish under \(V_i\mapsto-V_i\), and the fourth-order
remainder is \(Cg^{-4}\norm{V_i}^4\). A sphere net and scalar Gaussian
tails give
\begin{equation}\label{uvg:eq:tail}
 \Prob(\norm{V_i}>R)\le2\exp(2W\log9-WR^2/8),
 \qquad \E\norm{V_i}^p\le C_p
\end{equation}
for every fixed \(p\). On the complement of this cutoff, \eqref{uvg:eq:global},
boundedness of \(m_k,f_k\), and sufficiently high fixed moments
give another \(O(g^{-4})\) error.
\end{proof}

\medskip\noindent\textbf{The one-step fluctuations.}
For the fluctuation bounds we differentiate with the incoming state fixed.
For a Hermitian variation $\dot V_i$, differentiating
\eqref{lg:eq:matrix-update} and~\eqref{lg:eq:update} gives
\begin{equation}\label{lg:eq:gradients}
 \begin{aligned}
 dc_i&=\Re[ih_\theta g^{-1}v_i^*F(C_{i-1},V_i)^{-*}
             (I+C_{i-1}^*)\dot V_i v_i],\\
 dC_i&=ih_\theta g^{-1}(C_i-e^{-2i\theta}I)\dot V_i
             (I+C_{i-1})F(C_{i-1},V_i)^{-1}.
 \end{aligned}
\end{equation}
The derivative of $v_i$ is the derivative of $F^{-*}v_{i-1}$ followed
by differentiation of its normalization in \eqref{lg:eq:update}.

For a bounded real $Q$ put $Q_i=Q(C_i,v_i)$ and $Y_i=c_i-Q_i+Q_{i-1}$
as in \cite[Lemma~1]{SBanomaly}, \cite[Section~3]{SSB}. Then
\[
 \ell_n=Q_n-Q_0+\sum_{i=1}^n\E_{i-1}Y_i+
                    \sum_{i=1}^n(1-\E_{i-1})Y_i .
\]

\subsection{The quadratic correction for \texorpdfstring{$g\ge C_\kappa\sqrt W$}{g at least C kappa times sqrt W}}
\label{sec:uniform-large-g}

Fix $0<E_0<2$. For $E_0=2-\kappa$, denote the coupling threshold
constant below by $C_\kappa$.

\begin{theorem}[Negative transfer moments]\label{uvg:thm:transfer}
There exist \(C_*,W_0,C,c>0\) and \(s\in(0,1)\) such that, for
every \(W\ge W_0\), finite \(g\ge C_*\sqrt W\), integer \(n\ge0\),
deterministic \(|E|\le E_0\), and deterministic
\(\Phi\in\mathcal F_W\),
\[
 \E\norm{D_n(\Phi;E)^{-1}}^s
 \le CW\exp\{-cn/(Wg^2)\}.
\]
The initial frame may depend on random variables independent of
\(V_1,\ldots,V_n\). The exponent and width threshold are independent
of \(g,n\).
\end{theorem}

\medskip\noindent\textbf{The linear and quadratic corrections.}
Write $h=h_\theta$, $X=(C,v)$ for the incoming state, and $X'=(C',v')$
for the next state. Denote the increment $c_i$ by $c(X,V)$ and its
conditional mean by $P_gc$. On the bulk interval,
$h,h^{-1}$ and $(1+e^{-2i\theta})^{-1}$ are uniformly bounded.

Under the free step, $f_1$ acquires the phase $-e^{-2i\theta}$.
To cancel $h^2g^{-2}\Re f_1$ in \eqref{lg:eq:c2}, set
\[
 Q=-g^{-2}\Re\left[\frac{h^2}{1+e^{-2i\theta}}f_1\right].
\]

Define the corrected moments
\begin{align}
 M&=m_1-h^2g^{-2}\frac{e^{-2i\theta}+m_2+2m_1^2}{1+e^{-2i\theta}},
 \label{uvg:eq:M}\\
 F&=f_1-h^2g^{-2}\frac{4m_1f_1+2f_2-f_1^2-e^{-2i\theta}f_1f_{-1}}
 {1+e^{-2i\theta}}.\label{uvg:eq:F}
\end{align}
Thus \(M=m_1+O(g^{-2})\), \(F=f_1+O(g^{-2})\).

\begin{lemma}[Second-moment drifts]\label{uvg:lem:drifts}
For constants \(a,C>0\), uniformly in the state,
\begin{align}
 (P_g-I)|M|^2&\le-ag^{-2}|M|^2+Cg^{-2}(W^{-2}+g^{-4}),\label{uvg:eq:Mdrift}\\
 (P_g-I)|F|^2&\le-ag^{-2}|F|^2+Cg^{-2}|M|^2
                         +Cg^{-2}(W^{-1}+g^{-4}).\label{uvg:eq:Fdrift}
\end{align}
\end{lemma}

\begin{proof}[Proof of Lemma~\ref{uvg:lem:drifts}]
Write $A=I-ihg^{-1}V(I+C)$ in the graph update. Expanding with the cutoff used
for \eqref{lg:eq:c2} gives, uniformly in the incoming state,
\begin{align}
 P_gm_1&=-e^{-2i\theta}m_1+e^{-2i\theta}h^2g^{-2}(1+m_1)(1+2m_1+m_2)+O(g^{-4}),
 \label{uvg:eq:mmean}\\
 P_gf_1&=-e^{-2i\theta}f_1+e^{-2i\theta}h^2g^{-2}(B_o+B_e)+O(g^{-4}+g^{-2}/W),
 \label{uvg:eq:fmean}
\end{align}
where
\[
 \begin{aligned}
 B_o&=(3+m_2)f_1+2m_1f_2-f_1(m_1f_1+m_{-1}f_{-1}),\\
 B_e&=4m_1f_1+2f_2-f_1^2-f_1f_{-1}.
 \end{aligned}
\]
To prove \eqref{uvg:eq:mmean}--\eqref{uvg:eq:fmean},
expand the graph update and use
$\E VBV=W^{-1}\Tr(B)I$. To see the directional term, set
$Y=v^*A^{-1}A^{-*}v$ and
$N=-e^{2i\theta}v^*A^{-1}C'A^{-*}v$. With $a=h/g$, write
$Y=1+ay_1+a^2y_2+\cdots$ and $N=f_1+an_1+a^2n_2+\cdots$.
The contractions give
\[
 \begin{aligned}
 \E y_2&=-[(1+m_1)f_1+(1+m_{-1})f_{-1}],\\
 \E n_2&=-(3+4m_1+m_2)f_1-2(1+m_1)f_2,\\
 \E(n_1y_1)-f_1\E y_1^2&=O(W^{-1}).
 \end{aligned}
\]
Inserting these identities in
$\E(N/Y)=f_1+a^2\E[n_2-n_1y_1+f_1(y_1^2-y_2)]+O(g^{-4})$
gives \eqref{uvg:eq:fmean}. The normalized trace of the graph
expansion gives \eqref{uvg:eq:mmean}.

The zero-noise step multiplies \(m_k,f_k\) by \((-e^{-2i\theta})^k\).
Applying it to the correction polynomials yields
\begin{equation}\label{uvg:eq:means}
 \begin{gathered}
 P_gM=\alpha M+O(g^{-4}),\qquad
 P_gF=\alpha F+O(g^{-2}(|M|+g^{-2}+W^{-1})),\\
 \alpha=-e^{-2i\theta}[1-h^2g^{-2}(3+m_2)].
 \end{gathered}
\end{equation}
Since \(\Re(3+m_2)\ge2\), \(|\alpha|^2\le1-a_1g^{-2}\)
for some fixed \(a_1>0\), when \(g^{-2}\) is small.

The Gaussian Poincare inequality in real Hilbert--Schmidt
coordinates of variance \(W^{-1}\) gives
\(\E|Z-\E Z|^2\le W^{-1}\E\norm{\nabla Z}_{\HS}^2\).
By \eqref{lg:eq:gradients} and the normalization in \eqref{lg:eq:update},
the vector derivative has norm at most
\(Cg^{-1}(1+g^{-1}\norm V)^r\norm{dV}_{\HS}\)
for a fixed \(r\). Product differentiation then gives
\[
 \norm{\nabla M(C')}_{\HS}\le
 Cg^{-1}W^{-1/2}(1+g^{-1}\norm V)^r,\quad
 \norm{\nabla F(C',v')}_{\HS}\le
 Cg^{-1}(1+g^{-1}\norm V)^r.
\]
The normalized trace supplies \(W^{-1/2}\).
Thus conditional innovation variances are at most \(Cg^{-2}/W^2\)
and \(Cg^{-2}/W\), respectively. In
\(P_g|Z|^2=|P_gZ|^2+\E|Z-P_gZ|^2\), insert \eqref{uvg:eq:means}.
Young's inequality absorbs small multiples of \(g^{-2}|M|^2,g^{-2}|F|^2\).
The \(Cg^{-4}\) mean remainder costs \(Cg^{-6}\); the other costs
\(Cg^{-2}(|M|^2+g^{-4}+W^{-2})\). This proves both bounds.
\end{proof}

\begin{lemma}[Uniform positive corrected drift]\label{uvg:lem:positive}
There are fixed \(A,B,a_0,C_*,W_0>0\) such that for
\(W\ge W_0\), \(g\ge C_*\sqrt W\), the bounded real function
\[
 \mathcal Q=Q+R,\qquad R=A|M|^2+B|F|^2
\]
satisfies
\begin{equation}\label{uvg:eq:positive}
 P_gc-(P_g-I)\mathcal Q\ge a_0g^{-2}/W.
\end{equation}
Its supremum norm is independent of \(W,g\). For a fixed sufficiently
small $c_0>0$, on \(\{\norm V\le c_0g\}\), the derivative in $V$ of
\(c-\mathcal Q(C',v')\) is at most \(C/g\) in Hilbert--Schmidt norm.
\end{lemma}

\begin{proof}[Proof of Lemma~\ref{uvg:lem:positive}]
By \eqref{lg:eq:c2} and \eqref{uvg:eq:fmean}, the linear correction leaves
\[
 d:=P_gc-(P_g-I)Q=\frac{h^2g^{-2}}{2W}
 +h^2g^{-2}\left[\Re(m_1f_1)-\frac{\Re(f_1^2)}{2W}\right]+O(g^{-4}).
\]
Let \(h^2/2\ge a_*>0\). For any fixed \(\varepsilon>0\),
Young's inequality and sufficiently large \(W\) give
\[
 h^2|m_1f_1|+\frac{h^2}{2W}|f_1|^2
 \le K_\varepsilon|M|^2+\varepsilon|F|^2
                                     +K_\varepsilon g^{-4}.
\]
Choose \(B>0\) so small that the noise term \(CB/W\)
from \eqref{uvg:eq:Fdrift} is at most \(a_*/(8W)\).
Take \(\varepsilon=aB/2\), then choose \(A\) so large that
\(aA-CB\ge K_\varepsilon+1\).
The coefficients of \(|M|^2,|F|^2\) in
\(d-(P_g-I)R\) are nonnegative. The remaining lower bound is
\[
 g^{-2}\left[\frac{a_*}{W}-C\frac A{W^2}-C\frac B W
              -Cg^{-2}-C(A+B+K_\varepsilon)g^{-4}\right].
\]
Increase \(C_*\) to make \(Cg^{-2}\le a_*/(8W)\) for
\(g^{-2}\le C_*^{-2}W^{-1}\), then \(W_0\) to absorb the \(W^{-2}\)
terms. This proves \eqref{uvg:eq:positive}, with \(a_0=a_*/2\).
Boundedness and the derivative bound follow from the fixed
polynomials, the preceding derivative estimates, and
\eqref{uvg:eq:global} on the cutoff.
\end{proof}

Set $Y=c(X,V)-\mathcal Q(X')+\mathcal Q(X)$.

\begin{lemma}[One-step exponential contraction]\label{uvg:lem:mgf}
There exist fixed \(s,c>0\) such that, uniformly in the state and in the range of
Lemma~\ref{uvg:lem:positive},
\(\E(e^{-sY}\mid X)\le\exp(-cg^{-2}/W)\).
\end{lemma}

\begin{proof}[Proof of Lemma~\ref{uvg:lem:mgf}]
Fix the incoming state $X$.

On the convex set \(\mathcal K_g=\{\norm V\le c_0g\}\),
\(Y\) is \(C/g\)-Lipschitz with the same \(C\) for every state. Set
\[
 \widetilde Y(V)\coloneqq
 \inf_{H\in\mathcal K_g}\{Y(H)+(C/g)\|V-H\|_{\rm HS}\}.
\]
This extension has the same Lipschitz constant and is measurable in
the incoming state.
At zero noise the uncorrected squared moduli do not change,
so \(Y(0)=O(g^{-2})\), whence
\[
 |\widetilde Y(V)|\le Cg^{-2}+C\norm V_{\HS}/g.
\]
Gaussian concentration gives
\(\E e^{z(\widetilde Y-\E\widetilde Y)}
\le\exp(Cz^2g^{-2}/W)\).

By \eqref{uvg:eq:tail},
\(\Prob(\mathcal K_g^c)\le2e^{-cWg^2}\), while
\(|Y|\le C+\log(1+C\norm V/g)\).
For every fixed \(|z|\le1\), the exponential moments of
\(|Y|,|\widetilde Y|\) of order \(4|z|\) are uniformly bounded.
For the extension use \(\E\norm V_{\HS}\le\sqrt W\),
\(g\ge C_*\sqrt W\), and Gaussian concentration of the norm.
For \(Y\) use its logarithmic bound and fixed norm moments.
Cauchy--Schwarz on the exceptional set yields
\[
 |\E Y-\E\widetilde Y|\le Ce^{-cWg^2},\qquad
 \E e^{z(Y-\E Y)}
 \le\exp(Cz^2g^{-2}/W)+Ce^{-cWg^2}.
\]
Reducing \(c\) absorbs the centering error. Since
\(\E Y\ge a_0g^{-2}/W\), choose fixed \(s>0\) with
\(Cs^2\le a_0s/4\). Increase \(W_0\) so that
\(Ce^{-cWg^2}\le a_0s\,g^{-2}/(4W)\) uniformly over
\(g\ge C_*\sqrt W\). This uses the decrease of \(u e^{-cu}\)
for large \(u=Wg^2\). Combining the inequalities proves the claim.
\end{proof}

\begin{proof}[Proof of Theorem~\ref{uvg:thm:transfer}]
For \(\ell_n=\sum_{i=1}^n c_i\), iterate
Lemma~\ref{uvg:lem:mgf} and telescope the bounded corrector:
\begin{equation}\label{uvg:eq:scalarfraction}
 \E e^{-s\ell_n}\le Ce^{-cn/(Wg^2)}.
\end{equation}
For $0<s<1$, \eqref{lg:eq:scalar-reduction} and subadditivity give
$\|D_n(\Phi;E)^{-1}\|^s\le C\sum_{a=1}^W e^{-s\ell_n(C_0,e_a)}$.
Applying \eqref{uvg:eq:scalarfraction} to each term proves the bound
for all $n$. For an independent random frame, condition on the frame.

\end{proof}

\subsection{Higher-order correction near zero energy}
\label{lg:sec:drift-to-growth}

Fix $0<\epsilon<1$ and $W^\epsilon\le g<W$.

\begin{lemma}[Transfer growth near zero]\label{up:lem:center-transfer}
Fix $0<\epsilon<1$, $\kappa_0\in(0,2)$, and $D'>0$. There exist $0<\kappa_1<(2-\kappa_0)/2$ and constants
$C_1,W_0$ such that, for every deterministic
$W^\epsilon\le g<W$ and $W\ge W_0$,
\eqref{up:eq:input} holds with
\[
 J=[-2\kappa_1,2\kappa_1],\qquad \mathfrak b=Wg^2,\qquad
 n_0=\lceil C_1Wg^2\log W\rceil.
\]
One may take any sufficiently large fixed $C_1$.
The window and the lower threshold for $C_1$ depend only on
$\epsilon,\kappa_0,D'$; the width threshold may also depend on $C_1$.
The length is independent of $E,\Phi$
and the potentials, and the suprema in \eqref{up:eq:input} remain
outside the probability.
\end{lemma}

Fix $D>0$ and a polynomial horizon $W^{B_1}$ with $B_1>3$.
Put
\begin{equation}\label{lg:eq:initial-time}
 n_*(g,W):=\lceil C_*g^2\log W\rceil,
\end{equation}
where $C_*$ is a fixed constant. For $n\ge n_*$,
\begin{equation}\label{lg:eq:drift-sum}
 \ell_n-\ell_{n_*}=\sum_{i=n_*+1}^n\E_{i-1}c_i
       +\sum_{i=n_*+1}^n(c_i-\E_{i-1}c_i).
\end{equation}

For a deterministic function $Q_g(C,v)$, telescoping its increments gives
\begin{equation}\label{lg:eq:corrected-growth}
 \begin{aligned}
 \ell_n-\ell_{n_*}={}&Q_g(C_n,v_n)-Q_g(C_{n_*},v_{n_*}) +\sum_{i=n_*+1}^n\E_{i-1}
       [c_i-Q_g(C_i,v_i)+Q_g(C_{i-1},v_{i-1})]\\
 &+\sum_{i=n_*+1}^n(1-\E_{i-1})
       [c_i-Q_g(C_i,v_i)].
 \end{aligned}
\end{equation}

\begin{lemma}[Initial interval and martingale fluctuations]\label{lg:lem:fluctuations}
Fix $0<\epsilon<1$, $D>0$, $B_1>3$, and a deterministic family of real
functions $Q_g=Q_{g,\theta,W}$ on $U(W)\times\{v:\|v\|=1\}$ satisfying,
with a common constant and for real tangent directions,
\begin{equation}\label{lg:eq:Q-regularity}
 \begin{aligned}
 \|Q_g\|_\infty&\le Cg^{-2},\\
 |(\dot C\cdot\nabla_C+\dot v\cdot\nabla_v)Q_g(C,v)|
 &\le Cg^{-2}(\|\dot C\|_{\rm op}+\|\dot v\|_2).
 \end{aligned}
\end{equation}
For all sufficiently large $W$ the bounds below hold uniformly for
$W^\epsilon\le g<W$, $|\theta|\le\delta$ with $\delta>0$ fixed sufficiently small,
and all deterministic initial states.
\begin{enumerate}
\item For any fixed $C_*>0$ and $n_*$ as in \eqref{lg:eq:initial-time},
\begin{equation}\label{lg:eq:initial-segment}
 \Prob\{|\ell_{n_*}|>C\log W\}\le W^{-D}.
\end{equation}
\item For each deterministic $n_*\le n\le W^{B_1}$, with probability
at least $1-W^{-D}$,
\begin{equation}\label{lg:eq:centerMGF}
 \max_{n_*\le j\le n}\left|
       \sum_{i=n_*+1}^j(1-\E_{i-1})[c_i-Q_g(C_i,v_i)]\right|
       \le C\sqrt{ng^{-2}W^{-1}\log W}.
\end{equation}
\end{enumerate}
For a fixed constant in \eqref{lg:eq:Q-regularity}, the constants
are independent of $n$ and the choice of $Q_g$.
\end{lemma}

\begin{proof}[Proof of Lemma~\ref{lg:lem:fluctuations}]
The expansion \eqref{lg:eq:c2}, together with $|m_1|,|f_1|\le1$
and bounded $h_\theta$, gives, uniformly in the incoming state,
\begin{equation}\label{lg:eq:coarse-drift}
 |\E_{i-1}c_i|\le Cg^{-2}.
\end{equation}

\noindent\textbf{The martingale and initial interval.}
Fix the incoming state.
On the cutoff $\|V_i\|\le K$, the derivative formulas \eqref{lg:eq:gradients}
and the inverse bound \eqref{uvg:eq:global} give
\[
 \|\nabla_{V_i}c_i\|_\HS\le Cg^{-1},\qquad
 \|dC_i\|_{\rm op}+\|dv_i\|_2\le Cg^{-1}\|\dot V_i\|_\HS.
\]
The chain rule and \eqref{lg:eq:Q-regularity} give
$\|\nabla_{V_i}Q_g(C_i,v_i)\|_\HS\le Cg^{-3}$ on the same cutoff,
and hence $\|\nabla_{V_i}(c_i-Q_g(C_i,v_i))\|_\HS\le Cg^{-1}$.
Extend this function from the cutoff set with the same Lipschitz
constant. After conditional centering, Gaussian concentration gives
exponential moment bounds with variance parameter $Cg^{-2}/W$,
since the covariance of the normalized GUE is $W^{-1}$.
Summing the extended increments from $i=n_*+1$ and applying
the exponential supermartingale inequality at threshold
$C\sqrt{ng^{-2}W^{-1}\log W}$ gives failure probability at most
$\tfrac12W^{-D}$, with $C$ chosen depending on $D$.

To return to the true increments, the cutoff complement has
probability $Ce^{-cW}$, while
$|c_i|\le\log(1+2h_\theta g^{-1}\|V_i\|)$ and
$\|Q_g\|_\infty\le Cg^{-2}$.
The extension is bounded by $Cg^{-2}+C\|V_i\|_\HS/g$, so its
conditional $L^2$ norm is at most $C(g^{-2}+\sqrt W/g)$.
Cauchy--Schwarz bounds each conditional-mean error by this polynomial
factor times an exponential tail, hence by $Ce^{-cW}$ after decreasing $c$.
Summing the conditional-mean errors and taking a union bound over
at most $n$ inputs bounds each contribution by $Cne^{-cW}$. Since $n\le W^{B_1}$ and $g<W$, these
errors are absorbed into the bound in \eqref{lg:eq:centerMGF}
and the remaining $\tfrac12W^{-D}$ failure probability for large $W$.
This proves \eqref{lg:eq:centerMGF}.

For the initial interval up to $n_*$ in \eqref{lg:eq:initial-time},
combine \eqref{lg:eq:coarse-drift} with the same concentration argument
applied to $\sum_{i=1}^{n_*}(1-\E_{i-1})c_i$.
Using $\ell_{n_*}=\sum_{i=1}^{n_*}c_i$ from \eqref{lg:eq:additive-growth} gives
\[
 |\ell_{n_*}|\le Cg^{-2}n_*
       +C\sqrt{n_*g^{-2}W^{-1}\log W}\le C\log W
\]
with probability at least $1-W^{-D}$, proving
\eqref{lg:eq:initial-segment}.
\end{proof}

\subsubsection{First-order moment corrections at zero energy}
\label{lg:sec:state-estimates}
\label{lg:sec:concrete-Q}

Set $E=0$ and $t=(2g)^{-1}$. Here $M_j,F_j$ denote the
first corrections $M_j^{[1]},F_j^{[1]}$ of Section~\ref{sec:intro-direct-transfer}.
The free update and Gaussian symmetry give, for each fixed $k$,
\begin{equation}\label{lg:eq:example-leading-updates}
 P_gm_k=(-1)^km_k+O(g^{-2}),\qquad P_gf_k=(-1)^kf_k+O(g^{-2}).
\end{equation}

\noindent\textbf{First trace and directional moments.}
Take $M_1,F_1$ from \eqref{uvg:eq:M}--\eqref{uvg:eq:F} at $\theta=0$.
Equations \eqref{uvg:eq:mmean} and \eqref{uvg:eq:means} give
\begin{align}
 P_gm_1&=-m_1+t^2(1+m_1)(1+2m_1+m_2)+O(t^4),
 \label{lg:eq:example-raw-m}\\
 P_gM_1&=-[1-t^2(3+m_2)]M_1+O(t^4).
 \label{lg:eq:example-M-cancellation}
\end{align}
For $M_1$ the order-$t^2$ coupling to other moments is zero; for $F_1$
it is $O(t^2|M_1|)$ by \eqref{uvg:eq:means}.
For each of the fixed trace or directional polynomials $Z$ below,
\eqref{lg:eq:gradients} bounds its gradient in the new GUE block by
$C g^{-1}(1+g^{-1}\|V\|)^C$, with an additional $W^{-1/2}$ for a
normalized trace. Gaussian $L^{2p}$ Poincar\'e and \eqref{uvg:eq:tail}
give, conditionally on the incoming state, for fixed integers $p\ge1$,
\begin{equation}\label{lg:eq:example-fluctuations}
 \|(1-\E_{i-1})Z(C_i,v_i)\|_{2p,i-1}\le C_pg^{-1}
 \begin{cases}
  W^{-1},& Z\text{ involves only normalized traces},\\
  W^{-1/2},& Z\text{ also involves directional moments}.
 \end{cases}
\end{equation}
Expanding the $2p$th power around the conditional mean eliminates the
linear fluctuation. Young's inequality absorbs a fraction of the
$cg^{-2}|M_1|^{2p}$ damping from \eqref{lg:eq:example-M-cancellation},
leaving
\begin{equation}\label{lg:eq:example-moment-M}
 P_g|M_1|^{2p}\le(1-c_pg^{-2})|M_1|^{2p}
            +C_pg^{-2}(W^{-2p}+g^{-4p}).
\end{equation}
For $F_1$ the same expansion gives
\begin{equation}\label{lg:eq:example-moment-F}
 P_g|F_1|^{2p}\le(1-c_pg^{-2})|F_1|^{2p}
  +C_pg^{-2}\bigl(|M_1|^{2p}+W^{-p}+g^{-4p}\bigr).
\end{equation}
Iterating first for $M_1$ and then for $F_1$ gives, for $n_*\le i\le W^{B_1}$,
\[
 \|M_1(C_i)\|_{2p}\le C_p(W^{-1}+g^{-2}),\qquad
 \|F_1(C_i,v_i)\|_{2p}\le C_p(W^{-1/2}+g^{-2}).
\]
Here and below $C_*$ in \eqref{lg:eq:initial-time} is chosen for $p$:
the initial terms are bounded by $C_p(1+i/g^2)^2e^{-c_pi/g^2}$ before
taking the $2p$th root.

\noindent\textbf{Second trace moment.}
Write $H=I+C$ and $\operatorname{tr}_W(B)=W^{-1}\Tr B$.
In \eqref{lg:eq:matrix-update}, $C'=-C-itHVH+t^2HVHVH+O(t^3\|V\|^3)$
on the Taylor cutoff. The square of the linear term and the two cross
terms with $-C$ give
\[
 (P_g-I)m_2=-t^2\bigl[\operatorname{tr}_W(H^2)^2+
                  2\operatorname{tr}_W(H)\operatorname{tr}_W(CH^2)\bigr]+O(t^4)
\]
by $\E VBV=\operatorname{tr}_W(B)I$ and Gaussian symmetry. The odd part
of the bracket is $2(3m_1+m_3+4m_1m_2)$. Set
$M_2=m_2-t^2(3m_1+m_3+4m_1m_2)$ and $a=-3+\sqrt8$. Since
$1+6a+a^2=0$, \eqref{lg:eq:example-leading-updates} yields
\[
 P_g(M_2-a)=[1-t^2(6+m_2+a)](M_2-a)+O(t^2|M_1|+t^4).
\]
Here the remaining terms $6m_1^2+2m_1m_3$ are $O(|M_1|+t^2)$.
The multiplier has modulus at most $1-cg^{-2}$. Using
\eqref{lg:eq:example-fluctuations} and iterating as for $M_1$ gives
$\|M_2(C_i)-a\|_{2p}\le C_p(W^{-1}+g^{-2})$ for $i\ge n_*$.

\noindent\textbf{Second directional moment.}
Let $A=(I-itVH)^{-1}$. Expand $v^*A(C')^2A^*v/(v^*AA^*v)$,
contracting each pair by $\E VBV=\operatorname{tr}_W(B)I$.
At order $t^2$, contractions between two separate quadratic forms
carry $W^{-1}$. This gives
\begin{equation}\label{lg:eq:example-raw-f}
 \begin{aligned}
 P_gf_2={}&f_2+\tfrac14g^{-2}\bigl[
 -2(3+m_2)f_2-(m_3+6m_1)f_1-3m_1f_3\\
 &\hspace{12mm}+f_2(m_1f_1+m_{-1}f_{-1})
 -(3+4m_2)f_1-3f_3\\
 &\hspace{12mm}-8m_1f_2+(f_1+f_{-1})f_2\bigr]
 +O(g^{-4}+W^{-1}g^{-2}).
 \end{aligned}
\end{equation}

The free update changes the sign of each term in the bracket defining
\begin{equation}\label{lg:eq:example-corrections}
 F_2\coloneqq f_2-\tfrac12 t^2
       [(3+4m_2)f_1+3f_3+8m_1f_2-(f_1+f_{-1})f_2].
\end{equation}
Their increments cancel the odd terms in \eqref{lg:eq:example-raw-f}, leaving
\begin{equation}\label{lg:eq:example-updates}
 P_gF_2=[1-2t^2(3+m_2)]F_2
       +O(g^{-2}(|M_1|+|F_1|)+g^{-4}+W^{-1}g^{-2}).
\end{equation}
Using \eqref{lg:eq:example-fluctuations} as for $F_1$ gives
\begin{equation}\label{lg:eq:first-F2-moment}
 P_g|F_2|^{2p}\le(1-c_pg^{-2})|F_2|^{2p}
  +C_pg^{-2}\bigl(|M_1|^{2p}+|F_1|^{2p}+W^{-p}+g^{-4p}\bigr).
\end{equation}
Iteration gives
$\|F_2(C_i,v_i)\|_{2p}\le C_p(W^{-1/2}+g^{-2})$ for $i\ge n_*$.
In particular all four first corrections have the precision used in
Section~\ref{sec:intro-direct-transfer} when $g\asymp W^\epsilon$,
$0<\epsilon<1/4$.

\subsubsection{Wick expansion and the first even-mode correction}
\begin{lemma}[Uniform polynomial expansions]\label{lg:lem:ideal}
For fixed $0<\epsilon<1$, $W^\epsilon\le g<W$, and $|\theta|\le\delta$ with $\delta>0$ fixed sufficiently small,
let $p=p(m,f,W^{-1})$ be a fixed polynomial in finitely many
$m_k,f_k$ and in nonnegative powers of $W^{-1}$.
For every fixed integer $s\ge0$,
\[
 P_gp=\sum_{r=0}^{s}g^{-2r}\sum_{\ell\ge0}W^{-\ell}
             p_{r,\ell}(m,f;\theta)+O_{p,s}(g^{-2s-2}).
\]
Each $p_{r,\ell}$ is a polynomial, and only finitely many are
nonzero. Their coefficients and the remainder bound are uniform
in $W$, the energy, and the incoming state
$C\in U(W)$, $\|v\|_2=1$. For $0\le r\le s$, the coefficients
have the following properties.
\begin{enumerate}
\item[(i)] Let $\mathscr J\coloneqq(W^{-1},f_k:k\ne0)$ be the ideal generated
by $W^{-1}$ and the $f_k$ with $k\ne0$. Then
\[
 p\in\mathscr J\quad\Longrightarrow\quad
 \sum_{\ell\ge0}W^{-\ell}p_{r,\ell}\in\mathscr J.
\]
In particular, $p_{r,0}(m,0;\theta)=0$ for $p\in\mathscr J$:
each term with no explicit $W^{-1}$ retains at least one
nonzero-index $f_k$ factor.
\item[(ii)] For any two such polynomials $p,\widetilde p$,
\[
 (p\widetilde p)_{r,0}
   =\sum_{a+b=r}p_{a,0}\,\widetilde p_{b,0}.
\]
Here $(p\widetilde p)_{r,0}$ denotes the coefficient of
$g^{-2r}W^0$ in $P_g(p\widetilde p)$, with the analogous
notation for the two factors. Thus the formal $W^0$ expansion
preserves multiplication.
\item[(iii)] If $p$ depends only on the $m_k$'s, then
\[
 p_{r,\ell}=p_{r,\ell}(m;\theta),\qquad p_{r,1}=0.
\]
Thus these expansions contain only traces, and their first
possible finite-width correction has a factor $W^{-2}$.
\end{enumerate}
\end{lemma}
\begin{proof}
Isotropy and the shear/inverse-shear bounds give
$\|F(C_{i-1},V_i)^{-1}\|\le1+2h_\theta g^{-1}\|V_i\|$.
On $h_\theta g^{-1}\|V_i\|\le c_0$, expand the inverses and normalized vector
in \eqref{lg:eq:update}. Fixed-order derivatives are bounded by
$C\|V_i\|^a$, uniformly in the state and width; each coefficient at
$g^{-2r}$ increases degree by at most $2r$.
Gaussian symmetry eliminates the odd powers of $V_i$ after
expectation. A sphere-net estimate for GUE quadratic forms gives
uniform operator-norm moments of every fixed order. On the complement of the Taylor event, these moments and
the boundedness of the observables give any requested power of
$g^{-1}$. This proves the remainder estimate uniformly in the incoming state.

To verify the assertions about the $f_k$ factors, write $\alpha,\beta,\gamma,\delta_0$ for the blocks of the unrotated
shear at step $i$, so $F(C_{i-1},V_i)=\alpha+\beta C_{i-1}$.
For an auxiliary $w$, set
$\alpha_w\coloneqq \alpha-w\gamma$, $\beta_w\coloneqq \beta-w\delta_0$,
$K_w\coloneqq \alpha_w^{-1}\beta_w$,
$\mathcal S_w\coloneqq (I+K_wC_{i-1})^{-1}$.
The $J_{\rm sig}$-unitary block identities give exactly
\begin{equation}\label{lg:eq:poisson}
 [F(C_{i-1},V_i)-w(\gamma+\delta_0C_{i-1})]^{-1}
 F(C_{i-1},V_i)^{-*}=\mathcal S_w+\mathcal S_0^*-I.
\end{equation}
Indeed, $\alpha_w\alpha^*-\beta_w\beta^*=I$ and multiplication
reduces the claim to
$\mathcal S_w(I-K_wK_0^*)\mathcal S_0^*=\mathcal S_w+\mathcal S_0^*-I$.
If $Y\coloneqq\|F(C_{i-1},V_i)^{-*}v_{i-1}\|^2$, then
$Y=v_{i-1}^*(\mathcal S_0+\mathcal S_0^*-I)v_{i-1}$ and
\[
 v_i^*(I+we^{2i\theta}C_i)^{-1}v_i
   =1+v_{i-1}^*(\mathcal S_w-\mathcal S_0)v_{i-1}/Y.
\]
Expanding $\mathcal S_w$ shows that every initially marked word has a
nonzero same-sign $C_{i-1}$-power adjacent to its rank-one factor
$v_{i-1}v_{i-1}^*$;
negative powers are obtained by adjunction.

Contract each Gaussian pair using
$\E(V_i)_{ab}(V_i)_{cd}=W^{-1}\delta_{ad}\delta_{bc}$.
With $r$ pairs, $t$ initially normalized unmarked cycles and $u$
final unmarked cycles, the coefficient carries $W^{-(r+t-u)}$.
Each split increases $u$ by at most one, and a join decreases it.
Thus $r+t-u\ge0$. At exponent zero only splits creating unmarked
cycles are possible. Their cuts occur at $V_i$-letters, so the nonzero
same-sign $C_{i-1}$-power adjacent to $v_{i-1}v_{i-1}^*$ stays in its
marked cycle. With no joins, opposite-sign powers from other marked
cycles cannot cancel it, and a factor $f_k$ with $k\ne0$ survives.
Starting from an input
containing a nonzero-index $f_j$, a term with no remaining $f_k$
factor therefore carries $W^{-1}$. When the input contains only normalized traces,
$r+t-u$ is twice the number of joins, giving $W^{-2}$ as the first
width correction. Joins of two marked cycles cost $W^{-1}$, so
in the coefficient of $W^0$, contractions act separately on distinct
initial factors. This proves properties (i)--(iii); only finitely
many words and Gaussian pairings occur at each fixed order.
\end{proof}

\medskip\noindent\textbf{The linear $f_2$ coefficient.}
At $E=0$ set $t=(2g)^{-1}$ and $Q^{[1]}=-\tfrac12 t^2\Re f_1$.
We compute the $W^0$ coefficient linear in $f_2$ in
$P_gc-(P_g-I)Q^{[1]}$, in the convention of Section~\ref{sec:intro-direct-transfer}.
Let $w$ be a formal variable and use $\mathcal S_w$ from \eqref{lg:eq:poisson}.
Write $\E_0$ for the $W^0$ part of the Wick expansion,
$H=I+C$ and $H_w=H(I-wC)^{-1}$. Through any fixed Taylor order,
\[
 \E_0\mathcal S_w=(I+d_w C)^{-1},\qquad
 d_w=(s_w-w)/(1+s_w),\qquad
 s_w=t^2(1+w)^2\operatorname{tr}_W\!\left(H_w(I+s_w H_w)^{-1}\right).
\]
Indeed, putting $u=t(1+w)$ gives
$\mathcal S_w=(I-wC)^{-1}(I-iuVH_w)^{-1}(I-iuV)$.
For $K=\E_0(I-iuVH_w)^{-1}$, pairing the first $V$ with its partner gives
$K=I-u^2\operatorname{tr}_W(H_wK)H_wK$.
Only noncrossing pairs contribute by the splitting rule in
Lemma~\ref{lg:lem:ideal}. Applying the same pairing to the factor
$I-iuV$ on the right gives
$\E_0\mathcal S_w=(I-wC)^{-1}K(1+s_w)$, as claimed.

Set $d=d_0$ and $d'=\partial_w d_w|_{w=0}$.
In $Y=v^*(\mathcal S_0+\mathcal S_0^*-I)v$ the constant term in the
directional moments is one. Products factor at $W^0$, so the coefficient
of $\Re f_2$ in $-\tfrac12\E_0\log Y$ at real trace centers is $-d^2$.
For $P_g f_1$, extract $-[w]$ in the generating function following
\eqref{lg:eq:poisson}. The coefficient of $f_2$ is $-2dd'$ because
the numerator vanishes at $w=0$ and has no constant directional term.
Thus the two contributions in the corrected drift sum to $-d(d+t^2d')$.
Writing $a=-3+\sqrt8$, the trace substitution needed to this order is
$m_2=a+O(t^2)$ and $m_1=\tfrac12 t^2(1+a)+O(t^4)$: the corrected
odd trace has center zero and $M_1=m_1-\tfrac12t^2(1+m_2+2m_1^2)$
through its first correction. From the equation for $s_w$,
\[
 \begin{aligned}
 d&=t^2(1+m_1)-t^4\bigl[\operatorname{tr}_W(H)\operatorname{tr}_W(H^2)+
             \operatorname{tr}_W(H)^2\bigr]+O(t^6)
   =t^2-\tfrac12(3+a)t^4+O(t^6),\\
 d'&=-1+(3+a)t^2+O(t^4).
 \end{aligned}
\]
In particular the fourth-order constant cancels, and the sixth-order
coefficient and the corresponding term of the correction are
\begin{equation}\label{lg:eq:first-even-drift}
 -d(d+t^2d')=-\tfrac12(3+a)t^6+O(t^8),\qquad
 \frac{-\tfrac12(3+a)t^6}{-2(3+a)t^2}=\tfrac14 t^4.
\end{equation}
The divisor is the $f_2$ coefficient in \eqref{lg:eq:example-raw-f}
at the trace centers. This verifies the $\Re f_2$ term of
$Q^{[2]}-Q^{[1]}$.

\subsubsection{Higher-order moment estimates}
\label{lg:sec:general-state}

\begin{lemma}[Stability of $M_k$ and $F_k$]\label{lg:lem:dynamic}
Fix $0<\epsilon<1$, $D,\eta>0$, $B_1>3$, an integer $R\ge2$ with
$2\epsilon R>1$, and finitely many nonzero indices of the moments
\eqref{lg:eq:observables}. There exist $\delta,C_*>0$ such that,
for all sufficiently large $W$ and $W^\epsilon\le g<W$, the moments admit
polynomial corrections of the form
\[
 \begin{aligned}
 M_k&=m_k+\sum_{r=1}^{R+2}g^{-2r}\mathcal M_{k,r}(m;\theta), \quad 
 F_k =f_k+\sum_{r=1}^{R+2}g^{-2r}\mathcal F_{k,r}(m,f;\theta).
 \end{aligned}
\]
Here $\mathcal M_{k,r}$ is a polynomial in finitely many $m_j$'s,
and $\mathcal F_{k,r}$ is a polynomial in finitely many $m_j$'s and
$f_j$'s. Every monomial of $\mathcal F_{k,r}$ contains some $f_j$
with $j\ne0$.
The order $R$, polynomial degrees and index sets are independent of $W$.
The coefficient bounds are uniform on $|\theta|\le\delta$, with $\delta$ chosen for the fixed indices and order. In particular, $M_k=m_k+O(g^{-2})$ and
$F_k=f_k+O(g^{-2})$ uniformly in the state.
There are also bounded deterministic centers $a_k=a_k(g,\theta)$
such that, for each deterministic initial state and $|\theta|\le\delta$,
for all sufficiently large $W$, with probability at least
$1-W^{-D}$,
\begin{equation}\label{lg:eq:burnevent}
 |M_k(C_n)-a_k|\le W^{\eta-1},\qquad
 |F_k(C_n,v_n)|\le W^{\eta-1/2},
 \qquad n_*\le n\le W^{B_1},
\end{equation}
simultaneously for these indices and times. Here
$n_*$ is given by \eqref{lg:eq:initial-time}. The centers satisfy
\[
 a_{2j+1}=0,\qquad |a_2|<1/4+O(g^{-2}),\qquad
 a_{-k}=\overline{a_k}.
\]
\end{lemma}

\begin{proof}[Proof of Lemma~\ref{lg:lem:dynamic}]
\noindent\textbf{Polynomial corrections.}
For this fixed $R$, starting from the required indices,
include those needed by the polynomial corrections, inverse
substitutions, center corrections, and fixed-depth moment iterations
below, and by the finite Neumann construction and its final residual
in Section~\ref{lg:sec:cocycle-correction}.
Close each index set under lower nonzero indices. Only finitely many
indices are added at each step; the triangular inverses use only lower
indices. The resulting finite set is independent of $W$.

Write $[H]_{\le s,0}$
for the terms of orders $1,g^{-2},\ldots,g^{-2s}$ in $H$ with no
explicit factor $W^{-1}$, expressed in the coordinates used at that
step. For a monomial, set
\[
 p=\prod_j M_j^{\alpha_j}F_j^{\beta_j},\qquad
 \deg p=\sum_j j(\alpha_j+\beta_j),\qquad
 \alpha_j,\beta_j\in\mathbb Z_{\ge0}.
\]
We require
\begin{equation}\label{lg:eq:required-parity-updates}
 \begin{aligned}
 [P_gM_k]_{\le R+2,0}
 &=(-1)^ke^{-2ik\theta}M_k
   +\sum_{r=1}^{R+2}g^{-2r}
     \sum_{\deg p\equiv k\ (2)}c^M_{k,r,p}\,p(M),\\
 [P_gF_k]_{\le R+2,0}
 &=(-1)^ke^{-2ik\theta}F_k
   +\sum_{r=1}^{R+2}g^{-2r}
     \sum_{\deg p\equiv k\ (2)}c^F_{k,r,p}\,p(M,F).
 \end{aligned}
\end{equation}
Here $p(M)$ contains only trace moments.

Compute the coefficients by Lemma~\ref{lg:lem:ideal}, keeping the
phases at the specified energy. Initialize the recursion by
\[
 M_k^{(0)}=m_k,\qquad F_k^{(0)}=f_k\qquad(k\ne0),
 \qquad M_0^{(s)}=F_0^{(s)}=1.
\]
The superscript $s$ denotes the correction step.
At step $s$, $1\le s\le R+2$, the already constructed functions have
the form
\[
 M^{(s-1)}=m+\Delta^M(m),\qquad
 F^{(s-1)}=f+\Delta^F(m,f),\qquad
 \Delta^M,\Delta^F=O(g^{-2}).
\]
Here $m,f,M^{(s-1)},F^{(s-1)}$ denote the collections of coordinates,
and $\Delta^M,\Delta^F$ are the known correction polynomials at this
step. To express $m,f$ in the current coordinates, start with
$m^{[0]}=M^{(s-1)}$, $f^{[0]}=F^{(s-1)}$ and substitute recursively:
\begin{equation}\label{lg:eq:coordinate-reversion}
 \begin{aligned}
 m^{[r+1]}&=M^{(s-1)}-\Delta^M(m^{[r]}),\\
 f^{[r+1]}&=F^{(s-1)}-\Delta^F(m^{[r]},f^{[r]})
 \qquad(0\le r<R+2).
 \end{aligned}
\end{equation}
Truncate after $g^{-2R-4}$ at each substitution. Since the correction
polynomials start at $g^{-2}$, each substitution determines one more
order: $m=m^{[r]}+O(g^{-2r-2})$ and
$f=f^{[r]}+O(g^{-2r-2})$. Thus $r=R+2$ gives the required inverse substitution.

Now compute the conditional updates in $m,f$ by the expansion rule
above, substitute $m^{[R+2]},f^{[R+2]}$, and collect powers of
$g^{-2}$. For this fixed step $s$, this defines the polynomials
$B^M_{k,r},B^F_{k,r}$ by
\begin{equation}\label{lg:eq:current-coordinate-updates}
 \begin{aligned}
 [P_gM_k^{(s-1)}]_{\le R+2,0}
 &=(-1)^ke^{-2ik\theta}M_k^{(s-1)}
   +\sum_{r=1}^{R+2}g^{-2r}B^M_{k,r}(M^{(s-1)}),\\
 [P_gF_k^{(s-1)}]_{\le R+2,0}
 &=(-1)^ke^{-2ik\theta}F_k^{(s-1)}
   +\sum_{r=1}^{R+2}g^{-2r}B^F_{k,r}(M^{(s-1)},F^{(s-1)}).
 \end{aligned}
\end{equation}
At order $g^0$, the conditional update in the current coordinates multiplies
a monomial $p$ of degree $d$ by $(-1)^d e^{-2id\theta}$.
Corrections to this rule start at $g^{-2}$; after multiplication
by $g^{-2s}$, they affect only orders $g^{-2s-2}$ and higher.
For each monomial $p$, take its coefficients in these two polynomials:
\[
 b^M_{k,s,p}\coloneqq[p]B^M_{k,s},\qquad
 b^F_{k,s,p}\coloneqq[p]B^F_{k,s},
\]
where $[p]$ denotes the coefficient of $p$ in the current coordinates.
For the finitely many opposite-parity pairs $(k,d)$ occurring in
these corrections, choose a fixed $\delta>0$ such that
$|k-d|\delta\le\pi/3$ for every pair. Then
\begin{equation}\label{lg:eq:opposite-parity-gap}
 \inf_{|\theta|\le\delta}
 \left|(-1)^ke^{-2ik\theta}-(-1)^de^{-2id\theta}\right|
 =\inf_{|\theta|\le\delta}2|\cos((k-d)\theta)|\ge1.
\end{equation}
Both the window and this lower bound are independent of $W$.
For every monomial whose degree has parity opposite to $k$, make
the simultaneous updates
\begin{equation}\label{lg:eq:observable-correction}
 \begin{aligned}
 M_k^{(s)}&=M_k^{(s-1)}+g^{-2s}
 \sum_{\deg p\not\equiv k\ (2)}
 \frac{b^M_{k,s,p}\,p}
 {(-1)^ke^{-2ik\theta}-(-1)^{\deg p}e^{-2i(\deg p)\theta}},\\
 F_k^{(s)}&=F_k^{(s-1)}+g^{-2s}
 \sum_{\deg p\not\equiv k\ (2)}
 \frac{b^F_{k,s,p}\,p}
 {(-1)^ke^{-2ik\theta}-(-1)^{\deg p}e^{-2i(\deg p)\theta}}.
 \end{aligned}
\end{equation}
Here $p$ is evaluated at the coordinates from step $s-1$; each sum
contains only monomials with nonzero coefficients in that update.
Expand back in the original $m,f$ and discard powers beyond
$g^{-2R-4}$ after each step. Make no correction to same-parity terms.
Define negative-index corrections by conjugation.
Finally, the exact definitions are
\begin{equation}\label{lg:eq:corrected-observables}
 M_k\coloneqq M_k^{(R+2)},\qquad F_k\coloneqq F_k^{(R+2)}.
\end{equation}
For a term $b\,p$ with $\deg p=d$ of parity opposite to $k$, the
increment of the added polynomial in \eqref{lg:eq:observable-correction}
multiplies its coefficient by $(-1)^d e^{-2id\theta}-(-1)^k e^{-2ik\theta}$,
cancelling $b$ at order $g^{-2s}$.

By Lemma~\ref{lg:lem:ideal}, the part of the transition with no
explicit $W^{-1}$ preserves multiplication, polynomials in the
traces, and the ideal generated by the nonzero-index $f_k$'s.
The successive substitutions preserve these properties, and induction
in $s$ gives \eqref{lg:eq:required-parity-updates}.
Applying \eqref{lg:eq:coordinate-reversion} to these final
corrections gives the inverse polynomials $m_k^{[R+2]}(M)$ and
$f_k^{[R+2]}(M,F)$, with
uniform remainder $O(g^{-2R-6})$. Each inverse substitution
preserves the property that every monomial in $f_k^{[R+2]}$
contains a nonzero-index $F_j$.

\noindent\textbf{The leading terms of the corrected transition.}
We now compute the order-$g^{-2}$ terms that preserve parity.
They determine both the trace centers and the contraction used below.
For this coefficient calculation, set
$\mathcal F(w)\coloneqq \sum_{k\ge0}f_k w^k$,
$\mathcal M(w)\coloneqq \sum_{k\ge0}m_{k+1}w^k$.
Let $h_\theta^2\mathcal D_1$ denote the coefficient of $g^{-2}W^0$
in the expected unrotated transition. The shear gives
$K_w=-wI-ih_\theta g^{-1}(1+w)^2V_i+h_\theta^2g^{-2}(1+w)^3V_i^2+O(g^{-3})$.
Expanding \eqref{lg:eq:poisson} yields, for $k\ge1$,
\[
 \mathcal D_1f_k
 =-[w^k](1+w)^3[1+(1+w)\mathcal M(w)]\mathcal F'(w)
   +[(1+m_1)f_1+(1+m_{-1})f_{-1}]f_k.
\]
Retain only monomials whose degree has the same parity as $k$,
and denote the resulting operator by $\mathcal L$. For $k>0$ this
gives the \emph{second-order drift formula for the directional moments},
\begin{equation}\label{lg:eq:center-corner}
 \begin{aligned}
 \mathcal Lf_k={}&-3kf_k-\mathbf1_{k\ge3}(k-2)f_{k-2}
 -\sum_{n=1}^{k+1}nf_nm_{k+2-n}
 -6\sum_{n=1}^{k-1}nf_nm_{k-n}
 -\sum_{n=1}^{k-3}nf_nm_{k-2-n}\\
 &+f_k(m_1f_1+m_{-1}f_{-1}).
 \end{aligned}
\end{equation}
The first two coefficients of the update from $C_{i-1}$ before rotation,
in powers of $h_\theta g^{-1}$, are
$\delta C_1\coloneqq i(I+C_{i-1})V_i(I+C_{i-1})$,
$\delta C_2\coloneqq -(I+C_{i-1})V_i(I+C_{i-1})V_i(I+C_{i-1})$.
Contracting their contributions to $W^{-1}\Tr C_{i-1}^k$ and retaining
monomials of parity $k$ gives
\begin{equation}\label{lg:eq:center-trace}
 \mathcal Lm_k=-3k\,m_k-3k\sum_{a=1}^{k-1}m_am_{k-a}
 -\frac k2\sum_{a=0}^{k-2}m_am_{k-2-a}
 -\frac k2\sum_{a=1}^{k+1}m_am_{k+2-a}.
\end{equation}
Negative indices are conjugates. Empty sums vanish.
The only leading higher index is $k+1$, multiplied by $m_1$.
In particular,
$\mathcal Lm_2=-1-6m_2-m_2^2-6m_1^2-2m_1m_3$,
and $\mathcal Lf_1=-(3+m_2)f_1-2m_1f_2+
 f_1(m_1f_1+m_{-1}f_{-1})$.
The leading retained map is
$y_k\mapsto (-1)^ke^{-2ik\theta}[y_k+h_\theta^2g^{-2}\mathcal Ly_k]$;
the phase is the input phase, not separate phases of output monomials.

\noindent\textbf{Constructing the centers $a_k$.}
The updates of the $M_k$'s involve only the $M_j$'s. Define their
truncated expected increment by
\[
 \mathcal E_k(M)\coloneqq[(P_g-I)M_k]_{\le R+2,0},
 \qquad
 \mathcal E_k^{\rm lead}(M)\coloneqq[(P_g-I)M_k]_{\le1,0},
\]
where both polynomials are expressed in the final coordinates.
We choose $a$ as an approximate zero of this expected increment.
The concentration estimate for $M_k-a_k$ will follow from the
contraction of these updates.

First solve the leading equations
$\mathcal E_{2j}^{\rm lead}(a^{(0)})=0$ in increasing $j$,
with $a_0^{(0)}=1$, $a_{2j+1}^{(0)}=0$, and negative indices given
by conjugation. The first equation is
\begin{equation}\label{lg:eq:a2}
 (a_2^{(0)})^2+d a_2^{(0)}+1=0,\quad
 d\coloneqq6-g^2h_\theta^{-2}(1-e^{4i\theta}),\qquad
 |a_2^{(0)}|<1/4.
\end{equation}
The small root specifies $a_2^{(0)}$ uniquely. For $j\ge2$, the
$2j$th equation is linear in $a_{2j}^{(0)}$, with coefficient
\[
 e^{-4ij\theta}[1-2j h_\theta^2g^{-2}(3+a_2^{(0)})]-1;
\]
all its other nonzero terms involve already determined even indices.
Thus these equations are a recursive definition of every required
$a_{2j}^{(0)}$. At $E=0$ they give
$a_{2j}^{(0)}=(-3+\sqrt8)^j$.

To include the higher orders, let $\mathsf J$ be the Jacobian of
$\mathcal E^{\rm lead}$ at $a^{(0)}$ on the even coordinates.
Starting from $a^{(0)}$, perform $R$ corrections:
\begin{equation}\label{lg:eq:center-recursion}
 a^{(s+1)}=a^{(s)}-\mathsf J^{-1}\mathcal E(a^{(s)})
 \quad(0\le s<R),\qquad a\coloneqq a^{(R)}.
\end{equation}
Throughout the iteration, keep $a_0=1$, set the odd coordinates
to zero, and define the negative coordinates by conjugation. 

\noindent\textbf{Bounds for the trace centers and the inverse matrices.}
The trace formula \eqref{lg:eq:center-trace} gives
\eqref{lg:eq:a2}. For sufficiently large $g$,
$|d|^2=36+(2g^4h_\theta^{-4}-12g^2h_\theta^{-2})(1-\cos4\theta)\ge36$.
The contraction $a\mapsto-1/(d+a)$ on $|a|\le1/4$ proves existence
and uniqueness of $a_2^{(0)}$ uniformly in $g^2\theta$.
At $\theta=0$, $a_{2p}^{(0)}=(-3+\sqrt8)^p$, $a_{2p+1}^{(0)}=0$
solves all leading rows; substitution in \eqref{lg:eq:center-trace}
gives $-p^2(a_2^{(0)})^{p-1}[1+6a_2^{(0)}+(a_2^{(0)})^2]=0$.
The higher even equations are triangular. The diagonal coefficients
used in these equations have the form
\begin{equation}\label{lg:eq:center-denominators}
 d_k\coloneqq (-1)^ke^{-2ik\theta}[1-h_\theta^2g^{-2} k(3+a_2)]-1,\qquad
 |d_k|\ge c_kg^{-2}.
\end{equation}
For the leading trace equations use $a_2^{(0)}$ in place of $a_2$.
For even $k$, $e^{-2ik\theta}-1$ competes with $g^{-2}$, so here
$\theta=o(1)$ does not justify setting $\theta=0$.
The numerators in the higher-center equations have a factor $g^{-2}$;
hence the leading centers remain bounded at each fixed index.
The second-row Jacobian is
$e^{-4i\theta}[1-h_\theta^2g^{-2}(6+2a_2^{(0)})]-1$;
the higher rows are lower triangular with the diagonals just described.
All off-diagonal entries are $O(g^{-2})$, so successive substitution
on each fixed finite set gives $\|\mathsf J^{-1}\|=O(g^2)$.

The initial residual $\mathcal E(a^{(0)})$ is $O(g^{-4})$.
In \eqref{lg:eq:center-recursion}, a residual $O(g^{-2(s+2)})$
gives a correction $O(g^{-2(s+1)})$ and a new residual
$O(g^{-2(s+3)})$: the Jacobian differs from $\mathsf J$ by $O(g^{-4})$
and the leading second derivative is $O(g^{-2})$.
The finite iteration therefore gives
\begin{equation}\label{lg:eq:centers}
 a=a^{(0)}+O(g^{-2}),\quad
 \mathcal E(a)=O(g^{-2R-2}),\quad a_{2k+1}=0,
 \quad |a_2|<1/4+C_Rg^{-2}.
\end{equation}
The construction also gives $|M_k|,|F_k|\le C_{k,R}$ and
$|M_2|\le1+C_Rg^{-2}$. Their conditional means satisfy the corrected updates
in \eqref{lg:eq:required-parity-updates}, with remainders
\begin{equation}\label{lg:eq:widtherrors}
 O(g^{-2} W^{-2}+g^{-2R-2})\quad\hbox{for the }M_k\hbox{ updates},\qquad
 O(g^{-2} W^{-1}+g^{-2R-2})\quad\hbox{for the }F_k\hbox{ updates}.
\end{equation}
The powers of $W^{-1}$ follow from Lemma~\ref{lg:lem:ideal}(i), (iii)
after the truncated inverse substitution. The same lemma gives
a uniform bound on the Taylor remainder.
In the moment estimates, $M_{k,n}=M_k(C_n)$ and
$F_{k,n}=F_k(C_n,v_n)$ denote evaluation at time $n$.

\noindent\textbf{The one-step fluctuations.}
The goal is the conditional innovation bound \eqref{lg:eq:innovations}.
Differentiating $M_k(C_i)$ and $F_k(C_i,v_i)$ by the product rule,
using \eqref{lg:eq:gradients} and the derivative of the normalized
vector $v_i$, gives
\[
 \begin{aligned}
 \|\nabla_{V_i}M_k(C_i)\|_\HS
 &\le Cg^{-1}W^{-1/2}(1+g^{-1}\|V_i\|)^C,\\
 \|\nabla_{V_i}F_k(C_i,v_i)\|_\HS
 &\le Cg^{-1}(1+g^{-1}\|V_i\|)^C.
 \end{aligned}
\]
The Gaussian $L^{2p}$ Poincar\'e inequality therefore gives,
conditionally on the current state,
\begin{equation}\label{lg:eq:innovations}
 \|M_{k,n+1}-\E_nM_{k,n+1}\|_{2p,n}\le C_{k,p}g^{-1}\,W^{-1},
 \quad
 \|F_{k,n+1}-\E_nF_{k,n+1}\|_{2p,n}\le C_{k,p}g^{-1}W^{-1/2}.
\end{equation}
The dimension-independent $L^{2p}$ inequality follows from the
gradient commutation formula for the Gaussian Ornstein--Uhlenbeck
semigroup and H\"older's inequality; its constant depends only on $p$. The covariance contributes
$W^{-1/2}$ and uniform GUE norm moments bound the polynomial factors.

\noindent\textbf{Contraction of the conditional means.}
The goal is the moment recursion \eqref{lg:eq:momentrec}.
The factors $(-1)^ke^{-2ik\theta}$ have modulus one and therefore
do not affect these contraction estimates.
For the first trace row and the higher $M_k,F_k$ rows, the principal
multiplier has modulus $|1-h_\theta^2g^{-2}k(3+M_2)|$; for the centered
second trace it is $|1-h_\theta^2g^{-2}(6+M_2+a_2)|$.
Their squares are $\le1-c_kg^{-2}$, since $|M_2|\le1+O(g^{-2})$ and
$|a_2|<1/4+O(g^{-2})$. The remaining leading trace terms involve deviations at lower
indices or the product $M_1M_{k+1}$, which is bounded by $C_k|M_1|$
on the bounded state space. Terms in the $F_k$ update contain lower-index $F_j$'s or are bounded by
$C_k|M_1|$, including the nonlinear term in \eqref{lg:eq:center-corner}.
Choose positive weights in increasing order of the indices to absorb
these finitely many terms into the corresponding diagonal contractions.

For a finite index set $J$ containing, with each $k$, every $j$ with
$0<|j|\le|k|$, let
\[
 A_{J,p}(n)\coloneqq\sum_{k\in J}u_k\E|M_{k,n}-a_k|^{2p},\qquad
 B_{J,p}(n)\coloneqq\sum_{k\in J}v_k\E|F_{k,n}|^{2p},
\]
where the positive weights $u_k,v_k$ are chosen as just described
and are fixed independently of $W$. Expansion around the
conditional means, using the same conditional-moment inequality
as in Section~\ref{lg:sec:concrete-Q}, gives
\begin{equation}\label{lg:eq:momentrec}
 \begin{aligned}
 A_{J,p}(n+1)&\le(1-cg^{-2})A_{J,p}(n)
   +Cg^{-2}[W^{-2p}+g^{-4p}A_{J',p}(n)+g^{-4pR}],\\
 B_{J,p}(n+1)&\le(1-cg^{-2})B_{J,p}(n)
   +Cg^{-2}[W^{-p}+A_{J',p}(n)+g^{-4p}B_{J',p}(n)
                                       +g^{-4pR}].
 \end{aligned}
\end{equation}
Here $J'$ is a larger fixed set, also closed under lower nonzero indices,
containing the indices appearing in the higher-order update terms for $J$.
To justify the powers, use
\[
 g^{-4}|Y|^{2p-1}|U|\le cg^{-2}|Y|^{2p}
                       +Cg^{-2(2p+1)}|U|^{2p}.
\]
Centered innovations cancel the linear term. Their variance term
$g^{-2}W^{-2}|Y|^{2p-2}$ leaves $Cg^{-2}W^{-2p}$;
for the $F_k$'s replace $W^{-2}$ by $W^{-1}$. Higher innovation terms from
\eqref{lg:eq:innovations} are no larger after Young's inequality.
The residual $g^{-2R-2}$ leaves $g^{-4pR}$ after dividing
by the contraction rate. The finite-width errors in the conditional
means are smaller than
the variance contributions.

\noindent\textbf{Iteration and the probability bound.}
We now prove \eqref{lg:eq:momentout} and deduce
the simultaneous probability bound \eqref{lg:eq:burnevent}.
First enlarge the trace index sets to include every index needed
in the $F_k$ estimates. Iterate each contracting convolution at least $R$ times,
retaining a bounded terminal moment. The normalized geometric
kernel has total mass at most one; the coupling to the larger index set gains
$g^{-4p}$ at each iteration. Initial-state terms are bounded
by $C(1+g^{-2} n)^K e^{-cg^{-2} n}$ for a fixed finite $K$,
hence by $Ce^{-c'g^{-2} n}$. It follows that
\begin{equation}\label{lg:eq:momentout}
 A_{J,p}(n)\le C(e^{-c'g^{-2} n}+W^{-2p}+g^{-4pR}),\quad
 B_{J,p}(n)\le C(e^{-c'g^{-2} n}+W^{-p}+g^{-4pR}).
\end{equation}
Choose $R$ with $2\epsilon R>1$
and choose $p$ large enough to take a union bound over $W^{B_1}$
times with total failure probability at most $W^{-D}$. Then choose
$C_*$ in \eqref{lg:eq:initial-time} so that the exponential term is
also sufficiently small for every $n\ge n_*$.
Markov's inequality at thresholds
$W^{\eta-1}$ and $W^{\eta-1/2}$ proves \eqref{lg:eq:burnevent}
simultaneously for the finite set of indices used in the calculation.
The sets of indices and the moment order are fixed before the
union bound.

This completes the proof of Lemma~\ref{lg:lem:dynamic}.
\end{proof}

\subsubsection{Construction of the scalar correction and its drift}
\label{lg:sec:cocycle-correction}
\label{lg:sec:define-Q}
\label{lg:sec:prove-Q}

\begin{lemma}[A bounded corrector and positive drift]\label{lg:lem:centerQ}
Under the coupling assumptions of Lemma~\ref{up:lem:center-transfer}, fix $D>0$
and a polynomial horizon $W^{B_1}$ with $B_1>3$.
There are constants $\delta,C_*,C>0$ such that, for all sufficiently
large $W$ and each deterministic $|\theta|\le\delta$, $W^\epsilon\le g<W$,
there is a deterministic smooth real-valued function $Q_g$ on
$U(W)\times\{v:\|v\|_2=1\}$ with the following uniform bounds.
The notation $Q_g$ suppresses its dependence on $\theta,W$.
\begin{enumerate}
\item The size and derivative bounds \eqref{lg:eq:Q-regularity} hold.
\item Let $n_*$ be as in \eqref{lg:eq:initial-time}.
For each deterministic initial state, with probability at least $1-W^{-D}$,
simultaneously for $n_*\le i-1<W^{B_1}$,
\begin{equation}\label{lg:eq:corrected-drift}
 \left|\E_{i-1}[c_i-Q_g(C_i,v_i)+Q_g(C_{i-1},v_{i-1})]
       -\tfrac12h_\theta^2g^{-2}W^{-1}\right|
 \le\tfrac16h_\theta^2g^{-2}W^{-1}.
\end{equation}
\end{enumerate}
 
\end{lemma}

\begin{proof}[Proof of Lemma~\ref{lg:lem:centerQ}]
\noindent\textbf{The coefficient equation.}
Use the corrected coordinates $M_k,F_k$ in
\eqref{lg:eq:corrected-observables} and the centers $a_k$ in
\eqref{lg:eq:center-recursion}, constructed in Lemma~\ref{lg:lem:dynamic}
through order $g^{-2R-4}$, with $2\epsilon R>1$.

Lemma~\ref{lg:lem:ideal} and the polynomial substitutions give
the expansion of $(P_g-I)F_r$. For $\E_{i-1}c_i$, continue the
inverse and logarithm expansions used to prove \eqref{lg:eq:c2}.
The proof of Lemma~\ref{lg:lem:ideal} shows that every positive-order
expected coefficient of the squared norm contains a nonzero-index
$f_k$ factor or an explicit $W^{-1}$. The logarithm expansion and
the inverse substitution from $m,f$ to $M,F$ preserve this property,
with $F_k$ in place of $f_k$. The same cutoff argument gives a
uniform expected remainder at the chosen finite order.

Express both expansions in $M,F$. In the terms without an
explicit $W^{-1}$, set $M=a$ and extract the terms linear in $F$:
\begin{equation}\label{lg:eq:corrector-linear-expansions}
 \begin{aligned}
 \E_{i-1}c_i
 &=\tfrac12h_\theta^2g^{-2}W^{-1}
       +\sum_{k\ne0}b_kF_k+\mathcal R_0,\\
 (P_g-I)F_r&=\sum_{k\ne0}\mathsf A_{kr}F_k+\mathcal R_r.
 \end{aligned}
\end{equation}
With $[F_k]$ denoting the coefficient of exactly one factor
$F_k$ and no other $F_j$ factor, the coefficients are
\begin{equation}\label{lg:eq:Q-coefficients}
 \begin{aligned}
 b_k&\coloneqq[F_k]\left.\bigl[\E_{i-1}c_i\bigr]_{\le R+2,0}\right|_{M=a},\\
 \mathsf A_{kr}&\coloneqq[F_k]\left.\bigl[(P_g-I)F_r\bigr]_{\le R+2,0}\right|_{M=a}.
 \end{aligned}
\end{equation}
Here the brackets have the truncation convention of
Section~\ref{lg:sec:general-state}, applied after the polynomial
change of variables and before substituting $M=a$.
Positive and negative indices are separate formal variables.
The remainders collect terms containing a trace deviation
$M_j-a_j$ together with an $F$ factor, at least two nonzero-index
$F_j$ factors, an explicit $W^{-1}$, or a discarded Taylor order.
The constant positive drift is kept separate from $b$.

For $Q_g=\sum_rq_rF_r$, we obtain
\begin{equation}\label{lg:eq:corrector-coefficient-cancellation}
 \begin{aligned}
 \E_{i-1}c_i-(P_g-I)Q_g
 &=\tfrac12h_\theta^2g^{-2}W^{-1}
       +\sum_k[b_k-(\mathsf A q)_k]F_k +\mathcal R_0-\sum_rq_r\mathcal R_r.
 \end{aligned}
\end{equation}

\noindent\textbf{(1) Bounding the remainders.}
The goal is \eqref{lg:eq:corrector-remainder-bounds}.
Apply \eqref{lg:eq:burnevent} with $\eta/2$ in place of $\eta$,
where $0<\eta<\min(\epsilon,1/4)$. For the fixed finite set of
indices used below, it gives
\begin{equation}\label{lg:eq:corrector-remainder-bounds}
 |\mathcal R_0|=o(W^{-1}g^{-2}),\qquad
 \max_r|\mathcal R_r|\le C_RW^{-1+\eta}g^{-2}.
\end{equation}
Indeed, the terms in $\mathcal R_r$ start at $g^{-2}$ and contain
$(M_j-a_j)F_k$, two $F$ factors, or an explicit $W^{-1}$;
\eqref{lg:eq:burnevent} bounds each of these by $O(W^{-1+\eta})$.
For $\mathcal R_0$, \eqref{lg:eq:c2} and $a_1=0$ show that its
order-$g^{-2}$ terms are bounded by constant multiples of
$g^{-2}|(M_1-a_1)F_1|$ and $g^{-2}W^{-1}|F_1|^2$, both
$o(W^{-1}g^{-2})$ by \eqref{lg:eq:burnevent}.
The remaining terms have the extra factor $g^{-2}$ and total
size $O(g^{-4}W^{-1+\eta})=o(W^{-1}g^{-2})$.
The Taylor and inverse-substitution errors $O(g^{-2R-2})$ are
also negligible since $2\epsilon R>1$.
Thus, keeping the coefficient factor explicit,
\begin{equation}\label{lg:eq:centerdefect}
 \left|\mathcal R_0-\sum_rq_r\mathcal R_r\right|
 \le o(W^{-1}g^{-2})+C_RW^{-1+\eta}g^{-2}\sum_r|q_r|.
\end{equation}
\noindent\textbf{(2) Solving the coefficient equation.}
The goal is the coefficient and residual bound
\eqref{lg:eq:raw-Q-coefficient-bounds}.
By \eqref{lg:eq:required-parity-updates}, every retained monomial
in the update of $F_r$ has degree congruent to $r$ modulo two.
At $M=a$, a monomial containing an odd-index $M_j$ vanishes,
since $a_{2j+1}=0$. A surviving term linear in $F_k$ therefore
has $k\equiv r\pmod2$. Thus
\[
 \mathsf A_{kr}=0\quad(k\not\equiv r\pmod2),\qquad
 \mathsf A=
 \begin{pmatrix}\mathsf A_{\rm odd}&0\\0&\mathsf A_{\rm even}\end{pmatrix}.
\]

Let $\mathsf A_0$ consist of the terms of orders $1$ and $g^{-2}$ in $\mathsf A$:
\[
 (\mathsf A_0)_{kr}\coloneqq
 [F_k]\left.\bigl[(P_g-I)F_r\bigr]_{\le1,0}\right|_{M=a}.
\]
In \eqref{lg:eq:center-corner}, set $m=a$ and $f=F$.
Since $a_1=0$, the only higher-index term $-(k+1)a_1F_{k+1}$ vanishes;
the vanishing of all odd $a_j$ leaves only $F_k$ and lower indices
of the same parity, so each parity block of $\mathsf A_0$ is triangular.
The diagonal contribution is $-3kF_k-ka_2F_k=-k(3+a_2)F_k$,
where $-ka_2F_k$ comes from $n=k$ in the first sum.
Including the phase and subtracting $I$ gives the diagonal entry, for $k>0$,
\[
 (-1)^ke^{-2ik\theta}
       [1-kh_\theta^2g^{-2}(3+a_2)]-1,
\]
and the negative-index diagonals are the conjugates.
Shrink $\delta$ so that $|k|\delta\le\pi/3$ for every retained odd
index. Then $|1+e^{-2ik\theta}|\ge1$, so the odd diagonals remain
bounded away from zero after the $O(g^{-2})$ correction. The even
diagonals have modulus at least $c_kg^{-2}$ by
\eqref{lg:eq:center-denominators}. Off-diagonal entries are $O(g^{-2})$, so triangular
inversion gives $O(1)$ on the odd block and $O(g^2)$ on the even
block. In particular,
\[
 \|\mathsf A-\mathsf A_0\|\le C_Rg^{-4},\qquad
 \|\mathsf A_0^{-1}(\mathsf A-\mathsf A_0)\|\le C_Rg^{-2}.
\]
The second bound uses the exact energy dependence of the even
diagonals; it does not replace $\theta$ by zero.

Since $m_1=M_1+O(g^{-2})$, $f_1=F_1+O(g^{-2})$, and $a_1=0$,
\eqref{lg:eq:c2} gives
\[
 \sum_kb_kF_k=h_\theta^2g^{-2}\Re F_1
                        +\sum_kO(g^{-4})F_k.
\]
Thus the odd coefficients of $b$ start at $g^{-2}$, whereas the
even coefficients start at $g^{-4}$. The extra factor $g^{-2}$
in the even right-hand side compensates for its $O(g^2)$ inverse:
\[
 \| (\mathsf A_0^{-1}b)_{\rm odd}\|_1\le C_Rg^{-2},\qquad
 \| (\mathsf A_0^{-1}b)_{\rm even}\|_1
       \le C_Rg^2g^{-4}=C_Rg^{-2}.
\]
Define the successive coefficient corrections by
\begin{equation}\label{lg:eq:Q-recursion}
 q^{(0)}\coloneqq\mathsf A_0^{-1}b,\qquad
 q^{(s+1)}\coloneqq-\mathsf A_0^{-1}
       (\mathsf A-\mathsf A_0)q^{(s)},\quad 0\le s<R,
\end{equation}
and set
\begin{equation}\label{lg:eq:Q-definition}
 q\coloneqq\sum_{s=0}^{R}q^{(s)},\qquad
 Q_g(C_{i-1},v_{i-1})\coloneqq
       \sum_{k\ne0}q_kF_k(C_{i-1},v_{i-1}).
\end{equation}
Each correction gains a factor $g^{-2}$:
$\|q^{(s)}\|_1\le C_Rg^{-2s-2}$.
Telescoping the coefficient equations gives
$\mathsf A q-b=(\mathsf A-\mathsf A_0)q^{(R)}$, and hence
\begin{equation}\label{lg:eq:raw-Q-coefficient-bounds}
 \sum_k|q_k|\le C_Rg^{-2},\qquad
 \|\mathsf A q-b\|_1\le C_Rg^{-2R-6}.
\end{equation}
The recursion acts on the retained polynomial coefficients.

The finite-index convention in Section~\ref{lg:sec:general-state}
includes the support of $b$, the indices in all $R$ corrections
and the final residual, and the trace indices in their coefficients
and centers. Apply Lemma~\ref{lg:lem:dynamic} to this entire list;
no term is dropped because its index lies outside an initial list.

\noindent\textbf{(3) The norm and derivative bounds for $Q_g$.}
The goal is \eqref{lg:eq:Q-regularity}.
Conjugation of \eqref{lg:eq:Q-recursion} gives
$q_{-k}=\overline{q_k}$, so $Q_g$ is real.
For real tangent directions $\dot C,\dot v$,
\[
 \begin{aligned}
 |(\dot C\cdot\nabla_C)m_j|&\le |j|\|\dot C\|_{\rm op},\\
 |(\dot C\cdot\nabla_C+\dot v\cdot\nabla_v)f_j|
 &\le |j|\|\dot C\|_{\rm op}+2\|\dot v\|_2.
 \end{aligned}
\]
The finite polynomial definitions of $F_k$, their bounded
coefficients, and $|m_j|,|f_j|\le1$ therefore give
\[
 |F_k|\le C_R,\qquad
 |(\dot C\cdot\nabla_C+\dot v\cdot\nabla_v)F_k|
 \le C_R(\|\dot C\|_{\rm op}+\|\dot v\|_2).
\]
Combining these bounds with
\eqref{lg:eq:raw-Q-coefficient-bounds} proves both bounds in
\eqref{lg:eq:Q-regularity}.

\noindent\textbf{Completion of the proof.}
By \eqref{lg:eq:raw-Q-coefficient-bounds} and the bound $|F_k|\le C_R$,
the linear residual in \eqref{lg:eq:corrector-coefficient-cancellation}
has modulus at most $C_Rg^{-2R-6}$.
Substituting $\sum_r|q_r|\le C_Rg^{-2}$ into
\eqref{lg:eq:centerdefect} and using
\eqref{lg:eq:corrector-remainder-bounds} therefore gives
\begin{equation}\label{lg:eq:smalldefect}
 \begin{aligned}
 &\left|\E_{i-1}c_i-(P_g-I)Q_g
       -\tfrac12h_\theta^2g^{-2}W^{-1}\right|\\
 &\quad\le o(W^{-1}g^{-2})+C_RW^{-1+\eta}g^{-4}
       +C_Rg^{-2R-6}=o(W^{-1}g^{-2}),
 \end{aligned}
\end{equation}
since $g\ge W^\epsilon$, $\eta<\epsilon$, and $2\epsilon R>1$.
For sufficiently large $W$, the error
in the corrected drift is at most $h_\theta^2/(6Wg^2)$.
On the event \eqref{lg:eq:burnevent}, this holds simultaneously for
$n_*\le i-1<W^{B_1}$, with probability at least $1-W^{-D}$.
This proves \eqref{lg:eq:corrected-drift}.
\end{proof}

\subsubsection{From the corrected drift to transfer growth}

\begin{lemma}[Scalar growth]\label{lg:lem:scalar-growth}
Under the coupling assumptions of Lemma~\ref{up:lem:center-transfer}, let $\ell_n$
be defined by \eqref{lg:eq:additive-growth}. For every $D>0$, there are
$\delta>0$ and $C_D>0$ such that, for every fixed $C_1\ge C_D$ and all
sufficiently large $W$, the deterministic length
$n=\lceil C_1Wg^2\log W\rceil$ satisfies
\begin{equation}\label{lg:eq:scalar-growth}
 \sup_{|\theta|\le\delta}
 \sup_{\substack{C_0\in U(W)\\\|v_0\|_2=1}}
 \Prob\{\ell_{n}(C_0,v_0)<D\log W\}\le W^{-D},
\end{equation}
uniformly for $W^\epsilon\le g<W$. The constants $\delta,C_D$
depend only on $\epsilon,D$; the width threshold may also depend
on $C_1$.  
\end{lemma}

\begin{proof}[Proof of Lemma~\ref{lg:lem:scalar-growth}]
For $n\ge n_*$, take absolute values of the error terms in
\eqref{lg:eq:corrected-growth} to obtain
\begin{equation}\label{lg:eq:growth-lower-bound}
 \begin{aligned}
 \ell_n\ge{}&\frac{h_\theta^2g^{-2}(n-n_*)}{2W}
       -|\ell_{n_*}|-2\|Q_g\|_\infty
       -\left|\sum_{i=n_*+1}^n(1-\E_{i-1})[c_i-Q_g(C_i,v_i)]\right|\\
 &-\sum_{i=n_*+1}^n\left|\E_{i-1}
       [c_i-Q_g(C_i,v_i)+Q_g(C_{i-1},v_{i-1})]
       -\tfrac12h_\theta^2g^{-2}W^{-1}\right|.
 \end{aligned}
\end{equation}
Fix $D>0$ and $B_1>3$. Apply
Lemmas~\ref{lg:lem:centerQ} and~\ref{lg:lem:fluctuations}
with probability exponent $D+3$.
For their common $n_*$ from \eqref{lg:eq:initial-time} and each deterministic
$n_*\le n\le W^{B_1}$, the four bounds in those lemmas imply
\begin{equation}\label{lg:eq:integrated}
 \ell_n\ge h_\theta^2g^{-2}(n-n_*)/(3W)-C\log W
                        -C\sqrt{ng^{-2}W^{-1}\log W},
\end{equation}
with failure probability at most $3W^{-D-3}$.
Here the martingale sum is bounded directly by
\eqref{lg:eq:centerMGF}, and the endpoint cost $Cg^{-2}$ is absorbed
into $C\log W$.

Since $h_\theta^2\ge\tfrac14$, at $n=\lceil C_1Wg^2\log W\rceil$
equation \eqref{lg:eq:integrated} gives
\[
 \ell_n\ge(C_1/12-C\sqrt{C_1}-C)\log W,
\]
absorbing the initial-length contribution $O(C_*\log W/W)$.
Thus $\ell_{n}\ge D\log W$ for every fixed $C_1\ge C_D$,
with failure probability at most $W^{-D}$ for large $W$.
The choice of $B_1$ puts $n$ inside the stated horizon.
Uniformity in the deterministic initial state and energy gives
\eqref{lg:eq:scalar-growth}.
\end{proof}

\begin{proof}[Proof of Lemma~\ref{up:lem:center-transfer}]
Apply Lemma~\ref{lg:lem:scalar-growth} with $D>D'+2$ and choose
$0<\kappa_1<\min\{\sin\delta,(2-\kappa_0)/2\}$.
Applying \eqref{lg:eq:scalar-reduction} and taking a union bound over
the $W$ basis vectors gives \eqref{up:eq:input} with $J=[-2\kappa_1,2\kappa_1]$, $\mathfrak b=Wg^2$, and the
length $n_0=n$ from Lemma~\ref{lg:lem:scalar-growth}, uniformly
in the deterministic energy and initial frame.
\end{proof}

\subsection*{Completion of the upper bound}

\begin{proof}[Proof of Theorem~\ref{thm:large-upper}]
Set $J=[-2+\kappa,2-\kappa]$ and $\mathfrak b=Wg^2$.

\noindent\textbf{The range $g\ge W$.}
For sufficiently large $W$, we have $g\ge W\ge C_\kappa\sqrt W$,
so Theorem~\ref{uvg:thm:transfer} gives
\eqref{up:eq:negative-transfer-input}, uniformly in $g,n$.
Proposition~\ref{up:prop:conversion}\textup{(ii)} converts this to
\eqref{uva:eq:conversion-input} with $q=s/4$, and
Lemma~\ref{uva:thm:upper-audit}\textup{(ii)} gives
\eqref{eq:averaged-tail} with $R\le C_D\mathfrak b\log W$.
The constants and width threshold are independent of $g,L$.

\noindent\textbf{The range $W^\epsilon\le g<W$.}
This range is empty if $\epsilon\ge1$. For $0<\epsilon<1$, we have
$\mathfrak b<W^3$. Fix $0<q<1/4$ and the exponent $D'$ required by
Proposition~\ref{up:prop:conversion}\textup{(i)} with $t=2q$
and polynomial-range parameter $C_0=3$.
Lemma~\ref{up:lem:center-transfer} gives \eqref{up:eq:input}
on $[-2\kappa_1,2\kappa_1]$, with scale $\mathfrak b=Wg^2$.
The near-zero window is fixed after $D'$; neither depends on the
exceptional-count exponent $D$.
On $[\kappa_1,2-\kappa]$ and $[-2+\kappa,-\kappa_1]$,
Lemma~\ref{up:lem:density-contraction}\textup{(ii)}, followed by
Lemmas~\ref{up:lem:angular}, \ref{up:lem:operator-iteration}
and~\ref{up:lem:operator-reduction}, gives the same transfer input
with batch length $1000$.

On each interval, Proposition~\ref{up:prop:conversion}\textup{(i)}
gives \eqref{uv:eq:upper-volume-input}, and
Proposition~\ref{up:prop:fr-concentration} gives
\eqref{up:eq:fmr-pointwise} with exponent $q$ and scale
$\mathfrak b=Wg^2$. Choose the length parameter $C_1$ in the finite-volume
criterion to cover the resulting length constants, keeping $C_0=3$.
The lengths may differ between the three intervals; taking the larger
prefactor exponent and smaller decay rate gives
\eqref{up:eq:fmr-pointwise} on their union $J$.
All energy and frame suprema in the transfer inputs remain outside
the probability; there is no union over couplings or volumes.
Lemma~\ref{uva:thm:upper-audit}\textup{(i)} gives
\eqref{eq:averaged-tail} with $R\le C_D\mathfrak b\log W$,
with constants and width threshold uniform in this range and in $L$.

Taking the larger constants for the two ranges gives the same bound
uniformly for all $g\ge W^\epsilon$. For sufficiently large $W$,
$C_D\log W\le(\log W)^2$, proving \eqref{eq:large-radius} with a
width threshold independent of $g,L$. When the radius fills the cycle,
the tail is empty.
\end{proof}

\appendix
\section{Dyson solutions and deterministic kernel estimates}
\label{app:dyson-foundations}\label{app:deterministic}
This appendix proves the deterministic results of
Section~\ref{sec:deterministic-inputs}, together with the scalar
comparison facts used in Sections~\ref{sec:flow-overview} and
\ref{sec:input-proofs}. Each auxiliary result is proved before its use.

\subsection{Free Green functions and scalar identities}
\label{app:dyson-common}
To compare onsite variances at fixed hopping in
Lemma~\ref{da:lem:comparison-parameters}(a), allow a general variance
$\gamma>0$ and put $A_{g,\gamma}:=\sqrt\gamma\,gA_L$. The choice
$\gamma=\gamma_g$ gives $A_{g,\gamma_g}=A_g$ and recovers the notation of
\eqref{fo:eq:dyson-family}--\eqref{fo:eq:characteristic}.
For $w\in\C_+$, define
\begin{equation}\label{bl:eq:physical}
 \begin{split}
 f_g(w;\gamma)&\coloneqq \int_{-\pi}^{\pi}
       \frac{\mathrm dp/(2\pi)}{2\sqrt\gamma\,g\cos p-w}
       =-\frac1{\sqrt{w^2-4\gamma g^2}},\\
 f_{g,L}(w;\gamma)&\coloneqq L^{-1}\sum_{p\in(2\pi/L)\Z_L}
                   \frac1{2\sqrt\gamma\,g\cos p-w}.
 \end{split}
\end{equation}
The square root is asymptotic to $w$ at infinity and has positive
imaginary part on $\C_+$. At $g=0$, $f_0(w;\gamma)=-1/w$.
The scalar Dyson solutions and block matrix are
\begin{equation}\label{da:eq:dyson}
 \begin{gathered}
 m_g(z;\gamma)=f_g(z+\gamma m_g(z;\gamma);\gamma),\\
 m_{g,L}(z;\gamma)=f_{g,L}(z+\gamma m_{g,L}(z;\gamma);\gamma),\\
 M_B(z;\gamma)\coloneqq (A_{g,\gamma}-z-\gamma m_{g,L}(z;\gamma))^{-1},\qquad
 M(z;\gamma)\coloneqq M_B(z;\gamma)\otimes I_W.
 \end{gathered}
\end{equation}
We retain the Herglotz normalization and regular boundary values, write
$m\coloneqq m_{g,L}(z;\gamma)$, and suppress fixed variance arguments.
Let $\nu_g(\cdot;\gamma),F_g(\cdot;\gamma)$ be the measure and distribution
function of $m_g(\cdot;\gamma)$. At $\gamma=\gamma_g$,
$m_{g,L}(z;\gamma_g)=m_{g,L}(z)$ and $M_B(E;\gamma_g)=M_B(g,E)$,
as in
\eqref{fo:eq:dyson-family}--\eqref{fo:eq:characteristic}.

At a real regular energy $E$, take $m=m_{g,L}(E;\gamma)$ and
$M=M(E;\gamma)$. Translation invariance and \eqref{da:eq:dyson}
give $M_{aa}=mI_W$. With $\Lambda=A_{g,\gamma}\otimes I_W$,
$z_t=E+(1-t)\gamma m$ and
$\mathcal S[R]=\gamma W^{-1}\diag(\Tr R_{aa})_a\otimes I_W$,
direct substitution gives
\begin{equation}\label{da:eq:characteristic-general}
 \mathcal S[M]=\gamma mI_{WL},\qquad
 \Lambda-z_t-t\mathcal S[M]=\Lambda-E-\gamma m=M^{-1}.
\end{equation}

The spatial decay factor $\rho$ of the free Green function will be used
for finite-volume comparisons and the calculation of propagation scales.
For $w\in\C_+$ put
\begin{equation}\label{da:eq:free-decay-factor}
 \rho(w)\coloneqq \frac{2\sqrt\gamma\,g}{w+\sqrt{w^2-4\gamma g^2}},\qquad
 |\rho(w)|<1.
\end{equation}

\begin{lemma}[Free Green functions]\label{da:lem:free-green}
Write $R_{g,\infty}(x;w)$ and $R_{g,L}(x;w)$ for the $(0,x)$ entries
of the free resolvents on $\Z$ and $\Z_L$, respectively, both with
hopping $\sqrt\gamma\,g$. Thus $R_{g,L}(x;w)=((A_{g,\gamma}-w)^{-1})_{0x}$.
For $0\le x<L$, these entries satisfy
\begin{equation}\label{bl:eq:free}
 R_{g,\infty}(x;w)=f_g(w)\rho(w)^{|x|},\qquad
 R_{g,L}(x;w)=f_g(w)
       \frac{\rho(w)^x+\rho(w)^{L-x}}{1-\rho(w)^L}.
\end{equation}
At $x=0$ the numerator is $1+\rho^L$.
At $g=0$ the limiting row is $-w^{-1}\mathbf1_{x=0}$.
In particular, writing $M_B(x)=(M_B)_{0x}$, we have
$M_B(x)=R_{g,L}(x;z+\gamma m)$.\end{lemma}

\noindent\textbf{Proof of Lemma~\ref{da:lem:free-green} and \eqref{bl:eq:physical}.}
For $g>0$, apply the free Jacobi Green function formula
\cite[Section~1.3, (1.121)]{Teschl} at spectral parameter
$w/(2\sqrt\gamma\,g)$ and multiply by $(2\sqrt\gamma\,g)^{-1}$.
This gives the line formula, including $f_g(w)$ at $x=0$;
Fourier diagonalization gives the finite sum in \eqref{bl:eq:physical}.
Periodization and geometric summation yield
\[
 R_{g,L}(x;w)=\sum_{n\in\Z}R_{g,\infty}(x+nL;w)
 =f_g(w)\frac{\rho(w)^x+\rho(w)^{L-x}}{1-\rho(w)^L}.
\]
The $g=0$ formula is immediate.

\noindent\textbf{Proof of \eqref{da:eq:ward}.}
For the following identities, take $m=m_g(z;\gamma)$ with normalized
Fourier integration on the line, or $m=m_{g,L}(z;\gamma)$ with normalized
averaging over the finite grid, and denote the corresponding average by
$\langle\cdot\rangle$. Write
$R(p)\coloneqq (2\sqrt\gamma\,g\cos p-z-\gamma m)^{-1}$ and $z=E+i\eta$.
Taking imaginary parts in the Dyson equation gives
\begin{equation}\label{da:eq:ward}
 \langle|R|^2\rangle=\frac{\Im m}{\eta+\gamma\Im m},\qquad
 |m|\le\gamma^{-1/2},\qquad
 \Re(1-\gamma\langle R^2\rangle)
 \ge 2\gamma(\Im m)^2.
\end{equation}
The middle bound follows from $|m|^2\le\langle|R|^2\rangle$.
For the last inequality, use
$\Re R^2=|R|^2-2(\Im R)^2$ and Jensen's inequality.
For real $E$ with $\Im m_{g,L}(E;\gamma)>0$, Parseval gives the exact identity
$\gamma\sum_x|M_B(x)|^2=1$. For $\eta>0$, the difference of two
Dyson equations has coefficient of modulus less than one by
Cauchy--Schwarz and the first identity; this proves uniqueness.
The same identities control regular boundary values by continuity.

\subsection{Bulk regularity and the deterministic distribution}
\label{app:bulk-regularity}
We first prove the scalar regularity used in the distribution and
finite-cycle statements.

\begin{lemma}[Bulk regularity]\label{bl:lem:regular}
Fix $0<C<\infty$, $\kappa\in(0,1/2)$, and
$\max\{1,C\}^{-2}\le\gamma\le1$. Use the Dyson functions and their
probability measures in \eqref{da:eq:dyson}, suppressing the fixed
variance argument $\gamma$ throughout. There are
\(b_\kappa,a_{\rm strip},c,C'>0\), depending only on
\(C,\kappa\), such that, whenever
\begin{equation}\label{bl:eq:regular-domain}
 g\le C,\qquad
 \dist(x,\{E:F_g(-|E|;\gamma)\ge\kappa\})\le b_\kappa,\qquad
 z=x+i\eta,\qquad 0\le\eta\le1.
\end{equation}
we have \(|x|\le2C+2\) and
\begin{equation}\label{bl:eq:regular-bounds}
 |m_g(z)|\le\gamma^{-1/2},\qquad \Im m_g(z)\ge a_{\rm strip},\qquad
 \Re\bigl(1-\gamma f_g'(z+\gamma m_g(z))\bigr)
       \ge2\gamma(\Im m_g(z))^2.
\end{equation}
Every fixed mixed derivative of \(m_g\) in \(g\) and \(z\)
is uniformly bounded on this set. The free row defined in \eqref{bl:eq:free} satisfies
\[
 \sum_{n\in\Z}e^{c|n|}
       |R_{g,\infty}(n;z+\gamma m_g(z))|\le C'.
\]
The same conclusions hold with the bulk interval replaced by
interval \(\sqrt\gamma[-2+\kappa_0,2-\kappa_0]\),
\(\kappa_0\in(0,2)\), with constants depending on \(C,\kappa_0\).
\end{lemma}

\begin{proof}[\textbf{Proof of Lemma~\ref{bl:lem:regular}}]
First take $\gamma=1$. For the law $\tau_g$ of $2g\cos p$,
the boundary description in
\cite[Lemmas 2--4, Proposition 2 and Corollary 3, pp. 709--711]{Biane}
uses the negative of our Stieltjes transform. Since
$\int(u-y)^{-2}\tau_g(dy)$ is infinite on $[-2g,2g]$ and strictly
decreases to zero for $u>2g$, its set $U_1$ is $(-b_g,b_g)$;
at $g=0$, $b_0=1$. The boundary homeomorphism and density formula
therefore give one positive-density band $(-e_g,e_g)$, with
\[
 \int(b_g-y)^{-2}\tau_g(dy)=1,\quad b_g\le2g+1,\quad
 2\le b_g+b_g^{-1}\le e_g
 =b_g+\int(b_g-y)^{-1}\tau_g(dy)\le2g+2.
\]
Here Jensen gives the lower bound and Cauchy--Schwarz the upper bound.
The density is continuous and at most $1/\pi$ by \eqref{da:eq:ward},
so each interior quantile has distance at least $\pi\kappa$ from its
band edge. Weak continuity in $g$, interior quantile continuity, and
implicit stability in \eqref{da:eq:ward} give, by compactness, fixed
$b\le\pi\kappa/4$ and $c'>0$ with $\Im m_g\ge c'$ on the
$2b$-neighborhood of the quantile interval, uniformly for $g\le C$,
including $g=0$. For $x$ in its $b$-neighborhood and any fixed $H$,
Poisson integration on $[x-b,x+b]$ gives
\[
 \Im m_g(x+i\eta)\ge\frac{2c'}\pi\arctan(b/\eta)
 \ge\frac{2c'}\pi\arctan(b/H)>0\qquad(0<\eta\le H).
\]
The stability denominator in \eqref{da:eq:ward} is thus uniformly
positive. Implicit differentiation gives all fixed mixed derivatives;
shifting the free Fourier contour by $\pm2ic$, with
$2C\sinh(2c)$ smaller than half the imaginary gap, gives
$|R_{g,\infty}(n)|\le C'e^{-2c|n|}$ and the weighted row.
The fixed interval $[-2+\kappa_0,2-\kappa_0]$ is handled in the
same way, since $e_g\ge2$.

For general $\gamma$, use
\[
 m_g(z;\gamma)=\gamma^{-1/2}m_g(z/\sqrt\gamma;1),\qquad
 R_{g,\infty}(n;w;\gamma)
 =\gamma^{-1/2}R_{g,\infty}(n;w/\sqrt\gamma;1).
\]
Take $H=\max\{1,C\}$ and the physical buffer at most $b/H$.
These identities cover $\eta\le1$, rescale the quantiles and fixed
interval by $\sqrt\gamma$, and transfer the derivative and row bounds
uniformly over the stated variance interval; derivatives hold $\gamma$ fixed.
\end{proof}

\begin{proof}[\textbf{Proof of Lemma~\ref{da:lem:distribution}}]
In a fixed free realization, $\nu_g$ is the law of
$(g/\max\{1,g\})A+\max\{1,g\}^{-1}S$, where $A$ is arcsine and
$S$ is semicircular of variance one. Symmetry and the norm bound give
\eqref{eq:bounded-physical-bulk}; continuity of the coefficients on
$[0,\infty]$ gives weak continuity and the two endpoint laws.
Lemma~\ref{bl:lem:regular} supplies the single band and continuous
density for finite $g$.
Every law in this weakly continuous compact family is atomless,
including the arcsine endpoint. Compactness and the Portmanteau theorem
therefore give \eqref{rd:eq:noatoms}.

For \eqref{eq:bulk-bridge}, write $F_0,F_\infty$ for the semicircle
and arcsine distribution functions. In the same free realization,
$\|gA\|\le2g$ and $\|S/g\|\le2/g$, so spectral comparison gives
\[
 F_g(x)\le F_0(x+2g)\quad(g\le1),\qquad
 F_g(x)\le F_\infty(x+2/g)\quad(g\ge1).
\]
Take $\kappa_0=c_\kappa=\kappa^2/4$ and $C_\kappa=8/\kappa^2$.
For $g\le c_\kappa$ the first bound at $x=-2+\kappa_0$ is at most
$F_0(-2+3\kappa^2/4)<\kappa$; for $g\ge C_\kappa$ the second is at most
$F_\infty(-2+\kappa^2/2)<\kappa$, by the explicit endpoint densities.
Monotonicity of $F_g$ proves \eqref{eq:bulk-bridge}.
\end{proof}

\subsection{Finite-cycle regularity at bounded coupling}
\label{sec:bl-kernels}

\begin{lemma}[Finite-cycle Dyson regularity]\label{bl:lem:finite}
Fix $C<\infty$, $0<\kappa<1/2$, and
$\max\{1,C\}^{-2}\le\gamma\le1$. There are $b,c,C'>0$ such that, for
$g\le C$ and all sufficiently large $L$, uniformly on
\begin{equation}\label{da:eq:finite-domain}
 z=x+i\eta,\qquad 0\le\eta\le1,\qquad
 \operatorname{dist}\bigl(x,\{E:F_g(-|E|;\gamma)\ge\kappa\}\bigr)\le b,
\end{equation}
the line and finite-cycle solutions satisfy
\begin{equation}\label{da:eq:finite-regularity}
 \begin{gathered}
 |x|+|m_g(z)|+|m_{g,L}(z)|\le C',\qquad
 \Im m_g(z),\Im m_{g,L}(z)\ge c,\\
 |\partial_g^j\partial_z^r(m_{g,L}-m_g)(z)|
       \le C_{j,r}(1+L)^{C_{j,r}}e^{-cL},\qquad j,r\ge0\text{ fixed},\\
 \sum_{x\in\Z_L}e^{c\|x\|_L}|M_B(z)_{0x}|\le C'.
 \end{gathered}
\end{equation}
The fixed derivatives of each scalar solution are bounded. The same
conclusions hold on a fixed neighborhood of
$\sqrt\gamma[-2+\kappa_0,2-\kappa_0]$, $0<\kappa_0<2$.
Alternatively, for fixed $B_0,c'>0$, they hold locally at each
$(g,E)$ with $g\le C$, $|E|\le B_0$ and
$\Im m_{g,L}(E;\gamma)\ge c'$, with constants
depending only on $C,B_0,c'$ and the fixed variance interval.
All derivatives hold $\gamma$ fixed.
\end{lemma}

\begin{proof}[\textbf{Proof of Lemma~\ref{bl:lem:finite}}]
Use the free-row formula \eqref{bl:eq:free}. On $\Im w\ge c>0$
and $g\le C$, the identity
$w=\sqrt\gamma\,g(\rho+\rho^{-1})$ gives
$-\log|\rho|\ge\operatorname{arsinh}(c/(2C))$ for $g>0$.
Also $\rho=O(g)$ near zero.

Uniform exponential decay yields, at every fixed derivative order,
\begin{equation}\label{bl:eq:finiteerror}
 R_{g,L}(0;w)-f_g(w)=O((1+L)^{C'} e^{-cL}).
\end{equation}
The quantitative implicit function theorem, using
\eqref{bl:eq:regular-bounds}, gives the finite solution
$m_{g,L}(z)=R_{g,L}(0;z+\gamma m_{g,L}(z))$ near the line solution,
with the same error and all fixed derivatives. At real regular points
it has positive imaginary part, and its local upper continuation is
the unique finite Herglotz solution by the same Ward uniqueness
argument. No contraction $|f_g'|<1-c$ is assumed. Uniformly on
\eqref{bl:eq:regular-domain} and $0\le\eta\le1$, for $L\ge L_0$,
\begin{equation}\label{bl:eq:finitebulk}
 \Im m_{g,L}(x+i\eta)\ge a_{\rm strip}/2.
\end{equation}
The finite implicit comparison is also valid uniformly for
$g\le C$, $|E|\le B_0$ and $\Im m_{g,L}(E;\gamma)\ge c'$,
starting with its finite branch and finding
the nearby line branch. Free arguments there have imaginary part
at least $\gamma c'$ and bounded modulus; their stable denominator is
at least $2\gamma (c')^2$.

The derivative bounds and the weighted row in
\eqref{da:eq:finite-regularity} follow by differentiating the same
implicit equation and summing \eqref{bl:eq:free} with a smaller
exponential weight. Lemma~\ref{bl:lem:regular} supplies the stated
energy buffer and the alternative fixed interval.
\end{proof}

\begin{proof}[\textbf{Proof of Lemma~\ref{fo:lem:terminal-deterministic}(a)}]
Apply Lemma~\ref{bl:lem:regular} at $\gamma=\gamma_g$ and use
\eqref{bl:eq:finitebulk}. The energy bound in
\eqref{bl:eq:regular-domain} and the positive lower bound in
\eqref{bl:eq:finitebulk} give \eqref{fo:eq:terminal-bounded-regularity}.
The modulus bound in \eqref{da:eq:ward}, passed to real boundary values, gives
$|m_{g,L}(e)|\le\gamma_g^{-1/2}=\max\{1,g\}$.
\end{proof}

\paragraph{An internal small-coupling expansion.}
We prove \eqref{bl:eq:smallg}, which supplies the uniform semicircle
bulk margin used in the bounded-coupling transport calculation.
For $g\le C$, $\max\{1,C\}^{-2}\le\gamma\le1$ and sufficiently
large $L$ as in Lemma~\ref{bl:lem:finite}, assume $|E|\le B_0$ and
$\Im m_{g,L}(E;\gamma)\ge c'$, and write
$m\coloneqq m_{g,L}(E;\gamma)$. Ward's identity gives \(|m_{g,L}|\le\gamma^{-1/2}\).
The argument \(E+\gamma m_{g,L}\) of the free resolvent has imaginary part
at least \(\gamma c'\), so its exponential rows and implicit derivatives are
uniformly bounded. For sufficiently small $g$, the Neumann expansion
in $A_{g,\gamma}/(E+\gamma m)$ has no linear term since $\Tr A_{g,\gamma}=0$. Thus
\[
 m=-\frac1{E+\gamma m}
      +O\!\left(\frac{g^2}{|E+\gamma m|^3}\right),\qquad
 \gamma m^2+Em+1=O_{c',C}(g^2),
\]
where $|E+\gamma m|\ge\gamma c'$ controls the remainder uniformly.
Taking imaginary parts and dividing by $\Im m\ge c'$ gives the first
estimate below; substituting it into the real part gives the second:
\begin{equation}\label{bl:eq:smallg}
  \sqrt\gamma\,\Re m_{g,L}=-2^{-1}\gamma^{-1/2} E+O_{c',C}(g^2),\qquad
 \gamma(\Im m_{g,L})^2=1-4^{-1}\gamma^{-1} E^2+O_{c',C}(g^2).
\end{equation}
Thus, for such $(g,E)$ with \(g\) sufficiently small,
\(|E|\le2\sqrt\gamma\sqrt{1-\gamma (c')^2/2}<2\sqrt\gamma\). In particular, expansions around
zero hopping have a uniform semicircle bulk margin. For \(g\) bounded
away from zero, the positive imaginary part and the explicit free row
provide the corresponding uniform bounds. These facts give uniform bounds under the stated energy and
density conditions.

\subsection{Uniform finite-cycle comparison at large coupling}
\label{app:uniform-finite-comparison}

\begin{lemma}[Uniform periodic comparison]\label{uvp:prop:finite}
Fix $0<\kappa<2$ and $|E|\le2-\kappa$. There exist $C_\kappa\ge1$
and $R_\kappa,c,C>0$ with the following property.
Define
\[
 f(w)=-\frac1{\sqrt{w^2-4}},\qquad
 \rho(w)=\frac2{w+\sqrt{w^2-4}},\qquad w\in\C_+,
\]
where the square root is asymptotic to $w$ at infinity.
For every $g\ge C_\kappa$ and $L/g^2\ge R_\kappa$, the physical
Dyson solutions on the line and the cycle of $L$ blocks at variance $g^{-2}$ satisfy
\[
 \Im m_g(E),\Im m_{g,L}(E)\asymp1,\qquad
 |m_{g,L}(E)-m_g(E)|\le C e^{-cL/g^2}.
\]
Regard $P,Q$ from \eqref{fo:eq:PQ} as convolution kernels on $\Z_L$,
and let $P_0$ be the periodization of
$g^{-2}|f(E+g^{-2}m_g(E))\rho(E+g^{-2}m_g(E))^{|x|}|^2$.
Then $P,P_0$ are even nonnegative probability kernels and
\begin{equation}\label{uvp:eq:finiteerror}
 \norm{P-P_0}_1+g^2\norm{P-P_0}_\infty
 \le C(1+L/g^2)^C e^{-cL/g^2}=:\varepsilon_L.
\end{equation}
The first norm in \eqref{uvp:eq:finiteerror} also admits the weight
$e^{c'd(x,0)/g^2}$ for sufficiently small $c'>0$, where $d$ is cyclic
distance. The actual free row obeys
\begin{equation}\label{da:eq:uniform-free-row}
 |(M_B)_{0x}|\le Ce^{-cd(x,0)/g^2},\qquad
 \sum_x e^{c'd(x,0)/g^2}|(M_B)_{0x}|\le Cg^2
 \quad(0<c'<c).
\end{equation}
Every same-sign inverse $V=(I-tQ)^{-1}$ or its conjugate satisfies
\begin{equation}\label{uvp:eq:stable}
 \sum_x e^{c d(x,0)/g^2}|V(x)|\le C,\qquad 0\le t\le1.
\end{equation}
All constants are independent of $g,L$.
\end{lemma}

\begin{proof}[\textbf{Proof of Lemma~\ref{uvp:prop:finite}}]
We prove the comparison also at positive height for later use.
Keep hopping $A_L$ (or $A_{\mathbb Z}$), allow covariance $u/g^2$,
and, for fixed $C_0>0$, assume
\begin{equation}\label{uvp:eq:height-domain}
 |x|\le2-\kappa,\quad \tfrac12\le u\le1,\quad
 z=x+iy,\quad 0\le y\le C_0g^{-2},\quad
 g\ge C_{\kappa,C_0},\quad L/g^2\ge R_{\kappa,C_0}.
\end{equation}
Write $m_{u,\infty}(z)$ and $m_u(z)$ for the line and cycle solutions,
and $M_{B,u}(z)=(A_L-z-ug^{-2}m_u(z))^{-1}$.
For $\mu$ in a fixed small disk about $i/\sqrt{4-x^2}$, put
$w=z+ug^{-2}\mu$. The functions $f,\rho$ in the statement satisfy
\[
 f(w)=\frac{i}{\sqrt{4-x^2}}+O(g^{-2}),\qquad
 |ug^{-2}f'(w)|\le Cg^{-2},\qquad
 \Im w\asymp-\log|\rho(w)|\asymp g^{-2}.
\]
These estimates follow from analyticity on the upper branch over the
buffered bulk. The exact finite trace
$f_L(w)=f(w)(1+\rho(w)^L)/(1-\rho(w)^L)$ gives
\begin{equation}\label{uvp:eq:fixed-point-derivative}
 |(f_L-f)(w)|\le Ce^{-cL/g^2},\qquad
 |\partial_\mu(f_L-f)(w)|\le C(1+L/g^2)e^{-cL/g^2}.
\end{equation}
Increasing the two thresholds in \eqref{uvp:eq:height-domain} makes
both Dyson maps contractions on the disk. Their fixed points are the
Herglotz solutions by \eqref{da:eq:ward}, with boundary values obtained
as $y\downarrow0$. Subtracting the fixed-point equations gives
$|m_u-m_{u,\infty}|\le Ce^{-cL/g^2}$ and $\Im m_u,\Im m_{u,\infty}\asymp1$.

The exact periodic row is
$M_{B,u}(z)_{0r}=\sum_{n\in\mathbb Z}f(w_L)\rho(w_L)^{|r+nL|}$,
where $w_L=z+ug^{-2}m_u(z)$.
Its line summand has $\mu$-derivative bounded by
$Cg^{-2}(1+|r+nL|)e^{-c|r+nL|/g^2}$.
Replacing $m_u$ by $m_{u,\infty}$ therefore costs at most
$C(1+L/g^2)^Ce^{-cL/g^2}$ pointwise and $g^2$ times this in
exponentially weighted absolute row norm. Geometric summation and Ward
give, for a fixed sufficiently small $c'>0$,
\begin{equation}\label{uvp:eq:height-rows}
 \begin{gathered}
 \Im m_u(z)\asymp1,\qquad |M_{B,u}(z)_{0r}|\le Ce^{-cd(r,0)/g^2},\qquad
 \sum_r e^{c'd(r,0)/g^2}|M_{B,u}(z)_{0r}|\le Cg^2,\\
 \frac{u}{g^2}\sum_r|M_{B,u}(z)_{0r}|^2
 =\frac{ug^{-2}\Im m_u(z)}{y+ug^{-2}\Im m_u(z)}\le1.
 \end{gathered}
\end{equation}
When squaring the rows, the coefficient $u/g^2$ cancels their row
length $g^2$. The difference between the square of the periodized row
and the periodization of its square consists of distinct-image products,
for which $|r+nL|+|r+mL|\ge L$ if $n\ne m$. Their pointwise sum,
after multiplication by $u/g^2$, is at most
$Cg^{-2}(1+L/g^2)^Ce^{-cL/g^2}$, and their weighted absolute row
is at most $C(1+L/g^2)^Ce^{-cL/g^2}$.
These bounds apply to both modulus squares and complex squares.
At $u=1,y=0$, Ward gives row one for the mixed kernels;
this proves \eqref{uvp:eq:finiteerror} and \eqref{da:eq:uniform-free-row}.

For same-sign stability let $w_\infty=z+ug^{-2}m_{u,\infty}(z)$,
$R(p)=(2\cos p-w_\infty)^{-1}$, and
$Q_\infty(r)=ug^{-2}f(w_\infty)^2\rho(w_\infty)^{2|r|}$. Then
$\widehat Q_\infty(k)=ug^{-2}\int R(p)R(k-p)\,dp/(2\pi)$.
Since $\Im R(p)>0$,
\[
 \begin{aligned}
 \Re\widehat Q_\infty(k)
 &\le\frac{u}{g^2}\int(\Re R(p))^2\,\frac{dp}{2\pi}
 =\frac{ug^{-2}\int|R(p)|^2\,dp/(2\pi)+\Re\widehat Q_\infty(0)}2
 \le\tfrac12+O(g^{-2}),\\
 \widehat Q_\infty(0)&=ug^{-2}f'(w_\infty)=O(g^{-2}).
 \end{aligned}
\]
Here the first inequality uses $2ab\le a^2+b^2$, and the last uses
the line Ward mass, which is at most one as in \eqref{uvp:eq:height-rows}.
Thus $\Re(1-\tau\widehat Q_\infty)\ge c>0$ for $0\le\tau\le1$.
The profile $|Q_\infty(r)|\le Cg^{-2}e^{-c|r|/g^2}$ preserves this
gap in $|\Im k|\le c'/g^2$. Parseval on the edges of a narrower strip gives
\[
 \sum_r e^{2c''|r|/g^2}
 |[Q_\infty(I-\tau Q_\infty)^{-1}](r)|^2\le C/g^2.
\]
Weighted Cauchy--Schwarz gives a bounded exponential absolute row,
since $\sum_r e^{-c''|r|/g^2}=O(g^2)$.
The inverse $I+\tau Q_\infty(I-\tau Q_\infty)^{-1}$ has the same bound,
as does its periodization. For the actual kernel
$Q_L(r)=ug^{-2}M_{B,u}(z)_{0r}^2$, the image calculation above bounds
its difference from the periodized $Q_\infty$ by
$C(1+L/g^2)^Ce^{-cL/g^2}$ in weighted row norm.
The resolvent identity in the weighted convolution algebra is therefore
a convergent perturbation series for sufficiently large $R_{\kappa,C_0}$, giving
\begin{equation}\label{uvp:eq:height-stability}
 \sum_r e^{cd(r,0)/g^2}|[(I-\tau Q_L)^{-1}](r)|\le C,
 \qquad 0\le\tau\le1.
\end{equation}
Taking $u=1,y=0$ proves \eqref{uvp:eq:stable}; conjugation gives the other sign.
\end{proof}

\subsubsection{Deterministic terminal bounds}
\label{app:terminal-deterministic-large}

\begin{proof}[\textbf{Proof of Lemma~\ref{fo:lem:terminal-deterministic}(b)}]
Apply Lemma~\ref{uvp:prop:finite} on the buffered
interval $|e|\le2-3\kappa/4$. It gives $|m_{g,L}(e)|\le C$ and
$\Im m_{g,L}(e)\ge c$. At $u=1,y=0$, the fixed-point derivative
estimate \eqref{uvp:eq:fixed-point-derivative}, together with
$g^{-2}|f'(w)|\le Cg^{-2}$, bounds $g^{-2}|f_L'(w)|$, where
$w=e+g^{-2}m_{g,L}(e)$. Differentiating the Dyson equation gives,
after increasing both the coupling and $L/g^2$ thresholds,
\begin{equation}\label{fo:eq:terminal-anchor-derivative}
 g^{-2}|\partial_e m_{g,L}(e)|
 \le \frac{g^{-2}|f_L'(w)|}{1-g^{-2}|f_L'(w)|}
 \le C\bigl[g^{-2}+(1+L/g^2)e^{-cL/g^2}\bigr]\le\tfrac12.
\end{equation}
Thus $e\mapsto e+\delta g^{-2}\Re m_{g,L}(e)$ has derivative in
$[1/2,3/2]$ and displacement at most $Cg^{-2}$, uniformly for
$0\le\delta\le1$. Increasing $C_\kappa$ so that $Cg^{-2}\le\kappa/4$,
the endpoint images cover $[-2+\kappa,2-\kappa]$.
The intermediate value theorem gives $E_0$ in the claimed interval;
the bounds on $m_{g,L}(E_0)$ give \eqref{fo:eq:terminal-large-anchor}.
\end{proof}

\subsection{Line kernels and bounded-coupling transport}
\label{app:line-transport}

\begin{lemma}[Actual finite transport and same-sign bounds]\label{ba:lem:kernels}
Fix $C,B_0,c'>0$, and assume $L$ is sufficiently large depending only
on these constants. Suppose
\begin{equation}\label{bl:eq:intrinsic}
 g\le C,\qquad |E|\le B_0,\qquad \Im m_{g,L}(E)\ge c'.
\end{equation}
Use the finite Dyson matrix \eqref{fo:eq:characteristic} and the pair
kernels \eqref{fo:eq:PQ}. Uniformly for $0\le t<1$, with
\begin{equation}\label{ba:eq:controls}
 \ell_t=\min\left\{L,1+\frac{g}{\sqrt{\gamma_g(1-t)}}\right\},\quad
 \zeta_t=(\gamma_g W(1-t)\ell_t)^{-1},
\end{equation}
the actual mixed kernel satisfies exact row mass $(1-t)^{-1}$,
value comparable to $((1-t)\ell_t)^{-1}$ at zero, and upper exponential profile
$C'((1-t)\ell_t)^{-1}e^{-c\|x\|_L/\ell_t}$.
The free row, mixed kernel, and same-sign inverse $(I-tQ)^{-1}$ satisfy
\begin{equation}\label{ba:eq:weighted-rows}
 \begin{aligned}
 \sum_{x\in\Z_L}e^{c\|x\|_L}|M_B(0,x)|&\le C',\\
 \sum_{x\in\Z_L}e^{c\|x\|_L/\ell_t}
       |[P(I-tP)^{-1}]_{0x}|&\le \frac{C'}{1-t},\\
 \sum_{x\in\Z_L}e^{c\|x\|_L}|[(I-tQ)^{-1}]_{0x}|&\le C'.
 \end{aligned}
\end{equation}
The actual Fourier multipliers also satisfy
\begin{equation}\label{ba:eq:coercivity}
 1-\widehat P(k)\asymp g^2(1-\cos k),\qquad
 U_t(k)\coloneqq (1-t\widehat P(k))^{-1}
       \asymp((1-t)+g^2|k|_{\mathbb T}^2)^{-1}.
\end{equation}
With normalized averaging over the Fourier grid,
\begin{equation}\label{ba:eq:inverse-average}
 H\coloneqq\langle U_t\rangle
 =1+t[P(I-tP)^{-1}]_{00}\ge1,\qquad
 H\asymp((1-t)\ell_t)^{-1}.
\end{equation}
\end{lemma}
\paragraph{A common line kernel.}
On the line, at a real regular energy, write $m\coloneqq m_g(E;\gamma)$ and
$M_\infty(x)\coloneqq[(\sqrt\gamma\,gA_{\mathbb Z}-E-\gamma m)^{-1}]_{0x}$, and set
\[
 P_\infty(x)\coloneqq \gamma|M_\infty(x)|^2,\qquad
 Q_\infty(x)\coloneqq \gamma M_\infty(x)^2.
\]
On the cycle of $L$ blocks, use $P,Q$ from \eqref{fo:eq:PQ} with
$\gamma_g$ replaced by $\gamma$ and $M_B$ by $M_B(E;\gamma)$ from
\eqref{da:eq:dyson}. Products of kernels are convolutions, and all
Fourier averages are normalized. Write $\rho=e^{-q-i\theta}$,
$q>0$. Define the diffusivity, suppressing its dependence on $E,\gamma$, by
\begin{equation}\label{ba:eq:diffusivity}
 D_g\coloneqq \frac1{2(\cosh2q-1)}.
\end{equation}
\begin{lemma}[The line mixed kernel]\label{da:lem:line-kernel}
For every $0\le t<1$,
\begin{equation}\label{ba:eq:linekernel}
 \begin{gathered}
 \widehat P_\infty(k)=\frac1{1+2D_g(1-\cos k)},\qquad
 \frac{\widehat P_\infty(k)}{1-t\widehat P_\infty(k)}
 =\frac1{1-t+2D_g(1-\cos k)},\\
 [P_\infty(I-tP_\infty)^{-1}](x)=\frac{e^{-\lambda|x|}}{\sqrt{(1-t)(1-t+4D_g)}},\qquad
 \cosh\lambda=1+2^{-1}D_g^{-1}\bigl(1-t\bigr).
 \end{gathered}
\end{equation}
At $g=0$, $P_\infty(I-tP_\infty)^{-1}=(1-t)^{-1}\mathbf1_0$.
Periodizing \eqref{ba:eq:linekernel} gives a comparison kernel with
row mass $(1-t)^{-1}$ and block scale
\[
 \min\{L,1+\sqrt{D_g/(1-t)}\}.
\]
At bounded coupling in the regularity range of
Lemma~\ref{bl:lem:finite}, $D_g\asymp g^2$. At normalized large
coupling in the range of Lemma~\ref{uvp:prop:finite},
$D_g\asymp g^4$.
\end{lemma}

\noindent\textbf{Proof of Lemma~\ref{da:lem:line-kernel}.}
Write $\rho=e^{-q-i\theta}$ in \eqref{bl:eq:free}. Exact Ward gives
$\gamma|m|^2=\tanh q$. Summing
$\sum_{x\in\Z}e^{-2q|x|}e^{-ikx}$ proves the first Fourier identity
in \eqref{ba:eq:linekernel}; division by $1-t\widehat P_\infty(k)$ proves
the second. Inverse Fourier summation gives the displayed exponential
kernel and its periodization, with the stated limiting formula at $g=0$.
For $g>0$, put $w=E+\gamma m=2\sqrt\gamma\,g\cosh(q+i\theta)$.
Since $D_g=(4\sinh^2q)^{-1}$,
\[
 (\Im w)^2\le4\gamma g^2\sinh^2q\le|w|^2,\qquad
 \frac{\gamma g^2}{|w|^2}\le D_g\le\frac{\gamma g^2}{(\Im w)^2}.
\]
At bounded coupling, Lemma~\ref{bl:lem:finite} transfers the finite
density margin to $\Im m_g(E;\gamma)\ge c'/2$ for sufficiently large $L$.
Together with $|E|\le B_0$, $|m|\le\gamma^{-1/2}$ from
\eqref{da:eq:ward}, and the fixed positive variance interval, this gives
$|w|\le C'$, $\Im w\ge c$ and hence $D_g\asymp g^2$.
The uncapped block scale is therefore comparable to that in
\eqref{fo:eq:scale}, since $\gamma_g\asymp1$.
At large coupling, the proof of Lemma~\ref{uvp:prop:finite}
already gives $q\asymp g^{-2}$; substituting this into
\eqref{ba:eq:diffusivity} proves $D_g\asymp g^4$.

\begin{lemma}[Auxiliary increment and spatial-moment bounds]
\label{da:lem:transport-refinements}
Use the parameters, kernels and scales of Lemma~\ref{ba:lem:kernels}.
For $1-t\ge W^{-3}$, a fixed polynomial upper bound on $L$, and every accuracy $A$,
$L\ge C_A\log W$ gives, for $x,y\in\Z_L$,
\begin{equation}\label{gc:eq:gradient}
 \begin{split}
 &\bigl|[P(I-tP)^{-1}]_{0,x+y}-[P(I-tP)^{-1}]_{0,x}\bigr|\\
 &\quad\le\frac {C'}{(1-t)\ell_t}\min\{\|y\|_L/\ell_t,1\}
 \left(e^{-c\|x\|_L/\ell_t}+e^{-c\|x+y\|_L/\ell_t}\right)+W^{-A}.
 \end{split}
\end{equation}
For every fixed integer $j\ge0$, the difference of the $j$th time
derivatives satisfies \eqref{gc:eq:gradient} with
$C'/((1-t)\ell_t)$ replaced by $C'_j/((1-t)^{j+1}\ell_t)$ and
the same remainder $W^{-A}$, choosing $C_A$ also depending on $j$.
For every fixed $r>0$,
\begin{equation}\label{da:eq:stable-spatial-moments}
 \sum_{x\ne0}\|x\|_L^r|[(I-tQ)^{-1}]_{0x}|\le C'_r g^2.
\end{equation}
\end{lemma}

\begin{proof}[\textbf{Proof of Lemmas~\ref{ba:lem:kernels} and \ref{da:lem:transport-refinements}}]
We prove the estimates for every fixed
$\max\{1,C\}^{-2}\le\gamma\le1$, with $g\le C$, $|E|\le B_0$ and
$\Im m_{g,L}(E;\gamma)\ge c'$. Use the finite kernels defined above
and the scales \eqref{ba:eq:controls}, replacing $\gamma_g$ by $\gamma$
in $\zeta_t$; the constants are uniform over this variance interval.
The choice $\gamma=\gamma_g$ gives the two stated lemmas.
The free-row bound in \eqref{ba:eq:weighted-rows} is
\eqref{da:eq:finite-regularity} under these conditions.
Write $m\coloneqq m_{g,L}(E;\gamma)$ and set $Q(x)=\gamma M_B(x)^2$, $Q'\coloneqq Q-\gamma m^2\mathbf1_0$. Exact Ward gives
$\|Q'\|_1=1-\gamma|m|^2$. For $\Im m\ge a>0$,
\begin{equation}\label{ba:eq:gap}
 |1-t\gamma m^2|^2=(1-t\gamma|m|^2)^2+4t\gamma(\Im m)^2,
 \quad |1-t\gamma m^2|-t(1-\gamma|m|^2)\ge\gamma a^2/2.
\end{equation}
For $t\le1/2$ use the triangle inequality and a gap at least $1/2$;
for $t\ge1/2$ rationalize the square-root difference in
\eqref{ba:eq:gap}. Thus the complete, convergent expansion
\begin{equation}\label{ba:eq:stableinverse}
 (I-tQ)^{-1}=(1-t\gamma m^2)^{-1}
       \sum_{j\ge0}\left(\frac{tQ'}{1-t\gamma m^2}\right)^{*j}
\end{equation}
has absolute row $O((\gamma a^2)^{-1})$. A sufficiently small exponential
weight changes the gap by less than half. At small $g$, the
off-diagonal weighted row of $Q$ is $O(g^2)$, so the nonconstant
part of the full inverse is $O(g^2)$ too. The remaining compact
range bounded away from zero uses this proved weighted bound
with an enlarged fixed constant. This proves
\eqref{da:eq:stable-spatial-moments}. No fixed-depth truncation is used.

On the actual torus, for $L\ge L_0$,
\begin{equation}\label{ba:eq:coercivityinput}
 P(1)+P(-1)\ge cg^2,\quad
 \sum_{x\ne0}e^{c\|x\|_L}P(x)\le C'g^2,\quad
 \inf_k\widehat P(k)\ge c.
\end{equation}
For small $g$, the free Neumann and implicit expansions give
$m=m_0+O(g^2)$,
$M_B=m_0\mathbf1_0-\sqrt\gamma\,gm_0^2A_L(0,\cdot)+O_{\ell^1_c}(g^2)$,
$|m_0|=\gamma^{-1/2}$. The zero adjacency diagonal removes the linear scalar
term, and nearest-neighbor squares are at least $g^2/4$.
For $g\ge c>0$, line comparison and the explicit geometric
nearest-neighbor coefficient and positive symbol prove
\eqref{ba:eq:coercivityinput} by choosing a fixed larger $L_0$.
This separation avoids dividing an absolute comparison error by
an arbitrarily small $g$.
Evenness and row one now give exactly on the Fourier grid
\eqref{ba:eq:coercivity}.
Use nearest neighbors for the lower bound and
$1-\cos(kx)\le x^2(1-\cos k)$ for the upper bound. In particular
the zero mode is exact.

Lift $P$ to an even probability on minimal integer representatives,
splitting the Nyquist mass for even $L$. Its symbol $p^\sharp$ obeys
$\Re p^\sharp(k+iy)-p^\sharp(k)\le C'g^2y^2$.
With $|y|=c/(1+g/\sqrt{1-t})$ the denominator remains at least
$c((1-t)+g^2|k|^2)$ in real part. Contour shifting and
$\int dk/((1-t)+g^2k^2)\asymp[(1-t)((1-t)+g^2)]^{-1/2}$
give an exponential line-lift profile. Its convergent nonnegative
convolution series periodizes exactly to the actual $P(I-tP)^{-1}$;
summing translates proves the capped upper profile and weighted row.
The positive grid symbol in \eqref{ba:eq:coercivityinput} gives
\[
 \begin{gathered}
 [P(I-tP)^{-1}]_{00}
 \asymp L^{-1}\sum_k((1-t)+g^2|k|^2)^{-1}\\
 \asymp[(1-t)((1-t)+g^2)]^{-1/2}+(L(1-t))^{-1}
 \asymp((1-t)\ell_t)^{-1}.
 \end{gathered}
\]
Fourier inversion gives
$\langle U_t\rangle=1+t[P(I-tP)^{-1}]_{00}$.
For $t\ge1/2$ the displayed diagonal estimate proves
\eqref{ba:eq:inverse-average}; for $t\le1/2$, both $H$ and
$((1-t)\ell_t)^{-1}$ are comparable to one, since $g$ is bounded.

\noindent\textbf{Proof of \eqref{gc:eq:gradient} and its time derivatives.}
Compare with the periodized geometric
kernel of \eqref{ba:eq:linekernel}. Its discrete difference is
$C'((1-t)\ell_t)^{-1}\min\{\|h\|_L/\ell_t,1\}$ times the two endpoint
profiles. This also holds when $\ell_t\asymp1$, by applying the pointwise
bound at both endpoints.
The probability-kernel resolvent identity bounds the actual/comparison
row difference by
$(1-t)^{-2}\|P-\widetilde P_L\|_1$.
Here $\widetilde P_L$ is the periodized line square kernel, not
the square of the periodized line resolvent. Lemma~\ref{bl:lem:finite}
bounds their difference by $C'(1+L)^{C'}e^{-cL}$.
Choose $C_A$ after all polynomial cutoffs and derivative orders.
Finally $\partial_t^j[P(I-tP)^{-1}]=j![P(I-tP)^{-1}]^{j+1}$ has row
$j!(1-t)^{-j-1}$. Apply \eqref{gc:eq:gradient}, with accuracy
$A+3j+1$, to one convolution factor and the mixed weighted row bound
in \eqref{ba:eq:weighted-rows} to the other $j$ factors. Their weighted
row is at most $(C'/(1-t))^j$, which gives the stated derivative
prefactor. The error is at most $j!(1-t)^{-j}W^{-A-3j-1}\le W^{-A}$
for sufficiently large $W$, since $1-t\ge W^{-3}$.
Telescoping the same convolution identity gives comparison error
$C'_j(1-t)^{-j-2}\|P-\widetilde P_L\|_1$.
Exact row mass and Ward projection
always use the actual kernel.
\end{proof}

\subsection{Geometric kernels and uniform perturbations}
\label{app:geometric-kernels}
The following estimates separate the leading geometric mixed kernel
from its finite-cycle error. Here $a\ge1$ is an independent kernel
scale, $L\ge a$, $0<\delta\le1$, and
$\ell=\min\{L,a/\sqrt\delta\}$. All convolution norms use counting
measure on $\Z$ or $\Z_L$, and $\Delta f(x)=f(x+1)-f(x)$.

\begin{lemma}[A periodized geometric mixed kernel]\label{uvp:lem:geom}
Let $D\asymp a^2$ and
\[
 P^{\wedge}(q)=\frac1{1+2D(1-\cos q)},\quad t=1-\delta,
 \quad \Xi=tP(I-tP)^{-1}.
\]
On the line,
\[
 \Xi(x)=\frac{t e^{-\lambda|x|}}{\sqrt{\delta(\delta+4D)}},
 \qquad \cosh\lambda=1+\frac{\delta}{2D},
 \qquad \lambda\asymp\frac{\sqrt\delta}{a}.
\]
On $\Z$ and after periodization to $\Z_L$,
\begin{equation}\label{da:eq:geometric-kernel-bounds}
 \begin{gathered}
 \|\Xi\|_1\le C\delta^{-1},\qquad
 \|\Xi\|_\infty\le C(\delta\ell)^{-1},\\
 \|\Delta\Xi\|_1\le C\lambda\delta^{-1},\qquad
 \|\Delta\Xi\|_\infty\le C\lambda(\delta\ell)^{-1},\qquad
 \|\Delta^2\Xi\|_1\le C\lambda^2\delta^{-1}.
 \end{gathered}
\end{equation}
\end{lemma}

\begin{proof}[\textbf{Proof of Lemma~\ref{uvp:lem:geom}}]
The formula follows by solving the second-order difference equation
or summing the geometric Fourier series. Each row of $\Xi$ sums to $t/\delta$.
Periodization has supremum at most
$C(\delta^{-1}\lambda+\delta^{-1}/L)
 \le C/(\delta\ell)$.
A unit first difference is bounded in absolute value by $C\lambda$
times the sum of two translates of the original positive kernel.
This proves both first-difference estimates. Away from zero, a
second difference is $O(\lambda^2)$ times a translate of the kernel.
The two points adjacent to the cusp contribute at most
$C\lambda\Xi(0)\le C\lambda^2/\delta$ in total. Thus each absolute row
sum of $\Delta^2\Xi$ is $O(\lambda^2/\delta)$. Periodization preserves
the bound on the absolute row sums and commutes with differences. Notice that a pointwise
second-difference bound by $\lambda^2\norm\Xi_\infty$ is not asserted;
it fails at the cusp.
\end{proof}

\begin{lemma}[Small finite-kernel perturbations]\label{uvp:lem:perturb}
Let $P_0$ be the periodized kernel of Lemma~\ref{uvp:lem:geom}. Suppose
$P_1$ is an even nonnegative probability kernel and
\[
 \norm{P_1-P_0}_1\le\varepsilon,\qquad
 \norm{P_1-P_0}_\infty\le\varepsilon/a.
\]
If $\Xi_j=tP_j(I-tP_j)^{-1}$ and $\varepsilon\le c\delta^2$, then
$E=\Xi_1-\Xi_0$ satisfies
\begin{equation}\label{da:eq:abstract-kernel-error}
 \|E\|_1\le C\varepsilon\delta^{-2},\qquad
 \|E\|_\infty\le C\varepsilon(a\delta^2)^{-1}.
\end{equation}
\end{lemma}

\begin{proof}[\textbf{Proof of Lemma~\ref{uvp:lem:perturb}}]
The exact resolvent identity is
\[
 E=t(I-tP_1)^{-1}(P_1-P_0)(I-tP_0)^{-1}.
\]
Both inverse kernels are nonnegative and have row $\delta^{-1}$.
The row and pointwise Young inequalities give respectively
$\norm E_1\le\varepsilon\delta^{-2}$ and
$\norm E_\infty\le\varepsilon/(a\delta^2)$.
\end{proof}

\begin{lemma}[Mixed kernels uniformly in the coupling]\label{uvp:lem:actual-kernels}
Use the parameters and kernels of Lemma~\ref{uvp:prop:finite}.
Let $t=1-\delta$, $0<\delta\le1$, and
$\ell=\min\{L,g^2/\sqrt\delta\}$. If
$L/g^2\ge C(1+|\log\delta|)$, then the actual mixed inverse has the
decomposition
\begin{equation}\label{da:eq:actual-kernel-decomposition}
 (I-tP)^{-1}=I+\Xi+\mathcal E,\qquad
 \Xi=tP_0(I-tP_0)^{-1},\qquad \varepsilon=\varepsilon_L\le c\delta^2.
\end{equation}
Here $\varepsilon_L$ is defined in \eqref{uvp:eq:finiteerror},
$\Xi$ and $\mathcal E$ are even, and $d$ denotes cyclic distance.
The stable inverses satisfy \eqref{uvp:eq:stable}.
The actual mixed kernel satisfies
\begin{equation}\label{da:eq:uniform-mixed-row}
 0\le[P(I-tP)^{-1}](x)\le\frac{C}{\delta\ell},\qquad
 \sum_x[P(I-tP)^{-1}](x)=\delta^{-1}.
\end{equation}
With $\lambda\coloneqq\sqrt\delta/g^2$, the mixed part satisfies
\begin{align}
 \norm{\Xi}_1&\le C\delta^{-1},&
 \norm{\Xi}_\infty&\le C(\delta\ell)^{-1},\label{uvp:eq:xi}\\
 \norm{\Delta\Xi}_1&\le C\lambda\delta^{-1},&
 \norm{\Delta\Xi}_\infty&\le C\lambda(\delta\ell)^{-1},\label{uvp:eq:xi1}\\
 \norm{\Delta^2\Xi}_1&\le C\lambda^2\delta^{-1}.&&\label{uvp:eq:xi2}
\end{align}
The last bound also holds for any composition of two forward or
backward unit differences.
The error obeys
\begin{equation}\label{uvp:eq:error}
 \norm{\mathcal E}_1\le C\varepsilon\delta^{-2},\qquad
 \norm{\mathcal E}_\infty\le C\varepsilon(g^2\delta^2)^{-1},\qquad
 \varepsilon\le c\delta^2.
\end{equation}

In particular, $L\ge Cg^2\log W$ suffices when $\delta\ge W^{-3}$,
with no upper bound on $g,L$. 
The characteristic scales in \eqref{fo:eq:scale} satisfy
\begin{equation}\label{da:eq:primitive-scale-conversion}
 \ell\le\ell_t\le2\ell,\qquad
 \zeta_t\asymp\frac{g^2}{W\delta\ell},\qquad
 A\coloneqq\frac{g^2}{\delta\ell}\ge1.
\end{equation}
The actual and comparison inverses also satisfy
\begin{equation}\label{da:eq:weighted-actual-kernels}
 \begin{gathered}
 \sum_x e^{c d(x,0)/\ell}
    \bigl(|(I+\Xi)(x)|+|[(I-tP)^{-1}](x)|\bigr)\le C\delta^{-1},\\
 \sum_x e^{c d(x,0)/\ell}|\mathcal E(x)|\le C\varepsilon_L\delta^{-2}.
 \end{gathered}
\end{equation}
For every fixed $r,Q>0$, there is $C_{r,Q}$ such that, for sufficiently
large $W$,
\begin{equation}\label{da:eq:kernel-accuracy}
 \delta\ge W^{-3},\quad L\ge C_{r,Q}g^2\log W
 \quad\Longrightarrow\quad \delta^{-r}\varepsilon_L\le W^{-Q}.
\end{equation}
Finally, for any kernels satisfying
\eqref{uvp:eq:stable} and \eqref{uvp:eq:xi}--\eqref{uvp:eq:error}, with
$\varepsilon\le c\delta^2$, the local convolution bound is
\begin{equation}\label{uvp:eq:localconv}
 \|\mathcal E\|_1\le C,\qquad
 \sup_x\sum_y e^{-c d(y,0)/g^2}|K(x-y)|\le CA,
 \quad K\in\{I,\Xi,\mathcal E,\text{stable kernels}\}.
\end{equation}
\end{lemma}

\begin{proof}[\textbf{Proof of Lemma~\ref{uvp:lem:actual-kernels}}]
Writing $H=L/g^2$, absorb $(1+H)^C$ into $Ce^{cH/2}$ in
\eqref{uvp:eq:finiteerror}. The assumed lower bound
$H\ge C(1+|\log\delta|)$ then gives $\varepsilon_L\le c\delta^2$.
Use Lemma~\ref{da:lem:line-kernel} at $\gamma=g^{-2}$ and
Lemmas~\ref{uvp:lem:geom}--\ref{uvp:lem:perturb} with $a=g^2$,
$D\asymp g^4$, and $\varepsilon=\varepsilon_L\le c\delta^2$.
Then $(I-tP)^{-1}=I+tP(I-tP)^{-1}$ gives
\eqref{da:eq:actual-kernel-decomposition}, and the geometric and
perturbation bounds give \eqref{uvp:eq:xi}--\eqref{uvp:eq:error},
with decay rate $\lambda\asymp\sqrt\delta/g^2$.
The stable row is \eqref{uvp:eq:stable}.

\noindent\textbf{Proof of \eqref{da:eq:uniform-mixed-row} and
\eqref{da:eq:primitive-scale-conversion}.}
The convolution series of the probability kernel $P$ gives positivity
and row mass $\delta^{-1}$. For $t\ge1/2$, use
$P(I-tP)^{-1}=t^{-1}(\Xi+\mathcal E)$ and
$\|\mathcal E\|_\infty\le C\varepsilon_L/(g^2\delta^2)
\le C/g^2\le C/(\delta\ell)$.
For $t\le1/2$, the same series and $\|P\|_\infty\le C/g^2$ give
$\|P(I-tP)^{-1}\|_\infty\le C/(g^2\delta)$; here $\ell\asymp g^2$.
Also \eqref{fo:eq:scale} gives
$\ell\le\ell_t\le2\ell$, $\zeta_t\asymp g^2/(W\delta\ell)$,
and $\delta\ell\le g^2\sqrt\delta\le g^2$, hence $A\ge1$.

\noindent\textbf{Proof of \eqref{da:eq:weighted-actual-kernels}.}
In the weighted convolution norm
$\|K\|_{1,c/\ell}=\sum_xe^{cd(x,0)/\ell}|K(x)|$,
the explicit geometric formula in Lemma~\ref{uvp:lem:geom} gives
$\|R_0\|_{1,c/\ell}\le C\delta^{-1}$ for $R_0=(I-tP_0)^{-1}$,
by choosing $c/\ell$ below the decay rate when uncapped; at the cap
$\ell=L$, the weight is bounded and the row sum is $\delta^{-1}$.
Since $\ell\ge g^2$, \eqref{uvp:eq:finiteerror} gives
$\|P-P_0\|_{1,c/\ell}\le C\varepsilon_L$. Thus
\[
 (I-tP)^{-1}=R_0\sum_{j\ge0}[t(P-P_0)R_0]^{*j},\qquad
 \mathcal E=t(I-tP)^{-1}(P-P_0)R_0
\]
converges for $\varepsilon_L\le c\delta^2$ and yields weighted
inverse row $C\delta^{-1}$ and error row $C\varepsilon_L\delta^{-2}$.

\noindent\textbf{Proof of \eqref{da:eq:kernel-accuracy} and \eqref{uvp:eq:localconv}.}
For fixed $r,Q$, choose $C_{r,Q}$ sufficiently large in
$L\ge C_{r,Q}g^2\log W$; then
\[
 \delta^{-r}\varepsilon_L
 \le CW^{3r}(1+L/g^2)^Ce^{-cL/g^2}
 \le C'W^{3r}e^{-cL/(2g^2)}\le W^{-Q}.
\]
Thus $L/(g^2\log W)\to\infty$ permits every fixed accuracy,
without an upper bound on $g,L$.
For local convolution, the exponential summation costs $Cg^2$;
use $g^2\|\Xi\|_\infty\le CA$ and
$g^2\|\mathcal E\|_\infty\le C\varepsilon\delta^{-2}\le C$.
The identity and stable kernels use their bounded rows and $A\ge1$;
\eqref{uvp:eq:error} also gives $\|\mathcal E\|_1\le C$.
This argument uses only the stated abstract kernel bounds.
\end{proof}

\subsection{Scales and deterministic propagators}\label{app:main-propagators}
Under \eqref{da:eq:propagation-domain}, the following auxiliary bounds
control time ratios and propagator row sums in the continuation argument.
Use $1/2\le s\le u\le t<1$ and $1-t\ge W^{-3}$:
\begin{equation}\label{gc:eq:scale}
 \begin{gathered}
 1\le\ell_s\le\ell_u\le L,\qquad
 (W\eta_u)^{-1}=(\Im m)^{-1} \ell_u\zeta_u,\\
 c_*W^{-1}\le\zeta_s\le\zeta_u
       \le\frac{1-s}{1-u}\zeta_s.
 \end{gathered}
\end{equation}
\begin{equation}\label{gc:eq:propagator-input}
 \begin{gathered}
 \sum_b(U_{s,u}^{+-})_{ab}=\frac{1-s}{1-u},\qquad
 (U_{s,u}^{+-})_{ab}\ge0,\\
 \max_a\sum_b e^{c\|a-b\|_L/\ell_0}
       \bigl(|(U_{s,u}^{++})_{ab}|+|(U_{s,u}^{--})_{ab}|\bigr)\le C'.
 \end{gathered}
\end{equation}

The appendix arguments also use the following strengthened versions.
Here $\ell_0$ is the time-zero scale in \eqref{fo:eq:scale}; in the
stated parameter ranges, $\ell_0\asymp\max\{1,g^2\}$ and
$1\le\ell_0\le\ell_s$:
\begin{equation}\label{da:eq:weighted-propagation-dyson}
 \gamma_g\max_a\sum_b e^{c\|a-b\|_L/\ell_0}
       (|M_{B,ab}|+|M_{B,ab}|^2)\le C',\qquad
 \max_a\sum_b e^{c\|a-b\|_L/\ell_0}|[(I-uQ)^{-1}]_{ab}|\le C'.
\end{equation}
The mixed kernel has row mass $(1-u)^{-1}$ and, for every fixed $A>0$,
\begin{equation}\label{da:eq:mixed-profile-increment}
 \begin{aligned}
 &|[P(I-uP)^{-1}]_{0,x+h}-[P(I-uP)^{-1}]_{0,x}|\\
 &\quad\le\frac {C'}{(1-u)\ell_u}
      \min\{\ell_u^{-1}\|h\|_L,1\}
      (e^{-c\|x\|_L/\ell_u}+e^{-c\|x+h\|_L/\ell_u})
      +C_A\ell_u^{-1}W^{-A}.
 \end{aligned}
\end{equation}

Let $A_{\mathbb Z}$ be the nearest-neighbor adjacency operator on
$\mathbb Z$, and define the periodized line square kernel by
\[
 P^{\rm per}(x)=\gamma_g\sum_{j\in\mathbb Z}
 \left|\left[(\sqrt{\gamma_g}gA_{\mathbb Z}-E-
           \gamma_g m_g(E))^{-1}\right]_{0,x+jL}\right|^2.
\]
For $r_u=P(I-uP)^{-1}-P^{\rm per}(I-uP^{\rm per})^{-1}$ and every
prescribed $A>0$, uniformly for $s\le u\le t$,
\begin{equation}\label{gc:eq:kernel-remainder}
 \max_a\sum_b|(r_u)_{ab}|+\ell_u\max_{a,b}|(r_u)_{ab}|\le C_AW^{-A}.
\end{equation}
The same remainder bound holds for its spatial increments and
propagator version $(u-s)r_u$. Constants are independent of $g,L$; the width threshold
may depend on the rate of divergence in the volume lower bound.

\begin{proof}[\textbf{Proof of Lemma~\ref{da:lem:scales-propagators} and
\eqref{gc:eq:scale}--\eqref{gc:eq:kernel-remainder}}]
Use \eqref{da:eq:finite-regularity} with $\gamma_g\asymp\ell_0\asymp1$
at bounded coupling, and Lemma~\ref{uvp:prop:finite} together with
\eqref{da:eq:uniform-free-row} with
$\ell_0\asymp g^2$, $\gamma_gg^2=1$ at large coupling. They give
$\Im m\asymp1$, bounded entries, and
\[
 \gamma_g\max_a\sum_b e^{c\|a-b\|_L/\ell_0}|M_{B,ab}|\le C,
 \qquad |M_{B,ab}|^2\le C|M_{B,ab}|.
\]
Together with the stable rows \eqref{ba:eq:weighted-rows} and
\eqref{uvp:eq:stable}, these prove
\eqref{da:eq:weighted-propagation-dyson} and hence
\eqref{fo:eq:propagation-dyson}.
The scale definition also gives $1\le\ell_0\le\ell_s$.
Since $U_{s,u}^{++}=I+(u-s)u^{-1}((I-uQ)^{-1}-I)$ and $u\ge1/2$,
its weighted row is bounded. Ward gives $P\mathbf1=\mathbf1$,
so $U_{s,u}^{+-}\mathbf1=(1-s)(1-u)^{-1}\mathbf1$ with positive
entries, proving \eqref{gc:eq:propagator-input}.
The mixed profile \eqref{fo:eq:kernel-bounds} and the increment
bound \eqref{da:eq:mixed-profile-increment} follow from Lemmas~\ref{ba:lem:kernels} and
\ref{da:lem:transport-refinements} at bounded coupling, and from
the periodized geometric formula \eqref{ba:eq:linekernel} and
Lemma~\ref{uvp:lem:actual-kernels} at large coupling; the finite
remainder is checked below.
For $\theta_u=(1-u)P(I-uP)^{-1}$, direct differentiation gives
\begin{equation}\label{da:eq:profile-derivative}
 \partial_u\theta_u=P(I-uP)^{-1}\theta_u-(1-u)^{-1}\theta_u.
\end{equation}
For \eqref{gc:eq:scale}, the reciprocal scale
\[
 (\gamma_gW\zeta_u)^{-1}
 =\min\{L(1-u),(1-u)+g\gamma_g^{-1/2}\sqrt{1-u}\}
\]
is decreasing, whereas $\ell_u$ is increasing. Thus
$\zeta_s\le\zeta_u\le(1-s)(1-u)^{-1}\zeta_s$;
$\gamma_g(1-s)\ell_s\le C$ gives $\zeta_s\ge c/W$.
The identity for $\eta_u$ uses $\eta_u=(1-u)\gamma_g\Im m$.

\noindent\textbf{Proof of \eqref{gc:eq:commutator} and \eqref{da:eq:commutator-bound}.}
Write $K_v=P(I-vP)^{-1}$ and
$\mathcal Q_vA=A-(\operatorname{row}A)\prod_{j=2}^n\theta_v(a_1-a_j)$.
Symmetry and row mass give
\[
 \operatorname{row}\!\left(\sum_{i=1}^nK_v^{(i)}A\right)
 =K_v(\operatorname{row}A)+(n-1)(1-v)^{-1}\operatorname{row}A.
\]
Differentiate $\mathcal Q_vA$ using \eqref{da:eq:profile-derivative}:
the normalization terms cancel the last term above, and the
nonanchored generators cancel, leaving precisely the first-leg
commutator \eqref{gc:eq:commutator}.
Its bracket has row $1-1=0$; bounding each product profile by
$(C/\ell_v)^{n-1}$ and summing the kernel row $(1-v)^{-1}$ proves
\eqref{da:eq:commutator-bound}.

\noindent\textbf{Proof of \eqref{gc:eq:kernel-remainder}.}
At large coupling $\mathcal E=ur_u$, with $u\ge1/2$;
at bounded coupling use the inverse subtraction in the proof of
\eqref{gc:eq:gradient}. Thus
\[
 \max_a\sum_b|(r_u)_{ab}|+\ell_u\max_{a,b}|(r_u)_{ab}|
 \le\frac{C}{(1-u)^2}
 \begin{cases}
 (1+L)^{q+1}e^{-cL},&\text{bounded coupling},\\
 (1+\ell_u/g^2)\varepsilon_L,&\text{large coupling}.
 \end{cases}
\]
Here $(1-u)^{-1}\le W^3$, bounded-coupling $L$ is polynomial,
and $\ell_u/g^2\le CW^{3/2}$ at large coupling.
The volume lower bounds and \eqref{da:eq:kernel-accuracy} absorb these
costs at any fixed accuracy, chosen after the time schedule and loop
and moment orders. Spatial increments are differences of translates;
the propagator remainder is $(u-v)r_u$, $s\le v\le u$, so both
satisfy the same bound.
\end{proof}

\subsection{Spatial summation and weights}\label{app:scale-budgets}
\begin{proof}[\textbf{Proof of Lemma~\ref{da:lem:scale-budgets}}]
The spatial budget follows from polynomial $L$ at bounded coupling,
and from $\ell_u\ge g^2$, $L/g^2\le CW^2$ at large coupling.
The normalization identity in \eqref{fo:eq:spatial-sum-budget} follows
directly from the definition \eqref{fo:eq:scale}.

For \eqref{da:eq:weight-comparison} with $0<\theta\le1/2$, use the
cyclic triangle inequality,
$|x^\theta-y^\theta|\le|x-y|^\theta$, and the common nonnegative floor.
On mesocells
$n/\ell_u\le W^{-B}$, the ratio tends to one.
For \eqref{da:eq:weight-convolution}, retain the gain
$x^\theta+y^\theta-(x+y)^\theta\ge c_\theta\min(x,y)^\theta$.
Splitting the sum at the nearer endpoint leaves a summable factor of
mass $C_\theta\ell_u$. The floor terms are
$O(\ell_uW^{-F}+LW^{-2F})$.
For the square-root version use
$\sqrt{e^{-x}+W^{-F}}\le e^{-x/2}+W^{-F/2}$;
its floor product costs $LW^{-F}\le\ell_uW^{-F/2}$ by the stated
condition $LW^{-F/2}\le\ell_u$.

\end{proof}

\begin{lemma}[Weights and time integrals for the stopped pair equation]
\label{da:lem:stopped-pair-weights}
Under \eqref{da:eq:propagation-domain}, suppose
$L/\ell_v\le W^{B}$ for a fixed $B$ and $1/2\le s\le v\le u\le t\le1-W^{-3}$.
Fix $F>0$, and set $R=(1-s)/(1-t)$ and
$w_v(x)=e^{-c_{\rm p}(\|x\|_L/\ell_v)^\theta}+W^{-F}$,
where $0<\theta\le1/2$ and $c_{\rm p}>0$ are fixed.
Then
\begin{equation}\label{da:eq:short-weight}
 \sup_{\|x\|_L\le\ell_v(\log W)^{3/2}}w_v(x)^{-1/2}=W^{o(1)},
\end{equation}
and, for $U_{v,u}=U_{v,u}^{+-}$,
\begin{equation}\label{da:eq:propagated-weight}
 \sum_{a',b'}(U_{v,u})_{aa'}(U_{v,u})_{bb'}w_v(a'-b')
 \le C\left(\frac{1-v}{1-u}\right)^2 w_u(a-b).
\end{equation}
For fixed $q\ge p\ge0$,
\begin{equation}\label{da:eq:scale-integral}
 \int_s^u\frac{\zeta_v^q}{\zeta_u^p(1-v)}\,dv
 \le\zeta_u^{q-p}\log R.
\end{equation}
\end{lemma}
\begin{proof}[\textbf{Proof of Lemma~\ref{da:lem:stopped-pair-weights}}]
On the short halo the reciprocal weight is at most
$\exp(C(\log W)^{3\theta/2})=W^{o(1)}$, proving
\eqref{da:eq:short-weight}. Since $\ell_v\le\ell_u$, $w_v\le w_u$.
The triangle inequality gives
$w_u(x-y)\le w_u(x)e^{c_{\rm p}(\|y\|_L/\ell_u)^\theta}$.
The exponential part of \eqref{fo:eq:kernel-bounds}, in
$U_{v,u}=I+(u-v)P(I-uP)^{-1}$, therefore has weighted row norm
at most $C(1-v)/(1-u)$. Indeed, the resulting sum is bounded by
$C\ell_u^{-1}\sum_y e^{-c\|y\|_L/\ell_u+
 c_{\rm p}(\|y\|_L/\ell_u)^\theta}\le C$.
Apply this estimate successively in the two coordinates.
For the remainder, use $w_v\le2$, $w_u\ge W^{-F}$ and
\eqref{gc:eq:kernel-remainder}; choosing its accuracy after $B,F$
absorbs the polynomial factor $(L/\ell_u)W^F$.
This proves \eqref{da:eq:propagated-weight}.
Finally $\zeta_v\le\zeta_u$ by \eqref{gc:eq:scale}, so integrating
$(1-v)^{-1}$ proves \eqref{da:eq:scale-integral}.
\end{proof}

\begin{lemma}[Scalar gains used in loop improvement]
\label{da:lem:profile-improvement}
Let $R\ge1$, $0<\zeta_s\le\zeta_t\le R\zeta_s$, and
$R^{20}\zeta_t\le W^{-c}$ for fixed $c>0$.
For every fixed $C$, the following quantities have a strictly negative
power of $W$ as an upper bound:
\begin{equation}\label{gc:eq:ratios}
 \begin{gathered}
 (CR)^8\zeta_s\ll1,\quad
 CR\sqrt{\zeta_t}\,\zeta_s^{-1/4}\ll1,\quad
 CR\sqrt{\zeta_t}\,\zeta_s^{-1/8}\ll1,\\
 \zeta_t R^{10}\zeta_s^{-1/4}\ll1,\qquad R^2\zeta_s^{1/4}\ll1.
 \end{gathered}
\end{equation}
Also, for fixed $c_{\rm p},\theta>0$, $1\le R\le W^3$ and $x\ge0$,
there are fixed $c',C'>0$ such that
\begin{equation}\label{da:eq:profile-minimum}
 \min\{1,R^{10}e^{-c_{\rm p}x^\theta}\}
 \le e^{C'\sqrt{\log W}}e^{-c'x^{\theta/2}}.
\end{equation}
\end{lemma}
\begin{proof}[\textbf{Proof of Lemma~\ref{da:lem:profile-improvement}}]
For $b>d\ge0$ and $a+d\le20(b-d)$,
$R^a\zeta_t^b\zeta_s^{-d}\le R^{a+d}\zeta_t^{b-d}
\le W^{-(b-d)c}$. This gives the three terms with negative powers
of $\zeta_s$ in \eqref{gc:eq:ratios}; the other two follow from
$\zeta_s\le\zeta_t$. Fixed powers of $C$ preserve the strict margins.
For \eqref{da:eq:profile-minimum}, first suppose
$x^\theta\le C(1+\log R)$. Use the bound $1$ for the minimum;
$x^{\theta/2}\le C'\sqrt{\log W}$ gives the desired estimate.
If $x^\theta>C(1+\log R)$, choose $C$ sufficiently large so that
$R^{10}e^{-c_{\rm p}x^\theta}\le e^{-c_{\rm p}x^\theta/2}
\le C'e^{-c'x^{\theta/2}}$. These two cases prove the claim.
\end{proof}

For covariance-weighted sums at large coupling,
\eqref{fo:eq:spatial-sum-budget} gives total mass $L/g^2$ and halo mass
$\ell_t/g^2$. The restoring density $g^2\theta_t$ has size
$O(g^2/\ell_t)$, so $k-1$ relative-coordinate sums cost
$(L/\ell_t)^{k-1}$, a polynomial in $W$ uniformly in $g$.

\subsection{Deterministic comparison parameters and positive heights}
\label{app:matching-parameters}

\begin{proof}[\textbf{Proof of Lemma~\ref{da:lem:comparison-parameters}(a)}]
We construct the comparison energy and then check its regularity and scales.
Set $g'=g\sqrt{s/u}$ and
\[
 \xi=\sqrt{s\gamma_{g'}/(u\gamma_g)},\qquad
 a=u\gamma_g,\quad b=\gamma_g,\quad b'=(u/s)\gamma_g.
\]
Then $2^{-1/2}\le\xi\le1$, $a<b\le b'\le2a$,
$\xi^2b'=\gamma_{g'}$ and $A_{g'}=\xi A_g$.
Write $m(z;\beta)$ for the finite Dyson solution at fixed hopping
$A_g$ and variance $\beta$, and put
$x=\Re z_u(g,E)$, $\eta_0=\eta_u(g,E)$,
$\rho(y)=\Im m(x+iy;a)$.
Substitution gives $m(x+i\eta_0;a)=m(E;b)$ and
$\eta_0/\rho(\eta_0)=b-a$.
The Stieltjes formula
$\rho(y)/y=\int((\lambda-x)^2+y^2)^{-1}\,d\nu_a(\lambda)$
shows that $y/\rho(y)$ increases strictly to infinity. Thus for every
$\beta\in[b,b']$ there is a unique $\eta_\beta\ge\eta_0$ with
\begin{equation}\label{bl:eq:tube}
 \eta_\beta/\rho(\eta_\beta)=\beta-a,\qquad
 e_\beta=x-(\beta-a)\Re m(x+i\eta_\beta;a).
\end{equation}
Here $\eta_b=\eta_0$, $e_b=E$. The identity
$e_\beta+\beta m(x+i\eta_\beta;a)
=x+i\eta_\beta+a m(x+i\eta_\beta;a)$ makes
$m(x+i\eta_\beta;a)$ a positive-imaginary solution at $(e_\beta,\beta)$.
Ward at this covariance gives
$\Re(1-\beta\langle R^2\rangle)\ge2\beta\rho(\eta_\beta)^2>0$;
implicit continuation into $\C_+$ and Herglotz uniqueness identify the
root with $m(e_\beta;\beta)$. Covariance scaling now gives
\begin{equation}\label{bl:eq:exactlaw}
 E'=\xi e_{b'},\qquad m_{g',L}(E')=\xi^{-1}m(e_{b'};b'),\qquad
 w=\xi^{-1}z_s(g',E')=x+i\eta_{b'}.
\end{equation}
In particular, $\Re w=x$, $\Im w\ge\eta_0$, and
$\Im w/\eta_0=(u/s)(1-s)(1-u)^{-1}
\rho(\eta_{b'})/\rho(\eta_0)$.

\paragraph{Bounded coupling.}
Here $a\ge c>0$. Ward gives $|m(z;a)|\le a^{-1/2}$,
$|m(e_\beta;\beta)|\le\beta^{-1/2}$ and
$\eta_\beta\le(\beta-a)\beta^{-1/2}\le\sqrt2$.
Since $\|A_g\|\le2$ and $|x|\le C'+1$, imaginary parts of the
Dyson equation give $\rho(y)\ge cy$ on the segment. While
$\rho(y)\ge c'/2$, implicit differentiation and \eqref{da:eq:ward} give
$|\partial_zm(z;a)|\le2/(a^2(c')^2)$.
A first drop from $c'$ to $c'/2$ therefore needs a height increase
at least $a^2(c')^3/4$, after which $\rho(y)\ge cy$ supplies a fixed
positive lower bound. Hence $\rho\asymp1$ throughout, and
\eqref{bl:eq:exactlaw} yields
\begin{equation}\label{bl:eq:backbuffer}
 |E'|\le C'',\qquad g'=g\sqrt{s/u}\le g,\qquad
 \Im m_{g',L}(E')\ge(C'')^{-1},\qquad
 1\le\frac{\Im w}{\eta_u(g,E)}\le C''\frac{1-s}{1-u}.
\end{equation}
The same fixed energy and density bounds hold at every $e_\beta$.

\paragraph{Large coupling.}\label{app:matching-large}
Now $\xi=1$, $A_g=A_L$, $b=g^{-2}$, $b'=(g')^{-2}$.
Lemma~\ref{uvp:prop:finite} and \eqref{fo:eq:terminal-anchor-derivative},
at coupling $\beta^{-1/2}$ and margin $\kappa/2$, give
$|m(e;\beta)|\le C$, $\Im m(e;\beta)\ge c$ and
$\beta|\partial_em(e;\beta)|\le1/2$ for
$|e|\le2-\kappa/2$, $\beta\in[b,b']\subset[b,2b]$.
For $f_\beta(e)=e+(\beta-a)\Re m(e;\beta)$, differentiation of
Dyson gives $\partial_\beta m=m\partial_em$, hence
\[
 \partial_ef_\beta\in[1/2,3/2],\qquad
 \partial_\beta f_\beta=\Re m+(\beta-a)\Re(m\partial_em)=O(1).
\]
Since $f_b(E)=x$, inversion gives a unique root
$f_\beta(e)=x$ with $|e-E|\le C(\beta-b)$, still in the buffered
bulk for sufficiently large $g$. Its positive-height solution obeys
\eqref{bl:eq:tube}, so uniqueness identifies it with $e_\beta$. Thus
\begin{equation}\label{bl:eq:large-matching-bounds}
 |E'-E|\le C(b'-b),\qquad |E'|\le2-\kappa/2,\qquad
 1\le\frac{\Im w}{\eta_u(g,E)}\le C\frac{1-s}{1-u}.
\end{equation}

\paragraph{Scales and the comparison segment.}
Use $\zeta_t(g)^{-1}=W\min\{\gamma_gL(1-t),
\gamma_g(1-t)+g\sqrt{\gamma_g(1-t)}\}$.
Since $s(1-s)\ge u(1-u)$ and $\gamma_{g'}\ge\gamma_g$,
\begin{equation}\label{bl:eq:backdomain}
 \begin{gathered}
 \gamma_{g'}(1-s)+g'\sqrt{\gamma_{g'}(1-s)}
 \ge\gamma_g(1-u)+g\sqrt{\gamma_g(1-u)},\\
 \gamma_{g'}L(1-s)\ge\gamma_gL(1-u).
 \end{gathered}
\end{equation}
Thus $\zeta_s(g')\le\zeta_u(g)$, while $g'\asymp g$ and
$\gamma_{g'}\asymp\gamma_g$ give $\zeta_s(g')\asymp\zeta_s(g)$;
earlier comparison preserves both the time cutoff and scale smallness.
For $\eta_0\le y\le\eta_{b'}$, put
$\beta=a+y/\rho(y)$, $v=a/\beta\in[s,u]$,
$g_v=g\sqrt{v/u}$, $\xi_v=\sqrt{v\gamma_{g_v}/(u\gamma_g)}$
and $E_v=\xi_v e_\beta$.
Then $\xi_v^{-1}z_v(g_v,E_v)=x+iy$, with the regularity bounds just
proved. For the common Dyson matrix and its square kernel
$Q_\beta=\beta(M_{{\rm common},ab}^{\,2})_{a,b}$,
\begin{equation}\label{bl:eq:comparison-dyson}
 \begin{aligned}
 M_{\rm common}&=(A_g-x-iy-a m(x+iy;a))^{-1}\\
 &=(A_g-e_\beta-\beta m(e_\beta;\beta))^{-1}=\xi_vM_B(g_v,E_v),\\
 (I-(a/\beta)Q_\beta)^{-1}&=(I-vQ(g_v,E_v))^{-1},
 \qquad \beta\xi_v^2=\gamma_{g_v}.
 \end{aligned}
\end{equation}
\end{proof}

\begin{proof}[\textbf{Proof of Lemma~\ref{fo:lem:families}}]
We now pass from the deterministic comparison parameters to resolvents.
\paragraph{Resolvent comparisons under the Gaussian coupling.}
Lemma~\ref{da:lem:comparison-parameters}(a) gives $E',w,g',\xi$ and
their parameter bounds. Since $A_{g'}=\xi A_g$ and
$s\gamma_{g'}=\xi^2u\gamma_g$, Gaussian scaling gives,
under the coupling $H_s(g')=\xi H_u(g)$,
\begin{equation}\label{fo:eq:comparison-point}
 \begin{gathered}
 \Re w=\Re z_u(g,E),\qquad
 1\le\frac{\Im w}{\eta_u(g,E)}\le C''R,\\
 (H_u(g)-w)^{-1}=\xi G_s(g',E').
 \end{gathered}
\end{equation}
Write $G=G_u(g,E)$, $F=(H_u(g)-w)^{-1}$,
$\eta=\eta_u(g,E)$ and $y=\Im w$. Spectral calculus and the
resolvent identity give
\[
 \Im G\preceq\frac y\eta\Im F,\qquad
 G-F=i(\eta-y)GF.
\]
Using $GG^*=\eta^{-1}\Im G$ and $F^*F=y^{-1}\Im F$,
Cauchy--Schwarz on the intermediate block and orbital indices yields
\[
 \begin{aligned}
 \|G_{ab}-F_{ab}\|_{\op}
 &\le\frac{y-\eta}{\eta}
   \sqrt{\|\Im F_{aa}\|_{\op}\|\Im F_{bb}\|_{\op}},\\
 |\langle G-F\rangle_{ab}|
 &\le\frac{y-\eta}{\eta}
   \sqrt{\Im\langle F\rangle_{aa}\Im\langle F\rangle_{bb}}.
 \end{aligned}
\]
Taking maxima and using $F=\xi G_s(g',E')$ and
$y/\eta\le C''R$ proves the operator-norm and block-trace comparisons.
The first spectral-calculus inequality gives the imaginary-part comparison.
All three hold simultaneously under the same Gaussian coupling, which proves
\eqref{fo:eq:backward-matching} in the stated distributional sense.
\end{proof}

\subsubsection{Finite comparison chains}
\label{app:deterministic-iteration}
\begin{proof}[\textbf{Proof of Lemma~\ref{da:lem:comparison-parameters}(b)}]
Write $g_0$ for the target coupling and $g$ for a comparison coupling.
At bounded coupling, each application of part \textup{(a)} above
gives an energy bound and a positive lower bound on
$\Im m_{g',L}(E')$, depending only on the preceding fixed parameter
bounds, not on $W$. Since $J$ is fixed, taking the largest energy
bound and the smallest density lower bound over these finitely many
stages gives uniform constants depending only on the original bounds
and $J$; also $cg_0\le g\le C$.

The large-coupling calculation below needs only the lower-volume
condition of part \textup{(a)}, with no upper bound on $L$.
Each comparison preserves $ug^{-2}$, and moving to an earlier time
at fixed coupling decreases it. Since $u\ge1/2$ and $g'\le g$,
\begin{equation}\label{da:eq:comparison-chain-algebra}
 \begin{gathered}
 g_0^{-2}\le g^{-2}\le u^{-1}g_0^{-2}\le2g_0^{-2},\\
 \sum_{j=0}^{m-1}|E_{j+1}-E_j|
 \le C\sum_{j=0}^{m-1}(g_{j+1}^{-2}-g_j^{-2})
 =C(g_m^{-2}-g_0^{-2})\le Cg_0^{-2}.
 \end{gathered}
\end{equation}
The second line uses \eqref{bl:eq:large-matching-bounds} with margin
$\kappa/2$, up to any first exit from that buffered bulk.
Starting from $|E_0|\le2-3\kappa/4$, choose $C_\kappa$ so that
$Cg_0^{-2}<\kappa/4$ and $g_0/\sqrt2$ exceeds the input thresholds;
the displayed bound excludes the exit and proves
\eqref{fo:eq:large-input-domain}. Also
$L/(g^2\log W)\ge L/(g_0^2\log W)\to\infty$.
The scale inequality \eqref{bl:eq:backdomain} and monotonicity at
fixed coupling show that earlier comparisons preserve both
$\zeta\le W^{-c_0}$ and $1-u\ge W^{-3}$.
Finally, in the initializer range $L\le Cg_0^2W^2$, the first line of
\eqref{da:eq:comparison-chain-algebra} gives $L\le2Cg^2W^2$.
The initial and sampling inputs therefore apply with a buffered bulk
margin and twice the upper-volume constant.

\end{proof}

\subsection{Resolvent comparison and positive-height bounds}
\label{sec:bl-matching}

\begin{proof}[\textbf{Proof of Lemma~\ref{da:lem:positive-height}}]
We also prove the exponential-row refinement used below, with
$\ell_0\asymp\max\{1,g^2\}$:
\begin{equation}\label{da:eq:weighted-positive-stability}
 \begin{gathered}
 u\gamma_g\max_a\sum_b e^{c\|a-b\|_L/\ell_0}
       (|M_{B,u}(z)_{ab}|+|M_{B,u}(z)_{ab}|^2)\le C',\\
 \max_a\sum_b e^{c\|a-b\|_L/\ell_0}
 \left|\left[(I-u\gamma_g(M_{B,u}(z)_{ab}^{\,2})_{a,b})^{-1}
       \right]_{ab}\right|\le C'.
 \end{gathered}
\end{equation}
This holds on the comparison segments and on the direct large-coupling
domain in Lemma~\ref{da:lem:positive-height}.
On the characteristic, $M_{B,u}(z_u)=M_B(g,E)$; the real-energy
row estimates used in the proof of \eqref{da:eq:weighted-propagation-dyson}
therefore give \eqref{gc:eq:residual-row-bounds} also for $u<1/2$.
On a comparison segment, \eqref{bl:eq:comparison-dyson} identifies
$M_{B,u}=\xi_vM_B(g_v,E_v)$ and its stable inverse with
$(I-vQ(g_v,E_v))^{-1}$, where
$u\gamma_g\xi_v^2=v\gamma_{g_v}$, $\xi_v\asymp1$, and
$\ell_0(g_v)\asymp\ell_0(g)$.
The bounds \eqref{bl:eq:backbuffer} and
\eqref{bl:eq:large-matching-bounds} verify the same real-energy
hypotheses. This proves \eqref{da:eq:weighted-positive-stability}
and \eqref{gc:eq:positive-stability}, using $u\ge1/2$ to remove $u$
from the first row normalization.
For the direct large-coupling domain, apply
\eqref{uvp:eq:height-rows} and \eqref{uvp:eq:height-stability}
with $u=t$, $\tau=1$, bulk margin $\kappa/2$, and $C_0=C$;
the stated volume divergence supplies their fixed threshold.

The endpoint identity follows from \eqref{bl:eq:comparison-dyson} at
$v=s$. For \eqref{da:eq:height-derivative}, differentiate the scalar
Dyson equation at fixed $u,g,\Re z$:
\begin{equation}\label{da:eq:height-derivative-identities}
 \partial_r m_u=\frac{i\gamma_g L^{-1}\Tr M_{B,u}^2}
                         {1-u\gamma_g L^{-1}\Tr M_{B,u}^2},\qquad
 \partial_r M_{B,u}=\frac{i\gamma_g M_{B,u}^2}
                         {1-u\gamma_g L^{-1}\Tr M_{B,u}^2}.
\end{equation}
The second identity follows by substitution in the derivative of $M_{B,u}$.
The scalar denominator is the zero Fourier mode of the stable inverse
in \eqref{gc:eq:positive-stability}, so its reciprocal is bounded.
The absolute row bounds there give
\[
 \gamma_g\max_{a,b}|(M_{B,u}^2)_{ab}|\le C,\qquad
 \gamma_g^2\max_a\sum_b|(M_{B,u}^2)_{ab}|\le C,
\]
by convolution, proving the derivative bound. Uniqueness and the
characteristic identities follow from \eqref{da:eq:ward} and substitution
in \eqref{gc:eq:positive-dyson}, respectively.
\end{proof}

\subsubsection{Free rows and deterministic interval interpolation}\label{app:free-interpolation}
\begin{proof}[\textbf{Proof of Lemma~\ref{da:lem:free-interpolation}}]
For $w=E+m/g^2$, the assumed bounds on $m$ give
$\Im w\asymp g^{-2}$ and keep $\Re w$ in the buffered free bulk.
Hence $|\rho(w)|\le e^{-c/g^2}$ and $|f_g(w)|\le C$.
Since $L\ge4g^2$, the denominator in \eqref{bl:eq:free} satisfies
$|1-\rho(w)^L|\ge1-e^{-cL/g^2}\ge c'$. Substitution in
\eqref{bl:eq:free} gives the stronger pointwise estimate
\begin{equation}\label{da:eq:free-row-decay}
 |M_{B,bc}|\le Ce^{-cd(b,c)/g^2}.
\end{equation}
Its geometric sum proves \eqref{uvc:eq:free-row-bounds}.

We prove \eqref{da:eq:dyson-interpolation} through the general free
resolvent $R(w)=(A_L-w)^{-1}$, where $w\notin\operatorname{Spec}(A_L)$
and $|w-E|\le C/g^2$. With the coefficients defined in
\eqref{da:eq:cell-interpolation-coefficients}, for $a\in I$ and any $b$ we have
\begin{equation}\label{da:eq:free-interpolation}
 R(w)_{ab}=\sum_\alpha c_{a\alpha}R(w)_{\alpha b}+q_{ab}+e_{ab},\qquad
 |e_{ab}|\le C\frac{|I|}{g^2}(1+\max_\alpha|R(w)_{\alpha b}|).
\end{equation}
The coefficients in \eqref{da:eq:cell-interpolation-coefficients} solve
the free forced recurrence and are bounded because the bulk margin
bounds $\cos\theta$ below.
Definition~\ref{def:transfer-matrices}, with $W=1$, gives
$T(E)^j=S_ER_E^jS_E^{-1}$, whose norm is bounded uniformly for $j\in\Z$.
Variation of constants for its perturbation by $w-E$
on at most $|I|$ sites bounds the prefix error by
$C|I||w-E|e^{C|I||w-E|}\le C|I|/g^2$.
For $a\in I$, iterate the forced recurrence from $b_0,b_0-1$ to $a$;
at most one source occurs, with coefficient $q_{ab}$. This proves
\eqref{da:eq:free-interpolation}.
The entry bound \eqref{uvc:eq:free-row-bounds} then proves
\eqref{da:eq:dyson-interpolation}.
For a comparison matrix write
$M_{B,u}(z)=(A_L-w_u)^{-1}$ with $w_u=z+ug^{-2}m_u(z)$.
Lemma~\ref{da:lem:positive-height} gives $|w_u-E|\le C/g^2$ and bounded
entries. The identical free recurrence therefore proves the last assertion.
\end{proof}

\subsection{Polynomial-range derivatives and propagation scales}
\label{app:large-flow}

\begin{lemma}[Fourier inputs at polynomial coupling]
\label{ll:lem:deterministic}
Fix $B_g,B_L,A>0$ and $\kappa\in(0,2)$. Assume $1\ll g<W^{B_g}$,
$C_Ag^2\log W\le L\le W^{B_L}$, $|E|\le2-\kappa$, and
$0\le t\le1-W^{-3}$, with $C_A$ sufficiently large depending on
these fixed parameters. Use the Dyson solution
\eqref{fo:eq:dyson-family} and the scales \eqref{fo:eq:scale}.
For the finite solution $m=m_{g,L}(E)$, put
$f(p)=|(2\cos p-E-m/g^2)^{-1}|$ and
$U_t(p)=(1-t\widehat P(p))^{-1}$, with normalized finite-grid averages.
Then
\begin{equation}\label{ll:eq:forest-input}
 \|f\|_\infty\le Cg^2,\qquad \|f\|_1\le C(1+\log(g^2)),
 \qquad \sup_p\langle f(q)U_t(p-q)\rangle
       \le C(1+\log(g^2)+W\zeta_t).
\end{equation}
For a sufficiently small fixed $c>0$,
\begin{equation}\label{ll:eq:low-frequency-input}
 \begin{gathered}
 (1-t+g^4|p|_{\mathbb T}^2)U_t(p)\le C
       \quad(|p|_{\mathbb T}\le c/g^2),\qquad
 U_t(p)\le C\quad(|p|_{\mathbb T}\ge c/g^2),\\
 \langle\mathbf1_{\{|p|_{\mathbb T}\le c/g^2\}}U_t(p)\rangle
       \le CW\zeta_t/g^2,\qquad W\zeta_t\ge c.
 \end{gathered}
\end{equation}
\end{lemma}

\begin{lemma}[Auxiliary estimates in polynomial coupling ranges]
\label{da:lem:polynomial-derivatives}
Fix $B_g,B_L>0$ and $\kappa\in(0,2)$. Assume $1\ll g<W^{B_g}$,
$\log W\ll L\le W^{B_L}$, and $|E|\le2-\kappa$.
Use $\gamma=g^{-2}$, the diffusivity \eqref{ba:eq:diffusivity}, and
$\eta_t=(1-t)g^{-2}\Im m_{g,L}(E;g^{-2})$.
For sufficiently large $g$,
\[
 m_g(E;g^{-2})=\frac{i}{\sqrt{4-E^2}}+O(g^{-2}),\qquad
 \Im m_g(E;g^{-2})\asymp1,\qquad D_g\asymp g^4.
\]
Fix $A>0$ and restrict to $0\le t<1$, $1-t\ge W^{-3}$.
Then $L\ge C_Ag^2\log W$ gives the same density and
entry bounds for the actual solution on the finite cycle of $L$ blocks. Its free row has
pointwise bound $C$ and weighted row sum $Cg^2$ at spatial scale $g^2$.
Its mixed kernel has exact row mass $(1-t)^{-1}$, pointwise bound
$O(((1-t)\ell_t)^{-1})$, and the exponential profile and increment
bound \eqref{gc:eq:gradient}, where
\begin{equation}\label{ll:eq:normalization}
 \begin{gathered}
 \ell_t=\min\{L,1+g^2/\sqrt{1-t}\},\qquad
 \zeta_t=\frac{g^2}{W\ell_t(1-t)}
 \asymp\frac1{W\sqrt{1-t}}+\frac{g^2}{WL(1-t)},\\
 (W\eta_t)^{-1}=(\Im m)^{-1} \ell_t\zeta_t,\qquad
 (W/g^2)\ell_t\zeta_t=(1-t)^{-1}.
 \end{gathered}
\end{equation}
The same-sign inverse has a bounded weighted row on scale $g^2$,
including at $E=0$. The finite-cycle comparisons for the Dyson function and pair kernels
have error $O(W^{-A})$.
For each fixed $j\ge0$, the derivatives $\partial_E^j m_g$ are bounded and
$|\partial_E^jD_g|\le C_jg^4$. For every fixed accuracy $A$, the
finite-cycle comparisons for the Dyson solution and pair kernels have
error $O(W^{-A})$ after any fixed number of $E$-derivatives and, for
the pair kernels, $t$-derivatives, when the constant in
$L\ge C_Ag^2\log W$ is chosen after these orders. The Dyson row and
same-sign inverse bounds also hold on the comparison segments of
Lemma~\ref{da:lem:positive-height}.
\end{lemma}

\begin{proof}[\textbf{Proof of Lemmas~\ref{ll:lem:deterministic} and \ref{da:lem:polynomial-derivatives}}]
Differentiate the line Dyson contraction established in the proof of
Lemma~\ref{uvp:prop:finite}, on the buffered bulk
$|E|\le2-\kappa/2$. Its scalar derivative is $O(g^{-2})$, giving
the expansion and its differentiated form:
\begin{equation}\label{ll:eq:dyson-derivatives}
 \partial_E^j m_g(E;g^{-2})
 =\partial_E^j\!\left(\frac{i}{\sqrt{4-E^2}}\right)+O_j(g^{-2}),
 \qquad |\partial_E^j D_g|\le C_jg^4.
\end{equation}
The Ward relation $g^{-2}|m_g|^2=\tanh q$ and
\eqref{ba:eq:diffusivity} give
\[
 D_g=\frac{g^4}{4|m_g(E;g^{-2})|^4}-\frac14.
\]
The lower bound on $|m_g|$ and its fixed derivative bounds therefore
give $D_g\asymp g^4$ and the derivative bound in
\eqref{ll:eq:dyson-derivatives}. The zero-order mixed profile and increments
follow from \eqref{ba:eq:linekernel} with
$\lambda\asymp\sqrt{1-t}/g^2$; the actual finite-kernel bounds follow
from Lemma~\ref{uvp:lem:actual-kernels}. The same-sign row is
\eqref{uvp:eq:stable}, since $L\ge C_Ag^2\log W$ ensures its threshold
$L/g^2\ge R_\kappa$.

The finite comparison uses the same periodic free row
\eqref{bl:eq:free}, now with $-\log|\rho|\asymp g^{-2}$.
Its error is $O(W^Ce^{-cL/g^2})$ at every prescribed derivative order.
The scalar derivative is $O(g^{-2})$, so the implicit comparison
gives the actual finite Dyson solution and its derivatives with that error.
As in the bounded-coupling proof, compare $P,Q$ with the
periodized line square kernels. The mixed inverse costs at most
$(1-t)^{-2}$; stable inverses use a convergent weighted perturbation
series. Choose $C_A$ after these polynomial costs. Row masses always
use the exact finite Ward identity \eqref{da:eq:ward}.
On the vertical comparison segment, Lemma~\ref{da:lem:positive-height}
gives the Dyson rows and stable inverse, with the same bulk buffer.
The scale relations \eqref{ll:eq:normalization} follow directly from
\eqref{fo:eq:scale} and $\eta_t=(1-t)g^{-2}\Im m$.
Since $\ell_t\le2g^2/\sqrt{1-t}$, they also give
$W\zeta_t=g^2/((1-t)\ell_t)\ge1/2$.

\noindent\textbf{Proof of \eqref{ll:eq:forest-input} and
\eqref{ll:eq:low-frequency-input}.}
Write $\delta=1-t$. The two free poles have slopes bounded
away from zero on the buffered bulk, and $\Im m/g^2\asymp g^{-2}$.
Thus $\|f\|_\infty\le Cg^2$. Summing on the finite grid in dyadic
distance shells about these poles gives
\[
 \langle f\rangle\le C(1+\log(g^2)+g^2/L)\le C(1+\log(g^2)),
\]
since $L\ge C_Ag^2\log W$. The line multiplier
\eqref{ba:eq:linekernel}, with $D_g\asymp g^4$, and the finite-cycle
comparison give
\[
 U_t(p)\le\frac{C}{\delta+g^4|p|_{\mathbb T}^2}
 \quad(|p|_{\mathbb T}\le c/g^2),\qquad
 U_t(p)\le C\quad(|p|_{\mathbb T}\ge c/g^2).
\]
Here the comparison error is chosen smaller than $c\delta^2$;
all additional polynomial costs are absorbed by increasing $C_A$.
The same symbol and normalized grid summation give, including the zero mode,
\[
 \langle|U_t-1|\rangle+
 \langle\mathbf1_{\{|p|_{\mathbb T}\le c/g^2\}}U_t(p)\rangle
 \le C\left(\frac1{g^2\sqrt\delta}+\frac1{L\delta}\right)
 \le\frac{C W\zeta_t}{g^2}.
\]
This proves \eqref{ll:eq:low-frequency-input}. Finally,
\[
 \sup_p\langle f(q)U_t(p-q)\rangle
 \le\langle f\rangle+\|f\|_\infty\langle|U_t-1|\rangle
 \le C(1+\log(g^2)+W\zeta_t),
\]
which proves \eqref{ll:eq:forest-input}.
\end{proof}

\begin{lemma}[Early-time pair inverses]\label{da:lem:early-pair-inverses}
Assume \eqref{da:eq:propagation-domain} and use the time-zero scale
$\ell_0$ from \eqref{fo:eq:scale}.
For some fixed $c>0$, uniformly for $0\le u\le1/2$,
\begin{equation}\label{da:eq:early-pair-inverses}
 \max_a\sum_b e^{c\|a-b\|_L/\ell_0}
 \left(|[(I-uP)^{-1}]_{ab}|+|[(I-uQ)^{-1}]_{ab}|\right)\le C.
\end{equation}
\end{lemma}
\begin{proof}[\textbf{Proof of \eqref{da:eq:early-pair-inverses}}]
The exact Ward identity gives $P\mathbf1=\mathbf1$, and $|Q|\le P$.
Hence the unweighted inverse rows are at most $(1-u)^{-1}\le2$.
The exponential free-row bounds in Lemma~\ref{bl:lem:finite} and
Lemma~\ref{uvp:prop:finite} give
$\sum_b(e^{c\|a-b\|_L/\ell_0}-1)P_{ab}\le Cc$ for small fixed $c$.
Choose $c$ so that the weighted row norm of $P$ is at most $3/2$.
The Neumann series for both inverses then converge in that norm,
with sum at most $(1-3u/2)^{-1}\le4$.
\end{proof}

\begin{corollary}[Differentiated finite-cycle remainders]\label{da:cor:remainder-derivatives}
Under the assumptions of Lemma~\ref{da:lem:scales-propagators},
\eqref{gc:eq:kernel-remainder} also holds after any fixed finite
orders of characteristic-time differentiation. On a vertical comparison
segment, use coordinates $(u,r)$ with $z=x+i\gamma_g r$, holding $g,x$
fixed, and put $P(u,r)_{ab}=\gamma_g|M_{B,u}(z)_{ab}|^2$.
The same remainder bound holds after any fixed mixed partial derivatives
in $u,r$, comparing this kernel with its periodized line analogue.
\end{corollary}
\begin{proof}[\textbf{Proof of Corollary~\ref{da:cor:remainder-derivatives}}]
On the characteristic, $P$ is independent of time, so
$\partial_u^j(I-uP)^{-1}=j!P^j(I-uP)^{-j-1}$; its row bound costs only
$j!(1-u)^{-j-1}\le j!W^{3(j+1)}$.
For the segment derivatives, abbreviate $\gamma=\gamma_g$,
$M=M_{B,u}(z)$, $m=L^{-1}\Tr M$, $S=L^{-1}\Tr M^2$ and
$d=1-u\gamma S$. Differentiating the scalar Dyson equation and then
$M=(A_g-z-u\gamma m)^{-1}$ gives
\[
 \partial_u m=\frac{\gamma mS}{d},\qquad
 \partial_r m=\frac{i\gamma S}{d},\qquad
 \partial_u M=\frac{\gamma m}{d}M^2,\qquad
 \partial_r M=\frac{i\gamma}{d}M^2.
\]
The stability and weighted rows in \eqref{da:eq:weighted-positive-stability}
bound $d^{-1}$, all fixed-order derivatives of $m$, and, by repeated
convolution, the entries and covariance-weighted rows of the derivatives
of $M$. Thus all fixed-order derivatives of $P(u,r)$ have bounded
weighted row norms. Put $R=(I-uP(u,r))^{-1}$. The exact row mass and
the lower characteristic endpoint of the segment give
\[
 1-u\sum_bP(u,r)_{ab}
 =\frac{r}{r+u\Im m}\ge cW^{-3},\qquad \|R\|_{\infty\to\infty}\le CW^3.
\]
Here the density is bounded above and below, and the segment has
$r\ge c(1-u)\ge cW^{-3}$. The inverse derivatives are
\[
 \partial_uR=R(P+u\partial_uP)R,\qquad
 \partial_rR=R(u\partial_rP)R.
\]
Every fixed mixed derivative is a finite sum of such products and
therefore costs only a fixed power of $W$. The same identities apply
to the line kernels. Subtracting the finite and periodized-line products
leaves one differentiated finite-cycle error in every term; all other
factors have the polynomial bounds just proved.
Along $x(u)=\Re z_u(g,E)$, one has
$x'(u)=-\gamma_g\Re m(g,E)$ and
$\gamma_g\partial_xM=(\partial_rM)/i$, so the additional chain-rule
terms for $M$, and hence for $P$ by products, have the same bounds.

The error-absorption argument proving \eqref{gc:eq:kernel-remainder}
applies with polynomial powers increased by the fixed derivative orders.
Differentiating \eqref{bl:eq:free} and the finite-cycle comparisons in
Lemmas~\ref{bl:lem:finite} and~\ref{uvp:prop:finite} gives
$C(1+L)^q e^{-cL}$ at bounded coupling and
$C(1+L/g^2)^q e^{-cL/g^2}$ in row norm and in $g^2$ times
pointwise norm at large coupling.
The bounded variance interval makes
$\partial_r=\gamma_g\partial_{\Im z}$ harmless.
Pointwise normalization costs at most $\ell_u\le L\le W^{B_L}$,
or $\ell_u/g^2\le2W^{3/2}$, respectively.
The preceding inverse-derivative bounds add only fixed powers of $W$,
absorbed by the same volume lower bounds at any fixed accuracy.
At $g=0$, the finite and line errors vanish identically.
\end{proof}

\paragraph{Fixed derivative orders.}
In polynomial ranges, for every fixed derivative order the Dyson
matrices and stable inverses have polynomial bounds in the spectral,
coupling and time variables. For the Dyson solutions this follows
by implicit differentiation in Lemma~\ref{bl:lem:finite} and
\eqref{ll:eq:dyson-derivatives}; for the inverses differentiate
$(I-uQ)^{-1}$ and use the already proved weighted row bounds.
The general-variance formulas apply separately on $g\le1$ and $g\ge1$.

At large coupling, these conclusions also hold at the comparison couplings $g_v\asymp g$
in \eqref{bl:eq:comparison-dyson}: $L/(g_v^2\log W)\to\infty$, so
\eqref{da:eq:kernel-accuracy} and
Corollary~\ref{da:cor:remainder-derivatives} absorb the fixed polynomial
costs. In particular, they include the polynomial-coupling specialization
of \eqref{gc:eq:kernel-remainder}.

\subsection{Propagation of tensors with zero rows}
\label{sec:gc-projected}

Use the parameter ranges and deterministic scales of
Theorem~\ref{gc:thm:common}(i). Fix $s\le u\le v\le t$ and $n\ge2$.
For an $n$-tensor $A$ supported on relative diameter $H\ell_u$, $H\ge1$,
write $\mathcal U_{u,v}^{\boldsymbol\sigma}$ for the tensor product of
the $n$ external-leg propagators in \eqref{fo:eq:propagators}. Then
\begin{equation}\label{fo:eq:projected-propagation}
 \begin{aligned}
 \|(U_{u,v}^{+-})^{\otimes n}A\|_\infty
 &\le C_nH^n(\zeta_v/\zeta_u)^n\|A\|_\infty,
 &&\text{zero rows},\\
 \|\mathcal U_{u,v}^{\boldsymbol\sigma}A\|_\infty
 &\le C_nH^{n-1}(\zeta_v/\zeta_u)^{n-1}\|A\|_\infty,
 &&\text{at least one stable leg}.
 \end{aligned}
\end{equation}
Here zero rows means $\sum_{a_2,\ldots,a_n}A(a_1,\ldots,a_n)=0$
for every $a_1$. A same-sign neighboring pair gives a stable leg by
\eqref{fo:eq:propagation-dyson}.

Use the probability profile $\theta_v=(1-v)P(I-vP)^{-1}$ and its
derivative \eqref{da:eq:profile-derivative}.

\begin{proof}[\textbf{Proof of \eqref{fo:eq:projected-propagation}}]
Use the one-dimensional proof of
\cite[Lemma~7.7(I)--(II), Appendix~A.2]{TYY}, with
\eqref{gc:eq:propagator-input} and
\eqref{da:eq:mixed-profile-increment} as the kernel inputs.
Write $U=I+\Xi$, $\Xi\coloneqq(v-u)P(I-vP)^{-1}$, and
$\rho\coloneqq(1-u)/(1-v)$.
An identity or stable leg anchors one input and gives the factor
$C_n(H\zeta_v/\zeta_u)^{n-1}$.
For the all-$\Xi$ term, zero-row subtraction at $b_2=\cdots=b_n=b_1$
leaves one row sum, $n-1$ pointwise kernel bounds, and the increment
gain. Summing over the relative support gives
\[
 C_n\rho\left(\frac{\rho}{\ell_v}\right)^{n-1}
 (H\ell_u)^{n-1}\frac{H\ell_u}{\ell_v}\|A\|_\infty
 =C_n\left(\frac{\rho H\ell_u}{\ell_v}\right)^n\|A\|_\infty.
\]
These prove both bounds, including $\ell_u=1$ and $\ell_v=L$.
For tails in the zero-row case, set
$A_H=A\mathbf1_{\{\diam\le H\ell_u\}}$, $H=W^\nu$.
Its row is the negative of the discarded row, so restoration
against $\theta_v$ gives
\[
 \left\|\left(\sum_{b_2,\ldots,b_n}A_H(b_1,\ldots,b_n)\right)
       \prod_{j=2}^n\theta_v(b_j-b_1)\right\|_\infty
 \le C_n(L/\ell_v)^{n-1}\|A-A_H\|_\infty.
\]
The all-$\Xi$ subtraction for $A_H$ leaves its nonzero row; using
$n-1$ pointwise $\Xi$ bounds and one row sum bounds this term by
$((1-u)/(1-v))^n$ times the displayed right-hand side.
The discarded tensor and the restoration tensor satisfy the same
bound by the propagator row sums. Estimate these tensors separately:
the restoration has scale $\ell_v$, not support in $H\ell_u$.
Since $L/\ell_v\le L/\ell_u$ and the time ratio are polynomially
bounded, a tail bound $W^{-A}$ with $A$ exceeding these fixed polynomial
costs by the desired accuracy makes these errors negligible.
With a stable leg, the discarded tail is controlled directly by
the propagator row sums, without restoring a row.
\end{proof}

\subsection{A coarse stability bound in the upper half-plane}
\label{app:rough-stability}

\begin{lemma}[Coarse scalar stability]\label{da:lem:rough-stability}
Let $E\in\R$ and use the free transforms and
Dyson solutions in \eqref{bl:eq:physical}--\eqref{da:eq:dyson}.
Take the normalized parameters $\gamma=\max\{1,g\}^{-2}$ and
$z=E+i\eta$, $0<\eta\le1$. Suppose $u\in\C_+$, $|u|\le\eta^{-1}$,
and
\[
 u=f_g(z+\gamma u;\gamma)+r,\qquad |r|\le\eta/2.
\]
Then $|u-m_g(z)|\le C\eta^{-2}|r|$, uniformly in $g$ and $E$.
For every $w\in\C_+$ with $\Im w\ge\eta$, the same parameters satisfy
\[
 |f_{g,L}(w;\gamma)-f_g(w;\gamma)|\le CL^{-1}\eta^{-2}.
\]
\end{lemma}

\begin{proof}[\textbf{Proof of Lemma~\ref{da:lem:rough-stability}}]
Subtract the exact Dyson equation to obtain
$u-m_g=q(u-m_g)+r$, where
\[
 q\coloneqq \gamma\left\langle
 \frac1{(2\sqrt\gamma\,g\cos p-z-\gamma u)
          (2\sqrt\gamma\,g\cos p-z-\gamma m_g)}\right\rangle.
\]
Imaginary parts and Cauchy--Schwarz give
\[
 |q|^2\le
 \frac{\gamma(\Im u+|r|)}{\eta+\gamma\Im u}
 \frac{\gamma\Im m_g}{\eta+\gamma\Im m_g}.
\]
Since $\gamma\le1$, $|r|\le\eta/2$, and $\Im u\le\eta^{-1}$,
\[
 1-\frac{\gamma(\Im u+|r|)}{\eta+\gamma\Im u}
 =\frac{\eta-\gamma|r|}{\eta+\gamma\Im u}\ge\eta^2/4.
\]
The second fraction is also at most $1-c\eta^2$, since
$\Im m_g\le\eta^{-1}$. Thus $|1-q|\ge c\eta^2$, proving the stability estimate. For the Fourier-sum estimate,
the derivative of $(2\sqrt\gamma\,g\cos p-w)^{-1}$ has modulus
at most $2\eta^{-2}$, since $\sqrt\gamma\,g\le1$.
The elementary Riemann-sum bound proves the claim.
\end{proof}

\subsection{Terminal variance matching}\label{app:terminal-matching}

\begin{proof}[\textbf{Proof of Lemma~\ref{fo:lem:bounded-terminal-matching}}]
Set $z\coloneqq E+i\eta$. The definitions \eqref{bl:eq:terminal} give
\[
 \begin{gathered}
 A_{g_*}=c_*A_g,\qquad t_*\gamma_{g_*}=c_*^2\gamma_g,\\
 E_*+\gamma_{g_*}m/c_*=c_*(z+\gamma_gm).
 \end{gathered}
\]
The Dyson uniqueness and Gaussian scaling calculation proving
\cite[Lemma~3.3, (3.11)--(3.13)]{TYY}, with these covariance identities,
gives $m_{g_*,L}(E_*)=m/c_*$, $z_{t_*}(g_*,E_*)=c_*z$ and
$H_{t_*}(g_*)\stackrel d=c_*H$, hence
\eqref{fo:eq:bounded-terminal-resolvent}.
The bulk bounds of Lemmas~\ref{bl:lem:regular} and
\ref{bl:lem:finite} give $\Im m\asymp1$ and $|m|\le C$.
Thus $1-t_*\asymp\eta$, while
$\sqrt{t_*}\le c_*\le1$ follows directly from
$\gamma_g=\max\{1,g\}^{-2}$. Together with
$m_{g_*,L}(E_*)=m/c_*$ and \eqref{bl:eq:terminal}, these facts give
the parameter bounds \eqref{fo:eq:bounded-terminal-parameters}.
\end{proof}

\begin{lemma}[Terminal matching at large coupling]\label{fo:lem:matching}
Fix $\kappa\in(0,2)$ and assume $g\ge C_\kappa$ and
$L\ge C_*g^2\log W$, with $C_\kappa,C_*$ sufficiently large. Let
\[
 \tfrac12\le t_*<1,\qquad
 g_*\coloneqq\sqrt{t_*}g,\qquad
 \gamma_{g_*}=(t_*g^2)^{-1}.
\]
The flow hopping is $A_{g_*}=\sqrt{\gamma_{g_*}}g_*A_L=A_L$, and
\begin{equation}\label{fo:eq:large-matching}
 H_{t_*}(g_*)\stackrel d=H,\qquad
 t_*=(g_*/g)^2,\qquad t_*\gamma_{g_*}=g^{-2}.
\end{equation}
More explicitly, its terminal characteristic is
\[
 z_{t_*}(g_*,E_0)
       =E_0+(1-t_*)\gamma_{g_*}m_{g_*,L}(E_0).
\]
For sufficiently large $W$, its real and imaginary parts satisfy
\begin{equation}\label{fo:eq:large-match-domain}
 \begin{gathered}
 \forall E\in[-2+\kappa,2-\kappa]\quad
 \exists E_0\in[-2+3\kappa/4,2-3\kappa/4]:\\
 \Re z_{t_*}(g_*,E_0)=E,\\
 cg^{-2}\bigl(1-t_*\bigr)\le
       \Im z_{t_*}(g_*,E_0)
       \le Cg^{-2}\bigl(1-t_*\bigr),\\
 G\bigl(E+i\Im z_{t_*}(g_*,E_0)\bigr)
       \stackrel d=G_{t_*}(g_*,E_0).
 \end{gathered}
\end{equation}
Thus $(g_*,E_0)$ satisfies \eqref{gc:eq:initial-parameters}(b).
\end{lemma}

\begin{proof}[Proof of Lemma~\ref{fo:lem:matching}]
\noindent\textbf{Proof of \eqref{fo:eq:large-matching} and
\eqref{fo:eq:large-match-domain}.}
Here $A_{g_*}=A_g=A_L$ and $t_*\gamma_{g_*}=g^{-2}$, proving
\eqref{fo:eq:large-matching}. Since $g_*\ge g/\sqrt2$ and
$L/g_*^2\ge L/g^2$, Lemma~\ref{fo:lem:terminal-deterministic}(b)
applies at coupling $g_*$ and $\delta=1-t_*$ after enlarging the
fixed thresholds. It supplies $E_0$ in \eqref{fo:eq:large-match-domain}, with
\[
 \Im z_{t_*}(g_*,E_0)
 =(1-t_*)\gamma_{g_*}\Im m_{g_*,L}(E_0)
 \asymp\frac{1-t_*}{t_*g^2}\asymp(1-t_*)g^{-2}.
\]
The real characteristic map
$f(e)=e+(1-t_*)\gamma_{g_*}\Re m_{g_*,L}(e)$ has
$f'\in[1/2,3/2]$ by \eqref{fo:eq:terminal-anchor-derivative}, so its
inverse is differentiable. The resolvent identity in distribution
follows from the whole-matrix equality \eqref{fo:eq:large-matching}.
\end{proof}

\begin{lemma}[Terminal scale and energy curve in polynomial volumes]
\label{da:lem:polynomial-terminal}
Fix $B_g,B_L>0$, $\kappa\in(0,2)$ and $0<\nu<1$, and assume
$1\ll g<W^{B_g}$ and $g^2\log W\ll L\le W^{B_L}$.
Write $N=WL$ and $\ellstar=\min\{N,W^2g^2\}$, and define
\begin{equation}\label{ll:eq:endpoint}
 t_*\coloneqq1-\ellstar^{-1}g^2W^\nu,\qquad
 g_*\coloneqq\sqrt{t_*}g,\qquad \gamma_{g_*}=(t_*g^2)^{-1}.
\end{equation}
For sufficiently large $W$, $1/2\le t_*\le1-W^{-3}$ and
$(1-t_*)/(g^2/N)\ge W^\nu$. For the flow coupling $g_*$,
\begin{equation}\label{fo:eq:large-endbudget}
 \zeta_{t_*}\le C\left(
 \sqrt{W^{-2}g^{-2}\ellstar}\,W^{-\nu/2}
       +N^{-1}\ellstar W^{-\nu}\right)
 \le C(W^{-\nu/2}+W^{-\nu})\le W^{-c_0}
\end{equation}
for any fixed $0<c_0<\nu/2$. The inverse characteristic from
Lemma~\ref{fo:lem:matching} is differentiable, with
\begin{equation}\label{da:eq:polynomial-terminal-heights}
 \begin{gathered}
 f(e)=e+(1-t_*)\gamma_{g_*}\Re m_{g_*,L}(e),\quad
 E_0(E)=f^{-1}(E),\quad
 \eta(E)=(1-t_*)\gamma_{g_*}\Im m_{g_*,L}(E_0(E)),\\
 \tfrac12\le f'\le\tfrac32,\quad
 c\ellstar^{-1}W^\nu\le\eta(E)\le C\ellstar^{-1}W^\nu,\quad
 \eta(E)^{-1}\le W^{C_1},\quad |\eta'(E)|\le W^{C_2},\quad |E|\le2-\kappa.
 \end{gathered}
\end{equation}
Here $C_1,C_2$ are fixed; all constants are independent of $g,L$.
The prescribed height $\eta_L=W^\nu/\ellstar$ also satisfies
$\eta_L^{-1}\le W^{2+2B_g}$.
\end{lemma}
\begin{proof}[\textbf{Proof of Lemma~\ref{da:lem:polynomial-terminal}}]
The endpoint calculation is
\[
 \begin{gathered}
 1-t_*=\max\{W^{-1}L^{-1}g^2,W^{-2}\}W^\nu
 \le C\bigl((\log W)^{-1}W^{-1+\nu}+W^{-2+\nu}\bigr)=o(1),\\
 1-t_*\ge W^{-2+\nu}\ge W^{-3},\qquad
 \frac{1-t_*}{g^2/N}=\ellstar^{-1}NW^\nu\ge W^\nu.
 \end{gathered}
\]
Since $g_*^2\asymp g^2$, substitution into
\eqref{ll:eq:normalization} gives the first bound in
\eqref{fo:eq:large-endbudget}. The inequalities
$\ellstar\le W^2g^2$ and $\ellstar\le N$ give the second.
Lemma~\ref{fo:lem:matching} supplies the inverse map and the height
comparison in \eqref{da:eq:polynomial-terminal-heights}.
The derivative bounds in Lemma~\ref{da:lem:polynomial-derivatives}
and $|E_0'|\le2$ give
$|\eta'|\le2(1-t_*)\gamma_{g_*}
 |\partial_e m_{g_*,L}(E_0)|\le W^{C_2}$.
Finally $\eta^{-1}\le C\ellstar W^{-\nu}
\le CW^{2+2B_g-\nu}\le W^{C_1}$ for a fixed $C_1$;
the same bound without the constant gives the assertion for $\eta_L$.
\end{proof}

\subsection{Parameter control in polynomial volumes}\label{app:polynomial-iteration}
\begin{lemma}[Scale and parameter control in polynomial volumes]
\label{da:lem:polynomial-parameters}
Fix $B_g,B_L,\kappa>0$. For
$1\ll g<W^{B_g}$, $g^2\log W\ll L\le W^{B_L}$ and
$|E|\le2-\kappa$, the characteristic scales satisfy, for $0\le u<1$,
\begin{equation}\label{da:eq:polynomial-scale-budgets}
 \begin{gathered}
 \gamma_gL+L/\ell_u\le2W^{B_L},\qquad
 \gamma_gW\ell_u\zeta_u=(1-u)^{-1},\qquad \zeta_u\ge c/W,\\
 \ell_u\asymp g^2,\qquad\zeta_u\asymp W^{-1}
       \quad(0\le u\le1/2).
 \end{gathered}
\end{equation}
Starting with $|E_0|\le2-3\kappa/4$, every finite comparison chain
as in Lemma~\ref{da:lem:comparison-parameters}(b) satisfies
\eqref{fo:eq:large-input-domain}, preserves these polynomial ranges
with fixed buffered bounds, and cannot increase the scale $\zeta$.
\end{lemma}
\begin{proof}[\textbf{Proof of Lemma~\ref{da:lem:polynomial-parameters}}]
Here $\gamma_g=g^{-2}$ and
$\ell_u=\min\{L,1+g^2/\sqrt{1-u}\}\ge g^2$ for large $W$.
Thus $\gamma_gL+L/\ell_u\le2L/g^2\le2W^{B_L}$.
The normalization is the definition of $\zeta_u$; its lower bound
is \eqref{gc:eq:scale}. For $u\le1/2$, $L\gg g^2\log W$ leaves
the scale uncapped, giving the last line of
\eqref{da:eq:polynomial-scale-budgets}.
The volume-independent calculation \eqref{da:eq:comparison-chain-algebra}
and its buffered-bulk argument give \eqref{fo:eq:large-input-domain};
here $g_0\to\infty$ makes the displacement smaller than the fixed buffer.
The bounds $g_0/\sqrt2\le g\le g_0$ preserve
$g<W^{B_g}$, $L\le W^{B_L}$ and
$L/(g^2\log W)\ge L/(g_0^2\log W)\to\infty$.
Equation~\eqref{bl:eq:backdomain} preserves scale smallness.
\end{proof}

\subsection{The deformed one-block Dyson profile}
\label{app:deformed-dyson}

\begin{lemma}[A small-norm and small-rank deformation]\label{da:lem:deformed}
Fix $\kappa_0\in(0,2)$ and $C_*>0$. Let
$B=B_{\rm small}+B_{\rm bad}$ be deterministic Hermitian, with
$\|B_{\rm small}\|\le\varepsilon$ and $\rank B_{\rm bad}\le\varepsilon W$.
For sufficiently small fixed $\varepsilon$, the Herglotz solution
$M_B(z)=(-EI_W-B-z-m_B(z))^{-1}$, $m_B=W^{-1}\Tr M_B$, satisfies
\begin{equation}\label{us:eq:deformed-profile}
 \Im m_B(i\eta)\asymp1,\qquad
 \Im M_B(i\eta)\asymp(I_W+B^2)^{-1},\qquad
 |E|\le2-\kappa_0,\quad 0<\eta\le C_*.
\end{equation}
The constants depend only on the fixed bulk margin and height bound,
not on the exceptional eigenvalues of $B$. 
The comparison for $\Im M_B$ is an inequality of quadratic forms.
\end{lemma}

Under the hypotheses of Lemma~\ref{da:lem:deformed}, we also prove,
uniformly in its energy and height ranges,
\begin{equation}\label{da:eq:deformed-scalar-stability}
 |m_B(i\eta)|\le1,\qquad
 \left|1-W^{-1}\Tr M_B(i\eta)^2\right|\ge c.
\end{equation}

\begin{proof}[\textbf{Proof of Lemma~\ref{da:lem:deformed} and \eqref{da:eq:deformed-scalar-stability}}]
Write $M=M_B$, $m=m_B$, and $H=\max\{C_*,2\}$.
Imaginary parts and Cauchy--Schwarz give $|m|\le1$ as in
\eqref{da:eq:ward}. Min--max gives
$\#\{\lambda(B)\notin[-\varepsilon,\varepsilon]\}
\le\operatorname{rank}B_{\rm bad}\le\varepsilon W$.
On a fixed tube around the bulk semicircle branch, the imaginary gap
bounds the scalar inverses and their squares. The small eigenvalues
change each by $O(\varepsilon)$, and the exceptional fraction contributes
$O(\varepsilon)$ regardless of its eigenvalues. Thus
\[
 m=-(E+i\eta+m)^{-1}+O(\varepsilon).
\]
Start at $\eta=H$, where the scalar map is a contraction. The two
semicircle roots remain uniformly separated for
$|E|\le2-\kappa_0$, $0\le\eta\le H$; the displayed residual excludes
a first exit from the tube for small $\varepsilon$. Hence
$m=m_{\rm sc}(E+i\eta)+O(\varepsilon)$ and $\Im m\asymp1$.
Applying the same eigenvalue comparison to the squared inverses yields
$W^{-1}\Tr M^2=m_{\rm sc}(E+i\eta)^2+O(\varepsilon)$;
the bulk semicircle stability gap proves
\eqref{da:eq:deformed-scalar-stability}.
Finally, in each eigenvalue coordinate $b$ of $B$,
\[
 \Im M(i\eta)=\frac{\eta+\Im m}
 {(E+b+\Re m)^2+(\eta+\Im m)^2}\asymp(1+b^2)^{-1}.
\]
This proves the quadratic-form comparison uniformly in all exceptional
eigenvalues and as $\varepsilon\downarrow0$.
\end{proof}

\begin{proof}[\textbf{Proof of Lemma~\ref{da:lem:regular-dyson-profile}: \eqref{us:eq:mde}}]
Take $B=g^2R$:
its part on $\{|g^2R|\le g\}$ has norm at most $g$, and the
complementary rank is at most $Wg^{1/2}$. Both fractions tend to zero,
so \eqref{us:eq:deformed-profile} applies with
$\varepsilon=\max\{g,g^{1/2}\}$.
\end{proof}

\begin{lemma}[Residual stability and vector weights for the deformed profile]
\label{da:lem:deformed-stability}
Under the hypotheses of Lemma~\ref{da:lem:deformed}, write
$M=M_B(i\eta)$ and $m=m_B(i\eta)$. There are fixed $r_0,C>0$
such that, if $|q-m|\le r_0$, then
\begin{equation}\label{da:eq:deformed-residual-stability}
 r=q-W^{-1}\Tr(-EI_W-B-i\eta-q)^{-1}
 \quad\Longrightarrow\quad |q-m|\le C|r|.
\end{equation}
Uniformly in the same bulk and height ranges, for every vector $v$,
\begin{equation}\label{da:eq:deformed-vector-weight}
 \|M\|+\|\Im M\|\le C,\qquad \eta\le C\Im m,\qquad
 v^*\Im Mv=(\eta+\Im m)\|Mv\|_2^2\asymp\|Mv\|_2^2.
\end{equation}
\end{lemma}
\begin{proof}[\textbf{Proof of Lemma~\ref{da:lem:deformed-stability}}]
The positive imaginary gap in \eqref{us:eq:deformed-profile} bounds
the second derivative of the scalar Dyson map in a fixed neighborhood
of $m$. Taylor expansion gives
\[
 (1-W^{-1}\Tr M^2)(q-m)=r+O(|q-m|^2).
\]
The stability gap in \eqref{da:eq:deformed-scalar-stability} absorbs the
quadratic term when $r_0$ is sufficiently small, proving
\eqref{da:eq:deformed-residual-stability}. For a continuous approximate
solution starting close to $m$ at the upper height, this estimate
excludes a first exit from the neighborhood while the residual remains
small, uniformly throughout the height interval.
The same imaginary gap bounds $\|M\|$; $\Im m\asymp1$ gives
$\eta\le C\Im m$. Diagonalizing the Hermitian matrix $B$ gives
$\Im M=(\eta+\Im m)M^*M$, proving
\eqref{da:eq:deformed-vector-weight}.
\end{proof}

\section{Loop identities and primitive-loop estimates}\label{app:loop-identities}
We prove the loop identities and primitive estimates stated in
Section~\ref{sec:fo-loops}. The deterministic kernel inputs have already
been proved in Appendix~\ref{app:dyson-foundations}; the arguments below
are independent of stochastic propagation.

\subsection{Exact loop algebra and mixed Ward identity}\label{sec:gc-objects}
We use the flow, Dyson matrix, loops, and cut operations of
Section~\ref{sec:flow-overview}. Set $\gamma=\gamma_g$,
$\Lambda=A_g\otimes I_W$, and
$\mathcal S[R]=\gamma W^{-1}\diag(\Tr R_{aa})_{a\in\Z_L}\otimes I_W$.
The characteristic identity \eqref{da:eq:characteristic-general} gives
$\Lambda-z_t-t\mathcal S[M]=M^{-1}$ with $M$ independent of $t$.
Thus \cite[Assumption~4.1]{TYY} holds; its general tree representation
\cite[Lemma~4.16, (4.15)--(4.16)]{TYY} and primitive Ward identity
\cite[Corollary~4.21]{TYY} apply. These identities require no analytic
small-coupling hypothesis.

\begin{proof}[Proof of Lemma~\ref{fo:lem:mixed-ward}]
For actual loops, cyclically place the two opposite signs next to the
summed projector and use $\sum_cE_c=I/W$ and
$G_t-G_t^*=2i\eta_tG_tG_t^*$. This gives \eqref{gc:eq:ward}.
For primitives, the symmetric coordinate covariance is
$S_{ai,bj}=\gamma W^{-1}\mathbf1_{a=b}$.
Average each external coordinate of the general primitive in
\cite[Corollary~4.21]{TYY} against $P_a/W$, as in
\cite[Lemma~4.26, (4.38)]{TYY}.
The two new normalized projectors change the coordinate cut coefficient
$\gamma/W$ to $(\gamma/W)W^2=\gamma W$, so the averaged initial
values and evolution are exactly \eqref{fo:eq:primitive}.
The summed Ward vertex introduces one additional $W^{-1}$, giving
$1/(2iW\eta_t)$ in \eqref{gc:eq:ward}.
\end{proof}

\noindent\textbf{Proof of \eqref{fo:eq:hierarchy} and \eqref{gc:eq:bracket}.}
We start from the resolvent equation \eqref{fo:eq:ito}. The quadratic
contraction of two differentiated factors of a loop is its cut product.
A light insertion is $\gamma W$ times
$\langle\wt G^{\sigma}\rangle_{cc}$ and the $(k+1)$-loop obtained by
inserting the duplicated sign and $E_c$.
Let $\{T_\lambda\}_{\lambda=1}^{W^2}$ be a basis of Hermitian
$W\times W$ matrices with $\Tr(T_\lambda T_\mu)=\delta_{\lambda\mu}$,
and let $T_{c,\lambda}$ be its embedding into block $c$.
Then $d\mathcal B_t^{(g)}=\sqrt{\gamma/W}\sum_{c,\lambda}
T_{c,\lambda}\,d\beta_{c,\lambda}$ for independent standard real
Brownian motions $\beta_{c,\lambda}$.
Differentiating the $i$-th resolvent in \eqref{fo:eq:loops} and
cyclically moving $d\mathcal B_t^{(g)}$ to the front gives
$-\Tr(d\mathcal B_t^{(g)}B_i)$, with $B_i$ defined before
\eqref{gc:eq:bracket}.
Completeness of this basis gives
$\sum_\lambda\Tr(T_{c,\lambda}B_i)
\overline{\Tr(T_{c,\lambda}B_j')}
=\Tr(P_cB_iP_cB_j'^*)$, which proves \eqref{gc:eq:bracket}.
Replacing both $P_c$ by $WE_c$ makes this $\gamma W$ times
joined $(2k+2)$-loops with one free block, rather than a product
of unrelated loops.

\begin{proof}[\textbf{Proof of Lemma~\ref{uvc:lem:orbital-limit}: \eqref{uvc:eq:orbital-limit}}]
On $[0,T]$, the auxiliary model has
$\Im z_t\ge\eta_{\min}:=(1-T)\gamma_g\Im m>0$ and hence
$\|G_t^{(\wt W)}\|\le\eta_{\min}^{-1}$, uniformly in $\wt W$.
At fixed $t$, the directional derivative of a normalized trace word
in a Hermitian direction $V_c$ in block $c$ is a sum of $k$ terms
$\wt W^{-1}\Tr(YV_c)$. Each $Y$ is formed from resolvents and
unnormalized block projections, compressed to block $c$, so
$\|Y\|_{\mathrm{HS}}\le C\sqrt{\wt W}$.
Gaussian Poincare with the fixed-time covariance $t\gamma_g/\wt W$ gives
\[
 \sup_{0\le t\le T}
 \operatorname{Var}(\wt W^{k-1}\cL_k^{(\wt W)}(t))
 \le C_{k,g,L,E,T}\wt W^{-2}.
\]

Common orbital-unitary invariance defines $R_{\wt W}(t)\in\C^{L\times L}$ by
$\E(G_t^{(\wt W)})_{bc}=R_{\wt W}(t)_{bc}I_{\wt W}$.
By translation invariance,
$\bar m_{\wt W}(t):=R_{\wt W}(t)_{bb}=\E\cL_1^{(\wt W)}(t;b,+)$
is independent of $b$. Gaussian integration by parts gives
\[
 \bigl(A_g-(z_t+t\gamma_g\bar m_{\wt W}(t))I_L\bigr)R_{\wt W}(t)=I_L+o(1),
\]
uniformly on $[0,T]$ in $L\times L$ matrix norm. Each entry of the
residual matrix is a covariance term between a normalized diagonal trace
and a resolvent entry, hence is $o(1)$ by Poincare and the resolvent bound.
Since $\Im\bar m_{\wt W}(t)\ge0$, the inverse of
$A_g-(z_t+t\gamma_g\bar m_{\wt W}(t))I_L$ has norm at most $\eta_{\min}^{-1}$.
The uniform residual and boundedness imply that any subsequential limit
$\mu$ of $\bar m_{\wt W}(t)$ as $\wt W\to\infty$ and $t\to t_*\in[0,T]$ satisfies
\[
 \mu=L^{-1}\Tr\bigl(A_g-(z_{t_*}+t_*\gamma_g\mu)I_L\bigr)^{-1}.
\]
Scalar uniqueness in Lemma~\ref{da:lem:positive-height} and the characteristic
identity \eqref{gc:eq:positive-dyson} give $\mu=m$. Compactness therefore yields
$\sup_{t\le T}|\bar m_{\wt W}(t)-m|\to0$.
Together with Poincare, this gives
\[
 \sup_{0\le t\le T}\E|\cL_1^{(\wt W)}(t;b,+)-m|^2\longrightarrow0.
\]

Take expectations in \eqref{fo:eq:hierarchy} at width $\wt W$ and normalize.
The martingale has zero expectation; the centered insertion is $o(1)$
by the preceding $L^2$ convergence and the resolvent bound.
Poincare factorizes the normalized cut products up to $o(1)$, giving
\[
 \frac{d}{dt}\bigl(\wt W^{k-1}\E\cL_k^{(\wt W)}\bigr)
 =\gamma_g\wt W^k\sum_{i<j,c}
       \bigl(\E\operatorname{Cut}_{L,i,j}^c\cL_k^{(\wt W)}\bigr)
       \bigl(\E\operatorname{Cut}_{R,i,j}^c\cL_k^{(\wt W)}\bigr)+o(1),
\]
uniformly on $[0,T]$. For cut lengths $p+q=k+2$,
$\wt W^k\wt W^{1-p}\wt W^{1-q}=1$, so the normalized equation agrees
with \eqref{fo:eq:primitive} for $W^{k-1}\cK_k$, up to $o(1)$;
the initial data agree exactly.
For $k=2$, the resolvent bound and \eqref{fo:eq:pairs} bound both
$\wt W\E\cL_2^{(\wt W)}$ and the Riccati solution $W\cK_2$ on $[0,T]$.
Subtracting their integral equations and applying Gronwall gives uniform convergence.
For $k\ge3$, adjacent cuts are linear in the length-$k$ tensor and
all other cuts use shorter words; induction and Gronwall give
uniform convergence to $W^{k-1}\cK_k$.
There are finitely many label and sign tuples, proving
\eqref{uvc:eq:orbital-limit}.
\end{proof}

\subsection{Tree representation and stable molecule definitions}\label{sec:bl-primitives}
For a translation-invariant $k$-tensor, Fourier transform the
$k-1$ displacements from its first vertex; the momenta
$k_1,\ldots,k_k$ satisfy $\sum_i k_i=0$. Inverse Fourier transform
uses normalized averaging on each copy of $(2\pi/L)\mathbb Z_L$.
For a sign pair define
\[
 V_t^{+-}(q)=V_t^{-+}(q)\coloneqq (1-t\widehat P(q))^{-1},\qquad
 V_t^{++}(q)\coloneqq (1-t\widehat Q(q))^{-1},\qquad
 V_t^{--}(q)\coloneqq \overline{V_t^{++}(q)}.
\]
The mixed inverses may grow as $t\uparrow1$; the same-sign inverses
have the bounded weighted rows of Lemmas~\ref{ba:lem:kernels}
and~\ref{uvp:prop:finite}. For $k\ge3$, define
the Fourier core by
\[
 W^{k-1}\widehat{\cK}_k(t;\boldsymbol k;\boldsymbol\sigma)
 =C_t^{[k]}(\boldsymbol k;\boldsymbol\sigma)
       \prod_{i=1}^k V_t^{\sigma_i\sigma_{i+1}}(k_i),
 \qquad \sigma_{k+1}\coloneqq \sigma_1.
\]
The tree/core decomposition \cite[Lemmas~4.16, 4.22, 4.25]{TYY},
with the block reduction \cite[(4.38)--(4.42)]{TYY}, gives the
following polygon-dual formula over noncrossing dissections:
\begin{equation}\label{ba:eq:dissection}
 C_t^{[k]}=\sum_{\mathfrak d}
       \prod_{F\in\mathfrak d}C_0^{[|F|]}(F)
       \prod_{d\in\mathfrak d}\gamma_g tV_t(q_d).
\end{equation}
Faces have at least three sides, and each unordered nonadjacent diagonal
$d=\{i,j\}$ occurs once. It carries
$V_t(q_d)=V_t^{\sigma_i\sigma_j}(q_{ij})$ with
$q_{ij}=k_i+\cdots+k_{j-1}$, opposite momenta in the two children,
and the cut signs of \eqref{fo:eq:cuts}.
The block normalization in both coupling regimes is
$\gamma_g W W^{1-k_L}W^{1-k_R}/W^{1-k}=\gamma_g$,
where $k_L+k_R=k+2$.

Define a stable molecule $S_t^{[k]}$ as the spatial tensor whose Fourier
transform is the sum in \eqref{ba:eq:dissection} restricted to dissections
whose internal diagonals all carry same-sign inverses. Its external inverses are
still removed. Any unrestricted dissection is obtained by joining
these molecules along its mixed diagonals; the adjacency graph of
the molecules is a tree. For an anchored tensor, its weighted row
norm is the sum over the remaining vertices after multiplying by
$e^{c\operatorname{diam}(0,\boldsymbol x)}$, with periodic diameter.

\subsection{Core cancellation and bounded-coupling estimates}\label{app:bounded-primitives}
The next two estimates are uniform under the parameter conditions
\eqref{bl:eq:intrinsic}.

\begin{lemma}[Spatial moments and the alternating core]\label{ba:lem:moment}
Use the bounded-coupling parameters $C,B_0,c',g,E$ of
Lemma~\ref{ba:lem:kernels}, and let $S_t^{[k]}$ be the stable molecule defined after \eqref{ba:eq:dissection}.
Its norm $\|\cdot\|_{\ell^1_c}$ is the weighted row norm defined there,
with a sufficiently small fixed $c>0$.
For every fixed $k\ge3$, every external anchor, and every fixed
$a\ge1$, uniformly for $g\le C$, $0\le t\le1$,
\begin{equation}\label{ba:eq:moment}
 \|S_t^{[k]}\|_{\ell^1_c}\le C',\qquad
 \sum_{\boldsymbol x}\diam(0,\boldsymbol x)^a
       |S_t^{[k]}(0,\boldsymbol x)|\le C'g^2.
\end{equation}
For an alternating even molecule,
\begin{equation}\label{ba:eq:cancel}
 |\widehat S_t^{[k]}(\boldsymbol k)|
       \le C'\left((1-t)+g^2\sum_i|k_i|_{\mathbb T}^2\right).
\end{equation}
\end{lemma}
\begin{proof}
A term with $s$ stable diagonals has $s+1$ faces and $k+2s$
deterministic $M$-edges. Delete one edge in each face. The remaining
edges and covariance links form a spanning tree. Deleted edges are
pointwise bounded; all retained edges and stable inverses have bounded
weighted rows by Lemma~\ref{ba:lem:kernels}. Leaf summation, retaining
a small external exponential weight, proves the first estimate.

For small $g$, the free Neumann expansion and
\eqref{ba:eq:stableinverse} give weighted off-diagonal rows $O(g)$
for $M_B$ and $O(g^2)$ for each same-sign inverse, uniformly in $L,t$.
Expand each stable dissection into diagonal and off-diagonal edges.
A term with positive external diameter contains either an off-diagonal
stable edge, or at least two off-diagonal $M_B$-edges: if all stable
edges identify their endpoints, some closed face must leave a vertex
and return. In the preceding tree sum, each retained marked edge
keeps its weighted row gain, and each deleted marked edge keeps its
pointwise gain. Thus the weighted row restricted to positive diameter
is $O(g^2)$. Decreasing the weight proves every fixed positive spatial
moment in \eqref{ba:eq:moment}; for $g$ bounded away from zero, the
first bound gives the same conclusion after enlarging $C'$.

For an alternating molecule, use the signed Ward-deletion and subtraction
induction in \cite[Corollary~4.28 and proof of Lemma~4.29, (4.45)--(4.46)]{TYY}.
Here the stable weighted tree bounds pure-sign rows by $CW^{1-k}$,
each Ward deletion costs $C/[W(1-t)]$, and each normalized zero-mode
mixed join costs $C/(1-t)$, since $\gamma_g,\Im m\asymp1$.
These inputs give $|\widehat S_t^{[k]}(0)|\le C'(1-t)$.
The rows are signed; absolute summation before deletion is invalid.

Choose symmetric integer representatives for external displacements,
splitting each Nyquist coordinate when $L$ is even. The resulting
smooth Fourier interpolation agrees on the finite grid. Simultaneous
reflection gives zero gradient at the origin, and the second moment
in \eqref{ba:eq:moment} bounds its global Hessian by $C'g^2$.
Taylor's theorem in any anchored momentum coordinates proves
\eqref{ba:eq:cancel} globally on the fundamental domain. Circular
distances obey $|q|_{\mathbb T}\le\sum_{i\ne v}|k_i|_{\mathbb T}$
when the charges sum to zero, so aliasing does not remove this bound.
The preceding arguments apply for $t<1$. The uniform same-sign inverse
bound in \eqref{ba:eq:weighted-rows} implies that $I-Q$ is invertible.
These inverses and the finite stable-dissection formula extend continuously
to $t=1$; taking this limit proves \eqref{ba:eq:moment} and
\eqref{ba:eq:cancel} there.
\end{proof}

\begin{proposition}[Primitive-loop bounds]\label{ba:prop:primitive}
Under the parameter assumptions of Lemma~\ref{ba:lem:kernels}, use
$\ell_t,\zeta_t$ from \eqref{ba:eq:controls} and the primitive loops
\eqref{fo:eq:primitive}. Uniformly for $0\le t<1$, for every
fixed sign word of length \(k\ge2\) and \(g>0\),
\[
 |\cK_k(t)|\le C'\zeta_t^{k-1}.
\]
The constant is uniform for $0<g\le C$. For every fixed $B_L,\nu,A>0$, on
$1-t\ge W^{-3}$ and $L\le W^{B_L}$, the joint tail is $O(W^{-A})$
beyond $W^\nu\ell_t$.
The absolute pair row is at most \(C'/(\gamma_g W(1-t))\), with a
constant independent of \(g\), and the exact primitive Ward
identities hold.
\end{proposition}
\begin{proof}
For $k=2$ this follows from Lemma~\ref{ba:lem:kernels} and the
explicit two-loop formulas \eqref{fo:eq:pairs}. For $k\ge3$, put
$F_k\coloneqq W^{k-1}\widehat\cK_k(t)$ and use $H$ from
\eqref{ba:eq:inverse-average}; every integral below is normalized finite-grid averaging. Here
$\gamma_g$ is bounded above and below, so its fixed-order powers are
absorbed in the constants; the final estimate is stated using $\zeta_t$.
We prove the full and mixed-leg marginal estimates
\begin{equation}\label{ba:eq:norms}
 \|F_k\|_1\le C'H^{k-1},\qquad
 A_v(q)\coloneqq \int_{k_v=q}|F_k|
       \le C'H^{k-2}U_t(q).
\end{equation}
For an alternating molecule, \eqref{ba:eq:cancel} and
\eqref{ba:eq:coercivity} cancel one external mixed inverse in each
summand at bounded cost, since the numerator is bounded by
$\sum_i((1-t)+g^2|k_i|_{\mathbb T}^2)$.
Choose the cancelled momentum as dependent; all remaining $k-1$
independent inverse integrals cost $H$. For a fixed mixed leg
$q$, use the circular triangle inequality (or anchor at that leg)
to cancel a different leg; the remaining $k-2$ free integrals
cost $H$ and retain $U_t(q)$. For a nonalternating molecule there
is already a stable external leg; make it dependent. At most $k-1$
mixed integrations remain, or $k-2$ when one leg is retained.
The stable core is bounded by Lemma~\ref{ba:lem:moment}.
These single-molecule estimates imply \eqref{ba:eq:norms} for every $k\ge3$.

Mixed gluing of dressed children is exactly $\gamma_g tF_LF_R/U_t(q)$;
each child already contains the cut inverse. Since $\gamma_g t\le1$
and $k_L+k_R=k+2$,
\[
 \int U_t(q)^{-1} A_L(q)A_R(-q)\,dq
 \le\sup_qU_t(q)^{-1} A_L(q)\,\|F_R\|_1
 \le C'H^{k-1}.
\]
For a designated retained momentum $p$ on the left child, denote
its two-leg absolute marginal by $A_L(p,q)$. Bound the right child
by \(\sup_q A_R(-q)/U_t(q)\le
C'H^{k_R-2}\) before integrating the left child.
Since \(\int A_L(p,q)\,dq=A_L(p)\), the resulting bound is
\(C'H^{k-2}U_t(p)\).
The same argument applies with left and right interchanged.
This proves the second norm in \eqref{ba:eq:norms} without assuming
an additional two-leg estimate. Induction along the finite molecule tree gives
the bounds. Fourier inversion gives the stated pointwise bound on the primitive loops.
For tails, retain an external exponential weight along a connecting
tree in every dissection. Stable edges have bounded weighted rows;
mixed inverses have weighted row $C'/(1-t)$ at scale $\ell_t$.
For some fixed $C',c>0$ this gives
$|\cK_k(t)|\le C'W^{1-k}(1-t)^{-C'}
e^{-c\diam(\boldsymbol a)/\ell_t}$. Since $(1-t)^{-1}\le W^3$,
the tail is arbitrarily small past $W^\nu\ell_t$.
Ward and initial words are the exact algebra already verified.
\end{proof}

\noindent\textbf{Proof of Lemma~\ref{fo:lem:deterministic}(a).}
Its parameter assumptions imply those of Lemma~\ref{ba:lem:kernels},
and \eqref{ba:eq:controls} uses the scales \eqref{fo:eq:scale}.
Proposition~\ref{ba:prop:primitive} therefore proves
\eqref{fo:eq:primitivebound} in case \textup{(a)}, its length-two
bound, and the stated absolute pair-row bound. Its weighted estimate
also gives the tail at $L\le W^2$ and $1-t\ge W^{-3}$.

\subsection{The polynomial large-coupling auxiliary estimate}\label{app:polynomial-primitives}
\begin{proposition}[Large-coupling primitive bounds]\label{ll:prop:primitive}
Under the assumptions of Lemma~\ref{ll:lem:deterministic}, every fixed primitive word of length
$k\ge2$ satisfies
\[
 |\cK_k(t)|\le C(1+\log(g^2))^{6k}\zeta_t^{k-1}.
\]
It has arbitrarily small polynomial joint tails beyond
$W^\nu\ell_t$ for every fixed $\nu>0$. The absolute pair row is
$O(g^2/(W(1-t)))$ with no logarithmic factor.
\end{proposition}
\begin{proof}
Use the loops, primitive hierarchy and Ward identity of
\eqref{fo:eq:primitive} and \eqref{gc:eq:ward}, now with cut coefficient
\(W/g^2\), pair loops \(\cK_2^{+-}(t)=(g^2/W)P(I-tP)^{-1}\) and
\(\cK_2^{++}(t)=(g^2/W)Q(I-tQ)^{-1}\), and Ward denominator \(2iW\eta_t\).
All Fourier averages below are normalized finite-grid averages. Write
$m\coloneqq m_{g,L}(E)$ and $F_k\coloneqq W^{k-1}\widehat{\cK}_k(t)$.
Use $f,U_t$ and their Fourier bounds
\eqref{ll:eq:forest-input}--\eqref{ll:eq:low-frequency-input} from
Lemma~\ref{ll:lem:deterministic}.
The density, free-row and stable-inverse bounds are those of
Lemma~\ref{uvp:prop:finite}; the mixed row, pointwise
and weighted bounds are \eqref{da:eq:uniform-mixed-row} and
\eqref{da:eq:weighted-actual-kernels}. Their volume conditions follow
from $L\ge C_Ag^2\log W$ and $1-t\ge W^{-3}$ after choosing $C_A$
sufficiently large.

The common dissection formula \eqref{ba:eq:dissection} applies with
$\gamma_g=g^{-2}$: each internal diagonal carries $tV_t/g^2$.
Thus the core normalization and the retained-leg convention agree
with the bounded-coupling argument.

For the stable molecule cores $S^{[j]}$, repeat the weighted leaf
summation in the proof of Lemma~\ref{ba:lem:moment} at scale $g^2$.
A term with $s$ diagonals has $j+s-1$ free rows of cost $Cg^2$
and $s$ links of cost $C/g^2$, giving row bound $Cg^{2j-2}$;
a spatial moment of order $a=1,2$ costs a further $g^{2a}$.
The signed Ward deletion and mixed-reducible subtraction in the
same proof give the zero-momentum bound $O((1-t)g^{2k-2})$:
a normalized deletion costs $Cg^2/(1-t)$ and a zero-mode mixed join
costs $(g^2(1-t))^{-1}$, while $\sum_i j_i=k+2(b-1)$ for a tree
of $b$ molecules. The symmetric Fourier interpolation in that proof,
including the Nyquist convention, and reflection cancel the
first-order term. The second-moment bound therefore gives
\begin{equation}\label{ll:eq:alternating}
 |\widehat S^{[k]}(\boldsymbol k)|
 \le Cg^{2k-2}\left(1-t+g^{4}\sum_i|k_i|^2\right)
 \quad\text{on }|k_i|\le c/g^2.
\end{equation}

The momentum integration retains the factors of \(g^2\).
A stable term with \(s\) diagonals is a connected graph with
\(k\) degree-two external vertices, \(s\) degree-four internal
vertices, and \(k+2s\) edges. An \(s\)-edge forest joins each
internal vertex to an external root; its \(f\)-suprema cancel
\(g^{-2s}\). Its elimination has coefficients \(\pm1\) on the
finite grid. A spanning tree on the external roots has \(k-1\)
edges; outward integration uses \eqref{ll:eq:forest-input}, while
the \(s+1\) chords cost \(\|f\|_1\). With a bounded external
inverse this yields \(C(1+\log(g^2))^{s+1}
(1+\log(g^2)+W\zeta_t)^{k-1}\).
For a retained mixed leg $v$ with momentum $q$, remove its inverse
and choose instead a spanning forest rooted at all other external
vertices. Connectedness guarantees such a forest, with $s+1$ edges.
At $v$, the two incident momenta satisfy
\[
 p_{\rm in}-p_{\rm out}=q.
\]
Eliminating from leaves to roots therefore fixes the parent edge of
$v$ from its other incident momentum and $q$; at each internal vertex
the zero-charge conservation law fixes its parent edge in the same way.
Each solved momentum has coefficient $\pm1$, so this is a bijective
change of the independent finite-grid variables, with no volume factor.
The $s+1$ edge suprema leave $g^{2(s+1)}g^{-2s}=g^2$.
After contracting the forest, the graph has $k-1$ external roots and
$k+s-1$ edges: a spanning tree has $k-2$ edges and leaves $s+1$ chords.
With a bounded nonretained external inverse as the final root, the same
outward integration thus gives $k-2$ tree integrals and $s+1$ chord
integrals, uniformly in $q$.
For a nonalternating molecule a stable external leg is available.
For an alternating one a large external momentum provides it; on
the remaining region use \eqref{ll:eq:alternating} and the
low-frequency bounds \eqref{ll:eq:low-frequency-input}.
The full and one-leg absolute Fourier bounds are thus
\(C(1+\log(g^2))^{2k}(W\zeta_t)^{k-1}\) and
\(Cg^2(1+\log(g^2))^{2k}(W\zeta_t)^{k-2}U_t(q)\), respectively.
When the retained leg is mixed, its momentum is minus the sum of
the other momenta, so the canceled multiplier can be chosen among
the nonretained legs.

Repeat the mixed-gluing induction in the proof of
Proposition~\ref{ba:prop:primitive}, now with join
$tF_LF_R/(g^2U_t(q))$. The extra $g^2$ in the one-leg bound cancels
$g^{-2}$ at each join. For a retained leg, bound the opposite child's
one-leg marginal first, as there; no two-leg supremum is needed.
For a tree of $b$ molecules of orders $j_i$, we have
$\sum_i j_i=k+2(b-1)$ and $b\le k-2$, so the factors
$(1+\log(g^2))^{2j_i}$ contribute at most
$(1+\log(g^2))^{6k}$. Hence
\begin{equation}\label{ll:eq:primitive}
 |\cK_k(t)|\le C(1+\log(g^2))^{6k}\zeta_t^{k-1}.
\end{equation}
Stable edges have range \(g^2\) and mixed edges have range \(\ell_t\).
Retaining a small external exponential weight gives joint tails with
polynomial prefactors. At diameter \(W^\epsilon\ell_t\), these tails are
\(O_D(W^{-D})\) for every fixed \(D>0\). The pair row sums are
\(O(g^2/(W(1-t)))\), with no logarithmic loss.
Since \(g^2\le CW^{2B_g}\), every fixed \((1+\log(g^2))^{6k}\) is
\(W^{o(1)}\). Only the fixed-constant pair row
generates the linear evolution, so a logarithm is never exponentiated.
\end{proof}

\subsection{Uniform physical-space estimates at arbitrary coupling}\label{app:uniform-profile}
The tree and Ward identities used in Appendix~\ref{sec:bl-primitives}
remain exact for arbitrarily large coupling. Their Fourier estimates can,
however, lose a factor depending on $\log g$. We give a physical-space
estimate that removes this loss. The pointwise norm and the absolute row
norm with one external inverse removed close under the same tree gluing.
\subsubsection{An explicit deterministic contract}
Let $\Lambda=\Z/L\Z$ with counting measure, and let $d$ denote cyclic
distance. Fix an independent scale $a\ge1$, $L\ge a$, and $0<\delta\le1$, and put
\[
 \ell=\min\{L,a/\sqrt\delta\},\qquad
 A=\frac{a}{\delta\ell}\ge1.
\]
All constants below may depend on a fixed tensor order but not on
$a,L,\delta$. Convolution uses counting measure. The analytic kernel assumptions, stated in
Lemmas~\ref{uvp:prop:finite} and~\ref{uvp:lem:actual-kernels}, are
\eqref{uvp:eq:stable} and \eqref{uvp:eq:xi}--\eqref{uvp:eq:error}:
throughout this abstract argument, replace $g^2$ by $a$ in the kernel
bounds, replace $\mathcal E$ by $E$, and allow any even kernels satisfying them.
Different external legs may have different stable kernels; the mixed
kernel $U=I+\Xi+E$ is common. No difference estimate on $E$ is used.
The molecule cores are the stable parts of the primitive dissection
from Appendix~\ref{sec:bl-primitives}. Abstractly, a core $S_k:\Lambda^k\to\C$, $k\ge3$, is translation
invariant, invariant under simultaneous reflection, and satisfies
\begin{equation}\label{uvp:eq:coreprofile}
 |S_k(x_1,\ldots,x_k)|\le C
          e^{-c\diam(x_1,\ldots,x_k)/a}.
\end{equation}
Consequently, at every anchor $v$,
\begin{equation}\label{uvp:eq:moments}
 \sum_{\boldsymbol z}|S_k(0,\boldsymbol z)|
       (1+|\boldsymbol z|_1^2/a^2)\le C a^{k-1},
\end{equation}
where the other coordinates use minimal cyclic representatives.
A molecule is called alternating if all its external legs are mixed.
For such a molecule assume the signed row estimate
\begin{equation}\label{uvp:eq:mass}
 \left|\sum_{\boldsymbol z}S_k(0,\boldsymbol z)\right|
       \le C\delta a^{k-1}.
\end{equation}
Every nonalternating molecule has at least one stable external leg.
Dressing means convolution separately in all external coordinates:
\[
 F_k(\boldsymbol x)=\sum_{\boldsymbol y}S_k(\boldsymbol y)
                      \prod_{j=1}^k V_j(x_j-y_j),
 \qquad V_j\in\{U,\text{stable kernels}\}.
\]
For a mixed leg $v$, denote by $F_k^{\langle v\rangle}$ the same
expression with its factor $U$ replaced by $I$. Thus this notation
means removing that factor even if $U$ has not been assumed invertible.
Set
\begin{equation}\label{uvp:eq:norms}
 P(F_k)=\sup_{\boldsymbol x}|F_k(\boldsymbol x)|,
 \quad R_v(F_k)=\sup_{\boldsymbol x_{\ne v}}
       \sum_{x_v}|F_k^{\langle v\rangle}(\boldsymbol x)|.
\end{equation}

\subsubsection{Anchoring, Taylor cancellation, and the molecule estimate}
\begin{lemma}[Anchoring]\label{uvp:lem:anchor}
Suppose a tensor has the core profile \eqref{uvp:eq:coreprofile}.
In a dressed expression, if at least one external factor has absolute
row bounded by a constant, and the absolute value of each remaining
factor has convolution with $e^{-c d(\cdot,0)/a}$ bounded by $CA$
uniformly in the output coordinate, then its
pointwise norm is at most $CA^{k-1}$. If an additional leg is
undressed and summed in absolute value, the bound is
$CaA^{k-2}$.
\end{lemma}
\begin{proof}
The convolution hypothesis bounds the mass of each remaining factor
on every ball of radius $a$ by $CA$. Summing these bounds over annuli
of width $a$ gives the same convolution bound with $c$ replaced by
$c/(k-1)$, at a cost depending only on $k$.
Anchor at a leg $j$ with bounded row. The diameter bound implies
\[
 |S_k(\boldsymbol y)|\le C
  \prod_{i\ne j}e^{-c d(y_i,y_j)/((k-1)a)}.
\]
For a fixed $y_j$, each other convolution costs at most $CA$.
Summing the anchored kernel costs a constant. For the second assertion,
the undressed leg is distinct from $j$; its absolute sum costs
$\sum_x e^{-c d(x,y_j)/((k-1)a)}\le C_k a$. The other $k-2$ convolutions
cost $CA$ each. The argument is uniform in all fixed output coordinates.
\end{proof}

The row and convolution bounds \eqref{uvp:eq:stable} and
\eqref{uvp:eq:localconv}, with $\mathcal E=E$, verify the hypotheses of
Lemma~\ref{uvp:lem:anchor} for $I$, stable kernels, and $E$.
A nonalternating molecule is therefore controlled immediately.
An alternating molecule is expanded using $U=I+\Xi+E$ on each
external leg. Every term except the all-$\Xi$ term has an anchor.
For the undressed-leg row norm, expand the remaining $k-1$ factors;
again only the all-$\Xi$ term lacks an anchor.

\begin{lemma}[Discrete symmetric Taylor bound]\label{uvp:lem:taylor}
Let $B$ be a normed vector space on which translations in each of
$m$ cyclic coordinates are isometries. If $H(\boldsymbol z)$ is a
translated family and every mixed pair of unit differences has
$B$-norm at most $M$, then
\[
 \norm{H(\boldsymbol z)+H(-\boldsymbol z)-2H(0)}_B
      \le |\boldsymbol z|_1^2 M.
\]
Changing the universal constant accommodates a supremum over the
coordinate pairs. If $S(\boldsymbol z)=S(-\boldsymbol z)$, then
\[
 \norm{\sum_{\boldsymbol z}S(\boldsymbol z)H(\boldsymbol z)}_B
 \le \left|\sum S\right|\norm{H(0)}_B
      +C M\sum |S(\boldsymbol z)|\,|\boldsymbol z|_1^2.
\]
\end{lemma}
\begin{proof}
The operator identity
$T_z+T_{-z}-2I=-(T_z-I)(T_{-z}-I)$ is exact. Write each factor
$T_z-I$ as a telescoping sum of $|z|_1$ translated signed unit
differences along a shortest coordinate path. Products give the first
bound. Reflection pairs $z$ with $-z$ in the second assertion. The
proof works directly on the cyclic group, including its Nyquist points.
\end{proof}

\begin{theorem}[Single molecule]\label{uvp:thm:molecule}
Under the kernel assumptions \eqref{uvp:eq:stable} and \eqref{uvp:eq:xi}--\eqref{uvp:eq:error}
and the core assumptions \eqref{uvp:eq:coreprofile}--\eqref{uvp:eq:mass},
\[
 P(F_k)\le C A^{k-1},\qquad
 R_v(F_k)\le C a A^{k-2}
\]
for each mixed external leg $v$.
\end{theorem}

\begin{proof}[Proof of Theorem~\ref{uvp:thm:molecule}]
Only the all-$\Xi$ terms remain. Anchor the core at its first
coordinate and write $z_1=0$. For the pointwise norm define
\[
 H_{\boldsymbol z}(\boldsymbol x)=
       \sum_u\prod_{j=1}^k\Xi(x_j-u-z_j).
\]
Take $B$ to be the supremum norm in $\boldsymbol x$.
One kernel in absolute row norm and the others in pointwise norm give
\[
 \norm{H_0}_B\le C\delta^{-1}(\delta\ell)^{-(k-1)}.
\]
A second difference in a single $z_j$ places two differences on one
kernel, so \eqref{uvp:eq:xi2} gives
\[
 \norm{\Delta_i\Delta_j H_{\boldsymbol z}}_B
 \le C\lambda^2\delta^{-1}(\delta\ell)^{-(k-1)}.
\]
For differences on two distinct kernels use
$\norm{\Delta\Xi}_1\norm{\Delta\Xi}_\infty$
and $k-2$ undifferentiated pointwise factors. This gives the same bound.
Applying Lemma~\ref{uvp:lem:taylor}, \eqref{uvp:eq:mass}, and
\eqref{uvp:eq:moments} yields
\begin{align*}
 P(F_k)&\le C\left[
  \delta a^{k-1}\delta^{-1}(\delta\ell)^{-(k-1)}
  +a^{k+1}\lambda^2\delta^{-1}(\delta\ell)^{-(k-1)}\right]\\
 &\le C A^{k-1}.
\end{align*}

For the row norm, anchor at the undressed leg $v$, call its output
$u$, and put
\[
 H_{\boldsymbol z}(u,\boldsymbol x_{\ne v})
       =\prod_{j\ne v}\Xi(x_j-u-z_j).
\]
Now take $\norm H_B=\sup_{\boldsymbol x_{\ne v}}\sum_u|H(u,\boldsymbol
x_{\ne v})|$. The same row/pointwise calculation gives
\[
 \norm{H_0}_B\le C\delta^{-1}(\delta\ell)^{-(k-2)},\quad
 \norm{\Delta_i\Delta_jH}_B
   \le C\lambda^2\delta^{-1}(\delta\ell)^{-(k-2)}.
\]
The Taylor estimate is therefore at most
$Ca^{k-1}(\delta\ell)^{-(k-2)}=CaA^{k-2}$.
All other terms are controlled by Lemma~\ref{uvp:lem:anchor}.
\end{proof}

\subsubsection{Exact mixed gluing and closure}
Assume now that $U$ is an even invertible convolution kernel, and that
mixed gluing of children of orders $k_L,k_R$ is
\begin{equation}\label{uvp:eq:glue}
 F(\boldsymbol x_L,\boldsymbol x_R)
  =\frac{t}{a}\sum_z F_L(\boldsymbol x_L,z)
                F_R^{\langle\mathrm{cut}\rangle}(z,\boldsymbol x_R),
 \qquad 0\le t\le1.
\end{equation}
Here both cut legs are mixed and $k=k_L+k_R-2$.
The operation of undressing any retained external leg commutes with
gluing. Equivalently, at core level the join is $tU(y-z)/a$ and each
child has one cut-leg factor $U$. Formula \eqref{uvp:eq:glue} follows
by moving one even $U$ through the contraction and removing the other.

\begin{theorem}[Two norms close along mixed trees]\label{uvp:thm:tree}
If every molecule satisfies Theorem~\ref{uvp:thm:molecule}, every fixed
finite mixed tree satisfies
\[
 P(F_k)\le CA^{k-1},\qquad R_v(F_k)\le CaA^{k-2}.
\]
The constants depend on its order and on the fixed finite collection
of allowable tree topologies, but not on $a,L,\delta$.
\end{theorem}
\begin{proof}
For a join, \eqref{uvp:eq:glue} gives
\[
 P(F)\le a^{-1}P(F_L)R_{\mathrm{cut}}(F_R)
       \le C A^{k_L-1+k_R-2}=CA^{k-1}.
\]
If a retained mixed leg $v$ lies on the left, take its absolute row
before summing the cut variable. This gives
\[
 R_v(F)\le a^{-1}R_v(F_L)R_{\mathrm{cut}}(F_R)
          \le CaA^{k_L+k_R-4}=CaA^{k-2}.
\]
The right-retained case uses the symmetric version of the gluing
formula. No two-leg supremum or absolute Fourier integral occurs.
Induction proves the assertions.
\end{proof}

For the actual kernels, the abstract assumptions with $a=g^2$ are
verified by Lemmas~\ref{uvp:prop:finite} and~\ref{uvp:lem:actual-kernels}.

\subsubsection{Primitive cores and their signed sums}
We connect Theorem~\ref{uvp:thm:tree} and the finite-cycle kernel bounds of
Lemma~\ref{uvp:prop:finite} to primitive loops.
The external algebra is the tree representation and
Ward identity of Truong--Yang--Yin~\cite[Assumption 4.1, Lemma 4.16,
Corollary 4.21]{TYY}, with core decomposition
\cite[Definition 4.24 and Lemma 4.25, equation (4.34)]{TYY} and
block reduction \cite[Lemma 4.26, equation (4.38)]{TYY}.

For the block model, put $t=1-\delta$ and use $\gamma_g=g^{-2}$.
Equation~\eqref{da:eq:characteristic-general} verifies the external algebraic
hypothesis, and Lemma~\ref{uvp:prop:finite} supplies the density bounds. Under the block normalization,
$F_k=W^{k-1}\mathcal K_k(t)$ has cut coefficient $g^{-2}$.
Its mixed pair is $F_2=g^2P(I-tP)^{-1}$; its mixed external inverse
is $U=(I-tP)^{-1}$. Core joins are exactly $tU/g^2$, which is
the abstract formula \eqref{uvp:eq:glue} at $a=g^2$.

\begin{proposition}[Core hypotheses from a localized free row]\label{uvp:prop:cores}
Assume the free scalar row and stable inverses satisfy
\[
 |M(x)|\le Ce^{-c d(x,0)/a},\qquad
 \sum_xe^{c d(x,0)/a}|V^{++}(x)|
 +\sum_xe^{c d(x,0)/a}|V^{--}(x)|\le C,
\]
and that $\Im m$ is bounded above and below by positive constants.
Then the stable-only molecule cores in the exact primitive tree
algebra satisfy \eqref{uvp:eq:coreprofile}, reflection, and
\eqref{uvp:eq:mass}. Consequently, if the mixed kernel also satisfies
the bounds \eqref{uvp:eq:xi}--\eqref{uvp:eq:error}, then
\begin{equation}\label{uvp:eq:primitive-final}
 |\mathcal K_k(t)|\le C\left(\frac{a}{W\delta\ell}\right)^{k-1}.
\end{equation}
There is no factor depending on $\log a$ in this implication.
\end{proposition}
\begin{proof}
A stable-only dissection with $s$ internal diagonals has $k+2s$
free-resolvent edges and $s$ stable links, each with a factor $a^{-1}$.
Each stable link joins two internal vertices; the links are disjoint
as links in the expanded graph. Contract these $s$ links. Choose a
rooted forest of $s$ free edges connecting the resulting $s$ internal
vertices to the $k$ external vertices, with one external root per
component. Such a forest exists in every connected expanded graph:
contract all external vertices to a single root, choose a spanning
tree, and undo that contraction. The remaining free edges are bounded
pointwise. Summing the $s$ stable links costs $C^s$ and summing the
$s$ retained free edges costs $(Ca)^s$, canceled by $a^{-s}$.
Retain a smaller exponential weight before summation. Since the full
graph is connected, its total edge length dominates the external
diameter, so the result is \eqref{uvp:eq:coreprofile}. All fixed moments
follow by reducing the exponential constant. Evenness of the free
rows and stable kernels gives simultaneous reflection.

For the signed mass, repeat the Ward deletion and alternating-core
subtraction in the proof of Lemma~\ref{ba:lem:moment}, with
normalized cut coefficient $a^{-1}$ and $1-t=\delta$.
The preceding forest bound gives pure-sign rows $Ca^{k-1}$;
each normalized Ward deletion costs $Ca/\delta$, since $\Im m\asymp1$.
Thus the full signed row is at most $C(a/\delta)^{k-1}$, and removing
the $k$ external zero-mode inverses gives $C\delta a^{k-1}$ for the
complete alternating core. At a tree of $b$ molecules the scale is preserved:
\[
 \prod_i(\delta a^{k_i-1})(a\delta)^{-(b-1)}
 =\delta a^{k-1},\qquad \sum_i k_i=k+2(b-1).
\]
The same subtraction induction, starting at the quadrilateral,
therefore proves \eqref{uvp:eq:mass}.
Theorem~\ref{uvp:thm:tree} then gives \eqref{uvp:eq:primitive-final} for
$k\ge3$; the mixed and stable pair formulas give the $k=2$ case.
\end{proof}

\noindent\textbf{Proof of Lemma~\ref{fo:lem:deterministic}(b).}
Set $\delta=1-t$. For $L\ge C_*g^2\log W$ and
$\delta\ge W^{-3}$, Lemmas~\ref{uvp:prop:finite} and
\ref{uvp:lem:actual-kernels} verify the kernel hypotheses of
Proposition~\ref{uvp:prop:cores} at $a=g^2$; \eqref{da:eq:kernel-accuracy} supplies
the required error bound after choosing $C_*$ sufficiently large.
Its conclusion \eqref{uvp:eq:primitive-final}, with the scale conversion
\eqref{da:eq:primitive-scale-conversion}, is
\eqref{fo:eq:primitivebound} in case (b), including the stated bound at length two.
All constants are independent of $W,g,L,t$.

\subsection{Primitive joint tails}\label{app:primitive-tails}
\noindent\textbf{Proof of \eqref{fo:eq:primitive-tail} in case \textup{(a)}.}
The proof of Proposition~\ref{ba:prop:primitive} gives
$|\cK_k(t;\boldsymbol a)|\le C W^{1-k}(1-t)^{-C}
 e^{-c\diam(\boldsymbol a)/\ell_t}$, uniformly in the signs.
Since $(1-t)^{-C}\le W^{3C}$, this is at most $W^{-Q}$ beyond
$W^\epsilon\ell_t$ for sufficiently large $W$. If $\ell_t=L$,
the tail domain is empty. For arbitrarily large coupling we use the
following weighted estimate, also needed in later appendices.

\begin{corollary}[Coupling-independent joint tails]\label{uvp:cor:joint-tail}
Fix $0<\kappa<2$ and an integer $k\ge2$. Assume $g\ge C_\kappa$ and
$|E|\le2-\kappa$ as in Lemma~\ref{uvp:prop:finite}. Put
$\delta=1-t$, $0<\delta\le1$, and $\ell=\min\{L,g^2/\sqrt\delta\}$.
If $\delta\ge W^{-3}$ and $L\gg g^2\log W$, meaning
$L/(g^2\log W)\to\infty$, then for every $\epsilon,Q>0$ and sufficiently
large $W$,
\[
 \sup_{\diam(\boldsymbol x)>W^\epsilon\ell}
           |\mathcal K_k(t;\boldsymbol x)|\le W^{-Q}.
\]
All constants are independent of $g,L$; the width threshold is uniform
whenever $L/(g^2\log W)$ diverges uniformly.
More precisely the exact primitive satisfies the crude weighted bound
\[
 |W^{k-1}\mathcal K_k(t;\boldsymbol x)|
       \le C\delta^{-C}e^{-c\diam(\boldsymbol x)/\ell}
          +O(W^{-Q})
\]
with arbitrarily prescribed $Q$.  Sharp cancellation is not needed
for this tail estimate.
\end{corollary}
\begin{proof}[Proof of Corollary~\ref{uvp:cor:joint-tail} and
\eqref{fo:eq:primitive-tail} in case \textup{(b)}]
Use the weighted mixed rows and perturbation bound
\eqref{da:eq:weighted-actual-kernels}, the stable rows
\eqref{uvp:eq:stable}, and the free row \eqref{da:eq:uniform-free-row}.
Their costs are $C\delta^{-1}$, $C$, and $Cg^2$, respectively.
Equation~\eqref{da:eq:kernel-accuracy}, with its fixed exponent $r$
chosen after the number of inverse factors, makes the total finite-cycle
error smaller than $W^{-Q}$ for every prescribed $Q$.

In each fixed-order exact dissection, an internal covariance link
contributes $g^{-2}$.  Use the rooted forest from
Proposition~\ref{uvp:prop:cores} at $a=g^2$: one free-edge row $Cg^2$ for every
internal link cancels these factors.  Sum every stable or mixed
link and every external inverse in weighted absolute row.  Their
number is bounded in terms of $k$, so their combined cost is
$C\delta^{-C}$.  Bound the remaining free edges pointwise.
The connected graph has total edge length at least the external
diameter, so retaining a sufficiently small common exponential
weight gives the displayed bound.  The number of dissection
topologies is finite at fixed $k$.  If $\ell=L$, the displayed
joint-tail domain is empty for large $W$; otherwise
$\delta^{-C}\le W^{3C}$ and
$e^{-cW^\epsilon}$ beats every fixed power of $W$.
Finally, \eqref{da:eq:primitive-scale-conversion} gives
$\ell\le\ell_t\le2\ell$. Hence $\diam(\boldsymbol x)>W^\epsilon\ell_t$
implies $\diam(\boldsymbol x)/\ell>W^\epsilon$, proving
\eqref{fo:eq:primitive-tail} with its stated scale. All estimates are
uniform in the signs and $0\le t\le1-W^{-3}$; their constants are
independent of $g,L$, while the required error accuracy determines
how large $L/(g^2\log W)$ must be.
\end{proof}

\subsection{Primitive inputs for initialization and stopping}
\label{app:primitive-stochastic-inputs}
The following estimates supply the primitive terms in the initial error
equation, the stopped pair equation, and the polynomial-range time increments.

\begin{lemma}[Early-time primitive estimates]\label{b:lem:early-primitive}
Fix $k\ge2$.
\textup{(a)} Under the bounded-coupling conditions of
\eqref{da:eq:propagation-domain}, there are $C_k,c_k>0$ such that
\begin{equation}\label{b:eq:early-primitive}
 \sup_{0\le t\le1/2}\max_{\boldsymbol a,\boldsymbol\sigma}
 e^{c_k\diam(\boldsymbol a)}
 |\cK_k(t;\boldsymbol a,\boldsymbol\sigma)|\le C_kW^{1-k}.
\end{equation}
\textup{(b)} Fix $B_g,B_L>0$ and $\kappa\in(0,2)$, and assume
$C_\kappa\le g<W^{B_g}$, $g^2\log W\ll L\le W^{B_L}$ and
$|E|\le2-\kappa$, with $C_\kappa$ as in
Lemma~\ref{fo:lem:deterministic}\textup{(b)}. Then
\begin{equation}\label{b:eq:early-free-label}
 \sup_{0\le t\le1/2}\max_{\boldsymbol a_{\ne i},\boldsymbol\sigma}
 \sum_{a_i\in\Z_L}|\cK_k(t;\boldsymbol a,\boldsymbol\sigma)|
 \le C_kg^2W^{1-k},\qquad 1\le i\le k.
\end{equation}
The constants are independent of $W,g,L$ in their stated ranges; the
width threshold in \textup{(b)} may depend on the divergence rate of
$L/(g^2\log W)$.
\end{lemma}
\begin{proof}[\textbf{Proof of \eqref{b:eq:early-primitive} and
\eqref{b:eq:early-free-label}}]
For \textup{(a)} with $g>0$, use the weighted dissection estimate
in the proof of Proposition~\ref{ba:prop:primitive}:
$|\cK_k|\le C_kW^{1-k}(1-t)^{-C_k}
e^{-c_k\diam(\boldsymbol a)/\ell_t}$.
The bounded-coupling hypotheses imply \eqref{bl:eq:intrinsic};
on $t\le1/2$, both $(1-t)^{-1}$ and $\ell_t$ are bounded.
This proves \eqref{b:eq:early-primitive}. At $g=0$, all free rows and
external inverses are diagonal with bounded entries on $[0,1/2]$;
\eqref{ba:eq:dissection} and \eqref{fo:eq:pairs} give the same bound.

For \textup{(b)}, \eqref{da:eq:uniform-free-row} verifies the entry
bound used in the proof of \eqref{uvc:eq:primitive-coarse-row}.
Applying that estimate with $\gamma_g=g^{-2}$ and $t\le1/2$ gives
$\sum_{a_i}|\cK_k(t)|\le C_kg^2W^{1-k}$, uniformly in the remaining
labels and signs, proving \eqref{b:eq:early-free-label}.
\end{proof}

\begin{lemma}[Primitive pair bounds for continuation]
\label{b:lem:primitive-pair}
Assume \eqref{da:eq:propagation-domain}. For $0\le v<1$,
\begin{equation}\label{b:eq:primitive-pair-row}
 \max_{a,\sigma_1,\sigma_2}\sum_b
 |\cK_2(v;a,b;\sigma_1,\sigma_2)|
 \le\frac1{\gamma_gW(1-v)}.
\end{equation}
Fix $0<\theta\le1/2$, $c_{\rm p}>0$ and $F>0$, and put
$w_v(x)=e^{-c_{\rm p}(\|x\|_L/\ell_v)^\theta}+W^{-F}$.
For sufficiently large $W$, uniformly for $1/2\le v\le1-W^{-3}$,
\begin{equation}\label{b:eq:primitive-long-pair}
 |\cK_2(v;a,b;\sigma_1,\sigma_2)|
 \le C\zeta_v^2w_v(a-b)
 \quad\text{if }\|a-b\|_L\ge\ell_v(\log W)^{3/2}.
\end{equation}
\end{lemma}
\begin{proof}[\textbf{Proof of \eqref{b:eq:primitive-pair-row} and
\eqref{b:eq:primitive-long-pair}}]
Since $|Q|\le P$ and $P\mathbf1=\mathbf1$, the Neumann series give
$|Q(I-vQ)^{-1}|\le P(I-vP)^{-1}$. Formula
\eqref{fo:eq:pairs} and the mixed row mass $(1-v)^{-1}$ prove
\eqref{b:eq:primitive-pair-row} for all signs.
For \eqref{b:eq:primitive-long-pair}, set $r=\|a-b\|_L/\ell_v$.
The same pointwise comparison and \eqref{fo:eq:kernel-bounds} give
\[
 |\cK_2(v;a,b;\sigma_1,\sigma_2)|
 \le C\zeta_v e^{-cr}+C_A\zeta_vW^{-A}.
\]
For $r\ge(\log W)^{3/2}$, the ratio $e^{-cr+c_{\rm p}r^\theta}$ is
$O_D(W^{-D})$ for every fixed $D>0$.
By \eqref{gc:eq:scale}, $\zeta_v\ge c/W$; consequently the first
term is bounded by $C\zeta_v^2w_v$. Choose $A>F+2$ for the second
term and use $w_v\ge W^{-F}$. This proves the stated weighted bound.
\end{proof}

\section{Initialization and propagation of the loop estimates}
\label{app:gaussian-continuation}
\label{sec:fo-propagation-proof}\label{sec:gc-scheme}

We first prove the Gaussian and Haar moment estimates used in
Section~\ref{sec:input-averaging}, then the bounded-coupling initialization
of Lemma~\ref{fo:lem:bounded-initial} and the propagation theorem,
Theorem~\ref{gc:thm:common}.
The common five-step argument follows
\cite[Sections~5.1--5.7]{YY} and \cite[Sections~7.1--7.6]{TYY}.

In Appendices~\ref{sec:bounded-proof}--\ref{sec:gc-schedule},
$\gamma=\gamma_g$ and the parameter ranges are those of
Theorem~\ref{gc:thm:common}, including the earlier-time pairs required
in part \textup{(i)}. Spatial maxima and time suprema are inside
probability; suprema over deterministic parameters are outside,
as in that theorem. Appendix~\ref{sec:classical-random-matrix}
uses the hypotheses specified with its target estimates.

\subsection{Gaussian integration by parts and Haar estimates}
\label{sec:classical-random-matrix}

We give the Gaussian integration-by-parts and Haar calculations used
in Section~\ref{sec:input-averaging}. The residual estimates use either
parameter range of Lemma~\ref{gc:lem:residual}, including its comparison
points, and the deterministic Dyson rows
\eqref{gc:eq:residual-row-bounds}. The Haar estimate applies to an
arbitrary deterministic complex matrix.

\begin{proof}[\textbf{Proof of \eqref{gc:eq:residual-fixed-moment} and \eqref{fo:eq:residual-moments}}]
Fix deterministic $u,z$ as in Lemma~\ref{gc:lem:residual}, and abbreviate
$M_B=M_{B,u}(z)$ and $M=M_u(z)$. Use $R$ from
\eqref{gc:eq:residual-input} and $Z$ from
Section~\ref{sec:input-averaging}. We first prove
\eqref{gc:eq:residual-fixed-moment} under its moment hypotheses.
The assertion is immediate at $u=0$, since $R=0$; assume $u>0$.

Choose a real orthonormal basis $(T_\alpha)_{\alpha=1}^{W^2}$ of
Hermitian $W\times W$ matrices for the inner product $\Tr(AB)$,
and let $T_{c,\alpha}$ denote its embedding in block $c$. Thus
\[
 \mathcal B_u^{(g)}=\sum_{c,\alpha}\xi_{c,\alpha}T_{c,\alpha},\qquad
 \xi_{c,\alpha}\sim N(0,u\gamma_g/W)\ \text{independently},\qquad
 \sum_\alpha(T_\alpha)_{ij}(T_\alpha)_{k\ell}
       =\delta_{i\ell}\delta_{jk}.
\]
Write $\lambda=(c,\alpha)$ for a Gaussian coordinate. The identities
$\partial_\lambda G=-GT_{c,\alpha}G$ and
$\sum_\alpha T_\alpha A T_\alpha=(\Tr A)I_W$ imply
\begin{equation}\label{gc:eq:divergence}
 R_{ab}=\sum_\lambda\left(
       \xi_\lambda f_\lambda-W^{-1} u\gamma_g\partial_\lambda f_\lambda
                            \right),\qquad
 f_{c,\alpha}\coloneqq W^{-1} M_{B,ac}\Tr(T_\alpha G_{cb}).
\end{equation}
Indeed, the first term is $W^{-1}\Tr(M\mathcal B_u^{(g)}G)_{ab}$,
and the second is
\[
 -W^{-1} u\gamma_g\sum_{c,\alpha}\partial_{c,\alpha}f_{c,\alpha}
 =W^{-2} u\gamma_g\sum_cM_{B,ac}
       (\Tr G_{cc})(\Tr G_{cb})
 =u\gamma_g\sum_cM_{B,ac}\langle G\rangle_{cc}\langle G\rangle_{cb}.
\]
In particular, $\E R_{ab}=0$. More generally, for a smooth random
variable $F$, Gaussian integration by parts gives the exact cancellation
\[
 \E[R_{ab}F]=W^{-1} u\gamma_g
                  \sum_\lambda\E[f_\lambda\partial_\lambda F].
\]

We apply this cancellation to $\E|R_{ab}|^p$.
Label its $p$ factors by $1,\ldots,p$, the first $p/2$ being $R_{ab}$
and the others $\overline{R_{ab}}$, and write $f_{j,\lambda}$ for
$f_\lambda$ or $\overline{f_\lambda}$ accordingly.
Expanding \eqref{gc:eq:divergence} and integrating all Gaussian
coordinates by parts gives the following finite sum.
Choose $q$ disjoint pairs of factors; their Gaussian coordinates
are paired together. Each of the other $p-2q$ coordinates from a
factor $i$ differentiates a coefficient at a factor $j\ne i$.
For such a choice $\mathcal D$, there are $e=p-q$ freely summed
Gaussian labels, and its contribution is
\begin{equation}\label{gc:eq:residual-pairing}
 \left(W^{-1} u\gamma_g\right)^e
 \sum_{\boldsymbol\lambda\,\text{respecting the pairs}}
 \E\prod_{j=1}^p
     \left[\left(\prod_{i:\,i\to j}\partial_{\lambda_i}\right)
                    f_{j,\lambda_j}\right].
\end{equation}
The sum of these contributions is exactly $\E|R_{ab}|^p$.
To see why $j=i$ is excluded, a derivative acting on its own
coefficient occurs both from integration by parts and from the
negative term in \eqref{gc:eq:divergence}, with opposite signs.
These terms cancel even when other derivatives act on that coefficient,
since the derivatives commute. The number of terms depends only on $p$.

\emph{Block indices.}
An $h$-fold derivative of $f_{c,\alpha}$ is a sum of $h!$ terms
of the form
\[
 W^{-1} (-1)^hM_{B,ac}
 \Tr\bigl(T_\alpha G_{cd_1}T_{\alpha_1}
          G_{d_1d_2}\cdots T_{\alpha_h}G_{d_hb}\bigr),
\]
where $(d_1,\alpha_1),\ldots,(d_h,\alpha_h)$ are the incoming
Gaussian labels in one of their possible orders. Conjugated
coefficients give the reversed adjoint words. Thus every term has
$2e$ resolvent factors and $2e$ basis-matrix occurrences.
Each freely summed block label $c$ retains the factor $M_{B,ac}$
from its original coefficient, or two such factors when its two
original Gaussian coordinates are paired. No derivative acts on $M_B$.
Consequently \eqref{gc:eq:residual-row-bounds} gives
\[
 u\gamma_g\sum_c|M_{B,ac}|\le C,\qquad
 u\gamma_g\sum_c|M_{B,ac}|^2\le C.
\]
After applying the uniform moment bounds to the resolvent entries,
these estimates sum all $e$ block labels and cancel the factor
$(u\gamma_g)^e$ in \eqref{gc:eq:residual-pairing}.
The remaining normalization is $W^{-p-e}$.
In particular, no free block sum is weighted only by the random
array $(\langle G\rangle_{ab})_{a,b}$.

\emph{Orbital indices.}
Number the $2e$ basis-matrix occurrences. Let $\sigma$ give the
next occurrence in each original trace, and let $\pi$ pair the two
occurrences of each Gaussian label. If $F_h$ is the resolvent block
following occurrence $h$, the completeness identity for $(T_\alpha)$
reduces the orbital sum to
\begin{equation}\label{gc:eq:orbital-contraction}
 \sum_{i_1,\ldots,i_{2e}=1}^W
       \prod_{h=1}^{2e}(F_{\pi(h)})_{i_h,i_{\sigma\pi(h)}}.
\end{equation}
Each cycle of $\sigma\pi$ gives one trace of resolvent blocks.
There is no cycle of length one: $\sigma\pi(h)=h$ would pair
occurrence $h$ with its predecessor in the same original trace,
whereas every surviving Gaussian pair joins different original factors.
Thus each resulting trace contains at least two resolvent factors.

Only now substitute $G_{cd}=\langle G\rangle_{cd}I_W+Z_{cd}$, and likewise for
adjoints. If an expanded term has $x$ factors from $Z$ and $z_0$
traces containing only scalar factors, its orbital sum has at most
$W^{x+z_0}$ terms. Indeed, a trace with $x_F>0$ error factors has
at most $x_F$ free indices; an entirely scalar trace has one.
Since each entirely scalar trace contains at least two factors,
$x+2z_0\le2e$. H\"older's inequality uses at most $2e\le2p$
random factors and bounds each entry product by
$CK^{2e}\chi^x$. Hence every contribution is bounded by
\begin{equation}\label{gc:eq:facecount}
 CK^{2e}W^{x+z_0-p-e}\chi^x
 \le CK^{2p}W^{-p+x/2}\chi^x
 \le CK^{2p}\chi^{2p},
\end{equation}
where the last inequality uses $\chi\ge W^{-1/2}$.
Summing the finitely many terms proves
\eqref{gc:eq:residual-fixed-moment}. The expansion was performed
after every derivative, so it requires neither deterministic block traces
nor independence between the block traces and $Z$.

For \eqref{fo:eq:residual-moments}, repeat
\eqref{gc:eq:divergence}--\eqref{gc:eq:orbital-contraction} and
substitute $G_{cd}=M_{B,cd}I_W+\wt G_{cd}$ only after differentiation.
The scalar factors are deterministic and bounded, while part (iii)
controls the error entries in $L^{2p}$. Thus \eqref{gc:eq:facecount}
applies with fixed $K$, giving $C_p\chi^{2p}$ without a $W^\epsilon$ loss.
\end{proof}

\begin{proof}[\textbf{Proof of \eqref{gc:eq:haar-fixed-moment}}]
Let $U\in U(W)$ be Haar and $A\in\C^{W\times W}$ be deterministic.
For $i=j$, this is the centered quadratic-form moment bound for a
uniform unit vector, obtained by normalizing a complex Gaussian vector.
For $i\ne j$, condition on row $i$ of $U$: row $j$ is uniform on
its orthogonal complement. Its linear-form moment bound, followed
by averaging row $i$, gives the same estimate.
\end{proof}

\subsection{Initialization at bounded coupling}
\label{sec:bounded-proof}\label{sec:bl-seed}

\begin{proof}[\textbf{Proof of Lemma~\ref{fo:lem:bounded-initial}: \eqref{gc:eq:iteration-input}}]
Fix $g\le C$, $|E|\le B$ and
$\Im m_{g,L}(E)\ge c'>0$, with $\gamma=\gamma_g$, and take $0\le t\le1/2$.
The estimates below are uniform in these parameters. Write
$Z_{ab}\coloneqq G_{ab}-\langle G\rangle_{ab}I_W$ and use the residual of
\eqref{gc:eq:residual-input}. Here $\eta_t\ge\gamma c'/2$, so $\|G_t\|_\op\le2/(\gamma c')$
deterministically. Exponential spatial conjugation changes only
the bounded hopping; choose the weight small enough that its
operator perturbation is at most $\gamma c'/4$. A resolvent Neumann
series gives $\|(G_t)_{ab}\|_\op\le C'e^{-c\|a-b\|_L}$
for every onsite realization, without an onsite norm bound.
Thus every fixed actual word is $O(W^{1-k})$ with a joint
exponential tail. Haar concentration gives $Z\prec W^{-1/2}$,
bounded block traces, and the residual refinement gives $R\prec W^{-1}$.
Starting from $\wt G_0=0$, \eqref{gc:eq:scalar} and stable inversion
select the small quadratic branch by continuity:
$\max_{a,b}|\langle\wt G\rangle_{ab}|\prec W^{-1}$.

For every sign pair, Cauchy--Schwarz and Ward give
$\sum_b|\cL_2(a,b)|\le\Im\langle G\rangle_{aa}/(W\eta_t)\le C'/W$.
For the primitive loops, \eqref{b:eq:early-primitive} gives the same
$W^{1-k}$ bound and joint exponential tail at every fixed length.

In a cut difference write
$\cL_p\cL_q-\cK_p\cK_q=\cL_pD_q+D_p\cK_q$.
Its boundary cuts define a linear error operator of norm $C_k$,
by the pair rows and the coefficient $\gamma W$; interior cuts involve
shorter errors. A free block sum costs $C_k$, by exponential
localization. Induction on $k$ gives a drift source $O_\prec(W^{-k})$:
a light insertion is $O(\gamma W\cdot W^{-1}\cdot W^{-k})=O(W^{-k})$.
The raw bracket has size
$\gamma W\cdot W^{1-(2k+2)}=O(W^{-2k})$.
Raw martingale BDG and pathwise Gronwall yield $D_k\prec W^{-k}$.
The terminal linear operator may be random, so no anticipating
stochastic integral with its terminal propagator is used.
For $\diam\le C\log W$, multiplying the bound $D_k\prec W^{-k}$
by $e^{c\sqrt{\diam}}$ costs only $W^{o(1)}$. For $\diam>C\log W$,
the deterministic bound $|D_k|\le C_kW^{1-k}e^{-c'\diam}$ absorbs
the extra factor $W$ when $C$ is large. Thus
$|D_k|\prec W^{-k}(e^{-c\sqrt{\diam}}+W^{-F})$ for every fixed $F$.
The definition \eqref{fo:eq:scale} gives
$\ell_t\asymp1$ and $\zeta_t\asymp W^{-1}$ on $0\le t\le1/2$.

Uniform time follows from a polynomial mesh and Brownian maximal
moment estimates: with arbitrarily high probability, simultaneously
for all mesh intervals,
\[
 \sup_{t\in[t_j,t_{j+1}]}
 \|\mathcal B_t^{(g)}-\mathcal B_{t_j}^{(g)}\|_\op
 \le W^{1/2+\nu}(t_{j+1}-t_j)^{1/2}.
\]
Resolvent identities interpolate. Spectral and coupling meshes use
their deterministic polynomial derivative bounds. These meshes
do not condition any Gaussian integration-by-parts calculation.
\end{proof}

\subsection{One-step propagation}\label{sec:gc-one-step}

For arbitrarily large coupling, the additional verifications occur where
needed: Step~1 uses transfer moments and sampling to obtain spatial and
time suprema; Step~2 controls weighted sums and stopped estimates;
Steps~3--4 retain the covariance factors and the cancellation in the
projected hierarchy. These verifications allow the same five steps to
prove \eqref{gc:eq:one-step-output} in both parameter ranges.

\begin{proof}[\textbf{Proof of Theorem~\ref{gc:thm:common}(i)}]
Fix a deterministic step $[s,t]$ and a pair
$(g,E)$ satisfying \eqref{gc:eq:step-parameters} and the
input and time-step assumptions of part (i). Put
\[
 R\coloneqq\frac{1-s}{1-t},\qquad R^{20}\zeta_t\le W^{-c}.
\]
All estimates below hold for $s\le u\le t$, uniformly over these
pairs. In the step-goal displays, stochastic domination includes the
maximum over all indicated spatial labels and signs and the supremum
over $u\in[s,t]$ inside the probability; the supremum over deterministic
parameters stays outside. Intermediate estimates explicitly restricted
to deterministic labels or a fixed time retain that restriction. The previous-time estimates
are available at $(g,E)$ and at every earlier-time pair
$(g',E')$ specified in the input assumption of part (i).
The distribution comparisons \eqref{fo:eq:backward-matching} transfer
imaginary-part, block operator-norm and block-trace bounds from
$G_s(g',E')$ to $G_u(g,E)$, with time ratio
$(1-s)/(1-u)\le R$.

\paragraph{Step 1: rough loop bounds and a weak local law.}
The goal of this step is to obtain, for every fixed loop length $k$,
\begin{equation}\label{fo:eq:step1}
 \begin{aligned}
 &\sup_{s\le u\le t}\max_{a,b}
       \frac{\|(G_u)_{ab}\|_{\op}}{WR\zeta_u}\prec1,
 &&\sup_{s\le u\le t}\max_{a,b}
       \frac{\|(G_u)_{ab}\|_{\HS}^2}{W^2R\zeta_u}\prec1,\\
 &\sup_{s\le u\le t}\max_a
       \frac{|\langle\wt G_u\rangle_{aa}|}{R^2\zeta_u}\prec1,
 &&\sup_{s\le u\le t}\max_{\boldsymbol a,\boldsymbol\sigma}
       \frac{|\cL_k(u;\boldsymbol a,\boldsymbol\sigma)|}
       {R^{k-1}\zeta_u^{k-1}}\prec1.
 \end{aligned}
\end{equation}
All maxima and time suprema are inside the probability. We then deduce
\begin{equation}\label{fo:eq:step1-density}
 \sup_{s\le u\le t}\max_a\Im\langle G_u\rangle_{aa}\prec1.
\end{equation}

\emph{Starting bounds at a fixed time.}
Fix a deterministic $u\in[s,t]$ and let $(g',E')$ be supplied by
Lemma~\ref{fo:lem:families}. Write $G_s'=G_s(g',E')$. The time-$s$ positive-loop input gives,
for every fixed integer $p\ge1$,
\begin{equation}\label{fo:eq:starting-positive-loop}
 W^{-2p}\|(G_s')_{ab}\|_{\op}^{2p}
 \le W^{-2p}\Tr[((G_s')_{ab}(G_s')_{ab}^*)^p]
 =\Tr[(G_s'E_b(G_s')^*E_a)^p]
 \prec\zeta_s(g')^{2p-1}.
\end{equation}
The equality uses $E_a=W^{-1}P_a$. The bound is simultaneous in
$a,b$ because the positive-loop input already contains that maximum.
After taking the $2p$-th root and dividing by $\zeta_s(g')$, the
remaining factor is at most $CW^{1/(2p)}$, since
$\zeta_s(g')\ge c/W$. For a prescribed loss $W^\epsilon$, choose
$p>1/\epsilon$ and use a loss $W^{\epsilon/2}$ in the root estimate.
This proves $W^{-1}\max_{a,b}\|(G_s')_{ab}\|_{\op}\prec\zeta_s(g')$.
Together with the time-$s$ trace input and
$\zeta_s(g')\asymp\zeta_s$, \eqref{fo:eq:backward-matching} gives
\begin{equation}\label{fo:eq:rough-comparison}
 W^{-1}\max_{a,b}\|(G_u)_{ab}\|_{\op}\prec R\zeta_s,
 \qquad \max_{a,b}|\langle G_u\rangle_{ab}|\prec R.
\end{equation}
These spatial maxima are inherited from the input; they require no
new union over block labels.

The characteristic row and stability bounds are
\eqref{fo:eq:propagation-dyson} from
Lemma~\ref{da:lem:scales-propagators}; the corresponding comparison-segment
bounds are \eqref{gc:eq:positive-stability} from
Lemma~\ref{da:lem:positive-height}.

Let $w$ be the upper comparison point in \eqref{da:eq:comparison-parameters},
set $F=(H_u-w)^{-1}$, and write $G=(H_u-z)^{-1}$ for $z$ on the
vertical segment from $z_u$ to $w$. Thus
$1\le\Im w/\Im z\le CR$, and $F$ has the time-$s$ input bounds
under the Gaussian coupling. In particular,
$W^{-2}\|F_{ab}\|_{\HS}^2\prec\zeta_s$ follows from the length-two
positive loop, while \eqref{fo:eq:starting-positive-loop} gives
$W^{-1}\|F_{ab}\|_{\op}\prec\zeta_s$.
For equal real parts, spectral calculus gives
$\Im G\le(\Im w/\Im z)\Im F$. To obtain the off-diagonal block traces,
Cauchy--Schwarz and the two Hilbert--Schmidt Ward identities give
\[
 |\langle G-F\rangle_{ab}|
 \le |z-w|\sqrt{\frac{\Im\langle G\rangle_{aa}
                              \Im\langle F\rangle_{bb}}
                             {\Im z\Im w}}
 \le\frac{|z-w|}{\Im z}
       \sqrt{\Im\langle F\rangle_{aa}\Im\langle F\rangle_{bb}}
 \prec R.
\]
Adding the time-$s$ trace bound for $F$ proves
$\max_{a,b}|\langle G\rangle_{ab}|\prec R$ throughout the segment.
The resolvent identity and the block Ward identities give
\[
 G-F=(z-w)GF,\qquad
 \|P_aG\|_{\HS}^2=\frac{W\Im\langle G\rangle_{aa}}{\Im z},\qquad
 \|FP_b\|_{\op}^2=\frac{\|\Im F_{bb}\|_{\op}}{\Im w}.
\]
Since $W^{-2}\|F_{ab}\|_{\HS}^2\prec\zeta_s$ and
$W^{-1}\|\Im F_{bb}\|_{\op}\prec\zeta_s$, this yields
\begin{equation}\label{fo:eq:step1-hs-comparison}
 \begin{aligned}
 W^{-2}\|G_{ab}\|_{\HS}^2
 &\le 2W^{-2}\|F_{ab}\|_{\HS}^2
   +\frac{2|z-w|^2}{W\Im z\Im w}
      \Im\langle G\rangle_{aa}\|\Im F_{bb}\|_{\op}\\
 &\prec\zeta_s\bigl(1+R\Im\langle G\rangle_{aa}\bigr)
 \prec R^2\zeta_s.
 \end{aligned}
\end{equation}
The displayed spatial bounds are simultaneous at each fixed $(u,z)$;
the sharper centered-trace estimate has not yet been used.

\emph{Residual estimates before Dyson stability.}
Write $R_{ab}(u,z)$ for the residual in \eqref{gc:eq:residual-input}.
By \eqref{gc:eq:haar} and \eqref{fo:eq:step1-hs-comparison}, for each
fixed pair $a,b$,
\[
 \max_{i,j}|G_{ai,bj}-\delta_{ij}\langle G\rangle_{ab}|
       \prec R\sqrt{\zeta_s}.
\]
We use these pointwise entry estimates to bound residual moments under
the original Gaussian law. For each fixed moment order $p$, any scalar
field $X$ appearing here and any event $\Omega^c$ of probability at most
$W^{-M}$, the coarse moments in Lemma~\ref{fo:lem:sampling-continuity}
give
\[
 \E[|X|^p\mathbf1_{\Omega^c}]
 \le(\E|X|^{2p})^{1/2}\Prob(\Omega^c)^{1/2}
 \le C_pW^{C_p-M/2}.
\]
Choose $M$ after the required moment order and accuracy. The trace and
entry bounds therefore imply unconditional moment bounds with an
arbitrarily small power of $W$. Apply the fixed-pair moment implication
\eqref{gc:eq:residual-fixed-moment} with $K=CR$ and
$\chi=\max\{\sqrt{\zeta_s},W^{-1/2}\}\asymp\sqrt{\zeta_s}$,
allowing that small power in $K$. For every fixed even $p\ge2$ and
$\delta>0$, we obtain
\begin{equation}\label{fo:eq:step1-residual-moment}
 \max_{a,b}\E|R_{ab}(u,z)|^p
       \le C_{p,\delta}W^{\delta p}(R^2\zeta_s)^p.
\end{equation}
No Gaussian calculation is conditioned on a sampling or stopping event.

We next make \eqref{fo:eq:step1-residual-moment} simultaneous in space
and in the spectral heights needed for the stability argument. Fix a
desired failure exponent $D>0$ and a sufficiently large accuracy exponent
$Q$. At each fixed $(u,z)$, use the residual sample set from
Lemma~\ref{fo:lem:sampling-continuity}, with error $W^{-Q}$. Its size is at most $W^{C_{\rm S}}$. Markov's
inequality at those pairs and the sampling inequality give, for
$\epsilon>\delta$,
\[
 \Prob\!\left(\max_{a,b}|R_{ab}(u,z)|
                  >CW^\epsilon R^2\zeta_s+W^{-Q}\right)
 \le C_{p,\delta}W^{C_{\rm S}-p(\epsilon-\delta)}
       +\Prob(\text{sampling fails}).
\]
The sample-size exponent is fixed before choosing $p$. Thus the failure
exponent can absorb any fixed polynomial number of time and height
mesh points; this is where spatial sampling replaces a union over $L^2$.

For each deterministic time mesh point, use the closed comparison-height
consequence \eqref{fo:eq:comparison-height-interpolation}. Choose its
height mesh exponent $q_h>C_H+Q+1$ before the residual moment order.
It gives interpolation errors at most $W^{-Q}$ for $G,M_u(z),R$ on
the whole comparison segment, with arbitrary power probability.
Only the coarse inputs of Sections~\ref{app:uniform-compression} and
\ref{sec:sampling-propagation-extensions}, together with the deterministic
positive-height bounds, enter this step; the random Dyson
closure below has not yet been used.

We also prepare time interpolation, using the coarse increments rather
than the conclusion \eqref{fo:eq:step1}. Fix the required loop lengths
first, including a repeated positive loop of sufficiently large fixed
length for the final operator estimate. For a scalar field $X$ in
Lemma~\ref{fo:lem:sampling-continuity}(b) and a fixed normalization
power $j\ge0$, put $Y(u)=X(u)/\zeta_u^j$. Choose a fixed $C_0$,
independent of the moment order $p$, for the coarse moments and increments
of the finitely many fields needed below. Using $\zeta_u\ge c/W$ and
$|\partial_u\log\zeta_u|\le W^3$ from the definition \eqref{fo:eq:scale},
\[
 \begin{aligned}
 \|Y(v)-Y(u)\|_{L^p}
 &\le C_pW^{C_0+j}|v-u|^{1/2}
       +C_pW^{C_0+j+3}|v-u|\\
 &\le C_pW^{C_1}|v-u|^{1/2},\qquad C_1=C_0+j+3.
 \end{aligned}
\]
Request \eqref{fo:eq:weighted-sampling-input} with $w=1$, accuracy
exponent $Q+j+2$ and failure exponent $D+2$. Its sample set
$\mathcal S_X$ is independent of time and has at most $W^{C_{\rm S}}$
tuples. Choose a deterministic time mesh $\mathcal T\subset[s,t]$
with spacing $h=W^{-q}$, including the endpoints. On each mesh interval,
dyadic subdivision and Minkowski give, for $p>2$,
\[
 \left\|\sup_{u\in[v,v+h]\cap[s,t]}|Y(u)-Y(v)|\right\|_{L^p}
 \le C_pW^{C_1}\sum_{m\ge0}2^{m/p}(h2^{-m})^{1/2}
 \le C_pW^{C_1}h^{1/2}.
\]
Markov's inequality and the union over sampled tuples and intervals give
\begin{equation}\label{fo:eq:time-mesh-union}
 \begin{aligned}
 &\Prob\!\left(\max_{\boldsymbol a\in\mathcal S_X}
 \max_{v\in\mathcal T}\sup_{u\in[v,v+h]\cap[s,t]}
 |Y(u;\boldsymbol a)-Y(v;\boldsymbol a)|>W^{-Q}\right)\\
 &\qquad\le C_pW^{C_{\rm S}+(C_1+Q)p}h^{p/2-1}.
 \end{aligned}
\end{equation}
Choose $q>\max\{4,2(C_1+Q+1)\}$ first, and then
$p(q/2-C_1-Q)>C_{\rm S}+q+D+2$. The probability is at most
$W^{-D-1}$. The sampling error divided by $\zeta_u^j$ is
$O(W^{-Q-2})$, so on the sampling and continuity events,
\begin{equation}\label{fo:eq:step1-time-interpolation}
 \sup_{s\le u\le t}\max_{\boldsymbol a}|Y(u;\boldsymbol a)|
 \le C\max_{\substack{v\in\mathcal T\\
                      \boldsymbol a\in\mathcal S_X}}
             |Y(v;\boldsymbol a)|+CW^{-Q}.
\end{equation}
The $W^2$ orbital choices for entries are included in $C_{\rm S}$.
All mesh and sample-size exponents are fixed before the moment orders
and pointwise probability accuracies.

Apply this construction first to the residual, with $j=1$.
At the time and height mesh points choose the pointwise failure exponent
larger than $D+q+q_h+2$; the residual sample count was already included
in the Markov bound above. Height interpolation gives the second bound
below, while \eqref{fo:eq:step1-time-interpolation} on the characteristic
and $\zeta_s\le\zeta_u$ give the first:
\begin{equation}\label{fo:eq:step1-residual-uniform}
 \begin{gathered}
 \sup_{s\le u\le t}\max_{a,b}
        \frac{|R_{ab}(u,z_u)|}{R^2\zeta_u}\prec1,\\
 \max_{v\in\mathcal T}\ \sup_{\substack{\Re z=\Re z_v\\
                              \eta_v\le\Im z\le\Im w(v)}}
       \max_{a,b}\frac{|R_{ab}(v,z)|}{R^2\zeta_s}\prec1.
 \end{gathered}
\end{equation}
Here $w(v)$ is the upper comparison point at time $v$. The same finite
union and height interpolation make the starting trace bounds
simultaneous on these segments. No local law from Dyson closure has
been used to obtain \eqref{fo:eq:step1-residual-uniform}.

\emph{Dyson stability on the common event.}
Fix a mesh time $v$ and abbreviate $G=(H_v-z)^{-1}$ on its comparison
segment. Set
\[
 U(z)=\max_a|\langle G-M_v(z)\rangle_{aa}|,\qquad
 X(z)=\max_{a,b}|\langle G-M_v(z)\rangle_{ab}|.
\]
The two identities in \eqref{gc:eq:scalar}, with $v$ in place of $u$
and $M_v(z)$ in place of $M$, and the row bounds
\eqref{gc:eq:positive-stability} imply
\[
 X(z)\le\max_{a,b}|R_{ab}(v,z)|
       +C\max_{a,b}|\langle G\rangle_{ab}|\,U(z),\qquad
 U(z)\le C\max_{a,b}|R_{ab}(v,z)|+CU(z)X(z).
\]
For any sufficiently small fixed $\nu>0$, the preceding simultaneous
estimates give $\max|R_{ab}(v,z)|\le W^\nu R^2\zeta_s$ and
$\max|\langle G\rangle_{ab}|\le W^\nu R$ outside an event of arbitrarily
small polynomial probability. Substitution and absorption of the linear
term $CW^\nu R^2\zeta_sU$ yield
\begin{equation}\label{fo:eq:step1-dyson-closure}
 U(z)\le CW^\nu R^2\zeta_s+CW^\nu R\,U(z)^2,\qquad
 W^{2\nu}R^3\zeta_s\le W^{2\nu-c}=o(1).
\end{equation}
The last inequality uses $R^{20}\zeta_t\le W^{-c}$ and
$\zeta_s\le\zeta_t$. Thus the quadratic inequality excludes the interval
between a small branch of size $CW^\nu R^2\zeta_s$ and a large branch
of size $cW^{-\nu}R^{-1}$. At the upper point,
$(H_v-w(v))^{-1}=\xi G_s(g',E')$ and
$M_v(w(v))=\xi M(g',E')$ for the corresponding earlier-time pair,
under the Gaussian coupling of Lemma~\ref{fo:lem:families}, with the
Dyson endpoint identity from Lemma~\ref{da:lem:positive-height}.
The time-$s$ input starts $U(w(v))$ on the small branch. Continuity in
$\Im z$ therefore keeps $U$ on that branch down to $z_v$. Since $\nu$
can be arbitrarily small, this proves, simultaneously at all mesh times,
\[
 \max_{v\in\mathcal T}\max_a
        \frac{|\langle\wt G_v\rangle_{aa}|}{R^2\zeta_s}\prec1.
\]

Insert this bound into \eqref{fo:eq:step1-hs-comparison} at $z=z_v$.
Since $\Im\langle G_v\rangle_{aa}=\Im m+
\Im\langle\wt G_v\rangle_{aa}$ and $\Im m\asymp1$, it gives
$W^{-2}\max_{a,b}\|(G_v)_{ab}\|_{\HS}^2\prec R\zeta_s$.
The operator bound at mesh times is already \eqref{fo:eq:rough-comparison}.
For $k\ge2$, use two Hilbert--Schmidt factors and $k-2$ operator factors
in the cyclic trace:
\[
 \max_{\boldsymbol a,\boldsymbol\sigma}|\cL_k(v)|
 \le W^{-k}\left(\max_{a,b}\|(G_v)_{ab}\|_{\HS}\right)^2
                \left(\max_{a,b}\|(G_v)_{ab}\|_{\op}\right)^{k-2}
 \prec R^{k-1}\zeta_s^{k-1}.
\]
For $k=1$, use $\cL_1(v;a,+)=m+\langle\wt G_v\rangle_{aa}$
and its conjugate. These are spatially uniform bounds at the mesh times.

To pass to all times, apply \eqref{fo:eq:step1-time-interpolation} to
centered diagonal traces with $j=1$ and to loops with $j=k-1$.
The just-proved mesh bounds and $\zeta_s\le\zeta_v$ control its right-hand
side, giving the third and fourth bounds in \eqref{fo:eq:step1}.
No new spatial union is needed at mesh times; sampling reduces the
time interpolation to polynomially many tuples. Positive two-loops give
the Hilbert--Schmidt bound. For the operator norm, repeated positive
$2p$-loops give
\[
 \sup_{s\le u\le t}\max_{a,b}
 \frac{\|(G_u)_{ab}\|_{\op}}{WR\zeta_u}
 \prec\sup_{s\le u\le t}(R\zeta_u)^{-1/(2p)}
 \le CW^{1/(2p)}.
\]
For each requested domination loss $\epsilon>0$, choose a fixed
$p>1/\epsilon$ and use loss $\epsilon/2$ in the loop estimate.
This proves the first bound in \eqref{fo:eq:step1}; $p$ may depend on
$\epsilon$ but not on $W$. All deterministic $(g,E)$ suprema remain
outside the probability. This completes \eqref{fo:eq:step1}.

Finally, $\langle G_u\rangle_{aa}=m+\langle\wt G_u\rangle_{aa}$,
$\Im m\asymp1$ and the time-step condition \eqref{gc:eq:one-step},
together with monotonicity of $\zeta$, give
\[
 R^2\zeta_u\le R^2\zeta_t\le W^{-c},\qquad s\le u\le t.
\]
The third bound in \eqref{fo:eq:step1} therefore proves
\eqref{fo:eq:step1-density}.

\paragraph{Step 2: the sharp local law and a preliminary spatial bound.}
The goal of this step is to obtain the following sharp resolvent bounds
and preliminary two-loop spatial bound:
\begin{equation}\label{fo:eq:step2}
 \begin{gathered}
 |D_2(u;a,b;+,-)|\prec R^{10}\zeta_u^2
       \left(e^{-c(\|a-b\|_L/\ell_u)^\theta}+W^{-F}\right),\\
 W^{-2}\|(G_u)_{ab}\|_{\HS}^2\prec\zeta_u,\qquad
 \|\wt G_u\|_{\max}\prec\sqrt{\zeta_u},\qquad
 \max_{a,b}|\langle\wt G_u\rangle_{ab}|\prec\zeta_u.
 \end{gathered}
\end{equation}
Thus the resolvent bounds are sharp, while the spatial bound on $D_2$
still has the loss $R^{10}$.

We first control the free block sum using the final scale budget
\eqref{fo:eq:spatial-sum-budget} of Lemma~\ref{da:lem:scale-budgets}.

For the mixed pair, \eqref{fo:eq:error-hierarchy} becomes
\begin{equation}\label{fo:eq:pair-error}
 \begin{split}
 dD_2={}&\bigl[P(I-uP)^{-1}D_2+D_2P(I-uP)^{-1}\bigr]\,du\\
       &+\gamma_g W D_2^2\,du+\mathcal W_2\,du+d\mathcal M_2.
 \end{split}
\end{equation}
Write $w_v(a-b)=e^{-c_{\rm p}(\|a-b\|_L/\ell_v)^\theta}+W^{-F}$,
with a sufficiently small fixed $c_{\rm p}>0$. We may increase $F$,
since this strengthens the desired estimate; choose $F>2C$ with $C$
from \eqref{fo:eq:spatial-sum-budget}. The convolution bounds \eqref{da:eq:weight-convolution} apply to
$w_u$ and $\sqrt{w_u}$, since \eqref{fo:eq:spatial-sum-budget} and
$F>2C$ give $LW^{-F/2}\le\ell_u$.
Fix a polynomial barrier range $1\le J\le W^{C_J}$, where $C_J=32$ is
fixed before choosing the tail floor; the barrier used below has this
form since $R\le W^3$. Fix a sufficiently small $\nu>0$.
Stop at the first time $\tau_J$ when
$\max_{a,b}|D_2(v;a,b;+,-)|/(\zeta_v^2w_v(a-b))$ reaches $J$,
or when a normalized block operator norm, block Hilbert--Schmidt norm,
centered diagonal trace, or actual loop of length $3$ or $6$ exceeds
$W^\nu$ times its Step~1 envelope. All spatial maxima of this finite
list are included; these stopping tests depend only on the current
and past flow. Before the stop their bounds are ordinary inequalities.
Every subsequent small source loss denotes a fixed multiple of $\nu$;
choose $\nu$ so that their total is smaller than the eventual first-exit
margin. Before $\tau_J$, the quadratic term
therefore satisfies
\[
 \gamma_gW\sum_c|D_2(u;a,c)D_2(u;c,b)|
 \le C_\theta\gamma_gW\ell_u J^2\zeta_u^4w_u(a-b)
 =\frac{C_\theta J^2\zeta_u^3}{1-u}w_u(a-b).
\]
This is the required large-coupling normalization of the free sum.
The exact contraction \eqref{gc:eq:bracket} retains the same factor
$\gamma_gW$ in the joined six-loop bracket.
The estimate needed from the stopped equation is
\begin{equation}\label{gc:eq:spatialstop}
 \begin{aligned}
 J_{\rm out}&\coloneqq
 \sup_{s\le v\le t\wedge\tau_J}\max_{a,b}
 \frac{|D_2(v;a,b;+,-)|}{\zeta_v^2w_v(a-b)}\\
 &\prec R^6+R^2J^2\zeta_t
       +R^5J^{3/2}\zeta_t+R^3J^{3/2}\sqrt{\zeta_t}.
 \end{aligned}
\end{equation}
We prove this bound by estimating the light triangle and the joined
six-cycle before applying the deterministic propagators.
At time $u$, abbreviate $w_{ab}=w_u(a-b)$ and set
$h_u=\ell_u(\log W)^{3/2}$. An edge is long when its length exceeds
$h_u$. By \eqref{b:eq:primitive-long-pair} and the stopped pair bound,
\[
 W^{-1}\|(G_u)_{ab}\|_{\HS}
 \le C\sqrt J\,\zeta_u\sqrt{w_{ab}}\quad\text{on long edges}.
\]
On a short edge, \eqref{da:eq:short-weight} gives
$w_{ab}^{-1/2}\le W^\nu$ for sufficiently large $W$.
The stopped operator envelope therefore gives the ordinary inequality
\[
 W^{-1}\|(G_u)_{ab}\|_{\op}
 \le CRW^{2\nu}\zeta_u\sqrt{w_{ab}},\qquad u\le\tau_J.
\]
Here $J\ge1$ absorbs the primitive constant. These estimates are
available before any sharp local law or higher-loop tail estimate.

A distant light triangle has $d(a,b)>2h_u$. Its direct edge $ab$
and at least one of $ac,cb$ are long. Choose these two as the
Hilbert--Schmidt factors in its cyclic trace; bound the remaining
edge by its operator norm (or its Hilbert--Schmidt norm if long).
For every sign choice occurring in \eqref{fo:eq:light}, this gives
\[
 |\cL_3(u;a,c,b)|\prec
 CRJ^{3/2}\zeta_u^3\sqrt{w_{ab}w_{ac}w_{cb}}.
\]
The centered insertion has size $R^2\zeta_u$ by Step~1. Multiplying
by $\gamma_gW$, summing $c$ with the square-root convolution in
\eqref{da:eq:weight-convolution}, and using
$\gamma_gW\ell_u=[(1-u)\zeta_u]^{-1}$, we obtain
\begin{equation}\label{fo:eq:pair-far-light}
 |\mathcal W_{2,ab}(u)|\prec
 R^3J^{3/2}(1-u)^{-1}\zeta_u^3w_{ab},\qquad d(a,b)>2h_u.
\end{equation}

For the martingale bracket, \eqref{gc:eq:bracket} joins the two pairs
into the six-cycle $(a,b,c,b',a',c)$, up to reversal and adjoints,
with the same coefficient $\gamma_gW$ and one free label $c$.
When both $d(a,b)$ and $d(a',b')$ exceed $2h_u$, the direct edges
$ab,a'b'$ are long, and at least two of the four $c$-edges are long:
$c$ cannot be within $h_u$ of both ends of either distant pair.
Choose the direct edges as the two Hilbert--Schmidt factors.
At most two remaining edges are short, so the normalized six-word
is bounded by $(CR)^2J^3\zeta_u^6\prod_{\rm edges}\sqrt w$.
Cauchy--Schwarz followed by weight convolution yields
\[
 \sum_c\sqrt{w_{ac}w_{bc}w_{a'c}w_{b'c}}
 \le\Bigl(\sum_cw_{ac}w_{bc}\Bigr)^{1/2}
      \Bigl(\sum_cw_{a'c}w_{b'c}\Bigr)^{1/2}
 \le C\ell_u\sqrt{w_{ab}w_{a'b'}}.
\]
Writing $\mathcal V_{ab,a'b'}(u)$ for the density of
$d[\mathcal M_{2,ab},\overline{\mathcal M_{2,a'b'}}]_u$, we obtain
\begin{equation}\label{fo:eq:pair-far-bracket}
 |\mathcal V_{ab,a'b'}(u)|\prec
 R^2J^3(1-u)^{-1}\zeta_u^5w_{ab}w_{a'b'}
 \quad\text{when both pairs are distant}.
\end{equation}

For near endpoint pairs, choose a summation halo
$H_u=\ell_u(\log W)^{C_h}$ with $C_h\theta>1$.
Inside the union of the endpoint halos there are at most $CH_u$
possible free labels. The rough triangle and six-loop estimates in
\eqref{fo:eq:step1} have sizes $R^2\zeta_u^2$ and $R^5\zeta_u^5$.
Thus the near drift and bracket budgets are respectively
\[
 \begin{aligned}
 \gamma_gW H_u(R^2\zeta_u)(R^2\zeta_u^2)
   &=W^{o(1)}R^4(1-u)^{-1}\zeta_u^2,\\
 \gamma_gW H_uR^5\zeta_u^5
   &=W^{o(1)}R^5(1-u)^{-1}\zeta_u^4.
 \end{aligned}
\]
Outside these halos, a $c$-edge is longer than $H_u$. Use its stopped
positive-pair estimate as one Hilbert--Schmidt factor and rough bounds
for the remaining factors. Its exponential is smaller than every
power of $W$, since $C_h\theta>1$. The square root of its floor costs
$W^{-F/2}$, multiplied only by the fixed polynomial barrier and rough
powers and by the normalized sum $\gamma_gWL$. Increase $F$ after
those fixed degrees to make this contribution smaller than any
prescribed $W^{-A}$. This uses the stopped mixed pair itself, without
the joint-tail conclusion of Step~3.
For $d(a,b)\le2h_u$ one has $w_{ab}^{-1}=W^{o(1)}$, because
$\theta\le1/2$. Choose $A>6$ in the preceding remainder estimate.
Since $\zeta_u\ge c/W$, the near drift and bracket scales are at
least $cW^{-2}$ and $cW^{-4}$, respectively, before the near weights.
Thus the $W^{-A}$ remainder is absorbed into these near source bounds
on their near endpoint support. This absorption occurs before
propagation; only the resulting weighted near bounds are propagated.
Restoring the near weights gives
\begin{equation}\label{fo:eq:pair-near-sources}
 \begin{gathered}
 |\mathcal W_{2,ab}(u)|\prec R^4(1-u)^{-1}\zeta_u^2w_{ab}
       \quad\text{for a near pair},\\
 |\mathcal V_{ab,a'b'}(u)|\prec
       R^5(1-u)^{-1}\zeta_u^4w_{ab}w_{a'b'}
       \quad\text{for two near pairs}.
 \end{gathered}
\end{equation}
To include near--distant cross terms after propagation, fix a deterministic
terminal pair $(a,b)$ at time $v\ge u$ and put
$\mathcal U_{u,v}=(U_{u,v}^{+-})^{\otimes2}$. Its coefficients
$\mathcal U_{u,v}(a,b;x,y)=U_{u,v}^{+-}(a,x)U_{u,v}^{+-}(b,y)$
are nonnegative. Covariance Cauchy--Schwarz gives
\[
 \begin{aligned}
 &\sum_{x,y,x',y'}\mathcal U_{u,v}(a,b;x,y)
       \mathcal U_{u,v}(a,b;x',y')|\mathcal V_{xy,x'y'}(u)|\\
 &\quad\le\left(\sum_{x,y}\mathcal U_{u,v}(a,b;x,y)
                   \sqrt{\mathcal V_{xy,xy}(u)}\right)^2.
 \end{aligned}
\]
Split the inner sum into near and distant pairs. Their diagonal
square-root bounds have coefficients
$R^{5/2}(1-u)^{-1/2}\zeta_u^2$ and
$RJ^{3/2}(1-u)^{-1/2}\zeta_u^{5/2}$, respectively, each times
$w_u(x-y)$. By \eqref{da:eq:propagated-weight}, the sum of this weight
against $\mathcal U_{u,v}(a,b;\cdot,\cdot)$ is at most
$CR^2w_v(a-b)$. Applying $(x+y)^2\le2x^2+2y^2$ therefore bounds
the full propagated bracket by
\begin{equation}\label{fo:eq:propagated-pair-bracket}
 \frac{Cw_v(a-b)^2}{1-u}
       \bigl[R^9\zeta_u^4+R^6J^3\zeta_u^5\bigr],
\end{equation}
with the reserved small loss. The stronger floor has already absorbed
all remainders into the near source bounds before this calculation.

For a deterministic terminal time $v$, use the pair propagator just
defined in the stopped Duhamel expression
\[
 \begin{aligned}
 \widehat D_2(v)={}&\mathcal U_{s,v}D_2(s)
 +\int_s^v\mathbf1_{u\le\tau_J}\mathcal U_{u,v}
       (\gamma_gW D_2(u)^2+\mathcal W_2(u))\,du\\
 &+\int_s^v\mathbf1_{u\le\tau_J}\mathcal U_{u,v}\,d\mathcal M_2(u).
 \end{aligned}
\]
The stopped indicators are predictable, and
$\widehat D_2(v)=D_2(v)$ for $v\le\tau_J$.
The propagated-weight estimate \eqref{da:eq:propagated-weight} costs
at most $R^2$ on two labels and $R^4$ on the four bracket labels.
The initial packet contributes $O_\prec(R^2)$ after division by
$\zeta_v^2w_v(a-b)$. For the quadratic and light drifts, the preceding
source bounds and \eqref{da:eq:scale-integral} give the normalized
budgets
\begin{equation}\label{fo:eq:pair-drift-budgets}
 \begin{aligned}
 R^2J^2\int_s^v\frac{\zeta_u^3}{\zeta_v^2(1-u)}\,du
     &\le R^2J^2\zeta_t\log R,\\
 R^5J^{3/2}\int_s^v\frac{\zeta_u^3}{\zeta_v^2(1-u)}\,du
     &\le R^5J^{3/2}\zeta_t\log R,\\
 R^6\int_s^v\frac{\zeta_u^2}{\zeta_v^2(1-u)}\,du
     &\le R^6\log R.
 \end{aligned}
\end{equation}
The three lines correspond to the quadratic error, distant light
insertions, and near light insertions. For the bracket, division by
$\zeta_v^4w_v(a-b)^2$ gives
\begin{equation}\label{fo:eq:pair-bracket-budgets}
 \begin{aligned}
 R^6J^3\int_s^v\frac{\zeta_u^5}{\zeta_v^4(1-u)}\,du
     &\le R^6J^3\zeta_t\log R,\\
 R^9\int_s^v\frac{\zeta_u^4}{\zeta_v^4(1-u)}\,du
     &\le R^9\log R.
 \end{aligned}
\end{equation}
BDG takes the square root of these budgets, and therefore bounds the
normalized martingale by
$O_\prec(R^3J^{3/2}\sqrt{\zeta_t}+R^{9/2})$.
Since $\log R=O(\log W)$, the logarithmic and halo factors fit inside
the reserved arbitrarily small domination losses. Combining the initial
term, \eqref{fo:eq:pair-drift-budgets}, and
\eqref{fo:eq:pair-bracket-budgets} proves the right-hand side of
\eqref{gc:eq:spatialstop} at each deterministic pair and terminal time.

We now pass to the stopped supremum in \eqref{gc:eq:spatialstop}.
Choose the deterministic sample pairs from \eqref{fo:eq:weighted-sampling-input} with weight
$w_u$, accuracy exponent $Q+4$ and failure exponent $D+3$.
After division by $\zeta_u^2\ge cW^{-2}$, its additive error is
$O(W^{-Q-2})$. Apply \eqref{fo:eq:time-mesh-union} with $j=0$,
accuracy exponent $Q+F+4$ and failure exponent $D+3$.
Since $\zeta_u^2w_u\ge cW^{-2-F}$, its raw increment error also
becomes $O(W^{-Q-2})$ after normalization. Each exceptional event
has probability at most $W^{-D-2}$ for sufficiently large $W$.
Direct differentiation gives
$|\partial_u\log(\zeta_u^2w_u)|\le C_FW^3\log W$; on a sufficiently
fine mesh, the normalizers at adjacent times therefore differ by at most a factor two.
For $u\le\tau_J$ use the mesh point immediately to its left, which
also precedes $\tau_J$. On the sampling and continuity events,
\begin{equation}\label{fo:eq:stopped-sampling}
 J_{\rm out}\le
 2C\max_{\substack{v\in\mathcal T,\ v\le\tau_J\\
                         (a,b)\in\mathcal S_{D_2}}}
 \frac{|\widehat D_2(v;a,b)|}{\zeta_v^2w_v(a-b)}+W^{-Q}.
\end{equation}
Here $C$ is the fixed norming constant. Once the deterministic
pairs and times number at most $W^{C_1}$ with fixed $C_1$, apply the
fixed-terminal Duhamel and BDG bound with failure exponent $D+C_1+2$.
The total failure probability is at most
\[
 W^{C_1}W^{-D-C_1-2}+2W^{-D-2}\le W^{-D-1}.
\]
This proves the supremum form of \eqref{gc:eq:spatialstop}. The
sampling and continuity events enter only this last comparison; they
are not inserted into a stochastic integrand. Only the actual field
before $\tau_J$ is interpolated.

Take $J_0=W^\epsilon R^{10}$ with $0<\epsilon<\min\{c/10,1\}$.
Since $R\le W^3$, this barrier satisfies $J_0\le W^{31}$ and lies
in the polynomial range fixed before the tail floor.
Dividing the right-hand side of \eqref{gc:eq:spatialstop} by $J_0$
and using $R^{20}\zeta_t\le W^{-c}$ gives
\[
 \begin{aligned}
 &W^{-\epsilon}R^{-4}+W^\epsilon R^{12}\zeta_t
   +W^{\epsilon/2}R^{10}\zeta_t
   +W^{\epsilon/2}R^8\sqrt{\zeta_t}\\
 &\qquad\le W^{-\epsilon}+W^{\epsilon-c}
       +W^{\epsilon/2-c}+W^{\epsilon/2-c/2}=o(1).
 \end{aligned}
\]
Every term is polynomially small, so the stochastic-domination loss
can be absorbed. The initial bound in \eqref{gc:eq:one-step-input}
starts below $J_0$ with overwhelming probability. Continuity therefore
excludes the first exit through $J_0$; the rough-bound exits already
have arbitrarily small probability by Step~1. This proves the first
line of \eqref{fo:eq:step2}.

For the Hilbert--Schmidt bound, the loop definition \eqref{fo:eq:loops}
and the mixed pair formula \eqref{fo:eq:pairs}, together with
\eqref{fo:eq:kernel-bounds}, give
\[
 \begin{aligned}
 W^{-2}\|(G_u)_{ab}\|_{\HS}^2
 &=\cK_2(u;b,a;+,-)+D_2(u;b,a;+,-)\\
 &\le C\zeta_u+|D_2(u;b,a;+,-)|
 \prec\zeta_u+R^{10}\zeta_u^2\prec\zeta_u,
 \end{aligned}
\]
since $R^{10}\zeta_u\le R^{20}\zeta_t\le W^{-c}$.
We have proved the two-loop conclusions
\begin{equation}\label{gc:eq:weakpair}
 |D_2(u;a,b;+,-)|\prec R^{10}\zeta_u^2w_u(a-b),\qquad
 W^{-2}\|(G_u)_{ab}\|_{\HS}^2\prec\zeta_u.
\end{equation}
To make the Haar consequence simultaneous, apply
\eqref{gc:eq:haar-supremum} with $\Psi_u=\sqrt{\zeta_u}$ on $[s,t]$.
The scale is bounded above and below by fixed powers of $W$, and its
logarithmic derivative is at most $CW^3$ directly from
\eqref{fo:eq:scale}. The fixed-pair Haar moments are thus
used only at the deterministic traceless-entry samples and time mesh,
including the $W^2$ orbital choices, before spatial norming and time
interpolation. The result is
\begin{equation}\label{fo:eq:sharp-traceless}
 \sup_{s\le u\le t}\max_{ai,bj}
 \frac{|(G_u)_{ai,bj}-\delta_{ij}\langle G_u\rangle_{ab}|}
      {\sqrt{\zeta_u}}\prec1.
\end{equation}

We finish the scalar improvement in two applications of fluctuation
averaging. Set $X_u=\max_{a,b}|\langle\wt G_u\rangle_{ab}|$ and
$Y_u=\max_a|\langle\wt G_u\rangle_{aa}|$.
The rough trace bound \eqref{fo:eq:rough-comparison} is simultaneous
in time by \eqref{fo:eq:step1-time-interpolation} with $j=0$, since
centering subtracts the time-independent matrix $M$.
Apply \eqref{gc:eq:averaging-input} with $K=CR$ and
$\chi=\max\{\sqrt{\zeta_u},W^{-1/2}\}\asymp\sqrt{\zeta_u}$.
Its fixed-pair moment proof, residual sampling and the same time
interpolation give $\sup_{u,a,b}|R_{ab}(u)|/(R^2\zeta_u)\prec1$.
The first scalar identity in \eqref{gc:eq:scalar} and the normalized
Dyson row in \eqref{fo:eq:propagation-dyson} now give
\[
 X_u\le\max_{a,b}|R_{ab}(u)|
       +CY_u\max_{a,b}|\langle G_u\rangle_{ab}|
 \prec R^3\zeta_u\le W^{-c},
\]
where $Y_u\prec R^2\zeta_u$ comes from Step~1. The bounded Dyson
entries then imply $\max_{a,b}|\langle G_u\rangle_{ab}|\prec1$.
A second application of the same original-law moment implication,
now with $K=C$, improves the residual to
$\sup_{u,a,b}|R_{ab}(u)|/\zeta_u\prec1$.
Applying $(I-uQ)^{-1}$ and its absolute row bound in
\eqref{fo:eq:propagation-dyson} to the second scalar identity gives
\[
 Y_u\le C\max_{a,b}|R_{ab}(u)|+CY_uX_u.
\]
On the common event where $X_u$ is uniformly a small power of $W$,
absorb the last term to obtain $Y_u\prec\zeta_u$; the first scalar
identity then gives $X_u\prec\zeta_u$.
Finally,
\[
 \|\wt G_u\|_{\max}
 \le\max_{ai,bj}|(G_u)_{ai,bj}-\delta_{ij}\langle G_u\rangle_{ab}|+X_u
 \prec\sqrt{\zeta_u}.
\]
All these statements have the time supremum inside probability, by the
same residual sample and time-mesh construction used in Step~1.
This proves every bound in \eqref{fo:eq:step2}.

\paragraph{Step 3: sharp bounds on the actual loops.}
The goal of this step is to obtain, for every fixed $k$, the sharp
actual-loop bound
\begin{equation}\label{fo:eq:step3}
 |\cL_k(u)|\prec\zeta_u^{k-1},
\end{equation}
and the intermediate error estimate
\begin{equation}\label{fo:eq:step3-error}
 |D_k(u)|\prec\zeta_s^{-1/2}\zeta_u^k.
\end{equation}
We also need the joint spatial tails, for every fixed $k\ge2$ and
$\epsilon,A>0$:
\begin{equation}\label{fo:eq:step5-tails}
 \max_{\boldsymbol\sigma\in\{+,-\}^k}
 \max_{\boldsymbol a:\,\diam(\boldsymbol a)>W^\epsilon\ell_u}
 \bigl(|\cL_k(u)|+|D_k(u)|\bigr)\prec W^{-A}.
\end{equation}
We prove the tails first, since they localize every subsequent free label.
Fix the finite collection of lengths required below, and let $N_0$ be
its largest error length. The raw brackets have length at most
$M_0=2N_0+2$. Choose $B_*$ so that every crude normalized block norm
and $R^{10}\zeta_u^2$ is bounded by $W^{B_*}$ on the Step~1--2 events.
Let $P_*$ exceed the fixed polynomial costs of spatial sums, row
restoration and propagator bounds for these lengths. These are analytic
error costs; sample and mesh cardinalities will enter only the failure
probabilities. Given a final tail accuracy $A$, put $A_*=A+P_*+2$ and
choose the stronger Step~2 floor by
\begin{equation}\label{gc:eq:floor}
 F>2A_*+(2M_0+1)B_*+10.
\end{equation}
A word of length $m\le M_0$ and diameter larger than $W^\epsilon\ell_u$
has an edge of length at least $W^\epsilon\ell_u/m$.
The positive mixed pair on that edge is the sum of the primitive pair
and its error. For large $W$, the edge exceeds $W^{\epsilon/2}\ell_u$,
so \eqref{fo:eq:primitive-tail} applies with exponent $\epsilon/2$.
Together with \eqref{gc:eq:weakpair}, it gives
\[
 W^{-2}\|(G_u)_{ab}\|_{\HS}^2
 \le W^{-Q}+CW^{B_*}
       \left(e^{-c(\|a-b\|_L/\ell_u)^\theta}+W^{-F}\right)
\]
on the fixed Step~1--2 event; $B_*$ includes the reserved domination
loss. Choose $Q>2A_*+2(M_0-1)B_*+4$. Use the square root of this
bound as one Hilbert--Schmidt factor in the word; the remaining
$m-1$ normalized factors have crude Hilbert--Schmidt bounds $W^{B_*}$.
The primitive remainder and the error floor contribute, respectively,
\[
 W^{-Q/2+(m-1)B_*}\le W^{-A_*-2},\qquad
 W^{B_*/2-F/2+(m-1)B_*}\le W^{-A_*-1}.
\]
The stretched-exponential term is $O_D(W^{-D})$ for every fixed $D>0$.
This proves the actual-loop tail. The primitive tail in
\eqref{fo:eq:primitive-tail} then gives the error tail, proving
\eqref{fo:eq:step5-tails}. Every cut word shares its free vertex and
every bracket word is a single joined cycle, so this same argument
covers their forcing tensors. All summed tails keep the covariance
factor, and restored rows keep their profile factors $\ell_u^{-1}$;
\eqref{fo:eq:spatial-sum-budget} bounds the remaining costs by
$\gamma_gL$ and powers of $L/\ell_u$. Thus the discarded terms remain
at most $W^{P_*}W^{-A_*-1}\le W^{-A}$, including at arbitrary coupling.

The primitive size input used in this step is
\eqref{fo:eq:primitivebound}, together with its length-two bound.
At bounded coupling its constant is uniform for $0<g\le C$.
At large coupling the constants are independent of $g,L$;
$L/(g^2\log W)\to\infty$ meets every required fixed volume threshold.
The comparison pairs remain in these domains by
Lemma~\ref{da:lem:comparison-parameters}.

Use the normalized random controls
\begin{equation}\label{gc:eq:normalized-controls}
 \begin{gathered}
 A_j\coloneqq\sup_{s\le v\le t}
     \max_{\boldsymbol a,\boldsymbol\sigma}
     \frac{|D_j(v;\boldsymbol a,\boldsymbol\sigma)|}{\zeta_v^j},\\
 B_j(v)\coloneqq\max_{\boldsymbol a,\boldsymbol\sigma}
     \frac{|\cL_j(v;\boldsymbol a,\boldsymbol\sigma)|}{\zeta_v^{j-1}},
 \qquad B_j\coloneqq\sup_{s\le v\le t}B_j(v).
 \end{gathered}
\end{equation}
In particular, $A_1\prec1$ by Step~2, and the primitive input gives
$B_j(u)\prec1+\zeta_u A_j$. Localize the one free block label of each
source to $W^\nu\ell_u$ using
\eqref{fo:eq:step5-tails}, with arbitrarily small fixed $\nu>0$.
Using the normalization in \eqref{fo:eq:spatial-sum-budget} and
$p+q=n+2$, the preceding controls give
\begin{equation}\label{gc:eq:sourcetable}
\begin{array}{c|c}
 \text{shorter linear cuts}&(1-u)^{-1}\zeta_u^n\max_{j<n}A_j\\
 \text{quadratic errors}&(1-u)^{-1}\zeta_u^{n+1}
                         \max_{p+q=n+2}A_pA_q\\
 \text{light insertions}&(1-u)^{-1}\zeta_u^n B_{n+1}\\
 \text{cross bracket}&(1-u)^{-1}\zeta_u^{2n}B_{2n+2}.
\end{array}
\end{equation}
Here $B_j$ is the time supremum in \eqref{gc:eq:normalized-controls};
fixed-order constants and the arbitrarily small localization loss are
understood. Primitive--error cuts in the first line have shorter error
length; the adjacent pair terms are already in the linear generator.
All terms on the right will be replaced by previously proved
deterministic envelopes before the stopped estimate is applied.
When a floor is summed, retain $\gamma_g$ with the sum; when a row is
restored, retain every factor $\ell_u^{-1}$ in the restoring profile.
The polynomial budget is again \eqref{fo:eq:spatial-sum-budget}.

We use the tensor propagator bounds \eqref{fo:eq:projected-propagation}
from Appendix~\ref{sec:gc-projected}, for zero rows or at least one
stable leg. A same-sign neighboring pair gives a stable leg by
\eqref{fo:eq:propagation-dyson}.
For an alternating word, set
$h_u(a_1)=\sum_{a_2,\ldots,a_n}D_n(u;\boldsymbol a)$ and
$\theta_u=(1-u)P(I-uP)^{-1}$, a probability profile with
$\|\theta_u\|_\infty\le C/\ell_u$ by \eqref{fo:eq:kernel-bounds}.
Subtract the signed row by
\[
 \cQ_uD_n(\boldsymbol a)=D_n(\boldsymbol a)
       -h_u(a_1)\prod_{j=2}^n\theta_u(a_1-a_j).
\]
For every fixed $Q>0$, the needed row and commutator bounds are
\begin{equation}\label{fo:eq:ward-projection}
 \begin{gathered}
 \|h_u\|_\infty\prec\ell_u^{n-1}\zeta_u^n(A_{n-1}+W^{-Q}),\qquad
 \|D_n-\cQ_uD_n\|_\infty\prec\zeta_u^n(A_{n-1}+W^{-Q}),\\
 \|C_u(h_u)\|_\infty\prec(1-u)^{-1}\zeta_u^n(A_{n-1}+W^{-Q}),
 \qquad \operatorname{row}C_u(h_u)=0.
 \end{gathered}
\end{equation}
We verify the signed row before taking its absolute value.
By conjugation, take an alternating word beginning with $+$ and ending
with $-$. Subtract the $\cL$ and $\cK$ identities in
\eqref{gc:eq:ward}, and then sum the remaining $n-2$ labels. With
$\boldsymbol\sigma^+=(+,\sigma_2,\ldots,\sigma_{n-1})$ and
$\boldsymbol\sigma^-=(-,\sigma_2,\ldots,\sigma_{n-1})$, this gives
\[
 h_u(a_1)=\frac1{2iW\eta_u}\sum_{a_2,\ldots,a_{n-1}}
 \left[D_{n-1}(u;a_1,\ldots,a_{n-1};\boldsymbol\sigma^+)
       -D_{n-1}(u;a_1,\ldots,a_{n-1};\boldsymbol\sigma^-)\right].
\]
The localized absolute sum has at most
$C_n(W^\nu\ell_u)^{n-2}$ terms. Since
$(W\eta_u)^{-1}=\ell_u\zeta_u/(\Im m)$ and $\Im m\asymp1$,
\begin{equation}\label{gc:eq:row}
 |h_u(a_1)|\prec
 (W\eta_u)^{-1}(W^\nu\ell_u)^{n-2}
        \zeta_u^{n-1}A_{n-1}+\ell_u^{n-1}\zeta_u^nW^{-Q}
 \prec\ell_u^{n-1}\zeta_u^n(A_{n-1}+W^{-Q}).
\end{equation}
The discarded row tails have only fixed polynomial costs after
normalization, including powers of $L/\ell_u$ and $\zeta_u^{-1}$.
Choose their tail accuracy after $Q$ and these fixed costs to leave the
normalized remainder $W^{-Q}$. The raw random control $A_{n-1}$ need
not be bounded below, so this remainder is retained here.
This calculation uses the signed Ward identity before absolute summation.
Since $\|\theta_u\|_\infty\le C/\ell_u$, restoring the row costs
$\ell_u^{n-1}\zeta_u^n(A_{n-1}+W^{-Q})(C/\ell_u)^{n-1}
=C_n\zeta_u^n(A_{n-1}+W^{-Q})$.

The projected differential explains the remaining commutator.
Write the error hierarchy \eqref{fo:eq:error-hierarchy}, with the
generator $\Theta$ from \eqref{fo:eq:generator}, as
$dD_n=\Theta_{u,\boldsymbol\sigma}D_n\,du+F_n\,du+d\mathcal M_n$,
where $F_n$ contains the remaining cut and light sources.
For the alternating word every external-leg generator is
$P(I-uP)^{-1}$. Direct differentiation of the probability profile gives
\[
 \partial_u\theta_u=P(I-uP)^{-1}\theta_u-(1-u)^{-1}\theta_u.
\]
In differentiating $\cQ_uD_n$, the $j\ge2$ external-leg terms cancel
the corresponding profile derivatives: summing each such generator
contributes $(1-u)^{-1}$ to the signed row. The difference in the first
leg is exactly $C_u(h_u)$ in \eqref{gc:eq:commutator}. Therefore
\begin{equation}\label{fo:eq:projected-error-equation}
 d(\cQ_uD_n)=\Theta_{u,\boldsymbol\sigma}(\cQ_uD_n)\,du
       +[\cQ_uF_n+C_u(h_u)]\,du+\cQ_u\,d\mathcal M_n.
\end{equation}
The bound \eqref{da:eq:commutator-bound} with \eqref{gc:eq:row}
gives $(1-u)^{-1}\zeta_u^n(A_{n-1}+W^{-Q})$ for the commutator,
which proves the third line of \eqref{fo:eq:ward-projection}, including
its zero row.
The projected bracket is
$(\cQ_u\otimes\cQ_u)\mathcal V_n$, where $\mathcal V_n$ is the
unprojected bracket tensor. It has zero rows in both sets of labels.
For a source tensor $A$ supported on relative diameter $H\ell_u$,
$\|\theta_u\|_\infty\le C/\ell_u$ gives
\[
 \left\|(\operatorname{row}A)\prod_{j=2}^n\theta_u(a_1-a_j)\right\|_\infty
 \le C_n(H\ell_u)^{n-1}\|A\|_\infty\ell_u^{-(n-1)}
 =C_nH^{n-1}\|A\|_\infty.
\]
Apply this bound in each coordinate group for the bracket.
With $H=W^\nu$, projection thus retains the source-table scales up to
the reserved localization loss. Each restoring profile reconnects its
labels to its anchor; its exponential tails follow from
\eqref{fo:eq:kernel-bounds} and the definition of $\theta_u$.
The discarded polynomial floors are controlled next.
Apply \eqref{fo:eq:projected-propagation} to these zero-row tensors;
nonalternating words use its stable-leg bound. This proves the
projection mechanism needed for the recurrence.

For the truncation error use the tail argument in the proof of
\eqref{fo:eq:projected-propagation}, with discarded tensor $T$,
$\|T\|_\infty\le W^{-F}$. For $u\le v\le t$, restoration against $\theta_v$ costs at most
$C_n(L/\ell_v)^{n-1}W^{-F}\le C_n(L/\ell_u)^{n-1}W^{-F}$.
Propagate the discarded and restoration tensors separately; the latter
has scale $\ell_v$, not support in the truncated region.
Since $R\le W^3$ and \eqref{fo:eq:spatial-sum-budget} bounds the remaining
costs polynomially, choose $F$ after these costs and the required accuracy.
This controls the row-restoration and covariance-weighted cut remainders
without changing the main source normalizations.

The estimate supplied by the projected hierarchy is
\begin{equation}\label{gc:eq:recurrence}
 A_n\prec 1+\max_{j<n}A_j+
 \zeta_t\max_{2\le j\le n}A_jA_{n-j+2}
       +B_{n+1}+\sqrt{B_{2n+2}}.
\end{equation}
We use this estimate as a bound-improvement implication: replace the
controls on the right by deterministic polynomial bounds already proved
at the preceding induction stages, and stop predictably at those bounds.
Choose all these deterministic envelopes at least $1$; they absorb
the arbitrary normalized $W^{-Q}$ remainders in the row and commutator
bounds. The stopped estimate then gives the displayed deterministic
bound for $A_n$ in stochastic domination. To derive it, apply the deterministic-terminal Duhamel
formula to \eqref{fo:eq:projected-error-equation}, or to the stable-leg
hierarchy for a nonalternating word. A drift source of size
$(1-u)^{-1}\zeta_u^n C_u$ contributes at most
\[
 \zeta_v^{-n}\int_s^v(\zeta_v/\zeta_u)^n
            (1-u)^{-1}\zeta_u^n C_u\,du
 \le(\log R)\sup_{s\le u\le v}C_u.
\]
The stable-leg estimate is smaller since $\zeta_v\ge\zeta_u$.
The quadratic source has the additional factor $\zeta_u\le\zeta_t$.
For the bracket, use order $2n$ in
\eqref{fo:eq:projected-propagation}; its normalized integral is bounded
by $(\log R)B_{2n+2}$. BDG consequently gives
$(\log R)^{1/2}\sqrt{B_{2n+2}}$ for the normalized martingale.
At the initial time, \eqref{gc:eq:one-step-input} bounds both
$D_n(s)$ by $O_\prec(\zeta_s^n)$ and $D_{n-1}(s)$ by
$O_\prec(\zeta_s^{n-1})$. The signed-row calculation therefore gives
$\|\cQ_sD_n(s)\|_\infty\prec\zeta_s^n$; the initial tails have the
same floor budget. Propagation makes this a normalized contribution $1$.
The terminal row restoration contributes $A_{n-1}$.
These are exactly the five terms in
\eqref{gc:eq:recurrence}; the logarithms and localization powers are
absorbed in the reserved small losses.

At each induction node, stop at the already proved
shorter-length and preceding-round envelopes, together with a provisional
barrier for the current target. Apply the projected Duhamel formula at
deterministic sample tuples and terminal times. The passage to the
supremum is \eqref{fo:eq:stopped-sampling} with $D_2$ replaced by $D_n$
and normalization $\zeta_v^n$; the proof uses the actual error only
before the stop. Thus no estimate on the target at the current node
is used to justify its own improvement.

We now establish the same-sign two-loop base
\begin{equation}\label{fo:eq:same-pair}
 \sup_{s\le u\le t}\max_{\substack{a,b\in\Z_L\\\sigma\in\{+,-\}}}
       \zeta_u^{-2}|D_2(u;a,b;\sigma,\sigma)|\prec R^{10}.
\end{equation}
We adapt the Duhamel and martingale estimates from
\cite[proof of Lemma~5.11]{YY} and close the two-loop estimate
by the stopping argument below.
Write $D_2^{++}(u)\coloneqq(D_2(u;a,b;+,+))_{a,b\in\Z_L}$.
The exact equation is
\begin{equation}\label{gc:eq:same-pair-equation}
 \begin{split}
 dD_2^{++}={}&\bigl[Q(I-uQ)^{-1}D_2^{++}
                   +D_2^{++}Q(I-uQ)^{-1}\bigr]du\\
 &+\gamma W(D_2^{++})^2du+\mathcal W_2^{++}du
                         +d\mathcal M_2^{++}.
 \end{split}
\end{equation}
Here $Q_{ab}=\gamma(M_B)_{ab}^2$ as in \eqref{fo:eq:PQ}.
Its propagator has bounded absolute row and column sums by
\eqref{fo:eq:propagation-dyson} and symmetry. The initial packet gives
$\max_{a,b}|D_2^{++}(s;a,b)|\prec\zeta_s^2$.
The already proved local law gives $A_1\prec1$, whereas the rough bounds
\eqref{fo:eq:step1} give $B_3\prec R^2$ and $B_6\prec R^5$.
These are bounds on actual loops and do not require any improvement
of a centered higher loop.

Fix $0<\epsilon<c/4$ and put $J_0\coloneqq W^\epsilon R^{10}$.
Choose a total loss $0<\nu<\epsilon/4$, reserving smaller losses
for each estimate.
Stop at the first time
$\max_{a,b}|D_2^{++}(u;a,b)|/\zeta_u^2$ reaches $J_0$, also stopping
on the previously established rough, local-law and tail envelopes.
All tests use only the current and past flow. Before this stop, the
tails obtained above from the mixed positive pair localize the free
block label in $(D_2^{++})^2$ to $W^{\nu/10}\ell_u$ sites, up to an
arbitrarily small polynomial remainder. Thus
\[
 \gamma W\max_{a,b}|[(D_2^{++})^2]_{ab}|
 \le C W^{\nu/10}J_0^2(1-u)^{-1}\zeta_u^3+W^{-A},
\]
where the floor is chosen after $A$. The light term and the joined
six-loop bracket in \eqref{gc:eq:sourcetable} are bounded, respectively,
by $O_\prec(R^2(1-u)^{-1}\zeta_u^2)$ and
$O_\prec(R^5(1-u)^{-1}\zeta_u^4)$.

Apply Duhamel's formula with predictable stopped sources and the
deterministic propagator $U_{u,t'}^{++}$ at each deterministic terminal
time $t'$. Its bounded rows and BDG give the stopped version of
\eqref{gc:eq:recurrence} at length two:
\begin{equation}\label{gc:eq:same-pair-stop}
 J_{\rm out}\prec 1+R^2+R^{5/2}+\zeta_t J_0^2.
\end{equation}
Here $J_{\rm out}$ is the maximal normalized error up to the stop.
Time integration costs at most $1+\log R$, and the localization and
logarithmic losses are absorbed into an arbitrarily small domination
loss. The ordinary-tuple sampling and time-increment argument in the proof
of \eqref{fo:eq:stopped-sampling} gives the space-time supremum, also at
arbitrary coupling; no random terminal propagator is inserted into a
stochastic integral.
With the reserved total loss $W^\nu$, division by $J_0$ gives
\[
 J_0^{-1} J_{\rm out}
 \le C W^{\nu-\epsilon}R^{-15/2}
       +C W^{\nu+\epsilon}R^{10}\zeta_t=o(1)
\]
with arbitrarily high probability, by $R^{20}\zeta_t\le W^{-c}$.
The initial error is below $J_0$ with the same probability.
Continuity therefore excludes the first exit. Taking $\epsilon$
arbitrarily small proves the $(+,+)$ assertion in
\eqref{fo:eq:same-pair}; conjugation gives $(-,-)$.
No spatial profile for the same-sign error is needed here.
All constants used in this argument are the common propagation
constants, so it applies to both coupling cases of
Theorem~\ref{gc:thm:common}.

A direct induction on $n$ would not close, since the bracket contains
an actual loop of length $2n+2$. Put
$m_1=2\lfloor(n+1)/2\rfloor$ and
$m_2=2\lceil(n+1)/2\rceil$. Trace Cauchy--Schwarz gives
\begin{equation}\label{gc:eq:chain}
 B_{2n+2}(v)\le\zeta_v^{-1}B_{m_1}(v)B_{m_2}(v).
\end{equation}
This is the four-chain Cauchy--Schwarz calculation in
\cite[(A.74)--(A.77)]{TYY}; the present normalization contributes
$\zeta_v^{m_1+m_2-2-(2n+1)}=\zeta_v^{-1}$.

For the fixed interval $[s,t]$, use double induction in the loop
length $n$ and the improvement round $j\ge0$ to prove
\begin{equation}\label{fo:eq:finite-improvement}
 H_{n,j}\coloneqq\zeta_s^{-1/2}+(CR)^{n-1}\zeta_s^{j/4-1},\qquad
 A_n\prec H_{n,j}.
\end{equation}
Here $j$ counts improvements of the same random control, not time steps.
We verify every term of \eqref{gc:eq:recurrence} using only the already
proved controls assigned to the current induction node.
The rough bound \eqref{fo:eq:step1} and the primitive estimate give
$A_n\prec H_{n,0}$. The local law supplies $A_1\prec1$.
The mixed and same-sign two-loop bounds give, at every round,
\[
 A_2\prec R^{10}\le\zeta_s^{-1/2}\le H_{2,j},
\]
since $R^{20}\zeta_s\le R^{20}\zeta_t\le W^{-c}$.
For $j\ge1$ proceed in increasing $n\ge3$. Use shorter lengths
already proved in round $j$ for linear terms and interior quadratic
products. Use round $j-1$ for $B_{n+1}$, the target factor $A_n$ in
$A_2A_n$, and the two loop controls from \eqref{gc:eq:chain}.
Their lengths satisfy $m_1,m_2\le n+2$.

For a target length $k\ge3$, start with rough bounds through
$k+2(k+4)=3k+8$. At round $j=1,\ldots,k+4$, improve in increasing
order through length $k+2(k+4-j)$. The preceding round supplies every
length through $n+2$ required at a current node, and the last round
reaches $k$. Thus this construction has finitely many nodes and never
uses an unknown same-round target on the right.

The scale definition \eqref{fo:eq:scale} and the time-step condition
\eqref{gc:eq:one-step} give
\[
 c/W\le\zeta_s\le\zeta_u\le\zeta_t,\qquad
 \zeta_t/\zeta_s\le R,\qquad R^{20}\zeta_t\le W^{-c}.
\]
In particular $R\le\zeta_s^{-1/20}$ and $\zeta_s\le W^{-c}$.
The following strict power margins absorb all fixed-order constants
and the separately reserved losses at the finitely many nodes.

Use the scalar induction in \cite[Lemma~7.18, Appendix~A.8]{TYY}
with the assigned earlier envelopes. For the present time-ratio bound,
substitution of \eqref{fo:eq:finite-improvement} and \eqref{gc:eq:chain}
gives, uniformly for $s\le u\le t$,
\[
 \begin{aligned}
 B_{n+1}&\prec C\zeta_s^{1/2}H_{n,j},\\
 \sqrt{B_{2n+2}(u)}&\prec\zeta_s^{-1/2}
                       +C\zeta_s^{3/40}H_{n,j},\\
 \zeta_tR^{10}H_{n,j-1}&\le C\zeta_s^{1/5}H_{n,j},\\
 \zeta_tH_{p,j}H_{q,j}&\le C\zeta_s^{3/20}H_{n,j}
             \quad(p+q=n+2,\ 3\le p,q<n).
 \end{aligned}
\]
These control the light, bracket, boundary quadratic and interior
quadratic sources, respectively. Shorter linear terms are at most
$H_{n,j}$, so every source in \eqref{gc:eq:recurrence} has this bound.
Apply its adapted stopped implication at this node and continue through
the finite cone. Choosing smaller losses at its finitely many nodes
preserves the requested final domination loss.

At the final round $j=k+4$ for the target $k$,
\[
 (CR)^{k-1}\zeta_s^{(k+4)/4-1}
 \le C_k\zeta_s^{(4k+1)/20}.
\]
Thus $A_k\prec\zeta_s^{-1/2}$, proving
\eqref{fo:eq:step3-error}. This argument applies to each required fixed
length and the finite union of their cones.

Finally, $\zeta_u/\sqrt{\zeta_s}\le\sqrt{R\zeta_t}\le W^{-c/2}$
by the scale ratio and $R^{20}\zeta_t\le W^{-c}$.
Adding the primitive bound \eqref{fo:eq:primitivebound} and the
length-two estimate in Lemma~\ref{fo:lem:deterministic} therefore yields
\[
 |\cL_k(u)|\le|\cK_k(u)|+|D_k(u)|
 \prec\zeta_u^{k-1}\left(1+\frac{\zeta_u}{\sqrt{\zeta_s}}\right)
 \prec\zeta_u^{k-1},
\]
proving \eqref{fo:eq:step3}; cf.\
\cite[Section~7.4, (5.15)]{TYY} and \cite[Section~5.6, (2.77)]{YY}.

\paragraph{Step 4: sharp bounds on the loop errors.}
The goal of this step is to remove the factor $\zeta_s^{-1/2}$ and obtain
\begin{equation}\label{fo:eq:step4}
 |D_k(u)|\prec\zeta_u^k
\end{equation}
for every fixed $k$. Reinsert the sharp actual loops from Step~3
and the sharp one-loop error from Step~2 into the hierarchy.
The source bounds and projected propagation established in
Step~3 apply to these improved inputs. In particular, the four source
bounds in \eqref{gc:eq:sourcetable} give drift scale $(1-u)^{-1}\zeta_u^k$,
bracket scale $(1-u)^{-1}\zeta_u^{2k}$, and an extra factor
$\zeta_u$ for quadratic errors, also at arbitrarily large coupling.

With $B_j\prec1$ at every required length, the stopped implication
\eqref{gc:eq:recurrence} reduces to
\[
 A_n\prec1+\max_{j<n}A_j
       +\zeta_t\max_{2\le j\le n}A_jA_{n-j+2}.
\]
Use the length induction of \cite[Section~7.5, (7.66)--(7.67)]{TYY}
on this recurrence. With $A_1\prec1$, the successive input envelopes
$A_2\prec\zeta_s^{-1/2}$ from Step~3 and then $A_2\prec R$ give
\begin{equation}\label{fo:eq:two-loop-sharpening}
 A_2\prec1+\frac{\zeta_t}{\zeta_s}\lesssim R,\qquad
 A_2\prec1+\zeta_t R^2\lesssim1.
\end{equation}
Here $\zeta_t/\zeta_s\le R$ and $\zeta_tR^2\le W^{-c}$.
Each application stops at its already proved input envelope and then
uses the deterministic-sample and time-interpolation argument following
\eqref{gc:eq:recurrence}; the new output is not used in its own integrand.

For $n\ge3$, use shorter current-stage bounds $A_m\prec1$ in the
linear and interior quadratic terms, and the Step~3 envelope
$A_n\prec\zeta_s^{-1/2}$ in the boundary products $A_2A_n$. Then
\[
 A_n\prec1+\zeta_t+\zeta_t\zeta_s^{-1/2}\prec1,
\]
since $\zeta_t/\sqrt{\zeta_s}\le\sqrt{R\zeta_t}\le W^{-c/2}$.
All input controls precede this induction node.
The same joint-tail budget \eqref{gc:eq:floor} covers the localized
sums and stopped brackets. The finitely many loss and failure reserves
therefore give the spatial and time supremum in \eqref{fo:eq:step4}
for every fixed length, uniformly in the deterministic parameters.

\paragraph{Step 5: restoring the spatial profile and the next initial data.}
The goal of this step is to remove the spatial loss $R^{10}$:
\begin{equation}\label{fo:eq:step5}
 |D_2(u;a,b;+,-)|\prec\zeta_u^2
 \left(e^{-c(\|a-b\|_L/\ell_u)^{\theta/2}}+W^{-F}\right).
\end{equation}
Together with the joint tails \eqref{fo:eq:step5-tails} and the sharp
bounds from Steps~2--4, this restores every estimate in
\eqref{gc:eq:one-step-output}.

Put $x=\|a-b\|_L/\ell_u$ and apply Step~2 with a stronger floor
exponent $F'>F+32$. Because $R\le W^3$, its floor satisfies
$R^{10}W^{-F'}\le W^{-F-2}$. Taking the minimum with Step~4 gives
\[
 \frac{|D_2(u;a,b;+,-)|}{\zeta_u^2}
 \prec\min\{1,R^{10}e^{-c_{\rm p}x^\theta}\}+W^{-F-2}.
\]
The inequality $\min\{1,a+b\}\le\min\{1,a\}+b$ separates the floor.
Apply \eqref{da:eq:profile-minimum} with
$x=\|a-b\|_L/\ell_u$ and $1\le R\le W^3$.
Its prefactor $e^{C'\sqrt{\log W}}=W^{o(1)}$ is absorbed
in the reserved domination loss. Together with the
$W^{-F-2}$ remainder, this proves \eqref{fo:eq:step5}.
The constants are uniform in $a,b,u$ and the allowed parameters.
A fixed finite number of steps leaves a fixed positive profile exponent.

For the operator input at the next time, the sharp actual loops include
\begin{equation}\label{gc:eq:regenerate}
 \Tr[(G_t E_bG_t^*E_a)^p]
 =W^{-2p}\Tr[((G_t)_{ab}(G_t)_{ab}^*)^p]
 \prec\zeta_t^{2p-1}.
\end{equation}
Consequently,
\begin{equation}\label{fo:eq:operator-regeneration}
 W^{-1}\max_{a,b}\|(G_t)_{ab}\|_{\op}\prec\zeta_t.
\end{equation}
Repeat the Schatten-to-operator argument following
\eqref{fo:eq:starting-positive-loop}, now using
\eqref{gc:eq:regenerate} at time $t$.
For each $\epsilon,D>0$, choose a fixed integer $p>1/\epsilon$
before the loop-length cone, use loss $\epsilon/2$, and invoke
$\zeta_t\ge c/W$. The spatial maximum is inherited from the loop
estimate, so no $L^2$ union is needed. The same positive words supply
the fixed Schatten moments, including for off-diagonal blocks.
Together with the averaged law
and the comparison \eqref{fo:eq:backward-matching}, all inputs for
the next step are now available. All these estimates already include
the spatial and time suprema by \eqref{fo:eq:time-mesh-union} and
\eqref{fo:eq:stopped-sampling}. For each prescribed $\epsilon,D>0$,
reserve smaller losses at the finitely many uses above and failure
$W^{-D-2}$ for each step; their union has probability at most
$5W^{-D-2}\le W^{-D}$. Constants are uniform in deterministic $(g,E)$,
so their supremum remains outside the probability as required in
\eqref{gc:eq:one-step-output}. This proves part (i).
\end{proof}

\subsection{Finite time iteration}\label{sec:gc-schedule}

\begin{proof}[\textbf{Proof of Theorem~\ref{gc:thm:common}(ii)}]
We prove \eqref{gc:eq:iteration-output} with all
required finite orders and probability estimates fixed before the forward
iteration. Write $g_0$ for the coupling in the desired conclusion.
For the fixed $c_0>0$ and target $1/2\le t_*\le1-W^{-3}$ with
$\zeta_{t_*}\le W^{-c_0}$, set $J=\lceil240c_0^{-1}\rceil+1$ and
\begin{equation}\label{gc:eq:finite-schedule}
 s_j\coloneqq1-\max\{1-t_*,\tfrac12W^{-jc_0/80}\},\qquad 0\le j\le J.
\end{equation}
Since $Jc_0/80>3$, these times run from $s_0=1/2$ to $s_J=t_*$;
omit repeated endpoints. The definition and monotonicity of $\zeta$ give
$R_j\coloneqq(1-s_{j-1})/(1-s_j)\le W^{c_0/80}$ and
$R_j^{20}\zeta_{s_j}\le W^{-3c_0/4}$, while
$1-s_j\ge1-t_*\ge W^{-3}$. Thus part (i) applies with $c=c_0/2$;
$J$ is fixed independently of loop orders, losses, floors and failure exponents.

We must supply the starting estimates at every deterministic comparison
pair used within a step. Lemma~\ref{da:lem:comparison-parameters}(b)
applies to each backward comparison chain of length at most $J$.
At bounded coupling it gives $0<c_Jg_0\le g\le C$, bounded energies,
and a fixed positive lower bound on $\Im m$, preserving the input
range and the original $L\gg\log W$ assumption. At large coupling,
\eqref{fo:eq:large-input-domain} gives
$g_0/\sqrt2\le g\le g_0$ and $|E|\le2-\kappa/2$.
Consequently $L\gg g_0^2\log W$ implies $L\gg g^2\log W$, while
$L\le Cg_0^2W^2$ implies $L\le2Cg^2W^2$.
The preserved domain \eqref{fo:eq:large-time-domain} and its equality
$s\gamma_{g'}=u\gamma_g$ give precisely the earlier-time admissibility
needed in Step~1. The same lemma preserves $1-u\ge W^{-3}$ and the
endpoint scale bound $\zeta_u(g)\le\zeta_{t_*}(g_0)\le W^{-c_0}$.
For an earlier comparison pair required at time $s$, every preceding
time $v\le s$ satisfies
\[
 v\gamma_{g'}\le s\gamma_{g'}=u\gamma_g\le g_0^{-2}.
\]
Propagate that pair only through preceding intervals ending at $s$,
using the same schedule with deterministic admissibility truncation
when needed; propagation up to $t_*$ at the smaller $g'$ is not required.
These are deterministic parameter bounds; uniformity in them is outside
the probability and requires no union over an energy or coupling continuum.

Fix the desired final loop lengths, domination loss, spatial-tail
accuracy and failure exponent. For every operator loss $W^\delta$
needed at the next step, first select the fixed integer $p>1/\delta$
in \eqref{gc:eq:regenerate} and include its positive loop of length $2p$.
For an actual-loop output of length $n$, the finite cone in Step~3
requires rough lengths through $n+2(n+4)=3n+8$.
For an error output of length $k$, Step~4 needs actual loops through
$2k+2$, hence rough lengths through $6k+14$ suffice.
Take the finite union of these requirements and reserve smaller losses
at its nodes so that products and square roots leave the prescribed
output loss. The starting data at this step must supply the resulting
finite list at the original and the deterministic comparison pairs.

Now fix the analytic tail accuracies for this finite list. If $N_0$
is the largest required error length, its raw joined bracket has length
at most $M_0=2N_0+2$. Such a bracket uses the long-edge tail argument
of Step~3 and crude bounds on the remaining factors, not sharp errors
of length $M_0$. Fix their crude exponent $B_*$ and the polynomial
spatial-sum, row-restoration and propagator cost $P_*$. Given tail
accuracy $A$, choose $A_*=A+P_*+2$, the primitive remainder accuracy
$Q>2A_*+2(M_0-1)B_*+4$ from the long-edge calculation, and then the floor in \eqref{gc:eq:floor}.
For an output profile floor $W^{-F}$, choose the stronger input floor
$F'>F+32$ used in Step~5 as well. Increase the finitely many input
floor exponents to meet both requirements.
The discarded analytic errors are at most
$W^{P_*}W^{-A_*-1}\le W^{-A}$.

Only after these floor choices, fix each global norming and
interpolation event's failure probability at $W^{-D-2}$ for the desired
final exponent $D$. Invoke \eqref{fo:eq:weighted-sampling-input} at
this fixed failure exponent and the required raw accuracy, and choose
polynomial time meshes using \eqref{fo:eq:coarse-time-input}. For example, a weighted $n$-loop
normalization costs at most $C_nW^{F+n}$ because $\zeta_u\ge c/W$;
a raw error $W^{-Q-F-n-2}$ is therefore sufficient for normalized
accuracy $W^{-Q-1}$. Vertical comparisons use the closed
\eqref{fo:eq:comparison-height-interpolation}; its coarse envelope
exponent is fixed before its height mesh and moments. The Duhamel
integrals are exact, and interpolation compares a point with its
nearest or left mesh point and takes a supremum. No mesh cardinality
multiplies a deterministic tail error. Mesh sizes enter only the
probability union.

Let the combined spatial, orbital and time or height meshes at all
nodes have at most $W^{Q_{\rm mesh}}$ points. Their cardinality may
depend on the already fixed global failure exponent $D+2$.
The union of those global exceptional events has probability
$O(W^{-D-2})$. At the resulting fixed points, require probability
$W^{-D-Q_{\rm mesh}-2}$ only for the pointwise Gaussian, Haar and
stopped-martingale estimates. Their union costs another $O(W^{-D-2})$,
so the combined failure is at most $W^{-D}$ for large $W$.
Do not request the global sampling or interpolation input again with
this larger pointwise failure exponent. Choose the BDG and Markov
moment orders for these pointwise estimates last: a threshold gap
$W^\rho$ is handled by any fixed order $q$ with
$q\rho>D+Q_{\rm mesh}+2$.
Increasing this moment order does not enlarge the sharp-loop cone.
For a target of length $n$, the joined bracket word has the fixed
length $2n+2$; its deterministic envelope is imposed by an adapted
stop before BDG is applied. Higher moments are then taken of that
stopped martingale with its already bounded bracket.
Only the positive-loop order $p$ selected for the operator loss changes
the required loop lengths, and it was fixed before constructing the cone.
Brownian interpolation events are used after these adapted estimates;
no moment estimate is conditioned on a future interpolation event.

Assign this finite list of starting lengths, floors, losses and failure
accuracies to the preceding time step, and repeat backwards at most
$J$ times. Each individual step has finitely many improvement nodes,
and each comparison remains in the domains already verified. This
backward construction never requires one finite list to be closed
under repeated enlargement at the same time level. Its constants and
orders do not alter the deterministic schedule. Finally choose the
accuracy-dependent deterministic constants; the divergent volume
assumptions accommodate every one for sufficiently large $W$.
The bounded initial estimate \eqref{gc:eq:iteration-input} in
Lemma~\ref{fo:lem:bounded-initial}, or its
large-coupling version in Lemma~\ref{fo:prop:inputs}, supplies all these
fixed-order inputs at $s_0=1/2$, uniformly at the buffered pairs.

Apply part (i) forward over the chosen distinct intervals. The preceding
step provides the next starting packet and the operator/Schatten inputs
proved in Step~5, also at every comparison pair. The profile exponent
starts at $1/2$ and is halved at most $J$ times, so one may take
$\theta_*=2^{-J-1}>0$, with the minimum of the finitely many positive
profile constants. The time suprema are already inside the probability
on each interval; their finite union gives the full interval
$[1/2,t_*]$. For the original coupling $g=g_0$, the domain
$u\gamma_{g_0}\le g_0^{-2}$ imposes no additional restriction.
This proves \eqref{gc:eq:iteration-output} with its stated uniformity.
\end{proof}

\section{Orbital spreading and the rough density law}
\label{sec:elementary}
The rough density law at the end of this appendix is used in both parts
of Theorem~\ref{thm:main}.

\subsection{Common unitary conjugation}\label{subsec:pf_thm31}

\begin{proof}[\textbf{Proof of \eqref{eq:orbital-length} and \eqref{eq:haar-count}}]
Let $U$ be independent Haar unitary. The GUE law and scalar hopping give
$(I_L\otimes U)H(I_L\otimes U^*)\stackrel{\rm law}=H$.
Conditional on $H$, rotate all eigenvectors by $I_L\otimes U$;
each normalized nonzero block is uniform on the unit sphere of $\C^W$.
For a uniform unit vector $v\in\C^W$, the coordinate tail
$\Prob(|v_j|^2>s)=(1-s)_+^{W-1}$ for $s\ge0$, and a union bound
over the $N^2$ eigenvector-coordinate pairs, give with probability
at least $1-N^{-D}$
\begin{equation}\label{eq:orbital}
 |\psi_k(i)|^2\le C_DW^{-1}\log N\,\|\psi_{k,a}\|_2^2,
 \qquad i\in\mathcal I_a,\quad a\in\Z_L,\quad 1\le k\le N.
\end{equation}
The half-mass definition immediately gives \eqref{eq:orbital-length}.

For \eqref{eq:haar-count}, fix a normalized eigenvector and let
$r=\lceil c_0W\rceil<W$ with $c_0$ sufficiently small. A coordinate
interval of length $r$ meets at most two blocks. If it carries half
the mass, it contains mass at least $1/4$ from one of the at most
four blocks whose total mass is at least $1/4$.
For a fixed set $I$ of $r$ coordinates, the
$\operatorname{Beta}(r,W-r)$ tail gives
\[
 \Prob\left(\sum_{j\in I}|v_j|^2\ge1/4\right)\le e^{-cW}.
\]
A union bound over at most $W$ cyclic consecutive coordinate sets in each
of the four blocks bounds the conditional probability of
$\ell(\psi_k)<c_0W$ by $4We^{-cW}$. Averaging over the $WL$
eigenvectors and then over $H$ proves \eqref{eq:haar-count}.
\end{proof}

\subsection{The rough density law}
\label{sec:rough-density-proof}

\begin{proof}[\textbf{Proof of \eqref{eq:rough-counting}}]
Write $\gamma=\gamma_g$, so $H=A_g\otimes I_W+\sqrt\gamma V$,
with $0<\gamma\le1$ and $\|A_g\|\le2$. For $z=E+i\eta$,
$0<\eta\le1$, translation invariance gives
$\bar m:=\E\langle G(z)\rangle_{aa}=N^{-1}\E\Tr G(z)$.
Gaussian Poincar\'e and integration by parts yield
\begin{equation}\label{rd:eq:poincare}
 \E|\langle G\rangle_{aa}-\bar m|^2\le C\gamma W^{-2}\eta^{-4},
\end{equation}
and
\begin{equation}\label{rd:eq:mean}
 (A_g\otimes I_W-z-\gamma\bar m)\E G
 =I+\gamma\E\bigl[\diag((\langle G\rangle_{aa}-\bar m)I_W)G\bigr].
\end{equation}
Cauchy--Schwarz bounds the error on the right in operator norm by
$CW^{-1}\eta^{-3}$, independently of $L,g$.
After inversion and normalized trace, the Fourier-sum estimate in
Lemma~\ref{da:lem:rough-stability} gives
\begin{equation}\label{rd:eq:approximate}
 \bar m=f_g(z+\gamma\bar m;\gamma)
       +O\left(W^{-1}\eta^{-4}+L^{-1}\eta^{-2}\right).
\end{equation}
For fixed $\eta$ and sufficiently large $W,L$, the error is at most
$\eta/2$, so the same lemma gives
\begin{equation}\label{rd:eq:bias}
 |\bar m-m_g(z)|\le C\left(W^{-1}\eta^{-6}+L^{-1}\eta^{-4}\right).
\end{equation}
Gaussian concentration gives
\begin{equation}\label{eq:rough-density}
 \Prob\left(\left|N^{-1}\Tr G(E+i\eta)-\bar m\right|>u\right)
 \le4\exp(-cW^2L\eta^4u^2).
\end{equation}
The Gaussian second-moment bound
$\Prob(N^{-1}\Tr H^2>4)\le e^{-cW^2L}$ also gives empirical mass
at most $4/R^2$ outside $[-R,R]$, with the same uniformity in $g$.

Fix $\delta>0$. Choose $R>4$ large, then $h>0$ small and $\eta>0$
with $\eta/h$ small, so that
\[
 R^{-2}+\eta/h+\sup_{g,x}\nu_g([x-2h,x+2h])\ll\delta,
\]
using \eqref{rd:eq:noatoms}. On the second-moment event, Poisson
smoothing reduces the distribution-function comparison to resolvent
control at height $\eta$ on $[-R-1,R+1]$. Apply \eqref{rd:eq:bias}
and \eqref{eq:rough-density} on a finite energy mesh; the resolvent
derivative bound $2\eta^{-2}$ fills the gaps. Choose the mesh and
accuracy in terms of $\delta,R,h,\eta$, then take $W,L$ sufficiently
large. Removing the smoothing costs at most $C\eta/h$, the uniform
interval mass above, and the tails. Thus the distribution functions
differ by at most $\delta$, with failure probability
$C_\delta e^{-c_\delta W^2L}$. All choices depend only on $\delta$,
which proves \eqref{eq:rough-counting} uniformly in $g$.
\end{proof}

\section{Fractional moments and eigenvector concentration}
\label{app:strong-concentration}
We prove the fractional-moment conversion in
Lemma~\ref{uva:thm:upper-audit} and the finite-volume criterion in
Proposition~\ref{up:prop:fr-concentration}.

\subsection{From fractional moments to concentration}
\label{app:fr-concentration-proof}

\begin{proof}[\textbf{Proof of Lemma~\ref{uva:thm:upper-audit}\textup{(i)}}]
Define the block correlator
$\mathcal Q_{ab}(J)\coloneqq
\sum_{\lambda_k\in J}\|\psi_{k,a}\|_2\|\psi_{k,b}\|_2$.
Since $g\le W^{C_0}$, the conditional densities and $q$th moments
of the diagonal GUE entries of $H$ are bounded by powers of $W$.
The finite-volume spectral-averaging argument in
\cite[Appendix~A, proof of Theorem~A.1, (A.10)--(A.12)]{ASFH},
applied to the integral over $J$ of \eqref{up:eq:fmr-pointwise}
and summed over the $W^2$ coordinate pairs in blocks $a,b$, gives
\begin{equation}\label{up:eq:correlator-intro}
 \E\mathcal Q_{ab}(J)
 \le W^C e^{-c\|a-b\|_L/\mathfrak b},\qquad a\ne b.
\end{equation}

Choose a maximal-mass block $a_k$ for each eigenvector, and let
$B_D(R)$ be the number counted in \eqref{eq:averaged-tail}.
Since $\|\psi_{k,b}\|_2^2\le
\|\psi_{k,a_k}\|_2\|\psi_{k,b}\|_2$, summing
\eqref{up:eq:correlator-intro} gives
\begin{equation}\label{uv:eq:ordinary-count}
 B_D(R)\le W^D\sum_a\sum_{\|b-a\|_L>R}\mathcal Q_{ab}(J),\qquad
 \frac{\E B_D(R)}{WL}
 \le W^{D+C-1}\mathfrak b\,e^{-cR/\mathfrak b}.
\end{equation}
As $\mathfrak b\le W^{C_0}$, taking
$R=\lceil C_D\mathfrak b\log W\rceil$ with sufficiently large
$C_D$ proves \eqref{eq:averaged-tail}.
\end{proof}

\begin{proof}[\textbf{Proof of Lemma~\ref{uva:thm:upper-audit}\textup{(ii)}}]
To obtain \eqref{eq:averaged-tail}, we first control the mass correlations
between distant cells and then sum over their positions.
Use the assumptions of part~\textup{(ii)}. Since $g\ge W$, for
sufficiently large $W$ Lemma~\ref{uvs:lem:transfer-moments} applies
on $J$ with window length at most $2g^2$.
Set $n=\lfloor g^2\rfloor\asymp g^2$;
the constants and width thresholds below are uniform in $g,L$.
For $L\ge3n$, partition the cycle into consecutive
cells of lengths in $[n,2n]$. Write $n_A=|A|$ for each cell and let $d(A,B)$ be
the minimum block distance between two cells. Put
\[
 m_A(k)=\|P_A\psi_k\|^2,\qquad
 Q(A,B)=\sum_{\lambda_k\in J}\sqrt{m_A(k)m_B(k)},
\]
where $P_A$ projects onto all coordinates in $A$. Our target is
\begin{equation}\label{uva:eq:cell-correlator}
 \E Q(A,B)\le nW^C e^{-cd(A,B)/(Wg^2)}.
\end{equation}
For $A=B$, this follows from $Q(A,A)\le n_AW\le2nW$.
For $A\ne B$, set $t=q/4$ and
\[
 Q_t(A,B)=\sum_{\lambda_k\in J}m_A(k)^{1-t}m_B(k)^t.
\]
Cauchy--Schwarz over eigenstates gives
$Q(A,B)^2\le Q_t(A,B)Q_t(B,A)$, so it suffices to prove the same
decay bound for $\E Q_t(A,B)$.

To express the mass ratio $m_B/m_A$ through the resolvent on $A^c$,
we first compare the boundary values of an eigenvector with its cell mass.
For the transfer products in a cell $A$, number its blocks by
$1,\ldots,n_A$ in their natural order, and put
\[
 F_j(E)=T_j(E)\cdots T_1(E),\qquad
 K_A(E)=\max_{0\le j\le n_A}
             \max\{1,\|F_j(E)\|,\|F_j(E)^{-1}\|\}.
\]
Here $T_j(E)$ is defined in Definition~\ref{def:transfer-matrices}
and $F_0=I$. Since $n_A\le2g^2$,
Lemma~\ref{uvs:lem:transfer-moments} gives $\E K_A(E)^p\le C_p$
for every fixed $p>0$, uniformly on $J$.
Fix a normalized eigenpair $(\lambda,\psi)$, suppress its index in
$m_A,m_B$, and evaluate $K_A,K_B$ at $\lambda$.
Let $\psi_{A,\rm left},\psi_{A,\rm right}$ be the endpoint values on $A$,
and $x_-,x_+$ their respective neighbors outside $A$, the endpoints
of the Dirichlet path $\Lambda=A^c$.

The inverse-prefix bound, summed over the $n_A-1$ interior pairs,
controls the initial pair $(\psi_{A,\rm left},\psi_{x_-})$ by
$CK_A\sqrt{m_A/n_A}$. Forward propagation on $A,B$ then gives,
with $b_0,b_0+1$ the first two blocks of $B$,
\[
 \begin{aligned}
 \norm{\psi_{A,\rm left}}+\norm{\psi_{A,\rm right}}
 &\le CK_A^2\sqrt{m_A/n_A},\\
 m_B&\le Cn_BK_B^4
       (\norm{\psi_{b_0}}^2+\norm{\psi_{b_0+1}}^2).
 \end{aligned}
\]
Almost surely, $H$ and $H_\Lambda$ have no common eigenvalue:
condition on the exterior and all Gaussian directions in $A$ except
the common scalar shift. Every eigenvalue is strictly increasing in
this shift, since a solution vanishing on two consecutive blocks is zero.
Thus
\[
 \psi_b=-G_\Lambda(b,x_-;\lambda)\psi_{A,\rm left}
        -G_\Lambda(b,x_+;\lambda)\psi_{A,\rm right}.
\]
Since $n_B/n_A\le2$,
\begin{equation}\label{uva:eq:ratio-audit}
 \frac{m_B}{m_A}\le CK_A^4K_B^4
       \sum_{b\in\{b_0,b_0+1\},\ x\in\{x_-,x_+\}}
                    \norm{G_\Lambda(b,x;\lambda)}^2.
\end{equation}
For fixed $A,B$, define the exterior-measurable function
\[
 f(E,H_{A^c})=K_B(E)^{4t}
       \sum_{b\in\{b_0,b_0+1\},\ x\in\{x_-,x_+\}}
                    \|G_\Lambda(b,x;E)\|^{2t}.
\]
Since $4t=q<1$ and $K_A\ge1$, \eqref{uva:eq:ratio-audit} gives
\begin{equation}\label{uva:eq:averaging-reduction}
 Q_t(A,B)\le C\sum_{\lambda_k\in J}
        m_A(k)K_A(\lambda_k)f(\lambda_k,H_{A^c}).
\end{equation}
The input \eqref{uva:eq:conversion-input} controls an energy integral,
whereas \eqref{uva:eq:averaging-reduction} is evaluated at the random
eigenvalues. We therefore need the following averaging estimate, which
retains the transfer weight and is uniform in the exterior configuration.

\medskip\noindent\textbf{Claim (Weighted averaging over one cell).}
For each cell $A$ and every nonnegative measurable $f(E,H_{A^c})$,
\begin{equation}\label{uva:eq:average}
 \E_A\sum_{\lambda_k\in J}
 m_A(k)K_A(\lambda_k) f(\lambda_k,H_{A^c})
 \le C n_AW^2\int_J f(E,H_{A^c})\,dE.
\end{equation}
Here $\E_A$ averages over the blocks in $A$, with the exterior fixed.
We use classical spectral averaging and monotone boundary phases;
see \cite[Appendix~A]{ASFH} and
\cite[Section~4, Theorem~5 and Lemma~2]{SSB-averaging}.
The proof below keeps the transfer weight and the uniform constants
needed on cells of length comparable to $g^2$.
\begin{proof}[Proof of the claim]
Fix the exterior configuration. Since $\operatorname{Var}(\Tr V_j)=1$,
Gaussian orthogonal decomposition gives
\[
 u=n_A^{-1/2}\sum_{j\in A}\Tr V_j\sim N(0,1),\qquad
 V_j=V_j^\circ+\frac{u}{W\sqrt{n_A}}I_W,
\]
where $u$ is independent of the residual blocks $(V_j^\circ)_{j\in A}$.
Write $\E_u$ for expectation over $u$ at fixed residual blocks and
$\E_\circ$ for their expectation. Their joint average
$\E_A=\E_\circ\E_u$ restores the original independent GUE law on $A$.
We have $H(u)=H^\circ+\alpha uP_A$, with $\alpha=(gW\sqrt{n_A})^{-1}$;
$K_A(E,u)$ denotes the same transfer weight for these shifted blocks.
By the two-block uniqueness used above, the compression of $P_A$
to every eigenspace is positive definite.
The Rellich branches satisfy $\lambda_k'(u)=\alpha m_A(k)>0$,
with this compression diagonalized at crossings. The scalar area
formula, with standard Gaussian density $\varphi$, yields
\begin{equation}\label{uva:eq:coarea-audit}
 \E_u\sum_{\lambda_k\in J}m_A(k)K_A(\lambda_k,u)f(\lambda_k,H_{A^c})
 =\alpha^{-1}\int_J f(E,H_{A^c})
 \sum_{u:E\in\sigma(H(u))}\varphi(u)K_A(E,u)\,dE.
\end{equation}
In part~\textup{(i)}, the conditional densities arising in this conversion
from eigenvalue sums to energy integrals are bounded by powers of $W$
because $g\le W^{C_0}$. Here cell averaging and normalization remove
this coupling dependence.
Here and below, each root $u$ is counted with multiplicity
$\dim\ker(H(u)-E)$. It suffices to prove, for almost every $E\in J$,
\begin{equation}\label{uva:eq:weighted-root-sum}
 \E_\circ\sum_{u:E\in\sigma(H(u))}
       \varphi(u)K_A(E,u)\le CW.
\end{equation}

Fix $E\notin\sigma(H_{A^c})$.
For a solution of the cell recurrence put
\[
 v=(\psi_1,\psi_{n_A}),\quad w=(\psi_0,\psi_{n_A+1}),\qquad
 a=2^{-1/2}(v+iw),\quad b=2^{-1/2}(v-iw).
\]
Discrete Green's identity gives
$v^*\widetilde w-w^*\widetilde v=0$ for any two cell solutions.
Thus the boundary data form a Lagrangian space. The map from the
$2W$ initial data to $a\in\C^{2W}$ is injective: $a=0$ implies
$b=0$, hence $v=w=0$ and then the entire recurrence vanishes.
It is therefore invertible. We obtain a smooth unitary
$b=C_A(E,u)a$ for every real $u$.

The solution depends linearly on $a$; let
$\Psi_A:\C^{2W}\to\C^{n_AW}$ map $a$ to $(\psi_1,\ldots,\psi_{n_A})$.
Since the initial pair $(\psi_1,\psi_0)$ is part of $(v,w)$ and
$\norm v^2+\norm w^2=2\norm a^2$,
\begin{equation}\label{uva:eq:psi}
 \norm{\Psi_A}^2\le2n_AK_A^2.
\end{equation}
Differentiating the solution at fixed $a$ and summing by parts gives
\[
 v^*\partial_uw-w^*\partial_uv=-\alpha\norm{\Psi_Aa}^2.
\]
Together with $a^*da-b^*db=i(v^*dw-w^*dv)$, this quadratic-form
identity yields
\begin{equation}\label{uva:eq:generator-audit}
 -iC_A^*\partial_uC_A=\alpha\Psi_A^*\Psi_A>0,\qquad
 \norm{-iC_A^*\partial_uC_A}\le2\alpha n_AK_A^2.
\end{equation}
Strict positivity follows from the same two-block uniqueness property.

Eliminating the fixed exterior gives a unitary $C_{\rm ext}(E)$,
independent of $u$, in the same boundary coordinates. The full spectral
equation is $1\in\sigma(C_{\rm ext}^*C_A(E,u))$.
The positive generator in \eqref{uva:eq:generator-audit} makes all
lifted eigenphases strictly increasing, also at multiple phases by
compression to their eigenspaces. Let $N_I(E)$ count the roots in a
unit interval $I$, with multiplicity. The sum of the eigenphase increments
equals the integral of the trace of the generator, so
\begin{align}
 N_I(E)&\le2W+\frac1{2\pi}\int_I
                  \Tr(-iC_A^*\partial_uC_A)\,du\notag\\
 &\le2W+C\alpha n_AW\sup_IK_A^2
 \le CW(1+\sup_IK_A^2).\label{uva:eq:cross-audit}
\end{align}
Shared endpoints of adjacent intervals only overcount roots; the finite
set of excluded energies does not affect \eqref{uva:eq:coarea-audit}.

It remains to average weighted suprema. Product differentiation gives
\[
 \partial_uF_j=\sum_{k=1}^jF_jF_k^{-1}(\partial_uT_k)F_{k-1},\qquad
 \norm{\partial_uT_k}=\alpha.
\]
Using also the derivative of the inverse, the finite maximum $K_A$
is locally Lipschitz and satisfies almost everywhere
\[
 |\partial_uK_A|\le C\alpha n_AK_A^5.
\]
Since $\alpha n_A\le C/W$, applying
$\sup_I z\le\int_Iz+\int_I|z'|$ to $z=K_A^3$ on
$I_\ell=[\ell,\ell+1]$ gives
\begin{equation}\label{uva:eq:sobolev-audit}
 \sup_{I_\ell}K_A^3\le C\int_{I_\ell}K_A^7\,du.
\end{equation}
For $u\in I_\ell$, $\sup_{I_\ell}\varphi/\varphi(u)
\le e^{|u|+C}$.
Using \eqref{uva:eq:sobolev-audit}, Cauchy--Schwarz, and
Lemma~\ref{uvs:lem:transfer-moments},
\begin{align}
 \E_\circ\sum_\ell(\sup_{I_\ell}\varphi)
                         (\sup_{I_\ell}K_A^3)
 &\le C\E_A[e^{|u|}K_A^7]\notag\\
 &\le C(\E_A e^{2|u|})^{1/2}
              (\E_A K_A^{14})^{1/2}\le C.\label{uva:eq:weight-audit}
\end{align}
Since $K_A\ge1$, \eqref{uva:eq:cross-audit} and
\eqref{uva:eq:weight-audit} give \eqref{uva:eq:weighted-root-sum}.
Finally, $\alpha^{-1}=gW\sqrt{n_A}\le\sqrt2 n_AW$ proves
\eqref{uva:eq:average}.
\end{proof}

The averaging cost before using $n_A\asymp g^2$ is $CgW^2\sqrt{n_A}$.
After normalization by the cell dimension it is
\[
 \frac{CgW^2\sqrt{n_A}}{n_AW}
       =\frac{CgW}{\sqrt{n_A}}\asymp CW.
\]

We now apply \eqref{uva:eq:average} to the function $f$ in
\eqref{uva:eq:averaging-reduction}. Averaging the exterior and using
Cauchy--Schwarz jointly in energy and the exterior variables gives
\[
 \begin{aligned}
 \E Q_t(A,B)&\le CnW^2\sum_{b,x}\int_J
           \E[K_B(E)^{4t}\norm{G_\Lambda(b,x;E)}^{2t}]\,dE\\
 &\le nW^C e^{-cd(A,B)/(Wg^2)}.
 \end{aligned}
\]
The second inequality uses the $8t$th transfer moment from
Lemma~\ref{uvs:lem:transfer-moments}, the $4t=q$ moment in
\eqref{uva:eq:conversion-input}, and that the relevant path distances
are at least $d(A,B)-1$.
Interchanging $A,B$ and applying Cauchy--Schwarz in probability to
$Q(A,B)^2\le Q_t(A,B)Q_t(B,A)$ proves \eqref{uva:eq:cell-correlator}.

\medskip\noindent\textbf{Exceptional proportions.}
If the proposed block radius fills the cycle, the tail is empty;
otherwise $L\ge3n$ and we use the cell partition above.
For each eigenvector choose a maximal-mass cell $A_k$.
Since $m_B(k)\le\sqrt{m_{A_k}(k)m_B(k)}$, the number of eigenvectors
in $J$ with $\sum_{B:d(A_k,B)>r}m_B(k)>W^{-D}$ is at most
$W^D\sum_A\sum_{B:d(A,B)>r}\sum_{\lambda_k\in J}\sqrt{m_A(k)m_B(k)}$.
Using $O(L/n)$ cells, \eqref{uva:eq:cell-correlator}, and the geometric
tail sum gives a normalized exceptional proportion at most
\[
 C W^D\frac{L/n}{WL}\,nW^C\frac{\mathfrak b}{n}
       e^{-c'r/\mathfrak b}
 \le W^{D+C}e^{-c'r/\mathfrak b},
\]
since $\mathfrak b=Wg^2$ and $\mathfrak b/n\le CW$.
Choose $r=C_D\mathfrak b\log W$ with $c'C_D>2D+C+1$.
Every cell at distance at most $r$ from the maximal-mass cell lies
in a block neighborhood of radius at most $r+4n\le C_DWg^2\log W$, after increasing $C_D$, around a block
of the maximal-mass cell. This proves \eqref{eq:averaged-tail}
and the half-mass bound with the cap $WL$.
\end{proof}

\subsection{Finite-volume criterion for fractional-moment decay}
\label{app:finite-volume-criterion}

\begin{proof}[\textbf{Proof of Proposition~\ref{up:prop:fr-concentration}}]
Fix $E\in J$ and the length $n$ from \eqref{uv:eq:upper-volume-input}.
For a path $\Gamma=[l,b]$, $x<y$ with $y-x>n$, put
$\Lambda=[x,r]$, $r=x+n$. Conditional on the exterior of $\Lambda$,
the Schur complement gives
\[
 G_\Gamma(x,y;E)=-\frac{g}{\max\{1,g\}}
 G_\Lambda^{L_0,R_0}(x,r;E)G_{[r+1,b]}(r+1,y;E),
\]
where $L_0,R_0$ are Hermitian shifts determined by the exterior.
The last factor is fixed under this conditional average, so
\eqref{uv:eq:upper-volume-input} yields
\[
 \E\!\left[\|G_\Gamma(x,y;E)\|^s\mid (V_i)_{i\notin\Lambda}\right]
 \le W^{-2}\|G_{[r+1,b]}(r+1,y;E)\|^s.
\]
No independence between $R_0$ and the last factor is needed.
Iterating as in \cite[proof of Theorem~1.1, (2.23)--(2.24)]{ASFH},
using \eqref{eq:c-apriori} on the remainder and
$2\log W/(n+1)\ge c/\mathfrak b$, gives
\begin{equation}\label{eq:c-interval}
 \E\|G_\Gamma(x,y;E)\|^s\le W^C e^{-c|x-y|/\mathfrak b},
 \qquad E\in J,\quad x,y\in\Gamma.
\end{equation}
The constants are independent of the interval length; adjunction
handles $x>y$, and \eqref{eq:c-apriori} handles $x=y$.

On a cycle, deleting block $a$ leaves a Dirichlet path with resolvent
$G^{(a)}=(H_{\Z_L\setminus\{a\}}-E)^{-1}$. The Schur identity is
\begin{equation}\label{eq:c-cycle-Schur}
 G_{ab}=-\frac{g}{\max\{1,g\}}G_{aa}
          \bigl(G^{(a)}_{a+1,b}+G^{(a)}_{a-1,b}\bigr),\qquad a\ne b.
\end{equation}
Both path distances are at least $\|a-b\|_L-1$. Cauchy--Schwarz,
\eqref{eq:c-apriori} and \eqref{eq:c-interval}, with $q=s/2$, give
\begin{equation}\label{eq:c-cycle}
 \E\|G_{ab}(E)\|^q\le W^C e^{-c\|a-b\|_L/\mathfrak b},
 \qquad a\ne b,\quad E\in J.
\end{equation}
Jensen's inequality lowers the exponent in \eqref{eq:c-interval}
to $s/2$; \eqref{eq:c-apriori} handles the diagonal. This proves
\eqref{up:eq:fmr-pointwise} on paths and cycles.
\end{proof}

\section{Circle identities and operator estimates}
\label{app:compact-circle}

We use the form $J_{\rm sig}$ and group $\mathcal G$ from
Definition~\ref{up:def:isotropic-group}, and the volume ratio $c$ from
\eqref{eq:c-volume-cocycle}.
The operator iteration is given in Section~\ref{sec:upper-inputs}.

\subsection{The quotient circle and its exact Jacobian}

\begin{lemma}[Signature quotient]\label{cc:lem:quotient}
Let $U\subset\C^{2W}$ be a deterministic
subspace of dimension $W-1$ on which the form
$x^*J_{\rm sig}y$ vanishes. Then
\[
 H_U\coloneqq (U+J_{\rm sig}U)^\perp
\]
is \(J_{\rm sig}\)-invariant, has complex dimension two, and has signature \((1,1)\).
Choose unit vectors \(e,f\in H_U\) with \(J_{\rm sig}e=e\) and \(J_{\rm sig}f=-f\). Then
the extensions are \(Y=U\oplus\mathbb C v_\theta\), where
\begin{equation}
 v_\theta\coloneqq (e+e^{i\theta}f)/\sqrt2,\qquad
 dm_U=d\theta/(2\pi).
 \label{cc:eq:circle}
\end{equation}
For \(\mathcal T\in\mathcal G\), set
\[
 A_{\mathcal T}:H_U\longrightarrow H_{\mathcal T U},\qquad
 A_{\mathcal T}v\coloneqq \operatorname{proj}_{H_{\mathcal T U}}\mathcal T v .
\]
This is the induced \(J_{\rm sig}\)-unitary quotient isomorphism. For the volume
ratio $c(\mathcal T;U,Y)$ in \eqref{eq:c-volume-cocycle}, the angular Jacobian
is \eqref{cc:eq:jac}; equivalently,
$e^{-2c(\mathcal T;U,Y)}=\|A_{\mathcal T}v_\theta\|^{-2}$.
The function \(c\) is a cocycle, and
\begin{equation}
 |c(\mathcal T;U,Y)|\le\log\|\mathcal T\|.
 \label{cc:eq:boundc}
\end{equation}
\end{lemma}
\begin{proof}[\textbf{Proof of \eqref{cc:eq:jac}}]
Isotropy makes \(U\perp J_{\rm sig}U\). Their sum is \(J_{\rm sig}\)-invariant and has
signature \((W-1,W-1)\), which proves the claims concerning \(H_U\).
Moreover \(U^{\perp_{J_{\rm sig}}}=U\oplus H_U\). If \(v\in H_U\), then
\(\mathcal T v\in(\mathcal T U)^{\perp_{J_{\rm sig}}}=\mathcal T U\oplus H_{\mathcal T U}\). The removed component lies in
the isotropic space \(\mathcal T U\) and is \(J_{\rm sig}\)-orthogonal to all of
\((\mathcal T U)^{\perp_{J_{\rm sig}}}\). Thus the projection preserves the quotient form.
The quotient map for \(\mathcal T^{-1}\) is its inverse; their norms are bounded
by \(\|\mathcal T\|\) and \(\|\mathcal T^{-1}\|\), respectively. Since
\(\mathcal T^{-1}=J_{\rm sig}\mathcal T^*J_{\rm sig}\), we have \(\|\mathcal T^{-1}\|=\|\mathcal T\|\), hence
\(\|\mathcal T\|^{-1}\le\|A_{\mathcal T}v_\theta\|\le\|\mathcal T\|\).

Choose unit $J_{\rm sig}$-eigenvectors $e',f'$ of $H_{\mathcal T U}$
with eigenvalues $+1,-1$. In the bases $e,f$ of $H_U$ and $e',f'$ of
$H_{\mathcal T U}$, write \(A_{\mathcal T}=\begin{psmallmatrix}a&b\\c&d\end{psmallmatrix}\).
On \(|z|=1\) its action is \(z'=(c+dz)/(a+bz)\).
Since \(A_{\mathcal T}\) preserves \(\diag(1,-1)\), the two output components of
\(A_{\mathcal T}(1,z)^t/\sqrt2\) have equal modulus and
\(\|A_{\mathcal T}v_\theta\|^2=|a+bz|^2\).
Differentiating the fraction, and using the form-preserving identities
and \(|\det A_{\mathcal T}|=1\), gives
\(d\theta'/d\theta=|a+bz|^{-2}>0\).
Finally Gram--Schmidt on \(\mathcal T Y\), beginning with an orthonormal basis of
\(\mathcal T U\), identifies \(\|A_{\mathcal T}v_\theta\|\) with the last volume ratio.
The cocycle identity is the multiplicativity of restricted volumes,
and the norm bounds on the quotient and its inverse give \eqref{cc:eq:boundc}.
\end{proof}

The reciprocal convention is useful to record:
\[
 \frac{d(\mathcal T_*m_U)}{dm_{\mathcal T U}}(\mathcal T Y)=e^{2c(\mathcal T;U,Y)}.
\]
This circle has real dimension one. The exponent two in
\eqref{cc:eq:jac} must not be doubled because the ambient vector spaces
are complex.

\subsection{The inverse matrix and the last-volume ratios}

\begin{lemma}[Omitted-column identity]\label{cc:lem:cofactor}
Let $\mathcal T\in\mathcal G$ for the group \eqref{eq:c-group}, fix $0<s<1$,
and let $\Phi\in\mathcal F_W$ from \eqref{up:eq:boundary-frame},
\(Y\coloneqq \operatorname{ran}\Phi\), and \(U_j\) the span of all its columns
except the \(j\)-th. Set $D\coloneqq\sqrt2\begin{pmatrix}I_W&0\end{pmatrix}\mathcal T\Phi$.
Then
\begin{equation}
 e^{-c(\mathcal T;U_j,Y)}
 =\sqrt{[(D^*D)^{-1}]_{jj}}
 =\|e_j^*D^{-1}\|.
 \label{cc:eq:cofactor}
\end{equation}
Consequently
\begin{equation}
 \|D^{-1}\|^s\le\sum_{j=1}^W e^{-sc(\mathcal T;U_j,Y)} .
 \label{cc:eq:inverse}
\end{equation}
\end{lemma}
\begin{proof}[\textbf{Proof of Lemma~\ref{cc:lem:cofactor} and \eqref{eq:c-inverse-cocycle}}]
The identity \((\mathcal T\Phi)^*J_{\rm sig}(\mathcal T\Phi)=0\) equates the Gram matrices
of the upper and lower blocks, including cross terms. Thus
\((\mathcal T\Phi)^*(\mathcal T\Phi)=D^*D>0\), which explains the factor \(\sqrt2\)
and proves that \(D\) is invertible.
The squared restricted \(W\)-volume is \(\det(D^*D)\), and the squared
\((W-1)\)-volume is its \(j\)-th principal minor. Their quotient is
the reciprocal of the \(j\)-th diagonal entry of \((D^*D)^{-1}\).
Also \((D^*D)^{-1}=D^{-1}D^{-*}\).
Finally,
\(\|D^{-1}\|^s\le\|D^{-1}\|_\mathrm{HS}^s
\le\sum_j\|e_j^*D^{-1}\|^s\), by concavity for \(s/2<1\).
\end{proof}

\subsection{Bounds for the good and bad operators}
\label{app:c-operator-proofs}

The following estimates apply to any smooth good component with the
stated support and derivative bounds. Use the Lie algebra
\eqref{up:eq:lie-algebra}. For a smooth density on $\mathcal G$, use the right derivative
\cite[(3)]{Drogin}
\begin{equation}\label{eq:c-right-derivative}
 X_Mf(\mathcal T)\coloneqq\left.\frac{d}{dt}f(\mathcal T e^{tM})\right|_{t=0},\qquad
 \|X_Mf\|_1\coloneqq\int|X_Mf(\mathcal T)|\,d\mathcal T,
 \qquad M\in\mathfrak g.
\end{equation}
Here $d\mathcal T$ is Haar measure on $\mathcal G$. We use the real
Hilbert--Schmidt inner product
$\langle A,B\rangle=\Re\Tr(A^*B)$ on $\mathfrak g$.

\begin{lemma}[Bounds for the good and bad operators]\label{lem:c-operators}
Use the group and operators defined in Section~\ref{sec:upper-inputs},
the Lie algebra in \eqref{up:eq:lie-algebra}, and the right derivatives
in \eqref{eq:c-right-derivative}.
Let $\mu=\widetilde{\pi}f\,d\mathcal T+\mu_{\rm bad}$ be a probability measure,
where $0<\widetilde{\pi}\le1$, $f$ is a smooth compactly supported
probability density, and $\mu_{\rm bad}$ is a positive measure satisfying
$\int\|\mathcal T\|^{1/2}\,d\mu_{\rm bad}(\mathcal T)\le1/2$. Set
\[
 K\coloneqq\sup_{\mathcal T\in\supp f}\max\{\|\mathcal T\|,\|\mathcal T^{-1}\|\}.
\]
There is an absolute constant $C<\infty$ such that every bounded
measurable $h$ satisfies
\begin{align}
 \|G_{1/2}h\|_\infty
 &\le CK\left(1+
 \sup_{\substack{M\in\mathfrak g\\\|M\|_{\rm HS}=1}}\|X_Mf\|_1
          \right)^{1/4}
 \times\sup_U\|h(U,\cdot)\|_{L^4(\mathcal C_U,dm_U)},
 \label{cc:eq:smoothing}\\
 \|B_{1/2}h\|_\infty&\le\tfrac12\|h\|_\infty.
 \label{eq:c-bad-operator}
\end{align}
\end{lemma}

\begin{proof}[Proof of Lemma~\ref{lem:c-operators}]
\medskip\noindent\emph{From angular $L^4$ to a fixed flag.}
Fix $U,Y$ and take the compact Lie direction $K_U$ which is
$\frac12\diag(-i,i)$ in an orthonormal $J_{\rm sig}$-eigenbasis of
$H_U$, ordered $+1,-1$, and zero on $U+J_{\rm sig}U$. It is skew-Hermitian, has bounded
Hilbert--Schmidt norm, and rotates the source circle with unit
angular speed while fixing its base. For a function $\psi$,
measurable in the target base and smooth on each circle,
\[
 X_{K_U}\psi(\mathcal T U,\mathcal T Y)=e^{-2c(\mathcal T;U,Y)}\partial_\theta\psi(\mathcal T U,\mathcal T Y).
\]
On the support of $f$, Lemma~\ref{cc:lem:quotient} gives
\[
 e^{2c(\mathcal T;U,Y)}\le K^2,\qquad
 |X_{K_U}(e^{2c(\mathcal T;U,Y)})|\le CK^2.
\]
Indeed, $e^{2c(\mathcal T;U,Y)}=\|A_{\mathcal T}v_\theta\|^2$, and differentiation of the
source unit vector gives the second bound from $\|A_{\mathcal T}\|\le K$.
Haar integration by parts therefore yields
\begin{equation}\label{cc:eq:vertical}
 \begin{aligned}
 \left|\int\partial_\theta\psi(\mathcal T U,\mathcal T Y)f(\mathcal T)\,d\mathcal T\right|
 &\le\|\psi\|_\infty\|X_{K_U}(e^{2c(\mathcal T;U,Y)}f(\mathcal T))\|_1\\
 &\le CK^2\left(1+
 \sup_{\substack{M\in\mathfrak g\\\|M\|_{\rm HS}=1}}
                \|X_Mf\|_1\right)\|\psi\|_\infty.
 \end{aligned}
\end{equation}
This step uses the derivative bound over all Lie directions, since
$K_U$ is skew-Hermitian.

Disintegrate the output law of $(\mathcal T U,\mathcal T Y)$ under $f(\mathcal T)\,d\mathcal T$ over
its base marginal $\bar\nu(dV)$. Equation~\eqref{cc:eq:vertical}
bounds the vertical distributional derivative. Choose a countable
dense family of smooth circle tests with supremum norm at most one.
In \eqref{cc:eq:vertical}, choose measurably at each base the largest
signed value among the first finitely many tests, then pass to the
limit. This bounds the integral of the conditional total variations,
so the conditional circle laws have BV densities $k_V$ with
\[
 \int\|\partial_\theta k_V\|_{\rm TV}\,d\bar\nu(V)
 \le CK^2\left(1+
 \sup_{\substack{M\in\mathfrak g\\\|M\|_{\rm HS}=1}}
                \|X_Mf\|_1\right).
\]
For clarity, a circle measure with a finite distributional derivative
has a BV density: convolution by circle mollifiers gives densities
with uniformly bounded variation, and BV compactness gives the
claim on almost every conditional circle. Since $\int k_V\,dm_V=1$,
\[
 \int\|k_V\|_\infty\,d\bar\nu(V)
 \le CK^2\left(1+
 \sup_{\substack{M\in\mathfrak g\\\|M\|_{\rm HS}=1}}
                \|X_Mf\|_1\right).
\]
The weight $e^{-c(\mathcal T;U,Y)/2}$ is at most $K^{1/2}$.
H\"older's inequality on each circle and then concavity give
\[
 \begin{aligned}
 |(G_{1/2}h)(U,Y)|
 &\le K^{1/2}\int\!\int |h(V,Z)|k_V(Z)\,dm_V(Z)\,d\bar\nu(V)\\
 &\le K^{1/2}\sup_V\|h(V,\cdot)\|_{L^4(\mathcal C_V,dm_V)}
           \int\|k_V\|_{4/3}\,d\bar\nu(V)\\
 &\le K^{1/2}\sup_V\|h(V,\cdot)\|_{L^4(\mathcal C_V,dm_V)}
           \left(\int\|k_V\|_\infty\,d\bar\nu(V)\right)^{1/4}.
 \end{aligned}
\]
Here $\widetilde{\pi}\le1$ was discarded and
$\|k_V\|_{4/3}\le\|k_V\|_\infty^{1/4}$.
This proves \eqref{cc:eq:smoothing}, uniformly for the fixed source
flag. No density estimate on the base marginal is required.

\medskip\noindent\emph{The bad component.}
Finally, $|c(\mathcal T;U,Y)|\le\log\|\mathcal T\|$ gives
\[
 |(B_{1/2}h)(U,Y)|
 \le\|h\|_\infty\int\|\mathcal T\|^{1/2}\,d\mu_{\rm bad}(\mathcal T)
 \le\tfrac12\|h\|_\infty,
\]
which is \eqref{eq:c-bad-operator}.
\end{proof}

\section{GUE estimates and the common density}
\label{app:upper-gue}

The relative local law, preparation estimates and Riccati product bounds
are used in Section~\ref{su:sec:construction-estimates}. The deterministic
Dyson profile is given by Lemma~\ref{da:lem:deformed} in
Appendix~\ref{app:deformed-dyson}. Structured
Gaussian coercivity then supplies the cutoff for the common density,
whose construction proves Lemma~\ref{up:lem:good-bad}. The proof
is completed by the $G_{1/2},B_{1/2}$ estimates of
Appendix~\ref{app:c-operator-proofs}.

\subsection{The relative local law and preparation of the incoming flag}
\label{app:relative-local-law-proof}
\label{su:sec:spectral-proofs}

\begin{lemma}[Local law for a deformed GUE]\label{us:lem:local-law}
Fix $\kappa_0\in(0,2)$, $C_*,K>0$, and sufficiently small fixed
$\varepsilon>0$. Let $B=B_{\rm small}+B_{\rm bad}$ be Hermitian with
$\|B_{\rm small}\|\le\varepsilon$, $\rank B_{\rm bad}\le\varepsilon W$
and $\|B\|\le W^K$. Let $V$ be normalized complex
GUE, and set
\[
 G\coloneqq(V-EI_W-B-i\eta)^{-1},\qquad
 M\coloneqq M_B(i\eta),\qquad m\coloneqq W^{-1}\Tr M,
\]
where $M_B$ is the Herglotz solution of Lemma~\ref{da:lem:deformed}
with $\Im m>0$.
For every fixed $\delta>0$, uniformly for
$|E|\le2-\kappa_0$ and $W^{-1+\delta}\le\eta\le C_*$,
\begin{equation}\label{us:eq:weighted-local-law}
 \left|W^{-1}\Tr G-m\right|\prec W^{-1}\eta^{-1},\qquad
 |u^*(G-M)v|\prec
 (W\eta)^{-1/2} \norm{Mu}_2\norm{Mv}_2
\end{equation}
for deterministic vectors $u,v$. In particular, for each deterministic
$v\ne0$,
\begin{equation}\label{us:eq:relative-local-law}
 \left|\frac{\Im(v^*Gv)}{v^*\Im Mv}-1\right|
 \prec(W\eta)^{-1/2},\qquad
 \Im(v^*Gv)\asymp v^*(I_W+B^2)^{-1}v
\end{equation}
with the comparison holding with arbitrarily high polynomial probability.
The constants and probability bounds are uniform over the deterministic
parameters just specified; the event need not be common to all vectors.
\end{lemma}

For the isotropic and anisotropic local-law framework see \cite{KYiso}
and \cite[Section~12, Theorem~12.2]{KYani}.

\begin{proof}[\textbf{Proof of \eqref{us:eq:weighted-local-law} and \eqref{us:eq:relative-local-law}}]
We use the averaged estimate in \cite[Theorem~1.9]{HKR}, with
deformation $-EI_W-B$ and spectral parameter $i\eta$.
Its bounded-deformation assumption can be relaxed as explained in
\cite[Remark~1.11]{HKR}: Theorem~1.5 there has no such assumption.
The deterministic inputs to its Section~4 proof are
\eqref{us:eq:deformed-profile}, \eqref{da:eq:deformed-scalar-stability}, and
\eqref{da:eq:deformed-vector-weight}.
The local stability bound \eqref{da:eq:deformed-residual-stability}
verifies \cite[Definition~4.9]{HKR}, uniformly along the height interval.
The norm and density bounds in \eqref{da:eq:deformed-vector-weight}
replace the uses of a bounded deformation in its
Lemmas~4.1--4.2. In its Lemma~4.6, $\norm B\le W^K$ only requires
a finer polynomial mesh for the resolvent-continuity estimates.
The proof of its Theorem~1.9 therefore gives the first estimate in
\eqref{us:eq:weighted-local-law}, with uniform constants. The complex
case is covered by \cite[Appendix~B]{HKR}.

We next retain the vector weights, using the averaged estimate in
\eqref{us:eq:weighted-local-law} and the Gaussian structure. By unitary invariance we may diagonalize
$B$.
Let $G^{(i)}$ be the resolvent with row and column $i$ deleted, and let
$h_i\coloneqq(V_{ki})_{k\ne i}$. Schur's formula gives
\[
 G_{ii}^{-1}
 =-E-B_{ii}-i\eta+V_{ii}-h_i^*G^{(i)}h_i.
\]
Conditional Gaussian quadratic-form estimates and Ward's identity give
\[
 \left|h_i^*G^{(i)}h_i-W^{-1}\Tr G^{(i)}\right|
 \prec\left(W^{-2}\eta^{-1} \Im\Tr G^{(i)}\right)^{1/2}
 \prec(W\eta)^{-1/2}.
\]
Here the averaged estimate in \eqref{us:eq:weighted-local-law} and the interlacing bound
$|\Tr G-\Tr G^{(i)}|\le C\eta^{-1}$ control the trace on the right.
Consequently,
\begin{equation}\label{us:eq:weighted-entry-law}
 G_{ii}^{-1}=M_{ii}^{-1}+O_\prec((W\eta)^{-1/2}),\qquad
 |G_{ii}-M_{ii}|\prec(W\eta)^{-1/2} |M_{ii}|^2.
\end{equation}
Since $|M_{ii}|\le C$, inversion is justified uniformly, even when
$|M_{ii}|$ is small. The same calculation for $G^{(i)}_{jj}$, using
the resolvent with $i,j$ deleted, gives
$\Im G^{(i)}_{jj}\prec |M_{jj}|^2$; deleting two indices changes
$W^{-1}\Tr G$ by $O((W\eta)^{-1})$.
For $i\ne j$, the off-diagonal Schur formula and conditional Gaussian
linear-form estimate now yield
\begin{equation}\label{us:eq:weighted-offdiagonal-law}
 G_{ij}=-G_{ii}\sum_{k\ne i}V_{ik}G^{(i)}_{kj},\qquad
 |G_{ij}|\prec |M_{ii}|
       \left(W^{-1}\eta^{-1} \Im G^{(i)}_{jj}\right)^{1/2}
 \prec(W\eta)^{-1/2} |M_{ii}M_{jj}|.
\end{equation}

To pass to deterministic vectors, introduce independent uniform unit
complex numbers $\omega_i$. With
$\Omega\coloneqq\diag(\omega_1,\ldots,\omega_W)$, the pair $(G,M)$
has the same distribution as $(\Omega G\Omega^*,M)$.
For every fixed integer $p\ge1$ and every array $(c_{ij})$, the
usual decoupling and Khintchine inequalities give
\[
 \E_\omega\left|\sum_{i\ne j}c_{ij}\omega_i\overline{\omega_j}\right|^{2p}
 \le C_p\left(\sum_{i\ne j}|c_{ij}|^2\right)^p.
\]
Apply this conditionally on $G$, with $c_{ij}=\overline{u_i}G_{ij}v_j$.
Equation \eqref{us:eq:weighted-offdiagonal-law} bounds the sum of squares by
$(W\eta)^{-1}\norm{Mu}_2^2\norm{Mv}_2^2$, in stochastic domination.
By \eqref{us:eq:weighted-entry-law}, the diagonal contribution is bounded by
$(W\eta)^{-1/2}\sum_i|u_iv_i||M_{ii}|^2$, again in stochastic domination.
Cauchy--Schwarz and arbitrarily high moments prove the second estimate
in \eqref{us:eq:weighted-local-law}. Finally, \eqref{da:eq:deformed-vector-weight} converts the weighted
estimate to \eqref{us:eq:relative-local-law}.
\end{proof}

\subsubsection{Preparation of the incoming flag}
\label{su:sec:preparation}

\begin{proof}[\textbf{Proof of Lemma~\ref{su:lem:preparation}}]
To justify chart entry for every incoming flag, write
\[
 SY=\ran\binom AB,\qquad A^*B=B^*A,\qquad \rank\binom AB=W.
\]
Isotropy and full rank imply that $iA-B$ has trivial kernel, so
$\det(zA-B)$ is not the zero polynomial. The first upper block is
therefore invertible almost surely under the fresh GUE perturbation.
The second block gives an inverse Riccati state, followed by two
further steps; in particular,
\[
 R_4=(V_4-EI-g^2R_3)^{-1}.
\]
The line spanned by $R_4^{-1}v_4$ is fixed before $V_4$ is drawn, and
$v_4$ is determined before $V_5$. The one-block bounds
\eqref{us:eq:wegner} give the inverse-norm and eigenvalue-count bounds.
In particular, taking $\delta=2$ there and applying Markov's inequality
to the count on $[-g,g]$ gives regularity \eqref{us:eq:regular} with
failure $O(W^{-2}+g^{1/2})$ for each outgoing state, uniformly in the past.
In regime (i), their total failure is
$O(W^{-2}+W^{-\alpha/2}+W^{-\varepsilon/800})$.
In regime (ii), the inverse-norm failure is $O(W^{-\theta})$, and the
count has mean at most $CWg^2W^\theta\le CW^{2\theta}$, so Markov's
inequality gives failure $O(W^{-\theta})$. The last bound in
\eqref{nc:eq:prepared} follows from $g^2\le W^{-1+\theta}$.

It remains to compare $\|Pv_4\|^2$ with $b_4$ in regime (i).
On a regular past the deformation has polynomial norm, while its small
part and exceptional rank fraction tend to zero uniformly. Since
$Wg^2\ge W^{1-2\beta}$, Lemma~\ref{us:lem:local-law} applies at
$cg^2,g^2,Cg^2$ for fixed positive $c,C$, with common constants and any
prescribed failure power. The tested line is fixed by conditioning on
the past; no common event over all vectors is needed.

Set
\[
 u\coloneqq R_4^{-1}v_4,\qquad
 b_4=g^{-2}\Im\!\left[u^*(R_4^{-1}-ig^2)^{-1}u\right].
\]
Its line is fixed before $V_4$, so the homogeneous local law applies
despite its random length. By \eqref{us:eq:mde}, the Poisson integrals
at the three scales are comparable to
$u^*(I+g^4R_3^2)^{-1}u$. In the integral at $g^2$, the contribution
from $|\lambda|\le cg^2$ is at most $2c$ times the integral at $cg^2$,
and that from $|\lambda|\ge Cg^2$ is at most $2/C$ times the integral
at $Cg^2$. Taking $c$ small and $C$ large leaves a fixed fraction on
the annulus between them, which is retained by $P$ for large $W$.
Comparing the inverse-square and Poisson kernels there gives
\[
 \|Pv_4\|^2\ge c'g^{-2}u^*(I+g^4R_3^2)^{-1}u\ge c''b_4.
\]
This proves $b_4\le C\|Pv_4\|^2$. On the regular past,
$(I+g^4R_3^2)^{-1}\ge cW^{-6}I_W$, and
$\|R_4^{-1}v_4\|\ge W^{-3}\|v_4\|$ since $\|R_4\|\le W^3$.
Thus $\|Pv_4\|^2\ge cW^{-12}\|v_4\|^2$, which gives the weaker stated
$\|Pv_4\|\ge W^{-10}\|v_4\|$ at a common threshold.

\end{proof}

\subsection{Gaussian and spectral inputs for inverse images}
\label{app:inverse-image-inputs}

\subsubsection{Gaussian estimates for the dependent projection}
\label{app:gaussian-projection-inputs}

We verify \eqref{us:eq:lowerK}--\eqref{us:eq:quadratic} under the
assumptions of Lemma~\ref{us:lem:self}. Use the notation from its proof,
\[
 Y_\perp=(I-P_F)RY,\qquad
 E_>=\ran\one_{\{|R|\ge\eta^{-1}\}}.
\]

\paragraph{Proof of \eqref{us:eq:lowerK}.}
The compression of $RY$ to $E_>\cap F^\perp$ has row covariance at least
$(W\eta^2)^{-1}I$ and at least $C_0d-d$ rows. Factoring its covariance
reduces \eqref{us:eq:lowerK} to the rectangular Gaussian
least-singular-value bound.
For a Gaussian matrix with $n$ rows and $r$ columns, the rectangular
bound here only needs $n$ to exceed a sufficiently large constant times $r$.
For a fixed unit column vector, the squared Gaussian image norm is a
scaled sum of independent exponential variables, with lower and upper
tails $e^{-cn}$. A fixed fine net in the column sphere has $e^{O(r)}$
points. Its upper bound controls the operator norm, and then its lower
bound extends to the sphere. With $n\ge C r$ this proves the asserted
least singular value and exponential failure. A general row covariance
bounded below is factored out before this argument.

\paragraph{Proof of \eqref{us:eq:Gaussnorm}.}
The bounds in \eqref{us:eq:Gaussnorm} follow from a constant net: the expected norm in a fixed
column direction is at most $\norm{R_<}_{\HS}/\sqrt W$, and the Gaussian
Lipschitz deviation is $\norm{R_<}/\sqrt W$; use a deviation of order
$\sqrt d$. The unweighted projection bounds are the same calculation.

\paragraph{Proof of \eqref{us:eq:quadratic}.}
Set
\[
 M\coloneqq(I-P_F)R_<.
\]
The spectral assumptions give
$|\Tr M|\le CW\log W+d/\eta$, $\norm M_{\HS}\le C\sqrt{W/\eta}$,
and $\norm M\le\eta^{-1}$. The centered Gaussian quadratic-form
estimate and a constant net in the column space give
\[
 \norm{Y^*MY}\le |\Tr M|/W+
 C\sqrt d\norm M_{\HS}/W+Cd\norm M/W\le C\log W.
\]
For non-Hermitian $M$, apply this to its Hermitian real and imaginary
parts. The scalar tail estimate follows by diagonalization and the
exponential moment of centered independent exponential variables; the
net has $\exp(O(d))$ points. All failures are $Ce^{-cd}$.

\subsubsection{Spectral hypotheses for the complementary inverse}
\label{app:inverse-image-spectrum}

\paragraph{Proof of \eqref{us:eq:complement-spectrum}.}
Apply the averaged local law \eqref{us:eq:weighted-local-law} to the
complementary GUE with its $o(1)$-norm deformation and bulk energy. Set
$\eta=Ad/W$ and compare the imaginary parts of the normalized resolvent
traces at heights $c\eta$ and $\eta$. Outside $[-\eta,\eta]$, the
first Poisson kernel is at most $2c$ times the second. Choose $c$ small
to obtain at least $c'W\eta$ eigenvalues in this interval, then choose
$A$ large. The upper count follows from the kernel at $\eta$.
The upper count is $\Tr\mathbf1_{\{|\Gamma_c|\le r\}}\le CWr$
for $\eta\le r\le1$. Thus, for $j=1,2$,
\[
 \Tr|R_<|^j
 \le W+C\sum_{0\le l\le\lceil\log_2(1/\eta)\rceil}
       (2^l\eta)^{-j}W2^l\eta
 \le\begin{cases}CW\log W,&j=1,\\ CW/\eta,&j=2.\end{cases}
\]
The remaining norm bound $\|R_<\|\le\eta^{-1}$ is the definition of the truncation.
These are precisely the two truncated trace bounds in \eqref{us:eq:spectrum}. All these scales are in the local-law
domain, since $W\eta\ge AW^\delta$. No bound on the full real inverse
is needed in Lemma \ref{us:lem:self}.
For a complement of dimension $W-d$, its entry variance is $1/W$,
so divide the matrix, energy, and spectral parameter by $\sqrt{(W-d)/W}$
before applying Lemma~\ref{us:lem:local-law} in dimension $W-d$.
Since $d=o(W)$, this changes only fixed constants and keeps the energy
in a fixed bulk interval.

\subsection{Moment and lower bounds for Riccati products}
\label{app:product-moment-proof}
\label{su:sec:product-proof}

For fixed $b>0$ assume
\begin{equation}\label{us:eq:preliminary-range}
 W^{-b}\le g=o(1),\qquad |E|\le2-\kappa_0.
\end{equation}

\begin{lemma}[Killed product moment]\label{us:lem:product}
Assume \eqref{us:eq:preliminary-range}. Fix an integer $h\ge1$ and a
deterministic Hermitian $R_4$ satisfying \eqref{us:eq:regular}.
Let $V_5,\ldots,V_{h+4}$ be independent normalized complex GUE matrices,
and define
\[
 R_j\coloneqq(V_j-EI-g^2R_{j-1})^{-1},\qquad 5\le j\le h+4.
\]
Let $0<p<2/(h+1)$. For every $e>0$
there is a scalar spectral good event of probability $1-o(1)$ on which
\begin{equation}
 \Tr|R_{h+4}\cdots R_5|^p\le W^{1+e}.
 \label{us:eq:product}
\end{equation}
This event is invariant under simultaneous conjugation of
$R_4,V_5,\ldots,V_{h+4}$. The probability statement is uniform over
the deterministic $E,g,R_4$ just specified.
\end{lemma}

\begin{proof}[\textbf{Proof of \eqref{us:eq:product}}]
Choose a small $0<\delta<1$. Kill at any step where the inverse norm exceeds
$W^{1+\delta}$, the outgoing state is not regular, or an inverse dyadic
count exceeds $CW^{1+\delta}2^{-l}$. Equation \eqref{us:eq:wegner} and a
dyadic union bound make the total failure $O_h((\log W)W^{-\delta}+g^{1/2})$.
On regular pasts, Lemma~\ref{us:lem:local-law} gives, by testing every
fixed vector and then taking expectation,
\begin{equation}
 \E[\one_{\{|\Gamma|\le t\}}\mid\text{past}]
 \le CW^{2\delta}tI,
 \qquad W^{-1-\delta}\le t\le1.
 \label{us:eq:projexp}
\end{equation}
For $t\ge W^{-1+\delta}$ use
$\one_{|x|\le t}\le2t\Im(x-it)^{-1}$ and the local law; below that
scale use monotonicity at $W^{-1+\delta}$. The exceptional probability
can have an arbitrarily large power, and the projection is bounded by $I$.

For $5\le j\le h+4$, split a killed $R_j$ into bands $2^{l_j}A_{j,l_j}$, where
$\norm{A_{j,l_j}}\le1$ and the $l_j=0$ band contains the part of norm
at most one. Each killing indicator records survival through its
current step. In successive conditioning, earlier indicators are
measurable with respect to the past and the current indicator is
bounded by one; no future survival event is conditioned on.
For $Y\coloneqq A_{h+4,l_{h+4}}\cdots A_{5,l_5}$, including these
indicators, \eqref{us:eq:projexp} gives
\[
 \E\norm{Y}_{\HS}^2\le C_hW^{1+2h\delta}2^{-\sum l_j},\qquad
 \rank Y\le CW^{1+\delta}2^{-\max_j l_j}.
\]
The power-mean inequality for singular values and concavity imply
\[
 \E\Tr|2^{\sum l_j}Y|^p
 \le C_h W^{1+\delta(1-p/2+hp)}
       2^{-[1-p(h+1)/2]\max_j l_j}.
\]
The sum over multi-indices converges up to its harmless polynomial
multiplicity in $\max_j l_j$. To sum matrix bands, for $M_*$ matrices use the
rank approximation inequality
$N_{\sum Y_i}(t)\le\sum_iN_{Y_i}(t/M_*)$ and integrate it, obtaining
$\Tr|\sum Y_i|^p\le M_*^p\sum_i\Tr|Y_i|^p$.
Here $M_*=O((\log W)^h)$, absorbed by an arbitrarily small power.
Choose $\delta$ small and apply Markov with remaining power slack.
On the all-good event the killed product is the original product.
\end{proof}

For the application to Lemma~\ref{up:lem:small-contraction}(i), take
$h=q-6$ and $e=\varepsilon/3200$, with $q,p$ from \eqref{su:eq:packet2}.
One uniform choice of the small exponent in the proof is
\begin{equation}\label{su:eq:packet3}
 \delta\coloneqq
 \min\left\{\frac{\varepsilon}{204800(q-4)},\frac{1-2\beta}{64},
             \frac\alpha{64},\frac{\varepsilon}{51200}\right\}.
\end{equation}
Apply Lemma~\ref{us:lem:local-law} with domain exponent $\delta/4$
and stochastic-domination slack $\delta/8$. The scales $W^{-1+\delta}$
and $cg^2,g^2,Cg^2$ lie in this fixed domain. The geometric summation
margin is $1-p(q-5)/2=\varepsilon/6400$, and $1-p/2+(q-6)p<3$.
Thus the powers of $W^{\delta}$ and the logarithmic factors consume
less than half the available $\varepsilon/3200$ slack. The remaining Markov
failure is at most $W^{-\varepsilon/6400}$, while the killing failure is at most
$C_{q-5}((\log W)W^{-\delta}+W^{-\alpha/2})$. All parameters were fixed
before $a\in[\alpha,\beta]$ and $W$, which proves the asserted
uniformity.

\subsubsection{The output in the critical window}
\label{app:buffered-output}

\begin{lemma}[Buffered product lower bound]\label{nc:lem:buffer}
Under the assumptions of Lemma~\ref{up:lem:small-contraction}(ii),
use its $\theta,\kappa_0$ and define
\[
 \upsilon\coloneqq W^{-1+2\theta}.
\]
The matrices $V_5,\ldots,V_{q-1}$ are independent normalized complex GUE blocks.
Define
\[
 R_j\coloneqq (V_j-EI-g^2R_{j-1})^{-1},\qquad 5\le j\le q-1,
\]
with deterministic initial $R_4$ satisfying $\norm{R_4}\le W^{1+\theta}$ and regularity
\eqref{us:eq:regular}. Fix a nonzero vector $v_4$ and $|E|\le2-\kappa_0$.
For fixed $q$ and every fixed $D'>0$, the probability that the chain is
regular and \eqref{nc:eq:profile} fails for this $v_4$ is at most $W^{-D'}$.
\end{lemma}

\begin{proof}[\textbf{Proof of \eqref{nc:eq:profile}}]
Set $F_j^{(\upsilon)}=(I+\upsilon^2R_j^2)^{-1}$ for
$4\le j\le q-1$, with the empty-product convention $R_{[4,5]}=I_W$.
At each step $5\le j\le q-1$, a vector $u$ measurable before $V_j$ is fixed after
conditioning on the past. The exact spectral identity and
Lemma~\ref{us:lem:local-law} give
\begin{equation}
 u^*R_jF_j^{(\upsilon)}R_ju
 =\upsilon^{-1}\Im u^*(R_j^{-1}-i\upsilon)^{-1}u
 \ge c\upsilon^{-1}u^*(I+g^{4}R_{j-1}^2)^{-1}u
 \ge c\upsilon^{-1}u^*F_{j-1}^{(\upsilon)}u.
 \label{nc:eq:bufferstep}
\end{equation}
The first comparison in \eqref{nc:eq:bufferstep} has failure at most $W^{-D'}$ for any prescribed
$D'$, uniformly on regular pasts: $W\upsilon=W^{2\theta}$ and the
relative Dyson profile in \eqref{us:eq:mde} is retained. The second is a
deterministic comparison of two functions of the same Hermitian matrix,
using $g^2\le\upsilon$. Apply this with $u=R_{[j-1,5]}v_4$, normalized before
each local-law application. A finite conditional union bound, intersected
with chain regularity, gives the iterated lower bound. Finally
$\upsilon\norm{R_4}\le W^{3\theta}$, so
$v_4^*F_4^{(\upsilon)}v_4\ge\tfrac12W^{-6\theta}\|v_4\|^2$, and
$F_{q-1}^{(\upsilon)}\le I$. This proves \eqref{nc:eq:profile}.
\end{proof}

\subsection{Coercivity of a structured map of GUE matrices}
\label{app:compact-gaussian}

We adapt the proof of \cite[Lemma~A.2]{Drogin} to complex GUE matrices.
The complex map has the additional kernel direction \(iI_W\), which
is removed by the restriction \(\Im\Tr P=0\) below. After this
adjustment, we follow the same decoupling, norm-convergence and
Gaussian concentration argument.

\begin{lemma}[Structured Gaussian coercivity]\label{lem:c-gaussian}
Let \(m\coloneqq 128\) and let \(\widehat V=(\widehat V_1,\ldots,\widehat V_m)\) consist of independent complex GUE
matrices with
\(\E(\widehat V_j)_{ab}(\widehat V_k)_{cd}
  =W^{-1}\delta_{jk}\delta_{ad}\delta_{bc}\).
On the real Hilbert space
\[
 \mathcal Q_W\coloneqq \{P\in\mathbb C^{W\times W}:\Im\Tr P=0\}
\]
define
\begin{equation}
 \mathcal L_{\widehat V}(P)\coloneqq (\widehat V_jP+P^*\widehat V_j)_{j=1}^{m}.
 \label{cg:eq:LG}
\end{equation}
There exist absolute \(\tau,c,C>0\) such that, for sufficiently large \(W\),
\begin{equation}
 \Prob\{s_{\min}(\mathcal L_{\widehat V}|_{\mathcal Q_W})<\tau\}\le Ce^{-cW}.
 \label{cg:eq:SG}
\end{equation}
For any fixed finite additional family of normalized GUE matrices
(possibly including matrices from which the \(\widehat V_j\) are formed as disjoint normalized
differences), a sufficiently large fixed operator-norm cutoff has the
same exponentially small failure bound. There is a \(C^\infty\) cutoff
\(0\le\chi_P\le1\), compactly supported in the input Hermitian variables,
equal to one with probability \(1-Ce^{-cW}\), supported where the
singular value and norm conditions hold with fixed margins, and
\begin{equation}
 \|\nabla\chi_P\|\le C .
 \label{cg:eq:cutoff}
\end{equation}
Here the gradient is taken in all input Hermitian variables with
their real Hilbert--Schmidt product inner product.
The number of input matrices is fixed before \(W\) varies.

\end{lemma}

\begin{proof}[\textbf{Proof of \eqref{eq:c-main-cutoff} and \eqref{eq:c-main-coercivity}}]
We follow \cite[Appendix~A.1, proof of Lemma~A.2]{Drogin}, recording
the complex-valued adjustment and the cutoff construction.
Put
\[
 S\coloneqq\sum_j\widehat V_j^2,\qquad
 \Phi(P)\coloneqq\sum_j\widehat V_jP^*\widehat V_j,\qquad
 S_W\coloneqq s_{\min}(\mathcal L_{\widehat V}|_{\mathcal Q_W}).
\]
The unrestricted map annihilates \(iI_W\). On the full real matrix
space, \(\mathcal L_{\widehat V}^*\mathcal L_{\widehat V}(P)=2SP+2\Phi(P)\), and
GUE contraction gives
\[
 \E\Phi(P)=W^{-1} m\overline{\Tr P}\,I_W,\qquad
 \langle P,\E\Phi(P)\rangle
 =W^{-1} m\Re(\Tr P)^2\ge0\quad(P\in\mathcal Q_W).
\]
Thus the restriction \(\Im\Tr P=0\) removes the kernel direction
and gives the analogue of \cite[(42)]{Drogin}:
\begin{equation}
 S_W^2\ge2\lambda_{\min}(S)-2\|\Phi-\E\Phi\|.
 \label{cg:eq:coercive}
\end{equation}

Let \(\widehat V'\) be an independent copy of \(\widehat V\). Jensen's inequality
bounds \(\E\|\Phi_{\widehat V}-\E\Phi_{\widehat V}\|\) by \(\E\|\Phi_{\widehat V}-\Phi_{\widehat V'}\|\).
Under the Gaussian rotation
\(\widehat V=(B+C)/\sqrt2\), \(\widehat V'=(B-C)/\sqrt2\),
\[
 (\Phi_{\widehat V}-\Phi_{\widehat V'})(P)=\sum_j(B_jP^*C_j+C_jP^*B_j).
\]
The two cross operators have the same distribution. The triangle
inequality and the real isometry \(P\mapsto P^*\) give, as in
\cite[Appendix~A.1]{Drogin},
\begin{equation}
 \E\|\Phi-\E\Phi\|
 \le2\E\left\|\sum_jB_j\otimes C_j^t\right\|.
 \label{cg:eq:decouple}
\end{equation}
Use \cite[Corollary~1.2 and Lemma~9.6]{CGV} for ordinary GUE
norm convergence and for the two independent tensor factors,
respectively. The normalization is exactly that of
\cite[Definition~2.8(a)]{CGV}; transpose preserves the GUE law.
For variance-one semicircular families \(s_j,t_j\) on the two
factors, the free bounds used in \cite[Appendix~A.1]{Drogin} imply
\begin{equation}
 \sum_js_j^2\ge(m-2\sqrt m)I,\qquad
 \left\|\sum_js_j\otimes t_j\right\|\le2(\sqrt m+1).
 \label{cg:eq:freebounds}
\end{equation}
The lower spectral bound follows from norm convergence applied to
\(cI-S\) for a sufficiently large fixed \(c\). Fixed GUE norm
moments give uniform integrability in both limits. Consequently,
\eqref{cg:eq:coercive}--\eqref{cg:eq:freebounds} yield
\begin{equation}
 \liminf_{W\to\infty}\E S_W^2
 \ge2(m-2\sqrt m)-8(\sqrt m+1)>0\quad(m=128).
 \label{cg:eq:positive}
\end{equation}
Also \(\E S_W^4\le C\), since
\(S_W\le2(\sum_j\|\widehat V_j\|^2)^{1/2}\). Paley--Zygmund gives
\(\E S_W\ge c_0>0\) for large \(W\). The map \(\widehat V\mapsto S_W\)
is \(2\)-Lipschitz in the product Hilbert--Schmidt norm, so Gaussian
concentration with coordinate covariance \(I/W\) proves
\eqref{cg:eq:SG}. Standard GUE operator-norm tails give the same
failure bound for the fixed additional family.

For the cutoff, disjoint normalized differences are an orthogonal
Gaussian transformation, so the same Lipschitz bound holds in the
original inputs. Multiply Lipschitz cutoffs of \(S_W\) and the
input norms, with fixed nested margins, and mollify in the real
Hermitian input space inside a sufficiently small fixed
Hilbert--Schmidt ball. This preserves the support margins and the
uniform gradient bound, and leaves \(\chi_P=1\) on a stronger event
of probability \(1-Ce^{-cW}\). The norm cutoff makes its support
compact. This proves \eqref{cg:eq:cutoff} and the lemma.
\end{proof}

\subsection{Construction of the common density}
\label{sec:c-component-proof}

Lemma~\ref{up:lem:good-bad} asks for a smooth good component of the
batch law, the polynomial estimates \eqref{up:eq:coarse-density},
and the two operator bounds \eqref{up:eq:GB-bounds}. We first split
the law by a cutoff in the GUE inputs. To show that the retained
measure has a smooth density, we prove that the product map has
surjective differential on the cutoff support. A controlled right
inverse of this differential then gives an auxiliary density-derivative bound
by Gaussian integration by parts. The operator bounds follow from
Lemma~\ref{lem:c-operators}.
This construction follows \cite[Proposition~4.1, Lemma~4.2 and
Section~4.3]{Drogin}; the GUE coercivity input is
Lemma~\ref{lem:c-gaussian} in Appendix~\ref{app:compact-gaussian}.

\begin{proof}[\textbf{Proof of \eqref{cc:eq:decomp}, \eqref{up:eq:coarse-density}, and \eqref{up:eq:GB-bounds}}]
Fix $B,E_0$ and deterministic $W^{-B}\le g\le W^B$, $|E|\le E_0$.
Let $q$ be the fixed multiple of 1000 in the lemma. In particular,
$q$ is even and independent of $W$; all constants below may depend
on $B,E_0,q$.

\medskip\noindent\emph{Step 1: split the batch law and control the discarded part.}
The goal of this step is to obtain \eqref{up:eq:coarse-density}
from the cutoff conditions \eqref{eq:c-main-cutoff}.
Use $T_i$ from \eqref{eq:c-transfer} and the fixed unitary matrix
$S$ from Definition~\ref{def:transfer-matrices}. Set
$F\coloneqq T_1\cdots T_q$. By reversal of the independent inputs,
$S^*FS$ has the batch law $\mu$. Use $C_i$ from \eqref{eq:c-transfer-block}
and put $h_g:=\max\{1,g^{-1}\}$.
For a cutoff $0\le\chi_P\le1$ to be chosen, put
$\widetilde\pi\coloneqq\E\chi_P$ and define the positive measure
$\mu_{\rm bad}$ by
\[
 \int u(\mathcal T)\,d\mu_{\rm bad}(\mathcal T)
 \coloneqq\E[(1-\chi_P)u(S^*FS)]
\]
for bounded measurable $u$. Then
$\mu-\mu_{\rm bad}$ is positive, has mass $\widetilde\pi$, and
$\mu_{\rm bad}(\mathcal G)=1-\widetilde\pi$.
We will prove in Step 2 that $\mu-\mu_{\rm bad}$ has a smooth density
$\widetilde\pi f$.

For the mass and support estimates, we require
\[
 \begin{gathered}
 \chi_P\in C_c^\infty,\qquad 0\le\chi_P\le1,\qquad
 \Prob(\chi_P=1)\ge1-Ce^{-cW},\\
 \supp\chi_P\subset\{\max_i\|V_i\|\le C\}.
 \end{gathered}
\]
Step 2 will construct such a cutoff. We first check the consequences
of these requirements. The probability bound gives $\widetilde\pi\ge1/2$;
its stronger exponential form will also control the weighted bad mass.

On $\supp\chi_P$, each $\|C_i\|\le Ch_g$, and hence
$\|F\|=\|F^{-1}\|\le Ch_g^q\le W^C$. Fixed GUE norm moments give
\[
 \E\|F\|\le Ch_g^q.
\]
Consequently,
\[
 \int\|\mathcal T\|^{1/2}\,d\mu_{\rm bad}(\mathcal T)
 \le\Prob(\chi_P\ne1)^{1/2}(\E\|F\|)^{1/2}
 \le Ce^{-cW}h_g^{q/2}\le\tfrac12
\]
for large $W$. Thus all estimates in \eqref{up:eq:coarse-density}
are available once the density is shown to exist.

\medskip\noindent\emph{Step 2: construct the cutoff and obtain a smooth density.}
The goal of this step is to obtain \eqref{cc:eq:decomp} by proving
\eqref{eq:c-main-left-inverse}.
\smallskip\noindent\textbf{Proof of \eqref{eq:c-main-left-inverse}.}
We choose $\chi_P$ so that $V\mapsto F(V)$ has surjective differential
on a neighborhood of $\supp\chi_P$; this will give the required smooth
density. Conjugation by $S$ preserves this property.

For a Hermitian input variation $H=(H_i)$, differentiate the product
as in \cite[Claim~4.5]{Drogin}:
\begin{equation}\label{eq:c-differential}
 \mathcal A(H)\coloneqq F^{-1}(H\cdot\nabla_V)F
 =g^{-1}\sum_{i=1}^qQ_i^{-1}
       \begin{pmatrix}0&0\\H_i&0\end{pmatrix}Q_i,
 \qquad Q_i\coloneqq T_{i+1}\cdots T_q,
\end{equation}
with $Q_q=I$. We must show that $\mathcal A$ is onto.
For a target $Z$, use the block form \eqref{eq:c-Lie}, with
$(A,B_1,B_2)$ denoted by $(P,U,R)$.
We do so by recovering every $Z$ from $\mathcal A^*Z$; the adjoint of this
reconstruction will also provide the input perturbation needed in
Step 3.

Write $Z_i\coloneqq Q_i^{-*}ZQ_i^*$ and denote its blocks as in
\eqref{eq:c-Lie} by $(P_i,U_i,R_i)$, so
$(P_q,U_q,R_q)=(P,U,R)$. The trace pairing and transfer recursion,
as in \cite[Proposition~4.4, pp.~34--35]{Drogin}, give
\begin{equation}\label{eq:c-observe}
 \mathcal A^*Z=g^{-1}(R_i)_{i=1}^q,\qquad
 \begin{aligned}
 P_{i-1}&=-R_iC_i-P_i^*,\qquad
 U_{i-1}=-R_i,\\
 R_{i-1}&=P_iC_i+C_iP_i^*+C_iR_iC_i-U_i,
 \end{aligned}
 \quad 2\le i\le q.
\end{equation}
Thus the data are $\mathbf R=(R_1,\ldots,R_q)$, and we need to recover
$P,U,R$. For any Hermitian data array define
$E_q\coloneqq0$ and $E_{i-1}\coloneqq-R_iC_i-E_i^*$ for $2\le i\le q$.
When the data come from $Z$ as above, the first recursion gives
\[
 P_i=\begin{cases}P+E_i,&i\text{ even},\\-P^*+E_i,&i\text{ odd}.
 \end{cases}
\]
At an interior odd index, substituting this into the remaining
recursions yields
\begin{equation}\label{eq:c-pair}
 C_iP+P^*C_i=F_i(\mathbf R),\qquad
 F_i(\mathbf R)\coloneqq E_iC_i+C_iE_i^*+C_iR_iC_i+R_{i+1}-R_{i-1}.
\end{equation}
To recover $P$ from these equations, subtract them at disjoint pairs
of interior odd indices. Set
\[
 (i_j,k_j)\coloneqq(4j+1,4j-1),\qquad 1\le j\le128,\qquad
 \widehat V_j\coloneqq2^{-1/2} g(C_{i_j}-C_{k_j})
       =2^{-1/2}\bigl(V_{k_j}-V_{i_j}\bigr).
\]
The largest selected index is $513$, and its neighbor $514$ is still
within the batch since $q\ge1000$. The common energy term cancels,
and the disjoint pairs make the $\widehat V_j$ independent normalized GUE matrices.
Since $\Im\Tr P=0$, use $\mathcal L=\mathcal L_{\widehat V}$ from
\eqref{cg:eq:LG} for these $\widehat V_j$.
With
\[
 \mathcal F\mathbf R\coloneqq
       (F_{i_j}(\mathbf R)-F_{k_j}(\mathbf R))_{j=1}^{128},
\]
we obtain $\mathcal LP=(g/\sqrt2)\mathcal F\mathbf R$.
Only these odd-index pairs are needed to recover $P$; the other
blocks remain in the product and in the data $\mathbf R$.
Each $F_i$ is linear in $\mathbf R$ and has degree at most two in $C$.

Thus it suffices to arrange $s_{\min}(\mathcal L)\ge c$ on the
cutoff support. Lemma~\ref{lem:c-gaussian}, applied also to the norms
of all $q$ input blocks, now supplies a cutoff satisfying the
requirements of Step 1 and, in addition,
\[
 \supp\chi_P\subset\{s_{\min}(\mathcal L)\ge c\},\qquad
 \|\nabla_V\chi_P\|\le C.
\]
The lower singular-value bound and the input norm bounds hold on a
neighborhood of $\supp\chi_P$, after adjusting the fixed constants.
This choice of $\chi_P$ depends on the input blocks.
On this neighborhood, $\mathcal L$ has a bounded left inverse.
Define, on the cutoff neighborhood,
\[
 \mathcal J\coloneqq(\mathcal L^*\mathcal L)^{-1}\mathcal L^*,
 \qquad \mathcal K_P\coloneqq2^{-1/2} g\mathcal J\mathcal F.
\]
Then $P=\mathcal K_P\mathbf R$. The terminal recursion recovers
$U$ from $P,\mathbf R$, while $R=R_q$ is already observed. Accordingly,
for arbitrary Hermitian data define
\begin{equation}\label{eq:c-reconstruction}
 \mathcal K_0\mathbf R\coloneqq
 \begin{pmatrix}P&U\\R_q&-P^*\end{pmatrix},\qquad
 P\coloneqq\mathcal K_P\mathbf R,\qquad
 U\coloneqq PC_q+C_qP^*+C_qR_qC_q-R_{q-1}.
\end{equation}
This matrix belongs to the Lie algebra. For the data associated with
$Z$ we have just proved $\mathcal K_0\mathbf R=Z$.
Since $\mathcal A^*Z=g^{-1}\mathbf R$, taking adjoints gives
\begin{equation}\label{eq:c-right-inverse}
 g\mathcal K_0\mathcal A^*=I,\qquad
 \mathcal A(g\mathcal K_0^*)=I.
\end{equation}
In particular $\mathcal A$ is onto on the cutoff neighborhood.

Cover the compact cutoff support by finitely many submersion
charts and choose a partition of unity. In each chart, integration
of the smooth compactly supported input density along the fibers
gives a smooth target density. Their sum is a smooth compactly
supported Haar density for $\mu-\mu_{\rm bad}$; no global
submersion or properness assumption is needed. Its normalization defines
$f$ by
\begin{equation}\label{eq:c-cutoff-pushforward}
 \begin{aligned}
 \widetilde\pi\int u(\mathcal T)f(\mathcal T)\,d\mathcal T&\coloneqq\E[\chi_P\,u(S^*FS)],\\
 \int u(\mathcal T)\,d\mu_{\rm bad}(\mathcal T)&=\E[(1-\chi_P)u(S^*FS)].
 \end{aligned}
\end{equation}
Therefore $f\ge0$, $\int f=1$, $f\in C_c^\infty(\mathcal G)$,
and $\mu=\widetilde\pi f\,d\mathcal T+\mu_{\rm bad}$, proving
\eqref{cc:eq:decomp}. This density and cutoff were chosen before any
flag or differentiation direction.

\medskip\noindent\emph{Step 3: bound the directional derivatives of this density.}
The goal of this step is to obtain the following auxiliary bound for
the smoothing estimate in Step 4:
\begin{equation}\label{eq:c-coarse-derivative}
 \sup_{\substack{M\in\mathfrak g\\\|M\|_{\rm HS}=1}}
       \|X_Mf\|_1\le W^C.
\end{equation}
The directions and derivatives are defined in
\eqref{up:eq:lie-algebra}--\eqref{eq:c-right-derivative}.
For $M$ with $S^*MS\in\mathfrak g$ and $\|M\|_{\HS}=1$,
take the field
\[
 X\coloneqq\sqrt2\,g\mathcal K_0^*M,\qquad
 \mathcal A(X)=\sqrt2\,M,\qquad D_XF=\sqrt2\,FM.
\]
The last two identities follow from \eqref{eq:c-right-inverse}.
Extend $\chi_P X$ by zero and let $\gamma$ be the Gaussian law of the
original inputs $V_i$, with covariance $W^{-1}I$ in real
trace-orthonormal Hermitian coordinates. The definition
\eqref{eq:c-cutoff-pushforward}, the chain rule, and Gaussian integration
by parts \cite[Lemma~4.2]{Drogin} give, for $u\in C_c^\infty(\mathcal G)$,
\[
 \begin{aligned}
 \sqrt2\,\widetilde\pi\int (X_{S^*MS}u)(\mathcal T)f(\mathcal T)\,d\mathcal T
 &=\E[\chi_P D_X[u(S^*FS)]]\\
 &=\E[u(S^*FS)\,\delta_\gamma(\chi_P X)].
 \end{aligned}
\]
Here $D_X=X\cdot\nabla_V$, and
$\E D_X\varphi=\E[\varphi\delta_\gamma X]$ defines the Gaussian divergence. For completeness, we give the second-moment calculation.
\smallskip\noindent\textbf{Proof of \eqref{up:eq:lift-divergence-reduction}.}
In real trace-orthonormal coordinates $v_i$ of variance $1/W$, put
$Y=\chi_P X$ and $\delta_\gamma Y=\sum_i(Wv_iY_i-\partial_iY_i)$.
Two Gaussian integrations by parts give
\[
 \E|\delta_\gamma Y|^2
 =W\E\sum_iY_i^2+\E\sum_{i,j}(\partial_iY_j)(\partial_jY_i)
 \le W\E\|Y\|^2+\E\|\nabla Y\|_{\rm HS}^2.
\]
The field is real on the Hermitian input space, and $Y$ has compact
support, so no boundary term occurs. This proves
\eqref{up:eq:lift-divergence-reduction}.\smallskip

The reconstruction formulas use a fixed number of matrix products
and inverses, with $\|C_i\|\le Ch_g$ and $s_{\min}(\mathcal L)\ge c$
on the cutoff neighborhood. Differentiating them, and summing over the
$qW^2$ input coordinates, therefore gives
\[
 \sup_{\|M\|_{\rm HS}=1}
 (\|X\|+\|\nabla_VX\|_{\rm HS})+\|\nabla_V\chi_P\|\le W^C.
\]
Haar integration by parts and duality over $\|u\|_\infty\le1$ now yield
\begin{equation}\label{eq:c-score-to-density}
 \|X_{S^*MS}f\|_1
 \le2^{-1/2}\widetilde\pi^{-1}\|\delta_\gamma(\chi_P X)\|_2\le W^C.
\end{equation}
Since $S$ is unitary and $\widetilde\pi\ge1/2$, this proves
\eqref{eq:c-coarse-derivative} uniformly over unit Lie directions.

\medskip\noindent\emph{Step 4: deduce the two operator bounds.}
The goal of this step is to obtain \eqref{up:eq:GB-bounds}.
The support bound in \eqref{up:eq:coarse-density} and the derivative
bound \eqref{eq:c-coarse-derivative}, both proved above, are the
inputs to the smoothing estimate \eqref{cc:eq:smoothing} in
Lemma~\ref{lem:c-operators}; it gives
\[
 \|G_{1/2}h\|_\infty
 \le W^C\sup_U\|h(U,\cdot)\|_{L^4(\mathcal C_U,dm_U)}.
\]
For the bad operator, $e^{-c(\mathcal T;U,Y)/2}\le\|\mathcal T\|^{1/2}$ and Step 1 give
\[
 \|B_{1/2}h\|_\infty
 \le\|h\|_\infty\int\|\mathcal T\|^{1/2}\,d\mu_{\rm bad}(\mathcal T)
 \le\tfrac12\|h\|_\infty.
\]
The decomposition of $\mu$ gives $P_{1/2}=G_{1/2}+B_{1/2}$, so
\eqref{up:eq:GB-bounds} follows.

All constants and width thresholds depend only on $B,E_0,q$,
as asserted in Lemma~\ref{up:lem:good-bad}.
\end{proof}

\phantomsection

\end{document}